\documentclass[11pt]{amsart}
\usepackage[colorlinks,plainpages,backref=page,urlcolor=blue]{hyperref}

\usepackage{amsmath,amsthm,verbatim}
\usepackage{amscd,mathtools}
\usepackage{txfonts,amssymb}
\usepackage{tikz}
\usetikzlibrary{cd,calc}
\usetikzlibrary{matrix,arrows,decorations.pathmorphing,calc,shapes}
\usepackage{graphicx,calrsfs}
\usepackage{bbm}
\usepackage{mathrsfs}
\usepackage{xstring}
\usepackage{enumitem}
\usepackage[all]{xy}
\usepackage{spectralsequences}
\usepackage{bm}
\usepackage{accents}
\usepackage[final]{microtype}

\newtheorem{theorem}{Theorem}[section]
\newtheorem{corollary}[theorem]{Corollary}

\newtheorem{lemma}[theorem]{Lemma}

\newtheorem{proposition}[theorem]{Proposition}

\newtheorem{thm}{Theorem}

\theoremstyle{definition}
\newtheorem{definition}[theorem]{Definition}
\newtheorem{example}[theorem]{Example}
\newtheorem{remark}[theorem]{Remark}
\newtheorem{conjecture}[theorem]{Conjecture}
\newtheorem{question}[theorem]{Question}
\newtheorem{problem}[theorem]{Problem}

\DeclareMathOperator{\im}{im}
\DeclareMathOperator{\ini}{in}
\DeclareMathOperator{\id}{id}

\DeclareMathOperator{\coker}{coker}
\DeclareMathOperator{\Supp}{Supp}
\DeclareMathOperator{\supp}{Supp}
\DeclareMathOperator{\Ann}{Ann}
\DeclareMathOperator{\Iso}{Iso}
\DeclareMathOperator{\ord}{ord}
\DeclareMathOperator{\Syz}{Syz}

\DeclareMathOperator{\annn}{ann}
\newcommand{\wtH}{\widetilde{H}}
\DeclareMathOperator{\Isom}{Isom}
\DeclareMathOperator{\rank}{rank}

\DeclareMathOperator{\Tor}{Tor}
\DeclareMathOperator{\ab}{ab}

\DeclareMathOperator{\Prim}{Prim}
\DeclareMathOperator{\gr}{gr}
\DeclareMathOperator{\Sym}{Sym}
\DeclareMathOperator{\Hilb}{Hilb} 
\DeclareMathOperator{\TC}{TC}

\DeclareMathOperator{\Hom}{Hom}
\DeclareMathOperator{\Aut}{Aut}
\DeclareMathOperator{\End}{End}
\DeclareMathOperator{\OS}{OS}
\DeclareMathOperator{\ch}{char}

\DeclareMathOperator{\IA}{IA}
\DeclareMathOperator{\ev}{ev}
\DeclareMathOperator{\Ev}{Ev}

\DeclareMathOperator{\Spec}{Spec}
\DeclareMathOperator{\mSpec}{mSpec}
\DeclareMathOperator{\spn}{span}
\DeclareMathOperator{\Conf}{Conf}
\DeclareMathOperator{\Fitt}{Fitt}
\DeclareMathOperator{\wt}{wt}
\DeclareMathOperator{\CE}{CE}

\DeclareMathOperator{\In}{In}
\DeclareMathOperator{\Tors}{Tors}
\DeclareMathOperator{\GL}{GL}

\DeclareMathOperator{\PD}{PD}

\newcommand{\dga}{\ensuremath{\textsc{dga}}}
\newcommand{\cdga}{\ensuremath{\textsc{cdga}}}
\newcommand{\cga}{\ensuremath{\textsc{cga}}}
\newcommand{\cgas}{\cga{}\textnormal{s}}
\newcommand{\cdgas}{\cdga{}\textnormal{s}}
\newcommand{\N}{\mathbb{N}}
\newcommand{\R}{\mathbb{R}}
\newcommand{\Q}{\mathbb{Q}}
\newcommand{\C}{\mathbb{C}}
\newcommand{\Z}{\mathbb{Z}}
\newcommand{\E}{\mathbb{E}}

\newcommand{\K}{\mathbb{K}}
\renewcommand{\P}{\mathbb{P}}
\newcommand{\CP}{\mathbb{CP}}
\newcommand{\RP}{\mathbb{RP}}
\renewcommand{\k}{\Bbbk}

\newcommand{\bL}{\mathbf{L}}

\newcommand{\bV}{\mathbf{V}}

\newcommand{\g}{\mathfrak{g}}
\newcommand{\h}{\mathfrak{h}}
\newcommand{\fA}{\mathfrak{A}}
\newcommand{\fa}{\mathfrak{a}}
\newcommand{\fb}{\mathfrak{b}}
\newcommand{\B}{\mathfrak{B}}

\newcommand{\fM}{\mathfrak{M}}
\newcommand{\m}{\mathfrak{m}}
\newcommand{\n}{\mathfrak{n}}
\newcommand{\fL}{\mathfrak{L}}
\newcommand{\fr}{\mathfrak{r}}

\newcommand{\p}{\mathfrak{p}}
\newcommand{\ff}{\mathfrak{f}}
\newcommand{\Lie}{\mathbb{L}}
\newcommand{\fw}{\mathfrak{w}}

\newcommand{\fE}{\mathfrak{E}}
\newcommand{\sol}{\mathfrak{sol}}

\newcommand{\A}{\mathcal{A}}
\newcommand{\RR}{\mathcal{R}}
\newcommand{\M}{\mathcal{M}}

\newcommand{\cM}{\mathcal{M}}

\newcommand{\cF}{\mathcal{F}}

\newcommand{\cE}{\mathcal{E}}

\newcommand{\VV}{\mathcal{V}}

\DeclareMathAlphabet{\pazocal}{OMS}{zplm}{m}{n}

\newcommand{\dR}{\scriptscriptstyle{\rm dR}}
\newcommand{\PL}{\scriptscriptstyle{\rm PL}}
\newcommand{\apl}{A_{\scriptscriptstyle{{\rm PL}}}}

\newcommand{\oE}{\overline{E}}
\newcommand{\oU}{\overline{U}}
\newcommand{\oA}{\bar{A}}
\newcommand{\oS}{\bar{S}}

\newcommand{\wC}{\,\widehat{\CE}}

\newcommand{\DS}{\displaystyle}

\newcommand{\abs}[1]{\left| #1 \right|}

\newcommand{\surj}{\twoheadrightarrow}
\newcommand{\inj}{\hookrightarrow}
\newcommand\isom{\xrightarrow{ \,\smash{\raisebox{-0.6ex}{\ensuremath{\scriptstyle\simeq}}}\,}}
\newcommand{\longisom}{\xrightarrow{\hspace{4pt}\smash{\raisebox{-0.6ex}{\ensuremath{\scriptstyle\simeq}}}\hspace{4pt}}}
\newcommand{\longsurj}{\relbar\joinrel\twoheadrightarrow}
\newcommand{\longinj}{\lhook\joinrel\longrightarrow}

\newcommand{\bwedge}{\mbox{$\bigwedge$}}
\newcommand{\bcup}{\mbox{$\bigcup$}}
\newcommand{\boplus}{\mbox{$\bigoplus$}}
\newcommand{\sbm}[1]{{\let\amp=&\left(\begin{smallmatrix}#1\end{smallmatrix}\right)}}
\newcommand{\qA}{\mathrm{q}A}

\newcommand{\ql}{\mathrm{ql}}
\newcommand{\qql}{\mathrm{qql}}

\newcommand{\arxiv}[1]
{\texttt{\href{http://arxiv.org/abs/#1}{arXiv:#1}}}

\newcommand{\temp}{}
\providecommand{\MR}[1]{}
\renewcommand{\MR}[1]{%
\StrBefore{#1 }{ }[\temp]%
\href{https://mathscinet.ams.org/mathscinet-getitem?mr=\temp}{MR#1}}

\definecolor{lime}{HTML}{A6CE39}
\DeclareRobustCommand{\orcidicon}{
	\begin{tikzpicture}
	\draw[lime, fill=lime] (0,0) 
	circle [radius=0.16] 
	node[white] {{\fontfamily{qag}\selectfont \tiny ID}};
	\draw[white, fill=white] (-0.0625,0.095) 
	circle [radius=0.007];
	\end{tikzpicture}
	\hspace{-2mm}
}
\foreach \x in {A, ..., Z}{\expandafter\xdef\csname orcid\x\endcsname{\noexpand\href{%
https://orcid.org/\csname orcidauthor\x\endcsname}{\noexpand\orcidicon}}}

\begin{document}

%\date{August 10, 2026}

\title[Koszul modules, holonomy, and resonance]%
{Koszul modules, holonomy Lie algebras, and resonance of groups and CDGAs}

\author[Alexander~I.~Suciu]{Alexander~I.~Suciu$^1$\!\!\orcidA{}}
\address{Department of Mathematics, Northeastern University,
Boston, MA 02115, USA}
\email{\href{mailto:a.suciu@northeastern.edu}{a.suciu@northeastern.edu}}
\urladdr{\href{https://suciu.sites.northeastern.edu}%
{suciu.sites.northeastern.edu}}

\thanks{$^1$Partially supported by the Simons
Foundation Collaboration Grant for Mathematicians \#693825 and
by the project ``Singularities and Applications'' - CF 132/31.07.2023 funded by
the European Union - NextGenerationEU - through Romania's National Recovery
and Resilience Plan.}



\begin{abstract}
We develop a Koszul-theoretic framework for comparing classical 
Alexander-type invariants with infinitesimal invariants arising from 
finite-type commutative differential graded algebra models.
The central mechanism is Koszul linearization, which replaces 
nonlinear equivariant constructions with functorial algebraic 
objects defined from a $\cdga$.

To a connected $\cdga$ $(A,d)$ with finite-dimensional $H^1(A)$ 
we associate Koszul modules $\B_i(A)$ over the symmetric algebra 
on $H_1(A)$. We prove that the first Koszul module $\B_1(A)$ is 
isomorphic to the infinitesimal Alexander invariant of the holonomy 
Lie algebra $\h(A)$, yielding explicit formulas for holonomy Chen ranks.

We establish a tangent cone theorem for resonance varieties, showing 
that cohomology controls their first-order behavior at the origin.
For finitely generated groups admitting $1$-finite $1$-models, 
we prove that classical Alexander invariants agree with Koszul 
invariants after completion and passage to associated graded 
objects, and that Chen ranks are determined by the model. 
Applications include Chen rank computations for $2$-step 
nilpotent Lie algebras and pure elliptic braid groups, partial 
formality results for Sasakian manifolds, and a general framework 
for detecting non-formality of spaces, groups, and maps.
\end{abstract}

\maketitle

\setcounter{tocdepth}{1}
\tableofcontents

\setcounter{section}{0}
\renewcommand{\thesection}{\roman{section}}

%%%%%%%%%%%%%%
\part*{Introduction}
%%%%%%%%%%%%%%

%%%%%%%%%%%%%%%%%%%%%%
\section{Basic objects and scope}
\label{sect:basic-scope}
%%%%%%%%%%%%%%%%%%%%%%%%

We study spaces, groups, and maps via finite-type commutative differential 
graded algebra models, with the goal of comparing invariants defined purely 
from cohomology with finer, infinitesimal invariants extracted from such models.
The guiding theme is that cohomology captures only first-order information 
at the origin of natural representation spaces, while higher-order structure 
persists in finite-type models. In particular, discrepancies 
between cohomological invariants and their model-theoretic 
analogues provide concrete and computable obstructions to formality.

The central mechanism is Koszul linearization, which replaces 
nonlinear or equivariant construc\-tions---such as universal abelian covers 
and Alexander invariants---with functorial algebraic objects associated to 
a $\cdga$. This procedure yields computable invariants 
that interpolate between cohomology and minimal models and serve as 
infinitesimal avatars of classical topological constructions.

Concretely, let $(A,d)$ be a commutative differential graded algebra over a
field $\k$ of characteristic $0$.  In most results we assume $A$ is connected 
($A^0 = \k$) with $\dim_\k H^1(A) < \infty$.  To such a $\cdga$ we associate 
Koszul chain and cochain complexes whose homology and cohomology define 
a family of graded and filtered modules. 

A recurring theme is the comparison between invariants computed from the 
cohomology algebra $H^{*}(A)$ and those computed from a finite-type 
model $A$. In formal situations, these invariants necessarily coincide; 
discrepancies therefore quantify the failure of formality and provide 
refined information invisible to cohomology alone.  The examples worked out 
below are nilmanifolds, configuration spaces of elliptic curves, Sasakian 
manifolds, and groups admitting $1$-finite $1$-models.

Resonance varieties provide the most visible instance of the phenomenon 
described above. For a space or $\cdga$ $A$, one may compare the classical 
resonance varieties defined from $H^{*}(A)$ with the resonance varieties 
defined from a finite-type model. In many non-formal examples---most notably 
nilmanifolds---these two invariants differ dramatically. Such discrepancies 
cannot occur in the formal category and therefore provide effective obstructions 
to formality. More generally, the tangent-cone theorem proved here shows that 
cohomology governs only the first-order behavior of resonance at the origin, 
while higher-order structure is detected by the model.

A conceptual novelty of this approach is that it applies not only to
spaces and groups, but also to maps. Functoriality of Koszul modules 
detects non-formality of morphisms that act trivially on cohomology, 
a phenomenon invisible to classical cohomological methods. Resonance 
varieties, holonomy Lie algebras, Alexander invariants, and Chen ranks 
all arise naturally from this single Koszul-theoretic construction.

What these constructions have in common is that they are functorial, they 
require only a finite-type model rather than a formal one, and they convert 
discrepancies between a model and its cohomology into computable algebraic 
data.

%%%%%%%%%%%%%%%%%%%%%%
\section{Models and algebraic invariants}
\label{sect:models-alginv}
%%%%%%%%%%%%%%%%%%%%%%%%

We work systematically with finite-type $\cdga$ models of spaces and groups, 
together with their associated Koszul complexes, holonomy Lie algebras, and 
linearized invariants. Classical objects---such as Alexander invariants, Malcev Lie 
algebras, and resonance varieties---admit infinitesimal counterparts within 
this framework, related by a precise dictionary.

\subsection*{Topological versus infinitesimal invariants}
On the topological side lie equivariant chain complexes, Alexander invariants, 
and representation-theoretic jump loci; on the infinitesimal side lie Koszul modules, 
holonomy Lie algebras, and resonance schemes defined from models. 
In favorable cases these worlds coincide, while in general their 
divergence encodes higher-order structure.

The basic objects studied in this paper are spaces, groups, and maps,
together with their algebraic avatars arising from rational homotopy
theory, Lie theory, and commutative algebra.
Although manifolds and their fundamental groups provide many of the
motivating examples, our methods and results apply more broadly to
arbitrary spaces and finitely generated groups admitting finite-type
$\cdga$ models.

A central structural thread running through the paper concerns the interplay between 
the first Koszul module $\B_1(A)$, the holonomy Lie algebra $\h(A)$, degree-$1$ 
resonance, and Chen ranks. While these objects have appeared separately in the 
literature, one of the main contributions of this work is to show that they are governed 
by a single algebraic mechanism, intrinsic to finite-type $\cdga$ models.

A key technical ingredient is an explicit presentation of $\B_1(A)$ as a finitely presented 
module over the polynomial ring $S=\Sym(H_1(A))$. This presentation depends only on 
the linear differential and quadratic multiplication data of the model,
and extends, in a conceptually transparent way,
the foundational constructions of Papadima--Suciu in the quadratic case. 
It provides the algebraic bridge between Koszul homology and the infinitesimal Alexander 
invariant of the holonomy Lie algebra.

Building on this presentation, we establish a natural isomorphism $\B_1(A)\cong \B(\h(A))$, 
placing Chen ranks and degree-$1$ resonance under direct homological control. 
A notable feature of our approach is that this identification is obtained by purely 
algebraic means, via Koszul homology and general properties of exact complexes, 
rather than by dimension counting or external group-theoretic input. As a consequence, 
classical results---such as the computation of Chen ranks for free groups---emerge 
formally from the theory, rather than serving as foundational inputs.

From this perspective, resonance, holonomy, and Chen ranks emerge as 
different manifestations of a single algebraic mechanism, here made 
explicit through the isomorphism $\B_1(A)\cong\B(\h(A))$.
This clarifies the role of earlier quadratic and formal constructions and 
extends them to $\cdgas$ with nonzero differential.

At a methodological level, the theory replaces nonlinear or equivariant constructions 
by linearized, functorial invariants defined directly from models. Universal abelian covers, 
Alexander invariants, and Malcev completions thus acquire Koszul-theoretic counterparts 
that retain first-order information while making higher-order structure
accessible to computation.

%%%%%%%%%%%%%%%%%%%%%%
\subsection*{Finiteness and quadraticity}
%%%%%%%%%%%%%%%%%%%%%%

The primary focus is on spaces and groups that are typically non-formal 
but admit finite-type or $1$-finite models. Distinctions between quadratic, 
filtered, and non-quadratic structures play a central role in determining 
which invariants are rigid and which detect genuinely higher-order phenomena.

Objects fall naturally into three broad classes. \emph{Formal} objects 
are those for which cohomology provides a complete model and all 
invariants agree. At the opposite extreme are objects admitting no 
finite-type models, where one is forced to work entirely at the 
cohomology level. The primary focus of this paper lies in the 
intermediate regime: objects that are typically non-formal yet admit 
finite-type $\cdga$ models. It is precisely in this last setting that 
the comparison between classical and infinitesimal invariants becomes 
meaningful and effective.

A further organizing principle concerns quadraticity.
At the level of Lie algebras, one may distinguish between quadratic,
homogeneous, and non-homogeneous structures; at the level of $\cdgas$,
between formal, partially formal, and non-formal objects.
These distinctions are reflected in the behavior of Koszul homology,
holonomy Lie algebras, and resonance varieties, and they play a central
role in the applications developed later in the paper.

\begin{table}[ht]
\centering
\renewcommand{\arraystretch}{1.3}
\begin{tabular}{|l|l|}
\hline
\textbf{Classical / Topological}
& \textbf{Infinitesimal / Algebraic} \\
\hline
Space $X$, group $G$
& Finite-type $\cdga$ model $A$ \\
\hline
Universal abelian cover $X^{\ab}$
& Degree-$1$ minimal model \\
\hline
Alexander invariant $B_1(G;\k)$
& First Koszul module $\B_1(A)$ \\
\hline
Malcev Lie algebra $\m(G)$
& Holonomy Lie algebra $\h(A)$ \\
\hline
Chen ranks $\theta_n(G)$ 
& Holonomy Chen ranks $\theta_n(A)$ \\
\hline
Resonance of $H^{*}(X)$
& Resonance of a model $A$ \\
\hline
Formality of spaces or groups
& Formality of $\cdga$ models \\
\hline
\end{tabular}
\vspace{6pt}
\caption{Classical and infinitesimal invariants: a parallel.}
\label{tab:basic-dictionary}
\end{table}
\renewcommand{\arraystretch}{1.0}

Table~\ref{tab:basic-dictionary} summarizes the main parallel objects considered
throughout this work, arranged so as to emphasize the dictionary between
topological, rational-homotopical, Lie-theoretic, and algebraic invariants.
Each block of the table corresponds to a major thematic component of the
paper and serves as a guide to the structure of the next several subsections.

Having outlined the two parallel worlds we study, we now 
turn to the main results connecting them, organized by theme.
The results establish precise correspondences between
classical invariants arising from abelian covers and their
infinitesimal counterparts from finite-type $\cdga$ models,
and may be grouped into five interrelated themes: 
Koszul linearization and spectral sequences; 
holonomy Lie algebras and the first Koszul module; 
topological linearization and Alexander invariants; 
algebraic models, Chen ranks, and nilpotent Lie algebras; 
and resonance varieties.

%%%%%%%%%%%%%%%%%%%%%%%%%%%%%%
\section{Main themes and results}
\label{sect:results}
%%%%%%%%%%%%%%%%%%%%%%%%%%%%%%

The main results are organized around a functorial Koszul mechanism 
linking commutative differential graded algebra, Lie theory, algebraic topology, 
and group theory.  They fall into five interrelated themes: the Koszul 
linearization of a $\cdga$ and the spectral sequences it carries; the 
identification of the first Koszul module with the infinitesimal Alexander 
invariant of the holonomy Lie algebra, together with its topological 
counterpart; the computation of Chen ranks, for groups and for Lie algebras; 
resonance varieties, cohomology jump loci, and their tangent cones; and 
applications, to nilmanifolds, to configuration spaces of elliptic curves, and 
to Sasakian and Seifert fibered spaces.  The subsections below take these up in 
turn.

Everything in what follows is stated for $\cdgas$, and for spaces and groups
admitting finite-type models, without formality hypotheses.  Removing those 
hypotheses is the point: when $A$ is formal the invariants studied here are 
determined by $H^{*}(A)$, and the theory reduces to the graded case of
\cite{PS-crelle, AFRS24}.  Concretely, the $\cdga$-level theory yields the
following.
\begin{itemize}[itemsep=3pt]
\item The tangent cone inclusion $\TC_0\RR^{i,s}(A)\subseteq\RR^{i,s}(H^{*}(A))$
for every $i$-finite $\cdga$ (Theorem~\ref{thm:tangent-cone}), which answers
\cite[Prob.~3.2]{Su-indam}.  For finite-type $\cdgas$ the inclusion was
obtained in \cite{BR} by deformation-theoretic means.
\item Invariance of the germ at the origin of the resonance loci under
$q$-quasi-isomorphism (Theorem~\ref{thm:germ-qiso}), over any algebraically
closed field of characteristic $0$.  Over $\C$ this recovers the
model-independence deduced in \cite{DP-ccm} from the comparison with
character varieties.
\item A determination of the first resonance variety
$\RR^1(\Conf(\E,n);\k)$, for every $n\ge2$
(Proposition~\ref{prop:res-elliptic}).  This description was asserted in
\cite[Ex.~10.1]{DPS-duke} on the basis of computations for small $n$; the
proof given here is, as far as we know, the first.
\item Chen ranks of any finitely generated group with a $1$-finite $1$-model,
computed by a single graded module (Theorem~\ref{thm:alex-completed-intro}).
For the pure elliptic braid groups this gives
$\theta_k(P_{1,n})=\binom{n}{2}(k-1)$ for $k\ge2$, and non-$1$-formality of
$P_{1,n}$ for $n\ge3$ (Theorem~\ref{thm:elliptic-intro}), without recourse to
the tangent cone theorem for $1$-formal groups used in \cite{DPS-duke}.
\item A multigraded refinement of the construction
(Theorem~\ref{prop:multigraded-koszul}).  If the weight grading of $A$ is
replaced by a grading by a finitely generated free abelian group $\Gamma$,
compatible with the product and the differential and positive in the sense of
Definition~\ref{def:multihomog}, then the Koszul modules and the Chen ranks
refine into $\Gamma$-graded invariants, still computed from the cohomology
algebra.  For $\Gamma=\Z$ this is the positive-weight theory; for
$\Gamma=\Z^2$ it contains the bi-graded Koszul modules of Raicu and
Sam~\cite{RS}, and for spaces and groups admitting such a model it answers the
second half of a question of Farkas \cite[Quest.~5.8]{Farkas}, which asks for a
bigraded theory of Chen invariants of topological objects.
\end{itemize}

The Koszul modules used here are the objects that appear in the resolution of
Green's generic syzygy conjecture for canonical curves
\cite{AFPRW19, AFPRW22} and in the proof of the Chen ranks conjecture for
$1$-formal groups \cite{AFRS24, AFRS25}; Section~\ref{sect:chen-res}
describes that connection, and \cite{Farkas} surveys it from the
algebro-geometric side.

%%%%%%%%%%%%%%%%%%%%%%%%%%%%%%
\subsection*{Koszul linearization and spectral sequences}
%%%%%%%%%%%%%%%%%%%%%%%%%%%%%%

We begin with the Koszul (co)homology of a finite-type $\cdga$ model. 
Let $(A,d)$ be a connected $\k$-$\cdga$ with $H^1(A)$ 
finite-dimensional. Set $S=\Sym(H_1(A))$, where $H_1(A)$ denotes 
the dual of $H^1(A)$, and let 
$\omega_A\in H^1(A)\otimes_{\k} H_1(A)$ be the canonical element 
(a cocycle, since $d$ is zero on degree-$0$ elements). 
The \emph{Koszul chain complex} of $A$ is the $S$-linear complex
\begin{equation}
\label{eq:KbulletA}
K_\bullet(A)=(A\otimes_{\k} S,\partial^A),
\qquad
\partial^A = d+ \omega_A,
\end{equation}
where $\omega_A$ acts as an $S$-linear derivation of 
$A\otimes_{\k} S$ via contraction on the $A$-factor. 
Its homology modules, $\B_i(A)\coloneqq H_i(K_\bullet(A))$, 
are called the \emph{Koszul modules} of $A$; the notion and the
terminology originate in \cite{PS-crelle}, where the Koszul module of a
pair $(V,K)$ with $K\subseteq\bwedge^2 V$ was introduced, together with its
resonance variety.  All of these 
constructions are functorial with respect to $\cdga$ morphisms, 
a feature that becomes essential in later comparisons with 
topological and group-theoretic invariants.

Koszul linearization leads to functorial spectral sequences that 
encode higher-order algebraic data in a controlled way. 
Our first main result establishes the existence and structure of 
these Koszul spectral sequences. In both homological and 
cohomological form, they are shown to be compatible with natural 
comodule or algebra structures, with differentials expressed 
in terms of iterated left cup-products (or contraction) by the 
canonical class $\omega_A$.

\begin{thm}
\label{thm:koszul-ss-intro}
Let $(A,d)$ be a connected $\k$-$\cdga$ with $\dim_\k H^i(A)<\infty$ for all $i$.
\begin{enumerate}[itemsep=3pt]
\item 
There exists a natural cohomological spectral sequence 
of bigraded $\cdgas$, with $E_1$-page the Koszul cochain complex 
of\/ $H^*(A)$ and $E_2$-page its Koszul cohomology:
\[
E_1^{p,q} \cong K^{p+q}(H^*(A))_p, \qquad
E_2^{p,q} \cong \B^{p+q}(H^*(A))_p
\;\Longrightarrow\; \gr^\m_p \B^{p+q}(A),
\]
where the subscript $p$ denotes the degree-$p$ piece in the natural
$S$-grading on $K^\ast(H^*(A))$ and $\B^\ast(H^*(A))$.
\item
There exists a natural second-quadrant homological spectral sequence 
of bigraded comodules over $\Sym(H^1(A))$, 
with $E^1$-page the Koszul chain complex of $H^*(A)$ and 
$E^2$-page its Koszul homology:
\[
E^1_{p,q} \cong K_{p+q}(H^*(A))_{-p}, \qquad
E^2_{p,q} \cong \B_{p+q}(H^*(A))_{-p}
\;\Longrightarrow\; \gr^\m_{-p} \B_{p+q}(A),
\]
\textup{the abutment being the associated graded of the filtration 
induced on $\B_\bullet(A)$, which the $\m$-adic filtration determines only 
when the latter is separated.}
\end{enumerate}
The higher differentials are expressed in terms of iterated left
cup-products (or contractions) by the canonical element $\omega_A$. 
In particular, both spectral sequences degenerate at $E^2$ if $A$ is formal.
\end{thm}

The cohomological spectral sequence arises from the $\m$-adic 
filtration on $K^\bullet(A)$ and occupies the region $p\ge0$, $p+q\ge0$; 
in a fixed total degree it has infinitely many nonzero columns, so its 
convergence requires an Artin--Rees argument.  The homological spectral 
sequence arises from the $\m$-adic filtration on $K_\bullet(A)$ and lives 
in the second quadrant ($p\le 0$ in the natural indexing). Precise multiplicative and 
comultiplicative refinements are established in 
Theorems~\ref{thm:koszul-ss-coh} and~\ref{thm:koszul-ss-hom}.

The naturality of these constructions is developed in 
Section~\ref{sect:naturality-koszul}: Koszul modules are 
contravariantly functorial in general, but become covariantly 
functorial---and preserve algebra and coalgebra structures 
(Theorem~\ref{thm:naturality-alg-structures})---under morphisms 
inducing isomorphisms on $H^1$.
For instance, a quasi-isomorphism 
$\bigwedge(a) \to \CE(\sol_2)$ creates a Koszul module $\B_1$ 
supported away from the origin, showing that Koszul homology 
detects cochain-level data invisible to cohomology 
(Example~\ref{ex:non-functorial-qiso}).

When $A$ admits positive weights, the $\m$-adic filtration on 
$K_\bullet(A)$ is refined by the total-weight filtration of 
Section~\ref{subsec:weights-koszul}. As shown in 
Theorem~\ref{thm:weight-ss}, the associated weight spectral 
sequence identifies at the $E^1$-page with the Koszul chain 
complex of the cohomology algebra $H^*(A)$ (with zero 
differential), and converges strongly and unconditionally to 
the weight-graded Koszul homology of $A$. Higher differentials 
in this spectral sequence measure the deviation of $A$ from 
formality within each weight, and their vanishing in low total 
degrees provides effective obstructions to partial formality; 
see Section~\ref{subsec:koszul-formal-obstruction}.
Thus, for positive-weight models, the first nontrivial page is 
determined entirely by $H^*(A)$, while the filtration on Koszul 
homology is induced by an honest grading by weight.

%%%%%%%%%%%%%%%%%%%%%%%%%%%%%%
\subsection*{Holonomy Lie algebras and the first Koszul module}
%%%%%%%%%%%%%%%%%%%%%%%%%%%%%%

We now strengthen the finiteness hypothesis slightly: we assume 
$\dim_\k A^1<\infty$ (rather than merely $\dim_\k H^1(A)<\infty$), 
as the holonomy Lie algebra $\h(A)$ is constructed directly from $A^1$. 
Associated to $A$ is its \emph{holonomy Lie algebra} $\h(A)$, defined as the
quotient of the free Lie algebra on $A_1=(A^1)^{\vee}$ by the ideal generated by 
the image of the dual multiplication map $\mu^\vee\colon (A^2)^\vee \to \bigwedge^2 A_1$
and by the image of the dual differential $d^\vee\colon (A^2)^\vee \to A_1$.
This is a finitely presented Lie algebra whose relations encode both the 
cup product structure and the linear part of the differential of $A$.

A second main result (proved in Theorem~\ref{thm:B-holo}) is 
a substantial generalization of \cite[Prop.~9.3]{PS-imrn04}. 
It identifies the first Koszul homology module $\B_1(A)$ with 
$\B(\h(A))\coloneqq \h'(A)/\h''(A)$, the infinitesimal Alexander 
invariant of the holonomy Lie algebra. 

\begin{thm}
\label{thm:B-holo-intro}
Let $(A,d)$ be a connected $\k$-$\cdga$ with $\dim_\k A^1<\infty$. 
Then there is a natural isomorphism of $S$-modules,
\begin{equation}
\label{eq:b1a-bha}
\B_1(A) \cong \B(\h(A)).
\end{equation}
\end{thm}

The proof relies on an explicit $S$-presentation of $\B_1(A)$, established in 
Theorem~\ref{thm:Bpres}, which extends the classical presentation of infinitesimal 
Alexander invariants for quadratic algebras to the setting of finite-type $\cdga$ models.
This presentation provides the structural mechanism underlying the isomorphism 
\eqref{eq:b1a-bha}, and clarifies how $\B_1(A)$ is governed by the quadratic-linear 
relations of $\h(A)$.

As a consequence of Theorem~\ref{thm:B-holo-intro}, the holonomy Chen ranks of 
$A$ may be computed directly from the Koszul module $\B_1(A)$. More precisely,
\begin{equation}
\label{eq:theta-hilb}
\sum_{k\ge 0} \theta_{k+2}(A)\, t^k = \Hilb\bigl(\gr^\m(\B_1(A)),t\bigr),
\end{equation}
a fact proved in Proposition~\ref{prop:holo-massey}. 
This identity provides an effective method for computing holonomy Chen
ranks using commutative algebra, bypassing explicit Lie-theoretic or
Massey-product calculations.

From this perspective, $\B_1(A)$ isolates the quadratic core of the Koszul
theory: although higher Koszul homology detects non-formal and genuinely
higher-order phenomena, the first Koszul module is entirely determined by
holonomy.
As a result, invariants extracted from $\B_1(A)$---such as resonance support
loci and holonomy Chen ranks---are controlled by the quadratic data encoded
in $\h(A)$, even when the differential of $A$ is not itself quadratic.
At the same time, the isomorphism $\B_1(A)\cong\B(\h(A))$ is an 
identification of invariants of a fixed $\cdga$; under morphisms, 
the isomorphism is compatible with induced morphisms.
Both induced maps, $\h(\varphi)\colon\h(\oA)\to\h(A)$ and 
$\B_1(\varphi)$, are determined by the truncation of $\varphi$ 
in degrees at most $2$; the naturality statement identifies them 
under the isomorphism of Theorem~\ref{thm:B-holo-intro}
(see Remark~\ref{rem:naturality-B-holo-subtlety} and 
Example~\ref{ex:nonhomog-tre}). 
Since neither functor is invariant under quasi-isomorphisms, 
both detect strictly more than the induced maps on cohomology. 
This sensitivity is exploited in 
Section~\ref{sect:functorial-models}, where the induced map on 
Koszul modules obstructs the compatibility of a group homomorphism 
with chosen formality data---and, for an endomorphism inducing the 
identity on cohomology, detects a nontrivial action on the Malcev 
completion---even when source and target are individually $1$-formal 
(Proposition~\ref{prop:koszul-obstruction-groups} and 
Corollary~\ref{cor:koszul-obstruction-endo}).

Together, Theorems \ref{thm:koszul-ss-intro} and \ref{thm:B-holo-intro} form the 
algebraic backbone of the paper: they show that resonance, holonomy, and 
infinitesimal Alexander invariants are governed by a single Koszul mechanism 
intrinsic to finite-type $\cdga$ models.

%%%%%%%%%%%%%%%%%%%%%%%%%%%%%%
\subsection*{Topological linearization and Alexander invariants}
%%%%%%%%%%%%%%%%%%%%%%%%%%%%%%

The algebraic theory developed above interacts functorially with
the topology of universal abelian covers.
Let $X$ be a connected CW-complex with fundamental group
$G=\pi_1(X)$, and set $\Lambda=\k[G_{\ab}]$, $I\subset\Lambda$
the augmentation ideal, and $S=\Sym(H_1(X;\k))$.
The Alexander invariants $B_i(X;\k)\coloneqq H_i(X^{\ab};\k)$,
viewed as $\Lambda$-modules, encode information about the topology of $X$ 
that is often difficult to compute directly.
A key role in this paper is played by the identification of the
equivariant chain complexes of the universal abelian cover with
Koszul complexes, via a linearization process.

At the $E_1$-level, the equivariant spectral sequence associated to the
$I$-adic filtration on $C^\bullet(X;\Lambda)$ is canonically identified
with the Koszul cochain complex $K^\bullet(H^*(X;\k))$, for any finite-type
space $X$ (Theorem~\ref{thm:equiv-ss-coh}).
Under an additional minimality hypothesis on the CW-structure of $X$, 
this identification upgrades to an isomorphism of associated graded 
cochain complexes $\gr^I C^\bullet(X;\Lambda)\cong K^\bullet(H^*(X;\k))$,
and under further hypotheses---the existence of a finite-type $\cdga$ 
model with positive weights---it upgrades further to a filtered 
quasi-isomorphism of completions, recovering the Alexander invariants 
themselves from Koszul homology.

\begin{thm}
\label{thm:top-linearization-intro}
Let $X$ be a connected, finite-type, minimal CW-complex.
Set $\Lambda=\k[G_{\ab}]$ with augmentation ideal $I$, and
$S=\Sym(H_1(X;\k))$ with maximal ideal $\m$.
\begin{enumerate}
\item \label{C1}
There are canonical isomorphisms
\[
\gr^I C^\bullet(X;\Lambda) \cong K^\bullet(H^*(X;\k))
\cong \gr^\m K^\bullet\bigl(H^*(X;\k)\bigr),
\]
where the first is an isomorphism of bigraded dg-$S$-modules, and the second 
is an isomorphism of bigraded $S$-modules \textup{(}the $\m$-adic filtration 
on $K^\bullet(H^*(X;\k))$ is split by its natural $S$-grading\textup{)},
so the two associated graded objects coincide.
\item \label{C2}
If, moreover, the equivariant spectral sequence of $X$ degenerates at
the $E^2$-page in total degrees $\le q$, then
\[
\gr^I B_i(X;\k) \cong \B_i\bigl(H^*(X;\k)\bigr)
\qquad\text{for all } i\le q .
\]
For $i=1$ no such hypothesis is needed: if $X$ admits a $1$-finite
$1$-model $(A,d)$, then
\[
\widehat{B}_i(X;\k)_I \cong \widehat{\B}_i(A)_\m
\quad\text{and}\quad
\gr^I B_1(X;\k) \cong \gr^\m \B_1(A).
\]
\end{enumerate}
\end{thm}

Part~\eqref{C1} is established in Theorem~\ref{thm:lin-koszul}
(due to Papadima--Suciu~\cite{PS-tams}), together with the observation that
the $\m$-adic filtration on $K^\bullet(H^*(X;\k))$ is split by the 
$S$-degree grading. The minimality hypothesis is essential for the 
first isomorphism: without it, $\gr^I C^\bullet(X;\Lambda)$ is 
isomorphic to $C^\bullet(X;\k)\otimes_{\k} S$ with differential 
$\mathrm{gr}_I(\delta^{\ab})=\delta_X\otimes\mathrm{id}_S$, which differs 
from $K^\bullet(H^*(X;\k))$ whenever the cellular coboundary of $X$ is 
nontrivial; see Example~\ref{ex:trefoil} for an explicit illustration.

Part~\eqref{C2} is proved in
Theorem~\ref{thm:lin-koszul-degenerate} and
Theorem~\ref{thm:completed-koszul-space}.  The two parts play different roles.
Part~\eqref{C1} is a chain-level statement: it identifies the
linearization of the equivariant complex with a Koszul complex, and
holds for every minimal CW-complex.  Part~\eqref{C2} concerns the
homology of that complex, and here a genuine hypothesis is unavoidable.
Passing from $\gr$ to homology requires the higher differentials of the
equivariant spectral sequence to vanish, which is a formality-type
condition on $X$ rather than a property of its cell structure: the
presentation $2$-complex of the integer Heisenberg group is finite and
minimal, yet $\gr^I B_1=\k$ while $\B_1(H^*)\cong S(-1)$
\textup{(}Remark~\ref{rem:degeneration-not-minimality}\textup{)}.
In degree~$1$ the hypothesis can be dispensed with, because the
comparison may be routed through Malcev completion instead of the
spectral sequence; see Theorem~\ref{thm:alex-completed}.

Part~\eqref{C2} identifies the completed $\Lambda$-modules
$B_i(X;\k)$ with the completed Koszul modules $\B_i(A)$,
entirely by the model $(A,d)$.
When $(A,d)$ satisfies Poincar\'{e} duality
(see Section~\ref{sect:koszul-pd}), this correspondence is
further refined: if $(A,d)$ is a $\PD_m$-$\cdga$, then for
each $i\ge 0$ there is a natural isomorphism of graded $S$-modules
(Theorem~\ref{thm:pd-koszul-modules})
\begin{equation}
\label{eq:pd-bia}
\B_i(A)\cong\sigma^{*}\B^{\,m-i}(A),
\end{equation}
where $\sigma$ is the involution induced by $a\mapsto -a$ on $H^1(A)$;
thus Poincar\'{e} duality exchanges the homological and cohomological
Koszul modules in complementary degrees, rather than dualizing either.
Together, these results show that Koszul linearization provides
a functorial bridge between the $\Lambda$-module structure arising
from universal abelian covers and the computable commutative
algebra of Koszul modules attached to finite-type $\cdga$ models.

%%%%%%%%%%%%%%%%%%%%%%%%%
\subsection*{Algebraic models for spaces and groups}
%%%%%%%%%%%%%%%%%%%%%%%%

The Koszul-theoretic framework developed here builds on the rational homotopy 
theories of Quillen \cite{Quillen69} and Sullivan \cite{Sullivan77}, 
emphasizing the non-simply connected setting where the fundamental group 
and its associated algebraic invariants play a central role.
In contrast with the simply connected case---where minimal models are concentrated 
in degrees $\ge 2$---the presence of nontrivial $\pi_1$ introduces degree-one 
generators, filtered Lie algebras, and infinitesimal invariants associated to 
abelian covers. From a homological algebra viewpoint, Koszul complexes and BGG 
duality provide the mechanism for linearizing this additional structure, in the 
spirit of the Avramov–Halperin correspondence between rational homotopy theory 
and local algebra \cite{Avramov-Halperin}.

Given a $1$-finite $1$-model $(A,d)$ for a finitely generated group $G$, 
the holonomy Lie algebra $\h(A)$ is a filtered Lie algebra attached to $A$,
refining the quadratic holonomy Lie algebra $\h(H^*(A))$, which depends 
only on the cup product. By \cite[Cor.~5.8]{PS-jlms}, the Malcev Lie algebra $\m(G;\k)$ 
identifies with the lower central series completion of $\h(A)$. Thus $\h(A)$ 
encodes the full pronilpotent completion of $G$, while formality governs whether 
this completion is determined by the quadratic cohomology algebra.

Koszul modules provide a bridge between these Lie-theoretic invariants and 
algebraic models: discrepancies between $\h(A)$ and $\h(H^\ast(A))$---that is, 
failures of formality---are reflected in Koszul homology and Chen ranks. 
This viewpoint isolates precisely when quadratic data suffice and when 
higher-order structure is present.

Although phrased in terms of groups, the framework applies equally to spaces 
admitting finite-type $\cdga$ models. Examples include compact nilmanifolds, 
configuration spaces of complex curves, and Sasakian manifolds. In these 
settings, holonomy Lie algebras and Koszul invariants capture geometric 
information not visible at the level of cohomology.

%%%%%%%%%%%%%%%%%%%%%%%%%%%%%%%%%%%
\subsection*{Alexander invariants, Koszul modules, and Chen ranks}
%%%%%%%%%%%%%%%%%%%%%%%%%%%%%%%%%%%

Another organizing theme is the comparison between classical 
Alexander invariants of groups and their Koszul-theoretic 
counterparts arising from algebraic models.

For a finitely generated group $G$, the Alexander invariant $B_1(G;\k)$ is the 
$\k[G_{\ab}]$-module $H_1(G';\k)$, equivalently the first homology of the universal 
abelian cover of a $K(G,1)$. Given a $1$-finite $1$-model $(A,d)$, one associates 
the first Koszul module $\B_1(A)$, defined as the degree-one homology of the 
Koszul complex of $A$.

For $1$-formal groups, the completed Alexander invariant is determined 
by the holonomy Lie algebra \cite{PS-imrn04, DPS-duke}. 
Theorem~\ref{thm:alex-completed} removes the formality 
hypothesis and places this correspondence in a uniform 
$\cdga$-theoretic framework.

\begin{thm}
\label{thm:alex-completed-intro}
Let $G$ be a finitely generated group. Suppose $G$ admits a $1$-finite
$1$-model $(A,d)$ over a field $\k$ of characteristic $0$. Then:
\begin{enumerate}[itemsep=3pt]
\item 
The completed Alexander invariant of $G$ is isomorphic to the completed
first Koszul module of $A$,
\[
\widehat{B_1(G;\k)} \cong \widehat{\B_1(A)},
\]
as filtered modules over $\widehat{\k[G_{\ab}]}$.
\item 
Taking associated graded modules yields an isomorphism over $\Sym(H_1(G;\k))$,
\[
\gr(B_1(G;\k)) \cong \gr^\m(\B_1(A)).
\]
\item 
The Chen ranks of $G$ coincide with the holonomy Chen ranks of the model:
\[
\theta_n(G)=\theta_n(A)\qquad\text{for all }n\ge1.
\]
Moreover, the generating series of Chen ranks satisfies
\[
\sum_{n\ge 0} \theta_{n+2}(G)\,t^n = \Hilb\bigl(\gr^\m(\B_1(A)),t\bigr).
\]
\end{enumerate}
\end{thm}

This result shows that Chen ranks are infinitesimal invariants: although defined
at the level of the group, they are completely determined by the Koszul complex
of any $1$-finite $1$-model.
The identification of the generating series of Chen ranks with the Hilbert series
of $\gr(\B_1(A))$ is the mechanism behind their connection with resonance
varieties and eventual vanishing phenomena developed later in the paper.

In Section~\ref{sect:gr-koszul}, we further identify the graded pieces of 
$\B_i\bigl(H^*(A)\bigr)$ with graded $\Tor$-groups over the exterior algebra via 
BGG duality (Theorem~\ref{thm:gr-Bi-tor} and 
Corollary~\ref{cor:gr-Bi-tor-cdga}). In particular, Corollary~\ref{cor:chen-eventual} 
shows that Chen ranks exhibit rigid homological behavior: once they vanish in 
sufficiently high degree, they remain zero thereafter. Thus Chen ranks either 
persist indefinitely or vanish eventually.

%%%%%%%%%%%%%%%%%%%%%%%%%%%%%%
\subsection*{Nilpotent Lie algebras, CE-models, and Chen ranks}
%%%%%%%%%%%%%%%%%%%%%%%%%%%%%%

For finitely generated torsion-free nilpotent groups $G$, the
Chevalley--Eilenberg $\cdga$ $A=\CE(\g)$ of the Malcev Lie algebra
$\g$ provides a canonical $1$-finite $1$-model, so that 
Theorem~\ref{thm:alex-completed-intro}
applies and gives $\theta_n(G)=\theta_n(A)$ for all $n\ge1$.
The Koszul module $\B_1(A)$ thus carries complete information about
the Chen ranks of $G$, reducing an intrinsically group-theoretic
problem to a computation in commutative algebra over the polynomial
ring $S=\Sym(H_1(G;\k))$.

This reduction is especially effective for $2$-step nilpotent Lie
algebras.  For such an algebra $\h$, with $m=\dim\h_1\ge2$ and
$r=\dim\h_2$, the Chevalley--Eilenberg $\cdga$ $A=\CE(\h)$ carries a natural
positive-weight decomposition
(Proposition~\ref{prop:CE-positive-weights}), making $\B_1(A)$ a
\emph{graded} $S$-module with an explicit presentation
(Proposition~\ref{prop:B1-2step-presentation}).  That presentation is as
small as it could be: the holonomy Lie algebra of $A$ is $\h$ itself,
$\B_1(A)\cong\B(\h)=\h_2$ is a module of finite length $r$ on which $S$ acts
through the augmentation, and therefore $\theta_2(A)=r$, $\theta_n(A)=0$ for
$n\ge3$, and $\Supp_S\B_1(A)=\{0\}$, irrespective of $\h$.  All growth thus
occurs on the cohomological side, and the gap between the two sides is a
failure of $1$-formality.  The following result, proved in
Proposition~\ref{prop:B1-2step-presentation} and
Theorem~\ref{thm:chen-growth-2step}, makes this explicit.

\begin{thm}
\label{thm:2step-intro}
Let\/ $\h$ be a finite-dimensional $2$-step nilpotent Lie algebra 
over $\k$, with $m=\dim\h_1\ge2$ and $r=\dim\h_2$, and set $A=\CE(\h)$.

\begin{enumerate}[itemsep=4pt]
\item \label{tF1}
Write $W=\h_2^\vee\subseteq\bwedge^2\h_1^\vee$ and
$\bar\theta_n=\theta_n(H^*(A))$.  Then
$\h(H^*(A))=\Lie(\h_1)/\langle W^\perp\rangle$, and
\[
\RR^1\bigl(H^*(A)\bigr)=\bigcup\,\bigl\{P\subseteq\h_1^\vee \;\big|\;
\dim_\k P=2,\ \bwedge^2P\subseteq W\bigr\},
\]
a union of $2$-planes indexed by the decomposable elements of $W$.  The
ranks $\bar\theta_n$ grow polynomially of degree
$\dim\RR^1(H^*(A))-1$; they vanish for $n\gg0$ precisely when $W$
contains no nonzero decomposable element, and attain the maximal degree
$m-1$ precisely for\/ $\h=\ff_{m,2}$, where
$\bar\theta_n=(n-1)\binom{n+m-2}{n}$ is Chen's formula for a free group of
rank $m$.

\item \label{tF2}
If $A=\CE(\h)$ is a $1$-model for a finitely generated group $G$, then
$\theta_n(G)=0$ for $n\ge3$, so $G$ fails to be $1$-formal as soon as
$\bar\theta_n\ne0$ for a single $n\ge3$.  This happens for every\/ $\h$
with $r\ge1$ and $m\ge3$, whether or not $\bar\theta_n$ grows.
\end{enumerate}
\end{thm}

The geometry lives entirely in part~\eqref{tF1}: it is the
\emph{degenerate}\/ quadratic data---the decomposable elements of
$W$---that produce growth.

%%%%%%%%%%%%%%%%%%%%%%%%%%%%%%
\subsection*{Resonance varieties, cohomology jump loci, and tangent cones}
%%%%%%%%%%%%%%%%%%%%%%%%%%%%%%

Resonance varieties furnish the most concrete manifestation of the 
comparison between cohomology-level and model-level invariants developed here.
These varieties also connect to the broader theory of cohomology jump loci, 
capturing the tangent-cone approximation to characteristic varieties.

Let $(A,d)$ be a connected $\cdga$ with $\dim_\k H^1(A) < \infty$, 
and set $S = \Sym(H_1(A))$. The \emph{resonance varieties} of $A$ are the 
jump loci where Aomoto cohomology dimension increases:
\begin{equation}
\label{eq:risa-jump}
\RR^{i,s}(A) \coloneqq
\Bigl\{ \m\in \mSpec(S)\ \Big|\ 
\dim_\k H^i\bigl(K_\bullet(A)\otimes_{S} S/\m\bigr)\ge s \Bigr\}, 
\end{equation}
where $\mSpec(S)$ may be identified canonically with the affine space $H^1(A)$. 
Assuming that $\dim_\k A^j < \infty$ for $j \le i$, each set $\RR^{i,s}(A)$ with $s\ge 1$ 
is a Zariski-closed subset of $H^1(A)$, endowed 
with its natural determinantal scheme structure.

In parallel, the \emph{resonance support loci} are defined module-theoretically:
\begin{equation}
\label{eq:risa-supp}
\RR_{i,s}(A)
\coloneqq
\Supp_S\bigl(\bwedge^s H_i(K_\bullet(A))\bigr),
\end{equation}
where $\B_i(A) = H_i(K_\bullet(A))$ are the Koszul modules. While closely 
related, $\RR^{i,s}(A)$ and $\RR_{i,s}(A)$ need not coincide as sets or 
schemes (Example~\ref{ex:sol2-resonance}).

Applying these constructions to $H^{*}(A)$ (with $d=0$) yields the 
coho\-mology-level resonance varieties $\RR^{i,s}(H^{*}(A))$ and 
$\RR_{i,s}(H^{*}(A))$. In formal situations, model-level and 
cohomology-level resonance coincide; discrepancies detect non-formality.

Our results fall into two complementary parts.
First, we prove a general tangent-cone theorem showing that 
cohomology-level resonance controls only the lowest-order 
approximation at the origin. Second, under positive-weight hypotheses, 
we establish strong global constraints on model-level resonance, 
including conicality properties.

\begin{thm}
\label{thm:tangent-cone-intro}
Let $(A,d)$ be a connected $q$-finite $\cdga$ over an algebraically 
closed field $\k$ of characteristic $0$. Then for all $i \le q$,
\begin{enumerate}[label=\textup{(\roman*)}, itemsep=2pt]
\item \label{tc-i}
$\TC_0\bigl(\RR^{i,s}(A)\bigr) \subseteq \RR^{i,s}(H^{*}(A))$ 
for all $s\ge1$,
\item \label{tc-ii}
$\TC_0\bigl(\RR_{i,1}(A)\bigr) \subseteq \RR_{i,1}(H^{*}(A))$.
\end{enumerate}
Whether the inclusion in \ref{tc-ii} also holds in depth $s\ge2$ is an 
open question.  If, however, $A$ is $q$-formal, via a zig-zag of 
$q$-quasi-isomorphisms through $q$-finite $\cdgas$, then both families of 
loci are determined by the cohomology algebra in all depths: for all 
$i\le q$ and all $s\ge1$, the germs at the origin agree,
\[
\RR^{i,s}(A)_{(0)} = \RR^{i,s}(H^{*}(A))_{(0)}
\qquad\text{and}\qquad
\RR_{i,s}(A)_{(0)} = \RR_{i,s}(H^{*}(A))_{(0)},
\]
and consequently
\[
\TC_0\bigl(\RR^{i,s}(A)\bigr) = \RR^{i,s}(H^{*}(A))
\qquad\text{and}\qquad
\TC_0\bigl(\RR_{i,s}(A)\bigr) = \RR_{i,s}(H^{*}(A)).
\]
\end{thm}

The inclusions are proved in Theorem~\ref{thm:tangent-cone}: the 
cohomological one by Gaussian elimination on the determinantal matrices 
cutting out the jump loci, which exhibits every minor of the Aomoto 
matrix of $H^*(A)$ as the initial form of a function vanishing on 
$\RR^{i,s}(A)$; the homological one by an elementary embedding of 
$\gr^\m\B_i(A)$ into the second page of the Koszul spectral sequence, 
whose entries assemble to $\B_i(H^*(A))$; the open depth-$s$ case of 
\ref{tc-ii} is Problem~\ref{prob:tcone-hom-depth}.  The equalities under 
$q$-formality are Theorem~\ref{thm:formal-tc}: after passing to $\m$-adic 
completions, 
the Koszul modules become genuine invariants of $q$-quasi-isomorphism 
(Theorem~\ref{thm:q-iso-koszul}\eqref{iso3}), so the germs at the origin of 
all four families of loci can be transported along the $q$-formality zig-zag 
to $(H^*(A),0)$, where they are cones.  The germ invariance itself is 
Theorem~\ref{thm:germ-qiso}; it gives an algebraic proof, valid over any 
algebraically closed field of characteristic $0$, of the 
model-independence of these germs, which over $\C$ was obtained by 
Dimca and Papadima \cite{DP-ccm} (and, for $i=q=1$, in \cite{DPS-duke}) as a 
consequence of their analytic comparison with characteristic varieties.
The tangent cone inclusions answer a question 
posed in \cite{Su-indam}, and complement the work of Budur and Rubi\'{o}~\cite{BR}.
It shows that the initial homogeneous part of model-level resonance 
is governed by cohomology, while higher-order structure 
depends on the full $\cdga$ model. In later sections, 
we show how failures of equality in the inclusion provide 
effective obstructions to formality and $1$-formality. 
For instance, for the Heisenberg nilmanifold 
(Example~\ref{ex:heisenberg-resonance}), the model-level 
resonance is trivial while cohomology-level resonance fills 
all of $H^1$, demonstrating how cohomology alone misses 
essential infinitesimal structure.
It remains an open question whether, under $q$-formality, the 
irreducible components of $\RR^{i,s}(A)$ are always linear 
subspaces of $H^1(A)$ defined over $\Q$ 
(Conjecture~\ref{conj:formal-tc}).

Positive weights provide a particularly transparent framework for
understanding resonance in finite-type models. 
If $(A,d)$ admits a positive weight decomposition compatible 
with the differential, the Koszul modules become genuinely 
graded rather than merely filtered. 
In this graded setting, the resonance loci $\RR_{i,s}(A)$ are invariant 
under a weighted $\k^*$-action, and are conical whenever $H^1(A)$ is 
concentrated in a single weight---as it is, for instance, for the canonical 
weights on a formal model.  They then often decompose as finite unions of 
linear subspaces in $H^1(A)$, reflecting strong rigidity imposed by weight 
considerations.  In particular, Theorem~\ref{thm:resonance-inclusion} shows 
that for such models all resonance and support loci are conical and hence, 
by Theorem~\ref{thm:tangent-cone-intro}, contained in the corresponding 
loci of $H^{*}(A)$.

The tangent cone inclusion provides an algebraic 
foundation for the topological theory of cohomology jump loci. 
For a space $X$, the main loci of interest are: the 
\emph{characteristic varieties} $\VV^{i,s}(X;\k)$, defined from 
rank-one local systems; the \emph{Alexander varieties} 
$\VV_{i,s}(X;\k)=\Supp(\bigwedge^s B_i(X;\k))$, defined from 
Alexander invariants; and the \emph{resonance varieties} 
$\RR^{i,s}(X;\k)=\RR^{i,s}(H^*(X;\k))$ and support loci 
$\RR_{i,s}(X;\k)$, defined from cohomology.

Libgober \cite{Li02} proved that the tangent cone at the identity to $\VV^{i,s}(X;\C)$ 
is included in $\RR^{i,s}(X;\C)$. Dimca--Papadima--Suciu \cite{DPS-duke} and 
Dimca--Papa\-dima \cite{DP-ccm} showed that for $q$-formal spaces, the analytic 
germs of $\VV^{i,s}(X;\C)$ at the identity and $\RR^{i,s}(X;\C)$ at the origin coincide 
for $i \le q$.  Theorem~\ref{thm:tangent-cone-intro} supplies the algebraic mechanism 
underlying this topological inclusion. When $X$ admits a finite-type model $(A,d)$ 
over $\C$, the analytic germ of $\VV^{i,s}(X;\C)$ at the identity is identified with 
the germ of $\RR^{i,s}(A)$ at the origin. The containment 
$\TC_0(\RR^{i,s}(A)) \subseteq \RR^{i,s}(H^{*}(A))$ then shows that 
cohomology governs the lowest-order infinitesimal structure of resonance, 
a conclusion that extends to any $i$-finite model by a finite-dimensional 
approximation theorem for resonance data (Theorem~\ref{thm:cdga-skeleton}).
Higher-order terms depend on the full $\cdga$ model.

%%%%%%%%%%%%%%%%%%%%%%%%
\subsection*{Chen ranks and resonance}
%%%%%%%%%%%%%%%%%%%%%%%%

The relationship between Chen ranks and resonance is the most fully
developed thread in this circle of ideas, and a principal motivation for the
present work. The \emph{Chen ranks conjecture}, formulated in
\cite{Su-conm01} for fundamental groups of complex hyperplane arrangements,
predicts that in sufficiently high degrees the Chen ranks $\theta_n(G)$ are
governed by the geometry of the first resonance variety
$\RR^1(G;\C)$---concretely, by the number and dimensions of its irreducible
components. For $1$-formal groups this is now a theorem: building on
Cohen--Schenck \cite{Cohen-Schenck15}, the conjecture was established for
groups with strongly isotropic resonance in \cite{AFRS24}, and sharpened to
the effective range $n\ge b_1(G)-1$ in \cite{AFRS25}
(Theorem~\ref{thm:chen-afrs}).

These advances are part of a broader pattern. The Koszul-module technology
organized in this monograph has, in its several incarnations, repeatedly
produced results of independent interest beyond homotopy theory. The same
finite-length mechanism that forces Chen ranks to vanish underlies the
resolution of Green's generic syzygy conjecture for canonical curves
\cite{AFPRW19, AFPRW22}, while resonance and Koszul-module methods have
yielded structural results on Torelli groups and the Johnson filtration
\cite{DP-ann, PS-jtop}.  In each case the mechanism is a statement about 
$\B_1$ of a graded algebra.  What is identified here is the object 
responsible for the Chen ranks/resonance correspondence in the absence of 
formality, where no such graded algebra is available.

That object is the Koszul module $\B_1(A)$. For a $1$-formal group the link
is mediated by the holonomy Lie algebra, through the epimorphism
$\h(G)/\h(G)''\surj\gr(G/G'')$ of Theorem~\ref{thm:ps-sw}, an isomorphism
exactly when $G$ is $1$-formal. Theorem~\ref{thm:alex-completed-intro}
dispenses with formality: for any finitely generated group admitting a
$1$-finite $1$-model $(A,d)$, the Chen ranks are computed by $\B_1(A)$ alone,
via $\theta_{n+2}(G)=\dim_\k\gr^\m_n\B_1(A)$. The structural conditions under
which this module is in turn controlled by resonance---isotropicity and
separability of the resonance components---are isolated in
\S\ref{subsec:separable-isotropic} and applied in
\S\ref{subsec:chen-ranks-conj}; their failure is genuine, as right-angled
Artin groups and the upper McCool groups show (Remark~\ref{rem:chen-other}).

Beyond formality the picture is no longer purely cohomological: resonance
acquires differential features, and strong isotropicity must give way to a
condition on $(A,d)$ involving both the multiplication and the differential.
We do not settle the resulting non-formal conjecture
(Question~\ref{q:chen-nonformal}), but we record its first confirmed
instance. The pure elliptic braid groups $P_{1,n}$ are non-$1$-formal, yet
their Koszul modules decompose exactly along resonance components, 
yielding closed-form Chen ranks.
This is the prototype the non-formal theory is built to capture, 
and the reason the passage from $\cgas$ to $\cdgas$ undertaken 
here is more than generalization for its own sake.

%%%%%%%%%%%%%%%%%%%%%%%%
\subsection*{Configuration spaces on elliptic curves}
%%%%%%%%%%%%%%%%%%%%%%%%

The Koszul module techniques developed here apply to any $\k$-$\cdga$ model 
carrying a positive-weight decomposition, whether or not that model is a 
Chevalley--Eilenberg complex.  The configuration spaces of an elliptic curve 
$\E$ furnish an example of the second kind.  The pure elliptic braid groups 
$P_{1,n}=\pi_1(\Conf(\E,n))$ are not $1$-formal, yet they admit an explicit 
finite-type $\k$-$\cdga$ model $A(n)$, due to Bibby~\cite{Bi16}, whose 
positive weights come from Deligne's mixed Hodge theory on quasi-projective 
varieties.  The model is of Morgan--Gysin type: its degree-$1$ generators are 
the torus classes $a_i,b_i$ of weight~$1$ and the hyperplane classes $e_{ij}$ 
of weight~$2$, with $d e_{ij}=(a_i-a_j)(b_i-b_j)$.

\begin{thm}
\label{thm:elliptic-intro}
Let $\E$ be a smooth elliptic curve over $\C$, let $n\ge2$, let 
$P_{1,n}=\pi_1(\Conf(\E,n))$, and let $A(n)$ be the Bibby model, with 
$S=\k[x_1,y_1,\dots,x_n,y_n]$.
\begin{enumerate}[itemsep=4pt]
\item \label{tG1}
The Chen ranks of $P_{1,n}$ are given by 
$\theta_k(P_{1,n})=\binom{n}{2}(k-1)$ for all $k\ge2$, whereas the holonomy 
Chen ranks $\theta_k(\h(P_{1,n};\k))$ grow polynomially of degree exactly 
$n-1$.  For $n\ge3$ the group $P_{1,n}$ is not $1$-formal.

\item \label{tG2}
The first resonance variety of $\Conf(\E,n)$ consists of the pairs 
$(x,y)\in\k^n\times\k^n$ with $\sum_{i} x_i=\sum_{i} y_i=0$ and 
$x_iy_j-x_jy_i=0$ for all $i<j$; that is, it is the affine cone over the 
Segre variety $\mathbb{P}^1\times\mathbb{P}^{n-2}$.  The resonance variety 
of the model, $\RR^1(A(n))$, is a union of $\binom{n}{2}$ linear subspaces of 
dimension~$2$, and for $n\ge3$ the two are not isomorphic as varieties.

\item \label{tG3}
As a graded $S$-module,
\[
\B_1(A(n))\cong\bigoplus_{1\le i<j\le n} S/I_{ij},
\qquad
I_{ij}=\bigl(x_\ell,y_\ell \mid \ell\ne i,j\bigr)+(x_i+x_j,\;y_i+y_j),
\]
each summand being a polynomial ring in two variables.
\end{enumerate}
\end{thm}

Part~\eqref{tG3} is Theorem~\ref{thm:BA-n}.  Part~\eqref{tG2} is established
in Proposition~\ref{prop:res-elliptic},
Corollary~\ref{cor:BA-n-consequences}\eqref{bn3}, and
Proposition~\ref{prop:conf-elliptic-comparison}\eqref{pn1}; part~\eqref{tG1}
in Proposition~\ref{prop:conf-elliptic-comparison}\eqref{pn2}--\eqref{pn3}
and Corollary~\ref{cor:p1n-nonformal}.

The contrast with the classical pure braid group is sharp.  Being an
arrangement group, $P_n$ is $1$-formal, so its Chen ranks and its holonomy
Chen ranks agree; by Cohen--Suciu \cite{CS95} they are
$\theta_k(P_n)=(k-1)\binom{n+1}{4}$ for $k\ge3$.  Both families thus have Chen
ranks growing linearly in $k$, with leading coefficients $\binom{n+1}{4}$ and
$\binom{n}{2}$.  What separates them is not that growth but the gap between the
two invariants, which vanishes for $P_n$ and, for $P_{1,n}$ with $n\ge3$, grows
polynomially of degree $n-1$.

Part~\eqref{tG3} is what drives the other two: the summands of $\B_1(A(n))$ 
correspond bijectively to the components of $\RR^1(A(n))$, so the Chen ranks 
of $P_{1,n}$ are read off componentwise.  A decomposition of this shape is 
what strong isotropicity delivers in the $1$-formal case 
(Theorem~\ref{thm:chen-afrs}); here it holds for a group that is neither 
nilpotent nor $1$-formal, and it is the first such instance we know of.  The 
resulting linear growth in $k$ matches that of the $2$-step nilpotent case 
whose resonance is a union of $2$-planes 
(Theorem~\ref{thm:2step-intro}\eqref{tF1} and 
Example~\ref{ex:pencils-2step}\eqref{pen2}), but arises by a different 
mechanism.  These are the prototypes that 
Question~\ref{q:chen-nonformal} is meant to capture; the details are in 
Section~\ref{sect:qp}.

%%%%%%%%%%%%%%%%%%%%%%%%
\subsection*{Sasakian manifolds and Seifert fibered spaces}
%%%%%%%%%%%%%%%%%%%%%%%%

A compact Sasakian manifold $M^{2n+1}$ carries a canonical $S^1$-action whose
orbit space is a compact K\"{a}hler orbifold $B$, and the Tievsky
model~\cite{Tievsky} realizes $M$ as a $1$-step Hirsch extension
\begin{equation}
\label{eq:tievsky-intro}
\cM(M)\simeq \bigl(H^*(B;\k)\otimes_{\k} \bwedge(c),\,d\bigr),
\qquad \deg(c)=1,\ d(c)=\omega\in H^2(B;\k),
\end{equation}
of the formal $\cdga$ $(H^*(B;\k),0)$, with a natural positive-weight
structure. The Hirsch-extension theory of Part~\ref{part:res-holo} then
applies directly. Theorem~\ref{thm:partial-formality-hirsch} shows that such
an extension is $q$-formal whenever the Euler class is $q$-regular, while all
Massey products of length $\ge 4$ of pure classes---in particular, of
degree-one classes---vanish (Theorem~\ref{thm:massey-hirsch-length}) 
and the triple products are governed by the explicit criterion of
Proposition~\ref{prop:triple-massey-hirsch}. Applied to the Tievsky model,
where Hard Lefschetz supplies $(n-1)$-regularity, this recovers by purely
algebraic means the partial formality of Sasakian manifolds of~\cite{PS-imrn19}, 
isolating the triple Massey product as the first possible nonformal behavior.

Beyond formality, the framework produces explicit invariants and
realizability obstructions for \emph{Sasakian groups}, a class studied by
Chen~\cite{XChen}, who asks which nilpotent groups occur
(\cite[Prob.~2]{XChen}). We compute the Malcev Lie algebra of a Sasakian
group in closed form (equation~\eqref{eq:sasaki-malcev}), as a one-dimensional
central extension of the holonomy Lie algebra of the K\"{a}hler base, reducing
a Malcev computation to cup-product data on $B$. Combined with the
Carlson--Toledo bound for nilpotent K\"{a}hler groups~\cite{CT95}, this answers
Chen's question below the rank-$8$ threshold: a nilpotent Sasakian group with
$b_1\le 7$ has Heisenberg-by-abelian Malcev Lie algebra
(Proposition~\ref{prop:sasaki-nilpotent}), so its Chen ranks are governed by
the $2$-step nilpotent computation of Theorem~\ref{thm:chen-growth-2step}.
In particular, the free $2$-step nilpotent group on $4$ generators is not a
Sasakian group, even though it satisfies every cohomological necessary
condition (Example~\ref{ex:f42-not-sasaki}). This is the infinitesimal
counterpart of the geometric classifications of Sasakian
nilmanifolds~\cite{CDMY} and solvmanifolds~\cite{Kasuya16} as finite quotients
of Heisenberg nilmanifolds.

For orientable Seifert fibered spaces $M\to\Sigma_g$ with $g>0$, the same
Hirsch-extension framework yields filtered formality of $\pi_1(M)$ via
Proposition~\ref{prop:positive-weights-ff}; see
Section~\ref{subsec:seifert-hirsch}.

%%%%%%%%%%%%%%%%%%%%%%%%
\section{Organization of the paper}
\label{intro:organization}
%%%%%%%%%%%%%%%%%%%%%%%%

This monograph is organized into four main thematic blocks, progressing from 
algebraic foundations to applications in topology and geometry.

\medskip
\noindent\textbf{Part~\ref{part:Koszul}: Algebraic foundations.}
We develop the Koszul-theoretic framework underlying the entire work.
After reviewing background from commutative algebra and the
Bernstein--Gelfand--Gelfand correspondence, we introduce finite-type
$\cdgas$ and construct the Koszul chain and cochain complexes associated
to a connected $\cdga$ $(A,d)$ with finite-dimensional $H^1(A)$.
We analyze duality properties, establish spectral sequences comparing
Koszul modules of $A$ and of $H^{*}(A)$, and relate associated graded
Koszul homology to $\Tor$ groups via BGG duality. This part provides the
homological and spectral tools used in all subsequent comparison theorems. 
We also develop foundational commutative algebra tools, including 
the identification of tangent cones to support varieties via 
associated graded modules and the formal germ of a support locus at the 
origin, used throughout Part~\ref{part:res-holo}.

\medskip
\noindent\textbf{Part~\ref{part:structural}: Structural properties.}
We investigate additional structures that refine the Koszul framework.
Positive weight decompositions are introduced and shown to promote Koszul
modules from filtered to genuinely graded objects, imposing strong
rigidity on resonance and support loci.
For spaces whose Sullivan minimal model carries positive weights---most
notably, nilmanifolds and quasi-projective varieties---the weight 
spectral sequence on the Koszul complex provides a computable 
obstruction to partial formality, with higher differentials encoding 
deviations from the cohomology algebra.
Nothing in this requires the grading group to be $\Z$: replacing it by a 
finitely generated free abelian group $\Gamma$ whose occurring multidegrees 
span a strictly convex cone refines the Koszul modules, and with them the 
Chen ranks, into $\Gamma$-graded invariants, and specializes to the 
bi-graded Koszul modules from~\cite{RS}.
We also study functoriality of Koszul modules under $\cdga$ morphisms,
showing that Koszul homology can detect non-formality of maps, even when
they induce the identity on cohomology.

\medskip
\noindent\textbf{Part~\ref{part:res-holo}: Resonance and holonomy.}
This component develops the theory of resonance varieties and 
holonomy Lie algebras and connects them systematically to Koszul modules.
We define resonance jump loci $\RR^{i,s}$ and support loci $\RR_{i,s}$,
establish the determinantal scheme structure of the former and the annihilator
scheme structure of the latter, compare the two, and prove a 
finite-dimensional approximation theorem allowing all resonance 
results to be stated under minimal finiteness hypotheses. 
We prove the tangent cone theorem, showing that cohomology 
controls only first-order behavior at the origin; we then show that a 
$q$-quasi-isomorphism preserves the entire germ of the resonance loci at the 
origin, and deduce the equality under $q$-formality. We analyze resonance 
under standard algebraic constructions and obtain vanishing results for 
Chevalley--Eilenberg complexes of nilpotent Lie algebras. This part also 
establishes the identification $\B_1(A)\cong \B(\h(A))$ and derives the 
corresponding holonomy Chen rank formulas.

\medskip
\noindent\textbf{Part~\ref{part:applications}: Spaces, groups, and applications.}
We apply the algebraic theory to topology and group theory.
After revisiting equivariant spectral sequences and Alexander invariants,
we establish our main comparison theorem for finitely generated groups
admitting $1$-finite $1$-models, identifying completed and graded
Alexander invariants with Koszul modules and equating Chen ranks.
This completes the comparison theorems stated in the Introduction
for finitely generated groups admitting $1$-finite $1$-models.
We also study formality of maps and of group homomorphisms via $\cdga$
models for homomorphisms, and record several classes---nilpotent groups,
compact K\"{a}hler manifolds and holomorphic maps, quasi-projective
manifolds---that carry functorial families of such models.
The monograph concludes with two classes of geometric applications. 
For configuration spaces on elliptic curves, we compute the Koszul 
module of Bibby's positive-weight model, determine the first resonance 
variety of $\Conf(\E,n)$ for every $n$, and obtain closed-form 
Chen rank formulas for the pure elliptic braid groups $P_{1,n}$, along with a 
proof that these groups are not $1$-formal.  This provides the first 
non-formal instances of the Chen ranks pattern established in the 
$1$-formal setting.
For Sasakian manifolds and Seifert fibered spaces, the Hirsch-extension
framework yields partial formality and sharp Massey-vanishing bounds via the
Koszul--Gysin sequence, an explicit central-extension presentation of their
Malcev Lie algebras, and realizability obstructions for nilpotent Sasakian
groups that answer a question from~\cite{XChen}.

\medskip

The four parts are designed to be read sequentially, though readers
interested primarily in applications may begin with
Part~\ref{part:applications} and refer back as needed.
The Koszul spectral sequences and weight decompositions underlie all
formality obstructions, while the structure of $\B_1(A)$ is central to
both the holonomy and Alexander comparisons. 

\medskip

\noindent {\bf General background.}
The material in this monograph draws on several foundational areas.
For rational homotopy theory, including Sullivan models, formality, 
and applications to topology, we refer to the standard sources \cite{FHT, Tanre, FOT}.
Homological algebra, derived functors, and spectral sequences are 
treated in \cite{Weibel, HS, McCleary}, while group cohomology and 
homological methods in group theory are covered in \cite{Brown}.
For commutative algebra, including graded algebras, resolutions, 
and homological techniques, we rely on \cite{Eisenbud95, Eisenbud05}.
These references provide general context and terminology; 
more specialized sources are cited at the point of use.

\setcounter{section}{0}
\renewcommand{\thesection}{\arabic{section}}
\numberwithin{equation}{section}

%%%%%%%%%%%%%%%%%
\part{Algebraic foundations}
\label{part:Koszul}
%%%%%%%%%%%%%%%%%

In this part we develop the algebraic framework underlying the entire paper.
Starting from connected, finite-type commutative differential graded algebras,
we introduce a collection of homological constructions---Koszul complexes,
Koszul modules, and their associated spectral sequences---that encode subtle
information beyond cohomology.

A recurring theme is the systematic use of duality.
Linear duality, Poincar\'{e} duality, and filtration-induced dualities interact
to constrain the structure of Koszul (co)homology and to produce invariants
that admit both algebraic and geometric interpretations.
Tools such as the BGG correspondence and homological support varieties allow
these invariants to be detected and organized effectively.

The results established here provide the algebraic foundations for the study of
resonance varieties, holonomy Lie algebras, and Alexander-type invariants in
subsequent parts, where these constructions are applied to spaces,
groups, and their infinitesimal models.

%%%%%%%%%%%%%%%%%%%%%%%%%%%%%%%%%%%%%%
\section{The BGG correspondence, supports, and jump loci}
\label{sect:BGG}
%%%%%%%%%%%%%%%%%%%%%%%%%%%%%%%%%%%%%%

We review here several algebraic constructions that underlie much of the
analysis in this paper: the Bernstein--Gelfand--Gelfand correspondence, the
notion of support for graded modules and complexes, and the resulting
homological support and jump loci. Together, these constructions identify 
homology jump loci with determinantal and support conditions, allowing 
resonance varieties to be computed scheme-theoretically via dual complexes. 
This algebraic framework serves as the foundation for the study of 
resonance varieties, Koszul complexes, and their topological counterparts 
developed in the sections that follow.

%%%%%%%%%%%%%%%%%%%%
\subsection{The BGG correspondence}
\label{subsec:bgg}
%%%%%%%%%%%%%%%%%%%%

Let $V$ be a finite-dimensional $\k$-vector space, and 
let $V^{\vee}$ be the dual $\k$-vector space. 
The Bernstein--Gelfand--Gelfand correspondence is an explicit 
equivalence of the bounded derived categories of graded modules 
over the exterior algebra $E=\bwedge V^{\vee}$ and over the symmetric 
algebra $S=\Sym(V)$; see \cite{Eisenbud05}. 
In this work, we use the BGG correspondence in the form developed by
Eisenbud--Fl{\o}ystad--Schreyer \cite{EFS}, which realizes this equivalence
via linear free resolutions and provides a natural bridge between
homological invariants and geometric support conditions.
We assign the nonzero elements of $V$ degree $1$, so that 
$S_{1}=V$, while $E_{-1}=E^{1}=V^{\vee}$ (in degree $-1$).

\begin{remark}
\label{rem:bgg-grading}
The exterior algebra $E = \bigwedge V^\vee$ in this section is 
\emph{negatively graded}: its generators $V^\vee$ sit in degree 
$-1$, so that $E^1 = V^\vee$ in the cohomological convention.
This is the natural grading for BGG duality, which relates 
positively graded $S$-modules to negatively graded $E$-modules.
When exterior algebras appear later in the paper in the context 
of Koszul complexes associated to a $\cdga$ $(A,d)$, they will carry 
a positive grading following the cohomological convention on 
$H^1(A)$; the distinction will be clear from context.
\end{remark}

In one direction, the BGG correspondence is realized by a functor $\bL$ from the 
category of graded $E$-modules to the category of linear free complexes
over $S$. This functor assigns to a graded $E$-module $P=\bigoplus P_i$ 
the chain complex
\begin{equation}
\label{eq:lp}
\begin{tikzcd}[column sep=24pt]
\bL(P): \quad \cdots \ar[r]& S \otimes_{\k} P_i  \ar[r, "\delta^P_i"]
&  S \otimes_{\k} P_{i-1}  \ar[r]& \cdots ,
\end{tikzcd}
\end{equation}
with $S$-linear differentials $\delta^P_i$ induced by left-multiplication 
by the canonical element $\omega\in V^{\vee}\otimes_{\k} V$ corresponding 
to $\id_V\in \Hom(V,V)$ under the adjoint correspondence. 
The degree conventions ensure that $\bL(P)$ is a linear complex,
that is, all maps $\delta^P_i$ are homogeneous of degree $0$. 

For a finitely generated graded $E$-module $P$, its \emph{graded dual}
is $\widehat{P}=\bigoplus (P_i)^{\vee}$.
If $P$ is a left $E$-module, then $\widehat{P}$ is naturally
a right $E$-module, and may also be viewed as a left $E$-module
by setting $(\widehat{P})_{-i}=\widehat{P}_i$.
Note that $\widehat{E}\cong E$ as \emph{ungraded}\/ left $E$-modules; the
gradings differ, $\widehat{E}$ being generated in degree~$1$ and $E$ in
degree~$-1$.
As observed in~\cite[Ex.~7F7]{Eisenbud05}, there is a natural identification
\begin{equation}
\label{eq:blp dual}
\bL(P)^{\sharp}=\bL(\widehat{P}),
\end{equation}
where $\bL(P)^{\sharp}=\Hom_S(\bL(P),S)$ is the $S$-dual complex.

Furthermore, by~\cite[Thm.~7.8]{Eisenbud05}, there is a natural duality
isomorphism of $\k$-vector spaces,
\begin{equation}
\label{eq:duality}
\big(H_q(\bL(P))_{i+q}\big)^{\vee}
  \cong \Tor^E_i(\widehat{P},\k)_{-i-q},
\end{equation}
where the right-hand side denotes the degree $-i-q$ component of
the $i$-th $\Tor$ group computed from a free $E$-resolution of~$\widehat{P}$.
By the Reciprocity Theorem~\cite[Thm.~3.7]{EFS}, the complex $\bL(P)$ is a 
free $S$-resolution of $H_0(\bL(P))$; hence, 
\begin{equation}
\label{eq:tor-s-e}
\Tor^S_i(H_0(\bL(P)), \k)_j \cong H_i(\bL(P))_j.
\end{equation}

%%%%%%%%%%%%%%%%%%%%
\subsection{Supports, annihilators, and tangent cones}
\label{subsec:supp-vars}
%%%%%%%%%%%%%%%%%%%%

Let $R$ be a commutative ring, and let $M$ be an $R$-module. 
The \emph{support} of $M$, denoted $\Supp(M)$, is the set of 
maximal ideals $\m \subset R$ for which the localization $M_\m$ 
is nonzero. Supports behave additively:
if $0\to M\to N\to P\to 0$ is an exact sequence of $R$-modules,
then $\Supp(N)=\Supp(M)\cup\Supp(P)$.
The annihilator ideal $\Ann_R(M)=\{f\in R\mid fM=0\}$
encodes the same information in the finitely generated case.
We use the Fitting ideals in the convention of \cite[\S20.2]{Eisenbud95}: if
$R^{m}\xrightarrow{\ \varphi\ }R^{n}\to M\to 0$ is a presentation, then
$\Fitt_j(M)\coloneqq I_{n-j}(\varphi)$ is the ideal generated by the
$(n-j)$-minors of $\varphi$, so that
$\bV(\Fitt_j(M))=\{\p : \mu_\p(M)\ge j+1\}$; in particular
$\bV(\Fitt_0(M))=\Supp(M)$.  The following lemma is well-known; see
e.g.~\cite[\S20.2]{Eisenbud95}.

\begin{lemma}
\label{lem:ann-supp}
If $M$ is a finitely generated $R$-module, then:
\begin{enumerate}[itemsep=2pt]
\item  \label{ann1} \
$\Supp(M)=\bV(\Ann_R(M))$;
\item  \label{ann2} 
$\bV\big(\Ann_R(\bwedge^k M)\big)=\bV(\Fitt_{k-1}(M))$ for every $k\ge1$.
\end{enumerate}
\end{lemma}

In particular, when $R=S=\Sym(V)$ is the symmetric algebra on a
finite-dimensional vector space $V$, as in Section~\ref{subsec:bgg}, the affine
scheme $\mSpec(S)$ can be identified with the dual space $V^{\vee}$.
For a finitely generated $S$-module $M$, the support $\Supp(M)=\bV(\Ann_S(M))$
is a Zariski-closed subset of $V^{\vee}$, with scheme structure defined by the
ideal $\Ann_S(M)$.

\begin{definition}
\label{def:tcone}
Let $S=\k[x_1,\dots,x_n]$ be a standard graded polynomial ring with
irrelevant maximal ideal $\m=(x_1,\dots,x_n)$.
For an ideal $J\subset S$, the tangent cone at $0$ of $\bV(J)$ is
\[
\TC_0(\bV(J)) \coloneqq \Spec\big(\gr^\m(S/J)\big),
\]
where $\gr^\m(S/J)$ is the associated graded ring with respect to the
$\m$-adic filtration.
\end{definition}

The following lemma is folklore; for completeness, we include a short proof.

\begin{lemma}
\label{lem:tc-support-general}
Let $S=\k[x_1,\dots,x_n]$ be a standard graded polynomial ring with
irrelevant maximal ideal $\m=(x_1,\dots,x_n)$, and let $M$ be a finitely
generated $S$-module equipped with its $\m$-adic filtration. Then
\[
\TC_0(\Supp_S M) = \Supp_S(\gr^\m M).
\]
In particular, if $J=\Ann_S(M)$, then
\[
\TC_0(\bV(J)) = \bV(\In(J)).
\]
\end{lemma}

\begin{proof}
Let $J = \Ann_S(M)$; since $M$ is finitely generated,
$\Supp_S(M) = \bV(J)$ by Lemma~\ref{lem:ann-supp}\eqref{ann1}.

\smallskip
\noindent\textit{Inclusion $\subseteq$.}
Let $f\in J$ be nonzero, of $\m$-adic order $e$, and write $f=\In(f)+f'$ with
$\ord(f')>e$.  For $m\in\m^pM$ we have $fm=0$, whence
$\In(f)m=-f'm\in\m^{p+e+1}M$; that is, $\In(f)$ carries $\m^pM$ into
$\m^{p+e+1}M$, and so annihilates $\gr^\m M$.  Hence
$\Ann_S(\gr^\m M)\supseteq\In(J)$ and therefore
$\Supp_S(\gr^\m M)\subseteq\bV(\In(J))$.

\smallskip
\noindent\textit{Inclusion $\supseteq$.}
It suffices to show that $\Ann_S(\gr^\m M)\subseteq\sqrt{\In(J)}$, for then
$\bV(\In(J))=\bV\bigl(\sqrt{\In(J)}\bigr)\subseteq
\bV\bigl(\Ann_S(\gr^\m M)\bigr)=\Supp_S(\gr^\m M)$.
As $\gr^\m M$ is a graded module, its annihilator is a homogeneous ideal, so
we may take $g\in\Ann_S(\gr^\m M)$ homogeneous of degree $d$.  Then
$g\cdot\m^kM\subseteq\m^{k+d+1}M$ for all $k\ge0$; in particular
$gM\subseteq\m^{d+1}M$.  Since $M$ is generated by $r$ elements, say, the
determinant trick \textup{(}Cayley--Hamilton, as in
\cite[Thm.~4.3]{Eisenbud95}\textup{)} applied to multiplication by $g$ and
the ideal $I=\m^{d+1}$ yields a relation
\[
g^r+a_1g^{r-1}+\dots+a_r\in\Ann_S(M)=J,
\qquad a_i\in I^i=\m^{i(d+1)} .
\]
Now compare $\m$-adic orders: the term $g^r$ has order exactly $rd$, being
homogeneous of that degree, whereas
$\ord\bigl(a_ig^{r-i}\bigr)\ge i(d+1)+(r-i)d=rd+i>rd$ for $i\ge1$.  Hence
the initial form of the displayed element of $J$ is $g^r$, so
$g^r\in\In(J)$ and $g\in\sqrt{\In(J)}$, as required.

\smallskip
Combining both inclusions, 
$\Supp_S(\gr^\m M) = \bV(\In(J)) = \TC_0(\Supp_S M)$.
\end{proof}

\begin{corollary}
\label{cor:tc-support-graded}
Let $S=\k[x_1,\dots,x_n]$ with maximal ideal $\m$, and let $M$ be a finitely
generated graded $S$-module equipped with its $\m$-adic filtration.
Then:
\[
\Supp_S(\gr^\m M)=\Supp_S(M),
\qquad\text{and}\qquad
\TC_0(\Supp_S M)=\Supp_S(M).
\]
In particular, the associated graded does not change the support.
\end{corollary}

\begin{proof}\sloppy
Since $M$ is graded, $\Ann_S(M)$ is a homogeneous ideal, so 
$\Supp_S(M) = \bV(\Ann_S(M))$ is a conical (i.e., $\k^\ast$-stable)
closed subset of $\mSpec(S)$.  By Lemma~\ref{lem:tc-support-general},
$\TC_0(\Supp_S M) = \Supp_S(\gr^\m M)$.  Since $\Supp_S(M)$ is already 
conical, $\TC_0(\Supp_S M) = \Supp_S(M)$, and the two equalities follow.
\end{proof}

The tangent cone is the leading term of a finer local invariant, which we
will also need.  Let $\widehat{S}=\k[[x_1,\dots,x_n]]$ denote the $\m$-adic
completion of $S$, with maximal ideal $\widehat{\m}$.

\begin{definition}
\label{def:germ}
The \emph{germ at the origin} of the closed subscheme
$\bV(J)\subseteq \Spec S$ cut out by an ideal $J\subseteq S$ is the closed
subscheme
\[
\bV(J)_{(0)}\coloneqq \bV\bigl(J\widehat{S}\bigr)\subseteq \Spec\widehat{S}.
\]
\end{definition}

\begin{lemma}
\label{lem:germ-basics}
Let $J\subseteq S$ be an ideal and let $M$ be a finitely generated
$S$-module.  Then:
\begin{enumerate}[itemsep=2pt]
\item \label{gb1}
$\In\bigl(J\widehat{S}\bigr)=\In(J)$ in
$\gr^{\widehat{\m}}\widehat{S}=\gr^\m S=S$.  Hence the tangent cone
$\TC_0(\bV(J))$ is determined by the germ $\bV(J)_{(0)}$, and indeed by its
underlying set alone, since $\sqrt{\In(J)}=\sqrt{\In\bigl(\sqrt{J}\bigr)}$.
\item \label{gb2}
$\Ann_{\widehat{S}}\bigl(\widehat{M}\bigr)=\Ann_S(M)\,\widehat{S}$, and
$\widehat{S}\otimes_{S}\bwedge^s_{S}M\cong\bwedge^s_{\widehat{S}}\widehat{M}$
for every $s\ge1$.  Consequently
\[
\bigl(\Supp_S{\textstyle\bwedge^s}M\bigr)_{(0)}
=\Supp_{\widehat{S}}\bigl({\textstyle\bwedge^s}\widehat{M}\bigr),
\]
so the germ at the origin of $\Supp_S(\bwedge^sM)$ depends only on the
$\widehat{S}$-module $\widehat{M}$.
\item \label{gb3}
If $\k=\C$, two closed subschemes of $\Spec S$ have the same germ at the
origin if and only if they have the same analytic germ at the origin.
\end{enumerate}
\end{lemma}

\begin{proof}
\eqref{gb1}
Completion along $\m$ leaves associated graded objects unchanged:
$\widehat{S}/J\widehat{S}$ is the completion of $S/J$, so
$\gr^{\widehat{\m}}\bigl(\widehat{S}/J\widehat{S}\bigr)=\gr^\m(S/J)$, and
therefore $\In(J\widehat{S})=\In(J)$.  The first assertion now follows from
Definition~\ref{def:tcone}.  For the second, $J\subseteq\sqrt{J}$ gives one
inclusion; conversely, if $f\in\sqrt{J}$, say $f^N\in J$, then
$\In(f)^N=\In(f^N)\in\In(J)$, because $S$ is a domain, whence
$\In\bigl(\sqrt{J}\bigr)\subseteq\sqrt{\In(J)}$.

\eqref{gb2}
Being finitely generated over the Noetherian ring $S$, the module $M$ is
finitely presented; hence both the formation of annihilators and that of
exterior powers commute with the flat base change $S\to\widehat{S}$.  The
displayed equality of germs then follows from
Lemma~\ref{lem:ann-supp}\eqref{ann1}.

\eqref{gb3}
Let $\mathcal{O}$ be the local ring of germs at $0$ of holomorphic functions
on $\C^n$.  Then $\widehat{S}$ is faithfully flat over $\mathcal{O}$, so
$I\widehat{S}\cap\mathcal{O}=I$ for every ideal $I\subseteq\mathcal{O}$.
Applying this to $I=J\mathcal{O}$ shows that the analytic germ, cut out by
$J\mathcal{O}$, is determined by $J\widehat{S}=(J\mathcal{O})\widehat{S}$,
and conversely.
\end{proof}

The next lemma frees the identification of
Lemma~\ref{lem:tc-support-general} from the specific choice of the
$\m$-adic filtration: any equivalent stable filtration computes the same
support.  This will be used to compare the $\m$-adic associated graded of
a Koszul module with the associated graded arising from the Koszul
spectral sequence, whose filtration is equivalent, but not equal, to the
$\m$-adic one.

\begin{lemma}
\label{lem:stable-filt-supp}
Let $M$ be a finitely generated module over $S=\k[x_1,\dots,x_n]$, and
let $F$ be a decreasing, exhaustive filtration of $M$ by $S$-submodules
which is \emph{equivalent}\/ to the $\m$-adic filtration: there is a
constant $c\ge0$ with
\[
\m^pM\subseteq F^pM\subseteq \m^{p-c}M \qquad\text{for all } p\ge c,
\]
and which satisfies $\m\cdot F^pM\subseteq F^{p+1}M$.  Then
\[
\Supp_S\bigl(\gr^F M\bigr)=\Supp_S\bigl(\gr^\m M\bigr)
=\TC_0\bigl(\Supp_S M\bigr).
\]
\end{lemma}

\begin{proof}
The second equality is Lemma~\ref{lem:tc-support-general}, whose proof
also shows that $\sqrt{\Ann_S(\gr^\m M)}$ equals $\sqrt{\In(\Ann_SM)}$.  
For the first, we compare radicals of annihilators.

Let $g\in\Ann_S(\gr^FM)$ be homogeneous of degree $d$, so that
$g\cdot F^pM\subseteq F^{p+d+1}M$ for all $p$.  Iterating and using the
equivalence, $g^kM=g^kF^0M\subseteq F^{k(d+1)}M\subseteq\m^{k(d+1)-c}M$
for $k(d+1)\ge c$.  The determinant trick, exactly as in the proof of
Lemma~\ref{lem:tc-support-general} but applied to multiplication by
$g^k$, produces an element of $\Ann_S(M)$ with initial form $g^{kr}$
once $k>c$; hence $g\in\sqrt{\In(\Ann_SM)}=\sqrt{\Ann_S(\gr^\m M)}$.
This gives $\Supp_S(\gr^\m M)\subseteq\Supp_S(\gr^FM)$.

Conversely, let $g\in\Ann_S(\gr^\m M)$ be homogeneous of degree $d$, so
that $g\cdot\m^pM\subseteq\m^{p+d+1}M$.  Then, for $p\ge c$,
\[
g^k\cdot F^pM \subseteq g^k\cdot\m^{p-c}M \subseteq \m^{p+k(d+1)-c}M
\subseteq F^{p+k(d+1)-c}M ,
\]
and for $k\ge c+1$ we have $k(d+1)-c\ge kd+1$, so $g^k$, which is
homogeneous of degree $kd$, maps $F^pM$ into $F^{p+kd+1}M$; that is,
$g^k\in\Ann_S(\gr^FM)$ (the finitely many filtration steps $p<c$ are
absorbed by enlarging $k$).  Hence
$\Supp_S(\gr^FM)\subseteq\Supp_S(\gr^\m M)$.
\end{proof}

%%%%%%%%%%%%%%%%%%%%
\subsection{Support varieties and homology jump loci}
\label{subsec:supp-jump}
%%%%%%%%%%%%%%%%%%%%

Let $\k$ be an algebraically closed field of characteristic $0$, and 
let $S$ be an affine $\k$-algebra, that is, $S=\k[x_1,\dots, x_n]/I$ 
for some ideal $I$. Since $S$ is a finitely generated $\k$-algebra, 
it is Noetherian. Moreover, $S/\m \cong \k$ for every maximal ideal $\m\in \mSpec(S)$, 
and there is a natural identification $\mSpec(S)=\Hom_{\text{$\k$-alg}} (S,\k)$ 
under which a maximal ideal $\m$ corresponds to the $\k$-algebra morphism  
$S \to S/\m \cong \k$.  The hypothesis that $\k$ be algebraically closed is in
force for the remainder of this section, and enters only through the
identification $S/\m\cong\k$, which is what permits the loci below to be
described by dimensions over $\k$.  Elsewhere in the paper $\k$ is assumed only
to have characteristic $0$, and this stronger hypothesis is restated in those
results that require it. 

We associate to a chain complex $C=(C_i,\partial_i)_{i\ge 0}$ of $S$-modules
two families of subsets of $\mSpec(S)$.

\begin{definition}
\label{def:supp-varieties}
For $i\ge 0$ and $s\ge 1$, the \emph{support loci} of $C$ are defined by
\[
\RR_{i,s}(C)\coloneqq \Supp_S\big(\bwedge^s H_i(C)\big)
   \subseteq \mSpec(S).
\]
They depend only on the chain-homotopy type of $C$.
\end{definition}

If $H_i(C)$ is finitely generated over $S$ (for example, if $C_i$ is
finitely generated), then so is $\bigwedge^s H_i(C)$ for all $s\ge 1$, and
hence the support loci $\RR_{i,s}(C)$ are Zariski-closed subsets of
$\mSpec(S)$, with scheme structure given by the ideal 
$\Ann_S\bigl(\bigwedge^s H_i(C)\bigr)$.

While support loci detect intrinsic properties of the homology modules, 
homology jump loci detect failures of exactness under base change.

\begin{definition}
\label{def:jump-loci}
The \emph{homology jump loci} measure where base change from $S$ to $S/\m$
fails to be exact. For $i\ge 0$ and $s\ge 1$ set
\begin{equation}
\label{eq:jumploci}
\RR^{i,s}(C) \coloneqq \Big\{\m\in\mSpec(S) :
       \dim_\k H_i(C\otimes_{S} S/\m)\ge s\Big\}.
\end{equation}
Again, these sets depend only on the chain-homotopy type of $C$.
\end{definition}

As noted in Remark~\ref{rem:free-necessary} below, the sets $\RR^{i,s}(C)$ 
need not be Zariski closed, even if $C$ is a chain complex of 
finitely generated $S$-modules. Nevertheless, under suitable 
finiteness and freeness assumptions on $C$, the homology 
jump loci turn out to be determinantal varieties. 

Throughout the paper, we will use the shorthand
\begin{equation}
\label{eq:jump-short}
\RR^i(C) \coloneqq \RR^{i,1}(C)
\quad\text{and}\quad
\RR_i(C) \coloneqq \RR_{i,1}(C),
\end{equation}
and refer to these sets simply as the (degree $i$) jump loci and support loci,
respectively.

Both families of loci will be compared through the number of local generators
of a module.  We record the relevant elementary facts here, as they will be
used repeatedly.

\begin{lemma}
\label{lem:supp-mu}
Let $M$ be a finitely generated $S$-module, let $\m\in\mSpec(S)$, and write
$\mu_\m(M)\coloneqq\dim_\k\bigl(M\otimes_{S} S/\m\bigr)$ for the minimal number
of generators of $M_\m$.  Then:
\begin{enumerate}[itemsep=2pt]
\item \label{sm1}
Nakayama's lemma gives the equivalences
\begin{equation}
\label{eq:nakayama-supp}
\m\in\Supp_S(M)
\iff M_\m\neq0
\iff M\otimes_{S} S/\m\neq0 .
\end{equation}
\item \label{sm2}
For every $s\ge1$,
\[
\m\in\Supp_S\bigl(\bwedge^s M\bigr)
\quad\Longleftrightarrow\quad
\mu_\m(M)\ge s .
\]
Equivalently, $\RR_{i,s}(C)=\{\m : \mu_\m(H_i(C))\ge s\}$ whenever $H_i(C)$
is finitely generated.
\end{enumerate}
\end{lemma}

\begin{proof}
Part~\eqref{sm1} is Nakayama's lemma.  For part~\eqref{sm2}, exterior powers
commute with base change, so that
$\bwedge^s M\otimes_{S} S/\m\cong\bwedge^s_{\k}\bigl(M\otimes_{S} S/\m\bigr)$, and
the latter space is nonzero precisely when $\mu_\m(M)\ge s$.  Since
$\bwedge^sM$ is again finitely generated, part~\eqref{sm1} applies to it.
\end{proof}

The next lemma is a minimality property of complexes of free modules over a
local ring: such a complex is exact at any spot where its fiber is.

\begin{lemma}
\label{lem:fiber-exact}
Let $R$ be a local ring with maximal ideal $\n_R$ and residue field
$\kappa=R/\n_R$, and let $D_\bullet$ be a chain complex of $R$-modules such
that $D_{i-1}$ and $D_i$ are finitely generated and free, and $D_{i+1}$ is
free.  If $H_i(D_\bullet\otimes_{R}\kappa)=0$, then $H_i(D_\bullet)=0$.
\end{lemma}

\begin{proof}
Write $\bar\partial$ for the differentials of $D_\bullet\otimes_{R}\kappa$, and
set $r=\rank_\kappa\bar\partial_{i+1}$; note that $r\le\rank D_i<\infty$.
Since $\bar\partial_{i+1}$ has rank $r$, some $r\times r$ minor of
$\partial_{i+1}$ is a unit in $R$.  Hence there are decompositions
$D_{i+1}=P\oplus P'$ and $D_i=W\oplus U$ into free submodules, with
$\rank P=\rank W=r$, such that $\partial_{i+1}$ maps $P$ isomorphically onto
$W$.  Replacing each $y\in P'$ by
$y-(\partial_{i+1}|_P)^{-1}\pi_W\partial_{i+1}(y)$, where $\pi_W$ denotes the
projection of $D_i$ onto $W$, we may further assume that
$\partial_{i+1}(P')\subseteq U$.

By hypothesis, $\ker\bar\partial_i=\im\bar\partial_{i+1}$; the latter
contains $\bar W$ and has dimension $r=\dim_\kappa\bar W$, so
$\ker\bar\partial_i=\bar W$.  Therefore $\bar\partial_i$ is injective on
$\bar U$, so that $\partial_i|_U$ has a unit maximal minor and is in
particular injective.  On the other hand,
$\partial_i(W)=\partial_i\partial_{i+1}(P)=0$.  Hence, if
$w+u\in\ker\partial_i$ with $w\in W$ and $u\in U$, then $\partial_i(u)=0$ and
so $u=0$; that is, $\ker\partial_i=W$.  Finally,
$\partial_{i+1}(P')\subseteq U\cap\ker\partial_i=U\cap W=0$, and thus
$\im\partial_{i+1}=\partial_{i+1}(P)=W=\ker\partial_i$: the hypothesis has
forced the three-term truncation at $i$ to be a direct sum of contractible
complexes.
\end{proof}

%%%%%%%%%%%%%%%%%%%%
\subsection{Cohomology jump loci}
\label{subsec:cjl}
%%%%%%%%%%%%%%%%%%%%

While homology jump loci are defined in terms of base change for chain complexes,
they admit a more flexible description via cohomology. 
The cohomology jump loci arise naturally by dualizing a finite truncation 
of the complex, reflecting the fact that homology jump loci in degree $i$ 
are local with respect to homological degree.
The following proposition provides a precise link between homology and
cohomology jump loci, under mild finiteness assumptions.

\begin{proposition}
\label{prop:jump-dual-local}
Let $S$ be a $\k$-algebra, and let
$C=(C_\bullet,\partial)$ be a chain complex of $S$-modules.
Fix $i\ge 0$, and assume that the $S$-modules
$C_{i-1}$ and $C_i$ are finitely generated and free.
Let $C_{\le i+1}$ denote the brutal truncation
\begin{equation}
\label{eq:brutal}
\begin{tikzcd}[column sep=24pt]
\cdots \arrow[r] & C_{i+1}
\arrow[r, "\partial_{i+1}"] & C_i
\arrow[r, "\partial_i"] & C_{i-1}
\arrow[r] & 0 .
\end{tikzcd}
\end{equation}
Then, for all $s\ge 1$,
\[
\RR^{i,s}(C) = \Bigl\{ 
\m\in \mSpec(S) \;\big|\;
\dim_\k H^i\bigl(\Hom_S(C_{\le i+1},S)\otimes_{S} S/\m\bigr) \ge s
\Bigr\}.
\]
In particular, the homology jump loci in degree\/ $i$ are completely determined
by the truncated complex $C_{\le i+1}$, and may be computed scheme-theoretically
from the dual complex $\Hom_S(C_{\le i+1},S)$.
\end{proposition}

\begin{proof}
By definition, the homology group $H_i(C\otimes_{S} S/\m)$ depends only on the
three-term complex
\begin{equation}
\label{eq:three-cc}
\begin{tikzcd}[column sep=20pt]
C_{i+1}\otimes_{S} S/\m \arrow[r] &
C_i\otimes_{S} S/\m \arrow[r] &
C_{i-1}\otimes_{S} S/\m .
\end{tikzcd}
\end{equation}
Since $C_{i-1}$ and $C_i$ are finitely generated and free over $S$,
base change to $S/\m$ commutes with $S$-linear dualization in these degrees.
Thus, the $\k$-linear dual of the above complex is canonically identified with
\begin{equation}
\label{eq:hom-complex}
\begin{tikzcd}[column sep=20pt]
\Hom_S(C_{i-1},S)\otimes_{S} S/\m \arrow[r] &
\Hom_S(C_i,S)\otimes_{S} S/\m \arrow[r] &
\Hom_S(C_{i+1},S)\otimes_{S} S/\m ,
\end{tikzcd}
\end{equation}which computes the cohomology of
$\Hom_S(C_{\le i+1},S)\otimes_{S} S/\m$ in degree $i$.

Since duality over $\k$ interchanges homology and cohomology for
finite-dimensional complexes, we obtain natural isomorphisms
\begin{equation}
\label{eq:hic-ssm}
H_i(C\otimes_{S} S/\m) \cong
\bigl( H^i(\Hom_S(C_{\le i+1},S)\otimes_{S} S/\m) \bigr)^\vee .
\end{equation}
The claim follows by comparing dimensions.
\end{proof}

\begin{remark}
\label{rem:derived-dual}
A more general version of Proposition~\ref{prop:jump-dual-local},
valid without freeness or finite-generation assumptions,
can be formulated using derived duals in the derived category $D(S)$.
Since homology jump loci depend only on a finite truncation of the complex,
the elementary argument above suffices for all applications in this work,
and we will not pursue the derived formulation further.
\end{remark}

\begin{corollary}
\label{corollary:jump-dual}
Let $C$ be a chain complex of finitely generated free $S$-modules, 
and let $C^{\#} = \Hom_S(C,S)$ be its $S$-linear dual, regarded 
as a cochain complex. Then, for all $i\ge 0$ and $s\ge 1$,
\[
\RR^{i,s}(C) = \Big\{\m\in\mSpec(S) :
\dim_\k H^i(C^{\#}\otimes_{S} S/\m)\ge s\Big\}.
\]
\end{corollary}

\begin{proof}
Since all terms of $C$ are finitely generated free, the hypotheses 
of Proposition~\ref{prop:jump-dual-local} are satisfied in every 
degree, and $\Hom_S(C_{\le i+1}, S)$ coincides with the truncation 
of $C^\#$ in degrees $\le i+1$. Since the cohomology of 
$C^\# \otimes_{S} S/\m$ in degree $i$ depends only on 
those terms, the conclusion follows from 
Proposition~\ref{prop:jump-dual-local}.
\end{proof}

%%%%%%%%%%%%%%%%%%%%%
\subsection{Jump loci, support loci, and scheme structures}
\label{subsec:jump-supp-scheme}
%%%%%%%%%%%%%%%%%%%%%
In this subsection, we clarify the relationship between several a priori
distinct ways of defining resonance varieties for complexes of free modules
over a polynomial ring. These include homological jump loci, cohomological
jump loci of dual complexes, and support loci of homology modules.
Under mild finiteness assumptions, the first two determine the same closed
subscheme, which admits a determinantal presentation; the third is a
genuinely different construction, agreeing with the others only in depth $1$
and away from an explicit exceptional locus.

As shown in 
Proposition~\ref{prop:jump-dual-local}, the homology jump loci
in degree~$i$ may be computed from the dual of the truncated complex
$C_{\le i+1}$. Under suitable finiteness hypotheses, this observation
leads to explicit determinantal equations for $\RR^{i,s}(C)$.

\begin{proposition}
\label{prop:determinantal-finite}
Let $S$ be a polynomial ring over $\k$, and let $C$ be a chain complex of free
$S$-modules such that $C_{i-1}$, $C_i$, and $C_{i+1}$ are finitely generated,
of ranks $r_{i-1}$, $r_i$, and $r_{i+1}$, respectively.
Then, for all $s\ge 1$,
\[
\RR^{i,s}(C)
=
\bV\bigl(I_{r_i-s+1}(\partial_i\oplus \partial_{i+1})\bigr),
\]
where $\partial_i\oplus\partial_{i+1}$ denotes the block-\emph{diagonal} sum
\[
C_i\oplus C_{i+1}\xrightarrow{\ \partial_i\oplus\partial_{i+1}\ }
C_{i-1}\oplus C_i ,
\]
and $I_t(\cdot)$ denotes the ideal generated by the $t\times t$ minors of the
indicated matrix.
\end{proposition}

\begin{proof}
After choosing bases, the differentials $\partial_i$ and $\partial_{i+1}$
are represented by finite matrices with entries in $S$. The condition
$\dim_\k H_i(C\otimes_{S} S/\m)\ge s$ is equivalent to a rank condition on the
specializations of these matrices at $\m$, which is precisely expressed by
the vanishing of the minors of size $r_i-s+1$.
\end{proof}

The previous proposition can be extended to situations where
$C_{i+1}$ has infinite rank, by passing to the linear dual.

\begin{proposition}
\label{prop:determinantal}
Let $S$ be a polynomial ring over $\k$, and let $C$ be a chain complex of free
$S$-modules such that $C_{i-1}$ and $C_i$ are finitely generated.
Then, for all $s\ge 1$,
\[
\RR^{i,s}(C) = \bV\bigl(I_{r_i-s+1}(\partial_i^{\#}\oplus \partial_{i+1}^{\#})\bigr),
\]
where $\partial_i^{\#}$ and $\partial_{i+1}^{\#}$ are the dual differentials in
the cochain complex $C^{\#}=\Hom_S(C,S)$.
In particular, $\RR^{i,s}(C)$ is a Zariski closed subset of $\mSpec(S)$.
\end{proposition}

\begin{proof}
Since $C_{i-1}$ and $C_i$ are finitely generated free, so are the dual modules
$C^{i-1}$ and $C^i$; hence $\partial_i^{\#}\colon C^{i-1}\to C^i$ is given by a
finite matrix.  The module $C^{i+1}$ need not be finitely generated, so
$\partial_{i+1}^{\#}\colon C^i\to C^{i+1}$ may have infinitely many rows.  Its
source has finite rank $r_i$, though, so all minors of size exceeding $r_i$
vanish, and for $t\le r_i$ the minors of
$\partial_i^{\#}\oplus\partial_{i+1}^{\#}$ generate a well-defined ideal
$I_t$ of the Noetherian ring $S$, necessarily finitely generated.

By Proposition~\ref{prop:jump-dual-local}, the homology jump loci of $C$ coincide
with the cohomology jump loci of $C^{\#}$.  The rank computation in the proof of
Proposition~\ref{prop:determinantal-finite} now applies verbatim to $C^{\#}$: it
uses only that the middle term has finite rank, the finite generation of the
outer terms serving there merely to ensure that the minor ideals are defined.
This yields the stated determinantal description.
\end{proof}

\begin{remark}
\label{rem:free-necessary}
The freeness hypothesis on $C_{i-1}$ and $C_i$ cannot be dropped in general.
Papa\-dima--Suciu \cite[Ex.~2.4]{PS-mrl} exhibit a two-term complex
\begin{equation}
\label{eq:non-free-cc}
\begin{tikzcd}[column sep=20pt]
0 \arrow[r] & \k[x] \arrow[r, "\varepsilon"] & \k \arrow[r] & 0,
\end{tikzcd}
\qquad
\varepsilon(f)=f(0),
\end{equation}
for which $C_1$ is free of finite rank but $C_0$ is not free; the resulting
jump locus $\RR^{1}(C)$ is the punctured affine line
$\mathbb{A}^1_\k\setminus\{0\}$, which is not Zariski closed.
\end{remark}

The above remark highlights that freeness assumptions are
essential in order to define jump loci via determinantal ideals.
The next result reinterprets the determinantal description of
Proposition~\ref{prop:determinantal} module-theoretically, and compares it
with the support of the cohomology of the dual complex.  The first reading
provides the bridge between the BGG perspective and the module-theoretic
approach used throughout the rest of the paper; the second is an equality only
generically.

\begin{theorem}
\label{thm:jump-support}
Let $S$ be a polynomial ring over\/ $\k$, and let $C$ be a chain complex of
free $S$-modules such that $C_{i-1}$, $C_i$, and $C_{i+1}$
are finitely generated. Set
\[
N^i\coloneqq\coker\bigl(\partial^{\#}_i\bigr),
\qquad
N^{i+1}\coloneqq\coker\bigl(\partial^{\#}_{i+1}\bigr),
\]
let $r_j=\rank C_j$, and let $\rho_{i+1}$ be the generic rank of
$\partial^{\#}_{i+1}$.  Then:
\begin{enumerate}[itemsep=3pt]
\item \label{js1}
For every $\m\in\mSpec(S)$,
\[
\dim_\k H_i\bigl(C\otimes_{S} S/\m\bigr)
=\mu_\m\bigl(N^i\bigr)+\mu_\m\bigl(N^{i+1}\bigr)-r_{i+1} ,
\]
and consequently, for all $s\ge1$,
\[
\RR^{i,s}(C)
=\Supp_S\Bigl(\bwedge^{\,s+r_{i+1}}\bigl(N^i\oplus N^{i+1}\bigr)\Bigr)
=\bV\Bigl(\Fitt_{s+r_{i+1}-1}\bigl(N^i\oplus N^{i+1}\bigr)\Bigr) .
\]
\item \label{js2}
$\Supp_S\bigl(H^i(C^{\#})\bigr)\subseteq\RR^{i,1}(C)$.
\item \label{js3}
Let $U_{i+1}\coloneqq\mSpec(S)\setminus
\bV\bigl(I_{\rho_{i+1}}(\partial^{\#}_{i+1})\bigr)$ be the open locus where
$\partial^{\#}_{i+1}$ attains its generic rank.  Then, for all $s\ge1$,
\[
\RR^{i,s}(C)\cap U_{i+1}
=\Supp_S\bigl(\bwedge^s H^i(C^{\#})\bigr)\cap U_{i+1} .
\]
In particular, if $\partial^{\#}_{i+1}=0$, then $U_{i+1}=\mSpec(S)$ and
$\RR^{i,s}(C)=\Supp_S\bigl(\bwedge^s H^i(C^{\#})\bigr)$ for all $s\ge1$.
\end{enumerate}
\end{theorem}

\begin{proof}
\eqref{js1} By Corollary~\ref{corollary:jump-dual}, 
$\dim_\k H_i(C\otimes_{S} S/\m)=\dim_\k H^i(C^{\#}\otimes_{S} S/\m)$, and the 
latter equals $r_i-\rank\bar\partial^{\#}_i-\rank\bar\partial^{\#}_{i+1}$, 
where a bar denotes reduction modulo $\m$.  Since cokernels commute with base 
change, $\mu_\m(N^i)=r_i-\rank\bar\partial^{\#}_i$ and 
$\mu_\m(N^{i+1})=r_{i+1}-\rank\bar\partial^{\#}_{i+1}$, which gives the 
displayed identity.  Applying Lemma~\ref{lem:supp-mu}\eqref{sm2} to the 
finitely generated module $N^i\oplus N^{i+1}$ now yields
\[
\RR^{i,s}(C)=\bigl\{\m : \mu_\m\bigl(N^i\oplus N^{i+1}\bigr)\ge s+r_{i+1}\bigr\}
=\Supp_S\Bigl(\bwedge^{\,s+r_{i+1}}\bigl(N^i\oplus N^{i+1}\bigr)\Bigr),
\]
and the last expression equals 
$\bV\bigl(\Fitt_{s+r_{i+1}-1}(N^i\oplus N^{i+1})\bigr)$ by 
Lemma~\ref{lem:ann-supp}\eqref{ann2}.

\eqref{js2} Regard $C^{\#}$ as a chain complex by reindexing.  Under the 
stated hypotheses, $C^{i-1}$, $C^i$, and $C^{i+1}$ are finitely generated 
free, so Lemma~\ref{lem:fiber-exact} applies at each $\m$: if 
$H^i(C^{\#}\otimes_{S} S/\m)=0$, then $H^i(C^{\#})_\m=0$.  The contrapositive 
gives the stated containment, by Lemma~\ref{lem:supp-mu}\eqref{sm1}.

\eqref{js3} Fix $\m\in U_{i+1}$, put $R=S_\m$, and write $\rho=\rho_{i+1}$.  
Some $\rho\times\rho$ minor of $\partial^{\#}_{i+1}$ is a unit in $R$, so 
there are decompositions $C^i=U\oplus V$ and $C^{i+1}=W\oplus W'$ into free 
$R$-modules, with $\rank U=\rank W=\rho$, such that $\partial^{\#}_{i+1}$ 
maps $U$ isomorphically onto $W$; adjusting $V$ as in the proof of 
Lemma~\ref{lem:fiber-exact}, we may assume $\partial^{\#}_{i+1}(V)\subseteq W'$.
Since $\rho$ is the generic rank of $\partial^{\#}_{i+1}$, the restriction 
$\partial^{\#}_{i+1}|_V$ has generic rank $0$, and hence vanishes; therefore 
$\ker\partial^{\#}_{i+1}=V$ is a free direct summand of $C^i$ of rank 
$r_i-\rho$.  As $\im\partial^{\#}_i\subseteq\ker\partial^{\#}_{i+1}=V$, we 
obtain $H^i(C^{\#})_\m=\coker\bigl(\partial^{\#}_i\colon C^{i-1}\to V\bigr)$, 
whence
\[
\mu_\m\bigl(H^i(C^{\#})\bigr)
=(r_i-\rho)-\rank\bar\partial^{\#}_i
=r_i-\rank\bar\partial^{\#}_i-\rank\bar\partial^{\#}_{i+1}
=\dim_\k H^i\bigl(C^{\#}\otimes_{S} S/\m\bigr),
\]
the middle equality because $\rank\bar\partial^{\#}_{i+1}=\rho$ on $U_{i+1}$.
Lemma~\ref{lem:supp-mu}\eqref{sm2}, applied to $H^i(C^{\#})$, now gives the 
asserted equality of loci over $U_{i+1}$.  Finally, if 
$\partial^{\#}_{i+1}=0$, then $\rho_{i+1}=0$ and $I_0(\partial^{\#}_{i+1})=S$, 
so that $U_{i+1}=\mSpec(S)$.
\end{proof}

\begin{remark}
\label{rem:jump-support-strict}
The containment in part~\eqref{js2} may be strict, and the equality in 
part~\eqref{js3} genuinely requires the restriction to $U_{i+1}$.  Let 
$S=\k[x,y]$ and let $C$ be the Koszul complex on $(x,y)$, so that $C^{\#}$ is 
the cochain complex
\[
\begin{tikzcd}[column sep=32pt]
S \ar[r, "\sbm{ x \\ y}"]
&S^2 \ar[r, "\sbm{ -y \amp x}"]
&S .
\end{tikzcd}
\]
Then $H^1(C^{\#})=\ker\sbm{-y\amp x}/\im\sbm{x\\y}=0$, since both are 
generated by $(x,y)$; hence $\Supp_S\bigl(H^1(C^{\#})\bigr)=\emptyset$.  On 
the other hand, all differentials of $C$ vanish after specializing at the 
origin, so $\RR^{1,1}(C)=\{0\}$.  Here $\rho_2=1$ while 
$\rank\bar\partial^{\#}_2=0$ at the origin, so $U_2=\mSpec(S)\setminus\{0\}$, 
and the two loci do agree over $U_2$.  This complex is the Koszul complex of 
the $\cdga$ $H^{*}(T^2;\k)$; see Example~\ref{ex:torus-jump-support}.
\end{remark}

\begin{remark}
\label{rem:scheme-structures}
Even when defined over the same affine space, resonance varieties admit
distinct constructions and may differ both scheme-theoretically and
set-theoretically.

The homological resonance varieties $\RR^{i,s}(C)$ are defined by
determinantal ideals arising from the differentials of the complex
(Proposition~\ref{prop:determinantal}), while the support loci
$\RR_{i,s}(C)$ are defined by annihilator ideals of exterior powers of
homology modules (Definition~\ref{def:supp-varieties}).
Theorem~\ref{thm:jump-support}\eqref{js1} recasts the former as the support of
an exterior power of the pair of cokernels $N^i\oplus N^{i+1}$ of the dual
complex $C^{\#}$; the supports of exterior powers of the \emph{cohomology}
$H^i(C^{\#})$ form yet a third family, related to $\RR^{i,s}(C)$ only by a
containment in depth $1$ and by equality over the locus where
$\partial^{\#}_{i+1}$ has generic rank
(Theorem~\ref{thm:jump-support}\eqref{js2} and \eqref{js3}).

In general, the support loci $\RR_{i,s}(C)$ need not coincide with the jump
loci $\RR^{i,s}(C)$, even as sets.
These two constructions agree only in a filtered sense, as explained in
Theorem~\ref{thm:PS-mrl-refined}, and may differ in fixed degree, as shown
in Example~\ref{ex:sol2-resonance}.
\end{remark}

The preceding results apply to arbitrary chain complexes of free modules 
and belong to the general BGG framework. In later sections, we specialize 
to Koszul complexes associated to $\cdgas$, where additional structure, such 
as functoriality, support loci, and duality, allows for stronger conclusions.

%%%%%%%%%%%%%%%%%%%%%%%%%%%%%
\subsection{Comparing the support and jump loci}
\label{subsec:compare-jumps}
%%%%%%%%%%%%%%%%%%%%%%%%%%%%%

The next result (which compares the two types of loci in depth~$s=1$) 
sharpens \cite[Thm.~2.5]{PS-mrl}. As in the original proof, we will 
make use of the K\"{u}nneth spectral sequence associated with the
change-of-rings map $S\to S/\m$, whose $E^2$-page is $\Tor^S_p(H_q(C),M)$;
see \cite[Thm.~5.6.4]{Weibel} or \cite[Ch.~5]{McCleary}, and compare the
proof of \cite[Thm.~2.5]{PS-mrl}, where the same terminology is used.  It
should not be confused with the base-change spectral sequence for $\Tor$
along a ring map, which is \cite[Thm.~5.6.6]{Weibel}.

\begin{theorem}
\label{thm:PS-mrl-refined}
Let $S$ be an affine $\k$-algebra and let $C=(C_\bullet,\partial)$ be a chain
complex of free $S$-modules, bounded below. Fix an integer $q\ge 0$.
\begin{enumerate}[itemsep=2pt]
\item \label{ps1} 
If $H_r(C)$ is finitely generated over $S$ for all $r\le q$, then for every
maximal ideal $\m\subset S$ the change-of-rings spectral sequence
\[
E^2_{p,r}(\m)=\Tor_p^S\bigl(H_r(C),S/\m\bigr)
  \;\Longrightarrow\; H_{p+r}(C\otimes_{S} S/\m)
\]
has finite-dimensional $\k$-vector spaces on its $E^2$-page in total degrees
$\le q$, and therefore converges in this range.
\item \label{ps2} 
If $C_r$ is finitely generated over $S$ for all $r\le q$, then
\[
\bigcup_{i=0}^j \RR_{i,1}(C) = \bigcup_{i=0}^j \RR^{i,1}(C)
  \qquad\text{for all }j\le q.
\]
\end{enumerate}
\end{theorem}

\begin{proof}
\eqref{ps1} For each $r\le q$, $H_r(C)$ is a finitely generated
module over the Noetherian ring $S$.  Hence for every maximal 
ideal $\m$ and all $p\ge0$, the $\Tor$-groups
$\Tor^S_p(H_r(C),S/\m)$ are finite-dimensional $S/\m$-vector spaces.  If
$p+r\le q$, then necessarily $r\le q$; thus all entries of the $E^2$-page in
total degree $\le q$ are finite-dimensional.  Furthermore, for fixed
$n=p+r\le q$, only the finitely many pairs $(p,r)$ with $0\le r\le n$ contribute.
Hence the $E^2$-page is finite in degrees $\le q$, and the standard
convergence criterion for first quadrant spectral sequences yields~\eqref{ps1}.

\eqref{ps2} Since $C_r$ is finitely generated over the Noetherian ring $S$ 
for all $r\le q$, and $\ker\partial_r\subseteq C_r$, the homology modules 
$H_r(C)$ are finitely generated for all $r\le q$.  Hence part~\eqref{ps1} 
applies: for every $\m\in\mSpec(S)$ the spectral sequence $E^\bullet(\m)$ 
converges in total degrees $\le q$, so that each $H_i(C\otimes_{S} S/\m)$ with 
$i\le q$ carries a finite filtration with quotients $E^\infty_{p,i-p}(\m)$, 
$0\le p\le i$.  Note also that $E^2_{p,r}(\m)=0$ for $r<0$, since $C$ is 
concentrated in non-negative degrees.

By Lemma~\ref{lem:supp-mu}\eqref{sm1} and the identification 
$\RR_{r,1}(C)=\Supp_S\bigl(H_r(C)\bigr)$, we have 
$\RR_{r,1}(C)=\{\m : H_r(C)\otimes_{S} S/\m\neq 0\}$ for all $r\le q$.  
Fix now $j\le q$ and $\m\in\mSpec(S)$.

\emph{The inclusion $\bigcup_{i\le j}\RR_{i,1}(C)\subseteq
\bigcup_{i\le j}\RR^{i,1}(C)$.}
Suppose $\m\in\RR_{i,1}(C)$ for some $i\le j$, and let $i_0$ be the 
\emph{smallest} index with $\m\in\RR_{i_0,1}(C)$, so that $i_0\le i$.  By 
minimality and \eqref{eq:nakayama-supp}, $H_r(C)_\m=0$ for all $r<i_0$; since 
$\Tor$ commutes with localization, it follows that 
\begin{equation}
\label{eq:tor-vanish-below}
E^2_{p,r}(\m)=\Tor^S_p\bigl(H_r(C),S/\m\bigr)=0 
\qquad\text{for all $p\ge 0$ and all $r<i_0$}.
\end{equation}
A differential $d^r$ with $r\ge 2$ has bidegree $(-r,r-1)$.  No such 
differential emanates from $E^r_{0,i_0}$, since its target 
$E^r_{-r,\,i_0+r-1}$ vanishes; and none lands in $E^r_{0,i_0}$, since its 
source $E^r_{r,\,i_0-r+1}$ is a subquotient of $E^2_{r,\,i_0-r+1}$, which 
vanishes by \eqref{eq:tor-vanish-below} because $i_0-r+1<i_0$.  Therefore 
\[
E^\infty_{0,i_0}(\m)=E^2_{0,i_0}(\m)=H_{i_0}(C)\otimes_{S} S/\m \neq 0 ,
\]
and since $E^\infty_{0,i_0}(\m)$ is a subquotient of 
$H_{i_0}(C\otimes_{S} S/\m)$, we conclude that $H_{i_0}(C\otimes_{S} S/\m)\neq 0$; 
that is, $\m\in\RR^{i_0,1}(C)$ with $i_0\le j$.

\emph{The inclusion $\bigcup_{i\le j}\RR^{i,1}(C)\subseteq
\bigcup_{i\le j}\RR_{i,1}(C)$.}
Suppose $\m\in\RR^{i,1}(C)$ for some $i\le j$, so that 
$H_i(C\otimes_{S} S/\m)\neq 0$.  By the finite filtration recorded above, there 
is an index $p$ with $0\le p\le i$ and $E^\infty_{p,i-p}(\m)\neq0$, whence 
$E^2_{p,i-p}(\m)=\Tor^S_p\bigl(H_{i-p}(C),S/\m\bigr)\neq0$.  Localizing at 
$\m$ shows that $H_{i-p}(C)_\m\neq0$, that is, $\m\in\RR_{i-p,1}(C)$ with 
$i-p\le i\le j$.

Taking unions over $i\le j$ gives the stated equality.
\end{proof}

\begin{remark}
\label{rem:ps-compare}
The hypothesis in Theorem~\ref{thm:PS-mrl-refined} improves on
\cite[Thm.~2.5]{PS-mrl}, which assumed that each $C_i$ is finitely 
generated. Here, for part~\eqref{ps1} it suffices that the homology
groups $H_r(C)$ be finitely generated for $r\le q$, while for 
part~\eqref{ps2} we only require that $C_r$ be finitely generated 
for $r\le q$.

This mirrors the topological situation considered in 
\cite[Prop.~4.1 and Cor.~4.3]{PS-mrl}, in which $C$ is the 
equivariant chain complex of the maximal abelian cover of a 
connected CW-complex with finite $q$-skeleton.
Thus Theorem~\ref{thm:PS-mrl-refined} recovers 
\cite[Cor.~4.3]{PS-mrl} as a direct algebraic consequence of 
the spectral sequence argument.

Both proofs rely on Nakayama's lemma.  In \cite{PS-mrl} it is applied to 
$H_q(E)_\m$ to conclude that $H_q(E)_\m/\m H_q(E)_\m\neq0$; in the proof 
above it enters through \eqref{eq:nakayama-supp}, the point being that the 
minimality of the index $i_0$ is what forces the edge term $E^2_{0,i_0}$ to 
survive to $E^\infty$.  This step cannot be bypassed.  Indeed, the edge 
homomorphism in degree $i$ factors as 
\[
\begin{tikzcd}[column sep=22pt]
H_i(C)\otimes_{S} S/\m = E^2_{0,i} \ar[r, two heads] & E^\infty_{0,i} 
\ar[r, hook] & H_i(C\otimes_{S} S/\m),
\end{tikzcd}
\]
and in general it is neither injective nor surjective.  For instance, let 
$S=\k[x,y]$, let $\m=(x,y)$, and let $C$ be the complex
\[
C:\quad \begin{tikzcd}[column sep=32pt]
S \ar[r, "\sbm{ -xy \\ x^2}"]
&S^2 \ar[r, "\sbm{ x \amp y}"]
&S .
\end{tikzcd}
\]
Then $H_0(C)\cong\k$ and $H_1(C)\cong S/(x)$, the latter generated by the 
class of the Koszul cycle $(-y,x)$, which reduces to $0$ modulo $\m$.  Hence 
the edge map vanishes in degree $1$, even though 
$H_1(C)\otimes_{S} S/\m\cong\k$ and $H_1(C\otimes_{S} S/\m)\cong\k^2$; here 
$d^2\colon E^2_{2,0}\to E^2_{0,1}$ is an isomorphism.
\end{remark}

%%%%%%%%%%%%%%%%%%%%%%%%%%%%%
\subsection{Propagation of jump loci}
\label{subsec:propagation}
%%%%%%%%%%%%%%%%%%%%%%%%%%%%%

The following propagation property of jump loci was established 
in \cite[Thm.~1.2]{DSY17}, to which we refer for the proof.
Historically, this phenomenon was motivated by earlier work of 
Eisenbud, Popescu, and Yuzvinsky \cite{EPY03}, who observed 
propagation of (depth-$1$) resonance for complements of complex 
hyperplane arrangements as a consequence of a linearity property 
of the cohomology ring under the BGG correspondence.

\begin{proposition}[\cite{DSY17}]
\label{prop:propagation}
Let $S$ be a Noetherian $\k$-algebra, and let
$C=(C_\bullet,\partial)$ be a bounded-below chain complex of finitely generated
free $S$-modules.
Assume that the homology of $C$ is concentrated in a single degree $n$, that is,
\[
H_i(C)=0 \quad\text{for all } i\neq n.
\]
Then the homology jump loci satisfy the propagation property
\[
\RR^{i}(C)\subseteq \RR^{i+1}(C)
\qquad\text{for all } i<n.
\]
\end{proposition}

\begin{corollary}
\label{cor:propagation-BGG}
Let $P$ be a finitely generated graded $E$-module.
If the BGG complex $\bL(P)$ has homology concentrated in a 
single degree $n$, then the resonance varieties 
$\RR^i(P) = \RR^i(\bL(P))$ satisfy
\[
\RR^i(P) \subseteq \RR^{i+1}(P) \qquad \text{for all } i < n.
\]
\end{corollary}

\begin{proof}
Since $P$ is a finitely generated module over the finite-dimensional
algebra $E$, each $P_i$ is finite-dimensional and vanishes for all but
finitely many $i$; hence $\bL(P)$ is a bounded complex of finitely
generated free $S$-modules, by \eqref{eq:lp}.  Its homology being
concentrated in a single degree, Proposition~\ref{prop:propagation}
applies to $C=\bL(P)$.
\end{proof}

\begin{remark}
\label{rem:propagation-vs-support}
The propagation property of Proposition~\ref{prop:propagation} is specific 
to homology jump loci and does not extend to support loci in general.

First, the sets $\RR^i(C)$ and $\RR_i(C)$ need not coincide degreewise.  
Let $S=\k[x]$ and consider the chain complex
\[
C:\quad \begin{tikzcd}[column sep=32pt]
S \ar[r, "\sbm{ 0 \\ x-1}"]
&S^2 \ar[r, "\sbm{ x \amp 0}"]
&S .
\end{tikzcd}
\]
(This complex coincides with the Koszul complex $K_\bullet(A)$ of the 
$\cdga$ $A = \CE(\sol_2)$ of Example~\ref{ex:sol2-koszul} below,
computed in~\eqref{eq:toy}.) A direct computation 
shows that $\RR^1(C)=\{0,1\}$ while $\RR_1(C)=\{1\}$. 
Thus $\RR^1(C)\neq \RR_1(C)$, even though 
$\RR^0(C)\cup \RR^1(C)=\RR_0(C)\cup \RR_1(C)=\{0,1\}$, 
in accordance with Theorem~\ref{thm:resonance-comparison}.

Second, even when the jump loci propagate, the support loci need not.  
In the above example, the jump loci satisfy $\RR^0(C)\subseteq \RR^1(C)$, whereas
$\RR_0(C)=\{0\}$ while $\RR_1(C)=\{1\}$, so $\RR_1(C)\not\subseteq \RR_0(C)$.  
Thus propagation holds for jump loci but fails for support loci.
\end{remark}

%%%%%%%%%%%%%%%%%%%%%%%%%%%%
\section{Differential graded algebras}
\label{sect:cdga}
%%%%%%%%%%%%%%%%%%%%%%%%%%

In this section we recall the basic notions and conventions concerning
commutative differential graded algebras ($\cdgas$), which serve as the
primary algebraic framework for the constructions developed in this paper.
Working in the category of connected, finite-type $\cdgas$ allows us to
formulate Koszul modules, resonance varieties, and holonomy Lie algebras in
a uniform and functorial manner, while keeping close contact with the
cohomology algebras of spaces and groups.

%%%%%%%%%%%%%%%%%%%%
\subsection{Commutative differential graded algebras}
\label{subsec:cdga}
%%%%%%%%%%%%%%%%%%%%

We shall fix throughout a coefficient field $\k$ of characteristic $0$.  
A {\em graded $\k$-vector space}\/ is a vector space $V$ over $\k$, 
together with a direct sum decomposition, $V=\bigoplus_{n\ge 0} V^n$,  
into vector subspaces. An element $a\in V^n$ is said 
to be homogeneous; we write $\abs{a}=n$ for its degree.

A {\em graded $\k$-algebra}\/ is a graded $\k$-vector space, 
$A^{\bullet}=\bigoplus_{n\ge 0} A^n$, equipped with an associative 
multiplication map, $\cdot \colon A\times A\to A$, making $A$ 
into a $\k$-algebra with unit $1\in A^0$ such that 
$\abs{a\cdot b}=\abs{a}+\abs{b}$ for all homogenous 
elements $a,b\in A$.  A graded algebra $A$ is said to be 
{\em graded-commutative}\/ (for short, a $\cga$), if 
$a\cdot b = (-1)^{\abs{a}\abs{b}} b \cdot a$ for all 
homogeneous $a,b\in A$. A basic example is the free $\cga$ 
on a graded $\k$-vector space $V^{\bullet}$;  
this algebra, denoted by $\Lambda V$, is the tensor product 
of the symmetric algebra on $\bigoplus_{n\ge 0} V^{2n}$ with 
the exterior algebra on $\bigoplus_{n\ge 0} V^{2n+1}$. 

A graded $\k$-algebra $A$ is {\em connected}\/ if $A^0$ is the 
$\k$-span of the unit $1$ (and thus $A^0=\k$).  We say that $A$ is 
of {\em finite-type}\/ if all the graded pieces $A^n$ are finite-dimensional, 
and $A$ is {\em $q$-finite} (for some integer $q\ge 1$) if 
$\dim_{\k} A^n<\infty$ for $n \le q$. 
A morphism between two graded algebras is a $\k$-linear map 
$\varphi\colon A\to B$ that preserves gradings and satisfies 
$\varphi(a\cdot b)=\varphi(a)\cdot \varphi(b)$ for all $a,b\in A$. 
We say $\varphi$ is a {\em $q$-isomorphism}\/ if $\varphi^i\colon A^i\to B^i$ 
is an isomorphism for all $i\le q$.

A {\em commutative differential graded algebra}\/ (for short, a $\cdga$) over 
a field $\k$ is a $\k$-$\cga$, $A^{\bullet}$, equipped with a 
differential $d\colon A\to A$ of degree $1$ satisfying  
the graded Leibniz rule: $d(ab)=d(a)\cdot b+(-1)^{\abs{a}} a\cdot d(b)$ 
for all homogeneous $a, b\in A$. 
The cohomology of the underlying cochain complex, $H^{*}(A)$, 
inherits the structure of a $\cga$; we will let 
$b_i(A)=\dim_\k H^i(A)$ be its Betti numbers.
If $A$ is $q$-finite, then $H^i(A)$ is a subquotient of $A^i$ 
for each $i$, and so $b_i(A)<\infty$ for all $i\le q$; 
in particular, the cohomology algebra $H^*(A)$ is then $q$-finite as well.

A morphism between two $\cdgas$, $\varphi\colon A\to B$, is both 
an algebra map and a cochain map. Consequently, $\varphi$ induces a 
morphism, $\varphi^*\colon H^{*} (A)\to H^{*} (B)$, 
between the respective cohomology algebras.  
We say that $\varphi$ is a quasi-isomorphism if $\varphi^*$ is an 
isomorphism. Likewise, we say $\varphi$ is a $q$-quasi-isomorphism (for some 
$q\ge 1$) if $\varphi^*$ is an isomorphism in degrees up to $q$ 
and a monomorphism in degree $q+1$. 
Note that, unlike a $q$-isomorphism of $\cga$s, the 
$q$-quasi-isomorphism condition for $\cdgas$ includes a 
monomorphism requirement in degree $q+1$, needed to control 
cohomology at the boundary degree.

Two $\cdgas$ $(A,d_A)$ and $(B,d_B)$ are {\em weakly equivalent} 
(or only {\em $q$-equivalent}) if there is a finite zig-zag of 
quasi-isomorphisms (or $q$-quasi-isomorphisms) connecting $A$ to $\overline{A}$, 
\begin{equation}
\label{eq:zig-zag}
\begin{tikzcd}[column sep=20pt]
A\ar[r] & A_1 &A_2\ar[l] \ar[r]& \cdots & A_{\ell}\ar[r]\ar[l]  & B .
\end{tikzcd}
\end{equation}
In this case, we write $A\simeq B$ (or $A\simeq_q B$).  

A $\cdga$ $(A,d_A)$ is said to be {\em formal}\/ (or just {\em $q$-formal}) 
if it is weakly equivalent (or just $q$-equivalent) to its cohomology 
algebra equipped with the zero differential, $(H^{*}(A),0)$. 
Formality behaves well under field extensions. 
The classical result that a connected, finite-type $\cdga$ 
over $\k$ is formal if and only if its scalar extension to any 
field $\K\supset\k$ of characteristic $0$ is formal was 
established independently by Sullivan~\cite{Sullivan77}, 
Neisendorfer--Miller~\cite{NM78}, and Halperin--Stasheff~\cite{HS79}. 
This was generalized in \cite[Thm.~4.19]{SW-forum} to partial 
formality: under the assumption that $H^{\le q+1}(A)$ is 
finite-dimensional, the $\cdga$ $(A,d_A)$ is $q$-formal over 
$\k$ if and only if $(A\otimes_{\k}\K, d_A\otimes\id_\K)$ is 
$q$-formal over $\K$.
This generalization is used throughout the paper, particularly 
when passing between $\Q$, $\R$, $\C$, and an arbitrary field 
$\k$ of characteristic $0$ in formality arguments.

%%%%%%%%%%%%%%%%%%%%
\subsection{Maps of $\cdgas$}
\label{subsec:cdga-maps}
%%%%%%%%%%%%%%%%%%%%

Let $\Lambda (t, dt) = (\Lambda (t, u),d)$ be the free $\cdga$ 
generated by elements $t$ in degree $0$ and $u$ in degree $1$, 
with differential $d$ given by $d(t)=u$ and $d(u)=0$. (This is the 
algebraic analogue of the unit interval $I = [0,1]$ in topology.) 
For each $s\in \k$, let $\ev_s\colon  \Lambda (t, dt) \to \k$ be the 
$\cdga$ morphism which sends $t$ to $s$ and $dt$ to $0$.  
Given a $\cdga$ $A$, we let 
$\Ev_s \coloneqq \id\otimes \ev_s \colon A \otimes_{\k} \Lambda (t, dt) \to 
A \otimes_{\k} \k =A$. 
Two $\cdga$ maps, $\varphi_0, \varphi_1\colon A\to B$, are 
said to be {\em homotopic}\/ if there is a $\cdga$ map, 
$\Phi\colon A\to B \otimes_{\k} \Lambda (t, dt)$, 
such that $\Ev_s\circ \Phi = \varphi_s$ for $s=0,1$. 
As is well-known (see for instance \cite[Prop.~12.8(i)]{FHT}), 
homotopic $\cdga$ maps induce the same map in cohomology. 

Two $\cdga$ morphisms, $\varphi\colon A \to B$ and 
$\varphi'\colon A' \to B'$, are {\em weakly equivalent}\/ 
(written $\varphi\simeq \varphi'$) if there exist two zig-zags of 
quasi-isomorphisms $\psi_i$ and $\psi'_i$, and 
$\cdga$ maps $\varphi_1, \dots, \varphi_{\ell -1}$ such that 
the following diagram commutes, up to homotopy:
\begin{equation}
\label{eq:ziggy-zagg}
\begin{tikzcd}[row sep=22pt]
A \ar[d, "\varphi"] & A_1 \ar[d, "\varphi_1"] \ar[l, "\psi_1"']  \ar[r, "\psi_2"] 
& \cdots & A_{\ell-1}  \ar[d, "\varphi_{\ell-1}"] \ar[l]\ar[r, "\psi_{\ell}"] 
& A'\ar[d, "\varphi'"]\, \phantom{.}
\\
B& B_1 
 \ar[l, "\psi'_1"']  \ar[r, "\psi'_2"] & \cdots 
& B_{\ell-1} \ar[l]\ar[r, "\psi'_{\ell}"] & B' .  
\end{tikzcd}
\end{equation}
The notion of {\em $q$-equivalence}\/ (written 
$\varphi\simeq_q \varphi'$) is defined similarly. 

A $\cdga$ morphism $\varphi\colon A\to B$ is said to be {\em formal}\/ 
if there is a diagram of the form \eqref{eq:ziggy-zagg} connecting $\varphi$ 
to the induced homomorphism $\varphi^*\colon H^{*}(A)\to H^{*}(B)$ 
between cohomology algebras (viewed as $\cdga${}s with zero differentials). 
Likewise, $\varphi$ is said to be {\em $q$-formal}, for some $q\ge 0$, if 
$\varphi\simeq_q \varphi^*$.

We work exclusively with \emph{commutative} differential graded algebras 
throughout this monograph. While formality can be defined for non-commutative 
$\dga$s, and Saleh \cite{Sa17} showed that a $\cdga$ is formal as a $\dga$ 
if and only if it is formal as a $\cdga$---a result generalized by 
Campos--Petersen--Robert-Nicoud--Wierstra \cite{CPRW}, who proved that 
weak equivalence of $\cdgas$ is independent of whether one computes in 
the $\dga$ or $\cdga$ category---the Koszul complex construction underlying 
our theory is \emph{geometric} and requires commutativity.

Specifically, the verification that $\delta_A^2 = 0$ for the Koszul 
differential (see \S\ref{subsec:canonical-element}) uses that 
$e_i e_j = -e_j e_i$ for cocycles $e_i, e_j \in Z^1(A)$ of degree $1$. 
Without graded commutativity, the required cancellations fail and 
$\delta_A^2 \neq 0$, so one does not obtain a chain complex. 
While non-commutative generalizations involving Hochschild homology 
or $A_\infty$-structures exist, they represent different 
theories that lie beyond the scope of classical rational homotopy theory 
as developed by Sullivan and Quillen.

%%%%%%%%%%%%%%%%%%%%
\subsection{Massey products}
\label{subsec:massey}
%%%%%%%%%%%%%%%%%%%%

A well known obstruction to formality is provided by the higher-order 
Massey products, introduced by W.S.~Massey in \cite{Ma-58}, 
and studied for instance in \cite{Kraines, May}. 
Let $(A,d)$ be a $\k$-$\dga$ and let $u_1, \dots, u_n$ 
be elements in $H^{*}(A)$; without loss of generality, 
we may assume that $n\ge 3$ and each $u_i$ is homogeneous 
and of positive degree.  
Here and in what follows we write $\bar{a}\coloneqq(-1)^{\abs{a}}a$ for a
homogeneous element $a\in A$.
A {\em defining system}\/ for $u_1, \dots, u_n$ is a collection 
of elements $a_{i,j}\in A$ such that 
$d a_{i-1,i}=0$ and $[a_{i-1,i}]=u_i$ for $1\le i\le n$ and 
$d a_{i,j}=\sum_{i<r<j} \bar{a}_{i,r} a_{r,j}$ for 
$0\le i<j\le n$ and $(i,j)\ne (0,n)$.
It is readily verified that the element 
\begin{equation}
\label{eq:massey-prod}
\alpha\coloneqq \sum_{0<r<n} \bar{a}_{0,r} a_{r,n}
\end{equation}
is a cocycle, of degree $\abs{\alpha}=2-n+\sum_{i=1}^{n} \abs{u_i}$. 
The $n$-fold Massey product $\langle u_1, \dots, u_n \rangle$, 
then, is the subset of $H^{*}(A)$ consisting of the cohomology classes 
$[\alpha]$ corresponding to all possible defining systems for $u_1, \dots, u_n$.
We say that the Massey product is {\em defined}\/ if there is at least one such 
defining system, or, equivalently, $\langle u_1, \dots, u_n \rangle\ne \emptyset$,  
in which case the indeterminacy of the Massey product is the subset 
$\{u-v\mid u,v\in \langle u_1, \dots, u_n \rangle\}\subseteq H^*(A)$.
When a Massey product is defined, we say it {\em vanishes}\/ if it contains 
the element $0$; otherwise, it is non-vanishing.

The simplest Massey triple products are as follows. 
Let $u_1, u_2, u_3$ be elements in $H^1(A)$ 
such that $u_1u_2 = u_2 u_3 = 0$. We may then choose 
$1$-cocycles $a_{0,1}, a_{1,2}, a_{2,3}$  representing $u_1, u_2, u_3$  
and $1$-cochains $a_{0,2}$ and $a_{1,3}$ such that 
$da_{0,2} = - a_{0,1} a_{1,2}$ and $da_{1,3} = - a_{1,2} a_{2,3}$. 
The triple product $\langle u_1,u_2, u_3 \rangle$ is then  
the subset of $H^2(A)$ consisting of the cohomology classes 
$-[a_{0,1} a_{1,3} + a_{0,2}a_{2,3}]$, for all such choices of 
defining systems. Due to the ambiguity in the choices made, 
$\langle u_1, u_2, u_3\rangle$ may be viewed as a coset 
of $u_1 \cdot H^1(A)+H^1(A)\cdot u_3$ in $H^2(A)$.

\begin{example}
\label{ex:non-formal-cdga-again}
Let $(A,d)$ be the exterior algebra on generators $a_1,a_2, b$ in degree $1$ 
and differential given by $d a_i=0$ and $d b=a_1a_2$. Letting $u_i=[a_i]\in H^1(A)$, 
we have that the triple Massey products $\langle u_1, u_1, u_2\rangle=\{[a_1b]\}$ and 
$\langle u_1, u_2, u_2\rangle=\{[ba_2]\}$ are defined, have $0$ indeterminacy, 
and are non-van\-ishing; in fact, the two cohomology classes generate $H^2(A)$.
Therefore, $A$ is not formal.
\end{example}

Massey products enjoy the following (partial) functoriality properties.

\begin{proposition}[\cite{Kraines, May}]
\label{prop:massey-natural}
Let $\varphi\colon A\to B$ be a $\cdga$ morphism, and let 
$\varphi^*\colon H^*(A)\to H^*(B)$ be the induced morphism 
in cohomology; then 
\begin{equation}
\label{eq:massey-func}
\varphi^*(\langle u_1, \dots, u_n\rangle) 
\subseteq \langle \varphi^*(u_1), \dots, \varphi^*(u_n)\rangle .
\end{equation}
Moreover, if $\varphi$ is a quasi-isomorphism, then \eqref{eq:massey-func} 
holds as equality.
\end{proposition}

In particular, if $\langle u_1, \dots, u_n\rangle$ is defined, then 
$\langle \varphi^*(u_1), \dots, \varphi^*(u_n)\rangle$ is also defined; 
and if, in addition, $\langle \varphi^*(u_1), \dots, \varphi^*(u_n)\rangle$ is 
non-vanishing, then $\langle u_1, \dots, u_n\rangle$ 
is also non-vanishing. As another consequence, the following holds: 
if $A\simeq B$, then all Massey products in $H^*(A)$ vanish 
if and only if all Massey products in $H^*(B)$ vanish.

Finally, if the map $\varphi\colon A\to B$ is a $q$-quasi-isomorphism, 
for some $q\ge 0$, then \eqref{eq:massey-func} holds as equality in 
degrees up to $q+1$. Thus, if $A\simeq_{q} B$, then all Massey 
products in $H^{\le q+1}(A)$ vanish if and only if all Massey 
products in $H^{\le q+1}(B)$ vanish.

The vanishing of Massey products provides a well known 
obstruction to formality. An analogous statement 
holds for partial formality, see e.g.~\cite{Su-sullivan}.

\begin{proposition}
\label{prop:massey-formal}
Let $(A,d)$ be a $\k$-$\cdga$.  If $A$ is formal, then all Massey 
products in $H^{*}(A)$ vanish. Furthermore, if $A$ is $q$-formal, 
for some $q\ge 1$, then all Massey products in $H^{\le q+1}(A)$ vanish.
\end{proposition}

\begin{proof}
Suppose $A$ is formal, so that $A\simeq (H^{*}(A),0)$.  By the consequence of
Proposition~\ref{prop:massey-natural} recorded above, all Massey products in
$H^{*}(A)$ vanish precisely when all Massey products in the cohomology of
$(H^{*}(A),0)$ do.  The latter all vanish.  Indeed, let
$\langle u_1,\dots,u_n\rangle$ be defined in a $\cdga$ with zero differential,
and set $a_{i-1,i}=u_i$ for $1\le i\le n$ and $a_{i,j}=0$ for $j>i+1$.  The
conditions $d a_{i,j}=\sum_{i<r<j}\bar{a}_{i,r}a_{r,j}$ then reduce to
$\bar{u}_{i+1}u_{i+2}=0$, which holds because these products are coboundaries
in any defining system, hence zero when $d=0$.  This is a defining system, and
the corresponding element \eqref{eq:massey-prod} is $\alpha=0$.  The
$q$-formal case follows in the same way from the preceding paragraph, which
gives the equivalence in degrees $\le q+1$.
\end{proof}

In general, formality is stronger than the mere vanishing 
of all Massey products; it is equivalent to their {\em uniform}\/ vanishing, 
in the sense of Deligne--Griffiths--Morgan--Sullivan \cite[Thm.~4.1]{DGMS}: 
a single choice, made once and for all in a minimal model, must render the 
forms representing all Massey products---of every order---exact 
simultaneously, whereas vanishing of the individual products supplies such 
choices only one product at a time. 
In the coming sections we develop alternative obstructions to 
formality via resonance varieties of Koszul complexes and 
Hilbert series of Koszul modules, which in many cases detect 
non-formality more effectively than individual Massey products.

%%%%%%%%%%%%%%%%%%%%%%%%%%%
\section{The Koszul construction}
\label{sect:koszul}
%%%%%%%%%%%%%%%%%%%%%%%%%%%

We introduce the Koszul construction associated to a connected $\cdga$ 
$(A,d_A)$ with $\dim_\k H^1(A)<\infty$. This assigns to $A$ two canonically 
dual complexes---one cohomological, one homological---over polynomial 
rings encoding complementary aspects of $H^1(A)$. Their (co)homology 
defines the Koszul modules that play a central role throughout the paper.

%%%%%%%%%%%%%%%%%%%%
\subsection{The two coefficient rings: $S$ and $T$}
\label{subsec:S-and-T}
%%%%%%%%%%%%%%%%%%%%

Throughout this section, $A=(A^\bullet,d_A)$ denotes a connected $\cdga$
over a field $\k$ of characteristic $0$, and we assume $0< \dim_\k H^1(A)<\infty$.

Since $A$ is connected, the differential $d_A \colon A^0 \to A^1$ vanishes, 
and so we may identify $H^1(A) = Z^1(A) = \ker(d_A \colon A^1 \to A^2)$.
Set $E = \bwedge(H^1(A))$ and write $E^1 = H^1(A)$ and $E_1 = (E^1)^\vee$ 
for the dual vector space. 

We define the \emph{cohomological coefficient ring} as
\[
S \coloneqq \Sym(E_1) = \Sym(H_1(A)),
\]
the symmetric algebra on $E_1 = H_1(A)$, and the \emph{homological 
coefficient ring} as
\[
T \coloneqq \Sym(E^1) = \Sym(H^1(A)),
\]
the symmetric algebra on $E^1 = H^1(A)$.
The natural perfect pairing $E^1 \otimes_{\k} E_1 \to \k$ identifies $T$ with 
the graded dual of $S$; since $\k$ has characteristic $0$, that dual is 
again a polynomial ring. Thus $S$ and $T$ represent two complementary 
perspectives on one object, reflecting Eckmann--Hilton duality between 
cohomological and homological viewpoints.

Fix a basis $\{e_1, \ldots, e_n\}$ for $E^1 = H^1(A)$, and let 
$\{x_1, \ldots, x_n\}$ be the dual basis for $E_1 = H_1(A)$. 
Then $S = \k[x_1, \ldots, x_n]$ and $T = \k[e_1, \ldots, e_n]$ 
as polynomial rings, and $x_j \mapsto e_j$ defines an isomorphism 
$S \cong T$ of \emph{graded}\/ $\k$-algebras.
Let $\m = (x_1,\dots,x_n) \subset S$ and 
$\n = (e_1,\dots,e_n) \subset T$ be the respective augmentation ideals.

This isomorphism depends on the chosen basis; the correspondence of 
augmentation ideals that it induces does not.  Indeed, $\m$ and $\n$ are 
intrinsically the ideals of elements of positive degree, so \emph{every}\/ 
graded algebra isomorphism $S\cong T$ carries $\m^p$ onto $\n^p$, for every 
$p\ge0$.  The identification of the $\m$-adic and $\n$-adic filtrations used 
throughout the paper is therefore canonical, even though the identification 
of the two rings is not.

Both $S$ and $T$ carry natural graded Hopf algebra structures, 
with primitive generators in degree $1$: $S$ via the generators $x_j$ 
and $T$ via the generators $e_j$.
Under the isomorphism $S \cong T$, these Hopf algebra structures are 
identified, but the cohomological and homological gradings on $S$ and $T$ 
remain distinct, reflecting the dual roles they play in the Koszul 
construction developed in the following subsections.

%%%%%%%%%%%%%%%%%%%%
\subsection{The canonical element and the cohomological Koszul complex}
\label{subsec:canonical-element}
%%%%%%%%%%%%%%%%%%%%

Given a finite-dimensional $\k$-vector space $V$, we define the 
corresponding \emph{canonical element} to be the tensor 
$\omega_V \in V^\vee \otimes_{\k} V$ which corresponds to the identity 
automorphism of $V^\vee$ under the tensor-hom adjunction. 
In concrete terms, if we pick a basis $\{v_1, \ldots, v_m\}$ for $V^\vee$ 
and let $\{v_1^{\vee}, \ldots, v_m^{\vee}\}$ be the dual basis for $V$, then 
$\omega_V = \sum_{i=1}^m v_i \otimes v_i^{\vee}$.

For our connected $\cdga$ $(A,d_A)$ with $0<\dim_\k H^1(A)<\infty$, we set
\begin{equation}
\label{eq:omega-A}
\omega_A \coloneqq \omega_{H_1(A)} \in H^1(A) \otimes_{\k} H_1(A). 
\end{equation}
Since $A$ is connected, we have $H^1(A) = Z^1(A)$, and thus we may view 
$\omega_A$ as a (nonzero) element of $Z^1(A) \otimes_{\k} H_1(A)$. Therefore, 
left multiplication by $\omega_A = \sum_{j=1}^n e_j \otimes x_j$, 
defines an endomorphism of $A \otimes_{\k} S$ of bidegree $(1,1)$.  

Consider the map  
\begin{equation}
\label{eq:delta-a}
\delta_A \colon A \otimes_{\k} S \longrightarrow A \otimes_{\k} S, 
\quad \delta_A = \omega_A \cdot (-) + d_A \otimes \id_S.
\end{equation}
Then:
\begin{itemize}[itemsep=2pt]
  \item $\delta_A$ is $S$-linear.
  \item $\delta_A$ is a graded derivation of the bigraded algebra $A \otimes_{\k} S$, 
        of total degree $+1$.
  \item $\delta_A$ is bidegree-homogeneous only when $d_A=0$:
        the term $\omega_A \cdot (-)$ has bidegree $(1,1)$ while 
        $d_A\otimes\id_S$ has bidegree $(1,0)$.
\end{itemize}

\begin{proposition}
\label{prop:LA}
Let $(A^{*}, d_A)$ be a connected $\k$-$\cdga$ with 
$0<\dim_\k H^1(A) < \infty$. Then
\[
K^{\bullet}(A) \coloneqq (A^{\bullet} \otimes_{\k} S,\delta_A)
\]
is a cochain complex of free $S$-modules,
\[
\begin{tikzcd}[column sep=20pt]
 \cdots \ar[r] 
& A^{i}\otimes_{\k} S \ar[r, "\delta^{i}_A"] 
&[4pt] A^{i+1} \otimes_{\k} S \ar[r, "\delta^{i+1}_A"] 
&[4pt] A^{i+2} \otimes_{\k} S \ar[r] 
& \cdots ,
\end{tikzcd}
\]
and $K^{\bullet}(A)$ is naturally a bigraded $\k$-$\cdga$.
\end{proposition}

\begin{proof}
We verify directly that $\delta_A^2=0$.  Write
$\delta_A = \omega_A\cdot(-) + d_A\otimes \id_S$. Then
\[
\delta_A^2
= (\omega_A\cdot(-))^2
+ (\omega_A\cdot(-))\circ(d_A\otimes \id_S)
+ (d_A\otimes \id_S)\circ(\omega_A\cdot(-))
+ (d_A\otimes \id_S)^2 .
\]
The first and last terms vanish. 
Indeed, $\omega_A^2=0$ since $e_ie_j=-e_je_i$ in $A$ while
$x_ix_j=x_jx_i$ in $S$, and $d_A^2=0$ by assumption.

For the mixed terms, let $a\otimes s\in A\otimes_{\k} S$. Since
$\omega_A\in Z^1(A)\otimes_{\k} H_1(A)$ consists of cocycles, 
graded commutativity gives
\[
(d_A\otimes \id_S)\big(\omega_A\cdot(a\otimes s)\big)
= -\,\omega_A\cdot(d_A(a)\otimes s).
\]
Thus,
\[
(\omega_A\cdot(-))\circ(d_A\otimes \id_S) + 
(d_A\otimes \id_S)\circ(\omega_A\cdot(-)) = 0.
\]

Combining all four terms, we conclude that $\delta_A^2=0$.
Since $\delta_A$ is a graded derivation of total degree $+1$ 
on the bigraded algebra $A\otimes_{\k} S$, it follows that
$K^\bullet(A)=(A\otimes_{\k} S,\delta_A)$ is a bigraded $\cdga$.
\end{proof}

The differentials 
$\delta^i_A\colon A^i\otimes_{\k} S\to A^{i+1}\otimes_{\k} S$ 
take the explicit form
\begin{equation}
\label{eq:diff}
\delta^{i}_A(a \otimes s)= 
\sum_{j=1}^{n} e_j a \otimes x_j s + d_A(a) \otimes s ,
\end{equation}
for all $a\in A^{i}$ and $s\in S$.  
(Each $e_j$ is a cocycle in $A^1$, so $e_j a \in A^{i+1}$.)  
In particular, since $A^0=\k$ and $d_A(1)=0$, the map 
$\delta^0_A\colon S\to A^1\otimes_{\k} S$ is given by 
\begin{equation}
\label{eq:delta0}
\delta^0_A(1)=\sum_{j=1}^n e_j \otimes x_j .
\end{equation}

If the $\cdga$ $A$ has zero differential, then each $\delta^{i}_{A}$ 
is represented by a matrix with entries linear in the variables $x_1,\dots ,x_n$.  
In general, the matrices for $\delta_A^i$ may also have nonzero 
constant terms.

\begin{definition}
\label{def:aomoto}
Let $(A^\bullet,d_A)$ be a connected $\k$-$\cdga$ with
$\dim_\k H^1(A)<\infty$.  For $a\in H^1(A)$, the \emph{Aomoto complex of $A$
at $a$} is the cochain complex $(A,\delta_a)$ with the same underlying graded
vector space as $A$, and with differentials
$\delta^i(a)\colon A^i\to A^{i+1}$ given by
\[
\delta^i(a)(u) = a\cdot u + d_A(u) .
\]
That $\delta_a^2=0$ follows as in the proof of Proposition~\ref{prop:LA}, from
$a^2=0$ and $d_A(a)=0$. 
\end{definition}

The relationship between the Koszul cochain complex and the Aomoto complexes 
is given by the following  lemma, see e.g.~\cite{Su-indam, Su-bockres}. 
For completeness, we include a proof. 

\begin{lemma}
\label{lem:two aom}
The specialization of the cochain complex $K^{\bullet}(A)=(A\otimes_{\k} S,\delta_A)$ 
at an element $a\in H^1(A)$ coincides with the cochain complex $(A,\delta_{a})$. 
\end{lemma}

\begin{proof}
Let $a=\sum_{j=1}^{n} a_j e_j\in Z^1(A)=H^1(A)$, and let 
$\m_a=(x_1-a_1,\dots , x_n-a_n)$ be the corresponding maximal ideal. 
Specializing $K^\bullet(A)$ via the evaluation map
$\ev_a\colon S\to S/\m_a=\k$ yields the complex
$A(a)=K^{\bullet}(A)\otimes_{S} S/\m_a$, with differential induced by
$\delta_A$; the underlying $\cdga$ $A$ is not itself an $S$-module.
Specializing $\delta_A$ via the evaluation map
$\ev_a\colon S\to \k$ yields the differential
$\delta^i(a)(u)=a\cdot u+d_A(u)$, 
which agrees with Definition~\ref{def:aomoto}.
\end{proof}

Lemma~\ref{lem:two aom} shows that the Aomoto complexes of $A$ arise as the
fibers of the Koszul complex $K^\bullet(A)$ over points of $\mSpec S=H^1(A)$.
Consequently, the cohomology of $K^\bullet(A)$ governs, in a global and
algebraic manner, the behavior of the Aomoto cohomology groups
$H^i(A,\delta_a)$ as $a$ varies.

\begin{definition}
\label{def:Koszul-coho}
The \emph{Koszul cohomology modules} of $A$ are
\[
\B^i(A) \coloneqq H^i(K^\bullet(A)).
\]
Since $\delta^0_A(s)=\sum_{j}e_j\otimes x_js$ and the $e_j$ are linearly
independent in $A^1$, the map $\delta^0_A$ is injective; hence
$\B^0(A) = 0$.
The differential $\delta_A=d_A\otimes\id+\omega_A$ is the sum of a part of
$S$-degree $0$ and a part of $S$-degree $1$, so it is not homogeneous
unless $d_A=0$; accordingly, the modules $\B^i(A)$ are in general only
\emph{filtered}\/ by the $\m$-adic filtration of
Section~\ref{subsec:koszul-filtrations}, and need not admit a grading
compatible with their $S$-module structure---see
Example~\ref{ex:sol2-koszul}.  They are graded when $d_A=0$, and more
generally when $A$ has positive weights
\textup{(}Proposition~\ref{prop:koszul-separation}\textup{)}.
\end{definition}

%%%%%%%%%%%%%%%%%%%%
\subsection{The homological Koszul complex}
\label{subsec:hom-koszul}
%%%%%%%%%%%%%%%%%%%%

We now construct the \emph{Koszul chain complex} of $(A,d_A)$ as the 
$S$-linear dual of the cochain complex $K^{\bullet}(A)$ from 
Proposition~\ref{prop:LA}, with coefficients transported to $T$ via 
the graded isomorphism $S \cong T$ fixed in Section~\ref{subsec:S-and-T}. 

\begin{definition}
\label{def:koszul-chain}
The \emph{Koszul chain complex} of $(A,d_A)$ is the chain complex of free 
$T$-modules,
\[
K_{\bullet}(A) = (A_{\bullet} \otimes_{\k} T, \partial^A) \colon 
\begin{tikzcd}[column sep=22pt]
\cdots \arrow[r] & A_2 \otimes_{\k} T \arrow[r, "\partial_2^A"] & 
A_1 \otimes_{\k} T \arrow[r, "\partial_1^A"] & 
A_0 \otimes_{\k} T = T,
\end{tikzcd}
\]
where $A_i = (A^i)^\vee$ and $\partial_i^A = (\delta_A^{i-1})^\vee$ is the 
$S$-linear dual of $\delta_A^{i-1}$, transported to a $T$-linear map via 
the identification $S\cong T$.
\end{definition}

In terms of the pairing $\langle -, - \rangle \colon A^{i-1} \times A_{i-1} \to \k$, 
the map $\partial_i^A$ is characterized by
\begin{equation}
\label{eq:duality-pairing}
\langle b \otimes t,\; \partial_i^A(a^\vee \otimes s) \rangle 
= \langle \delta_A^{i-1}(b \otimes t),\; a^\vee \otimes s \rangle,
\end{equation}
for all $b \in A^{i-1}$, $a^\vee \in A_i$, $t \in T$, $s \in S$, where 
$\langle b \otimes t,\, a^\vee \otimes s \rangle = \langle b, a^\vee \rangle \cdot ts$.

Equivalently, we have
\begin{equation}
\label{eq:partial-i-A}
\partial_i^A = \omega_A^\vee \lrcorner (-) + d_A^\vee \otimes \id_T 
\colon A_i \otimes_{\k} T \longrightarrow A_{i-1} \otimes_{\k} T,
\end{equation}
where $\omega_A^\vee \in E_1 \otimes_{\k} E^1$ is the dual of the 
canonical element. Note that $(\omega_A^\vee)^2 = (\omega_A^2)^\vee = 0$.

To make this more concrete, recall that we fixed a basis $\{e_1, \ldots, e_n\}$ 
for $E^1 = H^1(A)$ with dual basis $\{x_1, \ldots, x_n\}$ for $E_1 = H_1(A)$. 
The dual canonical element is then
\[
\omega_A^\vee = \sum_{j=1}^n x_j \otimes e_j \in E_1 \otimes_{\k} E^1.
\]
The contraction operation $\omega_A^\vee \lrcorner (-) \colon 
A_i \otimes_{\k} T \to A_{i-1} \otimes_{\k} T$ is the $T$-linear map given by
\[
\omega_A^\vee \lrcorner (a^\vee \otimes t) 
= \sum_{j=1}^n (e_j \lrcorner a^\vee) \otimes e_j t,
\]
for $a^\vee \in A_i$ and $t \in T$, where the contraction 
$e_j \lrcorner a^\vee \in A_{i-1}$ is defined by 
\begin{equation}
\label{eq:contraction-def}
\langle e_j a, a^\vee \rangle = \langle a, e_j \lrcorner a^\vee \rangle
\end{equation}
for all $a \in A^{i-1}$. 

In the chosen basis, the boundary maps 
$\partial_i^A \colon A_i \otimes_{\k} T \to A_{i-1} \otimes_{\k} T$ 
take the explicit form
\begin{equation}
\label{eq:partial-explicit}
\partial_i^A(a^\vee \otimes t) = \sum_{j=1}^n 
(e_j \lrcorner a^\vee) \otimes e_j t  + (d_A^\vee a^\vee) \otimes t ,
\end{equation}
for all $a^\vee \in A_i$ and $t \in T$.
Here $d_A^\vee \colon A_i \to A_{i-1}$ denotes the differential dual to $d_A$,
and $e_j \lrcorner a^\vee$ is the contraction operator from 
\eqref{eq:contraction-def}.

Note that the generator $e_j$ occurs twice in \eqref{eq:partial-explicit}, in
two different roles: it contracts $a^\vee$, and it multiplies the coefficient
in $T$.  This is a consequence of writing the complex over $T$.  Written
instead over $S$, the same boundary map reads
\[
\partial_i^A(a^\vee \otimes s) = \sum_{j=1}^n
(e_j \lrcorner a^\vee) \otimes x_j s + (d_A^\vee a^\vee) \otimes s ,
\qquad s\in S,
\]
and now each of the two factors of $\omega_A^\vee = \sum_{j} x_j
\otimes_{\k} e_j$ is used exactly once: the $H^1(A)$-factor $e_j$ contracts
$a^\vee$, while the $H_1(A)$-factor $x_j$ multiplies the coefficient.  The two
formulas correspond under the identification $x_j \mapsto e_j$ fixed in
Section~\ref{subsec:S-and-T}.

\begin{definition}
\label{def:Koszul-homology}
The \emph{Koszul homology modules}\/ of $A$ are
\begin{equation*}
\label{eq:koszul-hom}
\B_i(A) \coloneqq H_i(K_\bullet(A)).
\end{equation*}
\end{definition}

These are $T$-comodules, with coaction induced by the Hopf 
comultiplication on $T$; as in the cohomological case, they are in general 
only filtered, and are graded when $d_A=0$ or, more generally, when $A$ has 
positive weights.
Via the identification $S\cong T$ fixed in Section~\ref{subsec:S-and-T},
we will freely regard $\B_i(A)$ as an $S$-module when discussing
graded pieces or Hilbert series, keeping in mind that its natural
structure is that of a $T$-comodule.

Since the image of $\partial_1^A\colon A_1\otimes_{\k} T \to A_0\otimes_{\k} T=T$ 
coincides with the maximal ideal $\n=(e_1,\dots,e_n)$, we have that 
$\B_0(A)=T/\n\cong \k$; in particular, $\B_0(A)\ne \B^0(A)$. 
 
\begin{example}
\label{ex:classical-koszul}
The most basic example is when 
$A$ is the exterior algebra $E = \bwedge(e_1,\dots,e_n)$ with $d=0$ and 
$S = \k[x_1,\dots,x_n]$.  Then $K_{\bullet}(E)=(E_\bullet \otimes_{\k} S, \partial^E)$ 
is the classical Koszul complex, with differentials
\[
\partial_i^E(e_{j_1}\wedge\cdots\wedge e_{j_i}) = \sum_{k=1}^{i} (-1)^{k-1}
e_{j_1}\wedge\cdots\wedge \widehat{e_{j_k}}\wedge\cdots\wedge e_{j_i}
\otimes x_{j_k}.
\]
Augmenting by $\varepsilon\colon S\to\k$, $\varepsilon(x_i)=0$, 
gives a free $S$-resolution of $\k$; hence, $\B_i(E)=0$ for all $i>0$.
\end{example}

%%%%%%%%%%%%%%%%%%%%%%%
\subsection{A tale of two Koszul complexes}
\label{subsec:ka=kha}
%%%%%%%%%%%%%%%%%%%%%%%

Having defined the Koszul chain complex 
\[
K_\bullet(A)= (A_\bullet \otimes_{\k} T , \partial^{A}),
\]
it is often useful to compare it with its ``linearized'' counterpart,
obtained by replacing the $\cdga$ $(A,d_A)$ with its cohomology algebra 
$H^\ast(A)$ equipped with zero differential:
\[
K_\bullet(H^\ast(A)) = \big(H_\bullet(A)\otimes_{\k} T , \partial^{H^\ast(A)}\big).
\]
The differential in $K_\bullet(A)$ decomposes as
$\partial^{A} = \omega_A^{\vee} + d_A^\vee$, 
where $\omega_A^{\vee}$ is contraction with the dual 
canonical element and $d_A^\vee$ is the dual differential.  
By contrast, $K_\bullet(H^\ast(A))$ retains only the linear component
$\omega_A^{\vee}$. Consequently, the Koszul homology modules 
\[
\B_i(A)=H_i(K_\bullet(A))
\quad\text{and}\quad 
\B_i(H^\ast(A))=H_i(K_\bullet(H^\ast(A)))
\]
may differ dramatically in positive degrees.

This discrepancy appears already in very low degrees.
For the Chevalley--Eilenberg $\cdga$ $A=\CE(\h(1))$ of the 
$3$-dimensional Heisenberg Lie algebra (see 
Section~\ref{subsec:heisenberg-1}), one has
$\B_1(A)\cong \k$, whereas $\B_1\bigl(H^\ast(A)\bigr)\cong S$ 
(as $S$-modules). Thus, although both $\cdgas$ have the same 
first cohomology, their Koszul homology in degree~$1$ differs maximally. 
This illustrates a general phenomenon:  
the homology of $K_\bullet(A)$ is sensitive to the full differential $d_A$,
while that of $K_\bullet(H^\ast(A))$ detects only the resonance operator.
The comparison between these two complexes---and the spectral sequences
linking them---will play a central role in the next section.

The two Koszul complexes of $A$ are related by $S$-duality.  For every
connected $\cdga$ $A$ with $\dim_\k H^1(A)<\infty$, the natural pairing
$A^i\otimes_{\k} A_i\to\k$ extends to a perfect $S$-bilinear pairing under
which $\partial^A$ is the transpose of $\delta_A$; equivalently,
\begin{equation}
\label{eq:EH-pairing}
K_\bullet(A) \cong \Hom_S\bigl(K^\bullet(A),S\bigr),
\end{equation}
so that $K_\bullet(A)$ is the Eckmann--Hilton dual of $K^\bullet(A)$.
Two dualities are in play here, and they should not be merged: the $S$-dual
$\Hom_S(-,S)$ of \eqref{eq:EH-pairing} returns coefficients in $S$, whereas
the graded $\k$-linear \textup{(}Kronecker\textup{)} dual returns
coefficients in the graded dual $S^\vee\cong T$, which is the form in which
$K_\bullet(A)$ was defined above.  Under the identification of
Section~\ref{subsec:S-and-T} the two complexes agree; the distinction between
the two dualities becomes essential only on \textup{(}co\textup{)}homology,
as we note next.  The relation \eqref{eq:EH-pairing} is already encoded in the
explicit formulas \eqref{eq:diff} and \eqref{eq:partial-explicit}.

It does not, however, descend to (co)homology: the functor
$\Hom_S(-,S)$ is not exact, so $\B_i(A)$ is not in general the $S$-dual
of $\B^i(A)$, nor its Kronecker dual.  A genuine comparison between
$\B_\bullet(A)$ and $\B^\bullet(A)$ requires additional input, supplied
by Poincar\'{e} duality through the orientation class rather than by the
pairing \eqref{eq:EH-pairing}; see
Theorem~\ref{thm:pd-koszul-modules}.  Because $\B_i(A)$ interacts more
directly with Alexander-type invariants and their filtrations, we work
primarily with Koszul homology in what follows.

%%%%%%%%%%%%%%%%%%%%%%%%%%%%%%%%%%%%
\subsection{Filtrations on Koszul modules}
\label{subsec:koszul-filtrations}
%%%%%%%%%%%%%%%%%%%%%%%%%%%%%%%%%%%%

We now introduce the natural filtrations on the Koszul (co)ho\-mology 
modules and address the fundamental question of their separation, 
which will play a recurring role throughout the paper.

Since several filtrations occur side by side---the $\m$-adic and $\n$-adic
ones, the weight filtration $W$, the filtration $F$ induced by a Koszul
spectral sequence, the $I$-adic filtration on a group ring, the lower
central series filtration $\gamma_\bullet$, and a general filtration
$\cF$ of a Lie algebra---we adopt throughout the convention
that the \emph{filtration}\/ appears as a superscript on $\gr$ and the
\emph{degree}\/ as a subscript.  Thus $\gr^\m M=\boplus_{p\ge0}\gr^\m_pM$
with $\gr^\m_pM=\m^pM/\m^{p+1}M$, and likewise $\gr^F_p$, $\gr^W_p$,
$\gr^{\gamma}_n$; the notation $\gr_n$ without a superscript refers to the
filtration in force at that point, which for Lie algebras and groups is
the lower central series $\gamma_\bullet$
\textup{(}see \eqref{eq:gr-lcs-convention}\textup{)}.  This
keeps expressions such as $\gr^\m_n\B_1(A)$ and $\gr^F_p\B_i(A)$
unambiguous.

The cohomological Koszul $S$-modules $\B^i(A) = H^i(K^\bullet(A))$ carry 
a natural descending $\m$-adic filtration, inherited from the 
$\m$-adic filtration on $K^\bullet(A)$:
\[
\B^i(A) \supseteq \m\B^i(A) \supseteq \m^2\B^i(A) \supseteq \cdots .
\]
Similarly, the homological Koszul $T$-modules $\B_i(A) = H_i(K_\bullet(A))$ 
carry a natural descending $\n$-adic filtration:
\[
\B_i(A) \supseteq \n\B_i(A) \supseteq \n^2\B_i(A) \supseteq \cdots .
\]
By Section~\ref{subsec:S-and-T}, every graded isomorphism $S\cong T$ carries 
$\m^p$ onto $\n^p$, so these two filtrations correspond canonically.  We use 
this correspondence to regard $\B_i(A)$ as a filtered $S$-module---which is in 
any case what the support and resonance loci of 
Section~\ref{sect:res-cdga} require---and accordingly write $\m$ for the 
augmentation ideal and $\gr^\m$ for the associated graded of both types of 
Koszul module.  The symbol $\m$ is retained only for the filtration 
$F^pK_\bullet(A)=A_\bullet\otimes_{\k}\n^p$ of the chain complex itself, where 
$T$ is genuinely the ring of coefficients.

If $A$ is $q$-finite, then for each $i\le q$ both $\B^i(A)$ and 
$\B_i(A)$ are finitely generated $S$-modules: they are quotients of 
submodules of the finitely generated free modules 
$A^i\otimes_{\k} S$ and $A_i\otimes_{\k} T$, and $S$ is Noetherian.
The separation of their filtrations---that is, whether 
$\bigcap_{k\ge 0}\m^k M = 0$---is a genuine condition, in the 
cohomological case just as much as in the homological one: it fails 
precisely when some nonzero element of the module is annihilated by a 
polynomial that does not vanish at the origin.

\begin{proposition}
\label{prop:separation-criterion}
Let $M$ be a finitely generated $S$-module, and let $M_\m$ denote its 
localization at the maximal ideal $\m$.  Then
\begin{equation}
\label{eq:sep-criterion}
\bigcap_{k\ge 0}\m^k M
=\ker\bigl(M\to M_\m\bigr)
=\bigl\{u\in M \mid su=0 \ \text{for some}\ s\in S\setminus\m\bigr\}.
\end{equation}
Consequently, the $\m$-adic filtration on $M$ is separated if and only if 
every associated prime of $M$ is contained in $\m$.  This holds, in 
particular, when $M$ has finite length, and when $M$ is a graded 
$S$-module.
\end{proposition}

\begin{proof}
Set $N=\bigcap_{k\ge 0}\m^kM$.  Each $M/\m^kM$ is a module over the 
$\m$-local ring $S/\m^k$, so the canonical maps 
$M/\m^kM\to M_\m/\m^kM_\m$ are isomorphisms, and $N$ is the kernel of 
the composite $M\to M_\m\to\varprojlim_{k} M_\m/\m^kM_\m$.  Since $M_\m$ 
is a finitely generated module over the Noetherian \emph{local}\/ ring 
$S_\m$, the Krull intersection theorem gives 
$\bigcap_{k}\m^kM_\m=0$, so the second map is injective and 
$N=\ker(M\to M_\m)$; this kernel is the set of elements whose 
annihilator meets $S\setminus\m$, which is the second equality 
in~\eqref{eq:sep-criterion}.

If $u\ne0$ lies in $N$, then any associated prime of $Su$---and hence 
of $M$---contains $\Ann(u)$, which meets $S\setminus\m$, so that prime 
is not contained in $\m$.  Conversely, if some associated prime $\p$ of 
$M$ satisfies $\p\not\subseteq\m$, write $\p=\Ann(u)$ with $u\ne0$ and 
choose $s\in\p\setminus\m$; then $su=0$, so $N\ne0$.  Finally, if $M$ 
has finite length, then $\m$ is its only associated prime; and if $M$ is 
graded, then every associated prime of $M$ is homogeneous, hence 
contained in the ideal $\m$ of elements of positive degree.
\end{proof}

Both filtrations can fail to be separated, and already do so for the 
same $\cdga$: for $A=\CE(\sol_2)$, treated in 
Section~\ref{subsec:sol2} below, one has 
$\B_1(A)\cong\k[x]/(x-1)$ and $\B^2(A)\cong\k[x]/(x-1)$, and in either 
case $x-1$ is a unit, so $\m$ acts invertibly, 
$\m^p\B=\B$ for every $p$, and the intersection 
in~\eqref{eq:sep-criterion} is the whole module.  Note that the Krull 
intersection theorem does not apply directly to $S$ itself: the ideal 
$\m$ is not contained in the Jacobson radical of $S$, which is zero.
Examples~\ref{ex:coho-nonseparated} and~\ref{ex:sol2-ss} record the 
consequences for the cohomological and the homological Koszul spectral 
sequence, respectively.

Separation is guaranteed in two important cases: when $A$ has positive 
weights, in which case the Koszul modules are graded 
(Proposition~\ref{prop:koszul-separation}), and when the module in 
question has finite length---for instance, when its support is 
contained in the origin.  Observe also that, by 
\eqref{eq:sep-criterion}, separation always holds after localizing at 
$\m$, and that $\gr^\m M=\gr^\m(M_\m)$ for every finitely generated 
$M$.  This is why those of our results that are stated at the level of 
associated graded modules, or of germs and tangent cones at the origin, 
require no separation hypothesis, whereas those that recover the modules 
themselves do.
Separation plays a role in spectral sequence convergence 
(Theorem~\ref{thm:koszul-ss-hom}\eqref{ssh4}), in isomorphism 
theorems under quasi-isomorphisms 
(Theorem~\ref{thm:q-iso-koszul}\eqref{iso1}), and in the 
comparison with Alexander invariants developed in 
Section~\ref{sect:equivariant-ss}.

%%%%%%%%%%%%%%%%%%%%
\subsection{A non-homogeneous Koszul complex}
\label{subsec:sol2}
%%%%%%%%%%%%%%%%%%%%

We conclude with a basic example showing that a nonzero linear component
in the differential can destroy homogeneity of the Koszul complex and
lead to Koszul homology that depends on the chosen $\cdga$ model,
even in a formal, degree-$1$-generated case.  This example recurs
throughout the paper; we record here, once and for all, the data that
will be needed later.

\begin{example}
\label{ex:sol2-koszul}
Let $A=\CE(\sol_2)=\bigl(\bigwedge\langle a,b\rangle,d\bigr)$ be the
Chevalley--Eilenberg algebra of the $2$-dimensional solvable,
non-nilpotent Lie algebra $\sol_2=\Lie(x,y)/([x,y]-y)$, so that
\begin{equation}
\label{eq:sol2-cdga}
da=0,\qquad db=-a\wedge b.
\end{equation}
Since $db$ has a nonzero linear part, $H^1(A)=\k\langle a\rangle$ and
$H^*(A)\cong\bigwedge(a)$.  Writing $S=T=\Sym(H_1(A))=\k[x]$ with
$x=a^\vee$, and $\omega_A=a\otimes x$, the Koszul chain complex
$K_\bullet(A)=(A_\bullet\otimes_{\k}S,\partial^A)$ and its
Eckmann--Hilton dual \eqref{eq:EH-pairing} take the form
\begin{equation}
\label{eq:toy}
K_\bullet(A)\colon\quad \begin{tikzcd}[column sep=32pt]
S \ar[r, "\sbm{ 0 \\ x-1}"] &S^2 \ar[r, "\sbm{ x \amp 0}"] &S ,
\end{tikzcd}
\end{equation}
\begin{equation}
\label{eq:toy-coh}
K^\bullet(A)\colon\quad \begin{tikzcd}[column sep=32pt]
S \ar[r, "\sbm{ x \\ 0}"]
&S^2 \ar[r, "\sbm{ 0 \amp x-1}"]
&S .
\end{tikzcd}
\end{equation}
A direct computation yields
\begin{equation}
\label{eq:sol2-modules}
\begin{gathered}
\B_0(A)\cong \k[x]/(x),\qquad
\B_1(A)\cong \k[x]/(x-1),\qquad
\B_i(A)=0 \quad (i\ge 2),\\
\B^1(A)\cong \k[x]/(x),\qquad
\B^2(A)\cong \k[x]/(x-1),\qquad
\B^i(A)=0 \quad (i\ne 1,2).
\end{gathered}
\end{equation}
By contrast, for the cohomology algebra $H^*(A)\cong\bigwedge(a)$---and
likewise for the sub-$\cdga$ $B=\bigwedge(a)\subset A$, whose Koszul
complexes are the truncations of \eqref{eq:toy} and
\eqref{eq:toy-coh} in which the generator $b$ is absent---one has
\begin{equation}
\label{eq:sol2-modules-coh}
\B_0\cong \k[x]/(x),\qquad
\B^1\cong \k[x]/(x),\qquad
\text{and all other Koszul modules vanish,}
\end{equation}
by Example~\ref{ex:classical-koszul} or by direct inspection.
In particular, $\B_1(A)\not\cong\B_1\bigl(H^*(A)\bigr)$ and
$\B^2(A)\not\cong\B^2\bigl(H^*(A)\bigr)$.

The ideal $(x)$ is homogeneous, and thus $\B_0(A)$ inherits a natural
grading from~$S$, as is always the case for $\B_0$.
By contrast, the relation $x=1$ is non-homogeneous, so neither
$\B_1(A)$ nor $\B^2(A)$ admits a grading compatible with its
$S$-module structure.  This failure of homogeneity is already visible at
the chain level: the boundary map $\partial^A_2$ contains a nonzero
constant term, so that, while the $\m$-adic filtration
$F^p=A_\bullet\otimes_{\k}\m^p$ is still a filtration by subcomplexes, the
standard grading on $S$ does not make $K_\bullet(A)$ a graded complex.
Finally, $A$ admits no positive-weight decomposition: weight-preservation
of $d$ would force $\wt(b)=\wt(db)=\wt(a)+\wt(b)$, which is impossible
for $\wt(a)>0$.
\end{example}

\begin{remark}
\label{rem:ce-sol2}
The $\cdga$ $A=\CE(\sol_2)$ is formal: the inclusion 
$H^*(A)\cong\bigwedge(a)\hookrightarrow A$ is a quasi-isomorphism 
of $\cdgas$. Thus $A\simeq H^*(A)$, yet $\B_1(A)\cong\k[x]/(x-1)$ 
while $\B_1(H^*(A))=0$. This shows that Koszul modules are 
\emph{not} quasi-isomorphism invariants: two weakly equivalent 
$\cdgas$ may have very different Koszul homology. 
The discrepancy arises because the Koszul differential 
$\partial^A = \omega_A^\vee\lrcorner(-) + d_A^\vee$ depends on 
the full differential $d_A$ of the $\cdga$ model, not merely on 
the induced maps on cohomology.
\end{remark}

%%%%%%%%%%%%%%%%%%%%%%%%
\subsection{Algebraic structures on Koszul modules}
\label{subsec:koszul-algebra-structure}
%%%%%%%%%%%%%%%%%%%%%%%%

The Koszul cochain complex $K^\bullet(A)$ is a dg-algebra over $S$ 
(Proposition~\ref{prop:LA}), and the Koszul chain complex $K_\bullet(A)$ 
is its dual dg-coalgebra over $T$. These structures pass to (co)homology, 
equipping the Koszul modules with additional algebraic structure beyond 
their bare $S$-module and $T$-comodule structures.

\begin{proposition}
\label{prop:koszul-alg-coalg}
Let $(A,d_A)$ be a connected, finite-type $\k$-$\cdga$. 
\begin{enumerate}[itemsep=2pt]
\item \label{ka1}
The cohomological Koszul modules $\B^\bullet(A) = H^\bullet(K^\bullet(A))$ 
form a graded-commu\-tative $S$-algebra, with product induced by 
the dg-algebra structure on $K^\bullet(A)$.
\item \label{ka2}
The homological Koszul modules $\B_\bullet(A) = H_\bullet(K_\bullet(A))$ 
form a graded-cocom\-mu\-tative $T$-coalgebra, with coproduct induced 
by the dg-coalgebra structure on $K_\bullet(A)$.
\end{enumerate}
\end{proposition}

\begin{proof}
Both statements follow from the fact that $\delta_A$ is a derivation with
respect to the product on $K^\bullet(A)$ \textup{(}Proposition~\ref{prop:LA}\textup{)},
so that the product descends to cohomology; dually, $\partial^A$ is a
coderivation for the coproduct on $K_\bullet(A)$.
The graded-commutativity of the induced product on $\B^\bullet(A)$ 
follows from the graded-commuta\-tivity of $K^\bullet(A)$ as a 
dg-algebra.
\end{proof}

In addition to the $S$-algebra structure on $\B^\bullet(A)$ and 
the $T$-coalgebra structure on $\B_\bullet(A)$, the cohomology 
algebra $H^*(A)$ acts on both Koszul modules via cap and cup 
products.

\begin{proposition}
\label{prop:hstar-action}
Let $(A,d_A)$ be a connected, finite-type $\k$-$\cdga$. 
The cohomology algebra $H^*(A)$ acts on both Koszul modules:
\begin{enumerate}[itemsep=2pt]
\item \label{ha1}
\emph{Cap product action on $\B_\bullet(A)$:} there is a left 
action
\[
\cap\colon H^p(A)\otimes_{\k} \B_i(A) \longrightarrow \B_{i-p}(A),
\qquad [\alpha]\cap[-] = [\alpha\lrcorner(-)],
\]
where $\alpha\lrcorner(-)$ denotes contraction of the $A$-factor 
by $\alpha\in Z^p(A)$. This action is compatible with the 
$T$-comodule structure on $\B_\bullet(A)$.
\item \label{ha2}
\emph{Cup product action on $\B^\bullet(A)$:} there is a left 
action
\[
\cup\colon H^p(A)\otimes_{\k} \B^i(A) \longrightarrow \B^{i+p}(A),
\qquad [\alpha]\cup[-] = [\alpha\wedge(-)],
\]
where $\alpha\wedge(-)$ denotes multiplication by $\alpha$ in
the $A$-factor of $K^\bullet(A)$. This action is compatible
with the $S$-algebra structure on $\B^\bullet(A)$.
\end{enumerate}
Moreover, the two actions are adjoint with respect to the $S$-bilinear
pairing $K_\bullet(A)\cong\Hom_S(K^\bullet(A),S)$ of \eqref{eq:EH-pairing}.
\end{proposition}

\begin{proof}
\eqref{ha1} Contraction by a cocycle $\alpha\in Z^p(A)$ 
defines a chain map 
$\alpha\lrcorner(-)\colon K_\bullet(A)\to K_{\bullet-p}(A)$ 
of degree $-p$, commuting with $\partial^A$ by the Leibniz rule 
and the cocycle condition. Passing to homology gives the cap 
product.

\eqref{ha2} Multiplication by $\alpha\in Z^p(A)$ in the 
$A$-factor defines a cochain map 
$\alpha\wedge(-)\colon K^\bullet(A)\to K^{\bullet+p}(A)$ 
of degree $+p$, commuting with $\delta_A$ by the dg-algebra 
structure on $K^\bullet(A)$. Passing to cohomology gives the 
cup product.  The adjointness of the two actions is immediate from
\eqref{eq:EH-pairing}, under which $\alpha\lrcorner(-)$ is the transpose of
$\alpha\wedge(-)$.  Since $\Hom_S(-,S)$ is not exact, this adjointness at the
level of complexes does not in general descend to a duality between
$\B_\bullet(A)$ and $\B^\bullet(A)$.
\end{proof}

\begin{remark}
\label{rem:koszul-alg-open-revised}
When $b_1(A)=1$, the cohomological Koszul modules $\B^i(A)$ are cyclic 
$\k[x]$-modules, and in all computed examples the nonzero modules 
$\B^i(A)$ and $\B^j(A)$ for $i\ne j$ have disjoint support in 
$\mSpec(\k[x])=\mathbb{A}^1$, forcing the cup product to vanish. 
Whether the multiplicative structure on $\B^\bullet(A)$ or the 
comultiplicative structure on $\B_\bullet(A)$ carries information 
beyond the $S$-module or $T$-comodule structures remains open.

The cohomology algebra $H^*(A)$ acts on both Koszul modules via 
cap and cup products. This action is used in the present paper 
primarily to identify the connecting morphism $\delta_e$ in the 
Koszul--Gysin sequence of a Hirsch extension as cap product 
with $[e]\in H^2(B)$ (Theorem~\ref{thm:koszul-hirsch-filtration} 
and Corollary~\ref{cor:supp-hirsch}). More generally, if $A$ carries 
a positive-weight decomposition (Section~\ref{subsec:weights-koszul}), 
then the $H^*(A)$-action preserves the total-weight filtration on 
$\B_\bullet(A)$ and induces the corresponding action on each associated 
graded piece. By Theorem~\ref{thm:koszul-weight-graded}, the associated
graded of $\B_\bullet(A)$ is a subquotient of $\B_\bullet(H^*(A))$, with
equality precisely when the weight spectral sequence degenerates at
$E^2$; this holds in particular when $A$ is formal
(Corollary~\ref{cor:koszul-completion-formal}), and in that case, for
weight-preserving morphisms (Section~\ref{subsec:wt-morphisms-bq}), the
compatibility implies multiplicative $\B_q$-formality.

Whether the $H^*(A)$-action carries deeper geometric or topological 
information beyond these formal properties---such as governing 
formality obstructions, controlling the decomposition of 
$\B_\bullet(A)$ along resonance components, or interacting with 
other filtrations---remains an interesting open question.
\end{remark}

%%%%%%%%%%%%%%
\subsection{K\"{u}nneth formulas for Koszul complexes and modules}
\label{subsec:koszul-kunneth}
%%%%%%%%%%%%%%

Let $(A,d_A)$ and $(B,d_B)$ be connected $\cdgas$ over $\k$ with
$\dim_\k H^1(A)<\infty$ and $\dim_\k H^1(B)<\infty$.
Set $S_A=\Sym(H_1(A))$, $S_B=\Sym(H_1(B))$, and let $C=A\otimes_{\k} B$
be the tensor product $\cdga$. Then $H^1(C)=H^1(A)\oplus H^1(B)$,
so $H_1(C)=H_1(A)\oplus H_1(B)$ and $S_C=\Sym(H_1(C))\cong S_A\otimes_{\k} S_B$
canonically as graded Hopf algebras.

\begin{theorem}
\label{thm:koszul-kunneth} 
There is a canonical isomorphism of dg-algebras
\[
K^\bullet(C) \cong K^\bullet(A) \otimes_{\k} K^\bullet(B),
\]
compatible with the $S_C$-module structures (via $S_C\cong S_A\otimes_{\k} S_B$).
Consequently, there is an isomorphism of dg-coalgebras
\[
K_\bullet(C) \cong K_\bullet(A) \otimes_{\k} K_\bullet(B),
\]
compatible with the natural $T_C$-comodule structures.
\end{theorem}

\begin{proof}
The graded vector space identification
\begin{equation}
\label{eq:csc-ab}
C\otimes_{\k} S_C = (A\otimes_{\k} B)\otimes_{\k} (S_A\otimes_{\k} S_B)
= (A\otimes_{\k} S_A)\otimes_{\k} (B\otimes_{\k} S_B) 
= K^\bullet(A)\otimes_{\k} K^\bullet(B)
\end{equation}
is immediate. The differential $\delta_C = d_C + \omega_C\cdot(-)$ 
expands to
\begin{equation}
\label{eq:dcdadb}
d_C = d_A\otimes 1 + 1\otimes d_B, \quad \omega_C\cdot (-) = 
\omega_A\cdot (-) + \omega_B\cdot (-),
\end{equation}
matching the standard tensor product differential.
This preserves the dg-algebra structure (multiplication on $A\otimes_{\k} S_A$ 
commutes with $\delta_A$). The chain version follows by graded duality.
\end{proof}

\begin{corollary}
\label{cor:koszul-kunneth-modules}
Let $(A,d_A)$ and $(B,d_B)$ be two connected, finite-type 
$\k$-$\cdgas$, and let $C=A\otimes_{\k} B$. Then for all $q\ge 0$:
\begin{enumerate}[itemsep=2pt]
\item \label{kk1}
$\B^q(C) \cong \bigoplus_{i+j=q} \B^i(A)\otimes_{\k} \B^j(B)$
as graded $S_C$-modules.
\item \label{kk2}
$\B_q(C) \cong \bigoplus_{i+j=q} \B_i(A)\otimes_{\k} \B_j(B)$
as graded $T_C$-comodules.
\end{enumerate}
\end{corollary}

\begin{proof}
By Theorem~\ref{thm:koszul-kunneth}, 
$K^\bullet(C)\cong K^\bullet(A)\otimes_{\k} K^\bullet(B)$
as dg-algebras over $S_C$.  The tensor product is taken over the field $\k$,
so every module in sight is flat and the K\"{u}nneth theorem holds in its
simplest form, as an isomorphism
$H^{r}\bigl(K^\bullet(A)\otimes_{\k}K^\bullet(B)\bigr)\cong
\bigoplus_{i+j=r}\B^i(A)\otimes_{\k}\B^j(B)$; no correction term arises.  This
gives~\eqref{kk1}, and the isomorphism is one of $S_C$-modules because the
identification of Theorem~\ref{thm:koszul-kunneth} is $S_C$-linear.
Part~\eqref{kk2} follows by graded duality, using 
$\B_0(A)\cong\B_0(B)\cong\k$.
\end{proof}

\begin{remark}
\label{rem:kunneth-alg-coalg}
The Künneth isomorphisms of Corollary~\ref{cor:koszul-kunneth-modules} 
are compatible with the algebraic structures of 
Proposition~\ref{prop:koszul-alg-coalg}: the algebra structure 
on $\B^\bullet(A\otimes_{\k} B)$ is the tensor product of those 
on $\B^\bullet(A)$ and $\B^\bullet(B)$, and similarly for the 
coalgebra structures. In particular, these structures carry no 
additional information beyond the $S$-module and $T$-comodule 
structures in the tensor product case.
\end{remark}

As we show in Section~\ref{subsec:tensor-coprod}, these K\"{u}nneth
isomorphisms determine both the jump loci $\RR^{q,s}$ and the support loci 
$\RR_{q,s}$ of a tensor product of $\cdgas$, in all degrees and depths, 
through the fiber dimensions of the two factors---though, for $s\ge2$, not 
as a union of products of the loci of the factors.

%%%%%%%%%%%%%%
\subsection{Koszul complexes of coproducts of $\cdgas$}
\label{subsec:coprod}
%%%%%%%%%%%%%%

Let $(A,d_A)$ and $(B,d_B)$ be two connected $\cdgas$ with
$\dim_\k H^1(A), \dim_\k H^1(B)<\infty$.  Their coproduct, $C=A\vee B$, is a 
new connected $\cdga$, whose underlying graded vector space in 
positive degrees is the direct sum $A^+\oplus B^+$, whose multiplication 
is given by $(a,b)\cdot (a',b') = (aa', bb')$, and whose differential is 
$d_{A\vee B}=d_A+d_B$. Note that $C^1=A^1\oplus B^1$; 
therefore, we may identify $H^1(A\vee B)$ with $H^1(A)\oplus H^1(B)$.

Set $S_A=\Sym(H_1(A))$, $S_B=\Sym(H_1(B))$, and 
$S_C=\Sym(H_1(C))\cong S_A\otimes_{\k} S_B$.
Since $C^q=A^q\oplus B^q$ for $q\ge1$ while $C^0=\k$, the Koszul cochain
complex of $C$ decomposes in positive degrees as
\begin{equation}
\label{eq:kq-aveeb}
K^q(A\vee B)\cong
\bigl(K^q(A)\otimes_{S_A}S_C\bigr)\oplus\bigl(K^q(B)\otimes_{S_B}S_C\bigr)
\qquad (q\ge1),
\end{equation}
with differential $\delta_A\oplus\delta_B$ in degrees $\ge1$: for $u\in A^{\ge1}$
one has $\omega_C\cdot u=\omega_A\cdot u$, since $A^{\ge1}\cdot B^{\ge1}=0$ in
$C$.  In degree $0$, by contrast, $K^0(C)=S_C$ receives only one copy of $S_C$
rather than two.  The two projections $C\to A$ and $C\to B$ therefore package
the comparison as a short exact sequence of cochain complexes,
\begin{equation}
\label{eq:kakb-otimes-ses}
\begin{tikzcd}[column sep=18pt]
0 \arrow[r] & K^\bullet(A\vee B)
\arrow[r] & \bigl(K^\bullet(A)\oplus K^\bullet(B)\bigr)\otimes_{\k} S_C
\arrow[r] & S_C[0]
\arrow[r] & 0,
\end{tikzcd}
\end{equation}
where $S_C[0]$ denotes a copy of $S_C$ concentrated in cohomological degree
$0$, and the base changes are along $S_A\to S_C$ and $S_B\to S_C$.
The extra copy of $S_C$ reflects the fact that the coproduct
has a single unit in degree $0$, while its degree-$1$ part is the direct sum
of those of $A$ and $B$, with no mixed products between them.

Dually, there is a short exact sequence of chain complexes
\begin{equation}
\label{eq:kakb-vee-ses}
\begin{tikzcd}[column sep=18pt]
0 \arrow[r] & S_C[0]
\arrow[r] & \bigl(K_\bullet(A)\oplus K_\bullet(B)\bigr)\otimes_{\k} S_C
\arrow[r] & K_\bullet(A\vee B)
\arrow[r] & 0,
\end{tikzcd}
\end{equation}
where $S_C[0]$ is concentrated in homological degree $0$, embedded
antidiagonally.  Concretely,
\begin{equation}
\label{eq:kq-aveeb-bis}
K_q(A\vee B)\cong
\bigl(K_q(A)\otimes_{S_A}S_C\bigr)\oplus\bigl(K_q(B)\otimes_{S_B}S_C\bigr)
\qquad (q\ge1),
\end{equation}
with boundary maps that restrict to those of $A$ and of $B$; in degree $1$,
both summands map to the same copy of $K_0(A\vee B)=S_C$.
As a consequence, we obtain the following description of the Koszul modules 
of a coproduct.

\begin{proposition}
\label{prop:koszul-coprod}
Let $C=A\vee B$, and write $\m_A,\m_B\subset S_C$ for the ideals generated by
$H_1(A)$ and by $H_1(B)$, respectively.  Then:
\begin{enumerate}[itemsep=2pt]
\item \label{kc1}
$\B_q(C)\cong \bigl(\B_q(A)\otimes_{S_A}S_C\bigr)\oplus
\bigl(\B_q(B)\otimes_{S_B}S_C\bigr)$, for all $q\ge 2$.
\item \label{kc2}
There is a short exact sequence
\[
\begin{tikzcd}[column sep=16pt]
0 \arrow[r] & \bigl(\B_1(A)\otimes_{S_A}S_C\bigr)\oplus
\bigl(\B_1(B)\otimes_{S_B}S_C\bigr) \arrow[r] & \B_1(C) \arrow[r] 
& \m_A\m_B \arrow[r] & 0.
\end{tikzcd}
\]
\end{enumerate}
\end{proposition}

\begin{proof}
Part~\eqref{kc1} follows from the decomposition 
\eqref{eq:kq-aveeb-bis}: for $q \ge 2$ the boundary maps of 
$K_\bullet(A \vee B)$ preserve the two summands, so homology splits 
accordingly; base change along the flat extensions $S_A\to S_C$ and 
$S_B\to S_C$ commutes with taking homology.

For part~\eqref{kc2}, note that $\partial_1^C$ restricts to $\partial_1^A$ and
to $\partial_1^B$ on the two summands of $K_1(C)$, both mapping into the same
copy of $S_C$; and that $\im\partial^C_2=\im\partial^A_2\oplus\im\partial^B_2$.
Hence there is an exact sequence
\[
0\longrightarrow \ker\partial^A_1\oplus\ker\partial^B_1
\longrightarrow \ker\partial^C_1
\longrightarrow \im\partial^A_1\cap\im\partial^B_1\longrightarrow 0 ,
\]
the third map being $(u,v)\mapsto\partial^A_1(u)=-\partial^B_1(v)$.  Now
$\im\partial^A_1=\m_A$ and $\im\partial^B_1=\m_B$; since $\m_A$ and $\m_B$ are
generated by disjoint sets of variables, $\m_A\cap\m_B=\m_A\m_B$.  Dividing by
$\im\partial^A_2\oplus\im\partial^B_2$ gives the stated sequence.
\end{proof}

\begin{remark}
\label{rem:coprod-extra-term}
The third term in \eqref{kc2} is isomorphic to $S_C$ precisely when
$b_1(A)=b_1(B)=1$; in general $\m_A\m_B$ requires $b_1(A)b_1(B)$ generators.
For instance, if $A=H^{*}(S^1;\k)$ and $B=H^{*}(S^1\vee S^1;\k)$, then
$\B_1(A)=0$ and $\B_1(B)\cong S_B$, while $C=A\vee B$ satisfies
$\B_1(C)=\ker\bigl(S_C^3\xrightarrow{(x,y_1,y_2)}S_C\bigr)$, which is
minimally generated by the three Koszul syzygies of $(x,y_1,y_2)$.  This is
consistent with \eqref{kc2}, whose outer terms need $0+1$ and $2$ generators,
respectively.
\end{remark}

We will see in Section~\ref{subsec:tensor-coprod} that these exact sequences
give rise to the corresponding coproduct formulas for support resonance varieties.

%%%%%%%%%%%%%%%%%%%%%%%%%%%%%%%%%%%%%
\section{Koszul complexes and Poincar\'{e} duality}
\label{sect:koszul-pd}
%%%%%%%%%%%%%%%%%%%%%%%%%%%%%%%%%%%%%

We analyze how Poincaré duality on a $\cdga$ interacts with the Koszul 
construction of Section~\ref{sect:koszul}. The key object is the \emph{PD-adjoint} 
differential, which differs from a naive matrix transpose and governs the 
relationship between the homological and cohomological Koszul complexes.

%%%%%%%%%%%%%%%%%%%%%
\subsection{Poincar\'{e} duality $\cdgas$}
\label{subsec:pd-cdga}
%%%%%%%%%%%%%%%%%%%%%

Throughout this subsection, $A=(A^\bullet,d)$ denotes a connected $\cdga$ over
a field $\k$ of characteristic $0$. We assume $A$ is of finite type when needed.

\begin{definition}
\label{def:pda}
A connected, finite-type graded-commutative $\k$-algebra
$A=\bigoplus_{i\ge 0} A^i$ is called a \emph{Poincar\'{e} duality algebra of formal
dimension $m$}\/ (briefly, a $\PD_m$-algebra) if there exists a linear functional
$\varepsilon\colon A^m \to \k$ such that, for all $0\le i\le m$, the bilinear
pairings
\begin{equation*}
\label{eq:PD-bilinear}
A^i \otimes_{\k} A^{m-i} \longrightarrow \k,
\qquad
a\otimes b \longmapsto \varepsilon(ab),
\end{equation*}
are nondegenerate.
\end{definition}

It follows that $A^i=0$ for $i>m$ and that multiplication induces canonical
isomorphisms of $\k$-vector spaces
\begin{equation}
\label{eq:pd-iso}
\PD^i\colon  A^i \longisom  (A^{m-i})^{\vee},
\qquad
\PD^i(a)(b)=\varepsilon(ab).
\end{equation}

Consequently, each element $a\in A^i$ has a {\em Poincar\'{e} 
dual}, $a^{*}\in A^{m-i}$, which is uniquely determined by the 
formula $\varepsilon (a a^{*})=1$. 
The \emph{orientation class} of $A$ is the element $\nu_A\in A^m$ uniquely
determined by $\varepsilon(\nu_A)=1$; equivalently, $\nu_A$ is the Poincar\'{e} dual
of $1\in A^0$.

\begin{definition}
\label{def:pdcdga}
A connected $\k$-$\cdga$ $(A^\bullet,d)$ is called a
\emph{Poincar\'{e} duality $\cdga$ of formal dimension $m$}\/ (briefly, a
$\PD_m$-$\cdga$) if:
\begin{enumerate}[itemsep=1pt]
\item the underlying graded algebra $A^\bullet$ is a $\PD_m$-algebra;
\item $d(A^{m-1})=0$, equivalently, $H^m(A)\cong \k$, with orientation induced by
$\varepsilon$.
\end{enumerate}
\end{definition}

If $(A,d)$ is a $\PD_m$-$\cdga$, then $H^{*}(A)$ is a $\PD_m$-algebra, 
with orientation induced by $\varepsilon$. The converse holds under 
suitable hypotheses: by a theorem of Lambrechts--Stanley \cite{LS-asens}, 
if $H^1(A)=0$ and $H^*(A)$ is a $\PD_m$-algebra, then $A$ is weakly 
equivalent to a $\PD_m$-$\cdga$.

Associated to a $\PD_m$ algebra over a field 
$\k$ there is an alternating $m$-form,
\begin{equation}
\label{eq:mu form}
\mu_A\colon \bwedge^m A^1 \longrightarrow \k, \quad
\mu_A(a_1\wedge \cdots \wedge a_m)= \varepsilon (a_1\cdots a_m).
\end{equation} 
This form encodes the top-degree multiplication in $A$ and will play a 
central role in the description of Koszul differentials in low formal dimension.

When $m=3$, the multiplicative structure of the underlying graded algebra of $A$
can be recovered from the $3$-form $\mu=\mu_A$ and the orientation $\varepsilon$,
as follows. 
Fix a basis $\{e_1,\dots ,e_n\}$ for $A^1$.  
Let $\{e_1^{*},\dots ,e_n^{*}\}$ be the Poincar\'{e} dual basis for $A^2$, 
and take as generator of $A^3=\k$ the class $\nu=1^{*}$.   
The multiplication in $A$, then, is given on basis elements by 
\begin{equation}
\label{eq:mult}
e_i  e_j=\sum_{k=1}^{n} \mu_{ijk}\,  e^{*}_k, \quad 
e_i e_j^{*} = \delta_{ij} \nu,
\end{equation} 
where $\mu_{ijk}=\mu(e_i\wedge e_j\wedge e_k )$ and $\delta_{ij}$ 
is the Kronecker delta.   An alternate way to encode this information 
is to let $A_i=(A^i)^{\vee}$ be the dual $\k$-vector space and to let $e_i^{\vee}\in A_1$ 
be the (Kronecker) dual of $e_i$.  Under the canonical identification between 
$(\bigwedge^3 A^1)^{\vee}$ and $\bigwedge^3 A_1$, we may view 
$\mu_A$ dually as a trivector, 
\begin{equation}
\label{eq:trivec}
\mu_A =\sum \mu_{ijk}\, e_i^{\vee} \wedge e_j^{\vee} \wedge e_k^{\vee} \in \bwedge^3 A_1.
\end{equation} 

%%%%%%%%%%%%%%%%%%%%
\subsection{Poincar\'{e} duality and Aomoto complexes}
\label{subsec:pd-aomoto}
%%%%%%%%%%%%%%%%%%%%

Let $(A^\bullet,d)$ be a $\PD_m$-$\cdga$ over a field $\k$ of characteristic $0$,
with orientation class $\nu_A\in A^m$. For each element $a\in A^1$, recall that
the associated Aomoto differential on $A^\bullet$ is given by
\begin{equation}
\label{eq:delta-ax}
\delta_a(x)=ax+d(x), \qquad x\in A^\bullet .
\end{equation}

Poincar\'{e} duality relates the Aomoto differentials associated to $a$ and $-a$ 
in complementary degrees. The following lemma appears in various forms in 
\cite{PS-imrn19, Su-edinb, Su-bockres}. For completeness, we include a proof.

\begin{lemma}
\label{lem:cd-uptosign}
Let $(A,d)$ be a $\PD_m$-$\cdga$ over a field $\k$ of characteristic $0$.  Then, 
for all $a\in A^1$ and all $0\le i\le m$, we have a commuting square,
\[
\begin{tikzcd}[column sep=42pt, row sep=28pt]
A^i \ar[r,"\delta_a^i"] \ar[d,"\PD^i"'] &
A^{i+1} \ar[d,"(-1)^{i+1}\PD^{i+1}"] \\
(A^{m-i})^{\vee} \ar[r,"(\delta_{-a}^{m-i-1})^{\vee}"] &
(A^{m-i-1})^{\vee} .
\end{tikzcd}
\]
\end{lemma}

\begin{proof}
Let $b \in A^{i}$ and $c \in A^{m-i-1}$.  Then
\begin{align*}
(-1)^{i+1}  \PD^{i+1}\circ \delta^{i}_{a}(b) (c)  
&=(-1)^{i+1}  \PD^{i+1}( a b + d(b))  (c) \\
&=(-1)^{i+1} \varepsilon(abc +d(b)c),
\end{align*}
while 
\begin{align*}
(\delta^{m-i-1}_{-a})^{\vee}\circ \PD^{i} (b) (c)
&=\PD^i (b) ( \delta^{m-i-1}_{-a} (c) )  \\
&=\PD^i(b) (-a c + d(c)) \\
&=\varepsilon(-bac +bd(c)). 
\end{align*}

Now, since $A$ is a $\cdga$, we have that $ab=(-1)^{i}ba$ 
and $d(b) c+(-1)^{i} b d(c)=d(b c)$. 
Moreover, since $A$ is a $\PD_m$-$\cdga$ and $bc\in A^{m-1}$, 
we have that $d(bc)=0$, and the claim follows.
\end{proof}

Lemma~\ref{lem:cd-uptosign} motivates the following terminology.
For $a \in A^1$ and $0 \le i \le m-1$, we say that 
$\delta_{-a}^{m-i-1}$ is the \emph{Poincaré-dual adjoint} of 
$\delta_a^i$, and write
\begin{equation}
\label{eq:dagger}
(\delta_a^i)^{\dagger_{\PD}} \coloneqq \delta_{-a}^{m-i-1}.
\end{equation}
In particular, the Aomoto differential $\delta_a$ is self-adjoint 
up to the involution $a \mapsto -a$ and degree reversal.

\begin{proposition}
\label{prop:pd-aomoto}
Let $(A^\bullet,d)$ be a $\PD_m$-$\cdga$. Then, for every $a\in A^1$ 
and every $i\ge 0$, there is a natural isomorphism
\[
H^i(A,\delta_a)^{\vee} \cong H^{m-i}(A,\delta_{-a}) .
\]
\end{proposition}

\begin{proof}
By the identity \eqref{eq:dagger}, Poincar\'{e}
duality identifies the cochain complex $(A^\bullet,\delta_a)$ with the graded
dual of $(A^{m-\bullet},\delta_{-a})$, up to the standard degree shifts and
Koszul signs. Passing to cohomology yields the asserted natural isomorphism.
\end{proof}

The involution $a\mapsto -a$ in Proposition~\ref{prop:pd-aomoto} 
is a chain-level artifact of graded commutativity: differentials 
in complementary degrees are related by PD-adjunction, not by 
naive transposition. For algebraic invariants insensitive to 
individual elements of $A^1$---such as Koszul modules and 
resonance varieties---this involution disappears, and Poincar\'{e} 
duality reduces to a palindromic symmetry in the formal 
dimension $m$.

%%%%%%%%%%%%%%%%%%%%%%%%%%%%
\subsection{Poincar\'{e} duality for Koszul homology}
\label{subsec:pd-koszul}
%%%%%%%%%%%%%%%%%%%%%%%%%%%%

We explain how Poincar\'{e} duality interacts with the Koszul complexes
introduced in Section~\ref{sect:koszul}.  The duality is induced by the
orientation class, and so, unlike the pairing of
Section~\ref{subsec:ka=kha}, it descends to Koszul (co)homology; it
requires no assumption that the graded algebra $A$ be generated in
degree~$1$.

Throughout this subsection, let $(A^\bullet,d)$ be a connected $\cdga$ over
$\k$, with $\dim_\k H^1(A)<\infty$, and set $S=\Sym(H_1(A))$.
Recall that the Koszul cochain complex of $A$ is
\[
K^\bullet(A)=(A^\bullet\otimes_{\k} S,\delta_A), \qquad
\delta_A=d+\omega_A,
\]
where $\omega_A\in Z^1(A)\otimes_{\k} H_1(A)$ is the canonical tensor.
Its cohomology is $\B^\bullet(A)$, and its $S$-dual $\Hom_S(-,S)$ is the
Koszul chain complex $K_\bullet(A)$, with homology $\B_\bullet(A)$.

Let $\sigma\colon S\to S$ be the $\k$-algebra involution determined by
$\sigma(x)=-x$ for $x\in H_1(A)$; on $\mSpec(S)=H^1(A)$ it is the
antipodal map $a\mapsto -a$.  For an $S$-module $M$ we write
$\sigma^{*}M$ for $M$ with its action twisted by $\sigma$, so that
$\Supp_S(\sigma^{*}M)=-\Supp_S(M)$.  Any subset of $H^1(A)$ stable under
scalar multiplication---in particular, the support of a graded
$S$-module---satisfies $-Z=Z$, in which case the twist is invisible.

Assume now that $A$ is a $\PD_m$-$\cdga$, with orientation class
$\nu_A\in A^m$.  Poincar\'{e} duality provides $\k$-linear isomorphisms
$\PD^i \colon A^i \isom (A^{m-i})^\vee=A_{m-i}$, which extend
$S$-linearly to isomorphisms of free graded $S$-modules
\begin{equation}
\label{eq:asa-amis}
K^i(A)=A^i\otimes_{\k} S \longisom A_{m-i}\otimes_{\k} S=K_{m-i}(A).
\end{equation}
These are isomorphisms of $S$-modules, not dualizations: the Kronecker
dual $(A^{m-i})^\vee$ is applied to the finite-dimensional vector space
$A^{m-i}$, before tensoring with $S$.  What relates $K^\bullet(A)$ to
$K_\bullet(A)$ by duality is the $S$-linear pairing
$K_\bullet(A)\cong\Hom_S(K^\bullet(A),S)$, valid for any connected
finite-type $\cdga$; the two mechanisms should not be conflated, since
$\Hom_S(-,S)$ does not commute with passing to (co)homology, whereas
\eqref{eq:asa-amis} does.

The following result establishes \eqref{eq:pd-bia} from the Introduction.

\begin{theorem}
\label{thm:pd-koszul-modules}
Let $(A^\bullet,d)$ be a $\PD_m$-$\cdga$ of finite type.  Then
\eqref{eq:asa-amis} induces an isomorphism of complexes of graded
$S$-modules
\[
K^{\bullet}(A) \longisom  \sigma^{*}K_{m-\bullet}(A),
\]
and hence, for each $i\ge 0$, a natural isomorphism of $S$-modules
\[
\B_i(A) \cong \sigma^{*}\B^{\,m-i}(A).
\]
\end{theorem}

\begin{proof}
By Lemma~\ref{lem:cd-uptosign}, for every $a\in A^1$ the maps $\PD^i$
intertwine $\delta^i_a$ with $(\delta^{m-i-1}_{-a})^{\vee}$ up to the
sign $(-1)^{i+1}$.  Replacing $\PD^i$ by $\epsilon_i\PD^i$, where the
signs $\epsilon_i$ are chosen recursively so that
$\epsilon_{i+1}(-1)^{i+1}=\epsilon_i$, we may assume that the squares
commute on the nose.

Now let $a$ vary universally.  The Koszul cochain differential is
$\delta_A=d+\omega_A$, and $\partial^A$ is its $S$-linear transpose, so
$(\delta^{m-i-1}_{-a})^{\vee}$ globalizes to the component
$\partial^A_{m-i}$ with $\omega_A$ replaced by $-\omega_A$.  Since
$\omega_A\mapsto-\omega_A$ is induced by $\sigma$ on the coefficient
ring, the $S$-linear extension of $\PD$ identifies $K^\bullet(A)$ with
$\sigma^{*}K_{m-\bullet}(A)$.  Passing to (co)homology gives the second
assertion.
\end{proof}

\begin{remark}
\label{rem:pd-not-kronecker}
The isomorphism of Theorem~\ref{thm:pd-koszul-modules} is not a duality
between $\B_i(A)$ and $\B^{m-i}(A)$, and no such duality holds.  Take
$A=H^{*}\bigl(\#^{n}(S^1\times S^2);\k\bigr)$, a $\PD_3$-algebra with
zero differential in which all products of degree-one classes vanish.
Writing $S=\k[x_1,\dots,x_n]$, the Koszul complexes give
\[
\B_1(A)=\B^{2}(A)=\ker\bigl(S^n\to S\bigr),
\qquad
\B_2(A)=\B^{1}(A)=\coker\bigl(S\to S^n\bigr),
\]
in accordance with the theorem; for $n=2$ the first of these is free of
rank one.  The Kronecker dual $\bigl(\B^{2}(A)\bigr)^{\vee}$, by
contrast, is not even finitely generated over $S$.
\end{remark}

\begin{corollary}
\label{cor:pd-res-symmetry}
Let $(A,d)$ be a $\PD_m$-$\cdga$ of finite type. Then for all
$i \ge 0$ and $s \ge 1$,
\[
\RR_{i,s}(A) = -\Supp_S\bigl(\bwedge^s\B^{\,m-i}(A)\bigr).
\]
If the support loci involved are conical, the antipodal twist may be
omitted.
\end{corollary}

\begin{proof}
Exterior powers commute with the twist, so
$\bwedge^s\B_i(A)\cong\sigma^{*}\bwedge^s\B^{m-i}(A)$ by
Theorem~\ref{thm:pd-koszul-modules}, and
$\Supp_S(\sigma^{*}M)=-\Supp_S(M)$.
\end{proof}

\begin{remark}
\label{rem:pd-vs-exterior}
There are two distinct mechanisms that give rise to duality between 
Koszul homology and cohomology modules of a $\cdga$ $(A,d)$.

\begin{enumerate}[itemsep=2pt]
\item If $A=\bigwedge V$ with $\dim_\k V=N$, duality is induced by the top
exterior power $\bigwedge^N V$ and yields a shift by $N$.

\item If $A$ is a $\PD_m$-$\cdga$, duality is induced by the orientation class
$\nu_A\in A^m$ and yields a shift by the formal dimension $m$.
\end{enumerate}

When both conditions hold, the two shifts coincide. In general, however, the
relevant duality shift is governed by the Poincar\'{e} duality dimension rather
than by $\dim_\k H^1(A)$ or by the number of generators in degree~$1$.
\end{remark}

%%%%%%%%%%%%%%%%%%%%%%%%%%%%%%%
\subsection{The Koszul complexes of a $\PD_3$-$\cdga$}
\label{subsec:koszul-pd3}
%%%%%%%%%%%%%%%%%%%%%%%%%%%%%%%

The description below generalizes the explicit construction of the Koszul
complex $K^\bullet(A)$ for $\PD_3$ algebras with zero differential given in
\cite[Sec.~7.1]{Su-edinb}, by incorporating the internal differential
$d$ of a $\PD_3$-$\cdga$ $(A,d)$.

Let $(A,d)$ be a $\PD_3$-$\cdga$ over a field $\k$ of characteristic $0$.
Set $n=\dim_\k H^1(A)$, and choose a basis $\{e_1,\dots,e_n\}$ 
for $H^1(A)=Z^1(A)\subset A^1$. Extend this to a basis $\{e_1,\dots,e_r\}$ 
with $r\ge n$ of the full vector space $A^1$, and let 
$\{e_1^*,\dots,e_r^*\}$ be the corresponding Poincar\'{e} dual basis of $A^2$. 
Let $H_1(A)=(H^1(A))^{\vee}$, and identify
$S=\Sym(H_1(A))\cong \k[x_1,\dots,x_n]$, 
where $x_i$ is the (Kronecker) dual of $e_i$ for $1\le i\le n$.

Let $\nu\in A^3$ be the orientation class, and let
$\mu=\mu_A\in \bwedge^3 A_1$ be the associated alternating $3$-form, with structure
constants $\mu_{ijk}=\varepsilon(e_i e_j e_k)$, skew-symmetric in the indices.

Consider the Koszul cochain complex
\begin{equation}
\label{eq:koszul 3mfd}
K^{\bullet}(A)\colon
\begin{tikzcd}
A^0\otimes_{\k} S \ar[r, "\delta^{0}_A"]  &
A^1\otimes_{\k} S \ar[r, "\delta^{1}_A"]  &
A^2\otimes_{\k} S \ar[r, "\delta^{2}_A"]  &
A^3\otimes_{\k} S .
\end{tikzcd}
\end{equation}
In the Poincar\'{e} dual bases fixed above, the differentials are given by
\begin{itemize}[itemsep=3pt]
\item $\delta^0_A(1)=\sum_{j=1}^{n} e_j \otimes x_j$.

\item $\delta^1_A(e_i)=\sum_{j=1}^{n}\sum_{k=1}^{r} \mu_{jik}\, e_k^{*}\otimes x_j
+ d(e_i)\otimes 1$, for $1\le i\le r$, where $d(e_i)=0$ for $1\le i\le n$.

\item $\delta^2_A(e_i^{*}) =\nu\otimes x_i$ for $1\le i\le n$, and
$\delta^2_A(e_i^{*}) =0$ for $n+1\le i\le r$.
\end{itemize}

Observe that, with respect to the Poincar\'{e} dual bases fixed above, the first
and third differentials have matrices
\[
\delta^0_A=\big(x_1\ \cdots\ x_n\ 0\ \cdots\ 0\big),
\qquad
\delta^2_A=(\delta^0_A)^{\top}.
\]
This equality reflects the fact that $\delta^2_A$ is the
Poincar\'{e}-dual adjoint of $\delta^0_A$, and that in Poincar\'{e} dual bases
the adjoint is represented by the ordinary transpose.

Conceptually, the middle differential 
$\delta^1_A\colon A^1\otimes_{\k} S \to A^2\otimes_{\k} S$
decomposes canonically as
\begin{equation}
\label{eq:delta-1a}
\delta^1_A=\delta^1_{\mu}+\delta^1_{d},
\end{equation}
reflecting the two independent structures carried by a 
$\PD_3$-$\cdga$. The summand $\delta^1_\mu$ encodes the 
linear terms: for each basis element $e_i \in A^1$,
\[
\delta^1_\mu(e_i) = \sum_{j=1}^n \sum_{k=1}^r \mu_{jik}\, 
e_k^* \otimes x_j,
\]
where the coefficients $\mu_{jik} = \varepsilon(e_j e_i e_k)$ 
are skew-symmetric in $j$ and $k$.
The summand $\delta^1_d$ encodes the constant terms: 
$\delta^1_d(e_i) = d(e_i) \otimes 1$, which is nonzero 
only for $i > n$ (i.e., for basis elements outside $Z^1(A)$).
Thus $\delta^1_A$ is completely determined by the pair $(\mu, d)$.

The Koszul chain complex $K_{\bullet}(A)=(A_{\bullet}\otimes_{\k} T,\partial^A)$
is obtained from the cochain complex \eqref{eq:koszul 3mfd} by combining
$\k$-linear duality with Poincar\'{e} duality on $A$.
Since $A$ is a $\PD_3$-$\cdga$, Poincar\'{e} duality provides canonical
identifications $A_i \cong (A^{3-i})^{\vee}$ for $0\le i\le 3$. 
Using also the graded duality between $S=\Sym(H_1(A))$ and
$T=\Sym(H^1(A))$, we obtain natural identifications
\begin{equation}
\label{eq:pd-ast-iso}
(A^{3-i}\otimes_{\k} S)^{\vee} \cong A_i\otimes_{\k} T .
\end{equation}
Under these identifications, the boundary maps
\[
\partial^A_{i+1}\colon  A_{i+1}\otimes_{\k} T \longrightarrow A_i\otimes_{\k} T
\]
are the \emph{Poincar\'{e}-dual adjoints} of the cochain differentials
\[
\delta_A^{2-i}\colon  A^{2-i}\otimes_{\k}S \longrightarrow  A^{3-i}\otimes_{\k} S .
\]

In particular, the extremal maps are adjoint in the sense that
$\partial^A_1$ is dual to $\delta^0_A$ and $\partial^A_3$ is dual to
$\delta^2_A$. With respect to the Poincar\'{e} dual bases fixed above, this yields
\[
\partial^A_1=(\delta^0_A)^{\top},
\qquad
\partial^A_3=(\delta^2_A)^{\top},
\]
while the middle boundary map $\partial_2$ is the graded transpose of
$\delta^1_A$ after applying Poincar\'{e} duality and degree reversal. 
Thus, although the first and third differentials are related by duality,
$\partial_3$ is not the transpose of $\partial_1$ in general.
Rather, $\delta^2_A$ is the Poincar\'{e}-dual adjoint of $\delta^0_A$,
and this symmetry is what ultimately governs the shape of the Koszul
chain complex for a $\PD_3$-$\cdga$.

The passage from the Koszul cochain complex $K^{\bullet}(A)$ to the Koszul
chain complex $K_{\bullet}(A)$ follows the general Poincar\'{e}-duality
formalism of Section~\ref{subsec:pd-aomoto}.
Explicitly, the boundary maps $\partial^A_{i+1}$ are obtained from the
cochain differentials $\delta_A^{2-i}$ by first taking Poincar\'{e}-dual
adjoints in complementary degree and then applying graded $\k$-linear
duality.
In particular, while $\delta_A^2$ is the Poincar\'{e}-dual adjoint of
$\delta_A^0$, the corresponding chain maps $\partial^A_3$ and $\partial^A_1$
are not related by naive matrix transposition.

%%%%%%%%%%%%%%%%%%%%%%%%%
\subsection{The Heisenberg Lie algebras}
\label{subsec:heisenberg}
%%%%%%%%%%%%%%%%%%%%%%%%%

Let $\h(n)$ denote the $(2n+1)$-dimensional Heisenberg Lie algebra, with basis
$\{x_1,\dots,x_n,\, y_1,\dots,y_n,\, z\}$ and nontrivial brackets $[x_i,y_i]=z$.
Its Chevalley--Eilenberg algebra $A(n)=\CE(\h(n))$ is the $\cdga$
\[
A(n)=\bwedge(e_1,\dots,e_n,f_1,\dots,f_n,g),
\qquad
dg=-\sum_{i=1}^n e_i\wedge f_i,
\quad
de_i=df_i=0.
\]

\begin{lemma}
\label{lem:heisenberg-pd}
For each $n\ge1$, the Chevalley--Eilenberg algebra $A(n)=\CE(\h(n))$ is a
$\PD_{2n+1}$-$\cdga$. Moreover, its cohomology algebra $H^*(A(n))$ is a 
$\PD_{2n+1}$-algebra, even though it is not generated in degree~$1$.
\end{lemma}

\begin{proof}
Let $V^\vee=\spn\{e_1,\dots,e_n,f_1,\dots,f_n\}$ and set
$h=e_1\wedge f_1\wedge\cdots\wedge e_n\wedge f_n$.
Then $A(n)^{2n+1}=\k\cdot\nu$, where $\nu=h\wedge g$.

We first verify Poincar\'{e} duality at the level of graded algebras.
As a graded algebra,
$A(n)=\bwedge(V^\vee\oplus \k g)$, and multiplication induces pairings
\[
A(n)^i\otimes_{\k} A(n)^{2n+1-i} \longisom \k\cdot\nu .
\]
Writing $A(n)^k=\bwedge^k V^\vee \oplus g\wedge\bwedge^{k-1}V^\vee$,
wedging with $h$ identifies $\bwedge^k V^\vee$ with the dual of
$g\wedge\bwedge^{2n-k}V^\vee$, while wedging with $g$ identifies
$\bwedge^{k-1}V^\vee$ with the dual of $\bwedge^{2n+1-k}V^\vee$.
Thus multiplication by $\nu$ yields a nondegenerate pairing
$A(n)^i\otimes_{\k} A(n)^{2n+1-i}\to\k$ for all $i$.

It remains to check that $d(A(n)^{2n})=0$.
Since $dg=-\sum_{i} e_i\wedge f_i$, we have $d\nu=0$.
Moreover, $A(n)^{2n}$ is spanned by $h$ and by elements of the form
$g\wedge\alpha$, with $\alpha\in\bwedge^{2n-1}V^\vee$.
The class $h$ is closed, while
$d(g\wedge\alpha)=(dg)\wedge\alpha=0$ for degree reasons.
Hence $d(A(n)^{2n})=0$ and $H^{2n+1}(A(n))\cong\k$, proving 
that $A(n)$ is a $\PD_{2n+1}$-$\cdga$.

That $H^{*}(A(n))$ is then a $\PD_{2n+1}$-algebra, with orientation induced by
$[\nu]$, is the general fact recorded after Definition~\ref{def:pdcdga}.
\end{proof}

Geometrically, the $\PD_{2n+1}$-algebra structure of $H^*(A(n))$ reflects 
the fact that $A(n)$ models the cohomology of the compact Heisenberg nilmanifold 
of dimension $2n+1$.

\begin{remark}
\label{rem:macinic}
The Heisenberg $\cdgas$ $A(n)=\CE(\h(n))$ exhibit subtle formality properties.
It is shown by M\u{a}cinic~\cite{Mc10} that $A(n)$ is $(n-1)$-formal but not
$n$-formal, while $A(1)=\CE(\h(1))$ is not $1$-formal.
Thus, among the Heisenberg Lie algebras, $\h(1)$ is the only case in which
non-formality already appears at the first stage. We will come back to 
this subject in Section~\ref{subsec:nilpotent-groups}.
\end{remark}

%%%%%%%%%%%%%%%%%%%%%%%%%%
\subsection{The Koszul modules of $\CE(\h(1))$}
\label{subsec:heisenberg-1}
%%%%%%%%%%%%%%%%%%%%%%%%%%

Among the Heisenberg $\cdgas$, $\CE(\h(1))$ is the simplest example in which
the Koszul construction distinguishes a $\cdga$ model from its cohomology algebra,
despite agreement in $H^1$ and Poincaré duality.
We now focus on this case.

Let $A=\CE(\h(1))=\bwedge(e_1,f_1,g)$, with $dg=-e_1\wedge f_1$ and $de_1=df_1=0$.
Its cohomology algebra is
\[
H^\ast(A)\cong \frac{\bwedge(e_1,f_1,b_1,b_2)} {(e_1f_1,\; e_1b_2+f_1b_1)},
\]
where $b_1=[e_1\wedge g]$ and $b_2=[f_1\wedge g]$ arise as Massey products.
In particular, $H^1(A)=\langle e_1,f_1\rangle$, and we set
$S=\Sym(H_1(A))\cong \k[x_1,x_2]$, with $x_1,x_2$ dual to $e_1,f_1$.
Writing $\omega_A=e_1\otimes x_1+f_1\otimes x_2$, and writing the Koszul chain
complex with coefficients in $S$ as in Section~\ref{subsec:koszul-filtrations},
$K_\bullet(A)=(A_\bullet\otimes_{\k} S,\partial^A)$ has the form
\begin{equation*}
\label{eq:cc-heis-3}
\begin{tikzcd}[column sep=20pt]
0 \arrow[r] & S 
   \arrow[r, "\partial_3"] & S^3
   \arrow[r, "\partial_2"] & S^3
   \arrow[r, "\partial_1"] & S.
\end{tikzcd}
\end{equation*}
The boundary maps $\partial_i$ below are the transposes of the cochain
differentials $\delta_A^{i-1}$ of Section~\ref{subsec:koszul-pd3}, as in
\eqref{eq:partial-explicit}.  We compute them in the Kronecker-dual bases
$\{e_1^{\vee},f_1^{\vee},g^{\vee}\}$ of $A_1$ and
$\{(e_1f_1)^{\vee},(e_1g)^{\vee},(f_1g)^{\vee}\}$ of $A_2$, rather than in
the Poincar\'{e} dual bases used in Section~\ref{subsec:koszul-pd3}:
\[
\partial_1 = \begin{pmatrix} x_1 & x_2 & 0 \end{pmatrix},\qquad
\partial_2 = 
\begin{pmatrix}
  -x_2 & 0 & 0 \\
  x_1 & 0 & 0 \\
  -1   & x_1 & x_2
\end{pmatrix},\qquad
\partial_3 = 
\begin{pmatrix}
 0\\ -x_2\\ x_1
\end{pmatrix}.
\]
Note that $\mu_A$ restricts to zero on $\bwedge^3 H^1(A)$, since
$\dim_\k H^1(A)=2$; the linear part of $\delta_A^1$ therefore involves only the
generator $g\notin Z^1(A)$.  The differential component $d\colon A^1\to A^2$ is
nontrivial as well, and it is this that contributes the constant entry of
$\partial_2$.
As a result, the middle Koszul differential encodes genuinely non-minimal
behavior, while the outer ones reflect the PD-adjunction between $\delta^0_A$
and $\delta^2_A$ noted in Section~\ref{subsec:koszul-pd3}.

A direct computation yields
\[
\B_0(A)\cong \B_1(A)\cong S/(x_1,x_2)\cong\k,
\qquad
\B_i(A)=0\quad(i\ge 2).
\]

By contrast, the Koszul differential for $H^*(A)$ vanishes in degree~$2$.
Consequently,
\[
\B_0\big(H^\ast(A)\big)\cong \k,\qquad
\B_1\big(H^\ast(A)\big)\cong S,\qquad
\B_2\big(H^\ast(A)\big)\cong S^2/\langle(-x_2,x_1)\rangle,
\]
and $\B_i\big(H^\ast(A)\big)=0$ for $i\ge3$. 
Finally, another direct computation gives $\B^{0}(A)=\B^{1}(A)=0$ and 
$\B^{2}(A)\cong\B^{3}(A)\cong \k$, while
\[
\B^1\big(H^\ast(A)\big)\cong S^2/\langle(-x_2,x_1)\rangle,\qquad
\B^2\big(H^\ast(A)\big)\cong S,\qquad
\B^3\big(H^\ast(A)\big)\cong \k.
\]

This example illustrates both the sensitivity of Koszul homology to the
differential on chains and the intrinsic Poincar\'{e} duality relating
$\B_i(A)$ and $\B^{3-i}(A)$.

%%%%%%%%%%%%%%%%%%%%%%%%%%%%%%%%%%%%%%%
\section{Spectral sequences for Koszul (co)homology}
\label{sect:koszul-ss}
%%%%%%%%%%%%%%%%%%%%%%%%%%%%%%%%%%%%%%%

Filtering the Koszul (co)chain complexes by symmetric degree produces 
two dual spectral sequences interpolating between the Koszul complexes 
of $A$ and those of $H^*(A)$. The higher differentials encode iterated 
Massey products along the universal Aomoto class, making precise how 
non-formality of $A$ manifests in Koszul (co)homology. These spectral 
sequences establish Theorem~\ref{thm:koszul-ss-intro} of the Introduction.

%%%%%%%%%%%%%%
\subsection{The cohomological Koszul spectral sequence}
\label{subsec:coho-ss}
%%%%%%%%%%%%%%

Let $(A,d_A)$ be a connected $\k$-$\cdga$ over a field $\k$ of characteristic $0$, 
and fix $q\ge1$ such that $A$ is $q$-finite, that is, $\dim_\k A^j<\infty$ for
all $j\le q$.  
Let $S=\Sym(H_1(A))$ and let $\m\subset S$ be the augmentation ideal.
We filter the Koszul cochain complex $K^\bullet(A)=(A^\bullet\otimes_{\k} S,\delta_A)$
by the descending $\m$-adic filtration
\[
F^p K^\bullet(A) = A^\bullet \otimes_{\k} \m^p, \qquad p\ge 0.
\]
In the resulting spectral sequence, $p\ge 0$ denotes the filtration
(symmetric) degree, $q$ the complementary degree, and $p+q$ the total
cohomological degree.

By a standard abuse of notation, we write $[\omega_A]\cup(-)$ 
for the $S$-linear cup product action of the canonical element 
$\omega_A \in H^1(A)\otimes_{\k} H_1(A)$ 
on $H^*(A)\otimes_{\k} S$, with the $H_1(A)$-factor acting by 
multiplication in $S$.

\begin{theorem}
\label{thm:koszul-ss-coh}
The filtration $F^p K^\bullet(A)=A^\bullet\otimes_{\k} \m^p$ is a filtration 
by dg-ideals of the dg-algebra $S$-module $K^\bullet(A)$, and induces a 
cohomological spectral sequence of bigraded $\cdgas$ 
$(E_r^{\bullet,\bullet},d_r)$, whose terms vanish outside the region 
$p\ge0$, $p+q\ge0$, with the following properties:

\begin{enumerate}[itemsep=2.5pt]
\item  \label{ssc1}
The $E_0$-page is $E_0^{p,q}=A^{p+q}\otimes_{\k} \Sym^p(H_1(A))$, 
with multiplication inherited from $A^\bullet\otimes_{\k} S$, 
and $d_0 = d_A\otimes\id_S \colon E_0^{p,q}\to E_0^{p,q+1}$.

\item \label{ssc2}
The $E_1$-page is the Koszul cochain complex of the cohomology algebra:
  \[
  (E_1^{\bullet,\bullet},d_1)\cong K^\bullet\bigl(H^{*}(A)\bigr)
  \quad\text{as dg-algebras over $S$},
  \]
 where $E_1^{p,q}=H^{p+q}(A)\otimes_{\k}  \Sym^p(H_1(A))$ and 
 $d_1(\alpha)=[\omega_A]\cup\alpha$.
 
 \item \label{ssc3}
For $r\ge 2$, the differential $d_r\colon E_r^{p,q}\to E_r^{p+r,q-r+1}$ 
is given by the $(r+1)$-fold Massey product 
$\langle[\omega_A],\ldots,[\omega_A],\alpha\rangle$ 
along the canonical class $[\omega_A]\in H^1(A)$.
In particular, $d_r$ is a derivation with respect to the 
bigraded algebra structure on $E_r^{\bullet,\bullet}$.

\item \label{ssc4}
The spectral sequence converges in total degrees at most $q$: for every
$i\le q$ there is an isomorphism of bigraded $S$-modules
  \[
  E_\infty^{p,\,i-p} \cong \gr^F_p \B^{i}(A),
  \]
where $F$ is the filtration induced on $\B^\bullet(A)$ by the filtration
of $K^\bullet(A)$; the filtration $F$ is $\m$-stable, hence equivalent to
the $\m$-adic filtration.
\end{enumerate}
\end{theorem}

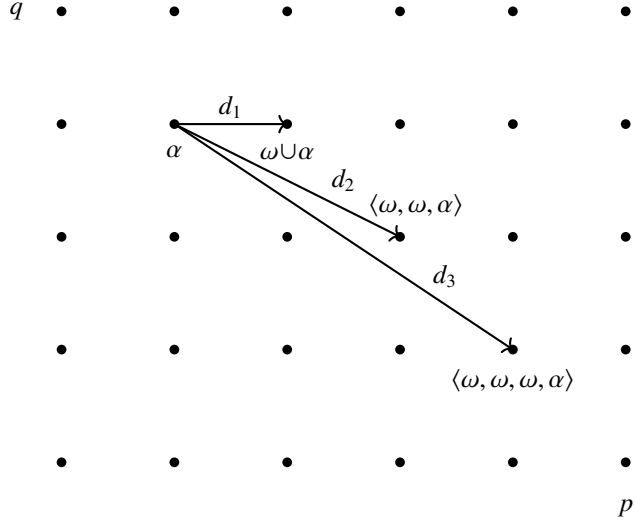
\begin{figure}[ht]
\centering
\begin{tikzpicture}[font=\small, x=1.55cm, y=1.05cm,
   dot/.style={circle,fill=black,inner sep=1.2pt}]
% Dots: the region $p\ge0$, $p+q\ge0$ in which the spectral sequence lives
\foreach \y in {0,1,2,3,4}
  \foreach \x in {0,1,2,3,4,5}
    \node[dot] at (\x,\y) {};
\foreach \x in {1,2,3,4,5} \node[dot] at (\x,-1) {};
\foreach \x in {2,3,4,5}   \node[dot] at (\x,-2) {};
\foreach \x in {3,4,5}     \node[dot] at (\x,-3) {};
\foreach \x in {4,5}       \node[dot] at (\x,-4) {};
\node[dot] at (5,-5) {};
% Boundary line $p+q=0$ of the region
\draw[gray!55, dashed] (-0.4,0.4) -- (5.45,-5.45);
\node[gray!70, font=\footnotesize] at (3.6,-4.75) {$p+q=0$};
% Axes
\draw[->, gray!60] (0,-0.45) -- (0,4.5) node[above, black] {$q$};
\draw[->, gray!60] (-0.45,0) -- (5.5,0) node[right, black] {$p$};
% The class $\alpha$ and its images
\node at (1,3) [left=3pt] {$\alpha$};
\node at (2,3) [above=4pt] {$\omega\!\cup\!\alpha$};
\node at (3,2) [above right=1pt] {$\langle \omega,\omega,\alpha\rangle$};
\node at (4,1) [below=4pt] {$\langle \omega,\omega,\omega,\alpha\rangle$};
% Differentials
\draw[->, thick] (1,3) -- (2,3);
\draw[->, thick] (1,3) -- (3,2);
\draw[->, thick] (1,3) -- (4,1);
\node at (1.5,3) [above=2pt] {$d_1$};
\node at (2.3,2.55) {$d_2$};
\node at (2.95,1.5) {$d_3$};
\end{tikzpicture}
\caption{Cohomological Koszul spectral sequence: 
differentials $d_1, d_2, d_3$ acting on a class $\alpha \in E_r^{1,3}$.}
\end{figure}

\begin{proof}
We filter the Koszul cochain complex $K^\bullet(A)=(A^\bullet\otimes_{\k} S,\delta_A)$ 
by the descending $\m$-adic filtration
\begin{equation}
\label{eq:fpk-amp}
F^p K^\bullet(A) = A^\bullet \otimes_{\k} \m^p.
\end{equation}
We take $p$ as the filtration degree (the symmetric degree on $S$) and $q$ as 
the internal degree, so that the total cohomological degree is $n = p+q$.
The differential is $\delta_A = \omega_A \cdot (-) + d_A \otimes \id_S$, where 
$\omega_A \in Z^1(A) \otimes_{\k} H_1(A)$ has $S$-degree $1$ while 
$d_A \otimes \id_S$ has $S$-degree $0$.

\smallskip
\noindent\emph{The $E_0$- and $E_1$-pages.}
The associated graded complex $\gr^\m(K^{\bullet}(A))$ has differential induced 
by the $S$-degree-zero component of $\delta_A$, which is $d_A \otimes \id_S$. 
Since $\Sym^p(H_1(A))$ is a free $\k$-module, the K\"{u}nneth formula 
gives $H^{p+q}(A^{\bullet}\otimes_{\k} \Sym^p(H_1(A)), d_A\otimes\id_S) 
\cong H^{p+q}(A)\otimes_{\k} \Sym^p(H_1(A))$ with no correction term, 
yielding
\begin{equation}
\label{eq:epq-hpqa}
E_1^{p,q} = H^{p+q}(A) \otimes_{\k} \Sym^p(H_1(A)).
\end{equation}
The $d_1$ differential arises from the $S$-degree-one component of $\delta_A$ 
(which increases filtration degree by $1$), namely multiplication by $\omega_A$. 
After passing to cohomology (to get the $E_1$-page), this becomes cup product 
with $[\omega_A] \in H^1(A)$. 
The resulting complex is precisely the Koszul complex of the cohomology algebra:
\begin{equation}
\label{eq:e1-d1-kb}
(E_1^{*,*}, d_1) \cong \bigl(H^{*}(A)\otimes_{\k} S, [\omega_A] \cup (-)\bigr) = 
K^{\bullet} \bigl(H^{*}(A)\bigr).
\end{equation}

\smallskip
\noindent\emph{Higher differentials.}
Let $\alpha \in E_r^{p,q}$ be a class surviving to page $r \geq 2$. 
Since $d_s(\alpha) = 0$ for all $1 \leq s < r$, the class $\alpha$
can be lifted step by step through the filtration: there exist elements
$a_0 \in A^{p+q}$ and $a_1, \ldots, a_{r-1} \in A^{p+q}$ satisfying
\begin{equation}
\label{eq:massey-system}
d_A(a_0) = 0, \qquad
d_A(a_{j+1}) = -\omega_A \cup a_j
\quad \text{for } j = 0, \ldots, r-2,
\end{equation}
with $[a_0] = \alpha$ on the $E_1$-page.
The differential $d_r(\alpha)$ is represented by the cocycle
$\omega_A \cup a_{r-1}$; this is indeed a cocycle, since
\[
d_A(\omega_A \cup a_{r-1})
= -\omega_A \cup d_A(a_{r-1})
= \omega_A \cup \omega_A \cup a_{r-2}
= 0,
\]
where the last equality uses $\omega_A^2 = 0$ in $A^2 \otimes_{\k} S$.
The system \eqref{eq:massey-system} is precisely a defining system
for the $(r+1)$-fold Massey product
$\langle [\omega_A], \ldots, [\omega_A], \alpha \rangle$
($r$ copies of $[\omega_A]$).

Different choices of lifts $a_0, \ldots, a_{r-1}$ satisfying
\eqref{eq:massey-system} yield different defining systems
(in the sense of \S\ref{subsec:massey}) for the $(r+1)$-fold
Massey product $\langle [\omega_A], \ldots, [\omega_A], \alpha \rangle$.
By the indeterminacy formula of \S\ref{subsec:massey}, any two
such choices yield cohomology classes in $H^{p+q+1}(A)$ differing by an
element of the indeterminacy of that Massey product.  For $r=2$ the latter is
the coset
\begin{equation}
\label{eq:massey-indet}
[\omega_A] \cdot H^{p+q}(A) + H^{1}(A) \cdot \alpha
\;\subseteq\; H^{p+q+1}(A),
\end{equation}
and for $r\ge3$ it is a larger subspace of $H^{p+q+1}(A)$, receiving
contributions also from the intermediate lifts.  In every case, at the
cochain level the ambiguity is represented by a cocycle
in $A^{p+q+1} \otimes_{\k} \m^{p+r+1} = F^{p+r+1}K^{\bullet}(A)$, and hence
maps to zero in the subquotient
\[
E_r^{p+r,\,q-r+1} = 
\frac{F^{p+r}K^{\bullet}(A) \cap \ker \delta_A}
     {F^{p+r+1}K^{\bullet}(A) \cap \ker \delta_A 
      + \delta_A\bigl(F^{p+r}K^{\bullet}(A)\bigr)}.
\]
Hence $d_r(\alpha)$ is a well-defined element of
$E_r^{p+r,\,q-r+1}$, independent of all choices:
the indeterminacy \eqref{eq:massey-indet}, which may be
nontrivial in $H^{p+q+1}(A)$, is entirely absorbed into
$F^{p+r+1}$ and vanishes in the $E_r$-subquotient.

\smallskip
\noindent\emph{Multiplicative structure.}
Observe that the descending $\m$-adic filtration 
$F^p K^\bullet(A) = A^\bullet \otimes_{\k} \m^p$ consists of dg-ideals of the 
bigraded dg-algebra $K^\bullet(A)$. Hence the induced spectral sequence 
is multiplicative from the $E_0$-page onward: each $E_r$ is a bigraded 
commutative algebra, every differential $d_r$ is a derivation, and the 
$E_\infty$-page is canonically isomorphic to the associated graded of 
$\B^\bullet(A) = H^{*}(K^\bullet(A))$ in the category of 
graded $S$-modules.

\smallskip
\noindent\emph{Convergence.}
Note first that the spectral sequence is not first-quadrant: in a fixed
total degree $n$ the terms $E_r^{p,n-p}$ can be nonzero for every
$p\ge0$, so the filtration on the abutment has infinitely many steps and
convergence is not formal.  It does hold in total degrees $\le q$, which is
all we shall need.
Indeed, fix $i\le q$; then $K^{i-1}(A)$ and $K^{i}(A)$ are finitely
generated free modules over the Noetherian ring $S$, so by the
Artin--Rees lemma the filtration
$F^p\B^i(A)=\im\bigl(H^i(F^pK^\bullet(A))\to\B^i(A)\bigr)$ is
$\m$-stable; the coboundaries arriving from $K^{i+1}(A)$, which need not
be finitely generated, are handled by exhaustiveness of the filtration
alone.  Hence the spectral sequence converges in those degrees, with
$E_\infty^{p,i-p}\cong\gr^\m_p\B^i(A)$; see \cite[pp.~22--24]{Serre}.  
The $\m$-adic filtration on $\B^i(A)$ need not be separated, however, 
so $\B^\bullet(A)$ is not in general recoverable from $E_\infty$; 
see Example~\ref{ex:coho-nonseparated}.
\end{proof}

\begin{example}
\label{ex:coho-nonseparated}
Let $A=\CE(\sol_2)$ be as in Example~\ref{ex:sol2-koszul}, so that
$S=\k[x]$ and, by \eqref{eq:toy-coh} and \eqref{eq:sol2-modules},
$\B^2(A)\cong\k[x]/(x-1)$, on which $\m=(x)$ acts invertibly.
Hence $\m^p\B^2(A)=\B^2(A)$ for every $p$, the $\m$-adic filtration is
not separated, and $\gr^\m\B^2(A)=0$ although $\B^2(A)\ne0$.  Thus the
cohomological Koszul spectral sequence, while convergent, carries no
information about $\B^2(A)$ in this case.  Contrast
Proposition~\ref{prop:koszul-separation}: positive weights exclude this
behaviour, and $\CE(\sol_2)$ has none.
\end{example}

Theorem~\ref{thm:koszul-ss-coh} provides the multiplicative refinement of 
Theorem~\ref{thm:koszul-ss-intro} announced in the Introduction. 
The corresponding degeneration result under formality or $q$-formality, 
for both the cohomological and homological spectral sequences simultaneously, 
is recorded in Corollary~\ref{cor:ss-formality} below.

%%%%%%%%%%%%%%%%%%%%%%%%%%%%%%
\subsection{The homological Koszul spectral sequence}
\label{subsec:hom-ss}
%%%%%%%%%%%%%%%%%%%%%%%%%%%%%%

We now present the homological counterpart of 
Theorem~\ref{thm:koszul-ss-coh}.
We filter the Koszul chain complex 
$K_\bullet(A)=(A_\bullet\otimes_{\k} T,\partial^A)$ 
by the descending $\n$-adic filtration
\[
F^p K_\bullet(A) = A_\bullet \otimes_{\k} \n^p, \qquad p\ge 0,
\]
where $\n = \bigoplus_{i\ge 1} T^i$ is the augmentation ideal of $T$,
and re-index by setting $F_{-p}=F^p$ following the convention 
of~\cite{PS-tams}.
The resulting spectral sequence lives in the second quadrant: 
$p\le 0$ denotes the filtration degree, $q\ge 0$ the complementary 
degree, and $p+q$ the total homological degree.
By a standard abuse of notation, we write $\omega_A^\vee\lrcorner(-)$ 
for the $T$-colinear contraction action of 
$\omega_A^\vee \in H_1(A)\otimes_{\k} H^1(A)$ on 
$H_\bullet(A)\otimes_{\k} T$, with the $H^1(A)$-factor acting 
by multiplication in $T$.

\begin{theorem}
\label{thm:koszul-ss-hom}
The filtration $F^p K_\bullet(A) = A_\bullet\otimes_{\k} \n^p$ is a 
filtration by dg-coideals of the right $T$-dg-comodule $K_\bullet(A)$, 
and induces a \emph{second-quadrant} homological spectral sequence of 
bigraded right $T$-comodules $(E^r_{\bullet,\bullet}, d^r)$, 
with $p\le 0$ and $q\ge 0$, and with the following properties:

\begin{enumerate}[itemsep=2.5pt]
\item \label{ssh1}
The $E^0$-page is $E^0_{p,q} = A_{p+q}\otimes_{\k} \Sym^{-p}(H^1(A))$,
with $T$-coaction inherited from $K_\bullet(A)$, 
and $d^0 = d_A^\vee\otimes\id_T \colon E^0_{p,q}\to E^0_{p,q-1}$.

\item \label{ssh2}
The $E^1$-page is canonically isomorphic, as a bigraded right 
$T$-comodule chain complex, to the Koszul chain complex of the 
cohomology algebra:
\[
(E^1_{\bullet,\bullet},d^1) \cong K_\bullet\bigl(H^{*}(A)\bigr),
\]
where $E^1_{p,q} = H_{p+q}(A)\otimes_{\k} \Sym^{-p}(H^1(A))$ and
$d^1 = \omega_A^\vee\lrcorner(-)$.
Passing to the $E^2$-page gives
\[
E^2_{p,q} \cong \bigl(\B_{p+q}(H^*(A))\bigr)_{-p},
\]
the piece of internal ($\Sym$-)degree $-p$ of the $(p+q)$-th Koszul 
homology module of $H^*(A)$; summing over the diagonal $p+q=n$ 
recovers the entire graded module $\B_n(H^*(A))$.

\item \label{ssh3}
For $r\ge 2$, the differential 
$d^r\colon E^r_{p,q}\to E^r_{p-r,\, q+r-1}$ is given by the 
$(r+1)$-fold homological Massey product along the 
$H^1(A)$-components of $\omega_A^\vee \in H_1(A)\otimes_{\k} H^1(A)$.
In particular, $d^r$ is a $T$-comodule morphism.

\item \label{ssh4}
Let $F$ denote the filtration induced on $\B_\bullet(A)$ by the filtration
of $K_\bullet(A)$.  The spectral sequence has a well-defined
$E^\infty$-page and there are natural isomorphisms of bigraded right
$T$-comodules
\[
\gr^F_{-p} \B_{p+q}(A) \cong E^\infty_{p,q}.
\]
The filtration $F$ is $\m$-stable, hence equivalent to the $\m$-adic
filtration.
\end{enumerate}
\end{theorem}

\begin{proof}
We filter the Koszul chain complex 
$K_\bullet(A) = (A_\bullet \otimes_{\k} T, \partial^A)$ by the 
descending $\n$-adic filtration
\[
F^p K_\bullet(A) = A_\bullet \otimes_{\k} \n^p, \qquad p \ge 0,
\]
where $\n = \bigoplus_{i \geq 1} T^i$ is the augmentation ideal of $T$,
and we re-index by setting $F_{-p} = F^p$ as in~\cite[\S 3.1]{PS-tams}.
The resulting increasing filtration $\{F_{-p}\}_{p \le 0}$ is bounded 
above by $F_0 = K_\bullet(A)$.

We decompose the Koszul differential as
\[
\partial^A = d_A^\vee \otimes \id_T + \omega_A^\vee \lrcorner (-),
\]
where $d_A^\vee$ has $T$-degree $0$ (preserves filtration degree) and 
$\omega_A^\vee \lrcorner (-)$ has $T$-degree $+1$ (increases filtration 
degree by $1$, hence maps $F^p$ into $F^{p+1} \subset F^p$).
Thus $\partial^A(F^p) \subset F^p$ and the filtration is preserved.

\smallskip
\noindent\emph{The $E^0$- and $E^1$-pages.}
The associated graded pieces are
\[
E^0_{p,q} = F^p K_{p+q} / F^{p+1} K_{p+q} 
= A_{p+q} \otimes_{\k} (\n^{-p}/\n^{-p+1})
\cong A_{p+q} \otimes_{\k} \Sym^{-p}(H^1(A)),
\]
with $d^0$ induced by $d_A^\vee \otimes \id$, going 
$E^0_{p,q} \to E^0_{p,q-1}$. Since $\Sym^{-p}(H^1(A))$ is a 
free $\k$-module, the K\"{u}nneth formula gives
\[
E^1_{p,q} = H_{p+q}(A) \otimes_{\k} \Sym^{-p}(H^1(A)).
\]
The $d^1$ differential is induced by $\omega_A^\vee \lrcorner (-)$,
which increases $T$-degree by $1$ (shifts $p$ by $-1$) and 
decreases chain degree by $1$, giving 
$d^1 \colon E^1_{p,q} \to E^1_{p-1,q}$.
This is precisely the differential of the Koszul chain complex
$K_\bullet(H^{*}(A))$, establishing the identification
\[
(E^1_{\bullet,\bullet}, d^1) \cong K_\bullet(H^{*}(A))
\]
as right $T$-comodule chain complexes. Taking homology yields 
$E^2_{p,q} \cong \bigl(\B_{p+q}(H^*(A))\bigr)_{-p}$, in agreement with 
item~\eqref{ssh2}.

\smallskip
\noindent\emph{Higher differentials.}
This step is formally dual to the corresponding step in the proof of
Theorem~\ref{thm:koszul-ss-coh}, with $(A^\bullet,d_A,\omega_A\cdot(-))$
replaced throughout by $(A_\bullet,d_A^\vee,\omega_A^\vee\lrcorner(-))$.
A class $\alpha\in E^r_{p,q}$ surviving to page $r\ge 2$ admits lifts
$a_0,\dots,a_{r-1}\in A_{p+q}$ satisfying
\begin{equation}
\label{eq:massey-system-hom}
d_A^\vee(a_0) = 0, \qquad
d_A^\vee(a_{j+1}) = -\omega_A^\vee \lrcorner\, a_j
\quad \text{for } j = 0, \ldots, r-2,
\end{equation}
with $[a_0]=\alpha$ on the $E^1$-page; the cycle
$\omega_A^\vee\lrcorner\,a_{r-1}$ represents $d^r(\alpha)$, and
\eqref{eq:massey-system-hom} is a defining system for the $(r+1)$-fold
homological Massey product
$\langle\omega_A^\vee,\dots,\omega_A^\vee,\alpha\rangle$ ($r$ copies of
$\omega_A^\vee$).  By the dual of the indeterminacy formula of
\S\ref{subsec:massey}, two such systems yield homology classes in
$H_{p+q-1}(A)$ differing by an element of the corresponding indeterminacy,
which for $r=2$ is
$\omega_A^\vee\lrcorner\,H_{p+q}(A)+H^1(A)\lrcorner\,\alpha$ and for $r\ge3$
is larger.  It is represented at the chain level by a cycle in
$A_{p+q-1}\otimes_{\k}\n^{-p+r+1}=F^{-p+r+1}K_\bullet(A)$, which
vanishes in the subquotient $E^r_{p-r,\,q+r-1}$; hence $d^r(\alpha)$ is
well defined.

\smallskip
\noindent\emph{Comultiplicative structure.}
The Koszul chain complex $K_\bullet(A) = A_\bullet \otimes_{\k} T$ 
carries a canonical right $T$-comodule structure via the coproduct 
$\Delta(e_i) = e_i \otimes 1 + 1 \otimes e_i$ on generators. 
Since $\partial^A(F^p) \subset F^p$ and 
$\rho(F^p) \subset F^p \otimes_{\k} T$, the filtration consists 
of dg-coideals, and the induced spectral sequence is comultiplicative: 
each $E^r$ is a bigraded right $T$-comodule and every $d^r$ is 
a $T$-comodule morphism.

\smallskip
\noindent\emph{Convergence.}
Since the filtration is bounded above ($F^0 = K_\bullet(A)$) 
and exhaustive, and each $E^1_{p,q}$ is finite-dimensional 
by the finite-type hypothesis, the spectral sequence has a 
well-defined $E^\infty$-page.  Since the filtration of $K_\bullet(A)$ is
separated, the general inclusion
$\gr^\m_{-p} \B_{p+q}(A) \subseteq E^\infty_{p,q}$
of~\cite[Prop.~5.1]{PS-tams} becomes the isomorphism of
item~\eqref{ssh4} for the induced filtration $F$.
\end{proof}

In homological degree $1$, the relation between the filtration $F$ of
item~\eqref{ssh4} and the $\m$-adic filtration of the abstract module
$\B_1(A)$ can be pinned down exactly: the two agree up to a shift by one.
This shift reflects the fact that $\B_1(A)$ is generated by syzygy-type
classes of internal degree $1$, and it is what reconciles the spectral
sequence, which computes $\gr^F$, with the Chen-rank formulas of
Section~\ref{sect:infalex-chen}, which are stated for the $\m$-adic
associated graded module.

\begin{lemma}
\label{lem:filtration-shift-B1}
Let $(A,d_A)$ be a connected, $1$-finite $\k$-$\cdga$, and let $F$ be the
filtration induced on $\B_1(A)$ by the $\m$-adic filtration
$F^pK_\bullet(A)=A_\bullet\otimes_{\k}\m^p$ of the Koszul chain complex, with
coefficients written in $S$.  Then
\[
F^0\B_1(A)=F^1\B_1(A)=\B_1(A),
\qquad
F^p\B_1(A)=\m^{p-1}\B_1(A) \quad\text{for all } p\ge1 .
\]
Consequently, $\gr^F_0\B_1(A)=0$, there are natural isomorphisms
$\gr^F_p\B_1(A)\cong\gr^\m_{p-1}\B_1(A)$ for $p\ge1$, and
\[
\Hilb\bigl(\gr^F\B_1(A),t\bigr)=t\cdot\Hilb\bigl(\gr^\m\B_1(A),t\bigr).
\]
For the $\cdga$ $(H^*(A),0)$, the module $\B_1(H^*(A))$ is graded and
generated in internal degree $1$, the filtration $F$ is the filtration by
internal degree, and $\gr^\m_k\B_1(H^*(A))$ is the piece of internal
degree $k+1$.
\end{lemma}

\begin{proof}
Since $A$ is connected, $d_A$ vanishes on $A^0=\k$, so $B^1(A)=0$ and the
cocycle representatives $\tilde e_1,\dots,\tilde e_n$ of a basis of
$H^1(A)$ span $Z^1(A)$.  Choose a complement, $A^1=Z^1(A)\oplus C^1$, and
decompose dually
\[
A_1 = \annn(C^1)\oplus\annn(Z^1(A)),
\qquad
\annn(C^1)\cong H_1(A).
\]
Because $d^0=0$, the Koszul boundary $\partial_1=\omega_A^\vee\lrcorner(-)$
is homogeneous of $S$-degree $+1$; in the coordinates above it is given by
$\partial_1(\alpha\otimes f)=\sum_{j}\alpha(\tilde e_j)\,x_jf$, so it
vanishes on $\annn(Z^1)\otimes_{\k} S$ and restricts on
$H_1(A)\otimes_{\k} S$ to the first Koszul differential of the regular sequence
$x_1,\dots,x_n$.  Hence the cycle module is graded and decomposes as
\begin{equation}
\label{eq:Z1-decomp}
Z_1(A)=\Syz\,\oplus\,\bigl(\annn(Z^1)\otimes_{\k} S\bigr),
\end{equation}
where $\Syz\subseteq H_1(A)\otimes_{\k} S$ is the module of Koszul syzygies of
$x_1,\dots,x_n$: since the Koszul complex resolves $\k$, this module is
generated in internal degree $1$, so that $\Syz_{\ge p}=\m^{p-1}\Syz$ for
all $p\ge1$.

\smallskip
\noindent (\emph{$F^p\subseteq\m^{p-1}$}.)  By \eqref{eq:Z1-decomp} and
gradedness,
\[
Z_1(A)\cap\bigl(A_1\otimes_{\k}\m^p\bigr)
=\Syz_{\ge p}\oplus\bigl(\annn(Z^1)\otimes_{\k}\m^p\bigr)
=\m^{p-1}\Syz\oplus\bigl(\annn(Z^1)\otimes_{\k}\m^p\bigr)
\subseteq\m^{p-1}Z_1(A),
\]
and $F^p\B_1(A)$ is by definition the image of the left-hand side.

\smallskip
\noindent (\emph{$\m^{p-1}\subseteq F^p$}.)  Since
$\im\bigl((d^1)^\vee\bigr)=\annn\bigl(\ker d^1\bigr)=\annn(Z^1)$, every
$\gamma\in\annn(Z^1)$ may be written $\gamma=(d^1)^\vee\beta$ for some
$\beta\in A_2$.  For any $f\in S$,
\[
\partial_2(\beta\otimes f)
=\gamma\otimes f+\sum_{j}(\tilde e_j\lrcorner\beta)\otimes x_jf ,
\]
so, modulo boundaries, $\gamma\otimes f$ is represented by the cycle
$-\sum_{j}(\tilde e_j\lrcorner\beta)\otimes x_jf$, which lies in
$A_1\otimes_{\k}\m^p$ whenever $f\in\m^{p-1}$.  Hence the image of
$\annn(Z^1)\otimes_{\k}\m^{p-1}$ in $\B_1(A)$ lies in $F^p\B_1(A)$; and the
image of $\m^{p-1}\Syz=\Syz_{\ge p}$ lies there directly.  As
$\m^{p-1}Z_1(A)$ is the sum of these two pieces,
$\m^{p-1}\B_1(A)\subseteq F^p\B_1(A)$.

\smallskip
Taking $p=1$ gives $\B_1(A)=\m^0\B_1(A)\subseteq F^1\B_1(A)$, whence
$F^0=F^1=\B_1(A)$ and $\gr^F_0=0$; the remaining assertions follow.  For
$(H^*(A),0)$ the complement $C^1$ is zero, so $Z_1=\Syz$ and
$\B_1(H^*(A))$ is a quotient of a module generated in internal degree
$1$, with $\B_1(H^*(A))_0=0$; for such a module $M$ one has
$\m^kM=M_{\ge k+1}$, which is the last assertion.
\end{proof}

The proof of Lemma~\ref{lem:filtration-shift-B1} is specific to
homological degree $1$: it rests on $B^1(A)=0$, forced by connectivity, and
on the syzygy module of the maximal ideal being generated in degree $1$.  In
higher degrees the shift is governed by where the cohomological Koszul
module begins, and by whether the associated graded module is generated
there.

\begin{proposition}
\label{prop:filtration-shift-higher}
Let $(A,d)$ be a connected, $q$-finite $\k$-$\cdga$, let $i\le q$, and
assume $\B_i(H^*(A))\ne0$.  Set
\[
\delta_i \coloneqq \min\bigl\{p\ge0 \bigm| \B_i(H^*(A))_p\ne0\bigr\},
\]
and let $F$ be the filtration induced on $\B_i(A)$ by that of
$K_\bullet(A)$.  Then:
\begin{enumerate}[itemsep=2pt]
\item \label{sh1}
$\gr^F_p\B_i(A)=0$ for $p<\delta_i$; consequently
$F^{\delta_i}\B_i(A)=\B_i(A)$, and
\[
\m^{p-\delta_i}\B_i(A)\subseteq F^p\B_i(A) \qquad\text{for all }
p\ge\delta_i .
\]
\item \label{sh2}
If moreover $\gr^F\B_i(A)$ is generated, as an $S$-module, in degree
$\delta_i$, then
\[
F^p\B_i(A)=\m^{p-\delta_i}\B_i(A)+\bigcap_{k\ge0}F^k\B_i(A)
\qquad\text{for all } p\ge\delta_i ;
\]
if in addition the $\m$-adic filtration on $\B_i(A)$ is separated, then
$F^p\B_i(A)=\m^{p-\delta_i}\B_i(A)$ and
\[
\Hilb\bigl(\gr^F\B_i(A);t\bigr)
=t^{\delta_i}\cdot\Hilb\bigl(\gr^\m\B_i(A);t\bigr).
\]
\item \label{sh3}
For $i=1$ one has $\delta_1=1$, and both hypotheses of \eqref{sh2} are
automatically satisfied: this is Lemma~\ref{lem:filtration-shift-B1}.
\end{enumerate}
\end{proposition}

\begin{proof}
\eqref{sh1} By Theorem~\ref{thm:koszul-ss-hom}\eqref{ssh2} and
\eqref{ssh4}, $\gr^F_p\B_i(A)\cong E^\infty_{p,i-p}$ is a subquotient of
$E^2_{p,i-p}\cong\B_i(H^*(A))_p$, which vanishes for $p<\delta_i$.  Hence
$F^p\B_i(A)=F^{p+1}\B_i(A)$ for $p<\delta_i$, and since
$F^0\B_i(A)=\B_i(A)$, the filtration is constant up to step $\delta_i$.
The inclusion follows because $\m F^k\subseteq F^{k+1}$ for all $k$, the
map $H_i(\m^kK_\bullet(A))\to\B_i(A)$ being $S$-linear.

\eqref{sh2} Generation in degree $\delta_i$ means
$\gr^F_{p+1}=S_1\cdot\gr^F_p$ for every $p\ge\delta_i$, that is,
$F^{p+1}\subseteq\m F^p+F^{p+2}$.  Iterating and using
$\m F^{p+j}\subseteq\m F^p$ gives $F^{p+1}\subseteq\m F^p+F^k$ for every
$k$, whence $F^{p+1}\subseteq\m F^p+\bigcap_{k}F^k$; the reverse inclusion
is clear.  Starting from $F^{\delta_i}=\B_i(A)$ and inducting on $p$ yields
the displayed equality.  If the $\m$-adic filtration is separated then so
is $F$, the two being equivalent by the Artin--Rees lemma, so the
intersection vanishes; the identity of Hilbert series then follows from
$\gr^F_p\B_i(A)\cong\gr^\m_{p-\delta_i}\B_i(A)$.

\eqref{sh3} For $i=1$, $\B_1(H^*(A))$ is a quotient of the module of Koszul
syzygies of $x_1,\dots,x_n$, hence $\delta_1=1$; the rest is
Lemma~\ref{lem:filtration-shift-B1}, whose conclusion is the equality of
\eqref{sh2} without any separation hypothesis.
\end{proof}

\begin{remark}
\label{rem:filtration-shift-higher}
The single-degree hypothesis in
Proposition~\ref{prop:filtration-shift-higher}\eqref{sh2} cannot be
dispensed with, and it is the reason no shift is available in general.  If
$\gr^F\B_i(A)$ has generators in two distinct degrees, no power of $t$ can
relate the two Hilbert series: already for the graded module
$M=S\oplus S(-1)$ over $S=\k[x,y]$ one has
$\Hilb(M;t)=(1+t)/(1-t)^2$, whereas $\gr^\m M\cong S\oplus S$ has
Hilbert series $2/(1-t)^2$.  When $A$ is $q$-formal, with graded Koszul
modules along a zig-zag as in
Corollary~\ref{cor:koszul-formal-separated}, one has
$\gr^F\B_i(A)\cong\B_i(H^*(A))$ with its internal grading, so the
hypothesis becomes the concrete requirement that $\B_i(H^*(A))$ be
generated in a single internal degree.  In all degrees, and with no
hypotheses at all, the supports of $\gr^F\B_i(A)$ and $\gr^\m\B_i(A)$
agree, by Lemma~\ref{lem:stable-filt-supp}.
\end{remark}

\begin{corollary}
\label{cor:ss-collapse-formal}
If $d_A=0$, then both Koszul spectral sequences of $A$ collapse at $E_2$,
and the edge isomorphisms are the identity of\/
$\B^\bullet(A)=\B^\bullet(H^*(A))$, respectively of\/
$\B_\bullet(A)=\B_\bullet(H^*(A))$.
\end{corollary}

\begin{proof}
For $d_A=0$ we have $H^*(A)=A$, so $E_1=K^\bullet(A)$ with $d_1$ the
Koszul differential and $E_2=\B^\bullet(A)$.  The modules
$\B^\bullet(A)$ are then graded, so the $\m$-adic filtration is the
grading and $\gr^\m\B^\bullet(A)=\B^\bullet(A)$; comparing graded
dimensions in each degree forces $d_r=0$ for $r\ge2$.  The homological
case is dual.
\end{proof}

\begin{remark}
\label{rem:convergence-asymmetry}
Strong convergence may fail for either Koszul spectral sequence, and 
for the same reason in both cases: the $E^\infty$-page computes the 
associated graded of the Koszul module, which determines the module 
itself only when the relevant filtration is separated.  This is what the 
convergence statements of Theorems~\ref{thm:koszul-ss-coh}\eqref{ssc4} 
and~\ref{thm:koszul-ss-hom}\eqref{ssh4} do and do not provide: neither 
$\B^\bullet(A)$ nor $\B_\bullet(A)$ is in general recoverable from 
$E_\infty$, as Examples~\ref{ex:coho-nonseparated} and~\ref{ex:sol2-ss} 
show.  By 
Proposition~\ref{prop:separation-criterion}, separation holds if and 
only if every associated prime of $\B_i(A)$---respectively of 
$\B^i(A)$---is contained in $\m$; it is guaranteed when the module has 
finite length, but not by Noetherianity alone.

When $d_A \neq 0$, the support of $\B_i(A)$ may contain components 
away from the origin, and any part of the module supported there is 
entirely invisible to the spectral sequence: if $M$ is a finitely 
generated $S$-module with $\m\notin\Supp_S(M)$, then $M_\m=0$, whence 
$\m \cdot M = M$, $\bigcap_{k\ge 0}\m^k M = M$, and $\gr^\m M = 0$.
Such an $M$ contributes nothing to $E^\infty$, 
even if the spectral sequence collapses at $E^2$.

This phenomenon is illustrated in detail by $A = \CE(\sol_2)$
in Example~\ref{ex:sol2-ss} below: the module
$\B_1(A) \cong \k[x]/(x-1)$ is supported at $\{x=1\}$, away 
from the origin, and is completely invisible to the homological 
spectral sequence despite the latter collapsing at $E^2$.

When $(A,d)$ has positive weights, each $\B_i(A)$ is a graded 
$T$-module, so its associated primes are homogeneous, hence contained 
in $\m$; by Proposition~\ref{prop:separation-criterion} the $\m$-adic 
filtration is separated (Proposition~\ref{prop:koszul-separation}) and 
strong convergence is restored.
\end{remark}

\subsection{Degeneration under formality and partial formality}
\label{subsec:ss-formality}
Passing to the $E^2$-page and applying the BGG correspondence 
formula \eqref{eq:tor-s-e} to $P = H^*(A)$ identifies
\begin{equation}
\label{eq:e2pq-tor}
E^2_{p,q} \cong \Tor^{S}_{-p}\bigl(H^{*}(A),\k\bigr)_{p+q}
\qquad (p\le 0).
\end{equation}
Under the re-indexing $p\mapsto -p$, this recovers the $E^2$-page 
identification stated in Theorem~\ref{thm:koszul-ss-intro}.
Theorem~\ref{thm:koszul-ss-hom} thus provides the comultiplicative 
refinement of Theorem~\ref{thm:koszul-ss-intro} announced in 
the Introduction.
Together, Theorems~\ref{thm:koszul-ss-coh} and~\ref{thm:koszul-ss-hom} 
yield the following degeneration result.

\begin{corollary}
\label{cor:ss-formality}
Let $(A,d_A)$ be a connected $\k$-$\cdga$ with $H^*(A)$ 
finite-dimensional in each degree.
\begin{enumerate}[itemsep=2pt]
  \item \label{cor-formal}
  If $A$ is formal, then $d_r=0$ for all $r\ge 2$ and 
  $d^r=0$ for all $r\ge 2$, and both spectral sequences 
  collapse at their respective $E_2$-pages, giving
  \begin{align*}
  \gr^\m_p \B^{i}(A) &\cong 
    \B^{i}\bigl(H^*(A)\bigr)^{\mathrm{weight}\,p}
    \quad (p\ge 0),\\
  \gr^\m_{-p} \B_{i}(A) &\cong 
    \bigl(\B_{i}\bigl(H^*(A)\bigr)\bigr)_{-p}
    \quad (p\le 0),
  \end{align*}
  as a bigraded $S$-module and as a bigraded right 
  $T$-comodule, respectively.
  \item \label{cor-qformal}
  If $A$ is $q$-formal for some $q\ge 1$, then 
  $d_r^{p,s}=0$ for all $r\ge 2$ and $p+s\le q$, and 
  $(d^r)_{p,s}=0$ for all $r\ge 2$ and $p+s\le q+1$.
  The respective isomorphisms of~\eqref{cor-formal} 
  hold for all cohomological degrees $\le q$ and all 
  homological degrees $\le q$.
\end{enumerate}
\end{corollary}

\begin{proof}
We start with the cohomological versions of the two statements. 
By Theorem~\ref{thm:koszul-ss-coh}\eqref{ssc3}, the differential 
$d_r$ on page $r\ge 2$ is given by an $(r+1)$-fold Massey product 
along $[\omega_A]\in H^1(A)$, which lies in $H^{p+s+1}(A)$ when 
acting on $E_r^{p,s}$.  Throughout this proof we write the second index of a
page as $s$, the letter $q$ being reserved for the formality parameter.

\eqref{cor-formal} 
Since $A$ is formal, Proposition~\ref{prop:massey-formal} gives 
$0 \in \langle [\omega_A], \ldots, [\omega_A], \alpha \rangle$ 
for every $\alpha \in H^k(A)$ and every $k \geq 0$; choosing 
a defining system realizing $0$ shows $d_r(\alpha)=0$ in 
$E_r^{p+r,\,s-r+1}$ for all $p,s$ and all $r \geq 2$. 
The collapse $E_2\cong E_\infty$ and the identification of the 
$E_\infty$-page with $\gr^\m_p \B^{i}(A)$ in total degree $i$ then follow from 
Theorem~\ref{thm:koszul-ss-coh}\eqref{ssc4}.

\eqref{cor-qformal}
If $A$ is $q$-formal, then by Proposition~\ref{prop:massey-formal} 
all Massey products in $H^{\le q+1}(A)$ vanish.
The differential $d_r^{p,s}$ lands in $H^{p+s+1}(A)$, which lies 
in $H^{\le q+1}(A)$ whenever $p+s\le q$.
Hence $d_r^{p,s}=0$ for all $r\ge 2$ and all $p+s\le q$, and 
the identification of $\gr^\m_p \B^i(A)$ with 
$\B^i(H^*(A))^{\mathrm{weight}\,p}$ for $i\le q$ follows as before.

The homological statements follow by the same argument, 
applied to Theorem~\ref{thm:koszul-ss-hom}\eqref{ssh3} 
and the dual version of Proposition~\ref{prop:massey-formal}, 
with $d^r$ in place of $d_r$.  The bound shifts by one, because $d^r$ lowers
the total degree where $d_r$ raises it: a map $H_n(A)\to H_{n-1}(A)$ vanishes
exactly when its dual $H^{n-1}(A)\to H^{n}(A)$ does, and the latter is the
cohomological Massey product valued in $H^{n}(A)$, which $q$-formality kills
as soon as $n\le q+1$.  Hence the condition $p+s\le q+1$ replaces
$p+s\le q$.  This is also what the last assertion requires: the isomorphism
in homological degree $i$ needs $d^r$ to vanish both out of degree $i$ and
into it from degree $i+1$, which gives $i\le q$, as stated.
\end{proof}

More generally, the higher differentials $d_r$ and $d^r$ ($r\ge 2$) 
in the two spectral sequences record successive Massey product 
obstructions to formality: the later the sequences collapse, the more 
severely the full nonlinear structure of $A$ departs from its 
linearization $H^*(A)$. In this sense, the pages of both spectral 
sequences provide a filtration of the gap between $A$ and its 
cohomology algebra, with the $E_2$-page capturing the purely 
linear (formal) approximation.

\begin{example}
\label{ex:sol2-ss}
We illustrate Theorems~\ref{thm:koszul-ss-coh} and~\ref{thm:koszul-ss-hom}
with the Chevalley--Eilenberg algebra $A = \CE(\sol_2)$ of
Example~\ref{ex:sol2-koszul}, whose Koszul complexes and modules are
recorded in \eqref{eq:toy}--\eqref{eq:sol2-modules}: here
$S=T=\k[x]$, $\omega_A=a\otimes x$, and
$\B_0(A)\cong\k[x]/(x)$, $\B_1(A)\cong\k[x]/(x-1)$.
The cohomology algebra is $H^*(A)\cong \bigwedge(a)$,
with $\B_1(H^*(A))=0$.
Since $A$ is formal (the inclusion 
$H^*(A)\hookrightarrow A$ is a quasi-isomorphism),
Corollary~\ref{cor:ss-formality}\eqref{cor-formal} applies.

\smallskip
\emph{Cohomological spectral sequence.}
The $E_1$-page is $K^\bullet(H^*(A)) = K^\bullet(\bigwedge(a))$,
the classical Koszul complex of the exterior algebra on one generator.
Since $\B^i(H^*(A))=0$ for $i\ge 1$, the $E_2$-page has 
$E_2^{p,q} = 0$ for $p+q \ge 1$, and $d_r = 0$ for all $r\ge 2$
by Corollary~\ref{cor:ss-formality}.
The spectral sequence collapses at $E_2$, giving
\[
\gr^\m_p \B^{p+q}(A) \cong \B^{p+q}(H^*(A))^{\mathrm{weight}\,p},
\]
which is nonzero only at $(p,q) = (0,0)$ and $(p,q) = (1,0)$, giving
$\gr^\m \B^1(A) \cong \k[x]/(x)$ and $\gr^\m \B^2(A) = 0$.
This is consistent with $\B^1(A) \cong \k[x]/(x)$ (already homogeneous)
and $\B^2(A) \cong \k[x]/(x-1)$ (non-homogeneous, so 
$\gr^\m \B^2(A) = 0$ since the filtration on $\k[x]/(x-1)$ 
is not separated: $\bigcap_{p} \m^p \B^2(A) = \B^2(A)$).

\smallskip
\emph{Homological spectral sequence.}
The $E^1$-page is $K_\bullet(H^*(A)) = K_\bullet(\bigwedge(a))$,
where $\B_0(H^*(A)) \cong \k$ and $\B_i(H^*(A))=0$ for $i\ge 1$.
By Corollary~\ref{cor:ss-formality}, $d^r=0$ for all $r\ge 2$,
and the spectral sequence collapses at $E^2$, giving
\[
\gr^\m_{-p} \B_{p+q}(A) \subseteq E^\infty_{p,q} 
= E^2_{p,q} \cong \bigl(\B_{p+q}(H^*(A))\bigr)_{-p}.
\]
This is nonzero only at $(p,q) = (0,0)$, where it equals $\k$.
Hence $\gr^\m \B_0(A) = \k$ (since $\B_0(A) = \k[x]/(x)$ 
is already homogeneous).
For $\B_1(A) \cong \k[x]/(x-1)$: one has $E^\infty_{-1,2} = 0$,
yet $\B_1(A) \neq 0$.
The spectral sequence does not converge strongly, since 
the filtration on $\B_1(A)$ is not separated:
$\m\B_1(A) = \B_1(A)$ (as $x-1$ is a unit in $\k[x]/(x-1)$,
so $x\cdot\B_1(A) = \B_1(A)$), whence 
$\bigcap_{p} \m^p\B_1(A) = \B_1(A) \neq 0$.
This illustrates the point of item~\eqref{ssh4}: the $E^\infty$-page 
computes the associated graded module, not the module itself, and here 
both $\gr\B_1(A)$ and $E^\infty$ vanish in total degree~$1$ while 
$\B_1(A)\ne0$:
the module $\B_1(A)$ is invisible to the spectral sequence 
despite being nonzero.

\smallskip
\emph{Summary.}
For $A = \CE(\sol_2)$, both spectral sequences collapse at $E_2$
as predicted by formality (Corollary~\ref{cor:ss-formality}),
yet neither recovers the full Koszul module $\B_1(A) \cong \k[x]/(x-1)$:
the cohomological spectral sequence misses the non-homogeneous part of
$\B^2(A)$, and the homological spectral sequence fails to converge
strongly to $\B_1(A)$.
This confirms that $E_2$-collapse is a property of the 
\emph{associated graded} of the Koszul modules,
not of the modules themselves.
\end{example}

%%%%%%%%%%%%%%%%%%%%%%%%%%%%%
\section{Associated graded Koszul modules and Hilbert series}
\label{sect:gr-koszul}
%%%%%%%%%%%%%%%%%%%%%%%%%%%%%

We introduce the associated graded Koszul modules and establish 
the key linearization result (Proposition~\ref{prop:koszul-linearization}): 
the associated graded of the Koszul complex of $A$ has homology 
$H^*(A)\otimes_{\k} S$, with induced differential the Koszul differential of 
$H^*(A)$.  In other words, the Koszul complex of the cohomology algebra 
appears as the first page of the spectral sequence of the $\m$-adic 
filtration---one page later than the associated graded complex itself.  
This identification underlies the tangent cone 
theorem and encodes the size of Koszul modules through Hilbert 
series, which we analyze in the final subsection.

%%%%%%%%%%%%%%%%%%%%%%%%%%%%%%%
\subsection{Linearization of the Koszul complex}
\label{subsec:koszul-linearization}
%%%%%%%%%%%%%%%%%%%%%%%%%%%%%%%

Let $(A,d_A)$ be a connected, finite-type $\k$-$\cdga$ and let 
$\m$ be the augmentation ideal of $S=\Sym(H_1(A))$. 
The following proposition is the key linearization result.  One might 
hope that the associated graded of the Koszul complex of $A$ is the Koszul 
complex of $H^*(A)$; this is not so---the two differentials live on 
different pages---but it becomes true after one homology step.  The 
proposition is the counterpart, valid with no weight or multigrading 
hypotheses, of Theorem~\ref{prop:multigraded-koszul}\eqref{mg3}, and 
underlies the tangent cone theorem for resonance varieties 
(Theorem~\ref{thm:tangent-cone}).

\begin{proposition}
\label{prop:koszul-linearization}
Let $(A,d_A)$ be a connected $\k$-$\cdga$ with $\dim_\k A^j<\infty$ for 
$j\le i+1$, for some $i\ge 0$.  Then:
\begin{enumerate}[itemsep=2pt]
\item \label{kl1}
The associated graded complexes of the $\m$-adic and $\n$-adic 
filtrations are the complexes with \emph{constant}\/ differential,
\[
\gr^\m K^\bullet(A) \cong \bigl(A\otimes_{\k} S,\, d_A\otimes\id\bigr),
\qquad
\gr^\n K_\bullet(A) \cong \bigl(A_\bullet\otimes_{\k} T,\, 
d_A^\vee\otimes\id\bigr);
\]
the Koszul part $\omega_A$ of the differential raises filtration degree 
and so induces the zero map on the associated graded.
\item \label{kl2}
Taking homology of \eqref{kl1} in $\cdga$-degrees $\le i+1$ and 
transporting the differential induced by $\omega_A$ yields canonical 
isomorphisms of complexes of bigraded $S$-modules and, respectively, 
bigraded right $T$-comodules,
\[
H^{\bullet}\bigl(\gr^\m K^\bullet(A)\bigr)
\cong K^\bullet\bigl(H^*(A)\bigr),
\qquad
H_{\bullet}\bigl(\gr^\n K_\bullet(A)\bigr)
\cong K_\bullet\bigl(H^*(A)\bigr),
\]
in degrees $\le i$.  Equivalently, the first pages of the Koszul spectral 
sequences of Theorems~\ref{thm:koszul-ss-coh} and~\ref{thm:koszul-ss-hom} 
are the Koszul complexes of the cohomology algebra: this is 
\eqref{ssc2} and \eqref{ssh2}.
\end{enumerate}
\end{proposition}

\begin{proof}
Both statements follow from the decomposition of the Koszul differential:
\[
\delta_A = d_A\otimes\id_S + \omega_A\cdot(-), \qquad
\partial^A = d_A^\vee\otimes\id_T + \omega_A^\vee\lrcorner(-),
\]
where in each case the first term has $S$-degree $0$ (preserving 
$\m$-adic filtration degree) and the second has $S$-degree $1$ (raising 
it by $1$).  Hence the induced differential on 
$\gr^\m_p K^\bullet(A)=A\otimes_{\k}\Sym^p(H_1(A))$ is $d_A\otimes\id$, 
which proves \eqref{kl1}; note that this is \emph{not}\/ the Koszul 
complex of $H^*(A)$, whose differential is induced by $\omega$, not 
by $d$.  Taking homology with respect to $d_A\otimes\id$ and using the 
K\"{u}nneth formula---valid since $\Sym^p(H_1(A))$ is a free 
$\k$-module---gives 
$H^i(\gr^\m_p K^\bullet(A))=H^i(A)\otimes_{\k}\Sym^p(H_1(A))$, and the 
differential induced by the filtration-degree-$1$ part $\omega_A$ of 
$\delta_A$ is $[\omega_A]\cup(-)$, the Koszul differential of 
$K^\bullet(H^*(A))$.  This is precisely the passage from the $E_0$-page 
to the $E_1$-page in Theorem~\ref{thm:koszul-ss-coh}, whose part 
\eqref{ssc2} records the identification.  The homological case is dual, 
via \eqref{ssh2}.
\end{proof}

\begin{remark}
\label{rem:linearization-one-page}
The distinction between the two pages in 
Proposition~\ref{prop:koszul-linearization} matters: identifying 
$\gr^\m K^\bullet(A)$ itself with $K^\bullet(H^*(A))$ would identify 
$\gr^\m\B^\bullet(A)$ with $\B^\bullet(H^*(A))$ unconditionally, which is 
false---the two differ whenever the Koszul spectral sequence supports a 
nonzero differential $d_r$ with $r\ge2$, and those differentials are 
Massey-type operations invisible to the cohomology ring 
\textup{(}Theorem~\ref{thm:koszul-ss-coh}\eqref{ssc3}\textup{)}.
\end{remark}

The graded isomorphism $S\cong T$ of Section~\ref{subsec:S-and-T} 
allows us to rewrite the chain complex $K_\bullet(A)$ with coefficients in 
$S$, and to regard its homology modules $\B_{\bullet}(A)$---naturally right 
$T$-comodules---as graded $S$-modules; by 
Section~\ref{subsec:koszul-filtrations} the augmentation-ideal filtrations 
correspond, so no ambiguity arises in writing $\m$ for both.
We denote by $S_k=\Sym^k(H_1(A))$ the degree-$k$ graded piece of $S$.

\begin{proposition}
\label{prop:koszul-ss-S}
Let $(A,d_A)$ be a connected $\k$-$\cdga$ with $H^*(A)$ 
finite-dimensional in each degree. 
The $\m$-adic filtration on $K_\bullet(A) = A_\bullet \otimes_{\k} S$
induces a homological spectral sequence
$(E^r_{p,q}, d^r)$ with $E^r_{p,q}=0$ unless $p\ge0$ and $p+q\ge0$
(after the reindexing $F_{-p}=F^p$, this is the second-quadrant
convention of Theorem~\ref{thm:koszul-ss-hom}), with
\[
E^1_{p,q} = H_{p+q}(A) \otimes_{\k} S_p, \qquad
E^2_{p,q} \cong \bigl(\B_{p+q}(H^*(A))\bigr)_{p},
\]
and $E^\infty_{p,q}\cong\gr^F_p\B_{p+q}(A)$ for the filtration $F$
induced on $\B_\bullet(A)$, as in
Theorem~\ref{thm:koszul-ss-hom}\eqref{ssh4}.
\end{proposition}

\begin{proof}
Via the isomorphism $S \cong T$ and the re-indexing $p \mapsto -p$,
this is the $S$-module counterpart of Theorem~\ref{thm:koszul-ss-hom}.
\end{proof}

%%%%%%%%%%%%%%%%%%%%%%%%%%%%%%%
\subsection{Associated graded Koszul modules}
\label{subsec:gr-koszul-mod}
%%%%%%%%%%%%%%%%%%%%%%%%%%%%%%%

Each Koszul module $\B_i(A)=H_i(K_{\bullet}(A))$ carries a natural 
descending $\m$-adic filtration,
\begin{equation}
\label{eq:m-adic-B}
\B_i(A)\supseteq \m\B_i(A)\supseteq \m^2\B_i(A)\supseteq \cdots,
\end{equation}
whose associated graded $S$-module records the infinitesimal structure 
of $\B_i(A)$ near the origin:
\begin{equation}
\label{eq:grB}
\gr\B_i(A) = \bigoplus_{k\ge 0} \m^k\B_i(A)/\m^{k+1}\B_i(A).
\end{equation}
The corresponding $\m$-adic completions,
\begin{equation}
\label{eq:hat-B}
\widehat{\B_i(A)} = \varprojlim_{k} \B_i(A)/\m^k\B_i(A),
\end{equation}
are filtered modules over the power series ring $\widehat{S}$, with 
$\gr\widehat{\B_i(A)}\cong \gr \B_i(A)$.

\begin{remark}
\label{rem:two-filtrations-gr}
Two filtrations on $\B_i(A)$ are in play throughout the paper, and any
statement about associated graded modules must say which one is meant.  The
first is the $\m$-adic filtration \eqref{eq:m-adic-B}, with associated graded
\eqref{eq:grB}.  The second is the filtration $F$ induced by the Koszul
spectral sequence, for which $E^\infty_{p,q}\cong\gr^F_p\B_{p+q}(A)$
\textup{(}Theorem~\ref{thm:koszul-ss-hom}\eqref{ssh4}\textup{)}.  The two are
equivalent as filtrations, and so have the same support and the same completion
\textup{(}Lemma~\ref{lem:stable-filt-supp}\textup{)}; they are not, however,
equal, and their graded pieces sit in different degrees.  For $i=1$ the
discrepancy is exactly one,
$\gr^F_p\B_1(A)\cong\gr^\m_{p-1}\B_1(A)$
\textup{(}Lemma~\ref{lem:filtration-shift-B1}\textup{)}; for $i\ge2$ the shift
is by the initial degree $\delta_i$ of $\B_i(H^*(A))$, and then only under a
generation hypothesis
\textup{(}Proposition~\ref{prop:filtration-shift-higher}\textup{)}.
Accordingly, an unadorned $\gr$ always denotes the $\m$-adic associated
graded \eqref{eq:grB}, and the superscript is written whenever the other
filtration is meant: $\gr^F$ for the one induced by the spectral sequence.  A third
indexing, by the total Koszul degree $i+p$ rather than by the filtration degree
$p$ alone, arises naturally on the BGG side of the correspondence; where it is
used it will be said so explicitly.
\end{remark}

\begin{example}
\label{ex:nonhomog-bis}
The $\cdga$ $A = \CE(\sol_2)$ of Example~\ref{ex:sol2-koszul} has 
$\B_1(A)\cong \k[x]/(x-1)$.
This module is not a graded $S$-module: the filtration is non-homogeneous,
since the relation $x=1$ is not homogeneous.
Moreover $\m$ acts invertibly on it, so
$\gr^\m\B_1(A)=0$ and $\widehat{\B_1(A)}=0$, while
$\B_1\bigl(H^*(A)\bigr)=0$ by \eqref{eq:sol2-modules-coh}.
Thus the non-homogeneous part of the differential---invisible to 
cohomology---affects the Koszul module;
the associated graded and completed Koszul modules detect only 
the linearized behavior encoded by $H^*(A)$.
\end{example}

The next theorem extends to all $i$ results of Fröberg--Löfwall \cite{FL} 
and of \cite{PS-mathann, SS06}, proved there by different methods for $i=1$.  
It expresses the graded pieces of the Koszul modules of a graded algebra in 
terms of Tor groups over the exterior algebra $\cE = \bigwedge A^1$.  The 
hypothesis that the differential vanish is essential; see 
Remark~\ref{rem:gr-Bi-tor-d-nonzero}. 

\begin{theorem}
\label{thm:gr-Bi-tor}
Let $A$ be a connected graded $\k$-algebra with $\dim_\k A^1<\infty$,
regarded as a $\cdga$ with zero differential.
Let $\cE=\bigwedge A^1$ be the exterior algebra on $A^1$ (with $A$ viewed
as a graded $\cE$-module in the natural way) and let $S=\Sym(A_1)$.
Then each $\B_i(A)$ is a graded $S$-module, and for all integers 
$i\ge 0$ and $k\ge 0$ there is a canonical isomorphism
of graded $\k$-vector spaces
\begin{equation}
\label{eq:gr-Bi-tor}
\bigl(\B_i(A)_k\bigr)^\vee \cong \Tor^\cE_{k-i}(A,\k)_{\,k},
\end{equation}
where $(-)^\vee$ denotes graded $\k$-linear dual, the subscript $k$ on the 
left denotes the component of total Koszul degree $k$---that is, the image of 
$A_i\otimes_{\k}S_{k-i}$---and the subscript $k$ on the right denotes the 
component of internal degree $k$.
\end{theorem}

\begin{proof}
Let $P=\widehat{A}$ be the graded dual (completed) $\cE$-module and form the
BGG complex $\bL(P)$ over $S$.  By \eqref{eq:lp} we have
$\bL(P)_i=S\otimes_{\k}A_i$, with differential the contraction by the
canonical element of $A^1\otimes_{\k}A_1$.  Since $d_A=0$, we have
$H^1(A)=A^1$ and $S=\Sym(A_1)$, and this is exactly the Koszul chain complex
$K_\bullet(A)$ of Section~\ref{subsec:hom-koszul}.  Hence
\begin{equation}
\label{eq:hi-bpl}
H_i(\bL(P))_{k}\cong \B_i(A)_k,
\end{equation}
the subscript $k$ being the total Koszul degree on either side.
Applying the BGG duality formula~\eqref{eq:duality} with $(q,i)\mapsto(i,k-i)$ yields
\begin{equation}
\label{eq:hi-bpl-tor}
\bigl(H_i(\bL(P))_{k}\bigr)^{\vee}
\cong
\Tor^{\widehat{\cE}}_{k-i}(A,\k)_{-k},
\end{equation}
where the negative internal degree reflects the grading convention on the
graded dual $P=\widehat A$.

Using the identification between polynomial degree on $S$ and the opposite
internal degree on $\cE$ (polynomial degree $k$ on $S$ corresponds to internal 
degree $-k$ on the negatively graded exterior algebra, 
see Remark~\ref{rem:bgg-grading}), we reindex to obtain
\begin{equation}
\label{eq:grb-tor}
\bigl(\B_i(A)_k\bigr)^\vee
\cong
\Tor^{\widehat{\cE}}_{k-i}(A,\k)_{k}.
\end{equation}
Since $\cE$ is finite-dimensional, completion does not affect Tor groups in
fixed internal degree (as $\Tor$ commutes with direct limits in fixed degree), 
and the isomorphism descends to $\Tor^{\cE}_{k-i}(A,\k)_{k}$.
\end{proof}

\begin{corollary}
\label{cor:gr-Bi-tor-cdga}
Let $(A,d_A)$ be a connected $\k$-$\cdga$ with $\dim_\k H^1(A)<\infty$, and
set $\cE=\bigwedge H^1(A)$.  Then, for all $i,k\ge0$,
\[
\bigl(\B_i(H^*(A))_k\bigr)^\vee \cong \Tor^\cE_{k-i}\bigl(H^*(A),\k\bigr)_{\,k},
\]
and $\gr^F_{k-i}\B_i(A)$ is a subquotient of $\B_i(H^*(A))_k$, for every $k$.
Equality holds for all $i$ and $k$ precisely when the Koszul spectral sequence
of Theorem~\ref{thm:koszul-ss-hom} degenerates at $E^2$.
\end{corollary}

\begin{proof}
The cohomology algebra $H^*(A)$ carries the zero differential, so
Theorem~\ref{thm:gr-Bi-tor} applies to it and gives the displayed isomorphism,
$\cE=\bigwedge H^1(A)$ being the relevant exterior algebra.  The second
assertion is the convergence of the Koszul spectral sequence: its $E^2$-page is
$\B_i(H^*(A))$ and its $E^\infty$-page is $\gr^F\B_i(A)$, by
Theorem~\ref{thm:koszul-ss-hom}\eqref{ssh2} and \eqref{ssh4}, so each
$E^\infty$-term is a subquotient of the corresponding $E^2$-term, with equality
throughout exactly when all higher differentials vanish.
\end{proof}

The degeneration hypothesis is met, in particular, whenever $A$ is formal, by
Corollary~\ref{cor:ss-formality}.

\begin{remark}
\label{rem:gr-Bi-tor-d-nonzero}
Theorem~\ref{thm:gr-Bi-tor} does not extend to $\cdgas$ with $d_A\ne0$ by
replacing $\B_i(A)$ with $\gr\B_i(A)$.  The obstruction is structural: the BGG
complex $\bL(\widehat A)$ has as differential the contraction by
$\omega_{A^1}$, with no $d_A^\vee$ term, whereas $\gr^\m K_\bullet(A)$ has
differential $d_A^\vee$ and no $\omega$ term
\textup{(}Proposition~\ref{prop:koszul-linearization}\eqref{kl1}\textup{)}.
The two complexes carry complementary halves of $\partial^A$, and agree
exactly when $d_A=0$.  For a concrete failure, take $A=\CE(\h(1))$ of
Section~\ref{subsec:heisenberg-1}: here $A^1$ is spanned by $e_1,f_1,g$, so
$\cE=\bigwedge A^1$ is $A$ itself and $A$ is free of rank one over $\cE$;
hence $\Tor^\cE_j(A,\k)=0$ for $j\ge1$ while $\Tor^\cE_0(A,\k)=\k$ sits in
internal degree $0$, so that $\Tor^\cE_{k-1}(A,\k)_k=0$ for every $k$.  Yet
$\B_1(A)\cong\k$, so that $\gr^\m_0\B_1(A)=\k$, and no reindexing can
reconcile the two sides.  What survives for a general $\cdga$ is
Corollary~\ref{cor:gr-Bi-tor-cdga}, together with the Hilbert series
inequality of Theorem~\ref{thm:Hilb-ineq-gr}.
\end{remark}

\begin{corollary}
\label{cor:gr-Bi-vanish}
If $A$ is a connected graded $\k$-algebra with $\dim_\k A^1<\infty$,
then $\B_i(A)_k=0$ for all $k<i$.
\end{corollary}

\begin{proof}
Immediate from Theorem~\ref{thm:gr-Bi-tor}, since
$\Tor^{\cE}_{p}(A,\k)=0$ for $p<0$.
\end{proof}

\begin{corollary}
\label{cor:chen-eventual}
Let $A$ be a connected graded $\k$-algebra with $\dim_\k A^1<\infty$.
If $\B_1(A)_k=0$ for some $k\ge 1$, then
$\B_1(A)_n=0$ for all $n\ge k$.
\end{corollary}

\begin{proof}
By Theorem~\ref{thm:gr-Bi-tor} with $i=1$,
\[
\bigl(\B_1(A)_k\bigr)^\vee \cong \Tor^{\cE}_{k-1}(A,\k)_{k}
= \beta_{k-1,k}^{\cE}(A),
\]
so $\B_1(A)_k=0$ if and only if $\beta_{k-1,k}^{\cE}(A)=0$.

It therefore suffices to show that $\beta_{k-1,k}^{\cE}(A)=0$ implies 
$\beta_{p,p+1}^{\cE}(A)=0$ for all $p\ge k-1$.
Consider the linear strand of the minimal free $\cE$-resolution of $A$:
\begin{equation}
\label{eq:linear-strand}
\begin{tikzcd}[column sep=22pt]
\cdots \ar[r] &
\cE(-p-1)^{\beta_{p,p+1}} \ar[r, "\phi_p"] &
\cE(-p)^{\beta_{p-1,p}} \ar[r] &
\cdots \ar[r] &
\cE(-1)^{\beta_{0,1}}.
\end{tikzcd}
\end{equation}
This is a complex of free $\cE$-modules with linear differentials.
Suppose $\beta_{k-1,k}=0$, so that the term $\cE(-k)^{\beta_{k-1,k}}=0$.
The differential $\phi_k\colon \cE(-k-1)^{\beta_{k,k+1}}\to 
\cE(-k)^{\beta_{k-1,k}}=0$ must then be zero, so by minimality 
$\beta_{k,k+1}=0$.
Proceeding inductively, $\phi_p\colon \cE(-p-1)^{\beta_{p,p+1}}\to 
\cE(-p)^{\beta_{p-1,p}}=0$ forces $\beta_{p,p+1}=0$ for all $p\ge k-1$.
Applying Theorem~\ref{thm:gr-Bi-tor} again gives
$\B_1(A)_n=0$ for all $n\ge k$.
\end{proof}

\begin{remark}
\label{rem:chen-eventual-i1}
The restriction to $i=1$ in Corollary~\ref{cor:chen-eventual} is essential.
That case corresponds to the \emph{linear strand} of the minimal free
resolution of $A$ over $\cE$, and the linear strand cannot revive once it
terminates: a nonzero differential into a zero module is impossible, and
minimality then forces the source to be zero.  For $i\ge2$, by contrast,
$\B_i(A)_k$ is dual to $\Tor^{\cE}_{k-i}(A,\k)_k$, which lies on the
strand of slope $k/(k-i)>1$, where propagation of vanishing is not controlled
by minimality alone; see \cite{EFS} for the BGG perspective.
\end{remark}

\begin{remark}
\label{rem:gr-finite}
Theorem~\ref{thm:gr-Bi-tor} has a perhaps unexpected consequence for a graded
algebra $A$.  Even when $A$ is infinite-dimensional in degrees $>1$, so that
the Koszul modules $\B_i(A)$ need not be finitely generated over
$S=\Sym(A_1)$, they exhibit a strong \emph{degreewise finiteness} property.

More precisely, although $\B_i(A)$ may be neither finitely generated over $S$
nor finite-dimensional as a whole, Theorem~\ref{thm:gr-Bi-tor} implies that it
is locally finite with respect to the internal grading inherited from $A$:
for each internal degree $d$ and each total Koszul degree $k$, the vector space
$(\B_i(A)_k)_d$ 
is finite-dimensional and determined functorially by the
$\cE=\bigwedge A^1$-module structure of $A$, via the groups
$\Tor^{\cE}(A,\k)$.

Thus the failure of finite generation of $\B_i(A)$ over $S$ is a global
phenomenon, reflecting the accumulation of infinitely many finite-dimensional
pieces, rather than a breakdown of finiteness in any fixed internal degree.
\end{remark}

\begin{example}
\label{ex:non-finitely-gen}
Let $A = \bigwedge(a) \otimes_{\k} \k[b_1,b_2,\dots]$, with $\abs{a}=1$,
$\abs{b_j}=2$, and zero differential. Then $\dim_{\k} H^1(A)=1$, so
$S=\k[x]$ and $\omega_A=a\otimes x$.
The Koszul complex $K_{\bullet}(A)$ has the form
\begin{equation*}
\label{eq:nonfg-ses}
\begin{tikzcd}[column sep=24pt]
A_2\otimes_{\k} S \arrow[r, "\partial_2"] &
A_1\otimes_{\k} S \arrow[r, "\partial_1"] &
A_0\otimes_{\k} S,
\end{tikzcd}
\end{equation*}
where $\partial_1(a\otimes s)=1\otimes xs$ and $\partial_2=0$.
Consequently,
\[
\B_1(A)=0
\qquad\text{and}\qquad
\B_2(A)=A_2\otimes_{\k} S \cong \boplus_{j\ge1} S\cdot b_j,
\]
which is not finitely generated as an $S$-module.

The graded pieces reflect the same global non-finiteness: for each $k\ge0$,
\[
\B_2(A)_k \cong \boplus_{j\ge1} S_{k-2}\cdot b_j,
\]
an infinite-dimensional $\k$-vector space.
Nevertheless, Theorem~\ref{thm:gr-Bi-tor} identifies its internal graded pieces as
\[
(\B_2(A)_k)_d \cong \Tor^{\bigwedge(a)}_{k-2}(A,\k)_{k,d},
\]
which are finite-dimensional for each fixed internal degree $d$, since
$\bigwedge(a)$ has global dimension $1$.
This example illustrates how infinite generation of $\B_i(A)$ may coexist with
degreewise finiteness in each internal degree, as predicted by
Theorem~\ref{thm:gr-Bi-tor}.
\end{example}

The $\cdga$ in Example~\ref{ex:non-finitely-gen} is a rational model for the space
$S^1 \times \bigvee_{j=1}^{\infty} \CP^{\infty}$.
From this perspective, the example shows that infinite generation of Koszul
modules can already occur for very simple spaces, even when $H^1$ is
finite-dimensional and the fundamental group is abelian.

%%%%%%%%%%%%%%%%%%%%%%%%%%%
\subsection{Hilbert series comparisons}
\label{subsec:hilb}
%%%%%%%%%%%%%%%%%%%%%%%%%%%

The Koszul spectral sequences of Section~\ref{sect:koszul-ss} impose strong 
numerical constraints on Koszul (co)homology. In particular, they yield coefficientwise 
inequalities between the Hilbert series of the associated graded Koszul modules 
of a $\cdga$ and those of its cohomology algebra. 
The results of this subsection should be understood as genuinely
degreewise statements: even when a Koszul module $\B_i(A)$ or its associated
graded $\gr\B_i(A)$ fails to be finitely generated over $S$, the coefficientwise
Hilbert-series comparisons remain meaningful.
Example~\ref{ex:non-finitely-gen} illustrates this phenomenon concretely.

\begin{theorem}
\label{thm:Hilb-ineq-gr}
Let $(A,d)$ be a connected $q$-finite $\k$-$\cdga$, and set 
$S=\Sym(H_1(A))$.  For each $0\le i\le q$, let $F$ denote the filtration
induced on the Koszul module $\B_i(A)$ by the spectral sequence of
Proposition~\ref{prop:koszul-ss-S}, with associated graded $\gr^F\B_i(A)$;
recall from Theorem~\ref{thm:koszul-ss-hom}\eqref{ssh4} that $F$ is
equivalent to the $\m$-adic filtration.  Then, coefficientwise,
\[
\Hilb\bigl(\gr^F\B_i(A);t\bigr) 
\preccurlyeq \Hilb\bigl(\B_i(H^\ast(A));t\bigr).
\]
For $i=1$, both sides shift uniformly to the $\m$-adic normalization: by
Lemma~\ref{lem:filtration-shift-B1},
\[
\Hilb\bigl(\gr^\m\B_1(A);t\bigr) 
\preccurlyeq \Hilb\bigl(\gr^\m\B_1(H^\ast(A));t\bigr)
=t^{-1}\,\Hilb\bigl(\B_1(H^\ast(A));t\bigr).
\]
Moreover, equality holds in degree $i$ if and only if no differential $d^r$ 
with $r\ge2$ enters or leaves the terms $E^r_{p,i-p}$ of the $\m$-adic 
Koszul spectral sequence; in particular, equality holds for all $i\le q$ if 
and only if that spectral sequence degenerates at $E^2$ in total degrees 
$\le q$.  No separation hypothesis is needed: both sides of the inequality 
are computed from associated graded modules.
\end{theorem}

\begin{proof}
Filter the Koszul chain complex $K_\bullet(A)=A_\bullet\otimes_{\k} S$ by the
$\m$-adic filtration on the coefficient ring, $F^t K_\bullet(A)=A_\bullet\otimes_{\k}\m^t$,
viewed as a decreasing filtration whose associated graded computes
$\gr^\m \B_\bullet(A)$.
The resulting homological spectral sequence (Proposition~\ref{prop:koszul-ss-S})
has terms vanishing outside $p\ge0$, $p+q\ge0$, and satisfies
\begin{equation}
\label{eq:e1pq-omegaA-bis}
E^1_{p,q}=H_{p+q}(A)\otimes_{\k} \Sym^p(H_1(A)),
\qquad d^1=\text{contraction by }\omega_A^\vee.
\end{equation}
The $E^1$-page in the $p$-direction is the Koszul chain complex of the
graded $S$-module $H_\ast(A)\otimes_{\k} S$, and its homology yields 
the $E^2$-page:
\begin{equation}
\label{eq:e2pq-bpq}
E^2_{p,q} \cong \bigl(\B_{p+q}(H^\ast(A))\bigr)_{p},
\end{equation}
where the subscript on the right denotes the summand of $S$-degree $p$.

In total degrees $\le q$ the spectral sequence converges to the associated 
graded of the Koszul homology of $A$, for the induced filtration $F$:
\begin{equation}
\label{eq:einfty-bpq}
E^\infty_{p,i-p}\cong \gr^F_p\bigl(\B_{i}(A)\bigr),
\qquad i\le q.
\end{equation}
Only the finite generation of $K_{i-1}(A)$ and $K_i(A)$ enters here---both 
hold because $\dim_\k A^j<\infty$ for $j\le q$---while the boundaries 
arriving from $K_{i+1}(A)$ are handled by exhaustiveness of the filtration; 
see the convergence discussion in the proof of 
Theorem~\ref{thm:koszul-ss-coh}.  In particular, no finiteness assumption 
on $A^{q+1}$ is required.
For fixed total homological degree $i$ and polynomial degree $d$, each
group $E^\infty_{p,i-p,d}$ is a subquotient of $E^2_{p,i-p,d}$. Hence
\begin{equation}
\label{eq:grpm-bia}
\dim_\k \gr^F_p\bigl(\B_i(A)\bigr)_d
=\dim_\k E^\infty_{p,i-p,d}
\le \dim_\k E^2_{p,i-p,d}
=\dim_\k \bigl(\B_i(H^\ast(A))\bigr)_d.
\end{equation}
Summing over $p\ge0$ gives
\begin{equation}
\label{eq:dim-grbi}
\dim_\k \bigl(\gr^F\B_i(A)\bigr)_d
=\sum_{p\ge0}\dim_\k \gr^F_p\bigl(\B_i(A)\bigr)_d
\le \dim_\k \bigl(\B_i(H^\ast(A))\bigr)_d,
\end{equation}
which establishes the claimed coefficientwise inequality; the $\m$-adic
form for $i=1$ follows by Lemma~\ref{lem:filtration-shift-B1}, applied to
$A$ on the left and to $(H^*(A),0)$ on the right.

Finally, \eqref{eq:grpm-bia} is an equality for every $p$ and $d$ precisely 
when $E^\infty_{p,i-p}=E^2_{p,i-p}$ for all $p$, that is, when no 
differential $d^r$ with $r\ge2$ enters or leaves the line $p+q=i$.
\end{proof}

\begin{corollary}
\label{cor:support-dim-ineq}
With hypotheses as above, the Krull dimension satisfies
\[
\dim \Supp\bigl(\gr^\m\B_i(A)\bigr) 
= \dim \Supp\bigl(\gr^F\B_i(A)\bigr)
\le \dim \Supp\bigl(\B_i(H^\ast(A))\bigr),
\]
the first equality by Lemma~\ref{lem:stable-filt-supp}, since $F$ is
equivalent to the $\m$-adic filtration.  In particular, vanishing of
$\B_i(H^\ast(A))$ in degree $d$ implies vanishing of $\gr^F\B_i(A)$ in
degree $d$.
\end{corollary}

\begin{proof}
If a finitely generated graded $S$-module $M$ has graded pieces equal to $0$ 
in all sufficiently large degrees then $\dim\Supp(M)=-\infty$; otherwise the
leading polynomial growth of $\Hilb(M;t)$ determines the Krull dimension.
The coefficientwise inequality of Theorem~\ref{thm:Hilb-ineq-gr} therefore
implies the stated inequality on dimensions of support.
\end{proof}

\begin{corollary}
\label{cor:Euler-Hilb}
Let $(A,d)$ be a connected, finite-type $\k$-$\cdga$, and suppose that the
total $E^2$-page is finite-dimensional, that is,
$\dim_\k\boplus_{i\ge0}\B_i(H^*(A))<\infty$.  Equivalently,
$\RR_{i,1}(H^*(A))\subseteq\{0\}$ for every $i$---so that each
$\B_i(H^*(A))$ has finite length---\emph{and}\/ $\B_i(H^*(A))=0$ for
$i\gg0$.  Then
$\boplus_{i}\gr^F\B_i(A)$ is finite-dimensional as well, and the two Euler
characteristics agree:
\[
\sum_{i\ge 0} (-1)^i \dim_\k \gr^F\B_i(A) =
\sum_{i\ge 0} (-1)^i \dim_\k \B_i\bigl(H^\ast(A)\bigr).
\]
\end{corollary}

\begin{proof}
By Theorem~\ref{thm:koszul-ss-hom}\eqref{ssh2} the total $E^2$-page is
$\boplus_{i}\B_i(H^*(A))$, finite-dimensional by hypothesis; each subsequent
page is a subquotient of it, hence also finite-dimensional, and
$E^\infty=\boplus_{i}\gr^F\B_i(A)$ by \eqref{ssh4}.  Each differential $d^r$
lowers the total degree $i=p+q$ by one, so for every $r\ge2$,
\[
\sum_{i}(-1)^i\dim_\k E^{r+1}_i
=\sum_{i}(-1)^i\Bigl(\dim_\k E^{r}_i-\rank d^r|_{E^r_i}
-\rank d^r|_{E^r_{i+1}}\Bigr)
=\sum_{i}(-1)^i\dim_\k E^{r}_i ,
\]
the two rank terms cancelling in pairs because they occur with opposite
signs.  Since all sums are finite and the pages stabilize, the Euler
characteristic of $E^2$ equals that of $E^\infty$.
\end{proof}

The identity does not refine to an equality of Hilbert series, and the
finiteness hypothesis cannot be dropped.  The reason is that the
differential $d^r$ shifts the internal degree: on every page the internal
$S$-degree of $E^r_{p,q}$ equals $p$, while $d^r$ maps $E^r_{p,q}$ to
$E^r_{p-r,q+r-1}$, so a fixed internal degree does not span a subcomplex
and the alternating sum of dimensions in each degree separately is not
invariant.  A page-by-page count instead gives
\[
\sum_{i}(-1)^i\Hilb\bigl(E^{r+1}_i;t\bigr)
-\sum_{i}(-1)^i\Hilb\bigl(E^{r}_i;t\bigr)
\in (t^r-1)\cdot\Z[[t]] ,
\]
which vanishes at $t=1$---whence the corollary---but not identically.

\begin{example}
\label{ex:euler-hilb-fails}
The Heisenberg $\cdga$ $A=\CE(\h(1))$ shows the failure concretely.  The
computations of Section~\ref{subsec:heisenberg-1} give
$\B_0(A)\cong\B_1(A)\cong\k$ and $\B_i(A)=0$ for $i\ge2$, so that, by
Lemma~\ref{lem:filtration-shift-B1},
\[
\sum_{i\ge0}(-1)^i\Hilb\bigl(\gr^F\B_i(A);t\bigr)=1-t .
\]
On the other hand, the Koszul differential of $H^*(A)$ vanishes in
degree~$2$---the cup product $\bwedge^2H^1(A)\to H^2(A)$ being zero---and
the same computations give $\B_0(H^*(A))\cong\k$,
$\B_1(H^*(A))\cong S$ \textup{(}free of rank $1$, generated in internal
degree $1$\textup{)}, $\B_2(H^*(A))\cong S^2/\langle(-x_2,x_1)\rangle$, and
$\B_i(H^*(A))=0$ for $i\ge3$.  Hence
\[
\sum_{i\ge0}(-1)^i\Hilb\bigl(\B_i(H^*(A));t\bigr)
=1-\frac{t}{(1-t)^2}+\frac{2-t}{(1-t)^2}
=1+\frac{2}{1-t}
=3+2t+2t^2+\cdots ,
\]
a series that differs from $1-t$ in every degree.  Note that the
hypothesis of Corollary~\ref{cor:Euler-Hilb} fails here, since
$\B_1(H^*(A))$ and $\B_2(H^*(A))$ are infinite-dimensional: the higher
differentials of the Koszul spectral sequence of $A$ destroy all of
$\B_1(H^*(A))$ but a single class, and they do so across internal
degrees.
\end{example}

The two conditions in the hypothesis of Corollary~\ref{cor:Euler-Hilb} are
independent: the vanishing of all the loci $\RR_{i,1}(H^*(A))$ does not by
itself bound the number of nonzero Koszul modules.

\begin{example}
\label{ex:euler-hypotheses-independent}
For a $\cdga$ $A$ with $H^*(A)=\bwedge(a)\otimes_{\k}\k[u]$,
where $\abs{a}=1$ and $\abs{u}=2$, one has $S=\k[x]$, and the Koszul chain
complex of $H^*(A)$ splits as the direct sum over $k\ge0$ of the two-term
complexes $S\cdot(au^k)^\vee\xrightarrow{\,x\,}S\cdot(u^k)^\vee$, in
homological degrees $2k+1$ and $2k$.  Hence
$\B_{2k}(H^*(A))\cong\k$ and $\B_{2k+1}(H^*(A))=0$ for every $k\ge0$: all
the loci $\RR_{i,1}(H^*(A))$ reduce to the origin, while the total
$E^2$-page is infinite-dimensional.
\end{example}

%%%%%%%%%%%%%%%%%
\part{Structural properties}
\label{part:structural}
%%%%%%%%%%%%%%%%%

In this part, we study additional structures and functoriality properties of
Koszul modules that sharpen their algebraic and geometric behavior.
We first introduce positive weight decompositions on $\cdgas$, which, 
when present, promote Koszul modules from filtered to genuinely graded 
objects and impose strong linearity constraints on resonance and support loci.
We then analyze the naturality of Koszul modules with respect to $\cdga$ 
morphisms, clarifying the role of base change and identifying precise 
obstructions to functoriality. Together, these results provide the technical 
foundation needed for the comparisons with topological and group-theoretic 
invariants developed in Part~\ref{part:applications}. 

In what follows, we work with the identification $S\cong T$ of 
Section~\ref{subsec:S-and-T} and regard all Koszul modules as 
filtered $S$-modules; the distinction between $S$ and $T$ plays no 
further role.

%%%%%%%%%%%%%%%%%%%%%%%%
\section{Positive weights and the weight spectral sequence}
\label{sect:posWeights}
%%%%%%%%%%%%%%%%%%%%%%%%

Positive weights provide a natural framework in which the Koszul modules
of a genuine $\cdga$ acquire honest gradings, rather than merely filtrations.
In this setting, the algebraic structures introduced in Part~\ref{part:Koszul}
become more rigid and more transparent.

When a $\cdga$ $(A,d)$ admits positive weights, the resulting grading 
forces the support loci $\RR_{i,s}(A)$ and $\RR^{i,s}(A)$ to be invariant 
under a weighted $\k^*$-action, and conical when $H^1(A)$ is concentrated 
in a single weight.  
In favorable cases these loci decompose as finite unions of linear 
subspaces of $H^1(A)$ (coordinate subspaces when the fine multigrading 
of Definition~\ref{def:multihomog} is available), a phenomenon that plays 
a central role in the applications to resonance and holonomy developed 
later in the paper.

%%%%%%%%%%%%%%%%%%%%
\subsection{Positive weights on a $\cdga$}
\label{subsec:weights-cdga}
%%%%%%%%%%%%%%%%%%%%

We begin by recalling the notion of positive weights for commutative
differential graded algebras, which formalizes the idea that the algebra
admits a grading compatible with both multiplication and differential,
and in which all generators carry strictly positive weight.

\begin{definition}
\label{def:pos-weight}
A connected $\cdga$ $(A,d_A)$ over a field $\k$ of characteristic $0$ is said to
have \emph{positive weights} if there exist decompositions
\[
A^i=\boplus_{\alpha>0}\, A^{i,\alpha}
\]
for all $i>0$ such that
\[
A^{i,\alpha}\cdot A^{j,\beta}\subseteq A^{i+j,\alpha+\beta}
\qquad\text{and}\qquad
d_A(A^{i,\alpha})\subseteq A^{i+1,\alpha}
\]
for all integers $i,j> 0$ and $\alpha,\beta>0$. (By convention, $A^0=\k$ has weight $0$.)
In particular, every $\cdga$ $(A,0)$ with zero differential has positive weights: 
assign weight $i$ to $A^i$ for $i>0$ and weight $0$ to $A^0=\k$.
\end{definition}

The notion of positive weights appears in the literature in several
equivalent formulations~\cite{BMSS, Sullivan77, SW-forum}, all valid
over any field $\k$ with $\ch\k=0$.
In the minimal Sullivan model $(\Lambda V, d)$, one assigns a positive
integer weight to each generator and extends multiplicatively, requiring
$d$ to preserve weight; see~\cite[Def.~2.1]{BMSS}.
Equivalently, by~\cite[Props.~2.3 and~2.7]{BMSS},
$(A,d)$ has positive weights if and only if there exists a degree-$0$
derivation $\Psi\colon A\to A$ commuting with $d$ whose eigenvalues on
$A^{>0}$ are strictly positive integers, with the weight decomposition
$A^{i,\alpha}=\ker(\Psi-\alpha\cdot\id)|_{A^i}$ recovering
Definition~\ref{def:pos-weight}. By~\cite[Thm.~2.7]{BMSS}, 
having positive weights is independent
of the ground field: $\Q$-positive weights, $\C$-positive
weights, and $\k$-positive weights are all equivalent.

For a minimal Sullivan $\cdga$ generated in degree~$1$, positive weights
in the sense of Definition~\ref{def:pos-weight} are equivalent to the
existence of positive \emph{Hirsch weights}~\cite{SW-forum}: one assigns
weight $k$ to each generator adjoined at the $k$-th stage of the Hirsch
extension tower, extends multiplicatively, and requires the differential
to preserve total weight.
This formulation is the one most directly relevant to the present work:
the $\CE$ models of nilpotent Lie algebras and the $\cdga$ models of
Sasakian manifolds both arise naturally as Hirsch extensions, and the
weight grading on their Koszul complexes is induced by the Hirsch weights
on the generators.

\begin{remark}
\label{rem:pos-weight-cga}
The same definition applies verbatim to a connected graded-commutative
$\k$-algebra $A^\bullet$ with no differential: one requires decompositions
$A^i=\bigoplus_{\alpha>0}A^{i,\alpha}$ such that $A^{i,\alpha}\cdot
A^{j,\beta}\subseteq A^{i+j,\alpha+\beta}$. We refer to this as a
\emph{positive-weight decomposition} of $A^\bullet$; it amounts to the
case $d=0$ of Definition~\ref{def:pos-weight}.
\end{remark}

We now return to the setup of Definition~\ref{def:pos-weight}.
Assume in addition that $\dim_\k H^1(A)<\infty$.
Choose a weight-homogeneous basis $\{e_1,\dots,e_n\}$ of $H^1(A)$, with
$\wt(e_j)=\alpha_j>0$, and let $\{x_1,\dots,x_n\}$ be the dual basis of
$H_1(A)=H^1(A)^\vee$, with $\wt(x_j)=-\alpha_j$.
Then the symmetric algebra $S=\Sym(H_1(A))=\k[x_1,\dots,x_n]$ inherits a
weighted grading given by
\begin{equation}
\label{eq:wt-x1-k1}
\wt(x_1^{k_1}\cdots x_n^{k_n})=\sum_{j} k_j(-\alpha_j).
\end{equation}
We define the \emph{total weight} on the tensor product $A_\bullet\otimes_{\k} S$ by
\[
\wt(a\otimes s)=\wt(a)+\wt(s).
\]

The canonical element $\omega_A = \sum e_j\otimes x_j$ has total weight $0$, 
so the Koszul differential
$\partial^A = \omega_A^\vee\lrcorner(-) + d_A^\vee\otimes\id_S$ 
preserves total weight.  Hence $K_\bullet(A)=(A_\bullet\otimes_{\k} S,\partial^A)$ 
is a \emph{weight-homogeneous} chain complex of graded $S$-modules, 
and its homology $\B_i(A) = H_i(K_\bullet(A))$ is naturally a \emph{graded} 
$S$-module (not merely filtered).

\begin{proposition}
\label{prop:koszul-separation}
Let $(A,d_A)$ be a connected, finite-type $\k$-$\cdga$ with positive 
weights. Then for each $i\ge 0$, the Koszul module $\B_i(A)$ is a 
graded $S$-module whose $\m$-adic filtration is separated:
\[
\bigcap_{p\ge 0}\, \m^p \B_i(A) = 0.
\]
\end{proposition}

\begin{proof}
Under the finite-type hypothesis, $\B_i(A)$ is a finitely generated
module over $S=\Sym(H_1(A))$, graded by total weight.  The generators of
$\m$ have strictly negative weight, while the weights occurring in
$\B_i(A)$ are bounded above, say by $C$.  Hence every element of
$\m^p\B_i(A)$ has weight at most $C-p$, and so
$\bigcap_{p\ge0}\m^p\B_i(A)=0$.
\end{proof}

\begin{remark}
\label{rem:krull-not-available}
The Krull intersection theorem is not available here: $\m$ is not
contained in the Jacobson radical of $S$.  By
Proposition~\ref{prop:separation-criterion}, separation fails as soon as
some associated prime of the module avoids $\m$, and
Example~\ref{ex:coho-nonseparated} exhibits a finitely generated
$S$-module on which $\m$ acts invertibly.  It is the weight grading, not
Noetherianity alone, that forces separation.
\end{remark}

%%%%%%%%%%%%%%%%%%%%
\subsection{Weight filtrations and Koszul complexes}
\label{subsec:weights-koszul}
%%%%%%%%%%%%%%%%%%%%

The presence of positive weights has strong homological consequences for
Koszul complexes.  Two distinct structures are in play, and it is important
to keep them apart.  The total weight is a \emph{grading} preserved by
$\partial^A$; this is what makes the Koszul modules graded and separated, as in
Proposition~\ref{prop:koszul-separation}.  The linearization of $\partial^A$,
on the other hand, requires a \emph{filtration} that the two summands of
$\partial^A=d_A^\vee+\omega_A^\vee$ treat differently.  Since both summands
preserve the total weight, the total-weight grading cannot serve this purpose;
what does is the weight of the coefficient factor alone.

Accordingly, we filter $K_\bullet(A)=A_\bullet\otimes_{\k}S$ by the weight of
the coefficients:
\begin{equation}
\label{eq:weight-filtration}
W^p K_\bullet(A) = A_\bullet\otimes_{\k} S_{\le -p},
\qquad
S_{\le -p}\coloneqq \boplus_{\wt(s)\le -p} \k\cdot s ,
\end{equation}
the sum being over the monomials $s=x_1^{k_1}\cdots x_n^{k_n}$ of weight
$\wt(s)=-\sum_{j} k_j\alpha_j\le -p$.  We write $S^{(p)}$ for the component of
$S$ of weight exactly $-p$, so that $S_{\le-p}/S_{\le-p-1}=S^{(p)}$ and hence
$\gr^W_pK_\bullet(A)=A_\bullet\otimes_{\k}S^{(p)}$.  This filtration is
\begin{itemize}
\item decreasing, exhaustive, and bounded below, with $W^0K_\bullet(A)=K_\bullet(A)$;

\item Hausdorff: $\bigcap_{p\ge0}W^pK_\bullet(A)=0$;

\item compatible with the $S$-module structure and with the total-weight
grading, each $W^p$ being a total-weight-graded submodule;

\item preserved by $\partial^A$: the term $d_A^\vee\otimes\id_S$ acts on the
first factor and so preserves each $W^p$ exactly, while
$\omega_A^\vee\lrcorner(-)$ raises the coefficient weight by some
$\alpha_j\ge 1$ and therefore maps $W^p$ into $W^{p+1}$.
\end{itemize}
When all the weights $\alpha_j$ are equal---as happens whenever $H^1(A)$ is
weight-homogeneous---the filtration \eqref{eq:weight-filtration} is the
$\m$-adic filtration of Section~\ref{sect:gr-koszul}, reindexed by the common
weight; see Corollary~\ref{cor:weight-splitting}.

Thus the filtration induces a homological spectral sequence
of graded $S$-modules,
\begin{equation}
\label{eq:erpq-ss}
E^r_{p,q}
\;\Longrightarrow\;
\gr^W_p\,\B_{p+q}(A)
\qquad (p\ge 0,\ p+q\ge 0).
\end{equation}

Before analyzing this spectral sequence, we record the finiteness property
that governs its convergence.  It is a statement about the filtered complex
alone, and is used both in Theorem~\ref{thm:weight-ss} below and,
in the form of a splitting of the resulting spectral sequence, in
Corollary~\ref{cor:weight-splitting}.

\begin{lemma}
\label{lem:weight-filtration-finite}
Let $(A,d_A)$ be a connected, finite-type $\cdga$ with positive weights.
Then $K_\bullet(A)$ decomposes as a direct sum, over the total weight
$\gamma\in\Z$, of filtered subcomplexes
\[
K_\bullet(A)=\bigoplus_{\gamma\in\Z}K_\bullet(A)^{[\gamma]},
\qquad
W^pK_\bullet(A)=\bigoplus_{\gamma\in\Z}\bigl(W^pK_\bullet(A)\cap
K_\bullet(A)^{[\gamma]}\bigr),
\]
and for each homological degree $i$ and each $\gamma$ the induced filtration
on $K_i(A)^{[\gamma]}$ is finite.
\end{lemma}

\begin{proof}
Both summands of $\partial^A$ preserve the total weight: $d_A^\vee$ is dual
to a weight-preserving map, and $\omega_A^\vee\lrcorner(-)$ preserves it
because $\omega_A$ has total weight~$0$.  Each $W^p$ is spanned by
total-weight-homogeneous elements, being defined by a condition on the
coefficient weight alone.  This gives the displayed decomposition.

For the finiteness, write $K_i(A)=A_i\otimes_{\k} S$ and note that an element
of $A_i\otimes_{\k} S^{(p)}$ has total weight $-(w+p)$, where $-w$ is the weight
of the factor in $A_i$.  Since $A$ has positive weights and $A^i$ is
finite-dimensional, $w$ ranges over a finite set of positive integers.
Hence, for fixed $i$ and $\gamma$, only finitely many values of $p$ occur, so
$W^pK_i(A)\cap K_i(A)^{[\gamma]}=0$ for $p\gg0$; and $W^0K_\bullet(A)$ is
everything.
\end{proof}

\begin{theorem}[Weight spectral sequence]
\label{thm:weight-ss}
Let $(A,d_A)$ be a connected, finite-type $\cdga$ with positive weights.
The spectral sequence associated to the filtration $W^\bullet$ of
\eqref{eq:weight-filtration} on $K_\bullet(A)$ has
\[
E^0_{p,\bullet}=\bigl(A_\bullet\otimes_{\k}S^{(p)},\,d_A^\vee\otimes\id\bigr),
\qquad
E^1_{p,\bullet}=H_\bullet(A)\otimes_{\k}S^{(p)},
\]
and it converges strongly to the associated graded of the Koszul homology
of $A$,
\[
E^{\infty}_{p,q} \cong \gr^W_p\,\B_{p+q}(A),
\]
the filtration on $\B_\bullet(A)$ being separated by
Proposition~\ref{prop:koszul-separation}.

Assume in addition that $H^1(A)$ is weight-homogeneous, normalized so that
$\alpha_1=\cdots=\alpha_n=1$.  Then the $E_1$-page, equipped with $d_1$, is
canonically isomorphic as a complex of graded $S$-modules to the Koszul
chain complex of the cohomology algebra $H^{*}(A)$ (with zero
differential):
\[
(E^1_{p,\bullet},d_1) \cong K_\bullet \bigl(H^{*}(A)\bigr)^{\mathrm{weight }\, p},
\]
and consequently
\[
E^2_{p,q} \cong \B_{p+q}\bigl(H^{*}(A)\bigr)^{\mathrm{weight }\, p}.
\]
\end{theorem}

\begin{proof}
By construction, $\gr^W_p K_\bullet(A)=A_\bullet\otimes_{\k}S^{(p)}$, where
$S^{(p)}$ denotes the weight $-p$ component of $S$.  As recorded above, the
summand $\omega_A^\vee\lrcorner(-)$ of $\partial^A$ strictly raises the
filtration degree, and therefore induces the zero map on $\gr^W$; the summand
$d_A^\vee\otimes\id_S$ preserves the filtration degree, and is the induced
differential on $\gr^W$.  Hence
\[
E^0_{p,\bullet}=\bigl(A_\bullet\otimes_{\k}S^{(p)},\,d_A^\vee\otimes\id\bigr),
\qquad
E^1_{p,\bullet}=H_\bullet\bigl(A_\bullet,d_A^\vee\bigr)\otimes_{\k}S^{(p)}
= H_\bullet(A)\otimes_{\k}S^{(p)},
\]
and, under the additional hypothesis that $H^1(A)$ be weight-homogeneous of
weight~$1$, the differential $d_1$ induced by the next-order term of
$\partial^A$ is contraction with $\omega_A^\vee$.  That is precisely the
Koszul differential of $K_\bullet(H^{*}(A))$, restricted to the weight $-p$
part of the coefficients, which gives the third assertion; taking
$d_1$-homology gives the fourth.

For convergence, apply Lemma~\ref{lem:weight-filtration-finite}: the
spectral sequence splits as a direct sum over the total weight $\gamma$, and
within each summand the filtration is finite in every homological degree.  A
filtered complex with filtration finite in each degree gives a strongly
convergent spectral sequence, and a direct sum of strongly convergent
spectral sequences is strongly convergent.  The abutment is
$\gr^W\B_\bullet(A)$, and this filtration on $\B_\bullet(A)$ is separated by
Proposition~\ref{prop:koszul-separation}, since \eqref{eq:weight-filtration}
is cofinal with the $\m$-adic filtration: $\m^p\subseteq S_{\le -p}$ and
$S_{\le -p}\subseteq\m^{\lceil p/\alpha_{\max}\rceil}$, where
$\alpha_{\max}=\max_j\alpha_j$.
\end{proof}

\begin{remark}
\label{rem:weight-homog-needed}
The weight-homogeneity hypothesis cannot be dropped from the $d_1$ and
$E^2$ statements.  Contraction with the component of $\omega_A^\vee$ dual
to a class of weight $\alpha_j$ raises the $W$-filtration by $\alpha_j$, so
only the components of weight~$1$ contribute to $d_1$; when the $\alpha_j$
differ, the remaining components appear on later pages.  For
$A=\bwedge(a,b)$ with $d=0$, $\wt(a)=1$ and $\wt(b)=2$, the middle homology
of $(E^1,d_1)$ is nonzero, whereas $\B_1(H^*(A))=0$.  The $E^0$- and
$E^1$-page identifications as graded modules, and the convergence
statement, hold without the hypothesis.
\end{remark}

\begin{remark}
\label{rem:total-weight-not-filtration}
The filtration $W^\bullet$ is by the weight of the coefficients, and not by
the total weight.  Both summands of $\partial^A$ preserve the total weight,
since $\omega_A$ has total weight $0$ and $d_A$ preserves weights; the total
weight therefore grades $K_\bullet(A)$ rather than filtering it, and its role
here is the one described in Proposition~\ref{prop:koszul-separation} and
Corollary~\ref{cor:weight-splitting}.
\end{remark}

The identification of the $E^1$-page in Theorem~\ref{thm:weight-ss}
may be rephrased at the level of associated graded complexes.
We record this consequence explicitly, since it will be used in
subsequent comparison arguments.

\begin{theorem}
\label{thm:koszul-weight-graded}
Let $(A,d_A)$ be a connected, finite-type $\cdga$ with positive weights, and
assume that $H^1(A)$ is weight-homogeneous, normalized to weight~$1$.
Then the $E^1$-page of the weight spectral sequence, equipped with the
differential $d_1$, is canonically isomorphic to the Koszul complex of the
cohomology algebra:
\[
\bigl(E^1_{\bullet,\bullet},\, d_1\bigr) \cong K_\bullet(H^*(A)).
\]
\end{theorem}

\begin{proof}
This is the $E^1$-page assertion of Theorem~\ref{thm:weight-ss}, together
with the identification $H_\bullet(A_\bullet,d_A^\vee)=H_\bullet(A)$ of the
$E^0$-page homology.  The hypothesis on $H^1(A)$ is what makes $d_1$ the
Koszul differential of $K_\bullet(H^*(A))$; see
Remark~\ref{rem:weight-homog-needed}.
\end{proof}

Here $E^1_{\bullet,\bullet}$ is the homology of the associated graded complex
$\gr^W K_\bullet(A)$ with respect to the induced differential
$d_0=d_A^\vee\otimes\id$.  It is that homology, and not
$\gr^W K_\bullet(A)$ itself, that is identified with $K_\bullet(H^*(A))$: the
stronger identification would force $\gr^W\B_\bullet(A)\cong\B_\bullet(H^*(A))$
unconditionally, which is false, as Remark~\ref{rem:linearization-one-page}
records.

The identification of Theorem~\ref{thm:koszul-weight-graded} is a
statement about associated graded complexes.  It does not by itself
lift to a comparison of the completed Koszul complexes of $A$ and of
$H^*(A)$: such a comparison requires a filtered morphism between them,
and producing one is precisely a formality statement.  See
Corollary~\ref{cor:koszul-completion-formal} for the comparison that is
available, and Remark~\ref{rem:pw-not-enough} for why positive weights
do not substitute for formality.

Finally, we record how the total-weight grading interacts with the spectral
sequence, and when the latter coincides with its $\m$-adic counterpart.

\begin{corollary}
\label{cor:weight-splitting}
Let $(A,d_A)$ be a connected, finite-type $\cdga$ with positive weights.
\begin{enumerate}[itemsep=2pt]
\item \label{ws1}
The spectral sequence of Theorem~\ref{thm:weight-ss} splits as a direct sum
of spectral sequences indexed by the total weight $\gamma\in\Z$, each of
which is the spectral sequence of a filtered complex whose filtration is
finite in every homological degree.
\item \label{ws2}
If $H^1(A)$ is weight-homogeneous, say $\alpha_1=\cdots=\alpha_n=\alpha$,
then $W^{p}K_\bullet(A)=A_\bullet\otimes_{\k}\m^{\lceil p/\alpha\rceil}$, so
\eqref{eq:weight-filtration} is the $\m$-adic filtration up to reindexing, and
the spectral sequence of Theorem~\ref{thm:weight-ss} is the $\m$-adic Koszul
spectral sequence of Theorem~\ref{thm:koszul-ss-hom}, refined by the
total-weight splitting of~\eqref{ws1}.
\end{enumerate}
\end{corollary}

\begin{proof}
\eqref{ws1} This is Lemma~\ref{lem:weight-filtration-finite}, together with
the fact that the spectral sequence of a direct sum of filtered complexes is
the direct sum of their spectral sequences.

\eqref{ws2} If all $\alpha_j=\alpha$, then a monomial of $S$ has weight
$-\alpha k$ exactly when it has degree $k$, so $S_{\le -p}=\m^{\lceil
p/\alpha\rceil}$.
\end{proof}

The hypothesis of \eqref{ws2} is satisfied in all the situations considered
later in this paper: for the Chevalley--Eilenberg $\cdga$ of a nilpotent Lie
algebra $\g$ one has $H^1(\CE(\g))=(\g/[\g,\g])^\vee$, concentrated in
bracket length $1$; and for the Bibby models of
Section~\ref{subsec:conf-elliptic}, $H^1(A(n))$ is spanned by the classes
$[a_i],[b_i]$, all of weight~$1$.  Note that it is only $H^1(A)$, and not
$A^1$, that is required to be weight-homogeneous.

The filtration considered in Theorem~\ref{thm:weight-ss} is formally analogous 
to the weight filtration arising from a positively weighted Sullivan minimal model 
in rational homotopy theory. In that setting, the associated spectral sequence 
computes the graded Lie algebra of rational homotopy groups from the cohomology 
algebra. The present construction is the purely algebraic counterpart: the 
$E^1$-page depends only on $H^{*}(A)$, while higher differentials measure 
the deviation of $A$ from formality.

\begin{remark}
\label{rem:res-linear-open}
It is natural to ask whether the existence of positive weights alone 
imposes linearity constraints on resonance varieties.  
More precisely, if a finite-type $\cdga$ $(A,d_A)$ carries a system of 
positive weights (without necessarily arising from a geometric model), 
must each resonance variety $\RR^i(A)$ be a finite union of linear subspaces?  
To the best of our knowledge, this question remains open (see \cite{Su-indam}).
\end{remark}

%%%%%%%%%%%%%%%%%%%%
\subsection{Multigradings and conicality of support loci}
\label{subsec:multigrading}
%%%%%%%%%%%%%%%%%%%%

Let $(A,d_A)$ be a connected positive-weight $\cdga$ with $\dim_\k H^1(A)<\infty$.  
Fix a weight-homogeneous basis $\{e_1,\dots,e_n\}$ of $H^1(A)$ of weights 
$\alpha_1,\dots,\alpha_n>0$ and let $\{x_1,\dots,x_n\}$ be the dual basis of 
$H_1(A)$, of weights $-\alpha_1,\dots,-\alpha_n$.  The total weight then 
defines a single $\Z$-grading on $K_\bullet(A)=A_\bullet\otimes_{\k} S$ which is 
preserved by $\partial^A$; by Corollary~\ref{cor:weight-splitting}\eqref{ws1} 
it splits the weight spectral sequence, and it is all that is needed for the 
conicality statement below.

Finer gradings are available whenever the $\cdga$ itself carries one.  Nothing
in the construction of $K_\bullet(A)$ requires the grading group to be $\Z$,
and it is worth recording the general form, since the two instances that
arise later sit at opposite extremes: the fine $\N^n$-grading of a monomial
algebra, and the $\GL(W)\times\GL(U_0)$-weight grading of
\S\ref{subsec:b1-conf-structure}.

\begin{definition}
\label{def:multihomog}
Let $\Gamma$ be a finitely generated free abelian group.  A connected $\cdga$
$(A,d_A)$ is \emph{$\Gamma$-multihomogeneous} if $A^\bullet$ admits a
decomposition $A^i=\boplus_{\mathbf{a}\in\Gamma}A^i_{\mathbf{a}}$ such that
\[
A^i_{\mathbf{a}}\cdot A^{j}_{\mathbf{b}}\subseteq A^{i+j}_{\mathbf{a}+\mathbf{b}},
\qquad 
d_A\bigl(A^i_{\mathbf{a}}\bigr)\subseteq A^{i+1}_{\mathbf{a}},
\]
with $A^0=A^0_{\mathbf{0}}=\k$, and if there is a homomorphism
$\lambda\colon\Gamma\to\Z$ with
\begin{equation}
\label{eq:multiwt-positive}
\lambda(\mathbf{a})>0
\qquad\text{whenever } A^i_{\mathbf{a}}\ne0 \text{ for some } i\ge1 .
\end{equation}
The grading is \emph{normalized} if, in addition,
$\lambda(\mathbf{a})=1$ whenever $H^1(A)_{\mathbf{a}}\ne0$.
\end{definition}

Condition \eqref{eq:multiwt-positive} says that the multidegrees occurring in
$A^{\ge1}$ lie in an open half-space, so that the cone they generate is
strictly convex.  It has two consequences.  First, composing the grading with
$\lambda$ gives $A$ positive weights in the sense of
Definition~\ref{def:pos-weight}, and normalization is exactly the hypothesis
that $H^1(A)$ be weight-homogeneous of weight~$1$; the case $\Gamma=\Z$ with
$\lambda=\mathrm{id}$ recovers Definition~\ref{def:pos-weight} verbatim.  Second, it
makes the multigraded pieces of $S=\Sym(H_1(A))$ finite-dimensional: a
monomial of multidegree $\mathbf{b}$ has at most $\lambda(\mathbf{b})$
factors, so only finitely many monomials occur in each multidegree, and the
multigraded Hilbert function of a finitely generated multigraded $S$-module is
defined.  Without \eqref{eq:multiwt-positive} both fail.

At one extreme, $\Gamma=\Z$ and $\lambda=\mathrm{id}$ impose no condition beyond
positive weights.  At the other, when $\dim_\k H^1(A)=n$ one may ask for
$\Gamma=\Z^n$ and a basis $\{e_1,\dots,e_n\}$ of $H^1(A)$ with
$e_j\in A^1_{\mathbf{e}_j}$, where $\mathbf{e}_1,\dots,\mathbf{e}_n$ are the
standard basis vectors; we call this the \emph{fine} grading, and take
$\lambda(\mathbf{a})=\sum_{j}a_j$, which is normalized.  The basic examples of
the fine grading are the monomial quotients
$A=\bwedge(e_1,\dots,e_n)/J$ with $d_A=0$ and $J$ generated by
monomials, in particular the exterior Stanley--Reisner rings of
Example~\ref{ex:sr-ring}.

\begin{theorem}
\label{prop:multigraded-koszul}
Let $\Gamma$ be a finitely generated free abelian group and let $(A,d_A)$ be a
connected, finite-type, $\Gamma$-multihomogeneous $\cdga$, as in
Definition~\ref{def:multihomog}.  Grade $A_i=(A^i)^\vee$ by placing the
subspace dual to $A^i_{\mathbf{a}}$ in multidegree $\mathbf{a}$, and grade
$S=\Sym(H_1(A))$ by placing the dual of $H^1(A)_{\mathbf{a}}$ in multidegree
$\mathbf{a}$.  Then the following hold.
\begin{enumerate}
\item \label{mg1}
The Koszul differential $\partial^A$ on $K_\bullet(A)=A_\bullet\otimes_{\k}S$ 
is multihomogeneous of multidegree $\mathbf{0}$.
\item \label{mg2}
The Koszul homology modules $\B_i(A)$ are $\Gamma$-multigraded $S$-modules,
with multigraded pieces of finite dimension over $\k$.
\item \label{mg3}
If the grading is normalized, then the associated graded of the
coefficient-weight filtration \eqref{eq:weight-filtration}, equipped with the
induced differential $d_1$, satisfies
\[
\bigl(\gr^W_\bullet K_\bullet(A),\,d_1\bigr) \cong K_\bullet\bigl(H^{*}(A)\bigr)
\]
as $\Gamma$-multigraded chain complexes of $S$-modules.
\end{enumerate}
\end{theorem}

\begin{proof}
Since $A$ is connected, $A^0=\k$ gives $B^1(A)=0$, so
$H^1(A)=Z^1(A)=\ker\bigl(d_A|_{A^1}\bigr)$; as $d_A$ preserves multidegrees,
this is a multigraded subspace of $A^1$, and the grading on $S$ is well
defined.  Choose a basis of $H^1(A)$ \emph{adapted} to the decomposition, each
member lying in a single $H^1(A)_{\mathbf{a}}$; the dual basis of $H_1(A)$ is
then adapted to the dual decomposition.  The canonical element
$\omega_A\in H^1(A)\otimes_{\k} H_1(A)$ is independent of this choice, and in an
adapted basis it is a sum of terms $e\otimes x_e$ with $e\in A^1_{\mathbf{a}}$
and $x_e$ of multidegree $\mathbf{a}$.  Contraction by $e$ is dual to
multiplication by $e$, hence lowers multidegree by $\mathbf{a}$, while
tensoring with $x_e$ raises it by $\mathbf{a}$; so each term of
$\omega_A^\vee\lrcorner(-)$, and therefore the whole of it, has multidegree
$\mathbf{0}$.  Likewise $d_A^\vee$ preserves multidegrees, since $d_A$ does.
This proves \eqref{mg1}, and \eqref{mg2} follows by passing to homology, the
finite-dimensionality of the pieces being the second consequence of
\eqref{eq:multiwt-positive} noted above.

For \eqref{mg3}, composing the grading with $\lambda$ gives $A$ positive
weights, and normalization makes $H^1(A)$ weight-homogeneous of weight~$1$; the
coefficient-weight filtration \eqref{eq:weight-filtration} is the one induced
by $\lambda$.  Theorem~\ref{thm:koszul-weight-graded} therefore applies, and
both the identification and the induced differentials are multihomogeneous
by \eqref{mg1}.
\end{proof}

For $\Gamma=\Z$ and $\lambda=\mathrm{id}$, Definition~\ref{def:multihomog} is
Definition~\ref{def:pos-weight} and
Theorem~\ref{prop:multigraded-koszul}\eqref{mg3} is
Theorem~\ref{thm:koszul-weight-graded}; the content of the multigraded
statement is that the whole of $\Gamma$ is carried along.

\begin{remark}
\label{rem:multigrading-fails}
Every positive-weight $\cdga$ is $\Z$-multihomogeneous, so
Theorem~\ref{prop:multigraded-koszul} always applies with $\Gamma=\Z$; what is
rare is a \emph{fine} grading.  Two distinct mechanisms obstruct that one.

First, the multiplicative relations may fail to be multihomogeneous for the
fine grading, even when $d_A=0$.  This happens for the 
Orlik--Solomon algebra $A=\OS(\A)$ of any arrangement having a rank-$2$ flat 
of multiplicity at least $3$: for the braid arrangement $\A_3$, the flat 
$\{12,13,23\}$ yields the relation 
$e_{13}e_{23}=e_{12}e_{23}-e_{12}e_{13}$ in $A^2$, whose three terms have 
pairwise distinct multidegrees $\mathbf{e}_{13}+\mathbf{e}_{23}$, 
$\mathbf{e}_{12}+\mathbf{e}_{23}$, and $\mathbf{e}_{12}+\mathbf{e}_{13}$.  
Were $A^2$ multigraded, each of those terms would have to vanish separately.

Second, the differential may fail to be multihomogeneous.  For the Bibby 
model $(A(n),d)$ of $\Conf(\E,n)$ from Section~\ref{subsec:conf-elliptic}, 
the basis $\{[a_i],[b_i]\}$ of $H^1(A(n))$ is weight-homogeneous, yet 
$d e_{ij}=(a_i-a_j)(b_i-b_j)$ expands into four terms of pairwise distinct 
multidegrees.  Here the obstruction is visible in the Koszul module itself: 
by Theorem~\ref{thm:BA-n}, 
$\B_1(A(n))\cong\boplus_{i<j}S/I_{ij}$ with $x_i+x_j\in I_{ij}$, so the 
associated primes of $\B_1(A(n))$ are the ideals $I_{ij}$, none of which is 
generated by variables.  Since every associated prime of an 
$\N^n$-multigraded module over $\k[x_1,\dots,x_n]$ is generated by a subset 
of the variables, $\B_1(A(n))$ admits no $\N^n$-multigrading.  For $n\ge3$ 
this conclusion is independent of the chosen basis of $H^1(A(n))$: were 
some basis to render all the $I_{ij}$ monomial, each $2$-plane 
$\bV(I_{ij})$ would be a coordinate subspace, forcing the pairwise 
non-proportional vectors $x_1-x_2$, $x_2-x_3$, and $x_1-x_3$ to be distinct 
members of a single basis---impossible, as they are linearly dependent.  
Equivalently: the irreducible components of $\RR_1(A(n))$ are linear, but 
they are not coordinate subspaces.

Neither obstruction rules out a coarser grading.  For $A(n)$ the
$\GL(W)\times\GL(U_0)$-weight grading of \S\ref{subsec:b1-conf-structure} is
available on the cohomology side, with $\lambda$ the total degree in the
$W$-factor; and Corollary~\ref{cor:conical} below uses only the single
$\Z$-grading by total weight, so it holds for all positive-weight
$\cdgas$.
\end{remark}

\begin{remark}
\label{rem:bigraded-rs}
The case $\Gamma=\Z^2$ deserves separate mention.  Let $V=V_1\oplus V_2$ with
$\dim V_i\ge2$, let $K\subseteq\bwedge^2V$, and let $A$ be the quadratic
algebra $\bwedge V^\vee/(K^\perp)$ with $d_A=0$, bigraded by placing
$V_1^\vee$ in bidegree $(1,0)$ and $V_2^\vee$ in bidegree $(0,1)$, so that
$A^2=\bwedge^2V^\vee/K^\perp$.  The pairing between
\[
\bwedge^2V=\bwedge^2V_1\oplus(V_1\otimes_{\k} V_2)\oplus\bwedge^2V_2
\qquad\text{and}\qquad
\bwedge^2V^\vee=\bwedge^2V_1^\vee\oplus(V_1^\vee\otimes_{\k} V_2^\vee)
\oplus\bwedge^2V_2^\vee
\]
is perfect on each summand and zero across distinct summands.  Hence $K^\perp$
is bihomogeneous if and only if $K$ is, and in that case $A$ is
$\Z^2$-multihomogeneous, with $\lambda(m_1,m_2)=m_1+m_2$ normalized.

Raicu and Sam \cite{RS} work under the stronger hypothesis
$K\subseteq V_1\otimes_{\k} V_2$.  By the same block-diagonality, taking
annihilators inside $\bwedge^2V$ turns the two inclusions
$\bwedge^2V_1^\vee\subseteq K^\perp$ and
$\bwedge^2V_2^\vee\subseteq K^\perp$ into
$K\subseteq(V_1\otimes_{\k} V_2)\oplus\bwedge^2V_2$ and
$K\subseteq\bwedge^2V_1\oplus(V_1\otimes_{\k} V_2)$, whose intersection is
$V_1\otimes_{\k} V_2$; so $K\subseteq V_1\otimes_{\k} V_2$ holds if and only if $A^2$ is
concentrated in bidegree $(1,1)$, equivalently if and only if $V_1^\vee$ and
$V_2^\vee$ are isotropic for the cup product of $A$.  This is strictly
stronger than bihomogeneity of $K$, and it is what forces $\B_1(A)$ to have
infinite length once $\dim V_i\ge2$: any $a\wedge b$ with $a,b\in V_1^\vee$
independent is a nonzero decomposable lying in $\bwedge^2V_1^\vee\subseteq
K^\perp$.  Accordingly \cite{RS} replace finite length by the vanishing of the
associated sheaf on $\mathbb{P}V_1\times\mathbb{P}V_2$, which holds precisely
when $K^\perp$ contains no nonzero tensor of rank at most~$2$
\cite[Prop.~3.5]{RS}, and they make that vanishing effective
\cite[Thm.~1.8]{RS}; the Grassmannian $\mathrm{Gr}_2(V^\vee)$ of the
single-graded theory is thereby replaced by the secant variety of the Segre
embedding.  Farkas \cite[Quest.~5.8]{Farkas} asks for the
topological counterpart of this theory.  Theorem~\ref{prop:multigraded-koszul}
provides one at the level of $\cdga$ models: the Chen ranks of a group with a
$\Gamma$-multihomogeneous $1$-finite $1$-model refine into a $\Gamma$-graded
Hilbert function, computed by
Theorem~\ref{prop:multigraded-koszul}\eqref{mg3} from the cohomology algebra.
What is not addressed here is whether the specific pairs $(V,K)$ arising in
\cite{RS} are realized by $\cdga$ models of spaces.
\end{remark}

\begin{corollary}
\label{cor:conical}
Let $(A,d_A)$ be a connected positive-weight $\cdga$ with 
$\dim_\k H^1(A)<\infty$, and let $\alpha_1,\dots,\alpha_n>0$ be the weights 
of a weight-homogeneous basis $e_1,\dots,e_n$ of $H^1(A)$.  Then, for all 
$i,s$, the two support loci
\[
\Supp_S\bigl(\textstyle\bigwedge^s \B_i(A)\bigr)
\quad\text{and}\quad
\Supp_S\bigl(\textstyle\bigwedge^s \B^i(A)\bigr)
\]
are invariant under the weighted $\k^\ast$-action 
$\lambda\cdot e_j=\lambda^{\alpha_j}e_j$ on $H^1(A)$; equivalently, they are 
cut out by weighted-homogeneous ideals.  If moreover $H^1(A)$ is 
concentrated in a single weight, then these loci are conical.
\end{corollary}

\begin{proof}
The $\k^\ast$-action on $H^1(A)$ by $\lambda\cdot e_j = \lambda^{\alpha_j} e_j$
dualizes to an action on $H_1(A)$ by $\lambda\cdot x_j=\lambda^{-\alpha_j}x_j$,
which extends multiplicatively to $S=\Sym(H_1(A))$.
Since $\omega_A = \sum e_j \otimes x_j$ has total weight $0$ under the
combined action, this $\k^\ast$-action commutes with the differentials of
both the homological Koszul complex $K_\bullet(A)$ and its cohomological
dual $K^\bullet(A)$.  Hence the modules $\B_i(A)$ and $\B^i(A)$ are 
equivariant, so their annihilators are homogeneous and the two support loci 
are stable under the $\k^\ast$-action.  Finally, if all the weights 
$\alpha_j$ coincide, then the weighted action is the standard scaling 
action, weighted-homogeneous ideals are homogeneous, and the loci are 
cones.
\end{proof}

\begin{remark}
\label{rem:weighted-not-conical}
The hypothesis that $H^1(A)$ be concentrated in a single weight cannot be 
omitted from the last assertion: a weighted-homogeneous subvariety need not 
be a cone.  For instance, with weights $\alpha_1=1$ and $\alpha_2=2$ on a 
$2$-dimensional space, the curve $\bV(y-x^2)$ is invariant under 
$\lambda\cdot(x,y)=(\lambda x,\lambda^2y)$, yet it is not conical.  The 
hypothesis is automatic when $\dim_\k H^1(A)=1$, and for a graded algebra 
with zero differential and weight equal to cohomological degree, where 
$H^1$ has weight $1$.
\end{remark}

This weighted invariance, and the conicality that follows from it in the 
single-weight case, is extended to the jump loci $\RR^{i,s}(A)$ in 
Theorem~\ref{thm:resonance-inclusion}, once the determinantal description 
of those loci is available; there, and in 
Proposition~\ref{prop:generic-resonance-agreement}, it is 
combined with the tangent cone inclusion of 
Theorem~\ref{thm:tangent-cone} to obtain sharp constraints on 
resonance varieties of $\cdgas$ with positive weights.

The use of $\N^n$-multigraded $S$-modules in this context is inspired by the
theory of squarefree multigraded modules developed by Yanagawa
\cite{Yanagawa00}, where Alexander duality for Stanley--Reisner rings is
formulated in terms of fine multigradings.
That framework---and with it the multigraded results of \cite{AFRSS24} used 
below---applies to the $\cdgas$ carrying the fine grading of 
Definition~\ref{def:multihomog}, in particular to the exterior 
Stanley--Reisner rings of Example~\ref{ex:sr-ring}, rather than to 
positive-weight $\cdgas$ in general.

%%%%%%%%%%%%%%%%%%%%
\subsection{Weighted Hilbert series and degeneration}
\label{subsec:hilbert-series}
%%%%%%%%%%%%%%%%%%%%

The total-weight grading provided by positive weights allows one to refine the
coefficientwise Hilbert-series inequalities of
Section~\ref{sect:gr-koszul}, and to characterize
degeneration of the weight spectral sequence in purely numerical terms.

Retain the setting above.  By Theorem~\ref{thm:weight-ss}, the weight 
spectral sequence has
\[
E^2_{p,q} \cong \B_{p+q}\bigl(H^{*}(A)\bigr)^{\text{weight }p}
\qquad\text{and}\qquad 
E^\infty_{p,q} \cong \gr^W_p\,\B_{p+q}(A) .
\]
All these $\k$-vector spaces are finite-dimensional, since the weights of $A$ 
are positive, those of $S$ are nonpositive, and $\dim_\k H^1(A)<\infty$.  
Consequently, the \emph{weighted Hilbert series}
\[
\Hilb_W(M,t)\coloneqq \sum_{p} \dim_\k \bigl(\gr^W_p M\bigr)\, t^{p}
\]
is well defined for $M=\B_i(A)$ and for $M=\B_i(H^{*}(A))$; in the latter case
the filtration is split by the grading, so $\gr^W_p\B_i(H^*(A))=
\B_i(H^*(A))^{\text{weight }p}$.

\begin{theorem}
\label{thm:Hilbert-inequality}
Let $(A,d_A)$ be a connected positive-weight $\cdga$ with $\dim_\k H^1(A)<\infty$,
and assume $H^1(A)$ is weight-homogeneous, normalized to weight~$1$.
For every $i\ge 0$ there is a coefficientwise inequality of weighted Hilbert 
series,
\[
\Hilb_W\bigl(\B_i(A),t\bigr) \preccurlyeq 
\Hilb_W\bigl(\B_i(H^{*}(A)),t\bigr),
\]
with equality in degree $i$ if and only if no differential $d^r$ with $r\ge2$ 
enters or leaves the terms $E^r_{p,i-p}$.  In particular, equality holds for 
all $i$ if and only if the weight spectral sequence degenerates, that is, 
$E^2\cong E^\infty$.

If, in addition, $A$ is $\Gamma$-multihomogeneous in the sense of 
Definition~\ref{def:multihomog}, then the same statements hold for the 
$\Gamma$-graded Hilbert series $\Hilb(-,\mathbf{t})$ in place of 
$\Hilb_W(-,t)$.
\end{theorem}

\begin{proof}
Each $E^\infty_{p,\,i-p}$ is an iterated subquotient of $E^2_{p,\,i-p}$, so
\begin{equation}
\label{eq:einfty-compare}
\dim_\k \gr^W_p\,\B_i(A) = \dim_\k E^\infty_{p,\,i-p} \le
\dim_\k E^2_{p,\,i-p} =
\dim_\k \B_i(H^{*}(A))^{\text{weight }p},
\end{equation}
the last equality being the $E^2$-identification of
Theorem~\ref{thm:weight-ss}, for which the weight-homogeneity hypothesis is
needed.  Equality holds for every $p$ precisely when no higher differential
affects the line $p+q=i$.  Multiplying by $t^p$ and summing over $p$ gives the 
inequality.  When $A$ is $\Gamma$-multihomogeneous, 
Theorem~\ref{prop:multigraded-koszul} shows that the weight spectral 
sequence is one of $\Gamma$-multigraded $S$-modules, so the comparison 
\eqref{eq:einfty-compare} may be carried out separately in each multidegree 
$\mathbf{a}\in\Gamma$.
\end{proof}

\begin{corollary}
\label{cor:degeneration-criteria}
Let $(A,d_A)$ be a connected positive-weight $\cdga$ with $\dim_\k H^1(A)<\infty$
and $H^1(A)$ weight-homogeneous, normalized to weight~$1$.
The following are equivalent:
\begin{enumerate}[itemsep=2pt]
\item \label{wh1}
The weight spectral sequence degenerates, i.e., $E^2\cong E^\infty$.
\item \label{wh2}
$\Hilb_W(\B_i(A),t) = \Hilb_W(\B_i(H^{*}(A)),t)$ for all $i$.
\item \label{wh3}
$\gr^W_\bullet \B_i(A) \cong \B_i(H^{*}(A))$ as weight-graded 
$S$-modules for all $i$.
\end{enumerate}
In particular, if $d_A=0$ then $\B_i(A)=\B_i(H^*(A))$ on the nose, so all 
three conditions hold.  If moreover $A$ is $\Gamma$-multihomogeneous, then the 
isomorphisms in~\eqref{wh3} and the Hilbert series in~\eqref{wh2} may be 
taken to be $\Gamma$-multigraded.
\end{corollary}

\begin{proof}
\eqref{wh1} $\Rightarrow$ \eqref{wh3}: If the spectral sequence degenerates, then
$E^2_{p,q} = E^\infty_{p,q}$, so $\gr^W_p \B_{p+q}(A) \cong 
\B_{p+q}(H^*(A))^{\text{weight }p}$ for all $p,q$. Summing over $p$ 
gives item~\eqref{wh3}.

\eqref{wh3} $\Rightarrow$ \eqref{wh2}: An isomorphism of weight-graded $S$-modules 
implies equality of the corresponding weighted Hilbert series.

\eqref{wh2} $\Rightarrow$ \eqref{wh1} : By Theorem~\ref{thm:Hilbert-inequality}, equality 
of Hilbert series holds if and only if there are no nonzero differentials 
on any page $E^r$ with $r \ge 2$.
\end{proof}

\begin{example}
\label{ex:sr-ring}
Let $\Delta$ be a simplicial complex on vertex set $[n]$. 
Recall that the \emph{exterior Stanley--Reisner ring} of $\Delta$ is 
$A = \k\langle\Delta\rangle = \bigwedge(e_1,\dots,e_n)/J_\Delta$, 
where $J_\Delta$ is generated by the monomials $e_{j_1}\wedge\cdots
\wedge e_{j_s}$ corresponding to non-faces of $\Delta$; 
see \cite{Yanagawa00} for the squarefree $\N^n$-graded framework 
used here. Since $J_\Delta$ is a monomial ideal and $d_A=0$, the $\cdga$ $A$ 
carries the fine grading of Definition~\ref{def:multihomog}: the 
space $A^i$ is spanned by the monomials $e_\sigma$ with $\sigma\in\Delta$ of 
dimension $i-1$, and $e_\sigma$ has multidegree 
$\sum_{j\in\sigma}\mathbf{e}_j$. Thus 
Theorem~\ref{prop:multigraded-koszul} applies; moreover 
Corollary~\ref{cor:degeneration-criteria} 
gives $\B_i(A) \cong \B_i(H^{*}(A))$ as multigraded $S$-modules.

Writing $\Delta_V$ for the restriction of $\Delta$ to a subset 
$V\subseteq[n]$ and $\mathbf{t}^V=\prod_{j\in V}t_j$, the Koszul modules 
$\B_i(A)$ are squarefree in the sense of \cite{Yanagawa00}, and, for every 
$i\ge1$, their multigraded Hilbert series is computed in 
\cite[Thm.~1.2]{AFRSS24}:
\[
\Hilb\bigl(\B_i(A),\mathbf{t}\bigr) = 
\sum_{V\subseteq [n]}
\dim_\k \tilde{H}_{i-1}(\Delta_V;\k)\, 
\frac{\mathbf{t}^V}{\prod_{j\in V}(1-t_j)} .
\]
Specializing to the single grading by total multidegree gives
\[
\Hilb\bigl(\B_i(A),t\bigr) = 
\sum_{V\subseteq [n]} \dim_\k \tilde{H}_{i-1}(\Delta_V;\k)
\Bigl(\frac{t}{1-t}\Bigr)^{\abs{V}} ;
\]
note that $A_i$ sits in total multidegrees of size $i$, so that these series 
differ by a factor of $t^{i}$ from the ones computed with respect to the 
internal $S$-degree used elsewhere in this paper.

Finally, by \cite[Thm.~1.3]{AFRSS24} the support loci are reduced, with 
\[
\RR_{i,1}(A) = \bcup \, \k^{V}
\qquad (i\ge 1),
\]
the union being over those subsets $V\subseteq[n]$ that are maximal with 
$\tilde{H}_{i-1}(\Delta_V;\k)\ne0$. In particular, all irreducible components 
are coordinate subspaces---as the fine multigrading predicts, and in contrast 
with the situation of Remark~\ref{rem:multigrading-fails}.
\end{example}

Theorem~\ref{thm:Hilbert-inequality} is the algebraic counterpart of the upper 
semicontinuity of Betti numbers in flat families of quadratic algebras (e.g., 
the degeneration $A \rightsquigarrow H^{*}(A)$ induced by the 
total-weight filtration).

%%%%%%%%%%%%%%%%%%%%%%%%%%%%%%%
\section{Naturality properties of Koszul modules}
\label{sect:naturality-koszul}
%%%%%%%%%%%%%%%%%%%%%%%%%%%%%%%

The Koszul modules associated to a $\cdga$ are functorial only in a restricted
sense. Since their definition depends on the symmetric algebra
$S=\Sym(H_1(A))$, a morphism of $\cdgas$ typically induces a comparison map
only after a suitable change of scalars. In this section, we make precise the
naturality properties of Koszul modules under $\cdga$ morphisms, identify
situations in which stronger functoriality holds, and explain how failures of
naturality reflect obstructions to formality, both for algebras and for maps.

%%%%%%%%%%%%%%%%%%%%%
\subsection{Naturality of Koszul modules}
\label{subsec:naturality-koszul}
%%%%%%%%%%%%%%%%%%%%%

We begin by describing the most general functoriality enjoyed by Koszul
modules. A morphism of $\cdgas$ induces a natural comparison map between 
the corresponding Koszul complexes only after restriction of scalars along the
induced map on symmetric algebras. Although this construction is formally
functorial, it is often too weak for applications, motivating the stronger
notions considered in subsequent subsections.

Let $(A^{\bullet},d_A)$ be a connected, graded-commutative $\k$-$\cdga$ 
with $\dim_{\k} A^{1}<\infty$.  
Set $S=\Sym(H_1(A))$, where $H_1(A) = (H^1(A))^{\vee}$.
In Section~\ref{subsec:hom-koszul} we defined the Koszul modules 
of $A$ by $\B_{i}(A) \coloneqq H_{i}(K_{\bullet}(A))$, where 
$K_{\bullet}(A) = (A_\bullet \otimes_{\k} S, \partial^A)$ is the 
Koszul chain complex of that section.

Let $\varphi \colon A \to \oA$ be a morphism of connected 
$\k$-$\cdgas$ with finite-dimensional $H^1$.  
The map $\varphi$ induces linear maps
\begin{equation}
\label{eq:h1a-h1a-prime}
H^1(\varphi)\colon  H^1(A) \longrightarrow H^1(\oA),
 \qquad
(H^1(\varphi))^{\vee}\colon  H_1(\oA) \longrightarrow H_1(A),
\end{equation}
and hence a morphism of symmetric algebras 
$\varphi_S \colon \oS=\Sym(H_1(\oA)) \to S=\Sym(H_1(A))$.
Since $\varphi$ is a $\cdga$ map $A\to\oA$, it induces a map on 
Koszul complexes in the \emph{opposite} direction after restriction 
of scalars along $\varphi_S$: a class in $K_\bullet(\oA)$ is pulled 
back to $K_\bullet(A)$ via $\varphi$, giving a chain map of 
$S$-modules
\begin{equation}
\label{eq:kaprime-ka}
K_{\bullet}(\oA) \otimes_{\oS} S \longrightarrow K_{\bullet}(A).
\end{equation}
Thus, for each $i\ge 0$, the assignment $A \mapsto \B_{i}(A)$
is contravariantly functorial: a $\cdga$ morphism
$\varphi \colon A \to \oA$ induces an $S$-linear map
\begin{equation}
\label{eq:bia-s-bia}
\B_{i}(\varphi) \colon \B_{i}(\oA) \otimes_{\oS} S\longrightarrow \B_{i}(A),
\end{equation}
natural in $\varphi$. Since the coefficient ring $S = \Sym(H_1(A))$
varies with $A$, this functoriality necessarily involves a change of
base rings, and takes the form of restriction of scalars along
$\varphi_S \colon \oS \to S$ rather than a functor into a single
module category.

%%%%%%%%%%%%%%%%%%%%%%%%%%%%%%%%%%
\subsection{Functoriality under isomorphisms on $H^1$}
\label{subsec:functorial-iso}
%%%%%%%%%%%%%%%%%%%%%%%%%%%%%%%%%%

When a $\cdga$ morphism $\varphi \colon A \to \oA$ induces an isomorphism 
on $H^1$, the isomorphism $H^1(\varphi)$ can be inverted to yield a 
covariant chain map $K_\bullet(A) \to K_\bullet(\oA)$ over the common 
ground ring $S \cong \oS$, in contrast to the contravariant 
restriction-of-scalars map \eqref{eq:kaprime-ka} of 
Section~\ref{subsec:naturality-koszul}.  
Set $S = \Sym(H_1(A))$, $\oS = \Sym(H_1(\oA))$, and let
\begin{equation}
\label{eq:Phi-map}
\Phi \coloneqq \Sym\bigl(H_1(\varphi)^{-1}\bigr) \colon S \longrightarrow \oS
\end{equation}
be the induced isomorphism of symmetric algebras.

\begin{proposition}
\label{prop:koszul-functorial}
Let $\varphi \colon (A,d_A) \to (\oA, d_{\oA})$ be a morphism of 
connected $1$-finite $\k$-$\cdgas$ inducing an isomorphism 
$H^1(\varphi) \colon H^1(A) \isom H^1(\oA)$.
Then the map
\[
K^\bullet(\varphi) \coloneqq \varphi \otimes \Phi \colon 
K^\bullet(A) = (A \otimes_{\k} S,\, \delta^A) 
\longrightarrow 
K^\bullet(\oA) = (\oA \otimes_{\k} \oS,\, \delta^{\oA})
\]
is a morphism of filtered cochain complexes, and its degreewise dual
\[
K_\bullet(\varphi) \coloneqq \varphi^\vee \otimes \Phi^{-1} \colon 
K_\bullet(\oA) = (\oA_\bullet \otimes_{\k} \oS,\, \partial^{\oA}) 
\longrightarrow 
K_\bullet(A) = (A_\bullet \otimes_{\k} S,\, \partial^A) 
\]
is a morphism of filtered chain complexes.  Hence, after identifying $S$ 
and $\oS$ via $\Phi$, there are natural morphisms
\[
\B^\bullet(A) \longrightarrow \B^\bullet(\oA)
\qquad\text{and}\qquad
\B_\bullet(\oA) \longrightarrow \B_\bullet(A).
\]
\end{proposition}

\begin{proof}
Write $\delta^A = d_A\otimes\id_S + \omega_A\cdot(-)$ and similarly for
$\delta^{\oA}$.  We first show that $K^\bullet(\varphi) = \varphi \otimes \Phi$ 
commutes with the Koszul differentials.

Recall that $\omega_A \in H^1(A) \otimes_{\k} H_1(A)$ is the canonical 
element corresponding to the identity map $\id \colon H^1(A) \to H^1(A)$ 
under the identification $H^1(A) \otimes_{\k} H_1(A) \cong 
\End_{\k}(H^1(A))$.  Under the isomorphism $H^1(\varphi) \colon 
H^1(A) \isom H^1(\oA)$, the canonical element transforms as
\begin{equation}
\label{eq:omega-naturality}
\bigl(H^1(\varphi) \otimes H_1(\varphi)^{-1}\bigr)(\omega_A) = \omega_{\oA},
\end{equation}
since both sides correspond to the identity map on $H^1(\oA)$ under 
the identification $H^1(\oA) \otimes_{\k} H_1(\oA) \cong \End_{\k}(H^1(\oA))$.  
Consequently $\varphi\otimes\Phi$ carries multiplication by $\omega_A$ to
multiplication by $\omega_{\oA}$, so the two multiplication terms are
intertwined.

Since $\varphi$ is a $\cdga$ morphism, it also intertwines the 
differentials:
\[
(\varphi \otimes \Phi) \circ (d_A \otimes \id_S)
= (d_{\oA} \otimes \id_{\oS}) \circ (\varphi \otimes \Phi).
\]
Adding these two identities shows that $\varphi\otimes\Phi$ commutes with 
the full Koszul differentials $\delta^A$ and $\delta^{\oA}$.  Since 
$\Phi(\m)=\bar{\m}$, it also preserves the $\m$-adic filtrations.  
Dualizing degreewise---so that \eqref{eq:omega-naturality} becomes
$(\varphi^\vee\otimes\Phi^{-1})(\omega_{\oA}^\vee)=\omega_A^\vee$---gives the
assertion for $K_\bullet(\varphi)=\varphi^\vee\otimes\Phi^{-1}$.
\end{proof}

Thus the cochain Koszul modules are covariant in the $\cdga$, while the chain
Koszul modules are contravariant---as they must be, since
$A_\bullet=(A^\bullet)^\vee$.  Compare the holonomy Lie algebra $\h(-)$, which
is likewise contravariant in the $\cdga$ and covariant in the space.

\begin{corollary}
\label{cor:koszul-iso-criterion}
The map $K_\bullet(\varphi)$ is an isomorphism of filtered chain complexes, 
and hence induces isomorphisms $\B_i(A) \cong \B_i(\oA)$ for all $i$, 
in either of the following cases:
\begin{enumerate}[itemsep=2pt]
\item \label{av1}
$\varphi$ is a $\cdga$ isomorphism;
\item \label{av2}
$A$ and $\oA$ are generated in degree $1$ and 
$\varphi^1 \colon A^1 \to \oA^1$ is an isomorphism.
\end{enumerate}
\end{corollary}

\begin{proof}
In case~\eqref{av1}, $\varphi$ is invertible at every degree, so 
$\varphi \otimes \Phi$ is an isomorphism of chain complexes.  
In case~\eqref{av2}, since $A$ and $\oA$ are generated in degree~$1$, the degree-$2$ 
piece of $\varphi$ is the map induced by $\varphi^1$ on $A^1 \wedge A^1$; 
since $\varphi^1$ is an isomorphism, so is $\varphi^2 \colon A^2 \to \oA^2$.  
Together with $\varphi^1$ being an isomorphism, this ensures that 
$\varphi \otimes \Phi$ is an isomorphism at the chain level in degrees~$0$, 
$1$, and~$2$, hence an isomorphism of filtered chain complexes.
\end{proof}

\begin{example}
\label{ex:non-functorial-qiso}
Let $A = \bwedge\langle a \rangle$ with $da = 0$, and let 
$\oA = \CE(\sol_2)$ be the $\cdga$ of Example~\ref{ex:sol2-koszul}.
The inclusion $\iota \colon A \inj \oA$ is a quasi-isomorphism: 
both algebras have cohomology $\bwedge\langle a \rangle$, but $\iota^1 \colon 
A^1 \to \oA^1$ is not an isomorphism.  
Setting $S = \oS = \k[x]$, we have $\B_1(A) = 0$ while 
$\B_1(\oA) \cong S/(x-1)$, by \eqref{eq:sol2-modules-coh} 
and~\eqref{eq:sol2-modules}.  
Thus $\iota$ creates a nonzero Koszul module in degree $1$, supported at 
the point $1 \in H^1(\oA)$---which corresponds to the cohomology class 
$[a]$---and so away from the origin.  
This illustrates that Koszul homology is sensitive to the 
cochain-level data of $\iota$, not merely to the induced map on cohomology.
\end{example}

\begin{lemma}
\label{lem:koszul-surjective}
Let $\varphi \colon (A, d_A) \inj (\oA, d_{\oA})$
be a $\cdga$ morphism that is injective in degree $1$ and
induces an isomorphism $H^1(\varphi)\colon H^1(A)\isom H^1(\oA)$.
Then the induced map on first Koszul homology
\[
\B_1(\varphi)\colon \B_1(\oA) \longrightarrow \B_1(A)
\]
is surjective.
\end{lemma}

\begin{proof}
Since $H^1(\varphi)$ is an isomorphism, Proposition~\ref{prop:koszul-functorial}
provides a chain map $K_\bullet(\varphi)\colon \allowbreak 
K_\bullet(\oA)\to K_\bullet(A)$
over the common ring $S=\Sym(H_1(A))\cong\Sym(H_1(\oA))$.
As $\varphi^1\colon A^1\to\oA^1$ is injective, its transpose
$(\varphi^1)^\vee\colon\oA_1\to A_1$ is surjective---any injection of
$\k$-vector spaces splits---so
\[
K_1(\varphi)\colon K_1(\oA)=\oA_1\otimes_{\k} S
\longrightarrow A_1\otimes_{\k} S = K_1(A)
\]
is surjective; and $K_0(\varphi)$ is the identity of $S$, since $A$ and
$\oA$ are connected.
Now let $z\in K_1(A)$ be a cycle, and choose $\bar z\in K_1(\oA)$ with 
$K_1(\varphi)(\bar z)=z$.  Then $K_0(\varphi)(\partial_1^{\oA}\bar z)
=\partial_1^{A}z=0$, and $K_0(\varphi)$ is injective, so $\bar z$ is a 
cycle as well.  Its class in $\B_1(\oA)$ maps to the class of $z$, which 
proves surjectivity of $\B_1(\varphi)$.
\end{proof}

The comparison also propagates to the spectral sequences of
Section~\ref{sect:koszul-ss} and to the weight spectral sequence of
Theorem~\ref{thm:weight-ss}.

\begin{lemma}
\label{lem:koszul-ss-natural}
Let $\varphi\colon A\to\oA$ be a morphism of connected, $q$-finite
$\k$-$\cdgas$ inducing an isomorphism on $H^1$, and identify $S$ with $\oS$
via $\Phi$ as in Proposition~\ref{prop:koszul-functorial}.  Then:
\begin{enumerate}[itemsep=2pt]
\item \label{nat-coh}
The map $K^\bullet(\varphi)\colon K^\bullet(A)\to K^\bullet(\oA)$ is a
morphism of filtered cochain complexes for the $\m$-adic filtrations, and
hence induces a morphism between the cohomological Koszul spectral
sequences of Theorem~\ref{thm:koszul-ss-coh}.  Under the identifications
of\/ \eqref{ssc2}, this morphism is $K^\bullet(H^*(\varphi))$ on
$E_1$-pages and $\B^\bullet(H^*(\varphi))$ on $E_2$-pages; on
$E_\infty$-pages, in total degrees $\le q$, it is
$\gr^\m\B^\bullet(\varphi)$.
\item \label{nat-hom}
Dually, $K_\bullet(\varphi)\colon K_\bullet(\oA)\to K_\bullet(A)$ is a
morphism of filtered chain complexes for the $\m$-adic filtrations, and
induces a morphism between the homological Koszul spectral sequences of
Theorem~\ref{thm:koszul-ss-hom}, equal to $K_\bullet(H^*(\varphi))$ on
$E^1$-pages, to $\B_\bullet(H^*(\varphi))$ on $E^2$-pages, and to
$\gr^\m\B_\bullet(\varphi)$ on $E^\infty$-pages in total degrees $\le q$.
\item \label{nat-weight}
If $A$ and $\oA$ have positive weights and $\varphi$ is weight-preserving,
the same conclusions hold for the weight filtration
\eqref{eq:weight-filtration} and the spectral sequences of
Theorem~\ref{thm:weight-ss}.
\end{enumerate}
\end{lemma}

\begin{proof}
\eqref{nat-coh} That $K^\bullet(\varphi)=\varphi\otimes\Phi$ is a cochain
map preserving the $\m$-adic filtration is
Proposition~\ref{prop:koszul-functorial}, the essential point being the
naturality \eqref{eq:omega-naturality} of the canonical element.  A
filtered map induces a morphism of the associated spectral sequences, so
it remains to identify it page by page.

On $E_0$-pages, $\gr^\m_p K^\bullet(A)=\bigl(A\otimes_{\k}
\m^p/\m^{p+1},\,d_A\otimes\id\bigr)$ and the induced map is
$\varphi\otimes\gr_p\Phi$; passing to cohomology gives
$H^*(\varphi)\otimes\gr_p\Phi$ on $E_1$-pages.  By \eqref{ssc2} the
$E_1$-page is $K^\bullet(H^*(A))$, whose differential $d_1$ is
multiplication by $\omega_A$; since $H^*(\varphi)$ is an algebra map
carrying $\omega_A$ to $\omega_{\oA}$ by \eqref{eq:omega-naturality}, the
map $H^*(\varphi)\otimes\gr\Phi$ commutes with $d_1$ and is precisely
$K^\bullet(H^*(\varphi))$.  Taking $d_1$-cohomology gives
$\B^\bullet(H^*(\varphi))$ on $E_2$-pages.  Finally, the map of abutments
is $\B^\bullet(\varphi)$, which preserves the filtrations, and the induced
map on $E_\infty\cong\gr^\m\B^\bullet$ of\/ \eqref{ssc4} is $\gr^\m$ of
it; the restriction to total degrees $\le q$ is the range in which
\eqref{ssc4} identifies $E_\infty$.

\eqref{nat-hom} The transpose $\varphi^\vee\colon\oA_\bullet\to A_\bullet$
is a morphism of dg-coalgebras and $\Phi^{-1}(\bar{\m})=\m$, so
$K_\bullet(\varphi)$ preserves the filtrations
$F^pK_\bullet(-)=(-)_\bullet\otimes_{\k}\n^p$ of
Theorem~\ref{thm:koszul-ss-hom}.  The identification of the pages is
verbatim as above, with $H_\bullet(\varphi^\vee)=H^*(\varphi)^\vee$ in
place of $H^*(\varphi)$.

\eqref{nat-weight} The weight filtration is preserved because $\varphi$ is
weight-preserving, and the identifications of the $E^0$- and $E^1$-pages
in Theorem~\ref{thm:weight-ss} are natural in the same way; for the
$E^2$-page one uses in addition that $\varphi$ carries the
weight-homogeneous $H^1(A)$ isomorphically onto $H^1(\oA)$.
\end{proof}

When the hypothesis that $H^1(\varphi)$ be an isomorphism is dropped,
one can no longer compare the Koszul complexes over a common symmetric
algebra.  Nevertheless, functoriality survives in a canonical
base-change form, as we now explain.

\begin{proposition}
\label{prop:canonical-basechange}
With notation as above, let $\varphi\colon A\to \oA$ be a $\cdga$ morphism. 
There is a canonical $\oS$-linear homomorphism
\[
\beta(\varphi)\colon \B(A)\otimes_{S} \oS \longrightarrow \B(\oA),
\]
functorial in $\varphi$, without any hypothesis on $H_1(\varphi)$.
\end{proposition}

\begin{proof}
Let $S\to \oS$ be the algebra map induced by
the homomorphism $H_1(\varphi)\colon H_1(\oA)\to H_1(A)$.
The Koszul complex $K_{\bullet}(A)=(A^{\bullet}\otimes_{\k} S,\partial^A)$ 
whose homology produces $\B(A)$ is 
built functorially from the $S$-module $A\otimes_{\k} S$.  Applying
$\varphi\otimes\id$ yields a chain map
\[
(A\otimes_{\k} S)\otimes_{S} \oS \cong A\otimes_{\k} \oS
\longrightarrow \oA \otimes_{\k} \oS,
\]
where the left-hand side is the base-change of $A\otimes_{\k} S$ 
along $S\to \oS$.
Passing to homology (and noting the $\oS$-linearity) produces the desired
map $\beta(\varphi)$.  Functoriality is immediate from the functoriality
of the whole Koszul construction.
\end{proof}

%%%%%%%%%%%%%%%%%%%%
\subsection{Behavior under quasi-isomorphisms}
\label{subsec:cdga-qiso}
%%%%%%%%%%%%%%%%%%%%

We now examine how Koszul modules behave under quasi-isomorphisms.
Although Koszul homology is not invariant under arbitrary quasi-isomor\-phisms,
the presence of an isomorphism on $H^{1}$ allows one to compare the associated
Koszul complexes over a common polynomial ring. In this situation, the
$\m$-adic filtration gives rise to a natural spectral sequence that controls
the effect of a $q$-(quasi-)isomorphism on Koszul homology and on its 
associated graded modules.

\begin{theorem}
\label{thm:q-iso-koszul}
Let $\varphi \colon A \to \oA$ be a $\cdga$ morphism between connected 
$q$-finite $\k$-$\cdgas$ inducing an isomorphism on $H^1$.  With notation 
as above,
\begin{enumerate}[itemsep=2pt]
\item \label{iso1}
If $\varphi$ is a $q$-isomorphism and the $\m$-adic filtrations 
on $\B_i(A)$ and $\B_i(\oA)$ are separated for all $i\le q$, 
then $\B_i(\varphi)\colon \B_i(\oA) \to \B_i(A)$
is an $S$-linear isomorphism for all $i \le q$.
\item \label{iso1b}
If $\varphi$ is only a $q$-quasi-isomorphism, and the $\m$-adic filtration
on $\B_i(\oA)$ is separated, then $\B_i(\varphi)$ is injective and its
cokernel $C_i$ satisfies $\m C_i=C_i$; in particular, $(C_i)_\m=0$, so the
origin does not lie in $\Supp_S C_i$.  If, in addition, $\B_i(A)$ is a
graded $S$-module, or if $\Supp_S\B_i(A)\subseteq\{0\}$, then
$\B_i(\varphi)$ is an isomorphism.
\item \label{iso3}
If $\varphi$ is a $q$-quasi-isomorphism, then the induced map on
$\m$-adic completions,
\[
\widehat{\B_i(\varphi)} \colon \widehat{\B_i(\oA)}
\longrightarrow \widehat{\B_i(A)},
\]
is an $\widehat{S}$-linear isomorphism for all $i\le q$.
\item \label{iso2}
If $\varphi$ is a $q$-quasi-isomorphism, then the induced map on 
associated graded Koszul modules
\[
\gr(\B_i(\varphi)) \colon \gr^F \B_i(\oA) \longrightarrow \gr^F \B_i(A)
\]
is an $S$-linear isomorphism for all $i \le q$, where $F$ denotes on each
side the induced filtration of Theorem~\ref{thm:koszul-ss-hom}\eqref{ssh4}.
\end{enumerate}
\end{theorem}

\begin{proof}
Since $\varphi$ induces an isomorphism on $H^1$, it induces an
$S$-linear morphism of Koszul chain complexes
$K_\bullet(\varphi)\colon K_\bullet(\oA)\to K_\bullet(A)$ preserving the
$\m$-adic filtrations $F^pK_\bullet=A_\bullet\otimes_{\k}\m^p$, hence a
morphism of the spectral sequences of
Theorem~\ref{thm:koszul-ss-hom}.  On $E^1$-pages this is
$H_\bullet(\varphi)\otimes\id_S$.

\smallskip
\noindent\eqref{iso3}
As $\varphi$ is a $q$-quasi-isomorphism, $H_i(\varphi)$ is an
isomorphism for $i\le q$ and an epimorphism for $i=q+1$; hence
$E^1(\varphi)$ is an isomorphism in total degrees $\le q$ and an
epimorphism in total degree $q+1$.  An induction on pages, applying the
five lemma to $E^{r+1}=\ker d^r/\im d^r$, shows that $E^r(\varphi)$ has
the same property for every $r\ge1$, and therefore so does
$E^\infty(\varphi)$.

In each homological degree $i\le q$, the modules $K_i$ and $K_{i-1}$ are 
finitely generated free over the Noetherian ring $S$; by the Artin--Rees 
lemma, the filtration $F^p\B_i=\im\bigl(H_i(F^pK_\bullet)\to\B_i\bigr)$ 
induced on homology is $\m$-stable (the boundary terms coming from 
$K_{i+1}$, which need not be finitely generated, are handled by 
exhaustiveness of the filtration alone: any $x\in F^pK_i\cap\im\partial$ 
is a boundary of an element of $F^0K_{i+1}$).  Consequently the spectral 
sequence converges in degrees $\le q$, with 
$E^\infty\cong\gr^F\B_\bullet$; see \cite[pp.~22--24]{Serre}.  Being 
$\m$-stable, $F$ is equivalent to the $\m$-adic filtration, so the two 
have the same completion.  An isomorphism on associated graded objects of 
exhaustive filtrations of complete Hausdorff modules induces an 
isomorphism of the modules themselves; applying this to the 
$\m$-adic completions $\widehat{\B_i(\oA)}\to\widehat{\B_i(A)}$---which 
are complete and Hausdorff, being finitely generated over the complete 
local ring $\widehat{S}_\m$ in the relevant gradings---gives 
\eqref{iso3} in degrees $\le q$.

\smallskip
\noindent\eqref{iso2}
For a finitely generated module $M$ over a Noetherian ring one has
$\widehat{M}/\m^p\widehat{M}\cong M/\m^pM$, whence
$\gr^\m M\cong\gr^\m\widehat{M}$.  Now apply $\gr^\m$ to \eqref{iso3}.

\smallskip
\noindent\eqref{iso1} Since $\varphi$ is a $q$-isomorphism, the chain 
map $K_\bullet(\varphi)$ is an isomorphism in degrees $\le q$; hence 
$\B_i(\varphi)\colon\B_i(\oA)\to\B_i(A)$ is an isomorphism outright for 
$i\le q-1$, and for $i=q$ it identifies the cycle modules, 
$Z_q(\oA)\cong Z_q(A)$, so that $\B_q(\varphi)$ becomes the canonical 
surjection
\[
\B_q(\oA)=Z_q\big/\im\partial^{\oA}_{q+1}
\longsurj
Z_q\big/\im\partial^{A}_{q+1}=\B_q(A),
\]
the containment $\im\partial^{\oA}_{q+1}\subseteq\im\partial^{A}_{q+1}$ 
holding because $K_\bullet(\varphi)$ is a chain map.  It remains to prove 
injectivity at $i=q$.  The kernel $M=\im\partial^A_{q+1}/
\im\partial^{\oA}_{q+1}$ is a submodule of $\B_q(\oA)$, and by 
\eqref{iso2} the surjection above induces an isomorphism on associated 
graded modules; hence every element of $M$ has vanishing symbol in every 
filtration degree, that is, $M\subseteq\bigcap_{p} F^p\B_q(\oA)
\subseteq\bigcap_{p}\m^{p-c}\B_q(\oA)$ for some constant $c$, by the 
Artin--Rees lemma applied to $Z_q\subseteq K_q(\oA)$.  Separation of 
$\B_q(\oA)$ then gives $M=0$.  Note that only separation of the 
\emph{source}\/ $\B_i(\oA)$ is used.

\noindent\eqref{iso1b} Write $f=\B_i(\varphi)\colon M=\B_i(\oA)\to
N=\B_i(A)$, a filtration-preserving $S$-linear map, and recall from
\eqref{iso2} that $\gr(f)$ is an isomorphism.  If $f(m)=0$ with $m\ne0$,
then the symbol of $m$ vanishes in every filtration degree, so
$m\in\bigcap_{p}F^pM\subseteq\bigcap_{p}\m^{p-c}M$ for some constant $c$, by
Artin--Rees; the latter intersection is $0$ by the separation hypothesis,
so $f$ is injective.  Set $C_i=N/f(M)$.  Surjectivity of $\gr(f)$ gives
$N=f(M)+F^pN$ for every $p$, by successive approximation, whence
$C_i=F^pC_i\subseteq\m^{p-c'}C_i$ for some constant $c'$; therefore
$C_i=\bigcap_{p}\m^pC_i$, and $\m C_i=C_i$ by Artin--Rees.  As $C_i$ is
finitely generated, Nakayama's lemma over $S_\m$ gives $(C_i)_\m=0$, that
is, $0\notin\Supp_SC_i$.  Finally, $C_i$ is a quotient of $N$: if $N$ is
graded, then $C_i$ is a finitely generated graded module with
$\m C_i=C_i$, so $C_i=0$ by the graded Nakayama lemma; and if
$\Supp_SN\subseteq\{0\}$, then $\Supp_SC_i\subseteq\{0\}$ as well, which
together with $0\notin\Supp_SC_i$ forces $C_i=0$.
\end{proof}

\begin{remark}
\label{rem:completion-vs-separation}
Parts~\eqref{iso3} and~\eqref{iso2} require no separation hypothesis.
The reason is that completion does not distinguish between the
$\m$-adic filtration and the filtration induced by the spectral
sequence, whereas the passage from $\gr$ back to the module itself does;
Example~\ref{ex:circle-models} shows that \eqref{iso1} genuinely fails
without separation.  Part~\eqref{iso1b} isolates what separation alone
buys, namely injectivity with a cokernel supported away from the origin;
to recover surjectivity one needs a hypothesis, such as gradedness, that
rules out such a cokernel.
\end{remark}

\begin{corollary}
\label{cor:koszul-completion-formal}
Let $(A,d_A)$ be a connected, $1$-finite, finite-type $\k$-$\cdga$ which
is $q$-formal.  Then, for all $i\le q$, there are natural isomorphisms
\[
\widehat{\B_i(A)} \cong \B_i\bigl(H^*(A)\bigr)\otimes_{S}\widehat{S},
\qquad
\gr^F \B_i(A) \cong \B_i\bigl(H^*(A)\bigr),
\]
the first of $\widehat{S}$-modules and the second of graded $S$-modules,
with $F$ the induced filtration of
Theorem~\ref{thm:koszul-ss-hom}\eqref{ssh4}.
\end{corollary}

\begin{proof}
Apply Theorem~\ref{thm:q-iso-koszul}\eqref{iso3} and \eqref{iso2} to
each morphism in a zig-zag of $q$-quasi-isomorphisms joining $A$ to
$(H^*(A),0)$, and note that $\B_i(H^*(A))$ is a graded $S$-module, hence
coincides with its own $\m$-adic associated graded.
\end{proof}

In particular
$\Hilb\bigl(\gr^F\B_i(A),t\bigr)=\Hilb\bigl(\B_i(H^*(A)),t\bigr)$.  For $i=1$,
Lemma~\ref{lem:filtration-shift-B1} converts this into the $\m$-adic statement
$\Hilb\bigl(\gr^\m\B_1(A),t\bigr)=\Hilb\bigl(\gr^\m\B_1(H^*(A)),t\bigr)$.

The gradedness hypothesis in the next corollary is met, for instance,
whenever each $\cdga$ in the zig-zag admits positive weights
\textup{(}Proposition~\ref{prop:koszul-separation}\textup{)}.

\begin{corollary}
\label{cor:koszul-formal-separated}
Let $(A,d_A)$ be as above and $q$-formal, and suppose that every $\cdga$
in some zig-zag realizing $q$-formality has \emph{graded}\/ Koszul modules
in degrees $\le q$.  Then
$\B_i(A) \cong \B_i\bigl(H^*(A)\bigr)$ as $S$-modules, for all $i\le q$.
\end{corollary}

\begin{proof}
Each arrow of the zig-zag is a $q$-quasi-isomorphism between $\cdgas$ whose
Koszul modules are graded, hence separated
\textup{(}Proposition~\ref{prop:separation-criterion}\textup{)}; so
Theorem~\ref{thm:q-iso-koszul}\eqref{iso1b} applies to it, in whichever
direction it points, and shows that the induced map on $\B_i$ is an
isomorphism for $i\le q$.  Composing these isomorphisms along the zig-zag,
and using that $\B_i(H^*(A))$ is graded, gives the assertion.
\end{proof}

\begin{remark}
\label{rem:separated-not-enough-formal}
Gradedness cannot be relaxed to separation here.  As noted in
Remark~\ref{rem:completion-vs-separation}, separation of the source buys
only injectivity, with a cokernel that is invisible at the origin but need
not vanish: over $S=\k[x]$, the map
$S\xrightarrow{\,x-1\,}S$ is injective, induces an isomorphism on
associated graded modules---the initial form of $x-1$ being a unit---and
has separated source and target, yet its cokernel is $\k$.  Nor can
non-isomorphisms be composed along a zig-zag, whose arrows point in both
directions.  What survives without any such hypothesis are the associated
graded and completed statements of
Corollary~\ref{cor:koszul-completion-formal}; the framing theory developed
below may be run at either of those levels verbatim, replacing $\B_i$
throughout by $\gr^\m\B_i$ or by $\widehat{\B_i}$.
\end{remark}

Corollary~\ref{cor:koszul-formal-separated} produces isomorphisms
$\B_i(A)\cong\B_i(H^*(A))$ for $i\le q$ from a zig-zag of
$q$-quasi-isomorphisms realizing $q$-formality, but these depend on the
zig-zag.  The next proposition identifies the resulting ambiguity: the set
of available isomorphisms is a single orbit under a group that depends
only on the cohomology algebra.

\begin{definition}
\label{def:koszul-monodromy}
Let $H$ be a connected graded $\k$-algebra with $\dim_\k H^n<\infty$ for
$n\le q$, regarded as a $\cdga$ with zero differential.  For $i\le q$, let
$G_i(H)\subseteq\Aut_\k\bigl(\B_i(H)\bigr)$ be the set of automorphisms
obtained by composing the maps $\B_i(-)$ and their inverses along a zig-zag
of $q$-quasi-isomorphisms from $H$ to $H$, all of whose terms are
$q$-finite with separated Koszul modules.  We call $G_i(H)$ the {\em Koszul
monodromy group}\/ of $H$ in degree $i$.
\end{definition}

\begin{proposition}
\label{prop:formality-iso-orbit}
Let $A$ be a connected, $q$-finite $\cdga$ which is $q$-formal via
zig-zags to which Corollary~\ref{cor:koszul-formal-separated} applies, put
$H=H^*(A)$, and for $i\le q$ let
\[
\Iso_i(A)\subseteq \Isom\bigl(\B_i(A),\B_i(H)\bigr)
\]
be the set of isomorphisms arising from such zig-zags.  Then:
\begin{enumerate}[itemsep=2pt]
\item \label{orb1}
$G_i(H)$ is a group, and each of its elements is semilinear over the ring
automorphism $\Sym(\tilde g|_{H^1})$ of $S$, where $\tilde g$ is the
algebra automorphism induced on $H^{\le q}$; the assignment
$g\mapsto\tilde g$ is a group homomorphism $G_i(H)\to\Aut(H^{\le q})$.
\item \label{orb2}
$\Iso_i(A)$ is a single orbit: $\Iso_i(A)=G_i(H)\cdot u$ for every
$u\in\Iso_i(A)$.  In particular, the ambiguity in the choice of
identification depends only on $H$, and not on $A$ or on any zig-zag.
\item \label{orb3}
$\{\B_i(\alpha) \mid \alpha\in\Aut(H)\}$ is a subgroup of $G_i(H)$.
\item \label{orb4}
If the zig-zags in Definition~\ref{def:koszul-monodromy} may be taken with
all terms of positive weights, all arrows weight-preserving, and all
weight spectral sequences collapsing at $E^2$ in total degrees $\le q$,
then equality holds in \eqref{orb3}.
\end{enumerate}
\end{proposition}

\begin{proof}
\eqref{orb1} Zig-zags from $H$ to $H$ are closed under concatenation and
reversal and $\B_i(-)$ is functorial, so $G_i(H)$ is a group.  Each arrow
$\psi$ contributes a map that is semilinear over $\Sym(H^1(\psi))$ by
Proposition~\ref{prop:koszul-functorial}; composing, $g$ is semilinear
over $\Sym(\tilde g|_{H^1})$, and $g\mapsto\tilde g$ is multiplicative
because $H^*(-)$ is a functor.

\eqref{orb2} If $u,u'\in\Iso_i(A)$ arise from zig-zags $Z,Z'$ from $A$ to
$H$, then $u'\circ u^{-1}$ is the automorphism of $\B_i(H)$ induced by the
concatenation of $Z'$ with the reverse of $Z$, a zig-zag from $H$ to $H$ of
the required type; hence $u'\in G_i(H)\cdot u$.  Conversely, $g\cdot u$
arises from the concatenation of a zig-zag realizing $g$ with $Z$.

\eqref{orb3} An algebra automorphism $\alpha$ of $H$ is in particular a
$q$-quasi-isomorphism $H\to H$, and $K_i(\alpha)=\B_i(\alpha)$;
functoriality gives a group homomorphism.

\eqref{orb4} Under the stated hypotheses each term $C$ of a zig-zag
carries the canonical isomorphism
$w(C)\colon\B_i(C)\to\B_i(H^*(C))$ given by the edge maps of its weight
spectral sequence, and by
Lemma~\ref{lem:koszul-ss-natural}\eqref{nat-weight} every arrow
$\psi\colon C\to D$ satisfies $w(C)\circ K_i(\psi)=\B_i(H^*(\psi))\circ
w(D)$.  Composing along the zig-zag, and using $w(H)=\id$ by
Corollary~\ref{cor:ss-collapse-formal}, gives $g=\B_i(\tilde g)$.
\end{proof}

\begin{remark}
\label{rem:orbit-vs-aut}
We do not know whether the inclusion in
Proposition~\ref{prop:formality-iso-orbit}\eqref{orb3} is always an
equality, that is, whether every element of the Koszul monodromy group is
induced by an algebra automorphism of $H$.  The argument for \eqref{orb4}
requires a chosen identification $\B_i(C)\cong\B_i(H^*(C))$ at each
intermediate term $C$ of the zig-zag, and outside the canonically framed
situation such identifications are themselves subject to the ambiguity
under discussion.  Since $\B_i(-)$ is not a quasi-isomorphism invariant
\textup{(}Example~\ref{ex:non-functorial-qiso}\textup{)}, a self-zig-zag of
$H$ inducing the identity on $H^{\le q}$ need not obviously induce the
identity on $\B_i(H)$.  A counterexample would exhibit a $q$-formal
$\cdga$ whose formality identifications are strictly more ambiguous than
its cohomology automorphisms account for.
\end{remark}

\begin{remark}
\label{rem:pw-not-enough}
Positive weights alone do not suffice in
Corollary~\ref{cor:koszul-formal-separated}: what they supply is
separation, not a comparison morphism.  The Chevalley--Eilenberg
$\cdga$ $A=\CE(\h(1))$ has positive weights, yet $\B_1(A)\cong\k$ while
$\B_1(H^*(A))\cong S$ \textup{(}Section~\ref{subsec:ka=kha}\textup{)}.
Formality is what produces the filtered morphism to which the
comparison theorem for filtered complexes may be applied; no amount of
control on the filtrations can manufacture one.
\end{remark}

\begin{theorem}
\label{thm:q-iso-koszul-coh}
Let $\varphi \colon A \to \oA$ be a $\cdga$ morphism between connected 
$1$-finite $\k$-$\cdgas$ inducing an isomorphism on $H^1$.
\begin{enumerate}[itemsep=2pt]
\item \label{ciso1}
If $\varphi$ is a $q$-isomorphism and the $\m$-adic filtrations on 
$\B^i(A)$ and $\B^i(\oA)$ are separated for all $i\le q$, then
$\B^i(\varphi)\colon \B^i(A) \to \B^i(\oA)$
is an $S$-linear isomorphism for all $i \le q$.
\item \label{ciso2}
If $\varphi$ is a $q$-quasi-isomorphism, then the induced map 
on associated graded cohomological Koszul modules
\[
\gr(\B^i(\varphi)) \colon \gr^F \B^i(A) \longrightarrow \gr^F \B^i(\oA)
\]
is an $S$-linear isomorphism for all $i \le q$, where $F$ is the induced
filtration of Theorem~\ref{thm:koszul-ss-coh}\eqref{ssc4}.
\item \label{ciso1b}
If $\varphi$ is only a $q$-quasi-isomorphism, and the $\m$-adic filtration
on $\B^i(A)$ is separated, then $\B^i(\varphi)$ is injective with cokernel
$C^i$ satisfying $\m C^i=C^i$; and $\B^i(\varphi)$ is an isomorphism if in
addition $\B^i(\oA)$ is a graded $S$-module or has support contained in
$\{0\}$.
\end{enumerate}
\end{theorem}

\begin{proof}
The argument is identical to that of Theorem~\ref{thm:q-iso-koszul}, 
replacing the homological Koszul spectral sequence 
(Theorem~\ref{thm:koszul-ss-hom}) throughout with its cohomological 
counterpart (Theorem~\ref{thm:koszul-ss-coh}), whose $E_1$-page is 
$\B^{p+q}(H^*(A))_p$ and which converges to 
$\gr^\m_p\B^{p+q}(A)$.
For~\eqref{ciso1}, the separation hypothesis cannot be dispensed with: 
by Proposition~\ref{prop:separation-criterion} it is a genuine 
condition on the cohomological Koszul modules, and 
Example~\ref{ex:coho-nonseparated} exhibits a $\cdga$ for which it 
fails.  Part~\eqref{ciso1b} is proved exactly as
Theorem~\ref{thm:q-iso-koszul}\eqref{iso1b}, with \eqref{ciso2} in place
of \eqref{iso2}.  It holds, for instance, when $A$ and $\oA$ have positive 
weights, by the weight argument of 
Proposition~\ref{prop:koszul-separation}, which applies verbatim to 
$\B^\bullet$.
\end{proof}

Geometric consequences of Theorem~\ref{thm:q-iso-koszul}\eqref{iso3} for 
resonance varieties and tangent cones are developed in 
Section~\ref{sect:tcone}, culminating in Theorem~\ref{thm:formal-tc}: 
after completion, the Koszul modules become genuine invariants of 
$q$-quasi-isomorphism, and both the jump and support loci can be 
compared through their germs at the origin.

\begin{corollary}
\label{cor:q-iso-koszul-pw}
Let $\varphi \colon A \to \oA$ be a $q$-quasi-isomorphism between 
connected, $q$-finite $\k$-$\cdgas$ inducing an isomorphism on $H^1$, and 
suppose that $A$ and $\oA$ both admit positive weight decompositions.  
Then
\[
\B_i(\varphi) \colon \B_i(\oA) \longrightarrow \B_i(A)
\quad\text{and}\quad
\B^i(\varphi) \colon \B^i(A) \longrightarrow \B^i(\oA)
\]
are $S$-linear isomorphisms for all $i \le q$.
\end{corollary}

\begin{proof}
By Proposition~\ref{prop:koszul-separation}, positive weights make the 
Koszul modules $\B_i(A)$ and $\B_i(\oA)$ graded, hence separated; the same 
weight argument applies to $\B^i(A)$ and $\B^i(\oA)$.  Being graded on both 
sides, the two conclusions now follow from 
Theorem~\ref{thm:q-iso-koszul}\eqref{iso1b} and 
Theorem~\ref{thm:q-iso-koszul-coh}\eqref{ciso1b}, respectively.
\end{proof}

\begin{remark}
\label{rem:no-coho-analogue}
The example of Remark~\ref{rem:ce-sol2} serves on the cohomological side as
well.  For the inclusion $\varphi\colon\bigwedge(a)\hookrightarrow\CE(\sol_2)$
of two connected, finite-type $\cdgas$, a quasi-isomorphism inducing the
identity on $H^1$, one has $\B^2(\bigwedge(a))=0$ while
$\B^2(\CE(\sol_2))\cong\k[x]/(x-1)\ne0$
\textup{(}Example~\ref{ex:coho-nonseparated}\textup{)}, so $K^2(\varphi)$ is
not an isomorphism.  Again separation is what fails, and $\CE(\sol_2)$ admits
no positive weights.
\end{remark}

%%%%%%%%%%%%%%%%%%%%%%%%%%%%%%%%%%%%
\subsection{Naturality of multiplicative structures}
\label{subsec:naturality-mult}
%%%%%%%%%%%%%%%%%%%%%%%%%%%%%%%%%%%%

The $S$-module and comodule naturality established in 
Sections~\ref{subsec:functorial-iso} and~\ref{subsec:cdga-qiso} 
does not capture the full algebraic structure carried by the 
Koszul modules. In this subsection we show that the 
graded-commutative algebra structure on $\B^\bullet(A)$, the 
graded-cocommutative coalgebra structure on $\B_\bullet(A)$, 
and the $H^*(A)$-module structures of 
Proposition~\ref{prop:hstar-action} are all natural with respect 
to $\cdga$ morphisms inducing isomorphisms on $H^1$. Under 
$q$-quasi-isomorphisms, these multiplicative structures are 
preserved at the level of associated graded modules, and under 
the appropriate separation hypotheses they are preserved at 
the module level itself.

\begin{theorem}
\label{thm:naturality-alg-structures}
Let $\varphi\colon A\to\bar{A}$ be a $\cdga$ morphism inducing 
an isomorphism on $H^1$. Then:
\begin{enumerate}[itemsep=2pt]
\item \label{nat-alg1}
The induced map $\B^i(\varphi)\colon \B^i(A)\to\B^i(\bar{A})$ 
is a morphism of graded $S$-algebras, compatible with the 
graded-commutative algebra structures on $\B^\bullet(A)$ 
and $\B^\bullet(\bar{A})$.
\item \label{nat-alg2}
The induced map $\B_i(\varphi)\colon \B_i(\bar{A})\otimes_{\bar{S}} S
\to\B_i(A)$ is a morphism of graded $S$-coalgebras, compatible 
with the graded-cocommutative coalgebra structures on $\B_\bullet(A)$ 
and $\B_\bullet(\bar{A})$.
\item \label{nat-alg3}
The maps $K^\bullet(\varphi)$ and $K_\bullet(\varphi)$ are 
equivariant with respect to the $H^*(A)$-module structures 
of Proposition~\ref{prop:hstar-action}: for all 
$[\alpha]\in H^p(A)$ and $x\in\B^i(A)$,
\[
\B^{i+p}(\varphi)\bigl([\alpha]\cup x\bigr)
= H^p(\varphi)([\alpha])\cup \B^i(\varphi)(x),
\]
and for all $[\alpha]\in H^p(A)$ and $y\in\B_i(\bar{A})$,
\[
\B_{i-p}(\varphi)\bigl(H^p(\varphi)([\alpha])\cap y\bigr)
= [\alpha]\cap \B_i(\varphi)(y),
\]
where the cap products on the left and right are taken with 
respect to $H^*(\bar{A})$ and $H^*(A)$ respectively.
\end{enumerate}
\end{theorem}

\begin{proof}
\eqref{nat-alg1} The map $K^\bullet(\varphi)\colon K^\bullet(A)\to 
K^\bullet(\bar{A})$ of Proposition~\ref{prop:koszul-functorial} is a morphism 
of dg-algebras over $S$: it is multiplicative because $\varphi$ is, and it 
intertwines the two Koszul differentials by \eqref{eq:omega-naturality}.  
Passing to cohomology preserves the algebra structure.

\eqref{nat-alg2} Dually, the chain map on $K_\bullet$ is a 
morphism of dg-coalgebras over $S$, and passing to homology 
preserves the coalgebra structure.

\eqref{nat-alg3} For the cup product: given a cocycle 
$\alpha\in Z^p(A)$, the diagram
\[
\begin{tikzcd}[column sep=30pt, row sep=22pt]
K^\bullet(A) \ar[r, "\varphi\otimes\id_S"] 
             \ar[d, "\alpha\wedge(-)"'] &
K^\bullet(\bar{A}) \ar[d, "\varphi(\alpha)\wedge(-)"] \\
K^{\bullet+p}(A) \ar[r, "\varphi\otimes\id_S"'] &
K^{\bullet+p}(\bar{A})
\end{tikzcd}
\]
commutes at the chain level since $\varphi$ is a $\cdga$ map. 
Passing to cohomology gives the cup product equivariance.
For the cap product: given a cocycle $\alpha\in Z^p(A)$, 
the diagram
\[
\begin{tikzcd}[column sep=30pt, row sep=22pt]
K_\bullet(\bar{A})\otimes_{\bar{S}}S 
  \ar[r, "K_\bullet(\varphi)"] 
  \ar[d, "\varphi(\alpha)\lrcorner(-)"'] &
K_\bullet(A) \ar[d, "\alpha\lrcorner(-)"] \\
K_{\bullet-p}(\bar{A})\otimes_{\bar{S}}S 
  \ar[r, "K_\bullet(\varphi)"'] &
K_{\bullet-p}(A)
\end{tikzcd}
\]
commutes at the chain level by the same reasoning. Passing to 
homology gives the cap product equivariance.
\end{proof}

\begin{corollary}
\label{cor:qiso-alg-structures}
Let $\varphi\colon A\to\bar{A}$ be a $q$-quasi-isomorphism of 
connected, finite-type $\k$-$\cdgas$. Then for all $i\le q$:
\begin{enumerate}[itemsep=2pt]
\item \label{qiso-alg1}
The isomorphism $\gr^\m\B^i(A)\cong\gr^\m\B^i(\bar{A})$ of 
Theorem~\ref{thm:q-iso-koszul-coh}\eqref{ciso2} is an isomorphism 
of graded $S$-algebras.
\item \label{qiso-alg2}
The isomorphism $\gr^\m\B_i(A)\cong\gr^\m\B_i(\bar{A})$ of 
Theorem~\ref{thm:q-iso-koszul}\eqref{iso2} is an isomorphism 
of graded $S$-coalgebras.
\item \label{qiso-alg3}
Both isomorphisms are equivariant with respect to the 
$H^*(A)$-module structures, via the isomorphism 
$H^i(A)\cong H^i(\bar{A})$ induced by $\varphi$ for $i\le q$.
\end{enumerate}
If moreover $\varphi$ is a $q$-isomorphism and the $\m$-adic 
filtrations on $\B^i(A)$ and $\B_i(A)$ are separated 
(which holds for both when $A$ has positive weights, by 
Proposition~\ref{prop:koszul-separation}, and which by 
Proposition~\ref{prop:separation-criterion} is a genuine condition in 
either case), 
then the same conclusions hold for $\B^i(A)\cong\B^i(\bar{A})$ 
and $\B_i(A)\cong\B_i(\bar{A})$ themselves.
\end{corollary}

\begin{proof}
All three parts follow from 
Theorem~\ref{thm:naturality-alg-structures} applied to 
$\varphi$, together with the isomorphisms on associated graded 
modules provided by Theorems~\ref{thm:q-iso-koszul}\eqref{iso2} 
and~\ref{thm:q-iso-koszul-coh}\eqref{ciso2}. The passage from 
associated graded to the modules themselves under the separation 
hypothesis follows from the standard filtered modules argument 
used in the proofs of those theorems.
\end{proof}

%%%%%%%%%%%%%%%%%%%%%%%%%%%
\subsection{Koszul modules and $\B_q$-formality of morphisms}
\label{subsec:koszul-formality-obstruction}
%%%%%%%%%%%%%%%%%%%%%%%%%%%

We now explain how Koszul modules detect failures of formality at the level
of morphisms, and introduce a notion of formality for maps that is naturally
measured by Koszul homology.

\begin{definition}
\label{def:bq-framing}
Let $A$ be a connected, $q$-finite $\cdga$ which is $q$-formal, via 
zig-zags to which Corollary~\ref{cor:koszul-formal-separated} applies.  A 
{\em $\B_q$-framing}\/ of $A$ is a choice $u=(u_i)_{i\le q}$ of 
isomorphisms $u_i\in\Iso_i(A)$, in the notation of 
Proposition~\ref{prop:formality-iso-orbit}.  A {\em framed $\cdga$}\/ is a 
pair $(A,u)$.
\end{definition}

\begin{definition}
\label{def:bq-formal-framed}
Let $\varphi\colon A\to\oA$ be a $\cdga$ morphism inducing an isomorphism 
on $H^1$, between framed $\cdgas$ $(A,u)$ and $(\oA,\bar u)$.  We say that 
$\varphi$ is {\em $\B_q$-formal}\/ if for each $i\le q$ the square
\begin{equation}
\label{eq:bi-compatible}
\begin{tikzcd}[column sep=40pt, row sep=22pt]
\B_i(\oA) \ar[r, "\B_i(\varphi)"] \ar[d, "\bar u_i"'] &
\B_i(A) \ar[d, "u_i"] \\
\B_i\bigl(H^*(\oA)\bigr) \ar[r, "\B_i(H^*(\varphi))"] &
\B_i\bigl(H^*(A)\bigr)
\end{tikzcd}
\end{equation}
commutes, and that $\varphi$ is {\em multiplicatively $\B_q$-formal}\/ if 
in addition the vertical isomorphisms are compatible with the algebra and 
coalgebra structures of Proposition~\ref{prop:koszul-alg-coalg} and with 
the $H^*$-actions of Proposition~\ref{prop:hstar-action}.  When we say 
that $\varphi$ is $\B_q$-formal without reference to framings, we mean 
that some choice of framings makes it so; by 
Proposition~\ref{prop:framed-bq}\eqref{fbq2} this is a condition on 
$\varphi$ alone.
\end{definition}

\begin{proposition}
\label{prop:framed-bq}
With notation as above:
\begin{enumerate}[itemsep=2pt]
\item \label{fbq1}
Framed $\cdgas$ and $\B_q$-formal morphisms form a category: identities 
are $\B_q$-formal, and if $\varphi\colon(A,u)\to(\oA,\bar u)$ and 
$\psi\colon(\oA,\bar u)\to(\tilde{A},\tilde u)$ are $\B_q$-formal, so is 
$\psi\circ\varphi$.
\item \label{fbq2}
Fix framings $u$ of $A$ and $\bar u$ of $\oA$.  Then $\varphi$ is 
$\B_q$-formal, for some choice of framings, if and only if
\[
\B_i(H^*(\varphi))\in
G_i\bigl(H^*(A)\bigr)\cdot
\bigl(u_i\circ\B_i(\varphi)\circ\bar u_i^{-1}\bigr)\cdot
G_i\bigl(H^*(\oA)\bigr)
\qquad\text{for all } i\le q ,
\]
a condition independent of the framings fixed at the outset.
\item \label{fbq3}
If $A$ has positive weights and its weight spectral sequence collapses at 
$E^2$ in total degrees $\le q$, the edge isomorphisms constitute a 
canonical $\B_q$-framing $w(A)$.
\end{enumerate}
\end{proposition}

\begin{proof}
\eqref{fbq1} The squares for $\varphi$ and for $\psi$ share the framing 
$\bar u$ along their common edge, so they may be juxtaposed; the identity 
square commutes trivially.

\eqref{fbq2} By Proposition~\ref{prop:formality-iso-orbit}\eqref{orb2}, 
the admissible replacements for $u_i$ and $\bar u_i$ are $gu_i$ and 
$\bar h\bar u_i$ with $g\in G_i(H^*(A))$ and $\bar h\in G_i(H^*(\oA))$; 
these replace $u_i\circ\B_i(\varphi)\circ\bar u_i^{-1}$ by 
$g\bigl(u_i\circ\B_i(\varphi)\circ\bar u_i^{-1}\bigr)\bar h^{-1}$, whence 
the double-coset criterion, which is manifestly independent of the initial 
choices.

\eqref{fbq3} Under the stated hypotheses the weight spectral sequence has 
$E^2\cong E^\infty$ in total degrees $\le q$, while positive weights make 
$\B_i(A)$ graded, so $\gr^W\B_i(A)=\B_i(A)$; the edge maps are therefore 
isomorphisms $\B_i(A)\to\B_i(H^*(A))$, and they lie in $\Iso_i(A)$ by 
Corollary~\ref{cor:koszul-formal-separated}.
\end{proof}

By definition, $q$-formality implies $\B_q$-formality, but the converse 
need not hold in general: $\B_q$-formality is the Koszul-homological 
shadow of $q$-formality, and may be strictly weaker.

\begin{theorem}
\label{thm:koszul-obstructs-formality}
Let $\varphi \colon (A,d_A) \to (\oA,d_{\oA})$ be a $\cdga$ morphism 
between connected $1$-finite $\cdgas$ that are both $q$-formal, and 
that induces an isomorphism on $H^1$.
If $\varphi$ is $q$-formal, then $\varphi$ is $\B_q$-formal.
\end{theorem}

\begin{proof}
By definition, $\varphi$ is $q$-formal if it fits into a commutative 
diagram of the form \eqref{eq:ziggy-zagg}, in which every horizontal map 
is a $q$-quasi-isomorphism.  Restricting that diagram to its source and 
target columns yields zig-zags realizing the $q$-formality of $A$ and of 
$\oA$; let $u$ and $\bar u$ be the resulting $\B_q$-framings 
\textup{(}Definition~\ref{def:bq-framing}\textup{)}.
By Proposition~\ref{prop:koszul-functorial}, each arrow of 
\eqref{eq:ziggy-zagg} induces a contravariant morphism between the 
corresponding Koszul chain complexes over a common symmetric algebra, 
which by Theorem~\ref{thm:q-iso-koszul}\eqref{iso1b} is an isomorphism on 
$\B_i$ for $i\le q$---the Koszul modules of the $\cdgas$ appearing in
\eqref{eq:ziggy-zagg} being graded, by the hypothesis carried by
Definition~\ref{def:bq-framing}.  Since \eqref{eq:ziggy-zagg} commutes, composing 
these isomorphisms along the two columns yields the commutative square 
\eqref{eq:bi-compatible} for the framings $u$ and $\bar u$; compatibility 
with the multiplicative structures is 
Theorem~\ref{thm:naturality-alg-structures}\eqref{nat-alg3}.
This proves the assertion.
\end{proof}

\begin{corollary}
\label{cor:koszul-nonformal-mor}
Let $\varphi \colon A \to \oA$ be a $\cdga$ morphism between formal 
$\cdgas$, inducing an isomorphism on $H^{*}$.  If $\varphi$ fails to 
be $\B_q$-formal for some $q$, then $\varphi$ is not formal.
\end{corollary}

\begin{proof}
This is the contrapositive of
Theorem~\ref{thm:koszul-obstructs-formality}: a formal $\cdga$ is
$q$-formal for every $q$, and an isomorphism on $H^{*}$ is in particular
an isomorphism on $H^1$, so the hypotheses of that theorem hold for every
$q$; were $\varphi$ formal, it would be $q$-formal, hence
$\B_q$-formal, for all $q$.
\end{proof}

\begin{theorem}
\label{thm:CE-auto-formal}
Let $\g$ be a finite-dimensional $2$-step nilpotent Lie algebra over 
$\k$, set $A = \CE(\g)$, and let $\alpha \in \Aut(\g)$ with 
$\varphi = \CE(\alpha) \in \Aut(A)$.
\begin{enumerate}[itemsep=2pt]
\item \label{ff-auto} $\varphi$ is filtered-formal.
\item \label{B1-compat} If $A$ is $1$-formal, then $\varphi$ is 
$\B_1$-formal.
\end{enumerate}
\end{theorem}

\begin{proof}
For~\eqref{ff-auto}: since $\g$ is $2$-step nilpotent, we have 
$\gr \g = \g$ as graded Lie algebras, with $\g_1 = \g/Z(\g)$ and 
$\g_2 = Z(\g)$, so $A = \CE(\g)$ is already filtered-formal.
Any $\alpha \in \Aut(\g)$ preserves the lower central series 
filtration, hence $\varphi = \CE(\alpha)$ is compatible with the 
filtered structure and induces $\CE(\gr\alpha) = \varphi$ on the 
associated graded. Therefore $\varphi$ is filtered-formal.

For~\eqref{B1-compat}: since $A$ is $1$-formal, the $1$-formality 
zigzag provides $S$-module isomorphisms $\B_1(A) \cong \B_1(H^*(A))$.
For any $\alpha \in \Aut(\g)$, the action of $\alpha$ on $Z(\g)$ is 
determined by its action on $\g/Z(\g)$ via the bracket map 
$\bigwedge^2(\g/Z(\g)) \to Z(\g)$.
Since this bracket map is exactly the cup product on $H^*(A)$ (up to 
duality), the action of $\varphi^* = H^*(\varphi)$ on $H^*(A)$ sees 
the same data, and hence the square \eqref{eq:bi-compatible} commutes 
for $i = 1$ under the isomorphisms provided by the $1$-formality zigzag.
Thus $\varphi$ is $\B_1$-formal.
\end{proof}

Part~\eqref{B1-compat} shows that the obstruction of 
Theorem~\ref{thm:koszul-obstructs-formality} cannot detect 
non-$1$-formality of CE automorphisms in the $2$-step nilpotent setting:
every such automorphism is automatically $\B_1$-formal when $A$ is 
$1$-formal. This raises the question of whether $\B_q$-formality is 
strictly weaker than $q$-formality for morphisms of CE type.

\begin{question}
\label{q:CE-auto-1-formal}
Let $\g$ be a finite-dimensional $2$-step nilpotent Lie algebra with 
$A = \CE(\g)$ being $1$-formal, and let $\alpha \in \Aut(\g)$.  
Is $\varphi = \CE(\alpha)$ necessarily $1$-formal?
More generally, does $\B_q$-formality imply $q$-formality for 
morphisms of CE type, or does a gap between the two notions appear 
already for $3$-step nilpotent Lie algebras?
\end{question}

The following example shows that $\B_1$-formality is a 
non-trivial condition on morphisms: there exist CE automorphisms 
of $3$-step nilpotent algebras that fail to be $\B_1$-formal, 
in contrast to the automatic $\B_1$-formality established in 
Theorem~\ref{thm:CE-auto-formal}\eqref{B1-compat} for the 
$2$-step nilpotent case.

\begin{example}
\label{ex:ia-auto-nonformal}
Let $\g = \ff_2/\gamma_4\ff_2$ be the free $3$-step nilpotent Lie 
algebra on two generators $x_1, x_2$, with basis
\[
x_1, x_2, \quad y_{12} = [x_1,x_2], \quad 
y_{112} = [x_1,[x_1,x_2]], \quad y_{212} = [x_2,[x_1,x_2]],
\]
and set $A = \CE(\g)$, $S = \k[x_1, x_2]$.
Note that $A$ is not $1$-formal: the cup product 
$\bigwedge^2 H^1(A) \to H^2(A)$ vanishes, and $H^2(A)$ 
contains classes arising from triple Massey products.

Define $\alpha \in \Aut(\g)$ by
\[
\alpha(x_1) = x_1, \quad \alpha(x_2) = x_2 + y_{12},
\]
with bracket-forced actions
\[
\alpha(y_{12}) = y_{12} + y_{112}, \quad
\alpha(y_{112}) = y_{112}, \quad
\alpha(y_{212}) = y_{212}.
\]
Then $\alpha$ is an IA-automorphism: it acts trivially on 
$\g/\gamma_2\g$, so $\varphi = \CE(\alpha)$ induces the 
identity on $H^*(A)$ and $\varphi^* = \id$.

A Macaulay2 computation \cite{M2} yields
\[
\Ann_S \B_1(A) = \m^2 = (x_1^2, x_1 x_2, x_2^2),
\qquad
\B_1(A) \cong \k^3
\]
as a $\k$-vector space, with generators 
$[y_{12}^\vee]$, $[y_{112}^\vee]$, $[y_{212}^\vee]$ 
in $\B_1(A)/\m\B_1(A)$.
The dual action of $\varphi$ on $K_1(A) = A_1 \otimes_{\k} S$ 
sends $y_{12}^\vee \mapsto y_{12}^\vee - y_{112}^\vee$ 
with all other basis vectors fixed.
Since $y_{112}^\vee \otimes 1$ is a cycle and 
\[
(\varphi - \id)(y_{12}^\vee \otimes 1) 
= -y_{112}^\vee \otimes 1 \notin \im(\partial_2),
\]
the automorphism $\varphi$ acts nontrivially on $\B_1(A)$.
Since $\varphi^* = \id$, the maps $\B_1(\varphi)$ and 
$\B_1(\varphi^*) = \id$ are not $\B_1$-compatible, so 
$\varphi$ is not $\B_1$-formal.
\end{example}

%%%%%%%%%%%%%%%%%%%%%%%%%%%%%%%%%%%%
\subsection{Multiplicative $\B_q$-formality}
\label{subsec:mult-bq-formal}
%%%%%%%%%%%%%%%%%%%%%%%%%%%%%%%%%%%%

The question of whether $\B_q$-formality is strictly weaker 
than $q$-formality for morphisms of CE type suggests a 
multiplicative refinement of the notion itself. The naturality 
of cup and cap products established in 
Theorem~\ref{thm:naturality-alg-structures} provides the 
necessary framework.

\begin{definition}
\label{def:mult-bq-formal}
Let $\varphi\colon A\to\bar{A}$ be a $\cdga$ morphism between 
$q$-formal $\cdgas$ inducing an isomorphism on $H^1$. We say 
$\varphi$ is \emph{multiplicatively $\B_q$-formal} if 
it is $\B_q$-formal and moreover:
\begin{enumerate}[itemsep=2pt]
\item \label{mbf1}
For all $i,j\le q$ with $i+j\le q$, the formality isomorphisms 
$\B^i(A)\cong\B^i(H^*(A))$ intertwine the cup product maps, 
so that the diagram
\[
\begin{tikzcd}[column sep=30pt, row sep=22pt]
\B^i(A)\otimes_{S}\B^j(A) \ar[r, "\cup"] 
  \ar[d, "\cong"'] & 
\B^{i+j}(A) \ar[d, "\cong"] \\
\B^i(H^*(A))\otimes_{S}\B^j(H^*(A)) \ar[r, "\cup"'] & 
\B^{i+j}(H^*(A))
\end{tikzcd}
\]
commutes, and similarly for $\bar{A}$.
\item \label{mbf2}
For all $p,i\le q$ with $p+i\le q$, the formality isomorphisms 
intertwine the cap product action of $H^*(A)$ on $\B_\bullet(A)$:
\[
\begin{tikzcd}[column sep=30pt, row sep=22pt]
H^p(A)\otimes_{\k}\B_i(A) \ar[r, "\cap"] 
  \ar[d, "\cong"'] & 
\B_{i-p}(A) \ar[d, "\cong"] \\
H^p(H^*(A))\otimes_{\k}\B_i(H^*(A)) \ar[r, "\cap"'] & 
\B_{i-p}(H^*(A))
\end{tikzcd}
\]
commutes, and similarly for $\bar{A}$.
\end{enumerate}
\end{definition}

\begin{remark}
\label{rem:mbf-conditions}
Conditions~\eqref{mbf1} and~\eqref{mbf2} are a priori independent: 
\eqref{mbf1} concerns the $S$-algebra structure on cohomological 
Koszul modules, while \eqref{mbf2} concerns the $H^*(A)$-module 
$H^*(A)$-module structure on homological Koszul modules. When $A$ is a
$\PD_m$-$\cdga$, the isomorphism $\B_i(A)\cong\sigma^{*}\B^{m-i}(A)$ of
Theorem~\ref{thm:pd-koszul-modules} carries the cup product action of
$H^{*}(A)$ on $\B^{\bullet}(A)$ to its cap product action on
$\B_{\bullet}(A)$, by Proposition~\ref{prop:hstar-action}; hence
\eqref{mbf1} and \eqref{mbf2} are equivalent in that case.
In general, however, one could conceivably hold without the other.
\end{remark}

\begin{problem}
\label{prob:mult-bq-formal}
Does there exist a $\cdga$ morphism $\varphi\colon A\to\bar{A}$ between 
$q$-formal $\cdgas$ that is $\B_q$-formal but not 
multiplicatively $\B_q$-formal? By Remark~\ref{rem:mbf-conditions}, 
conditions~\eqref{mbf1} and~\eqref{mbf2} are equivalent for 
$\PD$-$\cdgas$, so any such morphism must arise outside the 
Poincaré duality setting. In particular, does the 
IA-automorphism of Example~\ref{ex:ia-auto-nonformal} 
yield such a counterexample for $q=1$?
\end{problem}

\begin{problem}
\label{prob:mult-bq-formal-bis}
Let $\varphi\colon A\to\bar{A}$ be a morphism between $q$-formal $\cdgas$.  
Does $q$-formality of $\varphi$ imply multiplicative $\B_q$-formality?  
Conversely, under what additional hypotheses, if any, does multiplicative 
$\B_q$-formality imply that $\varphi$ is $q$-formal?
\end{problem}

\subsection{Weight-preserving morphisms and $\B_q$-formality}
\label{subsec:wt-morphisms-bq}
We now identify a natural class of morphisms for which 
multiplicative $\B_q$-formality holds automatically for 
all $q$: those that preserve the positive-weight decompositions 
of Section~\ref{subsec:weights-cdga}. This provides a 
large supply of examples and clarifies when the open problems 
of Section~\ref{subsec:mult-bq-formal} are genuinely nontrivial.

\begin{definition}
\label{def:weight-preserving}
Let $(A,d_A)$ and $(\bar{A},d_{\bar{A}})$ be two connected $\cdgas$ 
with positive weights. A $\cdga$ morphism $\varphi\colon A\to\bar{A}$ 
is \emph{weight-preserving} if 
$\varphi(A^{i,\alpha})\subseteq \bar{A}^{i,\alpha}$ for all 
$i\ge 1$ and $\alpha>0$.
\end{definition}

\begin{theorem}
\label{thm:weight-preserving-formal}
Let $\varphi\colon A\to\bar{A}$ be a weight-preserving 
$q$-quasi-isomorphism between connected, finite-type $\cdgas$ 
with positive weights.
\begin{enumerate}[itemsep=2pt]
\item \label{wpf1}
The chain map $K_\bullet(\varphi)\colon K_\bullet(\bar{A})\to 
K_\bullet(A)$ of Proposition~\ref{prop:koszul-functorial} is compatible 
with the total-weight filtrations, and 
hence induces a morphism between the weight spectral sequences of 
Theorem~\ref{thm:weight-ss} for $A$ and for $\bar{A}$.  This morphism is 
$H_\bullet(\varphi)\otimes\id_S$ on $E^1$-pages and 
$\gr^W\B_\bullet(\varphi)$ on $E^\infty$-pages, and it is compatible with 
the algebra and coalgebra structures and with the $H^*$-actions.

\item \label{wpf2}
Assume in addition that $H^1(A)$ is weight-homogeneous, and that $A$ and 
$\bar{A}$ are $q$-formal, via zig-zags to which 
Corollary~\ref{cor:koszul-formal-separated} applies.  Then both spectral 
sequences collapse at $E^2$ in total degrees $\le q$, the morphism 
of\/~\eqref{wpf1} is $\B_\bullet(H^*(\varphi))$ on $E^2$-pages, and 
$\varphi$ is multiplicatively $\B_q$-formal.
\end{enumerate}
\end{theorem}

\begin{proof}
\eqref{wpf1}
Since $\varphi$ is weight-preserving, $K_\bullet(\varphi)$ carries 
$W^p K_\bullet(\bar{A})$ into $W^pK_\bullet(A)$, so it induces a morphism 
of the associated spectral sequences; this is 
Lemma~\ref{lem:koszul-ss-natural}\eqref{nat-weight}, which also supplies 
the identifications on $E^1$- and $E^\infty$-pages.  For the multiplicative structures: 
$\varphi\otimes\id_S$ is a morphism of weight-graded dg-algebras, the cup 
and cap product actions are weight-homogeneous by 
Proposition~\ref{prop:hstar-action}, and so the equivariance diagrams of 
Theorem~\ref{thm:naturality-alg-structures}\eqref{nat-alg3} hold at the 
weight-graded level.

\eqref{wpf2}
Since $\varphi$ induces a weight-preserving isomorphism on $H^1$, the 
hypothesis that $H^1(A)$ be weight-homogeneous passes to 
$H^1(\bar{A})$, so the second part of Theorem~\ref{thm:weight-ss} applies 
to both $\cdgas$: the $E^1$-pages are the Koszul chain complexes of 
$H^*(A)$ and $H^*(\bar{A})$, the induced map is 
$K_\bullet(H^*(\varphi))$, and taking homology gives 
$\B_\bullet(H^*(\varphi))$ on $E^2$-pages.

We claim that $q$-formality forces the two spectral sequences to collapse 
in total degrees $\le q$.  Indeed, by 
Corollary~\ref{cor:koszul-formal-separated}, $q$-formality provides 
$S$-module isomorphisms $\B_i(A)\cong\B_i(H^*(A))$ for $i\le q$; 
these are graded, and finite-dimensional in each weight because $A$ is of 
finite type.  Since 
$E^2_{p,q}\cong\B_{p+q}(H^*(A))^{\mathrm{weight}\,p}$ and 
$E^\infty_{p,q}\cong\gr^W_p\B_{p+q}(A)$, comparing dimensions weight by 
weight in total degrees $\le q$ leaves no room for a nonzero differential 
$d^r$ with $r\ge2$.  The same argument applies to $\bar{A}$.

Positive weights make $\B_i(A)$ and $\B_i(\bar{A})$ graded, so 
$\gr^W\B_i=\B_i$; the collapse therefore turns the edge maps into 
isomorphisms $\B_i(A)\cong\B_i(H^*(A))$ and 
$\B_i(\bar{A})\cong\B_i(H^*(\bar{A}))$ for $i\le q$, and by \eqref{wpf1} 
these carry $\B_i(\varphi)$ to $\B_i(H^*(\varphi))$.  This is precisely 
the commutativity of \eqref{eq:bi-compatible}, and the multiplicative 
statement follows from the last part of \eqref{wpf1}.
\end{proof}

\begin{remark}
\label{rem:wpf-collapse}
What part~\eqref{wpf2} uses is the collapse of the two weight spectral 
sequences in total degrees $\le q$, for which $q$-formality is sufficient; 
positive weights alone do not suffice, by 
Remark~\ref{rem:pw-not-enough}.  Under collapse, the edge isomorphisms 
constitute the canonical $\B_q$-framing of 
Proposition~\ref{prop:framed-bq}\eqref{fbq3}, so the conclusion 
of\/ \eqref{wpf2} is the commutativity of \eqref{eq:bi-compatible} for 
these framings---a statement involving no choice of zig-zag, and hence 
stronger than $\B_q$-formality in the sense of 
Proposition~\ref{prop:framed-bq}\eqref{fbq2}.
\end{remark}

\begin{corollary}
\label{cor:weight-preserving-geometric}
For each $q\ge 1$, any morphism in the following classes of $\cdgas$ that is a 
$q$-quasi-isomorphism satisfies the conclusions of 
Theorem~\ref{thm:weight-preserving-formal}\eqref{wpf1}; if moreover the 
$\cdgas$ involved are $q$-formal and have weight-homogeneous $H^1$, such a 
morphism is multiplicatively $\B_q$-formal.
\begin{enumerate}[itemsep=2pt]
\item \label{wpg1}
Chevalley--Eilenberg maps $\CE(f)\colon\CE(\bar{\g})\to\CE(\g)$ 
induced by Lie algebra morphisms $f\colon\g\to\bar{\g}$ between 
finite-dimensional nilpotent Lie algebras, equipped with their 
lower central series weight gradings (see Section~\ref{subsec:CE-weights}).
\item \label{wpg2}
Pullback maps $f^*\colon A(Y)\to A(X)$ between Morgan--Gysin 
models of smooth quasi-projective varieties, induced by 
morphisms $f\colon X\to Y$, with weights 
$\mathrm{wt}(A^{p,q})=p+2q$ (see Section~\ref{subsec:qproj-weights}).
\end{enumerate}
\end{corollary}

\begin{proof}
Both classes consist of weight-preserving morphisms, so
Theorem~\ref{thm:weight-preserving-formal}\eqref{wpf1} applies.  For the
last assertion it remains to identify the weight-homogeneous cases.  In
\eqref{wpg1}, $H^1(\CE(\g))$ is the dual of $\g/\g'$, hence
weight-homogeneous of weight~$1$.  In \eqref{wpg2}, $H^1$ is
weight-homogeneous precisely when $H^1(X;\k)$ is pure of weight~$1$.
\end{proof}

\begin{problem}
\label{prob:mult-bq-formal-CE}
Does the failure of multiplicative $\B_1$-formality 
for a CE automorphism $\varphi = \CE(\alpha)$ of a $1$-formal 
$A=\CE(\g)$---which by Theorem~\ref{thm:weight-preserving-formal} cannot 
occur when $\varphi$ is a weight-preserving quasi-iso\-morphism---have 
detectable consequences for the Hilbert series or support 
loci of $\B_1(A)$, beyond those already detected by the 
failure of $\B_1$-formality at the level of graded $S$-modules?
\end{problem}

%%%%%%%%%%%%%%%%%
\part{Resonance and holonomy}
\label{part:res-holo}
%%%%%%%%%%%%%%%%%

In this part, we develop the theory of resonance varieties and holonomy Lie algebras, 
systematically connecting them to the Koszul-theoretic framework established in the 
preceding parts. Resonance varieties arise as algebraic jump loci that control the 
behavior of Koszul complexes under specialization, providing the infinitesimal counterpart 
to the characteristic varieties of spaces and groups.

A primary focus of this part is the relationship between the geometry of these varieties 
and the underlying cohomology algebra. We establish tangent cone theorems that describe 
how cohomology governs the first-order structure of resonance at the origin, with formality 
forcing an exact agreement. We also examine how resonance is transformed under standard 
algebraic constructions, such as Hirsch extensions and tensor products, which frequently 
arise in the construction of Sullivan models.

The second half of the part centers on the first Koszul module, $\B_1(A)$. We provide 
an explicit presentation for this module,
and identify it with the infinitesimal Alexander invariant of the holonomy Lie algebra.
This identification bridges commutative algebra and Lie theory,
yielding closed-form formulas for holonomy Chen ranks that generalize
classical patterns from the $1$-formal setting.

%%%%%%%%%%%%%%%%%%%%%%%%%%%%%%%%%%%%%%
\section{Resonance varieties of $\cdgas$}
\label{sect:res-cdga}
%%%%%%%%%%%%%%%%%%%%%%%%%%%%%%%%%%%%%%

In this section we introduce resonance varieties associated to a 
$\cdga$ and develop their basic properties from the perspective of 
Koszul (co)homology. Resonance varieties arise as algebraic jump 
loci controlling the behavior of Koszul complexes under 
specialization, and provide a precise algebraic counterpart to the 
characteristic varieties of spaces and groups. They play a central 
role in relating homological invariants, support loci, and formality 
properties, and will later serve as the main algebraic input for the 
topological applications developed in Part~\ref{part:applications}.

After comparing the two families of loci (the determinantal jump loci 
$\RR^{i,s}(A)$ and the module-theoretic support loci $\RR_{i,s}(A)$), we 
examine the particularly rigid case in which resonance vanishes, and how 
that property behaves as the $\cdga$ structure varies on a fixed underlying 
graded algebra.  Their first-order comparison with the resonance of 
$H^*(A)$ is the subject of Section~\ref{sect:tcone}.

%%%%%%%%%%%%%%%%%%%%
\subsection{Resonance varieties and support loci}
\label{subsec:res}
%%%%%%%%%%%%%%%%%%%%

We begin by defining resonance varieties in terms of the Koszul complexes
associated to a $\cdga$. From this viewpoint, resonance arises as a homological
jump phenomenon: it records those linear forms for which the specialized
Koszul differential fails to be exact in a given degree.
Their relationship with module-theoretic support will be clarified below.

Let $(A,d_A)$ be a connected $\cdga$ over an algebraically closed field $\k$ 
of characteristic $0$. We assume $0<\dim_{\k} H^1(A)<\infty$, and write 
$S=\Sym(H_1(A))$, where $H_1(A)=H^1(A)^{\vee}$.
For $a \in H^1(A)$, let $\ev_a\colon S \to \k$ be the evaluation homomorphism 
induced by the canonical pairing $H_1(A)\otimes_{\k} H^1(A)\to \k$. We denote 
by $\m_a=\ker(\ev_a)\subset S$ the corresponding maximal ideal, and by 
$\kappa_a=S/\m_a$ the associated residue field. 
The Koszul cochain complex $K^\bullet(A)=(A^\bullet\otimes_{\k} S,\delta_A)$ of
Proposition~\ref{prop:LA} and its dual $K_\bullet(A)=(A_\bullet\otimes_{\k} S,
\partial^A)$, written with coefficients in $S$, are as in
Section~\ref{sect:koszul}.

\begin{definition}
\label{def:hom-res}
The \emph{resonance varieties} of $A$ are defined as the homology 
jump loci of the Koszul chain complex:
\[
\RR^{i,s}(A) \coloneqq
\bigl\{ a \in H^1(A) \mid
\dim_\k H_i(K_\bullet(A)\otimes_{S} \kappa_a) \ge s
\bigr\}.
\]
The \emph{support resonance varieties} of $A$ are
\[
\RR_{i,s}(A) \coloneqq \Supp_S\bigl(\bwedge^s\B_i(A)\bigr)
\subseteq  H^1(A).
\]
(We use the term \emph{variety} loosely here; that $\RR^{i,s}(A)$ 
is indeed a Zariski closed subset of $H^1(A)$ under suitable 
finiteness hypotheses is established in Proposition~\ref{prop:res-descriptions}\eqref{rd1}.) 
As usual we abbreviate $\RR^i(A)=\RR^{i,1}(A)$ and $\RR_i(A)=\RR_{i,1}(A)$.
\end{definition}

Note that $0 \in \RR^{i,s}(A)$ if and only if $s < b_i(A)$,
since $H_i(K_\bullet(A)\otimes_{S} \k_0)\cong H_i(A)$. 
If $b_j(A)<\infty$ for $j\le i$, then by the BGG truncation principle
(Proposition~\ref{prop:jump-dual-local}), the variety $\RR^{i,s}(A)$
depends only on the maps $\partial_i^A$ and $\partial_{i+1}^A$. 

\begin{example}
\label{ex:res-degree-zero}
For a connected $\cdga$ $(A,d)$, we have $\B_0(A)\cong\k$ as an
$S$-module, supported only at the origin.
Moreover, for any $a\in H^1(A)$,
$H_0(K_\bullet(A)\otimes_{S}\k_a)\cong\k$ for $a=0$
and $=0$ for $a\ne 0$.
Hence, although the two definitions differ---$\RR^{0,s}$ via
fiber dimensions of the Koszul complex and $\RR_{0,s}$ via
support of $\B_0(A)$---they yield the same answer:
\[
\RR_{0,s}(A) = \RR^{0,s}(A) =
\begin{cases} \{0\} & s=1, \\ \emptyset & s\ge 2.\end{cases}
\]
\end{example}

\begin{lemma}
\label{lem:resonance-truncation-dependence}
For any connected $\cdga$ $(A,d)$ and all $i\ge 0$, $s\ge 1$,
the resonance variety $\RR^{i,s}(A)$ depends only on the truncated $\cdga$
$A^{\le i+1}$.
\end{lemma}

\begin{proof}
Consider the Koszul complex $K^\bullet(A)=(A\otimes_{\k}S,\delta_A)$, 
where $S=\Sym(H_1(A))$. The Koszul differential $\delta_A^j$ depends 
only on $d^{\le j+1}$ and on the canonical tensor 
$\omega_A \in H^1(A)\otimes_{\k} H_1(A)$. 
Hence $\delta_A^j$ depends only on the $\cdga$ structure 
of $A$ in internal degrees $\le j+1$.

The resonance variety $\RR^{i,s}(A)$ is defined via the homology in degree
$i$ of the specializations $K^\bullet(A)\otimes_{S} S/\m_a$ at points $a\in H^1(A)$. 
For computing $H^i(K^\bullet(A)\otimes_{S} S/\m_a)$,
only the differentials $\delta_A^{i-1}$ and $\delta_A^{i}$ enter, hence only
$d^{\le i+1}$ and $\omega_A$ are relevant.
Therefore the truncated complexes $K^{\le i+1}(A)\otimes_{S} S/\m_a$
and $K^{\le i+1}(A^{\le i+1})\otimes_{S} S/\m_a$ coincide.
Since the $i$-th cohomology of a cochain complex depends only on its terms in
degrees $i-1$, $i$ and $i+1$, the groups
$H^i(K^\bullet(A)\otimes_{S} S/\m_a)$ depend only on
$A^{\le i+1}$, and hence so do the jump conditions defining
$\RR^{i,s}(A)$.
\end{proof}

For $a \in H^1(A)$, the fiber $K_\bullet(A)\otimes_{S} \kappa_a$
identifies, via Lemma~\ref{lem:two aom}, with the Aomoto complex
$(A,\delta_a)$. Hence resonance can be interpreted as the jump
locus of the Aomoto complex.

\begin{proposition}
\label{prop:res-descriptions}
Let $(A,d)$ be a connected $\cdga$ with $\dim_\k A^j < \infty$ 
for $j \le i$.  Then:
\begin{enumerate}[itemsep=3pt]
\item \label{rd1}
  The resonance variety $\RR^{i,s}(A)$ is Zariski closed,
  with determinantal scheme structure
  \[
  \RR^{i,s}(A) = \bV\bigl(I_{r_i-s+1}
    (\delta^{i-1}_A \oplus \delta^{i}_A)\bigr),
  \qquad r_i=\dim_\k A^i .
  \]
  If moreover $\dim_\k A^{i+1}<\infty$, the same formula
  holds with $\partial^A_i \oplus \partial^A_{i+1}$ in place
  of $\delta^{i-1}_A \oplus \delta^{i}_A$.
  
  \item \label{rd2}
  Let $N^i(A)=\coker(\delta^{i-1}_A)$.  Then
  \[
  \RR^{i,s}(A) = \Supp_S\Bigl(\bwedge^{\,s+r_{i+1}}
  \bigl(N^i(A)\oplus N^{i+1}(A)\bigr)\Bigr).
  \]
  Moreover $\Supp_S\bigl(\B^i(A)\bigr)\subseteq\RR^{i,1}(A)$, and
  \[
  \RR^{i,s}(A)\cap U_{i+1} = \Supp_S\bigl(\bwedge^s\B^i(A)\bigr)\cap U_{i+1}
  \]
  for all $s\ge1$, where $U_{i+1}\subseteq H^1(A)$ is the open locus on which
  $\delta^i_A$ attains its generic rank.  In particular, if $A^{i+1}=0$, then
  $\RR^{i,s}(A)=\Supp_S\bigl(\bwedge^s\B^i(A)\bigr)$ for all $s\ge1$.

\item \label{rd3}
$\RR^{i,s}(A)$ coincides with 
  the jump locus of the Aomoto complex in degree~$i$:
  \[
  \RR^{i,s}(A) = \bigl\{ a \in H^1(A) \mid 
  \dim_\k H^i(A, \delta_a) \ge s \bigr\}.
  \]
\end{enumerate}
\end{proposition}

\begin{proof}
Part~\eqref{rd1} follows from Proposition~\ref{prop:determinantal}
(and Proposition~\ref{prop:determinantal-finite} under the stronger
finiteness hypothesis). 

Part~\eqref{rd2} is Theorem~\ref{thm:jump-support}, applied to the chain
complex $C=K_\bullet(A)$, whose dual is $C^{\#}=K^\bullet(A)$; here
$\partial^{\#}_i=\delta^{i-1}_A$ and $r_{i+1}=\dim_\k A^{i+1}$.

Part~\eqref{rd3} follows from Lemma~\ref{lem:two aom} and
Corollary~\ref{corollary:jump-dual}, which identify the specialization
of the Koszul complex at $a$ with the Aomoto complex $(A,\delta_a)$.
Thus
\[
H_i(K_\bullet(A)\otimes_{S} \k_a) \cong H^i(A,\delta_a)^{\vee},
\]
and the claimed description follows.
\end{proof}

\begin{example}
\label{ex:torus-jump-support}
The restriction to $U_{i+1}$ in part~\eqref{rd2} cannot be removed.  For
$A=H^{*}(T^2;\k)=\bwedge(e_1,e_2)$ with $d_A=0$, the complex $K^\bullet(A)$
is the Koszul complex on $x_1,x_2$ over $S=\k[x_1,x_2]$, so that
$\B^1(A)=0$ and $\Supp_S\bigl(\B^1(A)\bigr)=\emptyset$.  All differentials of
the Aomoto complex $(A,\delta_a)$ vanish at $a=0$, however, so
$\dim_\k H^1(A,\delta_0)=2$ and $\RR^{1,1}(A)=\{0\}$.  The discrepancy is
confined to the origin, which is exactly the complement of $U_2$: the generic
rank of $\delta^1_A$ is $1$, while $\delta^1_A$ vanishes at $a=0$.  Here
$\B_1(A)=0$ as well, so $\RR_{1,1}(A)=\emptyset$ and the containment
$\RR_{1,1}(A)\subseteq\RR^{1,1}(A)$ is strict.
\end{example}

\begin{remark}
\label{rem:aomoto-pointwise}
As above, let $(A,d)$ be a connected, $i$-finite $\cdga$. 
Under Lemma~\ref{lem:two aom}, the fiber of the Koszul complex
at $a\in H^1(A)$ is the Aomoto complex $(A,\delta_a)$.
Thus, membership in $\RR^{i,s}(A)$ may be detected by the dimension
of $H^i(A,\delta_a)$, or equivalently, by applying rank--nullity to the
two-step complex
\[
\begin{tikzcd}[column sep=24pt]
A^{i-1}\ar[r, "\delta^{i-1}_a"] & A^i \ar[r, "\delta^{i}_a"] &  A^{i+1},
\end{tikzcd}
\]
one obtains
\begin{equation}
\label{eq:point-res}
a \in \RR^{i,s}(A)
\;\Longleftrightarrow\;
\rank \delta^{i-1}_a + \rank \delta^i_a
\le \dim_\k A^i - s,
\end{equation}
provided $\dim_\k A^{i+1}<\infty$.
\end{remark}

We now show that all constructions above admit a controlled finite 
model obtained from a finite truncation of the $\cdga$, without 
changing resonance or support data in low degrees.
This is a $\cdga$ analogue of \cite[Prop.~4.1]{PS-mrl}, where 
a similar Noetherian argument produces a finite CW-skeleton 
preserving characteristic varieties and Alexander invariants 
in low degrees; the role of the group ring is here played 
by the polynomial ring $S=\Sym(H_1(A))$.

\begin{theorem}
\label{thm:cdga-skeleton}
Let $(A,d)$ be a connected, $q$-finite $\k$-$\cdga$. There exist a connected, 
$(q+1)$-finite $\k$-$\cdga$ $(\oA, \bar{d})$ and a $\cdga$ morphism 
$\varphi \colon A \to \oA$ such that:
\begin{enumerate}[itemsep=3pt]
\item \label{sk1}
$\varphi^j \colon A^j \to \oA^j$ is an isomorphism 
for all $j \le q$, so in particular $\oA^j = A^j$ for $j\le q$.
\item \label{sk2}
$\varphi^* \colon H^*(A) \to H^*(\oA)$ is an isomorphism 
in degrees $\le q$; in particular, it induces an isomorphism 
$H^1(A) \isom H^1(\oA)$.  If, moreover, $\dim_\k H^{q+1}(A)<\infty$, then 
$\oA$ may be chosen so that $\varphi^*$ is also a monomorphism in degree 
$q+1$, and hence so that $\varphi$ is a $q$-quasi-isomorphism.
\item \label{sk3}
The morphism $\varphi$ induces isomorphisms $\B_i(\oA)\isom \B_i(A)$ 
of $S$-modules for all $i\le q$.
\item \label{sk4}
For all $i \le q$ and $s \ge 1$, 
\[
\varphi^*\colon \RR^{i,s}(\oA) \longisom \RR^{i,s}(A),
\qquad
\varphi^*\colon \RR_{i,s}(\oA) \longisom \RR_{i,s}(A).
\]
\item \label{sk5}
If $(A,d)$ is $q$-formal and $\dim_\k H^{q+1}(A)<\infty$, then 
$(\oA,\bar d)$ is also $q$-formal, and the morphism $A \to \oA$ can be 
chosen to fit into a $q$-formality zig-zag between $A$ and $H^*(A)$.
\end{enumerate}
\end{theorem}

\begin{proof}
\textit{Construction of $\oA$.}
Set $S=\Sym(H_1(A))$, a polynomial ring over $\k$, which is Noetherian.
The Koszul differential in degree $q$ takes the form
\[
\delta_A^q \colon K^q(A) = A^q\otimes_{\k} S 
\longrightarrow A^{q+1}\otimes_{\k} S,
\]
where $K^q(A)$ is a finitely generated free $S$-module (since $\dim_\k A^q<\infty$).
Since $S$ is Noetherian, the image
\[
Z_1 \coloneqq \im(\delta_A^q) \;\subseteq\; A^{q+1}\otimes_{\k} S
\]
is a finitely generated $S$-submodule.  Each generator of $Z_1$ involves
only finitely many basis vectors of $A^{q+1}$; let $V_1\subseteq A^{q+1}$
be the finite-dimensional $\k$-subspace spanned by all of them, so that
$Z_1\subseteq V_1\otimes_{\k} S$.

In order for $\oA^{q+1}$ to support a $\cga$ structure extending $A^{\le q}$,
we must also accommodate the images of all cup products landing in degree $q+1$.
Set
\[
V_2 \coloneqq \sum_{j=1}^{q} \mu\bigl(A^j\otimes_{\k} A^{q+1-j}\bigr)
\subseteq A^{q+1},
\]
where $\mu$ denotes the multiplication in $A$. 
This is a finite-dimensional $\k$-subspace of $A^{q+1}$, since $A^j$ is
finite-dimensional for all $j\le q$.  Finally, if
$\dim_\k H^{q+1}(A)<\infty$, choose a finite-dimensional subspace
$V_3\subseteq Z^{q+1}(A)$ of cocycles whose classes span $H^{q+1}(A)$, and
set $V_3=0$ otherwise.  Define $V = V_1 + V_2 + V_3\subseteq A^{q+1}$.
Then $V$ is finite-dimensional, $\im(\delta_A^q)\subseteq V\otimes_{\k} S$,
and $V$ is stable under all cup products from degrees $\le q$. 

Set
\[
\oA^j = 
\begin{cases}
A^j & j\le q, \\
V   & j = q+1, \\
0   & j\ge q+2,
\end{cases}
\]
with differential $\bar{d} = d|_{A^{\le q}}$ in degrees $\le q-1$,
and $\bar{d}^q\colon A^q\to V$ the map $d^q$ corestricted to $V$. 
This is well-defined since $\im(\delta_A^q)\subset V\otimes_{\k} S$, and writing
\[
\delta_A^q(a\otimes 1)=d^q(a)\otimes 1+\omega_A(a),
\quad \text{with } \omega_A(a)\in A^{q+1}\otimes_{\k} \m,
\]
shows that the $S_0$-component is $d^q(a)\in V$.
By construction, $V$ contains all products $A^j\cdot A^{q+1-j}$ for $1\le j\le q$.
Hence for such elements, the projection $A^{q+1}\surj V$ acts as the identity.
It follows that the multiplication on $\oA$ induced from $A$ is well defined 
and graded-commutative.
The Leibniz rule is inherited from $A$.  Thus $(\oA,\bar{d})$ is a connected
$(q+1)$-finite $\k$-$\cdga$. 

\smallskip
\textit{Construction of $\varphi$.}
Define a map $\varphi\colon A\to\oA$ by setting $\varphi^j = \id$ for $j\le q$ and
$\varphi^{q+1}$ the projection $A^{q+1}\surj V$ (and zero in degrees $\ge q+2$). 
This map is a $\cdga$ morphism, since $\bar{d}^q$ is $d^q$ followed by the 
projection, and the compatibility with multiplication in degrees $\le q+1$
follows from the fact that $V$ contains all products landing in degree $q+1$. 

\smallskip
\textit{Verification of \eqref{sk1}.}
By construction, $\varphi^j = \id_{A^j}$ for $j\le q$. 

\smallskip
\textit{Verification of \eqref{sk2}.}
For $j< q$: $Z^j(\oA) = Z^j(A)$ and $B^j(\oA) = B^j(A)$, so $H^j(\oA)\cong H^j(A)$.
For $j=q$: $Z^q(\oA)=\ker(\bar{d}^q)=\ker(d^q)=Z^q(A)$ and
$B^q(\oA)=\im(\bar{d}^{q-1})=\im(d^{q-1})=B^q(A)$, so $H^q(\oA)\cong H^q(A)$.
In degrees $\le q$, $\varphi^*$ is an isomorphism, as required.

In degree $q+1$, assume $\dim_\k H^{q+1}(A)<\infty$, so that
$V\supseteq V_3$ with the classes of $V_3$ spanning $H^{q+1}(A)$.  Since
$\oA^{q+2}=0$ we have $Z^{q+1}(\oA)=\oA^{q+1}=V$ and
$B^{q+1}(\oA)=\im(\bar d^q)$; moreover $\im(d^q)\subseteq V_1\subseteq V$
and $\pi|_V=\id$, where $\pi\colon A^{q+1}\surj V$ is the projection, so
that $\im(\bar d^q)=\pi(\im d^q)=\im(d^q)$.
Let $v\in Z^{q+1}(A)$ with $\varphi^*([v])=0$, that is,
$\pi(v)=d^q(a)$ for some $a\in A^q$.  Since the classes of $V_3$ span
$H^{q+1}(A)$, we may write $v=w+d^q(b)$ with $w\in V_3$ and $b\in A^q$.
Applying $\pi$ and using $\pi|_V=\id$ gives
$d^q(a)=\pi(v)=w+d^q(b)$, so $w=d^q(a-b)$ is a boundary; hence
$[v]=[w]=0$ in $H^{q+1}(A)$, and $\varphi^*$ is injective in degree $q+1$.

\smallskip
\textit{Verification of \eqref{sk3}.}
For $i \le q-1$, the differentials $\partial_i^{\oA}$ and $\partial_i^A$ 
coincide, since $\oA^j = A^j$ and $d_{\oA}^j = d_A^j$ for all $j \le q$, 
and the canonical element $\omega$ is the same under the identification 
$H^1(\oA)\cong H^1(A)$. Hence
$\B_i(\oA)=H_i(K_\bullet(\oA))=H_i(K_\bullet(A))=\B_i(A)$
as $S$-modules for all $i \le q-1$.

For $i = q$, we have
$\B_q(\oA)=\ker(\partial_q^{\oA})/\im(\partial_{q+1}^{\oA})$.
As above, $\partial_q^{\oA}=\partial_q^A$, so the kernels agree.
By construction of $\oA$, the image of
$\delta_A^q \colon A^q \to A^{q+1}\otimes_{\k} S$
is contained in $V\otimes_{\k} S = \oA^{q+1}\otimes_{\k} S$.
Dualizing, it follows that the images of
$\partial_{q+1}^{\oA}$ and $\partial_{q+1}^A$
coincide as submodules of $K_q(A)=K_q(\oA)$.
Therefore, $\B_q(\oA)\cong \B_q(A)$ as $S$-modules.

\smallskip
\textit{Verification of \eqref{sk4}.}
The identification $\B_i(\oA)\cong \B_i(A)$ implies $\RR_{i,s}(\oA)=\RR_{i,s}(A)$ 
via the description of support loci. 
Moreover, for each $a\in H^1(A)\cong H^1(\oA)$, the Aomoto complexes
$(A,\delta_a)$ and $(\oA,\delta_a)$ agree up to degree $q$, hence
$H^i(A,\delta_a)\cong H^i(\oA,\delta_a)$ for all $i\le q$. Therefore,
$\RR^{i,s}(\oA)=\RR^{i,s}(A)$. 

\smallskip
\textit{Verification of \eqref{sk5}.}
Suppose $A$ is $q$-formal and $\dim_\k H^{q+1}(A)<\infty$, that is,
there is a zig-zag of $q$-quasi-isomorphisms connecting $A$ to $(H^*(A),0)$.
Since $\varphi\colon A\to\oA $ is a $q$-quasi-isomorphism, the induced 
ring map $\varphi^*\colon (H^{*}(A),0)\to (H^{*}(\oA),0)$ is also 
a $q$-quasi-isomorphism. Therefore, 
\[
\begin{tikzcd}[column sep=18pt]
\oA&A\ar[l, "\ \varphi"'] \ar[r]& A_1 & \ar[l] \cdots & A_{\ell}\ar[r]\ar[l] 
& (H^*(A),0) \ar[r, "\varphi^*"] & (H^*(\oA),0)
\end{tikzcd}
\]
is a $q$-formality zig-zag for $\oA$, and so $\oA$ is $q$-formal.
\end{proof}

\begin{remark}
\label{rem:skeleton-mono}
The finiteness hypothesis in \eqref{sk2} and \eqref{sk5} cannot be
omitted.  Take $q\ge1$, let $Z$ be an infinite-dimensional $\k$-vector
space placed in degree $q+1$, and let
$A=\bigwedge(e)\otimes_{\k}(\k\oplus Z)$ with $\abs{e}=1$, $d=0$, and all
products of $Z$ with itself and with $e$ equal to zero.  Then $A$ is
$q$-finite, while $H^{q+1}(A)=A^{q+1}$ is infinite-dimensional; since
$\oA$ is $(q+1)$-finite by construction, no morphism $A\to\oA$ can be
injective on $H^{q+1}$.  Parts \eqref{sk1}, \eqref{sk3} and \eqref{sk4}
are unaffected: they concern only degrees $\le q$, where
$\varphi$ is the identity, and the rank computations of \eqref{sk4}
involve $\delta^q_A$, whose image lies in $V_1\otimes_{\k} S$ by
construction.
\end{remark}

%%%%%%%%%%%%%%%%%%%%
\subsection{Comparing resonance and support loci}
\label{subsec:res-comparison}
%%%%%%%%%%%%%%%%%%%%

We now compare the two constructions of resonance introduced above:
the homological jump loci $\RR^{i}(A)$ and the support loci $\RR_i(A)$.
Although defined by different procedures---base change versus Koszul 
module support---these loci are closely related. Under finite-type assumptions, 
they agree after taking unions over degrees, and differ in a controlled 
way in each degree.

\begin{theorem}
\label{thm:resonance-comparison}
If $\dim_\k A^j<\infty$ for all $j\le q$, then
\[
\bigcup_{i=0}^q \RR_i(A) =  \bigcup_{i=0}^q \RR^i(A).
\]
\end{theorem}

\begin{proof}
Apply Theorem~\ref{thm:PS-mrl-refined}\eqref{ps2} to the chain complex
$C=K_\bullet(A)$. The hypothesis $\dim_\k A^j<\infty$ for $j\le q$
ensures that $C_j=A_j\otimes_{\k} S$ is a finitely generated free $S$-module
for all $j\le q$. The conclusion yields the claimed equality. 
\end{proof}

The previous theorem identifies the resonance and support loci after
taking unions over all degrees up to $q$. We now refine this comparison
by examining the individual strata $\RR^{i,s}(A)$ and $\RR_{i,s}(A)$.
In general, these loci need not coincide degree by degree, but they are
closely related: in depth $s=1$ the support loci are contained in the
corresponding resonance loci, and in every depth the discrepancy is
controlled by the support in lower degrees.
This is made precise in Theorem~\ref{thm:jump-support-containment} below.

The comparison rests on two general facts established in 
Section~\ref{subsec:supp-jump}: Lemma~\ref{lem:supp-mu}, which reinterprets 
the support loci as the loci where the number of local generators of a Koszul 
module jumps, putting them on the same footing as the fiber dimensions that 
define $\RR^{i,s}(A)$; and Lemma~\ref{lem:fiber-exact}, which asserts that a 
complex of free modules over a local ring is exact at any spot where its fiber 
is.

\begin{theorem}
\label{thm:jump-support-containment}
Let $(A,d)$ be a connected $q$-finite $\k$-$\cdga$, and let $i\le q$.
\begin{enumerate}[itemsep=3pt]
  \item \label{jsc1}
  $\RR_{i,1}(A) \subseteq \RR^{i,1}(A)$.
  \item \label{jsc2}
  For every $s\ge 1$, the loci $\RR^{i,s}(A)$ and $\RR_{i,s}(A)$ agree
  outside the support loci in lower degrees:
  \[
  \bigl(\RR^{i,s}(A) \setminus \RR_{i,s}(A)\bigr)
  \cup
  \bigl(\RR_{i,s}(A) \setminus \RR^{i,s}(A)\bigr)
  \;\subseteq\;
  \bigcup_{r<i}\RR_{r,1}(A) .
  \]
  In particular, for $i=1$ we have
  $\RR^{1,s}(A)\setminus\{0\} = \RR_{1,s}(A)\setminus\{0\}$ for all $s\ge1$.
\end{enumerate}
\end{theorem}

\begin{proof}
\eqref{jsc1}
Fix $a\in H^1(A)\cong\mSpec S$.  By Lemma~\ref{lem:supp-mu} (with $s=1$) it
suffices to show that $H_i(K_\bullet(A)\otimes_{S}\kappa_a)=0$ implies
$\B_i(A)_a=0$.  Apply Lemma~\ref{lem:fiber-exact} to the local ring $R=S_a$
and to the complex $D_\bullet=K_\bullet(A)\otimes_{S} S_a$, whose homology in
degree $i$ is $\B_i(A)_a$ and whose reduction modulo the maximal ideal is
$K_\bullet(A)\otimes_{S}\kappa_a$.  The $q$-finiteness hypothesis guarantees
that $D_{i-1}$ and $D_i$ are finitely generated and free, while $D_{i+1}$ is
free; the lemma therefore gives $\B_i(A)_a=0$, as required.

\eqref{jsc2}
Suppose $a\notin\bigcup_{r<i}\RR_{r,1}(A)$, so that
$\B_r(A)\otimes_{S}\kappa_a=0$ for every $r<i$.  Since each such $\B_r(A)$ is
finitely generated over the Noetherian ring $S$, Nakayama's lemma
\eqref{eq:nakayama-supp} gives $\B_r(A)_a=0$, and hence, $\Tor$ being
compatible with localization,
\begin{equation}
\label{eq:tor-vanish-koszul}
E^2_{p,r} = \Tor_p^S\bigl(\B_r(A),\kappa_a\bigr) = 0
\qquad\text{for all $p\ge0$ and all $r<i$}
\end{equation}
in the change-of-rings spectral sequence
\[
E^2_{p,r} = \Tor_p^S\bigl(\B_r(A),\kappa_a\bigr)
\Longrightarrow H_{p+r}(K_\bullet(A)\otimes_{S} \kappa_a),
\]
which converges by Theorem~\ref{thm:PS-mrl-refined}\eqref{ps1}.  Two
consequences follow.  First, no differential $d^{r}$ with $r\ge2$ lands in
$E^{r}_{0,i}$, since its source $E^{r}_{r,\,i-r+1}$ is a subquotient of
$E^2_{r,\,i-r+1}=0$; and none emanates from $E^r_{0,i}$, for degree reasons.
Hence $E^\infty_{0,i}=E^2_{0,i}=\B_i(A)\otimes_{S}\kappa_a$.  Second,
$E^\infty_{p,i-p}=0$ for every $p\ge1$, being a subquotient of
$E^2_{p,i-p}=0$.  The finite filtration on
$H_i(K_\bullet(A)\otimes_{S}\kappa_a)$ thus collapses onto its bottom piece,
yielding a natural isomorphism
\[
\B_i(A)\otimes_{S} \kappa_a \cong
H_i(K_\bullet(A)\otimes_{S}\kappa_a).
\]
By Lemma~\ref{lem:supp-mu}, it follows that $a \in\RR^{i,s}(A)$ if and only
if $a\in\RR_{i,s}(A)$, which is the asserted containment.  For $i=1$ the
exceptional locus is $\RR_{0,1}(A)=\{0\}$, by
Example~\ref{ex:res-degree-zero}.
\end{proof}

The proof shows that the discrepancy between $\RR^{i,s}(A)$ and
$\RR_{i,s}(A)$ is governed by the failure of degeneration of the
change-of-rings spectral sequence at the $E^2$-page. In particular,
this failure can occur only at points where some lower Koszul module
$\B_r(A)$, with $r<i$, fails to vanish after specialization.  Since the
support loci do not propagate (Remark~\ref{rem:propagation-vs-support}), the
union in part~\eqref{jsc2} cannot in general be replaced by its top term
$\RR_{i-1,1}(A)$.  Outside that union, the two notions coincide in degree
$i$, a fact we now record.

\begin{proposition}
\label{prop:generic-resonance-agreement}
Let $(A,d)$ be a connected $q$-finite $\k$-$\cdga$, and fix $i \le q$, $s \ge 1$.
Then the complement
\[
U_{<i} \coloneqq H^1(A)\setminus \bigcup_{r<i}\RR_{r,1}(A)
\]
is a Zariski open subset of $H^1(A)$, and on $U_{<i}$ one has
\[
\RR^{i,s}(A)\cap U_{<i} = \RR_{i,s}(A)\cap U_{<i}.
\]
\end{proposition}

\begin{proof}
Each $\B_r(A)$ with $r<i\le q$ is a finitely generated $S$-module, so
$\RR_{r,1}(A)=\Supp_S\B_r(A)$ is Zariski closed; hence the finite union
$\bigcup_{r<i}\RR_{r,1}(A)$ is closed and $U_{<i}$ is open.  The asserted
equality is precisely the content of
Theorem~\ref{thm:jump-support-containment}\eqref{jsc2}.
\end{proof}

The following two examples show that the containment in
Theorem~\ref{thm:jump-support-containment}\eqref{jsc2} can be strict on the
exceptional locus $\bigcup_{r<i}\RR_{r,1}(A)$, and that
part~\eqref{jsc1} does not extend to depths $s\ge2$.

\begin{example}
\label{ex:sol2-resonance}
Let $A=\CE(\sol_2)$ be the $\cdga$ from Section~\ref{subsec:sol2}.
As computed in Remark~\ref{rem:propagation-vs-support},
\[
\RR^1(A)=\{0,1\} \quad\text{while}\quad \RR_1(A)=\{1\}, 
\]
so the jump and support loci differ in degree $1$.
Nevertheless,
\[
\RR^0(A)\cup \RR^1(A)=\RR_0(A)\cup \RR_1(A)=\{0,1\},
\]
in agreement with Theorem~\ref{thm:resonance-comparison}.
\end{example}

\begin{example}
\label{ex:free-group-depth}
Let $A=H^{*}(\bigvee_{4}S^1;\k)$, so that $A^0=\k$, $A^1=\k^4$, and
$A^{\ge2}=0$; equivalently, $A$ is the exterior Stanley--Reisner ring of four
isolated points, as in Example~\ref{ex:sr-ring}.  The Koszul chain complex
reduces to $T^4\to T$, with $\partial_1=\pm(e_1\ \cdots\ e_4)$, so $\B_1(A)$
is the module of syzygies of $(e_1,\dots,e_4)$, minimally generated by the
$\binom{4}{2}=6$ Koszul syzygies.  Hence
$\dim_\k\bigl(\B_1(A)\otimes_{S}\kappa_0\bigr)=6$, whereas
$\dim_\k H_1\bigl(K_\bullet(A)\otimes_{S}\kappa_0\bigr)=\dim_\k A^1=4$, since
all differentials vanish after specializing at the origin.  For $a\neq 0$, on
the other hand, the localization $\B_1(A)_a$ is free of rank $3$ and
$\dim_\k H_1(K_\bullet(A)\otimes_{S}\kappa_a)=3$.  Lemma~\ref{lem:supp-mu}
therefore gives
\[
\RR_{1,s}(A)=
\begin{cases}
H^1(A) & s\le 3,\\
\{0\} & 4\le s\le 6,\\
\emptyset & s\ge 7,
\end{cases}
\qquad
\RR^{1,s}(A)=
\begin{cases}
H^1(A) & s\le 3,\\
\{0\} & s=4,\\
\emptyset & s\ge 5 .
\end{cases}
\]
In particular, $\RR_{1,5}(A)=\{0\}\not\subseteq\emptyset=\RR^{1,5}(A)$, so
that the containment of part~\eqref{jsc1} fails in depth $s=5$.  The two
families do agree away from the origin, in accordance with
part~\eqref{jsc2}.
\end{example}

\begin{remark}
\label{rem:scheme-structures-res}
The resonance variety $\RR^{i,s}(A)$ carries two natural scheme 
structures:
\begin{enumerate}[itemsep=2pt]
\item the \emph{Fitting scheme structure}, defined by the Fitting ideal
$\Fitt_{s+r_{i+1}-1}\bigl(N^i(A)\oplus N^{i+1}(A)\bigr)
=I_{r_i-s+1}(\delta^{i-1}_A\oplus\delta^i_A)$, the second expression being the
ideal of $(r_i-s+1)$-minors of the Koszul coboundary block matrix;
\item the \emph{cohomological support scheme structure}, 
defined by $\Ann_S(\bwedge^s\B^i(A))$.
\end{enumerate}

These two do not in general define the same underlying set: by 
Proposition~\ref{prop:res-descriptions}\eqref{rd2}, the second is contained 
in $\RR^{i,1}(A)$ when $s=1$, and agrees with $\RR^{i,s}(A)$ only over the 
locus where $\delta^i_A$ has generic rank, the discrepancy being illustrated 
in Example~\ref{ex:torus-jump-support}.  What the Fitting structure does 
compute, by Proposition~\ref{prop:res-descriptions}\eqref{rd2} again, is the 
support of an exterior power of $N^i(A)\oplus N^{i+1}(A)$, the pair of 
cokernels of the Koszul coboundaries.
Even where the underlying sets agree, the two ideals may differ: the inclusion 
$I_{r_i-s+1}(\delta^{i-1}_A \oplus \delta^i_A) 
\subseteq \Ann_S(\bwedge^s \B^i(A))$ 
may be strict, as illustrated in \cite[Ex.~5.3]{AFRSS24}. 
The support locus $\RR_{i,s}(A)$, by contrast, 
does not arise from a base-change (jump locus) construction, 
and carries a single natural scheme structure, defined by 
$\Ann_S(\bwedge^s\B_i(A))$.
The relationship between the scheme structures on $\RR^{i,s}(A)$ 
and $\RR_{i,s}(A)$ will be revisited in Remark~\ref{rem:scheme-difference}, 
after a set-theoretic comparison of the two loci.
\end{remark}

\begin{remark}
\label{rem:scheme-difference}
Even when the underlying sets $\RR_i(A)$ and $\RR^i(A)$ coincide, 
their scheme structures may differ substantially, as the annihilator 
ideal $\Ann_S(\bwedge^s\B_i(A))$ and the Fitting ideal 
$\Fitt_{s+r_{i+1}-1}(N^i(A)\oplus N^{i+1}(A))$ define in general 
different subschemes of $H^1(A)$; see 
Example~\ref{ex:sol2-resonance} and \cite[Ex.~5.3]{AFRSS24}.

Set-theoretically, Theorem~\ref{thm:resonance-comparison} gives 
equality of unions $\bigcup_{i=0}^q\RR_i(A)=\bigcup_{i=0}^q\RR^i(A)$ 
under finite-type hypotheses, 
but $\RR^{i,s}(A)$ and $\RR_{i,s}(A)$ need not coincide degreewise
for $i\ge 2$. For $i=1$, Theorem~\ref{thm:jump-support-containment}
shows that the discrepancy is supported at the origin; accordingly, 
the annihilator scheme structures on $\RR^{1,s}(A)$ 
and $\RR_{1,s}(A)$ agree away from $0$, for every $s\ge 1$. As 
Example~\ref{ex:free-group-depth} shows, they may well differ at the 
origin, in every depth.

The determinantal scheme structure on $\RR^{i,s}(A)$ arises 
naturally from the linear algebra of the cochain complex and is 
well suited for detecting cohomology jump loci. The annihilator 
scheme structure on $\RR_{i,s}(A)$ encodes the support of the 
Koszul module $\B_i(A)$ and is the one relevant for questions 
concerning decomposition, separability, and asymptotic behavior 
of homological invariants. In particular, for problems related 
to Chen ranks and their infinitesimal counterparts, it is the 
annihilator-defined scheme on $\RR_1(A)$ that governs whether 
additive formulas hold in high degrees: under suitable hypotheses 
(see Sections~\ref{subsec:separable-isotropic} 
and~\ref{subsec:chen-ranks-conj}), reducedness of this scheme 
governs the extent to which $\B_1(A)$ surjects onto a direct 
sum of modules supported on the irreducible components of 
resonance, with the length of the kernel controlling the 
effective range of the Chen ranks formula~\cite{AFRS25}.

Since the results of this paper concern primarily the underlying 
varieties and their geometric properties---such as linearity, 
conicality, and disjointness of components---we will not pursue 
these scheme-theoretic distinctions further here.
\end{remark}

%%%%%%%%%%%%%%%%%%%%%
\subsection{Vanishing resonance and deformations of the differential}
\label{subsec:vanishing-res-cdga}
%%%%%%%%%%%%%%%%%%%%%

Particularly rigid---and algebraically signifi\-cant---is the situation where
resonance varieties vanish.  It is natural to ask how this property behaves as
the $\cdga$ structure varies on a fixed underlying graded algebra.  As we now
explain, the answer depends sharply on what is held fixed, and in the most
obvious formulation there is no room for genericity at all.

Let $A^\bullet$ be a finite-type, connected graded-commutative $\k$-algebra, and
let
\[
\mathcal{D}(A)=\bigl\{d\colon A^\bullet\to A^{\bullet+1} \bigm|
d \text{ a degree-}1 \text{ derivation with } d^2=0 \bigr\}
\]
be the affine cone of $\cdga$ structures carried by $A^\bullet$.

\begin{definition}
\label{def:diff-cga}
Fix an integer $q\ge 1$. Let $\mathcal{D}_q(A)\subseteq\mathcal{D}(A)$ denote the
set of differentials $d$ such that $(A,d)$ is $q$-quasi-isomorphic to $(A,0)$.
\end{definition}

\begin{proposition}
\label{prop:Dq-rigid}
The set $\mathcal{D}_q(A)$ consists exactly of those $d\in\mathcal{D}(A)$ with
$d|_{A^j}=0$ for all $j\le q$; in particular, it is a Zariski-closed subcone of
$\mathcal{D}(A)$.  Moreover, for every $d\in\mathcal{D}_q(A)$, every $i\le q$ and
every $s\ge1$,
\[
\B^i(A,d)=\B^i(A^\bullet)
\qquad\text{and}\qquad
\RR^{i,s}(A,d)=\RR^{i,s}(A^\bullet).
\]
\end{proposition}

\begin{proof}
Suppose $(A,d)$ is $q$-quasi-isomorphic to $(A,0)$, so that
$\dim_\k H^j(A,d)=\dim_\k A^j$ for all $j\le q$.  We argue by induction on
$j\le q$ that $d|_{A^j}=0$.  For $j=1$, $H^1(A,d)=\ker d|_{A^1}$ has dimension
$\dim_\k A^1$ only if $d|_{A^1}=0$.  If $d|_{A^{j-1}}=0$, then
$H^j(A,d)=\ker d|_{A^j}$, which again has dimension $\dim_\k A^j$ only if
$d|_{A^j}=0$.  Conversely, suppose $d$ vanishes in degrees $\le q$.  
The subspace $A^{>q}$ is an ideal of 
$A^\bullet$ with $d(A^{>q})\subseteq A^{>q}$, hence a differential
ideal of $(A,d)$ and of $(A,0)$ alike, and the quotient $A/A^{>q}$ carries the
zero differential.  Both projections $(A,d)\to A/A^{>q}$ and
$(A,0)\to A/A^{>q}$ induce the identity of $A^j$ on $H^j$ for every $j\le q$,
and are surjective in degree $q+1$, the target vanishing there; so both are
$q$-quasi-isomorphisms, and together they join $(A,d)$ to $(A,0)$.  Note that
the identity of $A$ itself is not available for this purpose: it is a $\cdga$
map $(A,0)\to(A,d)$ only if $d$ vanishes in all degrees.  Finally, the
conditions $d|_{A^j}=0$ are linear, so $\mathcal{D}_q(A)$ is closed in
$\mathcal{D}(A)$, and stable under scaling.

For the second assertion, the canonical element $\omega_A$ and the ring
$S=\Sym(H_1(A,d))=\Sym\bigl((A^1)^\vee\bigr)$ are the same for all
$d\in\mathcal{D}_q(A)$, and
$\delta^j_{(A,d)}=d|_{A^j}\otimes\id+\omega_A\cdot(-)=\omega_A\cdot(-)
=\delta^j_{(A,0)}$ for $j\le q$.  Since $\B^i$ and, by
Proposition~\ref{prop:res-descriptions}\eqref{rd1}, the loci $\RR^{i,s}$ in
degree $i$ depend only on $\delta^{i-1}$ and $\delta^{i}$, both agree with their
counterparts for $(A,0)$ when $i\le q$.
\end{proof}

Thus, on $\mathcal{D}_q(A)$ the resonance varieties in degrees $\le q$ do not
depend on the differential at all: vanishing resonance holds for every
$d\in\mathcal{D}_q(A)$, or for none, according to whether
$\RR^i(A^\bullet)=\{0\}$.  Any genuine genericity statement must therefore either
vary the multiplication instead of the differential, as the known results do,
or else allow $d$ to change the cohomology, a possibility that
seems not to have been studied.

\begin{remark}
\label{rem:generic-zero-res-cdga}
Genericity results of the first kind were obtained by Papadima and Suciu
\cite{PS-mrl, PS-crelle}, who fix the graded vector space
$A=\k\oplus A^1\oplus A^2$ and let the graded-commutative multiplication
$\mu\colon\bwedge^2A^1\to A^2$ vary.  The locus of multiplications with
$\RR^1(A_\mu)\subseteq\{0\}$ is Zariski open, by
\cite[Thm.~3.3]{PS-mrl}; and it is non-empty precisely when
\begin{equation}
\label{eq:generic-vanishing-numerical}
\dim_\k A^2 \ge 2\dim_\k A^1-3 ,
\end{equation}
by \cite[Thm.~1.1]{PS-crelle}, whose criterion $m\ge 2n-3$ for an
$m$-dimensional space of relations in $\bwedge^2$ of an $n$-dimensional space
takes the form \eqref{eq:generic-vanishing-numerical} here; see
\cite[Rem.~3.4]{PS-mrl}.  There the
differential is zero throughout, and what varies is the Aomoto complex itself;
the Koszul spectral sequence plays no role.
\end{remark}

To formulate the remaining case precisely, we stratify $\mathcal{D}(A)$ by Betti
numbers.  Fix an integer $i\ge1$, and for $d\in\mathcal{D}(A)$ let
\begin{equation}
\label{eq:betti-vector}
\mathbf{b}_{\le i+1}(d)
=\bigl(\dim_\k H^1(A,d),\,\dots,\,\dim_\k H^{i+1}(A,d)\bigr) \in\N^{i+1} .
\end{equation}
For $\mathbf{b}\in\N^{i+1}$, put
\begin{equation}
\label{eq:betti-stratum}
\mathcal{D}(A)^{\mathbf{b}}
=\bigl\{d\in\mathcal{D}(A) \bigm| \mathbf{b}_{\le i+1}(d)=\mathbf{b}\bigr\} .
\end{equation}
The degree range is recorded in the length of\/ $\mathbf{b}$, so the stratum
needs no further decoration; the superscript distinguishes these strata from the
subcones $\mathcal{D}_q(A)$ of Definition~\ref{def:diff-cga}.  Since
$\dim_\k H^j(A,d)=\dim_\k A^j-\rank\bigl(d|_{A^j}\bigr)
-\rank\bigl(d|_{A^{j-1}}\bigr)$, and since $d|_{A^0}=0$, prescribing
$\mathbf{b}$ is the same as prescribing the ranks of $d$ in degrees $\le i+1$;
each condition $\rank(d|_{A^j})=r_j$ being locally closed, the strata
$\mathcal{D}(A)^{\mathbf{b}}$ are locally closed subvarieties, and
$\mathcal{D}(A)$ is the disjoint union of the finitely many nonempty ones.
On a fixed stratum, the spaces $H^1(A,d)=\ker\bigl(d|_{A^1}\bigr)$ all have
dimension $b_1$, so they are the fibers of a rank-$b_1$ subbundle
$\mathcal{H}^{\mathbf{b}}$ of the trivial bundle
$\mathcal{D}(A)^{\mathbf{b}}\times A^1$; and the resonance varieties
$\RR^i(A,d)\subseteq H^1(A,d)$ are the fibers of a closed subscheme
$\mathcal{R}^{\mathbf{b}}\subseteq\mathcal{H}^{\mathbf{b}}$, cut out fiberwise by
the determinantal ideals of
Proposition~\ref{prop:res-descriptions}\eqref{rd1}.

\begin{problem}
\label{prob:generic-vanishing-diff}
With the notation above, is
\[
\bigl\{d\in\mathcal{D}(A)^{\mathbf{b}} \bigm| \RR^i(A,d)\subseteq\{0\}\bigr\}
\]
a nonempty Zariski-open subset of\/ $\mathcal{D}(A)^{\mathbf{b}}$; equivalently,
does $\mathcal{R}^{\mathbf{b}}$ meet the generic fiber of\/
$\mathcal{H}^{\mathbf{b}}$ only in the zero section?  By the tangent cone
inclusion of Theorem~\ref{thm:tangent-cone} below, it would suffice that
$\RR^i\bigl(H^*(A,d)\bigr)=\{0\}$ for generic $d$ in the stratum, together with
conicality of $\RR^i(A,d)$; both conditions are open questions.
\end{problem}

%%%%%%%%%%%%%%%%%%%%%%%%%%%%%%%%%%%%%%
\section{Tangent cones and positive weights}
\label{sect:tcone}
%%%%%%%%%%%%%%%%%%%%%%%%%%%%%%%%%%%%%%

The main structural results of this section are the tangent cone 
theorems stated as Theorem~\ref{thm:tangent-cone-intro} in the 
Introduction: the inclusion 
$\TC_0(\RR^{q,s}(A)) \subseteq \RR^{q,s}(H^*(A))$, 
valid for all $q$-finite $\cdgas$ 
(Theorem~\ref{thm:tangent-cone}), and its upgrade to an equality of 
germs at the origin under $q$-formality (Theorem~\ref{thm:formal-tc}). 
Together these results show that cohomology governs the 
first-order infinitesimal structure of resonance, while 
formality forces exact agreement.  In the presence of positive weights the 
loci are conical, and the inclusion becomes global.

%%%%%%%%%%%%%%%%%%%%%%%
\subsection{The tangent cone inclusion}
\label{subsec:tcone}
%%%%%%%%%%%%%%%%%%%%%%%

We turn to the tangent cone inclusion theorem.  The 
cohomological half is proved by Gaussian elimination on the determinantal 
matrices of Proposition~\ref{prop:res-descriptions}\eqref{rd1}: splitting 
off the constant part of the differential exhibits, for every minor of the 
Aomoto matrix of $H^*(A)$, a function on $H^1(A)$ vanishing on 
$\RR^{i,s}(A)$ whose lowest-degree term is that minor.  The homological 
half combines the tangent cone identification of 
Lemma~\ref{lem:tc-support-general} with an elementary embedding of 
$\gr^\m\B_i(A)$ into the second page of the Koszul spectral sequence.

\begin{theorem}
\label{thm:tangent-cone}
Let $(A,d_A)$ be a connected $\k$-$\cdga$ over an 
algebraically closed field $\k$ of characteristic $0$, with 
$\dim_\k A^j<\infty$ for all $j\le i$.  Then:
\begin{enumerate}[itemsep=3pt]
\item \label{tc-coh}
For every $s\ge 1$,
\begin{equation}
\label{eq:tcone-coh}
\TC_0\bigl(\RR^{i,s}(A)\bigr) \subseteq \RR^{i,s}\bigl(H^*(A)\bigr).
\end{equation}
At the level of ideals: writing $I=I_{r_i-s+1}
(\delta_A^{i-1}\oplus\delta_A^{i})$ for the determinantal ideal of 
Proposition~\textup{\ref{prop:res-descriptions}\eqref{rd1}} and 
$b_i=\dim_\k H^i(A)$, the ideal of initial forms of the localization 
$I\cdot S_\m$ contains the corresponding determinantal ideal of the 
cohomology algebra:
\begin{equation}
\label{eq:tcone-coh-scheme}
\In\bigl(I\cdot S_\m\bigr)
\supseteq
I_{b_i-s+1}\bigl(\delta_{H^*(A)}^{i-1}\oplus\delta_{H^*(A)}^{i}\bigr).
\end{equation}
\item \label{tc-hom}
For $s=1$,
\begin{equation}
\label{eq:tcone-hom}
\TC_0\bigl(\RR_{i,1}(A)\bigr) 
= \Supp_S\bigl(\gr^\m\B_i(A)\bigr)
\subseteq \RR_{i,1}\bigl(H^*(A)\bigr),
\end{equation}
and for every $s\ge1$,
$\TC_0\bigl(\RR_{i,s}(A)\bigr)
=\Supp_S\bigl(\gr^\m({\textstyle\bwedge^s}\B_i(A))\bigr)$.
\end{enumerate}
\end{theorem}

\begin{proof}
\eqref{tc-coh}  Choose splittings of the (finite-dimensional) graded 
pieces
\[
A^j = B^j\oplus \wtH^j\oplus C^j, \qquad j=i-1,\, i,
\]
where $B^j=\im(d_A^{j-1})$, the subspace $\wtH^j$ is a complement of 
$B^j$ in $\ker(d_A^j)$, so that $\wtH^j\cong H^j(A)$, and $C^j$ is a 
complement of $\ker(d_A^j)$ in $A^j$; choose likewise a decomposition 
$A^{i+1}=B^{i+1}\oplus W^{i+1}$ with $B^{i+1}=\im(d_A^{i})$, allowing 
$\dim_\k A^{i+1}=\infty$.  Then $d_A$ vanishes on $B^j\oplus\wtH^j$ and 
restricts to isomorphisms $D_j\colon C^j\isom B^{j+1}$; note that 
$\dim_\k B^{j+1}=\dim_\k C^j<\infty$ even for $j=i$.

Fix cocycle representatives $\tilde e_1,\dots,\tilde e_n\in Z^1(A)$ of a 
basis of $H^1(A)$, and write the matrix of $\delta_A^j=d_A\otimes\id+
\sum_{k} \tilde e_k\cdot(-)\otimes x_k$ in block form with respect to these 
decompositions, as $M_j = D_j' + L_j$, where $D_j'$ is the constant 
matrix of $d_A$---whose only nonzero block is the invertible 
$D_j\colon C^j\to B^{j+1}$---and $L_j$ has entries in $S_1$.  Three block 
entries of $L_j$ vanish identically, as a formal consequence of 
$\delta_A^2=0$.  Indeed, write $\delta_A=\delta_0+\delta_+$ with 
$\delta_0=d_A\otimes\id$ and $\delta_+=\omega_A\cdot(-)$; then 
$\delta_0^2=0$ and $\delta_+^2=0$, the latter because $\omega_A^2=0$, 
whence
\begin{equation}
\label{eq:anticomm}
\delta_0\delta_++\delta_+\delta_0=0 .
\end{equation}
If $z\in B^j$, say $z=\delta_0y$, then $\delta_+z=-\delta_0\delta_+y
\in\im\delta_0=B^{j+1}$, so the blocks $B^j\to\wtH^{j+1}$ and 
$B^j\to C^{j+1}$ of $L_j$ are zero.  If $z\in\wtH^j\subseteq
\ker\delta_0$, then $\delta_0\delta_+z=-\delta_+\delta_0z=0$, so 
$\delta_+z\in\ker\delta_0=B^{j+1}\oplus\wtH^{j+1}$ and the block 
$\wtH^j\to C^{j+1}$ of $L_j$ is zero as well.  Moreover, the block 
$\wtH^j\to\wtH^{j+1}$ of 
$L_j$ is, by construction, the matrix of the Aomoto differential 
$\delta_{H^*(A)}^{j}=[\omega_A]\cup(-)$ of the cohomology algebra: 
multiplication by $\tilde e_k$ on cocycles descends to multiplication by 
$[e_k]$ on cohomology.

The block $D_j+ (L_j)_{C^j\to B^{j+1}}$ is a square matrix over $S$ of 
finite size whose value at $0$ is the invertible matrix $D_j$; hence it 
is invertible over the local ring $S_\m$.  Performing Gaussian 
elimination over $S_\m$ with this pivot block---that is, multiplying 
$M_j$ on both sides by unipotent block matrices over $S_\m$---brings 
$M_j$ to the form, in which $\delta^{j}_{\mathrm{red}}$ denotes the 
resulting Schur complement, the {\em reduced differential}\/ in 
degree~$j$:
\[
M_j \;\sim\;
\begin{pmatrix}
\,\mathrm{Id}_{C^j} & 0\, \\
\,0 & \delta^{j}_{\mathrm{red}}\,
\end{pmatrix},
\qquad
\delta^{j}_{\mathrm{red}}\colon (B^j\oplus\wtH^j)\otimes_{\k} S_\m
\longrightarrow (\wtH^{j+1}\oplus W')\otimes_{\k} S_\m,
\]
where $W'=C^{i+1}$ or $W^{i+1}$ according to the degree.  By the vanishing 
of the three linear blocks recorded above, and since the elimination 
alters entries only by terms in $\m^2 S_\m$, the Schur complement has the 
form
\begin{equation}
\label{eq:schur-form}
\delta^{j}_{\mathrm{red}} =
\begin{pmatrix}
\,0 & \delta_{H^*(A)}^{j}\, \\
\,0 & 0\,
\end{pmatrix}
\;+\; Q_j,
\qquad \text{all entries of } Q_j \text{ in } \m^2 S_\m .
\end{equation}

Now let $I=I_{r_i-s+1}(M_{i-1}\oplus M_i)$.  Since Fitting ideals are 
invariant under multiplication by invertible matrices (by the 
Cauchy--Binet formula, any minor of the transformed matrix is an 
$S_\m$-linear combination of same-size minors of the original), and 
since deleting an identity block of total size 
$c=\dim C^{i-1}+\dim C^i$ lowers the relevant minor size by exactly $c$, 
while $r_i=c+b_i$, we obtain
\[
I\cdot S_\m = I_{b_i-s+1}\bigl(\delta^{i-1}_{\mathrm{red}} \oplus 
\delta^{i}_{\mathrm{red}}\bigr)\cdot S_\m .
\]
Let $h$ be any nonzero $(b_i-s+1)$-minor of 
$\delta_{H^*(A)}^{i-1}\oplus\delta_{H^*(A)}^{i}$, a homogeneous 
polynomial of degree $b_i-s+1$, and let $f$ be the minor of 
$\delta^{i-1}_{\mathrm{red}}\oplus\delta^{i}_{\mathrm{red}}$ on 
the same rows and columns, taken inside the $\wtH$-blocks.  
Expanding $f$ using \eqref{eq:schur-form}: the term in 
which every entry contributes its linear part is $h$, of degree 
$b_i-s+1$, and every other term involves at least one entry of $Q$ and 
so has degree $\ge b_i-s+2$; hence $\In(f)=h$.  Clearing the denominator 
of $f$ by a unit $u\in S_\m$ with $u(0)\neq0$ produces $g=uf\in I$ with 
$\In(g)=u(0)\,h$.  This proves \eqref{eq:tcone-coh-scheme}, and since 
initial forms of an ideal are insensitive to localization at $\m$, we 
conclude
\[
\TC_0\bigl(\RR^{i,s}(A)\bigr)
=\bV\bigl(\In(I)\bigr)
\subseteq
\bV\Bigl(I_{b_i-s+1}\bigl(\delta_{H^*(A)}^{i-1}\oplus
\delta_{H^*(A)}^{i}\bigr)\Bigr)
=\RR^{i,s}\bigl(H^*(A)\bigr),
\]
the last equality being 
Proposition~\ref{prop:res-descriptions}\eqref{rd1} for the $\cdga$ 
$(H^*(A),0)$, whose cochain groups have dimensions $b_j$.

\eqref{tc-hom}  The equality 
$\TC_0(\RR_{i,s}(A))=\Supp_S(\gr^\m(\bwedge^s\B_i(A)))$, for every 
$s\ge1$, is Lemma~\ref{lem:tc-support-general} applied to the finitely 
generated $S$-module $\bwedge^s\B_i(A)$.

It remains to prove the containment in \eqref{eq:tcone-hom} for $s=1$.  
By Lemma~\ref{lem:stable-filt-supp}, applied to the filtration $F$ induced
on $\B_i(A)$ by the filtration $F^pK_\bullet(A)=A_\bullet\otimes_{\k}\m^p$
of the Koszul complex---which is $\m$-stable by the Artin--Rees
lemma---we have $\Supp_S(\gr^\m\B_i(A))=\Supp_S(\gr^F\B_i(A))$, so it
suffices to bound the support of $\gr^F\B_i(A)$.  Consider the spectral
sequence $(E^r_{p,q},d^r)$ of Theorem~\ref{thm:koszul-ss-hom} (in
unreindexed form, with $F^p$ decreasing and $p\ge0$).  Write 
$Z_r^p=\{x\in F^pK_i : \partial x\in F^{p+r}K_{i-1}\}$ and 
$Z^p_\infty=F^pK_i\cap\ker\partial$, so that
\[
\gr^F_p \B_i(A) = Z^p_\infty\,\big/\,
\bigl((F^p K_i\cap\im\partial)+Z^{p+1}_\infty\bigr).
\]
Fix $p\ge0$ and let $r\ge p+1$.  We claim that the inclusion 
$Z^p_\infty\subseteq Z^p_r$ induces an \emph{injective}\/ $S$-linear map
\[
\theta_r\colon \gr^F_p \B_i(A)\longrightarrow 
E^r_{p} = Z^p_r\,\big/\,
\bigl(\partial Z^{p-r+1}_{r-1}+Z^{p+1}_{r-1}\bigr).
\]
The map is well-defined: if $x\in F^pK_i\cap\im\partial$, write 
$x=\partial y$; since the filtration on $K_{i+1}$ is exhaustive with 
$F^0=K_{i+1}$ and $p-r+1\le 0$, we have 
$y\in F^{p-r+1}K_{i+1}$ and $\partial y\in F^p$, so 
$y\in Z^{p-r+1}_{r-1}$ and $x\in\partial Z^{p-r+1}_{r-1}$; and clearly 
$Z^{p+1}_\infty\subseteq Z^{p+1}_{r-1}$.  For injectivity, suppose 
$x\in Z^p_\infty$ lies in the denominator: $x=\partial y+z$ with 
$y\in Z^{p-r+1}_{r-1}$ and $z\in Z^{p+1}_{r-1}$.  Applying $\partial$ 
and using $\partial x=0$ gives $\partial z=0$, so 
$z\in Z^{p+1}_\infty$; and then $\partial y=x-z\in F^pK_i\cap\im\partial$.
Hence $x$ lies in $(F^pK_i\cap\im\partial)+Z^{p+1}_\infty$, as claimed.

Since each page $E^{r+1}$ is a subquotient of $E^r$ as an $S$-module, 
and passing to submodules and to subquotients can only shrink supports, 
we conclude
\[
\Supp_S\bigl(\gr^F_p\B_i(A)\bigr)
\subseteq \Supp_S\bigl(E^r_{p}\bigr)
\subseteq \Supp_S\bigl(E^2_{p}\bigr).
\]
By Theorem~\ref{thm:koszul-ss-hom}\eqref{ssh2}, the direct sum of the 
$E^2$-terms in total degree $i$ is the graded $S$-module 
$\B_i(H^*(A))$; the displayed containment, summed over $p$, therefore 
gives $\Supp_S(\gr^F\B_i(A))\subseteq\Supp_S(\B_i(H^*(A)))
=\RR_{i,1}(H^*(A))$, which together with the first paragraph completes 
the proof.
\end{proof}

\begin{remark}
\label{rem:tcone-i-finite}
The proof of \eqref{tc-coh} uses only that $\dim_\k A^j<\infty$ for 
$j\le i$: the target $A^{i+1}$ of $\delta^i_A$ may be 
infinite-dimensional, since the pivot block $D_i$ is a square matrix of 
size $\dim_\k C^i\le\dim_\k A^i$, and minors of a matrix with finitely 
many columns make sense regardless of the number of rows.  Likewise, the 
proof of \eqref{tc-hom} requires no finiteness on $A_{i+1}$: 
the boundary term is handled by exhaustiveness of the filtration alone.  
In particular, no skeleton reduction 
(Theorem~\ref{thm:cdga-skeleton}) is needed.
\end{remark}

\begin{problem}
\label{prob:tcone-hom-depth}
Does the containment 
$\TC_0\bigl(\RR_{i,s}(A)\bigr)\subseteq\RR_{i,s}\bigl(H^*(A)\bigr)$ hold 
for $s\ge2$?  Equivalently, by \eqref{tc-hom} and 
Lemma~\ref{lem:tc-support-general}, is 
$\Supp_S\bigl(\gr^\m(\bwedge^s\B_i(A))\bigr)$ contained in 
$\Supp_S\bigl(\bwedge^s\B_i(H^*(A))\bigr)$?  The argument above breaks 
down for $s\ge2$: the embedding $\gr^\m\B_i(A)\inj\B_i(H^*(A))$ realizes 
the associated graded module only as a \emph{submodule}, and supports of 
exterior powers are not monotone under submodules; e.g., for the submodule 
$\m\subset S=\k[x,y]$ one has $\bwedge^2\m\cong\k$, supported at the 
origin, while $\bwedge^2 S=0$.  We know of no counterexample; 
the containment does hold (as an equality) whenever $A$ is $q$-formal 
via a zig-zag of $q$-quasi-isomorphisms through $q$-finite $\cdgas$, by 
Theorem~\ref{thm:formal-tc} below.
\end{problem}

This establishes the inclusion part of 
Theorem~\ref{thm:tangent-cone-intro} from the Introduction;
the equality under $q$-formality is 
Theorem~\ref{thm:formal-tc} below.
The tangent cone inclusions reflect the fact that the non-linear part of 
the Koszul differential occurs in higher order at the origin, and 
therefore cannot create new tangent directions there.

\begin{corollary}
\label{cor:tcone-union}
In the setting of Theorem~\ref{thm:tangent-cone}, for every $q\ge 0$ such that 
$\dim_\k A^j<\infty$ for $j\le q$,
\[
\TC_0 \left(\bigcup_{i=0}^q \RR^i(A)\right) \subseteq 
   \bigcup_{i=0}^q \RR^i \bigl(H^{*}(A)\bigr)
   \quad\text{and}\quad
\TC_0 \left(\bigcup_{i=0}^q \RR_i(A)\right) \subseteq
   \bigcup_{i=0}^q \RR_i \bigl(H^{*}(A)\bigr).
\]
\end{corollary}

\begin{proof}
The tangent cone of a union is contained in the union of the tangent cones.
\end{proof}

\begin{remark}
\label{rem:tcone-history}
For the resonance varieties $\RR^i(A)$, the tangent cone inclusion was first
proved by Budur and Rubi\'{o}~\cite{BR} by deformation-theoretic means, resolving
a problem of Suciu~\cite[Prob.~3.2]{Su-indam}: they identify the analytic germ
of the resonance variety at the origin with a formal deformation functor of
$L_\infty$-pairs, and recover the tangent cone by linearizing the universal
$L_\infty$ differential.
The elimination mechanism used above, which pivots on the constant part of the
differential and reads off initial forms of minors, goes back to
Libgober~\cite{Li02}, and before him to Green and Lazarsfeld~\cite{GL87}.
In the setting of differential graded modules over a $\cdga$, both
inclusions were taken up in unfinished joint work of the author with
Denham~\cite[Thms.~6.2 and~6.3]{DS20}, where the block decomposition, the
vanishing of the relevant blocks by the anticommutation
\eqref{eq:anticomm}, the identification of the minors of the Aomoto matrix
as initial forms of minors of $\delta^i$, and the use of the $\m$-adic
spectral sequence for the support loci were all introduced.

The proofs given above are self-contained and differ from those drafted
there.  Working directly with $K^\bullet(A)$ over the local ring $S_\m$
replaces the reduction, by a Mumford-type lemma, to a bounded-above complex
of finitely generated free $S$-modules, and lowers the finiteness
requirement to $i$-finiteness (Remark~\ref{rem:tcone-i-finite}), where
\cite{BR} assumes finite type.  The elimination produces a Schur complement
with entries in $S_\m$ rather than in $S$, and the resulting initial forms
are compared by means of Lemma~\ref{lem:germ-basics}\eqref{gb1}.  For the
support loci, the relevant page of the spectral sequence is the second, whose
terms assemble to $\B_i(H^*(A))$; the embedding of $\gr^\m\B_i(A)$ into it,
together with Lemma~\ref{lem:stable-filt-supp}, is what yields
\eqref{tc-hom}.  The identification
$\TC_0\RR_{i,s}(A)=\Supp_S\gr^\m(\bwedge^s\B_i(A))$ holds in every depth,
while the inclusion itself remains open for $s\ge2$; see
Problem~\ref{prob:tcone-hom-depth}.
\end{remark}

%%%%%%%%%%%%%%%%%%%%%%%
\subsection{Germs at the origin and formality}
\label{subsec:germs-formality}
%%%%%%%%%%%%%%%%%%%%%%%

Under additional hypotheses---most notably, the presence of a 
$q$-quasi-isomorphism or a formality assumption---the inclusion of 
Theorem~\ref{thm:tangent-cone} upgrades to an equality.  In fact more is 
true: what a $q$-quasi-isomorphism preserves is not merely the tangent cone 
but the entire germ at the origin, in the sense of 
Definition~\ref{def:germ}.  The two theorems of this subsection make this 
precise.  We begin with an elementary comparison lemma for fibers of free 
complexes, which will let us transport jump loci across quasi-isomorphisms 
after completion, where the Koszul modules themselves become isomorphic by 
Theorem~\ref{thm:q-iso-koszul}\eqref{iso3}.

\begin{lemma}
\label{lem:fiber-comparison}
Let $R$ be a commutative Noetherian ring of finite global dimension and 
let $f\colon C\to D$ be a morphism of bounded-below complexes of flat 
$R$-modules such that $H_j(f)$ is an isomorphism for all $j\le q$.  Then 
for every $R$-module $M$ and every $i\le q$, the map 
$H_i(f\otimes\id_M)\colon H_i(C\otimes_{R}M)\to H_i(D\otimes_{R}M)$ is an 
isomorphism.  In particular, for every prime $P\subset R$ and $i\le q$,
\[
\dim_{\kappa(P)}H_i\bigl(C\otimes_{R}\kappa(P)\bigr)
=\dim_{\kappa(P)}H_i\bigl(D\otimes_{R}\kappa(P)\bigr).
\]
\end{lemma}

\begin{proof}
Since $C$ is a bounded-below complex of flat $R$-modules, there is a 
universal-coefficient spectral sequence
\[
E^2_{p,j}=\Tor_p^R\bigl(H_j(C),M\bigr)
\;\Longrightarrow\;
H_{p+j}(C\otimes_{R}M),
\]
natural in $C$ and in $M$, with differentials 
$d^r\colon E^r_{p,j}\to E^r_{p-r,\,j+r-1}$; see
\cite[Def.~5.2.1]{Weibel}.  
Since $R$ has finite global dimension, $E^2_{p,j}$ vanishes for $p$ 
large, so the filtration on the abutment is finite in each total 
degree.  The map $f$ induces a morphism of these spectral 
sequences which, on $E^2$-pages, is 
$\Tor_p(H_j(f),M)$---an isomorphism whenever $j\le q$.

We claim that every differential, on every page, that touches a term of 
total degree $\le q$ has both source and target with second index 
$j\le q$.  Indeed, a differential leaving $E^r_{p,j}$ with $p+j\le q$ 
lands in total degree $p+j-1$, with second index $j+r-1\le(p+j)-1\le q-1$ 
since $p\ge r$; and a differential arriving at total degree $i\le q$ 
departs from $E^r_{p+r,\,i-p-r+1}$, whose second index satisfies 
$i-p-r+1\le i-1\le q-1$ because $p\ge0$ and $r\ge2$.  (The terms 
$E^r_{0,j}$ with $j>q$, on which $f$ is not controlled, support no 
differentials at all, since their targets would have negative first 
index.)  By induction on $r$, the morphism of spectral sequences is 
therefore an isomorphism on every page in all total degrees $\le q$, 
hence on $E^\infty$ there; comparing the finite filtrations on the 
abutments via the five lemma gives the first assertion, and the second 
follows by taking $M=\kappa(P)$.
\end{proof}

\begin{theorem}
\label{thm:germ-qiso}
Let $\varphi\colon A\to B$ be a $q$-quasi-isomorphism of connected 
$q$-finite $\k$-$\cdgas$, and identify $H^1(A)$ with $H^1(B)$, and hence 
$S=\Sym(H_1(A))$ with $\Sym(H_1(B))$, by means of $\varphi$.  Then, for 
all $i\le q$ and $s\ge 1$:
\begin{enumerate}[itemsep=2pt]
\item \label{tcq1}
$\RR^{i,s}(H^*(A))=\RR^{i,s}(H^*(B))$.
\item \label{tcq2}
The resonance support loci have the same germ at the origin, as closed 
subschemes of $\Spec\widehat{S}$:
\[
\RR_{i,s}(A)_{(0)}=\RR_{i,s}(B)_{(0)} .
\]
\item \label{tcq3}
The resonance varieties have the same germ at the origin, as sets: the 
germs $\RR^{i,s}(A)_{(0)}$ and $\RR^{i,s}(B)_{(0)}$ have the same 
underlying reduced subscheme of $\Spec\widehat{S}$.
\item \label{tcq4}
$\TC_0\RR_{i,s}(A)=\TC_0\RR_{i,s}(B)$ and 
$\TC_0\RR^{i,s}(A)=\TC_0\RR^{i,s}(B)$.
\end{enumerate}
\end{theorem}

\begin{proof}
\eqref{tcq1} follows from Lemma~\ref{lem:resonance-truncation-dependence}
and the fact that $\varphi$ induces isomorphisms on $H^{\le q}$.

\eqref{tcq2} By Theorem~\ref{thm:q-iso-koszul}\eqref{iso3}, $\varphi$ 
induces $\widehat{S}$-linear isomorphisms 
$\widehat{\B_i(B)}\cong\widehat{\B_i(A)}$ for $i\le q$.  Since 
$\RR_{i,s}(-)=\Supp_S\bigl(\bwedge^s\B_i(-)\bigr)$ as schemes, 
Lemma~\ref{lem:germ-basics}\eqref{gb2} shows that the germ of this locus at 
the origin depends only on the $\widehat{S}$-module $\widehat{\B_i(-)}$; 
the two germs therefore coincide.

\eqref{tcq3} The chain map 
$K_\bullet(\varphi)\otimes_{S}\widehat{S}\colon 
K_\bullet(B)\otimes_{S}\widehat{S}\to K_\bullet(A)\otimes_{S}\widehat{S}$ is 
a morphism of bounded-below complexes of free---hence 
flat---$\widehat{S}$-modules over the regular ring $\widehat{S}$, and it 
induces isomorphisms on $H_j$ for $j\le q$: indeed, 
$H_j(K_\bullet(-)\otimes_{S}\widehat{S})\cong\B_j(-)\otimes_{S}\widehat{S}
=\widehat{\B_j(-)}$ by flatness of $\widehat{S}$, and these are 
identified by Theorem~\ref{thm:q-iso-koszul}\eqref{iso3}.  For any 
prime $P\subset\widehat{S}$, Lemma~\ref{lem:fiber-comparison} therefore 
gives
\[
\dim_{\kappa(P)}H_i\bigl(K_\bullet(A)\otimes_{S}\kappa(P)\bigr)
=\dim_{\kappa(P)}H_i\bigl(K_\bullet(B)\otimes_{S}\kappa(P)\bigr),
\qquad i\le q .
\]
By the determinantal description of 
Proposition~\ref{prop:res-descriptions}\eqref{rd1}, whose rank 
computation is valid fiberwise over any Noetherian ring, the locus 
$\{P\in\Spec\widehat{S} : \dim H_i(K_\bullet(-)\otimes_{S}\kappa(P))\ge s\}$ 
is the vanishing locus of the extended determinantal ideal 
$I^{(-)}\widehat{S}$, where $I^{(-)}$ denotes the ideal of 
Proposition~\ref{prop:res-descriptions}\eqref{rd1} for $A$, 
respectively $B$.  The displayed equality of fiberwise dimensions says 
that these two loci coincide as subsets of $\Spec\widehat{S}$, that is, 
$\sqrt{I^{(A)}\widehat{S}}=\sqrt{I^{(B)}\widehat{S}}$, which is 
precisely \eqref{tcq3}.

\eqref{tcq4} By Lemma~\ref{lem:germ-basics}\eqref{gb1}, the tangent cone at 
the origin of a closed subscheme is determined by its germ there, and indeed 
by the underlying set of that germ; so \eqref{tcq4} follows from 
\eqref{tcq2} and \eqref{tcq3}.
\end{proof}

Composing along a chain of such maps gives the following.

\begin{corollary}
\label{cor:germ-qequiv}
The conclusions of Theorem~\ref{thm:germ-qiso} hold whenever $A$ and $B$ 
are connected $q$-finite $\k$-$\cdgas$ joined by a zig-zag of 
$q$-quasi-isomorphisms all of whose terms are $q$-finite.
\end{corollary}

\begin{remark}
\label{rem:germ-qiso-scope}
For $\k=\C$, the germs in \eqref{tcq2} and \eqref{tcq3} may equivalently be
read as analytic germs at the origin, by
Lemma~\ref{lem:germ-basics}\eqref{gb3}.

Only the local structure at the origin is preserved.  Away from it no 
invariance holds, even for genuine quasi-isomorphisms: 
Example~\ref{ex:non-functorial-qiso} exhibits a quasi-isomorphism with 
$\RR_{1,1}(A)=\varnothing$ and $\RR_{1,1}(B)=\{1\}$.
\end{remark}

\begin{remark}
\label{rem:germ-qiso-dp}
A statement of this kind appears in unfinished joint work of the author with
Denham \cite[Prop.~5.4]{DS20}, for quasi-isomorphisms of differential graded
modules over a $\cdga$ with locally finite cohomology, resting there on the
quasi-isomorphism invariance of the completed Koszul complex
\cite[Thm.~4.10]{DS20}, which is proved by a mapping-cone reduction and a
comparison of complete filtered complexes rather than by the fiberwise
argument used above.  The treatment here is independent of that one, and
differs from it in three respects: the equality of germs is on the nose,
once $S$ is identified by means of $H_1(\varphi)$; the
hypotheses are $q$-truncated, so that only a $q$-quasi-isomorphism between
$q$-finite $\cdgas$ is required; and for the support loci the equality of
germs holds as schemes, by
Theorem~\ref{thm:germ-qiso}\eqref{tcq2}.

On the topological side, the germ at the origin is known to be a homotopy
invariant for an entirely different reason.  Let $X$ be a $q$-finite space
admitting a $q$-finite $q$-model $A$ over $\C$.  Dimca and
Papadima \cite[Thm.~B]{DP-ccm} construct a local analytic isomorphism
between the germ at $1$ of the character variety of $\pi_1(X)$ and the germ
at $0$ of $H^1(A)$, carrying $\VV^{i,s}(X;\C)_{(1)}$ onto
$\RR^{i,s}(A)_{(0)}$ for all $i\le q$ and $s\ge 1$; in the case $q=i=1$ this
is due to Dimca, Papadima and Suciu \cite[Thm.~A]{DPS-duke}.  Since the germ
on the topological side is a homotopy invariant of $X$, model-independence
follows for those $\cdgas$ which arise as models of spaces.  What the
algebraic route contributes is that neither a space, nor analytic geometry,
nor deformation theory is needed, and that the conclusion holds over any
algebraically closed field of characteristic $0$.  What it does not recover
is the comparison with the characteristic varieties themselves: the
exponential identification of $\VV^{i,s}(X;\C)_{(1)}$ with
$\RR^{i,s}(A)_{(0)}$ lies outside the scope of these methods.
\end{remark}

\begin{theorem}
\label{thm:formal-tc}
Let $(A,d)$ be a connected, $q$-finite $\k$-$\cdga$ which is $q$-formal 
via a zig-zag of $q$-quasi-isomorphisms all of whose terms are 
$q$-finite.  Then, for all $i\le q$ and $s\ge 1$, the resonance loci of $A$ 
have the same germs at the origin as those of its cohomology algebra,
\[
\RR^{i,s}(A)_{(0)}=\RR^{i,s}\bigl(H^{*}(A)\bigr)_{(0)}
\qquad\text{and}\qquad
\RR_{i,s}(A)_{(0)}=\RR_{i,s}\bigl(H^{*}(A)\bigr)_{(0)},
\]
the first as sets and the second as schemes.  Since the loci on the right 
are cones, it follows that
\[
\TC_0\RR^{i,s}(A)= \RR^{i,s}\bigl(H^{*}(A)\bigr) 
\qquad\text{and}\qquad
\TC_0\RR_{i,s}(A)= \RR_{i,s}\bigl(H^{*}(A)\bigr).
\]
\end{theorem}

\begin{proof}
The zig-zag joins $A$ to $(H^*(A),0)$ through $q$-finite $\cdgas$, so 
Corollary~\ref{cor:germ-qequiv} yields the two germ equalities.  It remains 
to see that the loci of $(H^*(A),0)$ are cones, hence equal to their own 
tangent cones.  Since the differential is zero, the Koszul module 
$\B_i(H^*(A))$ is a graded $S$-module, hence so is 
$\bwedge^s\B_i(H^*(A))$, and Corollary~\ref{cor:tc-support-graded} gives 
$\TC_0\RR_{i,s}(H^*(A))=\RR_{i,s}(H^*(A))$.  For the jump loci, the 
Koszul coboundaries of $(H^*(A),0)$ have linear entries, so the 
determinantal ideal of 
Proposition~\ref{prop:res-descriptions}\eqref{rd1} is homogeneous and 
$\TC_0\RR^{i,s}(H^*(A))=\RR^{i,s}(H^*(A))$ as well.  Combining these two 
identifications with Lemma~\ref{lem:germ-basics}\eqref{gb1} gives the 
displayed equalities of tangent cones.
\end{proof}

\begin{remark}
\label{rem:formal-tc-hypotheses}
The finiteness hypothesis on the terms of the zig-zag is needed so that 
Theorem~\ref{thm:q-iso-koszul}\eqref{iso3} applies to each arrow; the 
definition of $q$-formality \textup{(}Section~\ref{subsec:cdga-qiso}\textup{)} 
imposes no finiteness on the intermediate $\cdgas$, and we do not know 
whether every $q$-formal, $q$-finite $\cdga$ admits a $q$-formality 
zig-zag through $q$-finite $\cdgas$.  The hypothesis is satisfied in the 
cases arising in this monograph: whenever $A$ maps to, or receives, a 
$q$-quasi-isomorphism from $(H^*(A),0)$ directly, as for $\CE(\sol_2)$ via 
the inclusion of $\bigwedge(a)$; and whenever the zig-zag can be 
taken through finite-type models.  Note also that no 
separation hypotheses appear: the proof runs entirely through 
$\m$-adic completions, where 
Theorem~\ref{thm:q-iso-koszul}\eqref{iso3} holds unconditionally.  This 
is in contrast to statements at the level of the modules $\B_i(A)$ 
themselves, or of the varieties $\RR^{i,s}(A)$: those are not invariant 
under quasi-isomorphism, and $A=\CE(\sol_2)$---which satisfies the 
hypotheses above and is formal---has $\RR^1(A)=\{0,1\}\supsetneq\{0\}=
\RR^1(H^*(A))$, so that only the tangent cones at the origin, not the 
varieties, can agree in general.
\end{remark}

Together with Theorem~\ref{thm:tangent-cone}, this completes 
the proof of Theorem~\ref{thm:tangent-cone-intro}.

\begin{conjecture}
\label{conj:formal-tc}
In the setting of Theorem~\ref{thm:formal-tc}, each irreducible 
component of $\RR^{i,s}(A)$ and $\RR_{i,s}(A)$ is a $\k$-linear 
subspace of $H^1(A)$ defined over $\Q$.
\end{conjecture}

Theorem~\ref{thm:formal-tc} and Conjecture~\ref{conj:formal-tc} are 
the algebraic counterparts of several deep results on resonance 
varieties of spaces and groups. A detailed comparison with the 
topological literature, and further motivation for the conjecture, 
will be given in Sections~\ref{subsec:res-spaces} 
and~\ref{subsec:res-zero-model}.

%%%%%%%%%%%%%%%%%%%%%%%
\subsection{Collapse, vanishing, and finite-length Koszul modules}
\label{subsec:tcone-collapse}
%%%%%%%%%%%%%%%%%%%%%%%

Theorem~\ref{thm:tangent-cone} gives a one-sided inclusion 
$\TC_0(\RR^{i,s}(A))\subseteq\RR^{i,s}(H^*(A))$ valid for all 
finite-type $\cdgas$, and in general this inclusion may be strict: 
the tangent cone of $\RR^i(A)$ at the origin may be a proper subset 
of $\RR^i(H^*(A))$. When the cohomological Koszul spectral sequence 
collapses early, however, the inclusion becomes an equality and one 
can recover the tangent cone exactly when the resonance of $H^*(A)$ 
vanishes; early collapse also forces finite-length behavior on the Koszul module.

\begin{proposition}
\label{prop:no-hidden-res}
Let $(A,d)$ be a connected $\cdga$ over a field $\k$ of characteristic $0$, 
with $\dim_\k A^j<\infty$ for $j\le i+1$.  
Assume that the cohomological Koszul spectral sequence of 
Theorem~\ref{thm:koszul-ss-coh} collapses at the $E_2$-page 
in total degree~$i$ (i.e., $E_2^{p,q}=E_\infty^{p,q}$ for $p+q=i$).
If the resonance support locus satisfies
$\RR^i(H^\ast(A))\subseteq\{0\}$, then
\[
\TC_0\bigl(\RR^i(A)\bigr)\subseteq\{0\},
\]
the origin is an isolated point of $\RR^i(A)$ \textup{(}or does not lie on
it\textup{)}, and the localization $\B^i(A)_\m$ has finite length over
$S_\m$.  Equivalently, $\gr^\m\B^i(A)$ has finite length over
$S=\Sym(H_1(A))$; and likewise $\gr^F\B^i(A)$, for the equivalent induced
filtration of Theorem~\ref{thm:koszul-ss-coh}\eqref{ssc4}.
\end{proposition}

\begin{proof}
Set $S=\Sym(H_1(A))$.  By Theorem~\ref{thm:koszul-ss-coh}\eqref{ssc2}, the
$E_1$-page of the cohomological Koszul spectral sequence is
$K^\bullet(H^*(A))$, so $E_2^{p,i-p}\cong\B^i(H^*(A))_p$; and by
\eqref{ssc4}, $E_\infty^{p,i-p}\cong\gr^F_p\B^i(A)$, for the induced
filtration $F$.  The collapse hypothesis therefore gives an isomorphism
of graded $S$-modules
\begin{equation}
\label{eq:gr-collapse}
\gr^F\B^i(A)\cong\B^i\bigl(H^*(A)\bigr).
\end{equation}
Now $\RR^i(H^*(A))\subseteq\{0\}$ means that
$\Supp_S\B^i(H^*(A))\subseteq\{0\}$, that is, $\B^i(H^*(A))$ has finite
length; by \eqref{eq:gr-collapse} so does $\gr^F\B^i(A)$, and by
Lemma~\ref{lem:stable-filt-supp}---the filtration $F$ being equivalent to
the $\m$-adic one---so does $\gr^\m\B^i(A)$, whence
\[
\TC_0\bigl(\RR^{i}(A)\bigr)
=\Supp_S\bigl(\gr^\m\B^i(A)\bigr)\subseteq\{0\},
\]
the first equality being Lemma~\ref{lem:tc-support-general}.  Since
$\gr^\m$ commutes with localization at $\m$, the module
$\gr^\m\bigl(\B^i(A)_\m\bigr)$ has finite length as well, and hence so
does $\B^i(A)_\m$: the $\m$-adic filtration on a finitely generated
module over the local ring $S_\m$ is separated
(Proposition~\ref{prop:separation-criterion}), so a filtration with
finitely many nonzero graded pieces, each of finite length, exhausts the
module.  Finally, $\TC_0(\Supp_SM)\subseteq\{0\}$ says precisely that the
origin is isolated in, or absent from, $\Supp_SM$.
\end{proof}

\begin{remark}
\label{rem:no-hidden-res-global}
The conclusion cannot be strengthened to $\Supp_S\B^i(A)\subseteq\{0\}$,
nor to finite length of $\B^i(A)$ itself when $\dim_\k H^1(A)\ge2$: the
associated graded module sees only what happens at the origin.  For
$A=\CE(\sol_2)$ the hypotheses hold in degree $i=2$---one has
$\B^2(H^*(A))=0$, so $\RR^2(H^*(A))=\varnothing$, and the spectral
sequence collapses---while $\B^2(A)\cong\k[x]/(x-1)$ is supported at the
point $1$, away from the origin
\textup{(}Example~\ref{ex:sol2-koszul}\textup{)}.  Only because
$\dim_\k H^1(A)=1$ is that module of finite length; in general a component
of $\Supp_S\B^i(A)$ missing the origin can have positive dimension, and is
then invisible both to $\gr^\m$ and to the localization at $\m$.
\end{remark}

%%%%%%%%%%%%%%%%%%%%%
\subsection{Positive weights and resonance}
\label{subsec:res-weights}
%%%%%%%%%%%%%%%%%%%%%

We now combine the tangent cone theorem with the conicality of resonance 
varieties arising from positive weights.

\begin{theorem}
\label{thm:resonance-inclusion}
Let $(A,d)$ be a connected, finite-type $\cdga$ over a field of characteristic $0$, 
admitting positive weights. Then, for all $i,s$,
\begin{enumerate}[itemsep=2pt]
\item \label{rr1}
The resonance varieties $\RR^{i,s}(A)$ and the support loci 
$\RR_{i,s}(A)$ are invariant under the weighted $\k^*$-action of 
Corollary~\ref{cor:conical}; if $H^1(A)$ is concentrated in a single 
weight, they are conical.

\item \label{rr2}
If $H^1(A)$ is concentrated in a single weight, then
$\RR^{i,s}(A) \subseteq \RR^{i,s}(H^{*}(A))$ for all $s\ge1$, and 
$\RR_{i,1}(A) \subseteq \RR_{i,1}(H^{*}(A))$.
\end{enumerate}
\end{theorem}

\begin{proof}
For \eqref{rr1}, the support loci
$\RR_{i,s}(A)=\Supp_S\bigl(\bwedge^s\B_i(A)\bigr)$ are covered by
Corollary~\ref{cor:conical}.  The same $\k^*$-equivariance of $\B_i(A)$
renders the determinantal ideals of
Proposition~\ref{prop:res-descriptions}\eqref{rd1} weighted-homogeneous, so
the jump loci $\RR^{i,s}(A)$ are invariant as well; and when $H^1(A)$ is
concentrated in a single weight, weighted-homogeneous ideals are homogeneous
and all these loci are cones.

By Theorem~\ref{thm:tangent-cone}, 
$\TC_0(\RR^{i,s}(A)) \subseteq \RR^{i,s}(H^{*}(A))$ for all $s\ge1$, 
and $\TC_0(\RR_{i,1}(A)) \subseteq \RR_{i,1}(H^{*}(A))$.
Under the hypothesis of \eqref{rr2}, the loci are conical by 
\eqref{rr1}, so each coincides with its tangent cone at the origin, 
i.e., $\RR^{i,s}(A)=\TC_0(\RR^{i,s}(A))$, and similarly for 
$\RR_{i,s}(A)$; claim \eqref{rr2} follows.
\end{proof}

The single-weight hypothesis in \eqref{rr2} cannot be weakened to weighted
invariance alone: a weighted-homogeneous variety need not equal its tangent
cone at the origin, as the curve $\bV(y-x^2)$ of
Remark~\ref{rem:weighted-not-conical} shows.

\begin{remark}
\label{rem:resonance-inclusion-depth}
The same argument gives the homological inclusion 
$\RR_{i,s}(A)\subseteq\RR_{i,s}(H^{*}(A))$ in every depth $s\ge2$ for which 
Problem~\ref{prob:tcone-hom-depth} has a positive answer for $A$, since 
under the hypothesis of \eqref{rr2} the locus $\RR_{i,s}(A)$ again coincides 
with its own tangent cone.
\end{remark}

Thus, failure of weighted invariance---or, when $H^1$ is concentrated in a 
single weight, of conicality or of the resonance inclusion---provides an 
obstruction to the existence of positive weights on a $\cdga$; see 
Example~\ref{ex:sol2} for a concrete way in which this criterion applies.

Previous results on linearity of resonance varieties under positive weights 
were obtained in a topological setting, assuming the existence of a finite-type 
CW-complex realizing the model and a comparison map preserving rational 
structures \cite{DP-ccm, MPPS}.
Theorem~\ref{thm:resonance-inclusion} shows that weighted invariance of 
resonance---and its conicality, in the single-weight case---is an intrinsic 
property of $\cdgas$ admitting positive weights, independent of any 
topological realization of the $\cdga$.

\begin{proposition}
\label{prop:weight-discrepancy}
Let $(A,d)$ be a connected, finite-type $\cdga$ admitting positive weights, 
with $H^1(A)$ concentrated in a single weight.
Then for all $i,s$, the difference
\[
\RR^{i,s}(A)\setminus \RR_{i,s}(A)
\]
is a conical subset of $H^1(A)$.
Moreover, if $a\in H^1(A)$ lies outside $\bigcup_{r<i}\RR_{r,1}(A)$, then the
spectral sequence of Theorem~\ref{thm:jump-support-containment}
degenerates at $E^2$ for all $\lambda a$ with $\lambda\in \k^*$.
Consequently, $\P\bigl(\RR^{i,s}(A)\bigr)$ and $\P\bigl(\RR_{i,s}(A)\bigr)$
agree outside $\P\bigl(\bigcup_{r<i}\RR_{r,1}(A)\bigr)$; in particular,
\[
\P\bigl(\RR^{1,s}(A)\bigr) = \P\bigl(\RR_{1,s}(A)\bigr), 
\]
as projective subvarieties of $\P(H^1(A))$.
\end{proposition}

\begin{proof}
All constructions (Koszul complexes, Koszul modules, and change-of-rings
spectral sequence) are compatible with the positive grading on $A$, hence
inherit a natural $\k^*$-action induced by scaling on $H^1(A)$.
Since $\RR^{i,s}(A)$ and $\RR_{i,s}(A)$ are $\k^*$-invariant by
Theorem~\ref{thm:resonance-inclusion}\eqref{rr1}---the action being the 
standard one, as $H^1(A)$ is concentrated in a single weight---their 
difference is a conical subset as well.

For the second statement, note that the conditions
$\B_r(A)\otimes_{S} \kappa_a=0$, $r<i$, are preserved under scaling 
$a\mapsto \lambda a$, which corresponds to the induced $\k^*$-action 
on $S$, under which $\m_a$ is sent to $\m_{\lambda a}$, and each $\B_r(A)$ 
is equivariant as a graded $S$-module. Hence the degeneration criterion in 
Proposition~\ref{prop:generic-resonance-agreement} is invariant along 
$\k^*$-orbits, and the identification of jump and support loci holds projectively.
\end{proof}

Without the single-weight hypothesis, the same conclusions hold with
$\P(H^1(A))$ replaced by the weighted projective space determined by the
weights of $H^1(A)$.

In the presence of positive weights, resonance and support loci are naturally 
$\k^*$-invariant and therefore descend to well-defined projective subvarieties 
of $\P(H^1(A))$. From this perspective, the comparison between jump and 
support loci is a statement about projective cones in $\P(H^1(A))$, 
governed by the degeneration behavior of the Koszul spectral sequence 
along $\k^*$-orbits. In particular, nontrivial resonance corresponds to the 
existence of nonempty projective subvarieties, while vanishing resonance 
amounts to the collapse of these projective cones to the empty projective 
locus.

%%%%%%%%%%%%%%%%%%%%%
\subsection{Examples and discussion}
\label{subsec:tcone-discuss}
%%%%%%%%%%%%%%%%%%%%%

We conclude this section with some examples and further discussion of the 
tangent cone theorem for $\cdgas$. 

\begin{example}
\label{ex:sol2}
Let $A=\CE(\sol_2)$ be the $\cdga$ of Example~\ref{ex:sol2-koszul},
with $H^{*}(A)\cong\bigwedge(a)$, $S=\k[x]$, and
$\B_1(A)\cong\k[x]/(x-1)$, $\B_1\bigl(H^{*}(A)\bigr)=0$.
Since the relation $x=1$ is non-homogeneous, the $S$-module $\B_1(A)$
carries no grading compatible with the standard grading on $S$.
Consequently, the associated resonance varieties
\[
\RR^1(A)=\{0,1\} \quad\text{and}\quad
\RR_1(A)=\{1\} \subseteq \Spec\k[x]=\mathbb{A}^1
\]
are not homogeneous. In particular,
\[
\RR^1(A)\not\subseteq \RR^1\bigl(H^{*}(A)\bigr)=\{0\},
\qquad
\RR_1(A)\not\subseteq \RR_1\bigl(H^{*}(A)\bigr)=\emptyset.
\]

Passing to tangent cones at the origin yields
\[
\TC_0\RR^1(A)=\{0\}, \qquad \TC_0\RR_1(A)=\emptyset,
\]
since $\RR_1(A)$ does not pass through $0$.
Thus both homological and support resonance satisfy the inclusions of
Theorem~\ref{thm:tangent-cone}, with equality in degree~$1$
by Theorem~\ref{thm:formal-tc}: $A$ is formal via the one-arrow zig-zag 
$(H^*(A),0)=\bigwedge(a)\inj A$ 
(Remark~\ref{rem:ce-sol2}), whose terms are of finite type, so the 
hypotheses of that theorem hold.

In particular, since $\RR^1(A)=\{0,1\}$ is not homogeneous, 
$A$ admits no positive-weight decomposition: as $\dim_\k H^1(A)=1$, any 
positive-weight decomposition would concentrate $H^1(A)$ in a single 
weight, so Theorem~\ref{thm:resonance-inclusion}\eqref{rr1} would force 
$\RR^1(A)$ to be conical, contradicting the presence of the point $\{1\}$.
Thus, even when resonance varieties fail to be homogeneous 
or invariant under quasi-isomorphisms, their tangent cones 
at the origin are controlled solely by the cohomology algebra,
as guaranteed by Theorem~\ref{thm:tangent-cone}.
\end{example}

\begin{example}
\label{ex:heisenberg-resonance}
Let $A=\CE(\h(1))=\bwedge(e_1,e_2,e_3)$ with
$de_3=-e_1\wedge e_2$ and $de_1=de_2=0$,
the Chevalley--Eilenberg $\cdga$ of the $3$-dimensional Heisenberg Lie algebra.
Then $H^1(A)=\langle e_1,e_2\rangle\cong \k^2$.
A direct computation (Section~\ref{subsec:heisenberg-1}) shows that
$\B_1(A)\cong \k$, so that $\RR_{1,1}(A)=\{0\}$.  By
Theorem~\ref{thm:jump-support-containment}\eqref{jsc2}, the jump and support
loci agree away from the origin, so the first resonance variety of $A$ is
trivial as well:
\[
\RR^1(A)=\{0\}\subseteq H^1(A).
\]

By contrast, the cup product on the cohomology algebra $H^*(A)$ vanishes in
degree~$1$, and hence
\[
\RR^1\bigl(H^*(A)\bigr)=H^1(A)\cong \k^2.
\]

Thus, the tangent cone theorem yields a strict inclusion,
\[
\TC_0\RR^1(A)=\{0\} \subsetneqq \RR^1\bigl(H^*(A)\bigr).
\]

This strict containment reflects the presence of nontrivial triple Massey
products in $A$, which suppress resonance that would otherwise be predicted by
the degree-$1$ cup product alone.  The example is the one used for this
purpose in \cite[Ex.~7.3]{DS20}.

Since $A=\CE(\h(1))$ is not formal, this strict containment 
shows that the formality hypothesis in 
Theorem~\ref{thm:formal-tc} cannot be dropped: equality of 
tangent cones and cohomology-level resonance may fail for 
non-formal $\cdgas$, even when the model admits positive 
weights.
\end{example}

%%%%%%%%%%%%%%%%%%%%%%%%%%%%%%%%%%%%%%
\section{Resonance under algebraic constructions}
\label{sect:res-constructions}
%%%%%%%%%%%%%%%%%%%%%%%%%%%%%%%%%%%%%%

Many $\cdga$ models arising in geometry and topology are obtained from simpler ones
by means of tensor products, coproducts, and Hirsch extensions.  
In the next three subsections, we record several results describing 
how resonance varieties behave under these constructions.  
These structural results will be used repeatedly in later sections, 
and they also provide an alternative approach to the vanishing phenomena 
for resonance observed in nilpotent Chevalley--Eilenberg models.

%%%%%%%%%%%%%%%%%%%%%%%
\subsection{Resonance of tensor products and coproducts}
\label{subsec:tensor-coprod}
%%%%%%%%%%%%%%%%%%%%%%%

The Künneth decompositions of Corollary~\ref{cor:koszul-kunneth-modules} 
and Proposition~\ref{prop:koszul-coprod} translate into formulas for the 
resonance varieties of tensor products and coproducts of $\cdgas$. For 
tensor products the answer is complete, in all degrees and depths, but it is 
not a union of products of the resonance varieties of the factors: the 
relevant fiber dimensions multiply across the two factors, whereas depth 
conditions add. Earlier results of this type, obtained via Aomoto 
complex methods and valid only at depth $s=1$ or in degree $q=1$, 
appear in \cite{PS-plms10, SW-mz, Su-edinb, PS-springer, Su-bockres}.

\begin{proposition}
\label{prop:prod-res}
Let $(A,d_A)$ and $(B,d_B)$ be two connected, finite-type $\cdgas$, let
$C=A\otimes_{\k} B$, and identify $H^1(C)=H^1(A)\oplus H^1(B)$.  For 
$a\in H^1(A)$ and $b\in H^1(B)$, put
\[
d_i(a)=\dim_\k H^i(A,\delta_a),
\qquad
\mu_i(a)=\dim_\k\bigl(\B_i(A)\otimes_{S_A}\kappa_a\bigr),
\]
and define $e_j(b)$ and $\nu_j(b)$ analogously for $B$.  Then for all 
$q\ge 0$ and $s\ge 1$:
\begin{enumerate}[itemsep=3pt]
\item\label{t1}
$\RR^{q,s}(C)=\bigl\{(a,b) \;:\; 
\sum_{i+j=q} d_i(a)\,e_j(b)\ge s\bigr\}$.
\item\label{t2}
$\RR_{q,s}(C)=\bigl\{(a,b) \;:\; 
\sum_{i+j=q} \mu_i(a)\,\nu_j(b)\ge s\bigr\}$.
\end{enumerate}
In depth $s=1$ these formulas become
\[
\RR^{q,1}(C)=\bigcup_{i+j=q}\RR^{i,1}(A)\times\RR^{j,1}(B),
\qquad
\RR_{q,1}(C)=\bigcup_{i+j=q}\RR_{i,1}(A)\times\RR_{j,1}(B).
\]
\end{proposition}

\begin{proof}
\eqref{t1} By Proposition~\ref{prop:res-descriptions}\eqref{rd3}, membership 
in $\RR^{q,s}(C)$ is detected by the Aomoto complex of $C$ at the point 
$c=(a,b)$.  Specializing the identification 
$K^\bullet(C)=K^\bullet(A)\otimes_{\k} K^\bullet(B)$ of 
Theorem~\ref{thm:koszul-kunneth} at $c$, and using the decomposition 
\eqref{eq:dcdadb} of $\delta_C$, gives an isomorphism of cochain complexes
\[
(C,\delta_c)\cong (A,\delta_a)\otimes_{\k} (B,\delta_b) .
\]
Since $\k$ is a field, the Künneth theorem yields 
$H^q(C,\delta_c)\cong\boplus_{i+j=q}H^i(A,\delta_a)\otimes_{\k} 
H^j(B,\delta_b)$, and hence 
$\dim_\k H^q(C,\delta_c)=\sum_{i+j=q}d_i(a)e_j(b)$.

\eqref{t2} By Lemma~\ref{lem:supp-mu}, 
$\RR_{q,s}(C)=\{c : \dim_\k \B_q(C)\otimes_{S_C}\kappa_c\ge s\}$.  By 
Corollary~\ref{cor:koszul-kunneth-modules}\eqref{kk2}, 
$\B_q(C)\cong\boplus_{i+j=q}\B_i(A)\otimes_{\k}\B_j(B)$, where each summand 
carries the external $S_C$-module structure coming from 
$S_C\cong S_A\otimes_{\k} S_B$.  Since 
$\m_c=\m_a\otimes_{\k} S_B+S_A\otimes_{\k}\m_b$, we get
\[
\bigl(\B_i(A)\otimes_{\k}\B_j(B)\bigr)\otimes_{S_C}\kappa_c \cong 
\bigl(\B_i(A)\otimes_{S_A}\kappa_a\bigr)\otimes_{\k}
\bigl(\B_j(B)\otimes_{S_B}\kappa_b\bigr),
\]
and therefore 
$\dim_\k\B_q(C)\otimes_{S_C}\kappa_c=\sum_{i+j=q}\mu_i(a)\nu_j(b)$.

Finally, a sum of products of non-negative integers is positive if and only 
if one of those products is; and $d_i(a)\ge1$ means precisely that 
$a\in\RR^{i,1}(A)$, while $\mu_i(a)\ge1$ means precisely that 
$a\in\RR_{i,1}(A)$, by Proposition~\ref{prop:res-descriptions}\eqref{rd3} 
and Lemma~\ref{lem:supp-mu}.  This gives the two depth-$1$ formulas.
\end{proof}

\begin{remark}
\label{rem:higher-depth-products}
In depths $s\ge2$ the loci $\RR^{q,s}(C)$ and $\RR_{q,s}(C)$ are \emph{not} 
unions of products $\RR^{i,r}(A)\times\RR^{j,t}(B)$ with $r+t=s$.  Two 
distinct phenomena obstruct such a formula.

First, exterior powers do not distribute over external tensor products.  For 
$S_C$-modules of the form $M\otimes_{\k} N$ as above, $\bwedge^s_{S_C}$ is not 
computed from the exterior powers of $M$ over $S_A$ and of $N$ over $S_B$: 
already for free modules of ranks $m$ and $n$, the module 
$\bwedge^s_{S_C}(M\otimes_{\k} N)$ has rank $\binom{mn}{s}$, whereas 
$\boplus_{r+t=s}\bwedge^r_{S_A}M\otimes_{\k}\bwedge^t_{S_B}N$ has rank 
$\binom{m+n}{s}$.  This is an instance of the failure of exterior powers to 
commute with restriction of scalars.

Second, and more fundamentally, the fiber dimensions multiply across the two 
factors, while depth conditions add.  The inequality 
$\sum_{i+j=q}d_i(a)e_j(b)\ge s$ of part~\eqref{t1} is therefore not a 
Boolean combination of the conditions $d_i(a)\ge r$ and $e_j(b)\ge t$ with 
$r+t=s$.  Two examples make this concrete.
\begin{enumerate}[itemsep=3pt]
\item
Let $A=B=H^{*}(S^1;\k)$, so that $C\cong H^{*}(T^2;\k)$.  Then 
$d_0(a)=d_1(a)=1$ for $a=0$ and $d_i(a)=0$ otherwise, so part~\eqref{t1} 
gives $\RR^{1,1}(C)=\{0\}$, in agreement with the classical computation 
$\RR^{1,1}(H^*(T^2))=\{0\}$.  By contrast, the union 
$\bigcup_{i+j=1,\,r+t=1}\RR^{i,r}(A)\times\RR^{j,t}(B)$ is the union of the 
two coordinate axes if the extreme loci $\RR^{i,0}$ are read as all of 
$H^1$, and is empty if the indices $r,t$ are required to be positive.
\item
Let $A=B=H^{*}(S^1\vee S^1;\k)$ and take $q=s=2$.  Here $d_1(a)=2$ for $a=0$ 
and $d_1(a)=1$ otherwise, while $d_i=0$ for $i\ge2$; hence 
$\sum_{i+j=2}d_i(a)e_j(b)=d_1(a)e_1(b)$ and part~\eqref{t1} gives 
$\RR^{2,2}(C)=\bigl(\{0\}\times H^1(B)\bigr)\cup
\bigl(H^1(A)\times\{0\}\bigr)$.  The corresponding union contains 
$\RR^{1,1}(A)\times\RR^{1,1}(B)=H^1(C)$, hence is all of $H^1(C)$.
\end{enumerate}
\end{remark}

For coproducts, the answer is additive rather than multiplicative, since the
Aomoto complex of $A\vee B$ splits in positive degrees.  Here a genuine
exception occurs in degree $q=1$, but it is not a uniform shift of the depth:
the correction term is supported away from the two coordinate subspaces.
Throughout, we use the convention that $\RR^{q,0}(A)=\RR_{q,0}(A)=H^1(A)$.

\begin{proposition}
\label{prop:rescoprod}
Let $C=A\vee B$ be the coproduct of two connected, finite-type 
$\cdgas$ with $b_1(A)>0$ and $b_1(B)>0$, and identify 
$H^1(C)=H^1(A)\times H^1(B)$.  With $d_q(a)$, $e_q(b)$, $\mu_q(a)$, and 
$\nu_q(b)$ as in Proposition~\ref{prop:prod-res}, the following hold.
\begin{enumerate}[itemsep=3pt]
\item\label{rc1}
For $q\ge2$ and $s\ge1$,
\[
\RR^{q,s}(C)=\bigl\{(a,b) : d_q(a)+e_q(b)\ge s\bigr\}
=\bigcup_{j+k=s} \RR^{q,j}(A)\times\RR^{q,k}(B) ,
\]
while for $q=1$,
\[
\RR^{1,s}(C)=\bigl\{(a,b) : d_1(a)+e_1(b)+\chi(a,b)\ge s\bigr\},
\qquad
\chi(a,b)=\begin{cases} 1 & a\ne0 \text{ and } b\ne0,\\ 0 &\text{otherwise.}
\end{cases}
\]
\item\label{rc2}
For $q\ge2$ and $s\ge1$,
\[
\RR_{q,s}(C)=\bigl\{(a,b) : \mu_q(a)+\nu_q(b)\ge s\bigr\}
=\bigcup_{j+k=s} \RR_{q,j}(A)\times\RR_{q,k}(B) ,
\]
while for $q=1$,
\[
\mu_1(a)+\nu_1(b)\;\le\;
\dim_\k\bigl(\B_1(C)\otimes_{S_C}\kappa_{(a,b)}\bigr)
\;\le\;\mu_1(a)+\nu_1(b)+\mu_{(a,b)}(\m_A\m_B).
\]
In particular, $\RR_{1,1}(C)=H^1(C)$.
\end{enumerate}
\end{proposition}

\begin{proof}
\eqref{rc1} Since $C^0=\k$ and $C^q=A^q\oplus B^q$ for $q\ge1$, and since 
$A^{\ge1}\cdot B^{\ge1}=0$ in $C$, the Aomoto complex of $C$ at a point 
$c=(a,b)$ is
\[
\begin{tikzcd}[column sep=22pt]
\k \arrow[r, "{(a,b)}"] & \bigl(A^{\ge1},\delta_a\bigr)\oplus
\bigl(B^{\ge1},\delta_b\bigr) .
\end{tikzcd}
\]
For $q\ge2$ this gives $H^q(C,\delta_c)\cong H^q(A,\delta_a)\oplus 
H^q(B,\delta_b)$, whence the first formula, by 
Proposition~\ref{prop:res-descriptions}\eqref{rd3}; the union description 
follows since $d_q(a)+e_q(b)\ge s$ if and only if $d_q(a)\ge j$ and 
$e_q(b)\ge k$ for some $j+k=s$ with $j,k\ge0$.

For $q=1$, write $\epsilon(x)=1$ if $x\ne0$ and $\epsilon(x)=0$ otherwise.  
Then
\[
\dim_\k H^1(C,\delta_c)
=\Bigl(\dim_\k\ker\delta^1_a+\dim_\k\ker\delta^1_b\Bigr)-\epsilon(a,b)
=d_1(a)+\epsilon(a)+e_1(b)+\epsilon(b)-\epsilon(a,b),
\]
since $\dim\ker\delta^1_a=d_1(a)+\dim\im\delta^0_a=d_1(a)+\epsilon(a)$, and 
similarly for $b$.  As $\epsilon(a)+\epsilon(b)-\epsilon(a,b)=\chi(a,b)$, the 
second formula follows.

\eqref{rc2} For $q\ge2$, Proposition~\ref{prop:koszul-coprod}\eqref{kc1} gives 
$\B_q(C)\cong\bigl(\B_q(A)\otimes_{S_A}S_C\bigr)\oplus
\bigl(\B_q(B)\otimes_{S_B}S_C\bigr)$, and 
$\bigl(\B_q(A)\otimes_{S_A}S_C\bigr)\otimes_{S_C}\kappa_{(a,b)}\cong
\B_q(A)\otimes_{S_A}\kappa_a$; hence the fiber dimensions add, and 
Lemma~\ref{lem:supp-mu}\eqref{sm2} gives the formula, the union description 
following as before.

For $q=1$, apply $-\otimes_{S_C}\kappa_{(a,b)}$ to the short exact sequence of 
Proposition~\ref{prop:koszul-coprod}\eqref{kc2}: right exactness gives the 
upper bound, while the surjection onto $\m_A\m_B\otimes_{S_C}\kappa_{(a,b)}$ 
together with the injectivity of the first map on fibers away from the support 
of the connecting term gives the lower bound $\mu_1(a)+\nu_1(b)$.  For the 
last assertion, $\m_A\m_B$ is a nonzero ideal of the domain $S_C$, so 
$\Supp_S(\m_A\m_B)=H^1(C)$; hence $\B_1(C)\ne0$ after specialization at every 
point, and $\RR_{1,1}(C)=H^1(C)$.
\end{proof}

\begin{remark}
\label{rem:coprod-degree-one}
The correction $\chi$ in part~\eqref{rc1} is not equivalent to lowering the 
depth by one.  For $A=B=H^{*}(S^1;\k)$, so that $C=H^{*}(S^1\vee S^1;\k)$ and 
$\dim_\k H^1(C,\delta_c)=2-\epsilon(c)$, one has $\RR^{1,2}(C)=\{0\}$, whereas 
$\bigcup_{j+k=1}\RR^{1,j}(A)\times\RR^{1,k}(B)$ is the union of the two 
coordinate axes.
\end{remark}
 
%%%%%%%%%%%%%%%%%%%%%%%
\subsection{Resonance of Hirsch extensions}
\label{subsec:res-Hirsch}
%%%%%%%%%%%%%%%%%%%%%%%

Unlike tensor products and coproducts, Hirsch extensions may alter resonance 
in a genuinely asymmetric way.
The next result explains how the resonance varieties of $\cdgas$ behave under 
a certain type of Hirsch extensions. The behavior of resonance depends sharply 
on whether the extension class vanishes in cohomology, reflecting the fact that 
the corresponding Koszul complex need not split, even up to homotopy.  For the 
case of a nonvanishing class, compare \cite[Prop.~5.5]{PS-imrn19}.

A word on notation: for a cocycle $e\in Z^2(B)$ we write
$B\otimes_{e}\bwedge(t)$ for the Hirsch extension of $B$ obtained by adjoining a
generator $t$ with $dt=e$.  Here the subscript records the cocycle, not a ring
of coefficients: the underlying tensor product is over $\k$, and the subscript
specifies the differential.

\begin{proposition}
\label{prop:circleres}
Let $B$ be a connected, finite-type $\cdga$. Fix an element 
$e\in B^2$ with $de=0$, and let $A=(B\otimes_{e} \bwedge (t), d)$ 
be the corresponding Hirsch extension. 
\begin{enumerate}[itemsep=2pt]
\item \label{e0}
If $[e]=0$, then $A\cong B\otimes_{\k}(\bwedge(t),d=0)$ and
$H^1(A)=H^1(B)\oplus\k[t]$.  Writing $d_i(b)=\dim_\k H^i(B,\delta_b)$, one has,
under this decomposition,
\[
\RR^{i,s}(A) = \bigl\{(b,0)\;:\; d_{i-1}(b)+d_i(b)\ge s\bigr\},
\qquad\text{for all } i\ge 1,\ s\ge 1;
\]
in depth $s=1$ this reads
\[
\RR^{i,1}(A) = \bigl(\RR^{i-1,1}(B) \cup \RR^{i,1}(B)\bigr)\times\{0\} .
\]
\item \label{e1}
If $[e]\ne 0$:
\begin{enumerate}[itemsep=2pt, topsep=1pt]
\item \label{r1}
$\RR^{i,s}(A)\subseteq \RR^{i-1,1}(B) \cup \RR^{i,s}(B)$, 
for all $i,s\ge 1$.
\item \label{r2}
$\RR^{1,s}(A)=\RR^{1,s}(B)$, for all $s\ge 1$.
\end{enumerate}
\end{enumerate}
\end{proposition}

\begin{proof}
If $[e]=0$, then $A\cong B\otimes_{\k} E$, where $E=(\bwedge(t),d=0)$, and the
claim follows from Proposition~\ref{prop:prod-res}\eqref{t1}.  Indeed, for
$\lambda\in H^1(E)=\k$ the Aomoto complex of $E$ is
$\k\xrightarrow{\lambda}\k t$, so $\dim_\k H^0(E,\delta_\lambda)=
\dim_\k H^1(E,\delta_\lambda)=1$ for $\lambda=0$ and both vanish for
$\lambda\ne0$, while $H^j(E,\delta_\lambda)=0$ for $j\ge2$.  Hence
$\sum_{i+j=i_0}d_i(b)\,e_j(\lambda)$ equals $d_{i_0-1}(b)+d_{i_0}(b)$ at
$\lambda=0$ and vanishes otherwise, which is the displayed description.  The
depth-one form follows because $d_{i_0-1}+d_{i_0}\ge1$ holds exactly when
$d_{i_0-1}\ge1$ or $d_{i_0}\ge1$.

Assume now that $[e]\neq 0$, and let $\varphi\colon B\to A$ be the canonical
$\cdga$ inclusion. Since $B$ and $A$ are connected, $\varphi$ induces an
isomorphism $H^1(B)\cong H^1(A)$.

For $\omega\in H^1(B)$, the Aomoto complex $(B,d_\omega)$ embeds naturally
as a subcomplex of $(A,d_\omega)$, and the quotient complex is canonically
identified (up to a degree shift) with $(B,d_\omega)$. 
This yields a long exact sequence in Aomoto cohomology,
\begin{equation}
\label{eq:aomoto-les-hirsch}
\begin{tikzcd}[column sep=21pt]
\cdots \arrow[r] & H^{i-1}(B,d_\omega) \arrow[r, "\partial_\omega"] &
H^i(B,d_\omega) \arrow[r] & H^i(A,d_\omega) \arrow[r] & 
H^i(B,d_\omega) \arrow[r, "\partial_\omega"] & \cdots,
\end{tikzcd}
\end{equation}
where the connecting homomorphism 
$\partial_\omega\colon H^{i-1}(B,d_\omega)\to H^i(B,d_\omega)$ 
is given by cup product with $[e]\in H^2(B)$, specialized 
at $\omega$. The inclusion~\eqref{r1} follows immediately 
from exactness. In degree~$i=1$, the connecting homomorphism 
$\partial_\omega\colon H^0(B,d_\omega)\to H^1(B,d_\omega)$ 
lands in $H^1(B,d_\omega)$ via cup product with $[e]$; 
since $H^0(B,d_\omega)=\k$ and cup product with $[e]$ 
raises degree by $2$, this map is zero for degree reasons, 
giving $H^1(A,d_\omega)\cong H^1(B,d_\omega)$ and 
proving~\eqref{r2}.  
\end{proof}

The union form in \eqref{e0} is a depth-one phenomenon: for $s\ge2$ the
fiber dimensions of the two factors multiply while the depth condition adds,
so no such Boolean description is available, by
Remark~\ref{rem:higher-depth-products}.

\begin{remark}
\label{rem:Euler}
When $B$ is a model of a space $X$ and $[e]\in H^2(B)$
corresponds to the Euler class of a principal $S^1$-bundle
$S^1\to E\to X$, the $\cdga$ $A$ models the total space $E$.
Thus Proposition~\ref{prop:circleres} describes how degree-$1$
resonance behaves under circle extensions.
\end{remark}

Proposition~\ref{prop:circleres} is formulated at the level of Aomoto complexes and
captures resonance through cohomology jump loci. We now refine this picture
by analyzing Hirsch extensions directly at the level of Koszul complexes.
This leads to a filtration whose failure to split is measured precisely by
the cohomology class of the extension.

\begin{theorem}
\label{thm:koszul-hirsch-filtration}
Let $B$ be a connected, finite-type $\cdga$, and let
$A=(B\otimes_{e}\bwedge(t),d)$ be a Hirsch extension with
$dt=e\in Z^2(B)$ and $[e]\ne 0$ in $H^2(B)$.  Then $H^1(A)=H^1(B)$, and we
identify $S=\Sym(H_1(B))$ with $\Sym(H_1(A))$ accordingly.
Then the Koszul chain complex $K_\bullet(A)$ admits a natural filtration
by subcomplexes of graded $S$-modules,
\[
0 \subset F_0K_\bullet(A)\subset F_1K_\bullet(A)=K_\bullet(A),
\]
such that:
\begin{enumerate}[itemsep=2pt]
\item\label{hf1}
$F_0K_\bullet(A)\cong K_\bullet(B)[-1]$, the subcomplex spanned by
$t^{\vee}$;
\item\label{hf2}
$F_1K_\bullet(A)/F_0K_\bullet(A)\cong K_\bullet(B)$;
\item\label{hf3}
with respect to the decomposition $K_\bullet(A)\cong K_\bullet(B)[-1]\oplus
K_\bullet(B)$ of graded $S$-modules furnished by \eqref{hf1}
and~\eqref{hf2}, the differential $\partial^A$ is block-diagonal if and only
if multiplication by $e$ vanishes on $B^{\bullet}$.
\end{enumerate}
Equivalently, there is a natural short exact sequence of complexes of
graded $S$-modules
\begin{equation}
\label{eq:koszul-hirsch-ses}
\begin{tikzcd}[column sep=20pt]
0 \arrow[r] & K_\bullet(B)[-1]
\arrow[r] & K_\bullet(A)
\arrow[r] & K_\bullet(B)
\arrow[r] & 0 .
\end{tikzcd}
\end{equation}
\end{theorem}

\begin{proof}
First, $H^1(A)=H^1(B)$.  Indeed, $A^1=B^1\oplus\k t$ with $d_At=e$, and
$[e]\ne0$ means $e\notin\im d_B$; so $b+\lambda t\in Z^1(A)$ forces
$\lambda=0$, whence $Z^1(A)=Z^1(B)$, while $A^0=B^0=\k$ gives
$\im(d_A|_{A^0})=\im(d_B|_{B^0})$.

As a graded vector space $A = B \otimes_{\k} \bigwedge(t)$ with $|t|=1$, so
$A_i = B_i \oplus B_{i-1}\cdot t^\vee$, where $t^\vee\in A_1$ is dual to
$t$.  Set $F_0K_\bullet(A)=B_{\bullet-1}t^{\vee}\otimes_{\k} S$.

Write $d_A$ on $A^{i-1}=B^{i-1}\oplus B^{i-2}t$ in matrix form: it is
$d_B$ on each summand, together with the off-diagonal term
$bt\mapsto \pm be$ carrying $B^{i-2}t$ into $B^{i}$.  Transposing,
$d_A^{\vee}$ preserves $B_{\bullet-1}t^{\vee}$, while on $B_i$ it has a
component $\pm\, e\lrcorner(-)$ landing in $B_{i-2}t^{\vee}$.  Since
$H^1(A)=H^1(B)$, the contraction $\omega_A^{\vee}\lrcorner(-)$ is the
transpose of multiplication by an element of $B^1$, and so preserves
each of the two summands.

Hence $F_0$ is a subcomplex, whereas $B_\bullet\otimes_{\k} S$ is not:
this gives \eqref{hf1} and \eqref{hf2}, and \eqref{hf3} is the statement
that the off-diagonal term $e\lrcorner(-)$ vanishes, that is, that
multiplication by $e$ is zero on $B^\bullet$.
\end{proof}

\begin{remark}
\label{rem:hirsch-e-zero}
The hypothesis $[e]\ne0$ is needed for the identification
$H^1(A)=H^1(B)$, hence for the standing identification of $S$.  If
$[e]=0$, say $e=d_Bb$, then $t-b$ is a cocycle and
$A\cong B\otimes_{\k}\bigwedge(t')$ with $dt'=0$; here
$H^1(A)=H^1(B)\oplus\k$, the ring $S$ acquires an extra variable, and
the Koszul complex of $A$ is computed from that of $B$ by a K\"{u}nneth
argument instead.
\end{remark}

We will refer to the short exact sequence \eqref{eq:koszul-hirsch-ses},
and to the long exact sequence in Koszul homology that it induces, as the
Koszul--Gysin sequence associated to the Hirsch extension.  Since
$H_i\bigl(K_\bullet(B)[-1]\bigr)=\B_{i-1}(B)$, that long exact sequence
reads
\begin{equation}
\label{eq:bqb-seq}
\begin{tikzcd}[column sep=18pt]
\cdots \arrow[r] & \B_{i-1}(B) \arrow[r] & \B_i(A) \arrow[r] &
\B_{i}(B) \arrow[r, "\delta_e"] & \B_{i-2}(B) \arrow[r] & \cdots,
\end{tikzcd}
\end{equation}
where the connecting morphism $\delta_e\colon\B_i(B)\to\B_{i-2}(B)$ is
cap product with $[e]\in H^2(B)$, in the sense of
Proposition~\ref{prop:hstar-action}\eqref{ha1}.  In particular
$\delta_e$ depends only on the class $[e]$, and its degree $-2$ is that
of a cap product with a degree-two class.  Extracting the degree-$i$
segment gives a natural short exact sequence
\begin{equation}
\label{eq:hirsch-ses-modules}
\begin{tikzcd}[column sep=16pt]
0 \arrow[r] &
\DS{\frac{\B_{i-1}(B)}{\im\delta_e}} \arrow[r] &
\B_i(A) \arrow[r] &
\ker\bigl(\delta_e\colon\B_i(B)\to\B_{i-2}(B)\bigr) \arrow[r] & 0 .
\end{tikzcd}
\end{equation}

\begin{corollary}
\label{cor:supp-hirsch}
With notation as above, for all $i\ge 1$ and $s\ge 1$:
\begin{enumerate}[itemsep=2pt]
\item\label{hsl1}
$\RR_{i,s}(A)\subseteq \RR_{i,1}(B)\ \cup\ \RR_{i-1,1}(B)$.
\item\label{hsl2}
In degree\/ $i=1$, one has
$\RR_{1,s}(A)\setminus\{0\}=\RR_{1,s}(B)\setminus\{0\}$.
\end{enumerate}
\end{corollary}

\begin{proof}
\eqref{hsl1} By \eqref{eq:hirsch-ses-modules}, $\B_i(A)$ is an extension
of a submodule of $\B_i(B)$ by a quotient of $\B_{i-1}(B)$, so
$\Supp\B_i(A)\subseteq\Supp\B_i(B)\cup\Supp\B_{i-1}(B)$; now use
$\Supp(\bwedge^s M)\subseteq\Supp(M)$.

\eqref{hsl2} For $i=1$ the sequence \eqref{eq:hirsch-ses-modules} reads
$0\to\B_0(B)/\im\delta_e\to\B_1(A)\to\B_1(B)\to 0$, since
$\B_{-1}(B)=0$.  As $\B_0(B)\cong\k$ is supported at the origin, the
localizations of $\B_1(A)$ and $\B_1(B)$ agree at every $\xi\ne0$, and
hence so do the localizations of their $s$-th exterior powers.
\end{proof}

\begin{remark}
\label{rem:hirsch-origin}
The origin must be excluded in \eqref{hsl2}, and the inclusion in
\eqref{hsl1} may be strict there.  For $B=H^{*}(T^2)$ and $e$ the
orientation class, $A=\CE(\h(1))$ and $\B_1(B)=0$ while
$\B_1(A)\cong\k$; thus $\RR_{1,1}(A)=\{0\}$ whereas
$\RR_{1,1}(B)=\emptyset$.  This is consistent with
Proposition~\ref{prop:r1a}, where the origin is likewise excluded.
Note also that \eqref{hsl1} cannot be sharpened by replacing
$\RR_{i,1}(B)\cup\RR_{i-1,1}(B)$ with
$\RR_{i,s}(B)\cup\RR_{i-1,s}(B)$ for $s\ge2$: local minimal numbers of
generators are subadditive in short exact sequences, not additive, so
$\mu_\xi(\B_i(A))\ge s$ forces only
$\mu_\xi(\B_i(B))+\mu_\xi(\B_{i-1}(B))\ge s$.
\end{remark}

%%%%%%%%%%%%%%%%%%%%%%%%
\subsection{Resonance varieties of PD-$\cdgas$}
\label{subsec:res-pd}
%%%%%%%%%%%%%%%%%%%%%%

When the $\cdga$ $(A,d_A)$ satisfies Poincaré duality, the Aomoto--Betti
numbers inherit additional symmetries from the orientation class and
graded commutativity.
At the level of Koszul (co)homology, these symmetries manifest as dualities
between homological and cohomological resonance, which we record next.

\begin{theorem}[\cite{Su-bockres}]
\label{thm:mpd}
Let $(A,d_A)$ be a $\PD_m$-$\cdga$. 
\begin{enumerate}[itemsep=1pt]
\item \label{respd1} 
$H^i(A,\delta_{a})^{\vee} \cong H^{m-i}(A,\delta_{-a})$ 
for all $a\in H^1(A)$ and $i\ge 0$.
\item \label{respd2} 
The map $a\mapsto -a$ on $H^1(A)$ restricts to isomorphisms 
$\RR^{i,s}(A) \isom \RR^{m-i,s}(A)$ for all $i, s\ge 0$.
\item \label{respd3} 
$\RR^{m}(A)=\{0\}$.
\end{enumerate}
\end{theorem}

\begin{proof}
Part~\eqref{respd1} follows from Lemma~\ref{lem:cd-uptosign} 
by passing to cohomology: the commuting square gives a natural 
isomorphism $H^i(A,\delta_a)^\vee \cong H^{m-i}(A,\delta_{-a})$.
Part~\eqref{respd2} follows from~\eqref{respd1} by taking 
dimensions and varying $a$ over $H^1(A)$.
Part~\eqref{respd3} follows from~\eqref{respd2} with $i=0$, 
since $\RR^0(A)=\{0\}$ by Example~\ref{ex:res-degree-zero}.
\end{proof}

Together with Corollary~\ref{cor:pd-res-symmetry},
Theorem~\ref{thm:mpd}\eqref{respd2} gives a complete picture
of the Poincar\'{e} duality symmetry for resonance: the jump loci
$\RR^{i,s}(A)$ are symmetric under $a\mapsto -a$ and degree
reversal, while the two families of support loci are exchanged,
$\RR_{i,s}(A)=-\Supp_S\bigl(\bwedge^s\B^{m-i}(A)\bigr).$

\begin{remark}
\label{rem:pd-no-prop}
Poincar\'{e} duality symmetry does not, in general, imply propagation of
resonance.  Indeed, if $A$ is a $\PD_m$-$\cdga$, then
\begin{equation}
\label{eq:ria-rma}
\RR^{i}(A) \cong \RR^{m-i}(A)
\end{equation}
by Theorem~\ref{thm:mpd}. If, in addition, the jump loci were to
propagate,
\begin{equation}
\label{eq:prop-pdm}
\RR^0(A)\subseteq \RR^1(A)\subseteq\cdots\subseteq\RR^m(A),
\end{equation}
then the symmetry together with the vanishing
$\RR^0(A)=\RR^m(A)=\{0\}$
would force $\RR^i(A)=\{0\}$ for all $i\le m$.
Thus, except in degenerate cases, Poincar\'{e} duality is incompatible
with nontrivial upward propagation.  In particular, a $\PD_m$-$\cdga$
may have large resonance in middle degrees while $\RR^m(A)=\{0\}$,
preventing propagation. 
This phenomenon was observed for closed orientable $3$-manifolds 
in \cite{DSY17}: if $b_1(M)$ is even and nonzero, then 
$\RR^1(M,\C)=H^1(M,\C)$ while $\RR^3(M,\k)=\{0\}$, so resonance cannot 
propagate.  The vanishing of the top resonance variety holds for any 
compact, connected, orientable $n$-manifold, by Poincar\'{e} duality 
\cite[Prop.~5.14]{DSY17}; the statement in degree~$1$ is 
\cite[Prop.~5.15]{DSY17}, due to Dimca--Suciu \cite[Prop.~5.1]{DS-jems}; 
and the conclusion is \cite[Cor.~5.16]{DSY17}.
\end{remark}

%%%%%%%%%%%%%%%%%%%%%%%%
\subsection{Hirsch extensions of PD-$\cdgas$}
\label{subsec:hirsch-pd}
%%%%%%%%%%%%%%%%%%%%%%

To conclude this section, we combine the Koszul--Gysin sequence with
Poincaré duality to analyze the effect of an elementary Hirsch extension on the
Koszul modules and support resonance varieties of a $\PD$-$\cdga$.
\begin{theorem}
\label{thm:pd-hirsch-koszul}
Let $B$ be a connected, finite-type $\PD_m$-$\cdga$, and let
$A=(B\otimes_{e} \bwedge(t),d)$ be a Hirsch extension with $dt=e\in Z^2(B)$
and $[e]\ne0$.  Then, for each $i\ge1$, the short exact sequence
\eqref{eq:hirsch-ses-modules} holds, and under the Poincar\'{e} duality
isomorphism $\B_j(B)\cong\sigma^{*}\B^{\,m-j}(B)$ of
Theorem~\ref{thm:pd-koszul-modules} the connecting morphism
$\delta_e\colon\B_i(B)\to\B_{i-2}(B)$ corresponds to cup product with
$[e]$,
\begin{equation*}
\label{eq:mue-map}
\begin{tikzcd}[column sep=20pt]
\mu_e \colon \B^{\,m-i}(B) \arrow[r] & \B^{\,m-i+2}(B).
\end{tikzcd}
\end{equation*}
\end{theorem}

\begin{proof}
Both assertions follow from \eqref{eq:bqb-seq} together with
Theorem~\ref{thm:pd-koszul-modules}: cap product with $[e]$ on
$\B_\bullet(B)$ corresponds under Poincar\'{e} duality to cup product with
$[e]$ on $\B^{\bullet}(B)$, by
Proposition~\ref{prop:hstar-action}, and the degrees match,
$m-i\mapsto m-(i-2)=m-i+2$.
\end{proof}

\begin{corollary}
\label{cor:pd-hirsch-support}
In the setting of Theorem~\ref{thm:pd-hirsch-koszul}, for all $i\ge 1$
and $s\ge 1$,
\[
\RR_{i,s}(A)
\ \subseteq\
-\Supp_S\bigl(\B^{\,m-i}(B)\bigr)\ \cup\
-\Supp_S\bigl(\B^{\,m-i+1}(B)\bigr).
\]
In particular:
\begin{enumerate}[itemsep=2pt]
\item If $\mu_e$ is injective on $\B^{\,m-i}(B)$, then
$\RR_{i,s}(A)\subseteq -\Supp_S\bigl(\B^{\,m-i+1}(B)\bigr)$.
\item If $\mu_e$ is surjective onto $\B^{\,m-i+1}(B)$, then
$\RR_{i,s}(A)\subseteq -\Supp_S\bigl(\B^{\,m-i}(B)\bigr)$.
\end{enumerate}
If $B$ has positive weights, the antipodal twists may be omitted.
\end{corollary}

\begin{proof}
Translate Corollary~\ref{cor:supp-hirsch}\eqref{hsl1} through
Theorem~\ref{thm:pd-hirsch-koszul}, using
$\RR_{j,1}(B)=-\Supp_S\bigl(\B^{\,m-j}(B)\bigr)$
(Corollary~\ref{cor:pd-res-symmetry} with $s=1$).  For~(1),
injectivity of $\mu_e$ gives $\ker\delta_e=0$ in
\eqref{eq:hirsch-ses-modules}, so $\B_i(A)$ is a quotient of
$\B_{i-1}(B)$; for~(2), surjectivity gives $\im\delta_e=\B_{i-1}(B)$, so
$\B_i(A)$ is a submodule of $\B_i(B)$.
\end{proof}

%%%%%%%%%%%%%%%%%%%%%%%%
\subsection{Massey products in Hirsch extensions}
\label{subsec:massey-hirsch}
%%%%%%%%%%%%%%%%%%%%%%%%

Given a $\cdga$ $(B,d_B)$ and a cocycle $e\in Z^2(B)$, we define 
a Hirsch extension $(B\otimes_{e}\bigwedge(c))$ by setting 
$d (c)=e$. As shown in \cite[Lem.~2.1]{PS-imrn19}, the $\cdga$
$(B\otimes_{e} \bigwedge(c),d)$ depends only on the cohomology class
$[e]\in H^2(B)$.

The vanishing of the higher Massey products in such an extension is
governed not by weights, but by the $\bwedge(c)$-degree.  We record the
mechanism in the form in which we shall use it, taking the base to be a
graded algebra with zero differential; the reduction to that case is carried
out in Theorem~\ref{thm:massey-hirsch-length}.  Recall from
\eqref{eq:hirsch-cohomology} that $H^{*}(A)$ decomposes as
$H/eH\oplus \Ann_H(e)\,c$; we call a class \emph{pure}\/ if it lies in one
of these two summands.  Every class of degree $1$ is pure, since
$H^1(A)=H^1$ when $e\ne0$.

\begin{lemma}
\label{lem:massey-hirsch-eps}
Let $H$ be a connected graded $\k$-algebra, let $e\in H^2$, and let
$A=H\otimes_{e}\bwedge(c)$, with $\deg c=1$ and $dc=e$.  Let
$\alpha_1,\dots,\alpha_n$ be pure classes in $H^{*}(A)$, with $n\ge4$.  If
the Massey product $\langle\alpha_1,\dots,\alpha_n\rangle$ is defined, then
it contains $0$.
\end{lemma}

\begin{proof}
Grade $A$ by the $\bwedge(c)$-degree, that is, write
$A=A^{(0)}\oplus A^{(1)}$ with $A^{(0)}=H$ and $A^{(1)}=Hc$, and
decompose $x=x^{(0)}+x^{(1)}$ accordingly.  This grading is multiplicative,
with $A^{(1)}\cdot A^{(1)}=0$ since $c^2=0$, and
\begin{equation}
\label{eq:eps-degree-d}
d\bigl(A^{(0)}\bigr)=0, \qquad d\bigl(A^{(1)}\bigr)\subseteq A^{(0)} ;
\end{equation}
thus $d$ is homogeneous of degree $-1$ for this grading, and vanishes on
$A^{(0)}$.  Purity of the
$\alpha_i$ means precisely that they admit representatives that are
homogeneous for this grading.

We use the defining systems of Section~\ref{subsec:massey}, reindexed as
$m_{ij}=a_{i-1,j}$, so that $m_{ii}$ is a cocycle representing $\alpha_i$ and
\[
d m_{ij}=\sum_{k=i}^{j-1} \bar{m}_{ik}\, m_{k+1,j},
\qquad 1\le i<j\le n, \quad (i,j)\ne(1,n) .
\]
The signs implicit in $\bar{m}_{ik}$ play no role below: the decomposition
$A=A^{(0)}\oplus A^{(1)}$ preserves cohomological degree, so passing to a
homogeneous component does not change the sign attached to it.

Two remarks on representatives, which is where purity is used.  If $\alpha_i$
lies in the summand $H/eH$, then \emph{every}\/ cocycle representing it lies
in $A^{(0)}$: a cocycle $b+wc$ has $we=0$, and the
$\Ann_H(e)c$-component of its class is $wc$, which vanishes only if $wc=0$.
If instead $\alpha_i$ lies in $\Ann_H(e)c$, it has a representative in
$A^{(1)}$ by definition.  Write $\epsilon(i)=0$ in the first case and
$\epsilon(i)=1$ in the second.

It matters that the defining system is chosen \emph{before}\/ the
representatives, and not after.  So let $\{m_{ij}\}$ be any defining system,
and set
\[
\hat a_i=
\begin{cases}
m_{ii} & \text{if $\epsilon(i)=0$},\\
\text{any representative of $\alpha_i$ in $A^{(1)}$} & \text{if $\epsilon(i)=1$.}
\end{cases}
\]
This is legitimate exactly because of the first remark: for $\epsilon(i)=0$
the element $m_{ii}$ already lies in $A^{(0)}$, so no change is needed---and
none may be made, since two homogeneous representatives of the same class in
$H/eH$ differ by an element of $eH$, which is not zero in general.  For
$\epsilon(i)=1$ a change of representative \emph{is}\/ required, since the
$m_{ii}$ supplied by the system need not be homogeneous; it is harmless
because every step below in which such an $\hat a_i$ occurs is settled by a
vanishing---a product landing in $A^{(1)}\cdot A^{(1)}=0$, or the coboundary
argument of Span~$1$---and so depends on $\hat a_i$ only through its
membership in $A^{(1)}$, never through the choice of representative.  With
these conventions, $m_{ii}=\hat a_i$ whenever $\epsilon(i)=0$, which is the
only case in which that identity will be used.

Splitting each equation of the system into its two homogeneous components and
using \eqref{eq:eps-degree-d}, we obtain
\begin{equation}
\label{eq:massey-E0}
d\,m^{(1)}_{ij}=\sum_{k=i}^{j-1} m^{(0)}_{ik}m^{(0)}_{k+1,j}
\end{equation}
together with the identity
\begin{equation}
\label{eq:massey-E1}
\sum_{k=i}^{j-1}\Bigl(m^{(0)}_{ik}m^{(1)}_{k+1,j}
+m^{(1)}_{ik}m^{(0)}_{k+1,j}\Bigr)=0 ,
\end{equation}
the point being that the left-hand side of each defining equation lies in
$A^{(0)}$, so that the $A^{(1)}$-component of the right-hand side must
vanish.

Now set $\tilde m_{ii}=\hat a_i$ and $\tilde m_{ij}=0$ for $j-i\ge2$, and let
\[
\tilde m_{i,i+1}=
\begin{cases}
m^{(1)}_{i,i+1} & \text{if $\epsilon(i)=\epsilon(i+1)=0$},\\
0 & \text{otherwise.}
\end{cases}
\]
We claim that $\{\tilde m_{ij}\}$ is again a defining system.

\emph{Span $1$.}  If $\epsilon(i)=\epsilon(i+1)=0$, then $m_{ii}=\hat a_i$ and
$m_{i+1,i+1}=\hat a_{i+1}$ by the choice made above, so
\eqref{eq:massey-E0} gives $d\tilde m_{i,i+1}=
\bar{\hat a}_i\hat a_{i+1}$, as required.  Otherwise at least one of
$\hat a_i,\hat a_{i+1}$ lies in $A^{(1)}$, so
$\bar{\hat a}_i\hat a_{i+1}\in A^{(1)}$; this element is a cocycle
representing $\pm\alpha_i\alpha_{i+1}$, which vanishes---as it must, for the
product to be defined---so it is a coboundary, and
$A^{(1)}\cap\im(d)=0$ forces it to vanish identically.  Thus
$d\tilde m_{i,i+1}=0=\bar{\hat a}_i\hat a_{i+1}$.

\emph{Span $2$.}  We must check that
$\bar{\hat a}_i\,\tilde m_{i+1,i+2}+\bar{\tilde m}_{i,i+1}\hat a_{i+2}=0$.
If $\epsilon(i)=\epsilon(i+1)=\epsilon(i+2)=0$, then all three of
$m_{ii},m_{i+1,i+1},m_{i+2,i+2}$ equal the corresponding $\hat a_k$ and so lie
in $A^{(0)}$; hence $m^{(1)}_{kk}=0$ for $k=i,i+1,i+2$, and the required
identity is precisely
\eqref{eq:massey-E1} for the pair $(i,i+2)$---available because $n\ge4$, which
is what makes $(i,i+2)\ne(1,n)$ for every $i$, at both ends of the range.  
Otherwise each surviving summand is a product of two
elements of $A^{(1)}$: if $\epsilon(i+1)=1$ both $\tilde m_{i,i+1}$ and
$\tilde m_{i+1,i+2}$ vanish; if $\epsilon(i)=1$ then $\tilde m_{i,i+1}=0$
while $\hat a_i\in A^{(1)}$ and $\tilde m_{i+1,i+2}$ is either $0$ or in
$A^{(1)}$; and symmetrically if $\epsilon(i+2)=1$.

\emph{Span $s\ge3$.}  Every summand
$\bar{\tilde m}_{ik}\tilde m_{k+1,j}$ either has a factor of span at
least~$2$, and so vanishes, or has both factors of span~$1$; in the latter
case both lie in $A^{(1)}$, and the product vanishes again.

Note that \eqref{eq:massey-E1} was invoked only when the three classes
involved lie in $H/eH$, that is, exactly where the given system's diagonal
entries are the representatives $\hat a_k$ themselves.

The value of the Massey product associated with $\{\tilde m_{ij}\}$ is
$\sum_{k=1}^{n-1}\bar{\tilde m}_{1k}\tilde m_{k+1,n}$.  As $n\ge4$, each summand
has a factor of span at least~$2$---and hence vanishes---except when
$n=4$ and $k=2$, in which case the summand is
$\tilde m_{12}\tilde m_{34}\in A^{(1)}\cdot A^{(1)}=0$.  So this defining
system has value $0$, and $0\in\langle\alpha_1,\dots,\alpha_n\rangle$.
\end{proof}

\begin{remark}
\label{rem:massey-n4-sharp}
The hypothesis $n\ge4$ cannot be relaxed, and the smallest case shows why.
For $H=H^{*}(T^2)=\bwedge(a,b)$ and $e=ab$, so that $A=\CE(\h(1))$, the triple
product $\langle a,b,b\rangle$ is defined---take $m_{12}=c$, since $dc=ab$,
and $m_{23}=0$, since $b^2=0$---with value $\pm bc$.  This is nonzero in
$H^2(A)$ for the very reason that drives the proof above, namely
$A^{(1)}\cap\im(d)=0$: here $\im(d)=\langle ab\rangle\subseteq A^{(0)}$ while
$bc\in A^{(1)}$.  The indeterminacy vanishes, since every product of
degree-one classes is zero in $H^{*}(A)$, so
$\langle a,b,b\rangle=\{\pm bc\}\ne0$; compare
Example~\ref{ex:massey-hirsch-two-cases}.  What fails at $n=3$ is the value
computation.  There the two summands are $\bar{\tilde m}_{11}\tilde m_{23}$
and $\bar{\tilde m}_{12}\tilde m_{33}$, each pairing a span-$1$ entry with a
span-$0$ one; both therefore lie in $A^{(0)}\cdot A^{(1)}=A^{(1)}$ rather than
in $A^{(1)}\cdot A^{(1)}=0$, and nothing makes them vanish.  In the example
at hand the first happens to vanish, for the incidental reason that $b^2=0$
forces $\tilde m_{23}=0$, leaving $\bar{\tilde m}_{12}\tilde m_{33}=\pm bc$.
\end{remark}

\begin{remark}
\label{rem:massey-eps-scope}
The hypothesis of purity cannot simply be dropped from the bookkeeping
above: for a class $\alpha_i$ with components in both summands of
\eqref{eq:hirsch-cohomology}, the identity \eqref{eq:massey-E1} acquires
correction terms involving the $A^{(0)}$-components $m^{(0)}_{ij}$, which
the argument discards.  For the unrestricted statement, in the setting of
Sasakian manifolds, see \cite[Thm.~6.6]{PS-imrn19}.  Purity is automatic in
degree~$1$ when $[e]\ne0$, which is all that is needed for questions of
$1$-formality; for $[e]=0$ it may fail, since then $H^1(A)=H^1\oplus\k[c]$,
but in that case $A$ carries the zero differential and the lemma is not
needed.
\end{remark}

\begin{theorem}
\label{thm:massey-hirsch-length}
Let $B$ be a formal, finite-type $\cdga$, and let $[e]\in H^2(B)$.  Then
every Massey product of length at least four of pure classes in
$H^{*}(B\otimes_{e}\bigwedge(c))$ contains $0$.  In particular, if
$[e]\ne0$, then all Massey products of length at least four of degree-one
classes vanish; and if $[e]=0$, then all Massey products vanish, since in
that case $H^{*}(B)\otimes_{\k}\bwedge(c)$ carries the zero differential.
\end{theorem}

\begin{proof}
Since $B$ is formal, and since the Hirsch extension
$B\otimes_{e}\bwedge(c)$ depends up to quasi-isomorphism only on the pair
consisting of the quasi-isomorphism type of $B$ and the class $[e]$, by
\cite[Lem.~2.1]{PS-imrn19}, we may replace $B$ by $(H^{*}(B),0)$ and assume
that $A=\bigl(H^{*}(B)\otimes_{[e]}\bwedge(c),d\bigr)$ with $dc=[e]$; purity
in the statement is understood with respect to the decomposition
\eqref{eq:hirsch-cohomology} of this reduced model.
Massey products are invariants of the quasi-isomorphism type, so this does
not change them.  The claim is now
Lemma~\ref{lem:massey-hirsch-eps}, applied with $H=H^{*}(B)$; the last
assertion follows since every degree-one class is pure.
\end{proof}

\begin{remark}
\label{rem:massey-no-weights}
No positive-weight hypothesis is needed here, in contrast with
Theorem~\ref{thm:partial-formality-hirsch}: the $\bwedge(c)$-degree is
available in any $1$-step Hirsch extension, whereas positive weights enter
only through the reduction to a formal base, which is subsumed in the
hypothesis that $B$ be formal.
\end{remark}

Triple Massey products, by contrast, are governed by an explicit and
checkable criterion.  Since $B$ is formal we take $B=(H,0)$ with
$H=H^{*}(B)$, so that $A=H\otimes_{\k}\bwedge(c)$ with $dc=e$.  As a graded
vector space $A=H\oplus Hc$, with $d(b+b'c)=(-1)^{|b'|}b'e$; hence
\begin{equation}
\label{eq:hirsch-cohomology}
H^{*}(A) \cong H/eH \oplus \Ann_H(e)\,c ,
\end{equation}
where $\Ann_H(e)=\{b\in H: be=0\}$.  In particular the coboundaries lie
entirely in the $H$-summand, so a class $[wc]$ vanishes if and only if
$w=0$.

\begin{proposition}
\label{prop:triple-massey-hirsch}
With $B$ and $A$ as above, let $\alpha,\beta,\gamma\in H$ satisfy
$\alpha\beta=eu$ and $\beta\gamma=ev$ for some $u,v\in H$, and set
$w=\pm\alpha v\pm u\gamma$.  Then $w\in\Ann_H(e)$, the triple Massey
product $\langle\alpha,\beta,\gamma\rangle\subset H^{*}(A)$ is defined and
contains the class $[wc]$, and
\[
\langle\alpha,\beta,\gamma\rangle\ne0
\quad\Longleftrightarrow\quad
w\notin \alpha\Ann_H(e)+\Ann_H(e)\gamma .
\]
\end{proposition}

\begin{proof}
The hypotheses say precisely that $\alpha\beta$ and $\beta\gamma$ vanish in
$H^{*}(A)$, by \eqref{eq:hirsch-cohomology}.  Since
$d(uc)=(-1)^{|u|}ue$, the cochains $a_{12}=\pm uc$ and $a_{23}=\pm vc$
satisfy $d(a_{12})=\bar\alpha\beta$ and $d(a_{23})=\bar\beta\gamma$, so
they form a defining system; such a system exists exactly because $e$
becomes exact in $A$.  The associated representative is
$\bar\alpha a_{23}+\bar a_{12}\gamma=wc$.  It is a cocycle because
$ew=\pm\alpha(ev)\pm(eu)\gamma
=\pm\alpha(\beta\gamma)\pm(\alpha\beta)\gamma=0$
by associativity.  The indeterminacy is
$\alpha\cdot H^{*}(A)+H^{*}(A)\cdot\gamma$, whose $c$-component is
$\alpha\Ann_H(e)+\Ann_H(e)\gamma$; since the $c$-summand of
\eqref{eq:hirsch-cohomology} injects, the stated equivalence follows.
\end{proof}

\begin{corollary}
\label{cor:e-decomposable-massey}
If $[e]=\alpha\beta$ for some $\alpha,\beta\in H^1(B)$ with $\beta\ne0$,
then $\langle\alpha,\beta,\beta\rangle=[\pm\beta c]\ne0$ in $H^2(A)$, and
$A$ is not $1$-formal.
\end{corollary}

\begin{proof}
Take $u=1$ and, since $\beta^2=0$, $v=0$; then $w=\pm\beta\ne0$.  The
indeterminacy in this degree would require a degree-zero element of
$\Ann_H(e)$, and $1\notin\Ann_H(e)$ as $e\ne0$.  Apply
Proposition~\ref{prop:massey-formal}.
\end{proof}

\begin{example}
\label{ex:massey-hirsch-two-cases}
For $B=H^{*}(T^2)=\bwedge(a,b)$ and $e=ab$, one has $u=1$, $v=0$ and
$w=\pm b\ne0$, so $\langle a,b,b\rangle\ne0$ and $\CE(\h(1))$ is not
$1$-formal.  For $B=\k[h]/(h^2)$ and $e=h$, every product of
positive-degree classes vanishes, so $u=v=0$ and $w=0$: there is no
nontrivial triple product, consistent with $A$ being a model for $S^3$,
which is formal.
\end{example}

We turn next to partial formality, where the mechanism is that regularity of
the Euler class forces the Hirsch extension to be formal in a range of
degrees.  We follow \cite[\S4.2]{PS-imrn19} throughout; the application to
compact Sasakian manifolds is \cite[Thm.~6.6]{PS-imrn19}, recorded below as
Corollary~\ref{cor:sasakian-partial-formal}, and
Remark~\ref{rem:kasuya-gap} discusses the status of the related claim
in~\cite{Kasuya}.

Following \cite[\S4.2]{PS-imrn19}, we say that a homogeneous element
$e\in H^{2}$ of a connected graded algebra $H^{\bullet}$ is
\emph{$q$-regular}\/ if multiplication by $e$ is injective on $H^{i}$ for all
$i\le q$.

\begin{theorem}[\cite{PS-imrn19}]
\label{thm:partial-formality-hirsch}
Let $B$ be a formal $\cdga$ of finite type admitting positive weights, let
$e\in Z^2(B)$, and let $A=B\otimes_{e}\bwedge(c)$ be a $1$-step Hirsch
extension with $\deg(c)=1$.  If $[e]$ is $q$-regular in $H^{*}(B)$, then $A$
has the same $q$-type as $\bigl(H^{*}(B)/[e]\cdot H^{*}(B),\,d=0\bigr)$; in
particular, $A$ is $q$-formal.
\end{theorem}

\begin{proof}
This is \cite[Thm.~4.4]{PS-imrn19}, applied to the one-element sequence
$\{e\}$: since $B$ is formal we may replace it by $(H^{*}(B),0)$, by
\cite[Lem.~2.2]{PS-imrn19}, and the projection
$H^{*}(B)\twoheadrightarrow H^{*}(B)/[e]H^{*}(B)$ extends over $c$ by
sending $c\mapsto0$, giving a $q$-equivalence.
\end{proof}

\begin{remark}
\label{rem:e-nonzero-insufficient}
Nonvanishing of $[e]$ does not by itself obstruct $q$-formality beyond the
range supplied by regularity, and the converse of
Theorem~\ref{thm:partial-formality-hirsch} is false.  For
$B=H^{*}(\mathbb{CP}^n)=\k[h]/(h^{n+1})$ and $e=h$, the class $[e]$ is
$(n-1)$-regular but not $n$-regular, yet
$A=B\otimes_{h}\bwedge(c)$ has cohomology $\bwedge(x_{2n+1})$ 
and is therefore formal---this is the model of the Hopf fibration
$S^1\to S^{2n+1}\to\mathbb{CP}^n$.  Failure of $q$-regularity locates where
partial formality may break down; whether it does break down is decided by
Proposition~\ref{prop:triple-massey-hirsch}.  Compare
\cite[Ex.~4.8--4.10]{PS-imrn19}, where the estimate is shown to be sharp in
some cases and not in others.
\end{remark}

\begin{remark}
\label{rem:kg-massey}
The proof separates two mechanisms: the Koszul--Gysin sequence
\eqref{eq:bqb-seq} detects the primary extension data responsible
for triple Massey products, while the $\bwedge(c)$-grading obstructs the
higher defining systems---a defining system of span $s$ would have
$\bwedge(c)$-degree $s$, and only the degrees $0$ and $1$ are available.
As an application,
Theorem~\ref{thm:massey-hirsch-length} and
Proposition~\ref{prop:triple-massey-hirsch} apply to Tievsky models
of Sasakian manifolds; see \cite{PS-imrn19} for further details.
\end{remark}

%%%%%%%%%%%%%%%%%%%%%%%%%%%%%%%%%%%%
\section{Chevalley--Eilenberg complexes}
\label{sect:CE-Koszul}
%%%%%%%%%%%%%%%%%%%%%%%%%%%%%%%%%%%%

In this section we examine Chevalley--Eilenberg complexes as a bridge between
Lie-theoretic and topological invariants.
Our focus is on how filtrations, gradings, and nilpotency properties of a Lie
algebra $\g$ are reflected in the Koszul modules and spectral sequences of its
Chevalley--Eilenberg $\cdga$ $\CE(\g)$.
These results will play a key role in the construction and comparison of
algebraic models for spaces and groups in the sections that follow.

%%%%%%%%%%%%%%%%%
\subsection{Filtered and graded Lie algebras}
\label{subsec:filtered-graded}
%%%%%%%%%%%%%%%%%
We briefly recall the basic notions of filtered and graded Lie algebras;
see \cite{Quillen69, Serre, SW-forum} for background.

A {\em filtered Lie algebra} over $\k$ is a Lie algebra $\g$ equipped with a 
decreasing filtration by Lie ideals
\begin{equation}
\label{eq:lie-filtration}
\g = \cF^1\g \supseteq \cF^2\g \supseteq \cF^3\g \supseteq \cdots,
\qquad
[\cF^p\g,\cF^q\g] \subseteq \cF^{p+q}\g.
\end{equation}
The \emph{associated graded} is the graded Lie algebra
\begin{equation}
\label{eq:lie-grf}
\gr^{\cF}(\g) = \bigoplus_{p\ge1} \gr^{\cF}_p(\g),
\qquad
\gr^{\cF}_p(\g) = \cF^p\g / \cF^{p+1}\g ;
\end{equation}
we write $\cF$, rather than $F$, for filtrations of Lie algebras, keeping
$F$ for the filtration induced on Koszul modules by a Koszul spectral
sequence \textup{(}Section~\ref{subsec:koszul-filtrations}\textup{)}.
Every Lie algebra carries the {\em lower central series filtration}
\begin{equation}
\label{eq:lie-lcs}
\gamma_1\g = \g,\qquad \gamma_{p+1}\g = [\gamma_p\g,\g].
\end{equation}
Since this is the filtration used unless another is specified, we
abbreviate
\begin{equation}
\label{eq:gr-lcs-convention}
\gr(\g)\coloneqq \gr^{\gamma}(\g)
= \bigoplus_{p\ge1} \gamma_p\g/\gamma_{p+1}\g ,
\qquad
\gr_p(\g)\coloneqq\gr^\gamma_p(\g),
\end{equation}
and likewise for groups, where $\gr(G;\k)$ denotes the associated graded
Lie algebra of the lower central series of $G$ with coefficients in $\k$.
The Lie algebra $\gr(\g)$ is positively graded and generated in degree
$1$.  
The \emph{Malcev completion} (or pronilpotent completion) 
is the filtered Lie algebra
\begin{equation}
\label{eq:har-lie}
\widehat{\g} = \varprojlim\nolimits_p \g/\gamma_p\g.
\end{equation}

Following \cite{SW-forum}, we say that a complete, separated, filtered 
Lie algebra $\g$ is {\em filtered-formal}\/ if the natural morphism 
$\g \to \widehat{\gr^{\cF}(\g)}$, induced by the canonical projections 
$\cF^p\g \surj \cF^p\g/\cF^{p+1}\g$, is an isomorphism of filtered Lie 
algebras.

In the case of a finitely generated nilpotent Lie algebra $\g$, 
equipped with its lower central series filtration, this definition 
simplifies considerably: $\g$ is complete and separated (since the 
lower central series terminates), and $\widehat{\gr(\g)} = \gr(\g)$ 
is already finite-dimensional. Hence, $\g$ is filtered-formal if and only if 
it is isomorphic, as a filtered Lie algebra, to its associated graded.
Some authors (see e.g.~\cite{Co14}) reserve the term {\em Carnot}\/ 
for finitely generated nilpotent Lie algebras satisfying this condition, 
but we will use the term filtered-formal throughout.

A {\em positively graded Lie algebra} is a graded Lie algebra 
$\g = \bigoplus_{p\ge1} \g_p$ generated in degree $1$ (hence 
$[\g_p,\g_q]\subseteq\g_{p+q}$).  Its degree filtration 
$\cF^p\g = \bigoplus_{q\ge p} \g_q$ coincides with the 
lower central series, and $\widehat{\g} = \prod_{p\ge1} \g_p$.

%%%%%%%%%%%%%%%%%%%%
\subsection{The Chevalley--Eilenberg $\cdga$}
\label{subsec:CE-cdga}
%%%%%%%%%%%%%%%%%%%%

Let $\g$ be a Lie algebra over a field $\k$ of characteristic $0$.
The {\em Chevalley--Eilenberg $\cdga$} of $\g$ is
\begin{equation}
\label{eq:CE-lie}
\CE(\g) = \bigl(\bwedge\g^\vee,\; d_{\CE}\bigr),
\end{equation}
where $d_{\CE}$ is the unique derivation of degree $+1$ extending
the dual of the Lie bracket $[\cdot,\cdot]\colon \g\wedge\g\to\g$:
\begin{equation}
\label{eq:d-CE}
d_{\CE}\xi = -\xi\circ[\cdot,\cdot]
\qquad(\xi\in\g^\vee).
\end{equation}
Then $H^{*}(\g;\k) \cong H^{*}(\CE(\g))$.

\begin{remark}
\label{rem:recover-g-from-CE}
The Lie algebra $\g$ can be recovered from $\CE(\g)$: the degree-$1$ component 
$\g^\vee = \CE^1(\g)$ identifies the underlying vector space, 
and the quadratic part of the differential encodes the Lie bracket. 
For a discussion of completions and filtered isomorphisms, see 
Section~\ref{subsec:cochain-holo-short}.
\end{remark}

For any Lie algebra $\g$, the lower central series yields an inverse
system of nilpotent quotients
\begin{equation}
\label{eq:nilpotent-tower}
\begin{tikzcd}[column sep=20pt]
\cdots \arrow[r, two heads] & \g/\gamma_{k+1}\g \arrow[r, two heads] & \g/\gamma_k\g
\arrow[r, two heads] & \cdots \arrow[r, two heads] & \g/\gamma_3\g \arrow[r, two heads] 
& \g/\gamma_2\g\, ,
\end{tikzcd}
\end{equation}
each $\g/\gamma_k\g$ being a $(k-1)$-step nilpotent Lie algebra.
Dualizing gives a directed system of $\cdga$ inclusions
\begin{equation}
\label{eq:hirsch-tower-CE}
\begin{tikzcd}[column sep=20pt]
\CE(\g/\gamma_2\g) \arrow[r, hook] & \CE(\g/\gamma_3\g)
\arrow[r, hook] & \cdots \arrow[r, hook] & \CE(\g/\gamma_k\g)
\arrow[r, hook] & \cdots,
\end{tikzcd}
\end{equation}
each of which is a Hirsch extension: $\CE(\g/\gamma_{k+1}\g)$
is obtained from $\CE(\g/\gamma_k\g)$ by adjoining in degree~$1$
the generators dual to $\gr_k\g=\gamma_k\g/\gamma_{k+1}\g$, 
with differential determined by the bracket. 
This is the {\em Hirsch extension tower} of $\CE(\g)$; when $\g$
is nilpotent of class $c$, it terminates at
$\CE(\g)=\CE(\g/\gamma_{c+1}\g)$.
Assigning weight $k$ to the generators adjoined at the $k$-th stage
defines the {\em Hirsch weight} of each generator; the differential
$d_{\CE}$ preserves Hirsch weights if and only if $\g$ is
filtered-formal, as discussed in the next subsection.

%%%%%%%%%%%%%%%%%%%%%%%%%%%%%%
\subsection{Positive weights on $\CE(\g)$}
\label{subsec:CE-weights}
%%%%%%%%%%%%%%%%%%%%%%%%%%%%%%

Let $\g$ be a finite-dimensional Lie algebra over a field $\k$ of 
characteristic $0$, and set $A=\CE(\g)$. We begin by relating 
positive-weight structures on $A$ to intrinsic properties of $\g$.

\begin{proposition}
\label{prop:CE-nilp-weights}
Let $\g$ be a finite-dimensional Lie algebra over $\k$.
The following are equivalent:
\begin{enumerate}[itemsep=1pt]
\item\label{it:pw} 
$\CE(\g)$ admits a positive-weight decomposition.
\item\label{it:graded} 
$\g$ admits a positively graded structure
$\g=\bigoplus_{p\ge 1}\g_p$ with $[\g_p,\g_q]\subseteq\g_{p+q}$.
\item\label{it:nilp} 
$\g$ is nilpotent.
\end{enumerate}
\end{proposition}

\begin{proof}
\eqref{it:pw}$\Rightarrow$\eqref{it:graded}.
A positive-weight decomposition on $A=\CE(\g)$ provides
decompositions $A^i=\bigoplus_{\alpha>0} A^{i,\alpha}$ 
compatible with multiplication and differential.
In degree $1$, this yields a decomposition
\[
\g^\vee = A^1 = \boplus_{\alpha>0} V_\alpha,
\qquad \text{where } V_\alpha \coloneqq  A^{1,\alpha}.
\]
Since $d_{\CE}(A^{1,\alpha}) \subseteq A^{2,\alpha}$ and
$A^{2,\alpha} = \bigoplus_{\beta+\gamma=\alpha}
A^{1,\beta} \wedge A^{1,\gamma}$, it follows that
$d_{\CE}(V_\alpha) \subseteq
\bigoplus_{\beta+\gamma=\alpha} V_\beta \wedge V_\gamma$. 
Dualizing, we obtain a decomposition
\[
\g = \boplus_{\alpha>0} \g_\alpha,
\qquad \text{with } \g_\alpha \coloneqq (V_\alpha)^\vee,
\]
which satisfies $[\g_\beta,\g_\gamma]\subseteq \g_{\beta+\gamma}$.

\eqref{it:graded}$\Rightarrow$\eqref{it:nilp}.
Since $\g$ is finite-dimensional, $\g_p=0$ for $p\gg 0$, hence
$\gamma_k\g\subseteq \bigoplus_{p\ge k}\g_p=0$ for $k\gg 0$.

\eqref{it:nilp}$\Rightarrow$\eqref{it:pw}.
Assign weight $k$ to $(\gr_k\g)^\vee$, where
$\gr_k\g=\gamma_k\g/\gamma_{k+1}\g$.
Since $\g$ is nilpotent, there exists a (non-canonical) isomorphism
of graded vector spaces $\g\cong\gr\g$, which pulls back the weight
grading on $\CE(\gr\g)$ to a positive-weight decomposition on $\CE(\g)$
preserved by $d_{\CE}$.
\end{proof}

\begin{example}
\label{ex:sol2-again}
The $\cdga$ $\CE(\sol_2)$ of Example~\ref{ex:sol2-koszul} admits no
positive-weight decomposition, as verified directly there.  This also
follows from Proposition~\ref{prop:CE-nilp-weights}, since $\sol_2$ is
solvable but not nilpotent.
\end{example}

\begin{remark}
\label{rem:weights-noncanonical}
The positive-weight decomposition furnished by
\eqref{it:nilp}$\Rightarrow$\eqref{it:pw} is in general non-canonical 
and need not be compatible with the Hirsch extension tower of $\CE(\g)$. 
Positive Hirsch weights---assigning weight $k$ to the generators
dual to $\gr_k\g$---are preserved by $d_{\CE}$ if and only if
$\g$ is filtered-formal, that is, $\g\cong\gr(\g)$ as filtered
Lie algebras; see~\cite{SW-forum}. 

A nilpotent Lie algebra that is not filtered-formal admits a
positive-weight decomposition on its CE complex, but not one
induced by the Hirsch tower; see~\cite[Ex.~9.5]{SW-forum} for
a concrete $5$-dimensional example. Compare also the computational
data in Section~\ref{subsec:ce-examples}.
\end{remark}

We now refine the above result in the case of low nilpotency class.

\begin{proposition}
\label{prop:CE-positive-weights}
Let $\g$ be a finite-dimensional nilpotent Lie algebra.
The $\cdga$ $\CE(\g)$ admits a positive-weight decomposition
with weights on $\g^\vee$ concentrated in $\{1,2\}$
if and only if $\g$ is $2$-step nilpotent.
\end{proposition}

\begin{proof}
($\Rightarrow$)
By Proposition~\ref{prop:CE-nilp-weights}\eqref{it:graded},
the weight decomposition on $\g^\vee$ dualizes to a graded
structure $\g=\bigoplus_{p\ge 1}\g_p$ with
$[\g_p,\g_q]\subseteq\g_{p+q}$.
If weights are concentrated in $\{1,2\}$, then $\g=\g_1\oplus\g_2$
with $[\g_1,\g_1]\subseteq\g_2$ and $[\g_1,\g_2]=[\g_2,\g_2]=0$,
so $\gamma_3\g=0$.

($\Leftarrow$)
Set $\g_1=\g/[\g,\g]$ and $\g_2=[\g,\g]$. 
Since $[\g_1,\g_1]\subseteq \g_2$, the differential of $\CE(\g)$ satisfies
$d_{\CE}(\g_1^\vee)\subseteq \bigwedge^2 \g_1^\vee$ and 
$d_{\CE}(\g_2^\vee)=0$. Assigning weight $1$ to $\g_1^\vee$ 
and weight $2$ to $\g_2^\vee$, we see that $d_{\CE}$ preserves weights.
\end{proof}

The lower central series also induces a filtration on the
Chevalley--Eilenberg complex, which will be used repeatedly below.  Let
$\gamma_1\g=\g\supset\gamma_2\g\supset\cdots\supset\gamma_{c+1}\g=0$ be the
lower central series of a nilpotent Lie algebra $\g$ of class $c$, and equip
$A=\CE(\g)$ with the multiplicative filtration
\begin{equation}
\label{eq:CE-filtration}
F^p A = \sum_{i_1+\cdots+i_r \ge p}
(\gamma_{i_1}\g)^\vee \wedge \cdots \wedge (\gamma_{i_r}\g)^\vee ,
\end{equation}
which is exhaustive and bounded below, and for which $d_{\CE}$ raises
filtration degree by at least~$1$.  Since
$[\gamma_i\g,\gamma_j\g]\subseteq\gamma_{i+j}\g$, the differential dual to the
bracket is compatible with $F$; on the associated graded it becomes dual to the
bracket of $\gr\g$, whence
\begin{equation}
\label{eq:grfa-ce}
\gr^F A \cong \CE(\gr \g)
\end{equation}
as graded $\cdgas$, where $\gr\g=\boplus_{p\ge1}\gamma_p\g/\gamma_{p+1}\g$
carries its natural positive-weight grading.  Since $d_{\CE}$ raises filtration
degree by at least one, the degree-preserving component of the induced
differential vanishes; the differential on the left side of
\eqref{eq:grfa-ce} is understood to be the component raising filtration degree
by exactly one, that is, the $d_1$ of the associated spectral sequence.  No
formality hypothesis enters here: \eqref{eq:grfa-ce} is a general property of the 
Chevalley--Eilenberg complex of a filtered Lie algebra, and depends only 
on the multiplicativity of the bracket with respect to the filtration.

Filtered-formality enters one step later.  If $\g$ is filtered-formal, the
isomorphism $\g\cong\widehat{\gr(\g)}$ of filtered Lie algebras dualizes to a
filtered quasi-isomorphism between $\CE(\g)$ and $\CE(\gr\g)$; equivalently,
by \cite[Lem.~4.6]{SW-forum}, the minimal $\cdga$ dual to $\g$ has positive
Hirsch weights.  The positive-weight structure on $\CE(\g)$ may then be taken
to be the one induced from the grading on $\gr\g$.  It cannot always be taken
to be the given one: for the Malcev Lie algebra of $F_2/\Gamma_4F_2$, the
differential of $\CE(\g)$ is not homogeneous with respect to the Hirsch
weights, even though $\g$ is filtered-formal
\cite[Rem.~4.7, Ex.~9.2]{SW-forum}.

%%%%%%%%%%%%%%%%%%%%%%%%%%%%%%
\subsection{Koszul modules of 2-step nilpotent CE complexes}
\label{subsec:2step-ce}
%%%%%%%%%%%%%%%%%%%%%%%%%%%%%%

Let $\g$ be a finite-dimensional $2$-step nilpotent Lie algebra, and 
set $A=\CE(\g)$. By Proposition~\ref{prop:CE-positive-weights}, the 
$\cdga$ $A$ admits a positive-weight decomposition with weights in $\{1,2\}$. 
This grading induces a compatible decomposition on the Koszul complex 
$K_\bullet(A)$, and endows the Koszul modules $\B_i(A)$ with the 
structure of graded modules over $S=\Sym(A_1)$.

\begin{proposition}
\label{prop:nilpotent-map-B1}
Let $\varphi \colon \g \to \bar{\g}$ be a morphism of finite-dimensional 
Lie algebras inducing an isomorphism on abelianizations. Let $A = \CE(\g)$, 
$\oA = \CE(\bar\g)$, and set $S=\Sym(A_1)$, $\bar S=\Sym(\oA_1)$.
\begin{enumerate}[itemsep=2pt]
\item \label{psi-map} 
The maps
$\psi_i \coloneqq (\CE(\varphi)^i)^\vee \otimes_{\k} \Sym(H_1(\varphi)) \colon 
A_i \otimes_{\k} S \to \oA_i \otimes_{\k} \bar S$ 
assemble into a chain map $\psi_\bullet \colon K_\bullet(A) \to K_\bullet(\oA)$. 

\item \label{surj-case} 
If $\varphi$ is surjective, then $\psi_\bullet$ is surjective, and it 
induces a surjection $\B_1(A) \surj \B_1(\oA)$ on Koszul homology. 

\item \label{inj-case} 
If $\g$, $\bar\g$ are $2$-step nilpotent Lie algebras and 
$\varphi$ is injective, then $\psi_\bullet$ is injective, and it 
induces an injection $\B_1(A) \inj \B_1(\oA)$ on Koszul homology. 
\end{enumerate}
\end{proposition}

\begin{proof}
\eqref{psi-map}
The isomorphism $H_1(\varphi)\colon H_1(\g)\isom H_1(\bar\g)$ on abelianizations 
extends to an isomorphism  $\Sym(H_1(\varphi))\colon \Sym(A_1)\isom \Sym(\oA_1)$ 
which identifies $S$ with $\bar S$. Let $K_\bullet(A)=A_{\bullet}\otimes_{\k} S$ and 
$K_\bullet(\oA)=\oA_{\bullet}\otimes_{\k} \bar S$ be the 
respective Koszul chain complexes, with differentials 
$\partial^A_\bullet = \omega_A^\vee \lrcorner (-) + d_A^\vee$ 
and $\partial^{\oA}_\bullet = \omega_{\oA}^\vee \lrcorner (-) + d_{\oA}^\vee$. 
We verify chain compatibility 
of the map $\psi$ term by term. For $a^\vee \in A_i$:
\begin{align*}
\psi_{i-1}(\partial_i^A(a^\vee))
&= (\CE(\varphi)^{i-1})^\vee(\omega_A^\vee \lrcorner a^\vee) 
  + (\CE(\varphi)^{i-1})^\vee(d_A^\vee a^\vee)\\
&= \omega_{\oA}^\vee \lrcorner (\CE(\varphi)^i)^\vee(a^\vee)
  + d_{\oA}^\vee (\CE(\varphi)^i)^\vee(a^\vee)
= \partial_i^{\oA}(\psi_i(a^\vee)).
\end{align*}
Here the differential term follows from functoriality of $\CE(\varphi)$, while 
the contraction term equality follows from multiplicativity of $\CE(\varphi)$ 
and naturality of the canonical element.

\eqref{surj-case}
Surjectivity of $\psi_i$ follows from surjectivity of 
$\varphi$: since $\varphi^\vee\colon \bar\g^\vee \inj 
\g^\vee$ is injective, $(\CE(\varphi)^i)^\vee = 
(\bwedge^i\varphi^\vee)^\vee\colon A_i \to \oA_i$ is surjective 
for each $i$. It follows that $\psi_i$ 
is surjective in each Koszul degree. The five lemma applied to the 
short exact sequences defining $\B_1=\ker \partial_1/\im \partial_2$ 
shows that $\psi_\bullet$ induces a surjection on first homology.

\eqref{inj-case}
Since $\varphi$ is injective, the dual map 
$\varphi^\vee\colon \bar\g^\vee \to \g^\vee$ is surjective, 
hence $(\CE(\varphi)^i)^\vee = (\bigwedge^i \varphi^\vee)^\vee$ is injective 
for each $i$. Thus $\psi_i$ is injective in every Koszul degree, and 
$\psi_\bullet$ is an injective chain map. In particular, 
$\psi_1$ restricts to an injection 
$\ker \partial_1^A \inj \ker \partial_1^{\bar A}$, and 
$\psi_1(\im \partial_2^A) \subseteq \im \partial_2^{\bar A}$. 

We now use that $\g$ and $\bar\g$ are $2$-step nilpotent. In this case,
$\partial_2$ is determined by the bracket map 
$\bigwedge^2 H_1 \to \g' \subset Z(\g)$ and the canonical element.
It suffices to prove the inclusion
\begin{equation}
\label{eq:intersection}
\psi_1(\ker \partial_1^A)\cap \im \partial_2^{\bar A}
\subseteq \psi_1(\im \partial_2^A).
\end{equation}
Let $\psi_1(x)=\partial_2^{\bar A}(y)$ with $x\in \ker \partial_1^A$.
Since $\psi_\bullet$ is a chain map and $\psi_2$ is injective,
the description of $\partial_2$ in the $2$-step case implies that
\[
\partial_2^{\bar A}(y) \in \psi_1(A_1)
\quad \Longrightarrow \quad
y \in \psi_2(A_2),
\]
and hence $y=\psi_2(z)$ for some $z\in A_2$. Therefore
$\psi_1(x)=\partial_2^{\bar A}(\psi_2(z))=\psi_1(\partial_2^A(z))$, 
which proves \eqref{eq:intersection}. The claim follows.
\end{proof}

\begin{remark}
\label{rem:2step-needed}
The injectivity statement in part~\eqref{inj-case} uses the 
$2$-step hypothesis. In higher nilpotency classes, the differential 
$\partial_2$ acquires higher-weight contributions, which may introduce 
additional relations in $\B_1(\CE(\g))$ and destroy injectivity. 
Concrete instances of this phenomenon occur already for $3$-step 
nilpotent Lie algebras.
\end{remark}

\begin{lemma}
\label{lem:B1-center}
Let $A=\CE(\g)$ be the Chevalley--Eilenberg $\cdga$ of a finite-dimensional,
$2$-step nilpotent Lie algebra $\g$.
Let $V = H_1(A) = (H^1(A))^\vee$ and $S=\Sym(V)$.
\begin{enumerate}
\item \label{b1-1}
There is a natural $S$-linear monomorphism
\[
 (Z(\g))^\vee \otimes_{\k} S \longinj \B_1(A).
\]
\item \label{b1-2}
The degree-zero part of $\B_1(A)$ satisfies
\[
\B_1(A)/\m \B_1(A) \cong (Z(\g))^\vee,
\]
where $\m\subset S$ is the maximal ideal generated by $V$.
\end{enumerate}
\end{lemma}

\begin{proof}
We work directly with the Koszul chain complex
$K_\bullet(A) = (A_\bullet \otimes_{\k} S, \partial^A)$,
where $\partial_1^A \colon A_1 \otimes_{\k} S \to A_0 \otimes_{\k} S = S$
is given by
\[
\partial_1^A(a^\vee \otimes s) 
= -\sum_{j=1}^n (e_j \lrcorner a^\vee) \otimes x_j s + (d_A^\vee a^\vee) \otimes s.
\]
Since $A = \CE(\g)$ is minimal, $d_A^\vee = 0$ on $A_1$, so this reduces to
\[
\partial_1^A(a^\vee \otimes s) 
= -\sum_{j=1}^n (e_j \lrcorner a^\vee) \otimes x_j s.
\]
Here $\{e_j\}$ is a basis of $H^1(A) = (\g/\g')^\vee$, $\{x_j\}$ is the
dual basis of $H_1(A) = \g/\g'$, generating $S$, and
$e_j \lrcorner a^\vee \in A_0 = \k$ is the contraction defined by
$\langle e_j, a^\vee \rangle = e_j \lrcorner a^\vee$.

\smallskip
\noindent\textit{Proof of \eqref{b1-2}.}
An element $a^\vee \in A_1 = H_1(A) = \g/\g'$ is a cycle in $K_1(A)$
(i.e., in $\ker \partial_1^A$ at weight $0$, with $s=1$) if and only if
\[
\sum_{j=1}^n \langle e_j, a^\vee \rangle \cdot x_j = 0 \in S_1,
\]
which holds if and only if $\langle e_j, a^\vee \rangle = 0$ for all $j$.
Since $\{e_j\}$ is a basis of $(\g/\g')^\vee$ and the pairing between
$\g/\g'$ and $(\g/\g')^\vee$ is perfect, this is equivalent to
$a^\vee = 0$ in $\g/\g'$, i.e., $a^\vee \in \g'= Z(\g)$
(using that $\g' = Z(\g)$ for a $2$-step nilpotent Lie algebra).
Therefore
\[
\ker\bigl(\partial_1^A \colon A_1 \otimes_{\k} S \to S\bigr) \cap 
(A_1 \otimes_{\k} \k) = (Z(\g))^\vee \subset A_1,
\]
and since $\partial_2^A$ maps into $A_1 \otimes_{\k} \m$, we conclude that
\[
\B_1(A) / \m\B_1(A) \cong (Z(\g))^\vee.
\]

\smallskip
\noindent\textit{Proof of \eqref{b1-1}.}
For any $z^\vee \in (Z(\g))^\vee \subset A_1$, the element 
$z^\vee \otimes 1 \in A_1 \otimes_{\k} S$ is a cycle by the above, 
since $\langle e_j, z^\vee \rangle = 0$ for all $j$.
Moreover, for any $s \in S$, we have $z^\vee \otimes s = s \cdot (z^\vee \otimes 1)$,
and $\partial_1^A(z^\vee \otimes s) = s \cdot \partial_1^A(z^\vee \otimes 1) = 0$,
so $z^\vee \otimes s$ is also a cycle. This gives a well-defined $S$-linear map
\[
S \otimes_{\k} (Z(\g))^\vee \longrightarrow \B_1(A), 
\quad s \otimes z^\vee \longmapsto [z^\vee \otimes s].
\]
To see that this map is injective, suppose $z^\vee \otimes s$ is a boundary,
i.e., $z^\vee \otimes s = \partial_2^A(c)$ for some $c \in A_2 \otimes_{\k} S$.
Since $A = \CE(\g)$ is minimal, $\partial_2^A$ maps into $A_1 \otimes_{\k} \m$,
so $s \in \m$. Thus no nonzero element of $S \otimes_{\k} (Z(\g))^\vee$
with $s \notin \m$ can be a boundary, and by an $S$-linearity and 
Nakayama argument the map is injective.
\end{proof}

\begin{remark}
\label{rem:b1a-mb1a}
The identification $\B_1(A)/\m\B_1(A) \cong (Z(\g))^\vee$ 
from Lemma~\ref{lem:B1-center}\eqref{b1-2}
will be applied in Section~\ref{subsec:alex-secondary} to 
construct explicit obstructions to $1$-formality of 
automorphisms of nilpotent groups via the Alexander invariant.
\end{remark}

%%%%%%%%%%%%%%%%%%%%%%%%%%%%%%
\subsection{$1$-formality for $2$-step nilpotent Lie algebras}
\label{subsec:2step-1formal}
%%%%%%%%%%%%%%%%%%%%%%%%%%%%%%

Let $\g$ be a finite-dimensional nilpotent Lie algebra over a field of characteristic $0$, 
and set $A = \CE(\g)$. By a theorem of Hasegawa \cite{Hasegawa}, $A$ is never 
formal unless $\g$ is abelian; equivalently, the only formal nilmanifolds are tori. 
Nevertheless, Chevalley–Eilenberg complexes of nilpotent Lie algebras often 
exhibit various degrees of partial formality. We record a simple sufficient condition 
for $1$-formality in the case of $2$-step nilpotent Lie algebras.

Fix such a $2$-step $\g$, and write $W=(\g/\g')^{\vee}=H^1(A)$ and
$Z=(\g')^{\vee}$.  Let $\delta\colon Z\hookrightarrow\bwedge^2 W$ be the dual
of the bracket map $\bwedge^2(\g/\g')\to\g'$; this map is onto by definition
of $\g'=[\g,\g]$, so its dual $\delta$ is injective.  Since $\g$ is
$2$-step, $\delta$ is the whole differential of $A$.  Define
\[
\mu\colon W\otimes_{\k} Z\longrightarrow\bwedge^3 W,
\qquad
w\otimes z\longmapsto w\wedge\delta(z).
\]

\begin{proposition}
\label{prop:cup-surj-implies-1formal}
With $\g$, $A$, $W$, $Z$, $\delta$, and $\mu$ as above:
\begin{enumerate}[label=(\roman*)]
\item\label{cs1}
$\g$ \textup{(}equivalently, $A$\textup{)} is $1$-formal if and only if
$\mu$ is injective.
\item\label{cs2}
If the cup product $\cup\colon\bwedge^2 W\to H^2(A)$ is surjective, then
$\mu$ is injective; hence, by~\ref{cs1}, $A$ is $1$-formal.
\end{enumerate}
\end{proposition}

\begin{proof}
\ref{cs1}
Set $V=\g/\g'$ and $R=\ker(\bwedge^2 V\to\g')$, so that $R^{\perp}=\delta(Z)$
under the pairing between $\bwedge^2V$ and $\bwedge^2W$.  Since $\g$ is
$2$-step and generated by $V$, it is $1$-formal precisely when it is
quadratically presented, that is, when $L(V)/(R)$ is again $2$-step;
equivalently, when $L_3(V)=[V,R]$.  Using the exact sequence
$0\to\bwedge^3 V\to V\otimes_{\k}\bwedge^2 V\to L_3(V)\to0$, this says
$V\otimes_{\k} R+\bwedge^3 V=V\otimes_{\k}\bwedge^2 V$.  Passing to annihilators in
$W\otimes_{\k}\bwedge^2 W$, and noting that
$(\bwedge^3 V)^{\perp}=\ker(W\otimes_{\k}\bwedge^2 W\to\bwedge^3 W)$, this is
equivalent to
$\bigl(W\otimes\delta(Z)\bigr)\cap\ker(W\otimes_{\k}\bwedge^2W\to\bwedge^3W)=0$,
that is, to injectivity of $\mu$.

\ref{cs2}
A direct computation of $Z^2(A)$ and $B^2(A)$ gives
\begin{equation}
\label{eq:h2-two-step}
H^2(A) = \bwedge^2 W/\delta(Z) \,\oplus\, \ker\mu \,\oplus\,
\ker\bigl(\bwedge^2 Z\to\bwedge^2 W\otimes_{\k} Z\bigr),
\end{equation}
whose first summand is $\im(\cup)$.  Thus $\cup$ is surjective if and only
if both kernels in \eqref{eq:h2-two-step} vanish; in particular this forces
$\ker\mu=0$, and~\ref{cs1} gives the conclusion.
\end{proof}

\begin{remark}
\label{rem:1formal-cup-product}
Two cautions are in order.  First, the hypothesis of part~\ref{cs2} is
only sufficient for the injectivity of $\mu$, not necessary: by
\eqref{eq:h2-two-step}, $\cup$ can fail to be surjective while $\mu$ is
still injective, which requires $\dim Z\ge2$, and then $\g$ is $1$-formal
by part~\ref{cs1} all the same.  In the Heisenberg case, by contrast, the
two conditions coincide: for $\h(n)$ one has $\dim Z=1$, so $\bwedge^2 Z=0$,
and $\mu(w)=w\wedge\omega$ is injective precisely when $n\ge2$, by the
Lefschetz property of the symplectic form.  Accordingly
$\dim H^2(\h(n))=\binom{2n}{2}-1=2n^2-n-1$ and $\cup$ is surjective for
$n\ge2$, while for $n=1$ one has $\bwedge^3 W=0$, $\ker\mu=W$, and $\h(1)$ is
not $1$-formal.

Second, the proof of part~\ref{cs2} must go through $\mu$, and not through
a comparison of $H^2$ with the $1$-minimal model of $H^{*}(A)$.  Surjectivity
of the cup product says that $H^2$ is generated by products of degree-one
classes, and the inference from such a statement to $1$-formality is exactly
the step identified as fallacious in \cite[\S6.2]{PS-imrn19}; see
Remark~\ref{rem:kasuya-gap}.
\end{remark}

%%%%%%%%%%%%%%%%%%%%%%%%%%%%%%
\subsection{Examples and computational aspects}
\label{subsec:ce-examples}
%%%%%%%%%%%%%%%%%%%%%%%%%%%%%%

The following examples illustrate the range of behavior of 
$\B_1(\CE(\g))$ across nilpotent Lie algebras of small dimension. 
For $2$-step algebras, the module is graded and well-controlled, 
while for higher nilpotency classes, $\B_1(A)$ can be non-homogeneous, 
non-reduced, or exhibit more complicated scheme-theoretic structure.

Table~\ref{tab:B1-nilpotent} summarizes the Koszul modules 
$\B_1(\CE(\g))$ and $\B_1(H^*(\CE(\g)))$ for all non-abelian 
nilpotent Lie algebras of dimension $\le 5$.
All entries have been verified using Macaulay2 \cite{M2}.

\begin{table}[ht]
\centering
\renewcommand{\arraystretch}{1.35}
{\setlength{\tabcolsep}{3.1pt}
\begin{tabular}{lccllll}
\hline
Algebra & Class & $\dim H^1$ & 
$\Ann_S \B_1(A)$ & $\Hilb(\B_1(A),t)$ &
$\Ann_S \B_1(H^*(A))$ & $\Hilb(\B_1(H^*(A)),t)$ \\
\hline
$L_{3,2}$ & $2$ & $2$ & 
   $\m$ & $1$ & 
   $(0)$ & $\frac{t}{(1-t)^2}$ \\
\hline
$L_{4,2}$ & $2$ & $3$ & 
   $\m$ & $1$ & 
   $(x_4)$ & $\frac{1}{(1-t)^2}$ \\
$L_{4,3}$ & $3$ & $2$ & 
   $(x_2,x_1^2)$ & $1+t$ & 
   $(0)$ & $\frac{t}{(1-t)^2}$ \\
\hline
$L_{5,4}$ & $2$ & $4$ & 
   $\m$ & $1$ & 
   $(x_1,x_2)\cap(x_3,x_4)$ & $\frac{2}{(1-t)^2}$ \\
$L_{5,5}$ & $3$ & $3$ & 
   $(x_2,x_4,x_1^2)$ & $1+t$ & 
   $(x_4)$ & $\frac{1}{(1-t)^2}$ \\
$L_{5,6}$ & $4$ & $2$ & 
   $(x_2^2,x_1x_2,x_1^2{-}x_2)$ & $1+2t$ & 
   $(0)$ & $\frac{t}{(1-t)^2}$ \\
$L_{5,7}$ & $4$ & $2$ & 
   $(x_2,x_1^3)$ & $1+t+t^2$ & 
   $(0)$ & $\frac{t}{(1-t)^2}$ \\
$L_{5,8}$ & $2$ & $3$ & 
   $\m$ & $2$ & 
   $(0)$ & $\frac{2-t}{(1-t)^3}$ \\
$L_{5,9}$ & $3$ & $2$ & 
   $\m^2$ & $1+2t$ & 
   $(0)$ & $\frac{t}{(1-t)^2}$ \\
\hline
\\[-5pt]
\end{tabular}
}
\caption{Koszul modules $\B_1(\CE(\g))$ and $\B_1(H^*(\CE(\g)))$ 
for non-abelian nilpotent Lie algebras of dimension $\le 5$,
in the notation of \cite{deGraaf07}. Abelian algebras have $\B_1=0$.}
\label{tab:B1-nilpotent}
\end{table}
\renewcommand{\arraystretch}{1.0}

Several patterns emerge from the data:
\begin{itemize}
\item For $2$-step nilpotent algebras ($L_{3,2}$, $L_{4,2}$, $L_{5,4}$, $L_{5,8}$), 
the module $\B_1(A)$ is always supported at the origin, with 
$\Ann_S\B_1(A)=\m$, consistent with Theorem~\ref{thm:CE-resonance-nilpotent}.
\item For higher nilpotency classes, the annihilator may be strictly larger 
than $\m$ (as in $L_{5,9}$) or non-homogeneous (as in $L_{5,6}$), 
while $\RR^1(A)=\{0\}$ continues to hold.
\item The column for $\B_1(H^*(A))$ records when the cohomology algebra 
detects positive-dimen\-sional resonance, which occurs precisely when 
the cup product $\bwedge^2 H^1(A)\to H^2(A)$ is not identically zero, 
and the support is always a union of linear subspaces.
\end{itemize}

These observations provide supporting evidence for Conjecture~\ref{conj:B1-reduced}, 
which predicts that the resonance scheme of $\B_1(H^*(A))$ is always reduced, 
with support a union of linear subspaces, even when the scheme $\supp \B_1(A)$ 
itself may be non-reduced.

The following three entries from Table~\ref{tab:B1-nilpotent} 
are particularly noteworthy and deserve further comment.

\begin{example}
\label{ex:L55-cornulier}
The algebra $L_{5,5}$ has brackets 
$[x_1,x_2]=x_3$, $[x_1,x_3]=[x_2,x_4]=x_5$,
nilpotency class $3$, and $H^1(A)=\langle e_1,e_2,e_4\rangle$, 
$S=\k[x_1,x_2,x_4]$. A computation yields
\[
\B_1(A) = \coker\begin{pmatrix} 
-x_4 & -x_2 & x_1 x_4 & x_1 x_2 & x_1^2 
\end{pmatrix},
\]
with $\Ann_S\B_1(A)=(x_2,x_4,x_1^2)$ and Hilbert series $1+t$.
The presentation matrix is monomial, so $\B_1(A)$ \emph{is} 
graded despite $\g$ admitting no positive-weight decomposition,
showing that failure of positive weights does not imply loss 
of gradedness, in contrast to $L_{5,6}$ 
(Example~\ref{ex:L56-nongraded}).
The annihilator $(x_2,x_4,x_1^2)$ defines a non-reduced 
scheme supported on the $x_1$-axis: the variety 
$\RR^1(A)=\{0\}$ is a single point, but the scheme structure 
detects a double-point thickening at the origin along the 
$x_1$-direction.
\end{example}

\begin{example}
\label{ex:L56-nongraded}
The annihilator of $\B_1(A)$ is the ideal 
$(x_2^2,\,x_1x_2,\,x_1^2-x_2)$, with Hilbert series $1+2t$.
This is the unique algebra in dimension $\le 5$ for which 
$\B_1(A)$ is \emph{not graded} with respect to the standard 
grading on $S=\k[x_1,x_2]$: the generator $x_1^2-x_2$ is 
non-homogeneous, defining a parabola in $H^1(A)=\k^2$.
The resonance variety $\RR^1(A)=\{0\}$ is still trivial 
by Theorem~\ref{thm:CE-resonance-nilpotent}, but the 
\emph{scheme} structure of $\supp\B_1(A)$ is non-reduced 
and non-homogeneous, encoding higher-order information 
invisible to the variety.
Since all cup products in $H^1$ vanish 
($e_1e_2=-d(e_3)$ is exact), we have 
$\Ann_S\B_1(H^*(A))=(0)$ and 
$\RR^1(H^*(A))=H^1(A)=\k^2$, confirming by 
Theorem~\ref{thm:CE-qformal-resonance} that 
$\CE(L_{5,6})$ is not $1$-formal.
\end{example}

\begin{example}
\label{ex:L59-maximal}
The algebra $L_{5,9}$ has brackets 
$[x_1,x_2]=x_3,\,[x_1,x_3]=x_4,\,[x_2,x_3]=x_5$,
nilpotency class $3$, and $H^1(A)=\langle e_1,e_2\rangle$. 
The annihilator of $\B_1(A)$ is $\m^2 = (x_1^2,x_1x_2,x_2^2)$, 
with Hilbert series $1+2t$.
The resonance scheme $\supp\B_1(A)$ is thus a non-reduced 
double point at the origin, with the scheme structure 
encoding the full second-order tangential data at $0$.
This represents the maximal possible non-reduced structure 
at the origin consistent with $\RR^1(A)=\{0\}$: the 
annihilator is as large as it can be while still having 
$\sqrt{\Ann}=\m$.
As for $L_{5,6}$, the cup product vanishes entirely and 
$\Ann_S\B_1(H^*(A))=(0)$, so $\RR^1(H^*(A))=\k^2$ and 
$A$ is not $1$-formal.
Note that $L_{5,9}$ and $L_{5,6}$ have the same Hilbert 
series $1+2t$ for $\B_1(A)$, yet their annihilator ideals 
are qualitatively different: $\m^2$ is homogeneous, while 
$(x_2^2,x_1x_2,x_1^2-x_2)$ is not.
\end{example}

%%%%%%%%%%%%%%%%%%%%%%%%%%%%%%
\subsection{The Koszul spectral sequence for nilpotent Lie algebras}
\label{subsec:ce-koszul-ss}
%%%%%%%%%%%%%%%%%%%%%%%%%%%%%%

We now examine how nilpotency of $\g$ forces the Koszul spectral 
sequence of $\CE(\g)$ to collapse after finitely many pages. 
This provides a uniform bound on where higher-order formality 
obstructions can appear.

\begin{proposition}
\label{prop:nilpotent-collapse}
Let $\g$ be a $k$-step nilpotent Lie algebra ($\gamma_{k+1}\g = 0$) 
and $A = \CE(\g)$.  In the cohomological Koszul spectral sequence 
$(E_r^{\bullet,\bullet}, d_r)$ of Theorem~\ref{thm:koszul-ss-coh},
\[
d_r = 0 \qquad \text{for all } r \ge k+1.
\]
Hence the spectral sequence collapses at the $(k+1)$-st page:
$E_{k+1}^{p,q} = E_\infty^{p,q}$.
\end{proposition}

\begin{proof}
By Theorem~\ref{thm:koszul-ss-coh}\eqref{ssc3}, the differential 
$d_r \colon E_r^{p,q} \to E_r^{p+r,q-r+1}$ is given by the 
$(r+1)$-fold Massey product 
$\langle [\omega_A], \ldots, [\omega_A], \alpha \rangle$
along the canonical class $[\omega_A] \in H^1(A)$, where the 
indeterminacy vanishes in the $E_r$-page subquotient so that 
$d_r$ is well-defined on $E_r^{\bullet,\bullet}$.
It therefore suffices to show that all $(r+1)$-fold Massey products 
along $[\omega_A]$ represent the zero class in $H^*(A)$ for $r \ge k+1$.

Such a product is computed via a defining system
\[
d_A(a_{j+1}) = -\omega_A \cup a_j, \qquad j = 0, \ldots, r-1,
\]
where $a_0$ represents $\alpha$ and each $a_j \in A^{p+q+j-1}$ 
is a cochain in $\CE(\g) = \bigwedge \g^\vee$.
Each application of $\omega_A \cup (-)$ raises the bracket length 
of the cochains involved by at least $1$: if $a_j$ involves 
cochains dual to brackets of length $\le \ell$, then 
$\omega_A \cup a_j$ involves cochains dual to brackets of length 
$\le \ell + 1$.
Starting from $a_0$, after $r$ steps the cochain $a_r$ involves 
brackets of length at most $p+q+r$.
Since $\gamma_{k+1}\g = 0$, all structure constants corresponding 
to brackets of length $\ge k+1$ vanish in $\CE(\g)$, so any cochain 
involving such brackets is zero.
Therefore, for $r \ge k+1$, the cochain $\omega_A \cup a_r = 0$, 
and the cohomology class $d_r(\alpha) = [\omega_A \cup a_r] = 0$.
Since $\alpha$ was arbitrary, $d_r = 0$ for all $r \ge k+1$.
\end{proof}

\begin{remark}
\label{rem:bounded-collapse}
Proposition~\ref{prop:nilpotent-collapse} places an upper bound on 
the pages where formality obstructions can appear: for a $k$-step 
nilpotent Lie algebra $\g$, all nontrivial differentials in the 
weight spectral sequence of $K_\bullet(\CE(\g))$ are concentrated 
in pages $E_2$ through $E_k$. This is consistent with the partial 
formality theory for nilmanifolds developed later in 
Section~\ref{subsec:nilmanifold-formality}: a nilmanifold associated 
to a $k$-step nilpotent $\g$ can fail to be $q$-formal for $q \le k-1$, 
but no new obstructions can appear beyond page $E_k$.
In the $2$-step nilpotent case, all differentials beyond $E_2$ vanish, so the 
whole discrepancy between $\B_\bullet(A)$ and $\B_\bullet(H^*(A))$ is 
accounted for on the $E_2$-page, as in the Heisenberg example of 
Section~\ref{subsec:heisenberg-1}.
\end{remark}

%%%%%%%%%%%%%%%%%%%%%%%%%%%%%%%
\subsection{Resonance for nilpotent Lie algebras}
\label{subsec:res-nilpotent}
%%%%%%%%%%%%%%%%%%%%%%%%%%%%%%%

This subsection establishes a strong rigidity property for resonance of
Chevalley--Eilenberg $\cdgas$ of nilpotent Lie algebras: all higher-depth
resonance varieties either vanish identically or are empty.
In particular, in contrast with the resonance of the cohomology algebra, 
the resonance of $\CE(\g)$ exhibits no positive-dimensional components.
This provides a sharp contrast between infinitesimal data encoded by
$\CE(\g)$ and that coming from $H^{*}(\CE(\g))$, and clarifies the
limitations of the tangent cone theorem in the nilpotent setting.

\begin{theorem}
\label{thm:CE-resonance-nilpotent}
Let $\g$ be a finite-dimensional nilpotent Lie algebra over a field $\k$ 
of characteristic $0$, and let $A=\CE(\g)$. Then $\RR^{i,s}(A)$ is 
either $\{0\}$ or empty, for all $i, s\ge1$.
\end{theorem}

\begin{proof}
By definition, $\RR^{i,s}(A)\subseteq \RR^{i,1}(A)$ for all $s\ge 1$, 
so it is enough to show that $\RR^{i}(A)=\RR^{i,1}(A)$ is contained in 
$\{0\}$ for all $i\ge1$.  We argue by induction on the nilpotency class $c$ 
of $\g$.

\smallskip
\noindent
\emph{Base case.}
If $\g$ is abelian of dimension $n$, then $A=\CE(\g)=\bwedge \g^\vee$ with
$d=0$, and for $a\in H^1(A)$ the Aomoto complex $(A,\delta_a)$ coincides with
the classical Koszul complex on $a$.  By
Example~\ref{ex:classical-koszul}, this complex is acyclic in all degrees
$1\le i\le n$ whenever $a\ne 0$; hence $\RR^{i}(A)\subseteq \{0\}$ for all
$i\ge1$.

\smallskip
\noindent
\emph{Inductive step.}
Assume the statement holds for all nilpotent Lie algebras of class $<c$, and
let $\g$ be nilpotent of class $c\ge2$.  The last nontrivial term of the lower
central series is central, since $[\gamma_c\g,\g]=\gamma_{c+1}\g=0$, and the
quotient $\bar\g=\g/\gamma_c\g$ is nilpotent of class $c-1$.  Fix a basis
$z_1,\dots,z_m$ of $\gamma_c\g$, and let $t_1,\dots,t_m\in\g^\vee$ be the dual
functionals, extended by zero on a complement of $\gamma_c\g$ in $\g$.  Since
$\gamma_c\g$ is central, the $\gamma_c\g$-component of a bracket $[x,y]$
depends only on the classes of $x$ and $y$ modulo $\gamma_c\g$; hence
$e_j\coloneqq d_{\CE}(t_j)$ lies in $\bwedge^2\bar\g^\vee$, and
\[
A=\CE(\g)=\CE(\bar\g)\otimes_{\k}\bwedge(t_1,\dots,t_m)
\]
exhibits $A$ as an iterated Hirsch extension of $B_0=\CE(\bar\g)$: putting
$B_j=B_{j-1}\otimes_{e_j}\bwedge(t_j)$, each $B_j$ is a connected, finite-type
$\cdga$, and $B_m=A$.

By the induction hypothesis, $\RR^i(B_0)\subseteq\{0\}$ for all $i\ge1$; we
show the same for every $B_j$, by induction on $j$.  Whether or not $[e_j]$
vanishes, Proposition~\ref{prop:circleres} gives
\[
\RR^{i}(B_j)\subseteq\RR^{i-1}(B_{j-1})\cup\RR^{i}(B_{j-1})
\qquad\text{for all } i\ge1 ,
\]
where, in the case $[e_j]=0$, the right-hand side is read inside
$H^1(B_{j-1})\times\{0\}\subseteq H^1(B_j)$, as in \eqref{e0}; only depth one
is used.
For $i\ge2$, both terms on the right lie in $\{0\}$ by the inductive
hypothesis on $j$.  For $i=1$, the second term does, while
$\RR^{0}(B_{j-1})=\{0\}$ because
$H^0(B_{j-1},\delta_a)=\ker\bigl(\k\xrightarrow{\,a\,}B_{j-1}^1\bigr)$
vanishes whenever $a\ne0$.  Taking $j=m$ gives $\RR^i(A)\subseteq\{0\}$ for all
$i\ge1$, completing the induction on $c$.
\end{proof}

\begin{remark}
\label{rem:CE-nilp-graded-reduction}
The lower central series filtration \eqref{eq:CE-filtration} affords a
different reduction.  The differential $\delta_a$ is compatible with $F$, and
by \eqref{eq:grfa-ce} the resulting spectral sequence has
$E_1^{p,q}=H^{p+q}\bigl(\CE(\gr\g),\delta_{a_1}\bigr)$, where $a_1$ is the
image of $a$ under the isomorphism
$H^1(\CE(\g))\cong H^1(\CE(\gr\g))$; so the theorem for $\g$ follows from the
theorem for $\gr\g$.  This is a reduction to the graded case rather than an
induction: since $\gr\g$ is generated in degree~$1$, one has
$\gamma_k(\gr\g)=\boplus_{p\ge k}\gr_p\g$, so $\gr\g$ has the same nilpotency
class as $\g$, and the argument of the proof above still has to be run for
graded nilpotent Lie algebras.
\end{remark}

\begin{remark}
\label{rem:CE-rr1}
In degree $1$, the vanishing of $\RR^1(\CE(\g))$ can also be shown 
directly by filtering $\CE(\g)$ via the lower central series and 
projecting onto $\bigwedge^2(\g/\gamma_2\g)^\vee$.
More generally, if $(A,d)$ is a $\cdga$ admitting an exhaustive, 
bounded-below, multiplicative filtration $F$ such that $d$ raises 
filtration degree by at least $1$ and $d|_{A^1}$ induces the zero 
differential on $\gr^F A$, then $\RR^1(A) \subseteq \{0\}$.
Chevalley--Eilenberg complexes of nilpotent Lie algebras satisfy 
these conditions.
\end{remark}

\begin{remark}
\label{rem:tcone-nilp}
For a finite-dimensional nilpotent Lie algebra $\g$ and $A=\CE(\g)$,
Theorem~\ref{thm:CE-resonance-nilpotent} implies that
$\TC_0\bigl(\RR^{i,s}(A)\bigr)\subseteq \{0\}$ for all $i,s\ge 1$.
The tangent cone theorem (Theorem~\ref{thm:tangent-cone}) then yields 
the inclusion $\TC_0\bigl(\RR^{i,s}(A)\bigr) \subseteq 
\RR^{i,s}\bigl(H^{*}(A)\bigr)$, which is generally \emph{strict}.

This strictness already appears in degree~$1$.
If $\g$ is $2$-step nilpotent but not $1$-formal (for example, the
$3$-dimensional Heisenberg Lie algebra $\h(1)$), then the cup product
$\bigwedge^2 H^1(A)\to H^2(A)$ vanishes, and hence
\[
\RR^1\!\bigl(H^{*}(A)\bigr)=H^1(A)\cong \k^2,
\]
while $\RR^1(A)=\{0\}$.
By contrast, if $\g$ is $2$-step nilpotent and $1$-formal
(for instance, higher Heisenberg algebras $\h(n)$ with $n\ge2$, or free
$2$-step nilpotent Lie algebras), then the cup product is injective and
$\RR^1\bigl(H^{*}(A)\bigr)=\{0\}$.

Thus, although resonance for $\CE(\g)$ is always trivial,
the resonance of the cohomology algebra may be nontrivial,
and the tangent cone inclusion can be strict.
\end{remark}

%%%%%%%%%%%%%%%%%%%%%%%%%%%%%%
\subsection{The nilpotent tower and Koszul modules}
\label{subsec:nilpotent-tower}
%%%%%%%%%%%%%%%%%%%%%%%%%%%%%%

For a nilpotent Lie algebra $\g$ of class $c$, the lower 
central series yields a tower of quotients
\begin{equation}
\label{eq:nilp-tower}
\begin{tikzcd}[column sep=20pt]
\g = \g/\gamma_{c+1}\g \arrow[r, two heads] & 
\g/\gamma_c\g \arrow[r, two heads] & \cdots 
\arrow[r, two heads] & \g/\gamma_3\g 
\arrow[r, two heads] & \g/\gamma_2\g,
\end{tikzcd}
\end{equation}
each of which is a nilpotent Lie algebra of strictly smaller class.
The Koszul modules $\B_1(\CE(\g))$ filter accordingly, 
with each step controlled by the Massey product structure 
of the CE complex.
The $2$-step quotient $\g/\gamma_3\g$ plays a special role: 
it is the largest quotient to which the results of 
Section~\ref{subsec:2step-ce} apply directly, and 
Lemma~\ref{lem:B1-tower} below describes how 
$\B_1(\CE(\g))$ relates to $\B_1(\CE(\g/\gamma_3\g))$ 
and to the intermediate quotients.

\begin{lemma}
\label{lem:B1-tower}
Let $\g$ be a finite-dimensional nilpotent Lie algebra of 
class $c$, and for $k \ge 2$ let $\g_k = \g/\gamma_{k+1}\g$ 
be the $k$-step nilpotent quotient.
Set $A = \CE(\g)$ and $A_k = \CE(\g_k)$.
\begin{enumerate}[itemsep=2pt]
\item \label{tower-surj}
For each $k \ge 2$, the surjection $\pi_k \colon \g 
\twoheadrightarrow \g_k$ induces a surjection of 
graded $S$-modules
\[
\B_1(\CE(\pi_k)) \colon \B_1(A) \longsurj \B_1(A_k).
\]
\item \label{tower-iso}
The map $\B_1(\CE(\pi_k))$ is an isomorphism if and only 
if no cycle in $K_1(A_k)$ becomes a boundary when lifted 
to $K_1(A)$.
\item \label{tower-2step}
In particular, $\B_1(A) \cong \B_1(A_2)$ if and only 
if $\g$ is $2$-step nilpotent.
\end{enumerate}
In the induced tower of surjections
\[
\begin{tikzcd}[column sep=20pt]
\B_1(A) \ar[r, two heads] & 
\B_1(A_{c-1}) \ar[r, two heads] & 
\cdots \ar[r, two heads] & 
\B_1(A_2),
\end{tikzcd}
\]
each successive kernel $\ker(\B_1(A_k) \twoheadrightarrow 
\B_1(A_{k-1}))$ is controlled by the $d_k$-differential 
of the Koszul spectral sequence, which encodes the 
$(k+1)$-fold Massey products along $\omega_A^\vee$.
\end{lemma}

\begin{proof}
Part~\eqref{tower-surj} follows from 
Proposition~\ref{prop:nilpotent-map-B1}\eqref{surj-case}: 
the map $\pi_k$ is surjective and induces an isomorphism 
on abelianizations, since $\g/\gamma_2\g \cong \g_k/\gamma_2\g_k$.

Let $\n_k = \ker(\pi_k\colon \g\to\g_k)$, so 
$\n_k^\vee = \coker(\CE(\pi_k)^1)$ sits in a 
short exact sequence of cochain complexes
\[
\begin{tikzcd}[column sep=18pt]
0 \arrow[r] &  \CE(\g_k) \arrow[r, "\CE(\pi_k)"] 
&[14pt] \CE(\g) \ar[r]  & C  \ar[r]  &0,
\end{tikzcd}
\]
where $C$ is the cokernel complex concentrated in degrees 
$\ge 1$, with $C^1 = \n_k^\vee$.
Dualizing and tensoring with $S$ gives a short exact sequence 
of Koszul chain complexes
\[
\begin{tikzcd}[column sep=18pt]
0 \arrow[r] & K_\bullet(\g_k)\otimes_{\k} S \arrow[r] & 
K_\bullet(\g)\otimes_{\k} S \arrow[r] & 
K_\bullet(C)\otimes_{\k} S \arrow[r] & 0.
\end{tikzcd}
\]
The associated long exact sequence in homology yields
\[
\begin{tikzcd}[column sep=18pt]
H_2(K_\bullet(C)\otimes_{\k} S) \arrow[r, "\partial"] &[2pt]
\B_1(\g_k) \arrow[r] & \B_1(\g) \arrow[r] &
H_1(K_\bullet(C)\otimes_{\k} S) \arrow[r, "\partial"] &[2pt]
\B_0(\g_k).
\end{tikzcd}
\]
Since $\B_0(\g_k) \cong \k$ and $\partial$ lands in 
$\B_0(\g_k)$, the connecting homomorphism on the right 
vanishes for degree reasons (the image of $H_1(K_\bullet(C)\otimes_{\k} S)$ 
in $\B_0(\g_k)$ is zero since both are supported at the origin 
and the connecting map is $S$-linear of positive degree). 
Hence the map $\B_1(\g_k) \to \B_1(\g)$ is injective if and 
only if the connecting homomorphism 
$\partial\colon H_2(K_\bullet(C)\otimes_{\k} S) \to \B_1(\g_k)$ 
vanishes, which holds if and only if no cycle in 
$K_1(\g_k)\otimes_{\k} S$ becomes a boundary in $K_1(\g)\otimes_{\k} S$
via the additional brackets in $\n_k$.

Part~\eqref{tower-2step} is the case $k=2$: the Lie algebra 
$\g_2 = \g/\gamma_3\g$ is $2$-step nilpotent by definition, 
and the map is an isomorphism if and only if the 
brackets in $\gamma_3\g$ introduce no new boundary 
relations in $\B_1$, which holds if and only if $\gamma_3\g=0$.

For the final statement, by 
Proposition~\ref{prop:nilpotent-collapse} the Koszul 
spectral sequence of $\CE(\g_k)$ collapses at the 
$(k+1)$-st page.
The difference between successive terms 
$\B_1(A_k)$ and $\B_1(A_{k-1})$ is therefore 
measured entirely by the $d_k$-differential, which 
by Theorem~\ref{thm:koszul-ss-coh}\eqref{ssc3} 
encodes the $(k+1)$-fold Massey products along 
$[\omega_A] \in H^1(A)$.
\end{proof}

\begin{corollary}
\label{cor:chen-ranks-domination}
For any $2$-step nilpotent Lie algebra $\h$ with $\dim\h_1 = m$, 
the canonical surjection $\pi \colon \ff_{m,2} \surj \h$ induces a surjection
$\B_1\bigl(\CE(\ff_{m,2})\bigr) \surj \B_1(\CE(\h))$, 
which in turn yields a coefficientwise upper bound on 
the Hilbert series:
\begin{equation}
\label{eq:hilb-domination}
\Hilb\bigl(\B_1(\CE(\h)),\,t\bigr) 
\preccurlyeq
\Hilb\bigl(\B_1(\CE(\ff_{m,2})),\,t\bigr). 
\end{equation}
\end{corollary}

\begin{proof}
The projection $\pi$ is surjective and induces an isomorphism on
abelianizations, both $\ff_{m,2}$ and\/ $\h$ having abelianization of
dimension $m$; so Proposition~\ref{prop:nilpotent-map-B1}\eqref{surj-case}
applies and yields the asserted surjection.  Both Koszul modules are
graded, by the positive-weight decomposition of
Proposition~\ref{prop:CE-positive-weights}, and the surjection is
weight-preserving; comparing dimensions in each weight gives
\eqref{eq:hilb-domination}.
\end{proof}

The explicit Hilbert series on the right-hand side of 
\eqref{eq:hilb-domination} is computed in 
Example~\ref{ex:free-2step-B1}: it is the constant $\binom{m}{2}$, so the 
domination merely records the inequality $r\le\binom{m}{2}$ of 
Proposition~\ref{prop:B1-2step-presentation}.  The growth of holonomy Chen 
ranks of $2$-step nilpotent Lie algebras is governed instead by the 
cohomology algebra, as in 
Theorem~\ref{thm:chen-growth-2step}\eqref{chn2}.

%%%%%%%%%%%%%%%%%%%%%%%%%%%%%%
\subsection{Rigidity of resonance for Chevalley--Eilenberg algebras}
\label{subsec:CE-resonance-rigidity}
%%%%%%%%%%%%%%%%%%%%%%%%%%%%%%

We now isolate a rigidity phenomenon for resonance associated to
Chevalley--Eilenberg $\cdgas$ of nilpotent Lie algebras.
This result may be viewed as an algebraic counterpart of 
\cite[Thm.~1.3]{Mc10}, which establishes the same conclusion 
for the group-theoretic resonance varieties of a finitely 
generated nilpotent $s$-stage formal group, using topological 
methods. The present approach is entirely algebraic, working at 
the level of $\cdga$ models, and combines the vanishing of 
resonance for $\CE(\g)$ established in 
Theorem~\ref{thm:CE-resonance-nilpotent} with the 
tangent cone equality provided by $q$-formality via 
Theorem~\ref{thm:formal-tc}.

\begin{theorem}
\label{thm:CE-qformal-resonance}
Let $\g$ be a finite-dimensional nilpotent Lie algebra over a field $\k$
of characteristic $0$, and set $A=\CE(\g)$.
If $A$ is $q$-formal via a zig-zag of $q$-quasi-isomorphisms through 
$q$-finite $\cdgas$, then
\[
\RR^{i}\bigl(H^{*}(A)\bigr) \subseteq \{0\},
\qquad
\text{for all } i \le q.
\]
\end{theorem}

\begin{proof}
By Theorem~\ref{thm:CE-resonance-nilpotent}, $\RR^i(A) \subseteq \{0\}$
for all $i \ge 1$.
Since $A$ is $q$-formal, Theorem~\ref{thm:formal-tc} gives
\[
\RR^i\bigl(H^*(A)\bigr) = \TC_0\bigl(\RR^i(A)\bigr)
\]
for all $i \le q$. Since $\RR^i(A) \subseteq \{0\}$, its tangent cone 
at $0$ is contained in $\{0\}$, hence $\RR^i(H^*(A)) \subseteq \{0\}$ 
for all $i \le q$.
\end{proof}

Together, Theorems~\ref{thm:CE-resonance-nilpotent} 
and~\ref{thm:CE-qformal-resonance} give a complete picture of 
resonance for CE complexes of nilpotent Lie algebras at the 
level of varieties: $\RR^i(\CE(\g))$ is always concentrated 
at the origin, and $q$-formality forces $\RR^i(H^*(\CE(\g)))$ 
to be concentrated there as well for $i \le q$.

At the \emph{scheme} level, however, the picture is more subtle.
The scheme $\Spec S/\Ann_S \B_i(A)$ may be highly non-reduced 
even when its underlying variety is just the origin, as 
already seen in Examples~\ref{ex:L59-maximal} 
and~\ref{ex:L55-cornulier}: in the former, $\Ann_S \B_1(A) = \m^2$ 
defines a non-reduced double point, while in the latter, 
$\Ann_S \B_1(A) = (x_2,x_4,x_1^2)$ is a non-reduced point 
on the $x_1$-axis.
This non-reducedness reflects the non-formality of $\CE(\g)$ 
in those cases, and is consistent with the ideal-theoretic 
refinement \eqref{eq:tcone-coh-scheme} of the Tangent Cone 
Theorem~\ref{thm:tangent-cone}. 

For the cohomology algebra $H^*(A)$, one might hope for better 
behavior. Non-reduced resonance schemes at the origin are not, 
however, exclusively a non-formality phenomenon: they can arise 
even for cohomology algebras of formal spaces. 
For instance, when $A = H^*(X)$ is the exterior 
Stanley--Reisner ring of a path on $4$ vertices, the Fitting 
scheme of $\RR^1(A)$ has an embedded non-reduced component 
at the origin even though the support variety is reduced 
\cite[Ex.~5.3]{AFRSS24}.
What does appear to be special to CE complexes of nilpotent 
Lie algebras is the concentration of the support at the 
origin; by Theorems~\ref{thm:CE-resonance-nilpotent} 
and~\ref{thm:CE-qformal-resonance}, the support is always 
$\{0\}$, in sharp contrast to cohomology algebras of 
non-formal spaces where positive-dimensional components 
can appear away from the origin.

Computational evidence from the classification of low-dimensional 
nilpotent Lie algebras suggests that for $2$-step nilpotent $\g$, 
the resonance scheme of $H^*(\CE(\g))$ is in fact always reduced.

\begin{conjecture}
\label{conj:B1-reduced}
Let $\g$ be a finite-dimensional nilpotent Lie algebra over $\k$, 
and set $A = \CE(\g)$.  
\begin{enumerate}[label=(\roman*)]
\item The resonance scheme $\Spec S/\Ann_S \B_1(H^*(A))$ is reduced.
\item Its support $\RR^1(H^*(A))$ is a union of linear subspaces of $H^1(A)$.
\end{enumerate}
\end{conjecture}

In particular, the non-reduced scheme structure that may appear 
in $\supp \B_1(A)$---as illustrated by the examples in 
Table~\ref{tab:B1-nilpotent}---never propagates to 
$\B_1(H^*(A))$.

%%%%%%%%%%%%%%%%%%%%%%%%%%%%%%%%%%%
\section{The first Koszul module and its presentation}
\label{sect:BAd}
%%%%%%%%%%%%%%%%%%%%%%%%%%%%%%%%%%%

In this section we focus on the first Koszul module $\B_1(A)=H_1(K_\bullet(A))$ 
of a connected $\cdga$ $(A,d_A)$ with $\dim A^1<\infty$.  
This module plays a distinguished role: it admits an explicit presentation
in terms of the linear and quadratic data of $A$, fits into an infinitesimal
analogue of Crowell’s exact sequence from group theory, and governs 
the degree-$1$ resonance varieties of $A$. These features make $\B_1(A)$ 
the basic algebraic invariant connecting Koszul homology, resonance, and 
holonomy Lie algebras. In particular, the results of this section provide the 
algebraic backbone for the identification of $\B_1(A)$ with the infinitesimal 
Alexander invariant of the holonomy Lie algebra and for the extraction of 
holonomy Chen ranks from Koszul homology modules.

%%%%%%%%%%%%%%%%%%%%%%
\subsection{A presentation for the $S$-module $\B_1(A)$}
\label{subsec:Bpres}
%%%%%%%%%%%%%%%%%%%%%%

We begin by giving an explicit presentation of $\B_1(A)$ as a module over
$S=\Sym(H_1(A))$, making transparent its dependence on the linear differential
$d_A$ and the quadratic multiplication on $A$. This presentation generalizes 
the one given in \cite{PS-imrn04, PS-crelle} in the case when $d_A=0$. 
Unlike the quadratic case, this presentation separates the contributions 
of the linear differential and the quadratic multiplication, making explicit 
how non-formal features enter $\B_1(A)$ through their interaction.
This presentation will serve as the starting point for the comparison with
holonomy Lie algebras and for the analysis of resonance in degree~$1$.

As before, set $E=\bwedge (H^1(A))$, and write $A_i=(A^i)^{\vee}$ and 
$E_i=(E^i)^{\vee}$.  Denote by
$d^1=d_A^1\colon A^1\to A^2$ the degree-one differential and by
$d_1=(d_A^1)^\vee\colon A_2\to A_1$ its $\k$-linear dual.

Since $A$ is connected, we may identify $E^1=H^1(A)$ with 
$Z^1=\ker\big(d_A\colon A^1\to A^2\big)$, and we set
\begin{equation}
\label{eq:U1-def}
U_1=\im(d_1)\subseteq A_1 .
\end{equation}
Both $E^1$ and $U_1$ are canonically attached to $A$, and they annihilate 
each other: since $\langle d_1(\gamma),a\rangle=\langle \gamma,d_A a\rangle$ 
for $\gamma\in A_2$ and $a\in A^1$, we have 
\begin{equation}
\label{eq:perp-pairs}
(U_1)^{\perp}=Z^1=E^1, \qquad U_1=(E^1)^{\perp} .
\end{equation}
Splitting off these subspaces, however, requires a choice. Fix a complement 
$U^1$ of $Z^1$ in $A^1$, and set $E_1\coloneqq (U^1)^{\perp}\subseteq A_1$.  
Then $E_1\cap U_1=(U^1)^{\perp}\cap (E^1)^{\perp}=(U^1+E^1)^{\perp}=0$, 
while $\dim E_1+\dim U_1=\dim A^1$; hence we obtain mutually dual 
$\k$-vector space decompositions
\begin{equation}
\label{eq:A1A1}
A_1 = E_1\oplus U_1,\qquad 
A^1 = E^1\oplus U^1,
\end{equation}
along with an identification $E_1\cong (E^1)^{\vee}=H^1(A)^{\vee}=H_1(A)$. 
We denote by $\pi_E\colon A_1\to E_1$ and $\pi_U\colon A_1\to U_1$ 
the projections onto the respective summands.  The decompositions 
\eqref{eq:A1A1} depend on the choice of $U^1$ and are \emph{not}\/ 
functorial; the effect of this on the naturality of the presentation 
below is recorded in Remark~\ref{rem:bphi-triangular}.

Since $A$ is graded-commutative, the multiplication map
$A^1\otimes_{\k} A^1 \to A^2$ descends to a $\k$-linear map, 
$\mu_A\colon A^1\wedge A^1 \to A^2$; we let $\mu_A^{\vee}\colon 
A_2 \to A_1\wedge A_1$ be the dual map.  We also let 
$\nu_A\colon E^1\wedge E^1 \to A^2$ be the 
restriction of $\mu_A$ to the summand 
$E^1\wedge E^1=E^2$, and denote by 
\begin{equation}
\label{eq:nu-vee}
\begin{tikzcd}[column sep=22pt]
\nu_A^{\vee}=(\pi_E\wedge \pi_E)\circ \mu_A^{\vee} 
\colon A_2 \arrow[r] & E_1\wedge E_1=E_2
\end{tikzcd}
\end{equation}
the $\k$-dual map.

Now let $S=\Sym(E_1)$ be the symmetric algebra on the (finite-dimensional) 
$\k$-vector space $E_1$. As noted previously, left-multiplication by the 
canonical element $\omega_A\in Z^1(A) \otimes_{\k} H_1(A)$ 
defines $S$-linear maps $A^i\otimes_{\k} S\to A^{i+1} \otimes_{\k} S$. 
Likewise, multiplication by the canonical element 
$\omega_E\in E^1\otimes_{\k} E_1$ defines $S$-linear maps 
$E^i\otimes_{\k} S\to E^{i+1} \otimes_{\k} S$. 
With these definitions, it is readily seen that 
diagram \eqref{eq:omega-ae} below commutes. 
\begin{equation}
\label{eq:omega-ae}
\begin{tikzcd}
A^1\otimes_{\k} S \ar[r, "\omega_A"] 
&
A^2 \otimes_{\k} S \phantom{.}
\\
E^1\otimes_{\k} S \ar[r, "\omega_E"] 
\arrow[hookrightarrow]{u}[swap]{\pi_E^{\vee}\otimes \id_{S}} 
&
E^2 \otimes_{\k} S . \ar[u, "\nu_A\otimes \id_S" ']
\end{tikzcd}
\end{equation}
The diagram expresses the compatibility of the canonical elements 
$\omega_A$ and $\omega_E$ under the morphisms $\pi_E^{\vee}$ 
and $\nu_A$. 
Taking duals, this means that $(\pi_E\otimes \id_S) \circ \omega_A^{\vee}= 
\omega_E^{\vee}\circ (\nu_A^{\vee} \otimes \id_S)$.  Finally, consider the composite 
$\beta_A\colon U^1\otimes_{\k} S\to A^1 \otimes_{\k} S \xrightarrow{\, \omega_A\,} 
A^2 \otimes_{\k} S$, that sends an element $u\in U^1$ to $\sum_{i=1}^n u e_i \otimes x_i$ 
via the canonical tensor $\omega_A$, and let 
\begin{equation}
\label{eq:beta-vee}
\beta_A^{\vee}=(\pi_U\otimes \id_S) \circ \omega_A^{\vee}\colon 
A_2\otimes_{\k} S\longrightarrow U_1\otimes_{\k} S
\end{equation}
be its $S$-dual. 

With notation in place, we are ready to state and prove the main result of this
section.  Modules are left modules throughout, so that a morphism of free
modules acts on the right of row vectors; accordingly, in the block matrices
displayed below, here and in \eqref{eq:bphi}, the rows are indexed by the
summands of the source and the columns by those of the target. 

\begin{theorem}
\label{thm:Bpres}
Let $(A,d_A)$ be a connected $\k$-$\cdga$ with $A^1$ finite-dimensional. 
Then the Koszul module $\B_1(A)$ has presentation as an $S$-module
\begin{equation}
\label{eq:Bmod-pres}
\begin{tikzcd}[column sep=20pt]
\big(E_3\oplus A_2\big) \otimes_{\k} S
\ar[r, "
\sbm{\partial_3^{E}  \amp 0 \\  
\nu_A^{\vee} \otimes\id_S\amp d_A^{\vee} \otimes\id_S+ \beta_A^{\vee}}
"]&[65pt] 
\big( E_2\oplus U_1\big)  \otimes_{\k} S
\ar[r]& \B_1(A)  \ar[r]& 0 .
\end{tikzcd}
\end{equation} 
\end{theorem}

\begin{proof}
By definition, $\B_1(A)$ is the first homology of 
the chain complex of free $S$-modules 
$K_{\bullet}(A)=(A_{\bullet}\otimes_{\k} S,\partial^A)$. 
As noted previously, the Koszul complex 
$K_{\bullet}(E)=(E_{\bullet}\otimes_{\k} S,\partial^E)$ is exact; 
in particular, $\ker(\partial^E_1)\cong\coker( \partial^E_3)$. 

In order to find a presentation for the $S$-module $\B_1(A)$, 
recall that $\partial^A_2=d_A^{\vee}\otimes \id_S +\omega_A^{\vee}$ 
and $\partial^A_1=\omega_A^{\vee}$, whereas 
$\partial^E_2=\omega_{E}^{\vee}$ and 
$\partial^E_1=\omega_{E}^{\vee}$. Consider 
now the following diagram, 
\begin{equation*}
\label{eq:chain bgg}
\begin{tikzcd}[column sep=26pt, row sep=38pt]
& 
& A_2  \otimes_{\k}   S  \ar[r, "\partial^A_2"]   
\ar[d, " \sbm{\nu_A^{\vee}\\[2pt]  d_A^{\vee}\otimes \id_S + \beta_A^{\vee}}" {yshift=5pt}]
&[6pt] A_1  \otimes_{\k}   S  \ar[r, "\partial^A_1"] \ar[d, equal]
&[6pt] A_0  \otimes_{\k}   S=S \phantom{.} \ar[d, equal] \\
& E_3  \otimes_{\k} S  \ar[r, "\sbm{\partial^E_3\\ 0}"]
&  (E_2 \oplus U_1) \otimes_{\k} S    \ar[r, "\sbm{\partial^E_2 \amp 0\\0\amp \id}"]
&  (E_1 \oplus U_1)  \otimes_{\k} S   \ar[r, "\sbm{ \partial^E_1 \amp 0}"]
&  E_0  \otimes_{\k} S =S .
\end{tikzcd}
\end{equation*}

The chain complex on the bottom is exact in degrees $1$ and $2$
(since it is chain-homotopic to $(E_{\le 3}\otimes_{\k} S,\partial^E)$).
The vertical arrows define a 
chain map between the two complexes. Indeed, using 
\eqref{eq:omega-ae} and the fact that 
$\pi_U\circ d_A^{\vee} = d_A^{\vee}$, we get 
\begin{align*}
\partial_2^E\circ (\nu_A^{\vee}\otimes \id_S)& 
=\omega_E^{\vee} \circ (\nu_A^{\vee}\otimes \id_S)
= (\pi_E\otimes \id_S) \circ \omega_A^{\vee}
= (\pi_E\otimes \id_S) \circ  \partial_2^A,
\\
d_A^{\vee} \otimes \id_S +\, \beta_A^{\vee}&= 
(\pi_U \otimes \id_S) \circ  (d_A^{\vee} \otimes \id_S) 
+ (\pi_U\otimes \id_S)\circ \omega_A^{\vee}=
(\pi_U\otimes \id_S)\circ  \partial_2^A .
\end{align*}

Under the identification $A_1=E_1\oplus U_1$, the map 
$\partial^A_1=\omega_A^\vee$ agrees with 
$\partial^E_1\circ \pi_E$, and thus 
$\ker(\partial^A_1)=\ker(\partial^E_1\circ \pi_E)$. 
Using the commutativity of the diagram in degree $2$, together with 
the exactness of the bottom complex in degrees $1$ and $2$, 
it follows that the quotient 
$\B_1(A)=\ker (\partial^A_1)/\im (\partial^A_2)$ 
is the cokernel of the map in \eqref{eq:Bmod-pres}.
\end{proof}

In particular, the presentation of $\B_1(A)$ depends only on the linear differential 
$d_A^1$, the induced quadratic multiplication on $H^1(A)$, and their interaction 
via the canonical tensor $\omega_A$ (through the map $\beta_A$). 
No higher-degree information from $A$ enters into the structure of 
$\B_1(A)$ as an $S$-module.

\begin{remark}
\label{rem:B-zero-d}
When $d_A=0$, the presentation \eqref{eq:Bmod-pres}  recovers the 
Koszul-type presentation of the infinitesimal Alexander invariant introduced 
by Papadima and Suciu in \cite[Thm.~6.2]{PS-imrn04} and later formalized  
in \cite[Def.~2.1]{PS-crelle}. Indeed, in this case $d_A^{\vee}=0$ and $U_1=0$ 
(hence $\beta_A^{\vee}=0$), while $\nu_A^{\vee}=\mu_A^{\vee}$. Therefore, 
\begin{equation}
\label{eq:b1a-mua}
\B_1(A)=\coker \left( \begin{smallmatrix}  \partial_E^3 \\ 
\mu_A^{\vee} \end{smallmatrix}\right).
\end{equation}
Hence, Theorem~\ref{thm:Bpres} 
recovers the aforementioned infinitesimal Alexander invariant  
as a special case, while simultaneously identifying the precise 
correction terms arising from a nontrivial differential.
In contrast to \cite{PS-imrn04, PS-crelle}, where this presentation 
serves as the starting point for the theory, here it is obtained 
as a consequence of the intrinsic definition of $\B_1(A)$ via 
homology of the Koszul complex $K_{\bullet}(A)$. 
\end{remark}

When $d_A=0$ and $\mu_A$ is onto, the presentation \eqref{eq:b1a-mua}
determines the space of minimal generators of $\B_1(A)$ completely, in terms
of the single subspace $K=\ker(\mu_A)$.

\begin{lemma}
\label{lem:b1-degree-zero}
Let $A$ be a connected $\k$-$\cdga$ with $d_A=0$ and $A^1$ finite-dimensional,
and suppose the cup product $\mu_A\colon\bwedge^2A^1\to A^2$ is surjective.
Set $K=\ker(\mu_A)\subseteq\bwedge^2A^1$.  Then there is a natural
isomorphism
\[
\B_1(A)/\m\B_1(A) \;\cong\; K^{\vee} .
\]
In particular, $\B_1(A)$ is minimally generated by $\dim_\k K$ elements.
\end{lemma}

\begin{proof}
Write $V=A_1=H_1(A)$ (using $d_A=0$ to identify $A^1=H^1(A)$), and let
$E=\bwedge(A^1)$ as before, so that $E_2=\bwedge^2 V$ canonically.  Since
$\mu_A$ is onto, its $\k$-dual $\mu_A^{\vee}\colon A_2\to E_2$ is injective,
with image $K^{\perp}$, the annihilator of $K$ under the canonical pairing
$\bwedge^2A^1\otimes_{\k}\bwedge^2V\to\k$; thus $A_2\cong K^{\perp}\subseteq E_2$,
and $E_2/K^{\perp}\cong K^{\vee}$ canonically, via restriction of linear
functionals.

By Remark~\ref{rem:B-zero-d}, $\B_1(A)=\coker(\partial_3^E;\mu_A^{\vee})$,
where $(\partial_3^E;\mu_A^{\vee})\colon E_3\otimes_{\k} S\oplus A_2\otimes_{\k} S
\to E_2\otimes_{\k} S$.  By definition, $\B_1(A)/\m\B_1(A)$ is the further
quotient of this cokernel by the submodule generated by $S_{\ge1}\cdot E_2$;
equivalently, it is the cokernel of the composite map to $E_2\otimes_{\k} S_0=E_2$
obtained by retaining only the internal-degree-$0$ part of the target and
the parts of the source landing there.  The map $\partial_3^E$, 
like $\partial^A$ in \eqref{eq:KbulletA}, is given by
contraction with a canonical class $\omega_E\in E^1\otimes_{\k} E_1$, and hence
raises internal $S$-degree by $1$; as its source $E_3\otimes_{\k} S$ therefore
contributes nothing in internal degree $0$, only the
map $\mu_A^{\vee}\colon A_2\to E_2$ survives at that degree, giving
\[
\B_1(A)/\m\B_1(A)\cong E_2/\im(\mu_A^{\vee})=E_2/K^{\perp}\cong K^{\vee} .
\qedhere
\]
\end{proof}

\begin{remark}
\label{rem:b1-degree-zero-shift}
In the internal-degree convention of Lemma~\ref{lem:filtration-shift-B1},
$\B_1(A)$ is generated in internal degree $1$, and
Lemma~\ref{lem:b1-degree-zero} identifies this generating piece; the shift
by one accounts for the discrepancy between $\B_1(A)/\m\B_1(A)$
\textup{(}internal degree $1$\textup{)} and what
Section~\ref{subsec:b1-conf-structure} below calls ``degree $0$'' of
$\B_1(H^{*}(\Conf(\E,n);\k))$, which is the standard reindexing making the
module generated starting in degree $0$.
\end{remark}

The presentation given in Theorem~\ref{thm:Bpres} is natural, in the 
following sense. Suppose $\varphi \colon A \to \oA$ is a morphism 
of connected $\k$-$\cdgas$ such that $H^1(A)$ and $H^1(\oA)$ are 
finite-dimensional. 
Let $\B_{1}(\varphi)\colon \B_{1}(\oA) \otimes_{\oS} S \to \B_{1}(A)$ 
be the induced morphism of $S$-modules from \eqref{eq:bia-s-bia}, and write 
$\varphi_i=(\varphi^i)^{\vee}$. Decomposing the restriction of $\varphi_1$ 
to $\oE_1$ according to \eqref{eq:A1A1}, we set 
\begin{equation}
\label{eq:phi-E-U}
\varphi_1^{E}=\pi_E\circ \varphi_1|_{\oE_1}\colon \oE_1\to E_1,
\qquad 
\varphi_1^{U}=\pi_U\circ \varphi_1|_{\oE_1}\colon \oE_1\to U_1 ,
\end{equation}
so that $\varphi_1|_{\oE_1}=\varphi_1^{E}+\varphi_1^{U}$, where 
$\varphi_1^{E}=H_1(\varphi)$ induces the base-change map 
$\Sym(\varphi_1^{E})\colon \oS\to S$.  On the other summand, no such 
correction occurs: dualizing $\varphi^2\circ d_A=d_{\oA}\circ \varphi^1$ 
gives $d_A^{\vee}\circ \varphi_2=\varphi_1\circ d_{\oA}^{\vee}$, whence
\begin{equation}
\label{eq:phi-U-canonical}
\varphi_1(\oU_1)\subseteq U_1 .
\end{equation}
Thus, the matrix of $\varphi_1$ with respect to \eqref{eq:A1A1} is 
triangular, with off-diagonal entry $\varphi_1^{U}$. Accordingly, set 
\begin{equation}
\label{eq:theta-phi}
\theta_{\varphi}=(\varphi_1^{U}\otimes \id_S)\circ \partial_2^{\oE}
\colon \oE_2\otimes_{\k} S \longrightarrow U_1\otimes_{\k} S .
\end{equation}
We then have the following commuting diagram of $S$-module 
presentations, where we suppress $\otimes\id_S$, $\otimes\id_{\oS}$, and 
the base change along $\Sym(\varphi_1^E)$ from the notation, for brevity.
\begin{equation}
\label{eq:bphi}
\begin{tikzcd}[column sep=24pt,row sep=36pt]
\big(\oE_3\oplus \oA_2\big) \otimes_{\k} S
\ar[r, "\sbm{\partial_3^{\oE}  \amp 0 \\  
\nu_{\oA}^{\vee} \amp d_{\oA}^{\vee} + \beta_{\oA}^{\vee}}"]
\ar[d, pos=0.4, "\bwedge^3\varphi_1^{E} \oplus 
\varphi_2 "]  &[40pt] 
\big( \oE_2\oplus \oU_1\big)  \otimes_{\k} S\ar[r] 
\ar[d, pos=0.4, "\sbm{\bwedge^2\varphi_1^{E} \amp \theta_{\varphi}\\ 
0 \amp \varphi_1}" ] 
& \B_{1}(\oA) \otimes_{\oS} S \ar[r] \ar[d, pos=0.4, "\B_1(\varphi)"]
& 0\phantom{\,.}
\\
\big(E_3\oplus A_2\big) \otimes_{\k} S
\ar[r, "\sbm{\partial_3^{E}  \amp 0 \\  
\nu_{A}^{\vee} \amp d_{A}^{\vee} + 
\beta_{A}^{\vee}}"]
&[40pt] 
\big( E_2\oplus U_1\big)  
\otimes_{\k} S\ar[r] 
& \B_1(A) \ar[r] & 0 .
\end{tikzcd}
\end{equation}
Recall that the rows of the displayed matrices are indexed 
by the summands of the source and the columns by those of the target; thus, 
the middle vertical arrow sends $\bar{e}\wedge \bar{f}\otimes 1$ to 
$\big(\varphi_1^{E}(\bar e)\wedge \varphi_1^{E}(\bar f)\otimes 1,\, 
\theta_{\varphi}(\bar e\wedge \bar f)\big)$ and $\bar{u}\otimes 1$ to 
$\big(0,\varphi_1(\bar u)\otimes 1\big)$, for 
$\bar e,\bar f\in \oE_1$ and $\bar u \in \oU_1$.

Commutativity of the right-hand square is the definition of 
$\B_1(\varphi)$.  For the left-hand square, note first that 
\eqref{eq:phi-U-canonical} gives 
$\pi_E\circ \varphi_1=\varphi_1^{E}\circ \pi_{\oE}$, so multiplicativity 
of $\varphi$, in the form 
$(\varphi_1\wedge \varphi_1)\circ \mu_{\oA}^{\vee}=
\mu_A^{\vee}\circ \varphi_2$, yields 
\[
\nu_A^{\vee}\circ \varphi_2=
(\pi_E\wedge \pi_E)\circ \mu_A^{\vee}\circ \varphi_2=
\big(\varphi_1^{E}\wedge \varphi_1^{E}\big)\circ \nu_{\oA}^{\vee} .
\]
Next, \eqref{eq:phi-U-canonical} also gives 
$\pi_U\circ \varphi_1-\varphi_1\circ \pi_{\oU}=
\varphi_1^{U}\circ \pi_{\oE}$.  Since 
$d_A^{\vee}\otimes \id_S+\beta_A^{\vee}=
(\pi_U\otimes \id_S)\circ \partial_2^A$ and $\varphi$ induces a chain map 
of Koszul complexes, we obtain 
\[
\big(d_A^{\vee}+\beta_A^{\vee}\big)\circ \varphi_2=
(\pi_U\otimes \id_S)\circ \varphi_1\circ \partial_2^{\oA}=
\varphi_1\circ \big(d_{\oA}^{\vee}+\beta_{\oA}^{\vee}\big)
+\big(\varphi_1^{U}\otimes \id_S\big)\circ 
(\pi_{\oE}\otimes \id_S)\circ \partial_2^{\oA} ,
\]
and the last term equals $\theta_{\varphi}\circ \nu_{\oA}^{\vee}$, by the 
first identity established in the proof of Theorem~\ref{thm:Bpres}. 
Together, these two computations say precisely that the left square 
commutes on the summand $\oA_2$.  On the summand $\oE_3$, the required 
identities are the functoriality of the Koszul differential, 
$\bwedge^2\varphi_1^E\circ \partial_3^{\oE}=
\partial_3^{E}\circ \bwedge^3 \varphi_1^E$, together with 
$\theta_{\varphi}\circ \partial_3^{\oE}=
(\varphi_1^{U}\otimes \id_S)\circ \partial_2^{\oE}\circ 
\partial_3^{\oE}=0$.

\begin{remark}
\label{rem:bphi-triangular}
The vertical maps in \eqref{eq:bphi} are triangular rather than block 
diagonal, and the correction term $\theta_{\varphi}$ cannot be dispensed 
with in general: it vanishes precisely when $\varphi_1^{U}=0$, that is, 
when $\varphi_1(\oE_1)\subseteq E_1$, or equivalently when 
$\varphi^1(U^1)\subseteq \oU^1$ (the case $\oE_2=0$ being trivial).  As the 
splittings \eqref{eq:A1A1} involve a choice, this compatibility can always 
be arranged for a single morphism, but not simultaneously for all 
morphisms.  In particular, $\theta_{\varphi}=0$ whenever 
$d_A=0$ and $d_{\oA}=0$, since then $U_1=0$ and $\oU_1=0$; in that case 
\eqref{eq:bphi} reduces to the block-diagonal naturality of the 
presentation of \cite{PS-imrn04, PS-crelle}.
\end{remark}

This naturality will be important in later comparisons with holonomy Lie algebras 
and group-level invariants, where $\cdga$ morphisms encode functoriality of 
Malcev and Chen constructions.

%%%%%%%%%%%%%%%%%%%%
\subsection{A Crowell-type exact sequence and degree $1$ resonance}
\label{subsec:crowell}
%%%%%%%%%%%%%%%%%%%%

The presentation above relates the Koszul module $\B_1(A)$ to the 
Alexander-type module $\fA(A)$ and leads naturally to an infinitesimal 
analogue of the classical Crowell exact sequence for Alexander 
invariants of groups \cite{Cr65}. This generalizes the construction from 
\cite[Thm.~15.1]{Su-pisa}, where the case $d_A=0$ was treated. 
Moreover, the sequence is functorial in $A$.

\begin{proposition}
\label{prop:inf-crowell}
Let $(A,d_A)$ be a connected $\k$-$\cdga$ with $\dim A^1<\infty$.
Let $\fA(A)=\coker (\partial_2^A)$ and
let $\m$ be the augmentation ideal of $S=\Sym(H_1(A))$.  Then there is a
natural exact sequence of $S$-modules
\[
\begin{tikzcd}[column sep=20pt]
0\ar[r]& \B_1(A)\ar[r]&\fA(A)\ar[r]&\m \ar[r]& 0  .
\end{tikzcd}
\]
\end{proposition}

\begin{proof}
Using the decomposition $A_1 = E_1 \oplus U_1$ from \eqref{eq:A1A1}, we have
\[
A_1\otimes_{\k} S = (E_1\otimes_{\k} S) \oplus (U_1\otimes_{\k} S), \qquad
\ker(\partial^A_1) \cong \ker(\partial^E_1) \oplus (U_1\otimes_{\k} S).
\]
It follows that
\[
\fA(A)/\B_1(A) \cong (E_1\otimes_{\k} S) / \ker(\partial^E_1) 
\cong \im(\partial^E_1) \cong \m,
\]
giving the exact sequence. Naturality follows as in the discussion of $\B_1(A)$.
\end{proof}

The degree-$1$ resonance varieties $\RR^{1,s}(A)$ are determined by 
the differential $d_A\colon A^1\to A^2$ and the multiplication 
$\mu_A\colon A^1\wedge A^1\to A^2$, or equivalently by their 
corestrictions to $\im d_A$ and $\im \mu_A$. These varieties 
are tightly related to the support loci 
$\RR_{1,s}(A) =\supp \bigl( \bwedge^s  \B_1(A)\bigr)$, as follows. 

\begin{proposition}
\label{prop:r1a}
Let $(A,d_A)$ be a connected $\cdga$ with $\dim A^1<\infty$. Then
\[
\RR^{1,s}(A)\setminus\{0\} = \RR_{1,s}(A)\setminus\{0\}
\qquad\text{for all } s\ge 1 .
\]
\end{proposition}

\begin{proof}
Immediate from 
Theorem~\ref{thm:jump-support-containment}\eqref{jsc2}, since 
$\RR_{0,1}(A)=\{0\}$.
\end{proof}

By Example~\ref{ex:free-group-depth}, the origin must be excluded here for
every $s\ge1$, and not merely for $s=1$.

\begin{corollary}
\label{cor:res-fitt}
Let $(A,d_A)$ be a connected $\k$-$\cdga$ with $\dim A^1<\infty$.
Then, for all $s\ge 1$,
\[
\RR^{1,s}(A)\setminus\{0\} = \bV\big(\Fitt_{s}(\fA(A))\big) \setminus\{0\}.
\]
\end{corollary}

\begin{proof}
Consider the Crowell exact sequence
$0 \to \B_1(A)  \to \fA(A)  \to  \m  \to  0$, 
with $\m \subset S$ the augmentation ideal.  Localizing at a maximal 
ideal $\m_\xi \ne \m$ gives $\m_{\m_\xi} = S_{\m_\xi}$, so the localized 
sequence has free quotient and therefore splits:
\[
\fA(A)_{\m_\xi}\cong \B_1(A)_{\m_\xi}\oplus S_{\m_\xi} .
\]
Hence the local numbers of generators satisfy
$\mu_\xi(\fA(A))=\mu_\xi(\B_1(A))+1$, and Lemma~\ref{lem:supp-mu} converts 
this into a shift of supports,
\[
\Supp\big(\bwedge^{s+1} \fA(A)\big) \setminus \{0\} = 
\Supp\big(\bwedge^{s} \B_1(A)\big) \setminus \{0\}.
\]
By Lemma~\ref{lem:ann-supp}\eqref{ann2} the left-hand side is
$\bV\big(\Fitt_{s}(\fA(A))\big)\setminus\{0\}$, while the right-hand side is
$\RR_{1,s}(A)\setminus\{0\}$; the claim now follows from
Proposition~\ref{prop:r1a}. 
\end{proof}

\begin{remark}
\label{rem:alt-aomoto}
Alternatively, let $a\in H^1(A)$. By \eqref{eq:point-res}, 
$a\in \RR^{1,s}(A)$ if and only if 
\[
\rank \partial_1^A(a) + \rank \partial_2^A(a) \le \dim_{\k} A^1 - s.
\]
Since $A_0 = \k$, the map $\partial_1^A(a)\colon A_1 \to \k$ has rank $0$
if $a=0$ and rank $1$ if $a\neq 0$. Thus, for $a\neq 0$, the condition 
reduces to $\rank \partial_2^A(a) \le \dim_{\k} A^1 - s - 1$. 
By definition, $\fA(A) = \coker(\partial_2^A)$, a module presented by a
matrix with target of rank $\dim_\k A^1$, so
$\Fitt_{s}(\fA(A))$ is generated by the
$(\dim_{\k} A^1 - s)$-minors of $\partial_2^A$.  These vanish precisely
when $\rank \partial_2^A(a) \le \dim_{\k} A^1 - s - 1$, which is the stated
condition.  This gives the equality away from $0$.
\end{remark}

\begin{remark}
\label{rem:r1s-vann}
Away from $0$, degree-$1$ resonance is completely governed by $\B_1(A)$:
\[
\RR^{1,s}(A)\setminus\{0\} = \bV\bigl(\Ann(\bwedge^s \B_1(A))\bigr)\setminus\{0\},
\]
so the Crowell exact sequence provides a very concrete algebraic model for these 
varieties. Even for $s=1$, the scheme structure at $0$ can differ from the 
support locus; see \cite[Ex.~5.3]{AFRSS24}.
\end{remark}

%%%%%%%%%%%%%%%%%%%%%%%%%%
\section{The holonomy Lie algebra of a $\cdga$}
\label{sect:holo-lie}
%%%%%%%%%%%%%%%%%%%%%%%%%%

In this section, we develop the structure theory of holonomy Lie
algebras associated to connected $\cdgas$, extending the classical
quadratic framework to the non-formal setting.
In particular, we express the universal enveloping algebra $U(\h(A))$
in terms of functorial constructions on quadratic and quadratic--linear
algebras, and analyze its associated graded algebra via Koszul duality
and the Poincar\'{e}--Birkhoff--Witt theorem.
We also show that the Lie algebra $\h(A)$ functorially determines the
$1$-minimal model of $(A,d)$, and thus encodes its infinitesimal
homotopy invariants.

%%%%%%%%%%%%%%%%%
\subsection{Holonomy Lie algebra}
\label{subsec:holonomy-def}
%%%%%%%%%%%%%%%%%

In what follows, $\Lie(V)$ will denote the free Lie algebra on a finite-dimensional 
$\k$-vector space~$V$, graded by bracket length. We identify $\Lie^1(V)=V$ and 
$\Lie^2(V)=V\wedge V$ via the standard correspondence
$[v,w] \leftrightarrow v\wedge w$. 

Let $(A,d_A)$ be a connected $\k$-$\cdga$ with $\dim_\k A^1<\infty$.
Set $A_i=(A^i)^\vee$, and denote by
\begin{equation}
\label{eq:muA-dA}
\mu_A \colon  A^1\wedge A^1  \longrightarrow A^2
\quad\text{and}\quad
d_A \colon  A^1  \longrightarrow A^2
\end{equation}
the multiplication and differential, respectively.
Their linear duals combine into a single map
\begin{equation}
\label{eq:nablaA}
\nabla_A \coloneqq d_A^\vee + \mu_A^\vee
\colon A_2 \longrightarrow A_1 \oplus (A_1\wedge A_1) ,
\end{equation}
which simultaneously encodes the linear and quadratic relations
defining $\h(A)$, yielding a quadratic--linear presentation. 
Following~\cite{MPPS,PS-jlms}, we make the following definition. 

\begin{definition}
\label{def:holonomy-lie}
The \emph{holonomy Lie algebra} of $A$ is the quotient
\[
\h(A) \coloneqq \Lie(A_1)\big/\big\langle\! \im(\nabla_A)\big\rangle,
\]
where $\langle \im(\nabla_A)\rangle$ is the Lie ideal generated by the image
of $\nabla_A$, viewed as a subspace of $\Lie(A_1)$ via the natural inclusions
$A_1 \hookrightarrow \Lie(A_1)$ and $A_1\wedge A_1 \hookrightarrow \Lie^2(A_1)$.
\end{definition}

By construction, $\h(A)$ is a finitely presented Lie algebra with generators
in degree~$1$ and relations in degrees~$1$ and~$2$.
The linear relations come from the image of $d_A^\vee$, while the quadratic
relations are dual to the multiplication in~$A$. Its abelianization satisfies
\[
\h(A)/\h'(A)\cong H_1(A),
\]
since the linear relations imposed by $d_A^\vee$ identify
$\h(A)_{\ab}$ with $\ker(d_A)^\vee = (H^1(A))^{\vee}=H_1(A)$. 

The free Lie algebra $\Lie(A_1)$ is equipped with the descending filtration by
bracket length,
\[
F^p \Lie(A_1) = \bigoplus_{n \ge p} \Lie^n(A_1).
\]
Since $\im(d_A^\vee) \subset A_1$ and $\im(\mu_A^\vee) \subset A_1\wedge A_1$,
the ideal $\langle \im(\nabla_A)\rangle$ is compatible with this filtration.
Thus, $\h(A)$ inherits a natural descending filtration
\[
F^p \h(A) = \im\big(F^p \Lie(A_1) \longrightarrow \h(A)\big),
\]
making $\h(A)$ into a filtered Lie algebra.
Passing to associated graded Lie algebras, the linear relations coming from
$\im(d_A^\vee)$ vanish, while the quadratic relations coming from
$\im(\mu_A^\vee)$ remain. Thus,
\begin{equation}
\label{eq:grhA-lie}
\gr(\h(A)) \cong \Lie(A_1)/\langle \im(\mu_A^\vee)\rangle .
\end{equation}

If $d_A=0$, then $\nabla_A=\mu_A^\vee$, and one recovers Chen’s classical
holonomy Lie algebra~\cite{Chen73}, which is quadratic and graded by bracket length.
In general, $\h(A)$ need not admit such a grading, but only the filtration above.

The construction is contravariantly functorial:
a $\cdga$ morphism $\varphi\colon A\to A'$ induces 
a Lie algebra morphism in the opposite direction
\[
\h(\varphi)\colon  \h(A') \longrightarrow \h(A).
\]

\begin{proposition}
\label{prop:holo-CE}
For any finite-dimensional Lie algebra $\g$ over $\k$, there is a
canonical isomorphism $\h(\CE(\g)) \cong \g$.
\end{proposition}

\begin{proof}
Set $A = \CE(\g) = (\bigwedge \g^\vee,\, d_{\CE})$.
Since $\g$ is finite-dimensional, $A_1 = (A^1)^\vee = (\g^\vee)^\vee \cong \g$
and $A_2 = (A^2)^\vee \cong \bigwedge^2\g$.
Under these identifications, we compute $\nabla_A = d_A^\vee + \mu_A^\vee$
on the basis elements $x \wedge y \in A_2$.

For $\xi \in \g^\vee$, the Chevalley--Eilenberg differential satisfies
$(d_{\CE}\xi)(x,y) = -\xi([x,y]_\g)$, so the dual map gives
$d_A^\vee(x\wedge y) = -[x,y]_\g \in \g = A_1$.
The cup-product dual satisfies 
$\mu_A^\vee(x\wedge y) = x\wedge y \in A_1\wedge A_1$, 
which corresponds to $[x,y]_{\Lie(\g)}$ in $\Lie^2(A_1)$.
Hence
\[
\nabla_A(x\wedge y) = [x,y]_{\Lie(\g)} - [x,y]_\g,
\]
and $\langle\im(\nabla_A)\rangle$ is the ideal in $\Lie(\g)$ 
generated by $\{[x,y]_{\Lie(\g)}-[x,y]_\g : x,y\in\g\}$.
This is precisely the kernel of the canonical surjection
$\pi\colon \Lie(\g)\surj\g$ induced by $\id_\g$,
so the first isomorphism theorem gives
$\h(A) = \Lie(\g)/\langle\im(\nabla_A)\rangle \cong \g$.
\end{proof}

%%%%%%%%%%%%%%%%%
\subsection{Holonomy Lie algebras under weight conditions}
\label{subsec:holo-weights}
%%%%%%%%%%%%%%%%%

In the presence of weight structures, the holonomy Lie algebra behaves in a 
more rigid fashion than in the general qua\-dratic--linear setting. In particular, 
positivity of weights allows the linear and quadratic components of the defining 
relations to be separated, leading to a compatible grading on $\h(A)$. Moreover, 
under suitable weight restrictions on the differential, holonomy is stable under 
passage from a $\cdga$ to its cohomology algebra in low weights. The present 
subsection records these two complementary phenomena.

\begin{proposition}
\label{prop:pos-weights-graded-holo}
Let $(A,d_A)$ be a connected $1$-finite $\cdga$ with positive weights.
Then the holonomy Lie algebra $\h(A)$ is graded by weight,
with grading induced from the positive weight decompositions of $A^1$ 
and $A^2$ via the map $\nabla_A = d_A^\vee + \mu_A^\vee$.
\end{proposition}

\begin{proof}
Since $(A,d_A)$ has positive weights, the spaces $A^1$ and $A^2$ 
decompose as $A^1 = \bigoplus_{\alpha>0} A^{1,\alpha}$ and 
$A^2 = \bigoplus_{\alpha>0} A^{2,\alpha}$, with both $\mu_A$ and 
$d_A$ weight-homogeneous. Dualizing, both $d_A^\vee$ and $\mu_A^\vee$ 
are homogeneous with respect to the weight grading, and hence 
so is their sum $\nabla_A$. Explicitly, for each weight $\alpha$,
\[
\nabla_A\colon A_{2,\alpha} \longrightarrow 
A_{1,\alpha} \oplus 
\boplus_{\beta+\gamma=\alpha} A_{1,\beta}\wedge A_{1,\gamma}.
\]
Hence $\im(\nabla_A)$ is a weight-homogeneous subspace of 
$\Lie(A_1)$, and the quotient 
$\h(A) = \Lie(A_1)/\langle\im(\nabla_A)\rangle$ 
inherits a well-defined grading by weight.
\end{proof}

This grading reflects the fact that, under positivity assumptions, 
the quadratic and linear components of the relations defining $\h(A)$ are compatible 
with the weight decomposition, so that no mixing between different weights occurs 
in the quotient.

\begin{example}
\label{ex:heis-bis}
Let $(A,d)$ be the Chevalley--Eilenberg complex of the Heisenberg Lie algebra $\h(1)$,
with $A=\bwedge(a_1,a_2,a_3)$ and $d(a_1)=d(a_2)=0$, $d(a_3)=-a_1\wedge a_2$, 
see Section~\ref{subsec:heisenberg-1}. 
Identifying $\Lie(A_1)$ with the free Lie algebra on generators $x_1,x_2,x_3$,
the ideal $\langle\im(\nabla_A)\rangle$ is generated by
$x_3-[x_1,x_2]$, $ [x_1,x_3]$, $[x_2,x_3]$.
Thus $\h(A)$ is generated by $x_1$ and $x_2$, with $x_3=[x_1,x_2]$, and hence
admits a grading by bracket length despite arising from a quadratic--linear presentation. 
By Proposition~\ref{prop:holo-CE}, we recover the intrinsic identification
$\h(A)\cong \h(1)$, and the above presentation shows that this
isomorphism is compatible with the induced weight grading.
\end{example}

\begin{lemma}
\label{lem:holo-to-cohomology}
Let $(A,d_A)$ be a connected $1$-finite $\k$-$\cdga$ equipped with a
positive grading such that $U_1=\im(d_A^\vee)$ is concentrated in 
weights $\ge 2$. There is an epimorphism of graded Lie algebras
\[
\Xi_A\colon \h(H^*(A)) \longsurj \h(A),
\]
functorial with respect to $\cdga$ morphisms, which is an isomorphism
in weights $1$ and~$2$.
\end{lemma}

\begin{proof}
Since $U_1 \subseteq \mathrm{weight} \ge 2$, proceed by induction on
weight. For each weight $w \ge 2$, consider the elements 
$u \in U_1$ of weight $w$, writing $u = d_A^\vee(a)$ for some $a \in A_2$
with $\mathrm{wt}(a) = w$. Decompose
\[
\mu_A^\vee(a) = (\pi_{E_1}\wedge\pi_{E_1})(\mu_A^\vee(a))
+ (\text{terms involving } U_1 \text{-elements of weight} < w),
\]
where the second summand involves only elements of $U_1$ with weight
strictly less than $w$, hence already expressed in terms of $E_1$ by
the inductive hypothesis. Substituting yields 
$u = -(\pi_{E_1}\wedge\pi_{E_1})(\mu_A^\vee(a)) + (\text{elements of }\Lie(E_1))$,
so $u$ is expressed entirely in $\Lie(E_1)$.
Since $U_1$ is finite-dimensional and each step strictly decreases the
minimum weight of uneliminated generators, the process terminates,
yielding a presentation $\h(A) = \Lie(E_1)/J_A$ with
$J_A \supseteq J_{H^*(A)}$.
The natural projection gives $\Xi_A$; functoriality is immediate.

Since any element of $U_1$ has weight $\ge 2$, the minimum weight of any
element in $J_A \setminus J_{H^*(A)}$ as a Lie bracket expression in $E_1$
is at least $1 + 2 = 3$ (one bracket with an element of weight $\ge 2$). 
Hence $J_A$ and $J_{H^*(A)}$ agree in weights $\le 2$, and $\Xi_A$ is an
isomorphism in those weights.
\end{proof}

The morphism $\Xi_A$ should be understood as a stability map: it expresses 
that, under the stated weight hypothesis, the linear part of the differential $d_A$ 
does not contribute new holonomy relations in weights $\le 2$. In particular, 
$\Xi_A$ is not induced by a universal functorial construction in the category 
of all $\cdgas$, but rather arises from a controlled elimination of high-weight 
linear terms.

%%%%%%%%%%%%%%%%%
\subsection{The universal enveloping algebra}
\label{subsec:UhA}
%%%%%%%%%%%%%%%%%

We now identify the universal enveloping algebra of the holonomy Lie algebra 
in terms of some functorial operations on quadratic algebras. 

Assume first that $A$ is a connected graded-commutative $\k$-algebra with
$A^0=\k$ and $\dim_\k A^1<\infty$. Every such algebra may be realized 
as the quotient $T(V)/I$ of the tensor algebra on $V=A^1$ by a 
homogeneous two-sided ideal $I$ generated in degrees $\ge 2$. 
Let $\qA=T(V)/\langle I^2\rangle$ denote the quadratic closure of $A$,
where $I^2$ is the degree-$2$ part of the ideal $I$, and let 
$(\qA)^{!}=T(V^{\vee})/I^{\perp}$ be its quadratic dual, where 
$I^{\perp}\subset T(V^{\vee})$ is the ideal generated by
$(I^2)^{\perp}\subset V^{\vee}\otimes_{\k} V^{\vee}$.

The next result was proved by Shelton--Yuzvinsky \cite[Lem.~4.1]{SY97} for 
Orlik--Solomon algebras, and generalized by Papadima--Yuzvinsky 
\cite[Lem.~4.1]{Papadima-Yuzvinsky}, with some improvements given 
by Suciu--Wang in \cite[Prop.~3.4]{SW-forum}.

\begin{proposition}
\label{prop:Papadima-Y}
Let $A$ be a graded, connected $\k$-algebra with $\dim_{\k} A^1<\infty$. Then 
\[
U(\h(A))\cong (\qA)^!.
\]
In particular, $U(\h(A))$ is quadratic, and depends only on the quadratic closure $\qA$.
\end{proposition}

For a general $\cdga$ $(A,d_A)$, the presence of a nontrivial differential
introduces linear (degree-$1$) relations in addition to quadratic ones.
Accordingly, the holonomy Lie algebra is governed by the quadratic-linear
dual of $A$ in the sense of Polishchuk--Positselski \cite{Polishchuk-Positselski}.

More precisely, let $V=A^1$, $A_1=V^\vee$, and $A_2=(A^2)^{\vee}$. 
Recall that the differential $d_A\colon A^1\to A^2$ induces a dual map
$d_A^\vee\colon A_2\to A_1$, and the multiplication 
$\mu_A \colon  A^1\wedge A^1  \to A^2$  induces a dual map 
$\mu_A^\vee\colon A_2 \to A_1 \wedge A_1$.  As in \eqref{eq:ha-nabla-a} below,
these combine into the map $\nabla_A=d_A^\vee+\mu_A^\vee\colon
A_2\to A_1\oplus(A_1\wedge A_1)$, whose image
\begin{equation}
\label{eq:lp-qp}
R\ \coloneqq\ \im(\nabla_A)\ \subset\ A_1\oplus (A_1\otimes_{\k} A_1)
\end{equation}
we regard as living in $A_1\oplus(A_1\otimes_{\k}A_1)$ via the antisymmetrization
$u\wedge v\mapsto u\otimes v - v\otimes u$.  This is a single subspace, not a
pair of independently chosen ones: a generic element of $R$ has both a
nonzero linear part and a nonzero quadratic part, and the two are locked
together by $\nabla_A$.  Consequently
\begin{equation}
\label{eq:quad-linear-closure}
U(\h(A)) \ =\ T(A_1)/\langle R\rangle ,
\end{equation}
an identity we record now and justify in the proof of
Proposition~\ref{prop:Uh-cdga}.  Following
\cite{Polishchuk-Positselski}, we call $\ql(A,d_A)$ the quadratic--linear
algebra, on generators $V$, whose Koszul dual (in the quadratic--linear
sense of \emph{loc.\ cit.}) is \eqref{eq:quad-linear-closure}; its existence
and uniqueness are part of that theory, which we do not otherwise need: every
statement below is phrased directly in terms of $(\ql(A,d_A))^!=T(A_1)/\langle
R\rangle$, which is computable from $A$ without reference to $\ql(A,d_A)$
itself.  (An earlier version of this closure attempted to describe
$\ql(A,d_A)$ by annihilating $L^\perp=\im(d_A^\vee)$ and
$Q^\perp=\im(\mu_A^\vee)$ \emph{independently}; this is not the same
construction as \eqref{eq:quad-linear-closure} unless one of the two
vanishes, and gives the wrong holonomy Lie algebra whenever $\nabla_A$
genuinely mixes degrees---already for $A=\CE(\sol_2)$ of
Example~\ref{ex:sol2-koszul}, it returns the abelian $1$-dimensional Lie
algebra in place of $\sol_2$ itself.)

\begin{proposition}
\label{prop:Uh-cdga}
Let $(A,d_A)$ be a connected $\k$-$\cdga$ with $\dim_\k A^1<\infty$. 
Then there is a natural isomorphism
\[
U(\h(A))\cong (\ql(A,d_A))^! .
\]
If $d_A^\vee=0$, then $\ql(A,d_A)=\qA$, and the identification 
above specializes to the quadratic case of Proposition~\ref{prop:Papadima-Y}.
\end{proposition}

\begin{proof}
Recall that
\begin{equation}
\label{eq:ha-nabla-a}
\h(A)=\Lie(A_1)\big/\langle \im(\nabla_A)\rangle,
\qquad
\nabla_A=d_A^\vee+\mu_A^\vee\colon 
A_2 \longrightarrow A_1\oplus(A_1\wedge A_1).
\end{equation}
Thus $\h(A)$ is the Lie algebra on generators $A_1$ whose defining ideal is
generated by the elements $d_A^\vee(\xi)+\mu_A^\vee(\xi)$, $\xi\in A_2$: each
such relation identifies the linear element $d_A^\vee(\xi)\in A_1$ with
$-\mu_A^\vee(\xi)\in\Lie_2(A_1)$, rather than setting the two separately to
zero.  Equivalently, the defining ideal is generated by the graph $R$ of
\eqref{eq:lp-qp}, embedded in $\Lie(A_1)$ via
$u\otimes v-v\otimes u\mapsto[u,v]$ on its quadratic part.

By the universal property of the enveloping algebra applied to the quotient
Lie algebra $\h(A)=\Lie(A_1)/\langle\im(\nabla_A)\rangle$, together with the
Poincar\'{e}--Birkhoff--Witt theorem for the free Lie algebra $\Lie(A_1)$
(which gives $U(\Lie(A_1))=T(A_1)$), the universal enveloping algebra
$U(\h(A))$ is the quotient of $T(A_1)$ by the two-sided ideal generated by
the same graph $R$, embedded via $u\wedge v\mapsto u\otimes v-v\otimes u$:
\[
U(\h(A)) = T(A_1)/\langle R\rangle .
\]
This is \eqref{eq:quad-linear-closure}, and by definition of $\ql(A,d_A)$ as
the quadratic--linear algebra dual to this presentation, $U(\h(A)) \cong
(\ql(A,d_A))^!$.

If $d_A^\vee=0$, every element of $R$ is purely quadratic, $R=\im(\mu_A^\vee)$,
so $T(A_1)/\langle R\rangle=(\qA)^!$ by Proposition~\ref{prop:Papadima-Y},
and hence $\ql(A,d_A)=\qA$.
\end{proof}

%%%%%%%%%%%%%%%%%
\subsection{Lower central series and enveloping algebra}
\label{subsec:lcs-short}
%%%%%%%%%%%%%%%%%

Although $\h(A)$ is defined from the cohomological grading of $A$, 
it need not inherit a grading compatible with bracket length. The 
discrepancy is measured by the lower central series filtration.

The \emph{lower central series} of a Lie algebra $\g$ is defined inductively by 
$\gamma_1\g=\g$ and $\gamma_{p+1}\g=[\gamma_p\g,\g]$.  It is readily seen 
that each $\gamma_p\g$ is an ideal of $\g$ and that 
$[\gamma_p \g, \gamma_q \g] \subset \gamma_{p+q} \g$. 
The associated graded Lie algebra 
$\gr(\g)=\bigoplus_{p\ge1}\gamma_p\g/\gamma_{p+1}\g$ 
is generated in degree~$1$, with grading given by bracket length.  
When $\g$ is already graded, $\gr(\g)=\g$; otherwise, 
the two may differ substantially.

\begin{example}
\label{ex:nonhomog-bis-short}
Let $A=\CE(\sol_2)$ be the $\cdga$ of Example~\ref{ex:sol2-koszul}.
Then $\h(A)=\sol_2=\Lie(x,y)/([x,y]-y)$, 
yet $\gr(\h)=\h/\h'=\k x$ is strictly smaller than $\h$, 
reflecting that the non-homogeneous relation identifies nontrivial 
commutators with lower-degree elements.
\end{example}

\begin{example}
\label{ex:not-carnot-short}
Let $A=(\bwedge(a_1,\dots,a_5),d)$ with 
$d(a_4)=a_1a_3$ and $d(a_5)=a_1a_4+a_2a_3$.  
As shown in~\cite[Ex.~10.5]{SW-forum}, the holonomy Lie algebra 
$\h(A)$ is not isomorphic to its associated graded $\gr(\h(A))$ 
under the LCS filtration.
\end{example}

The lower central series filtration on a Lie algebra $\g$ induces a
natural increasing filtration on its universal enveloping
algebra $U(\g)$, defined by
\begin{equation}
\label{eq:Ug-filtration}
\gamma_p U(\g) = \spn \{x_1\cdots x_r \mid x_i\in \gamma_{k_i}\g,\;
k_1+\cdots+k_r\le p\}.
\end{equation}
Since the LCS filtration on $\g$ is multiplicative, 
the induced filtration on $U(\g)$ is also multiplicative, that is, 
$\gamma_p U(\g) \cdot \gamma_q U(\g) \subset \gamma_{p+q} U(\g)$ for all 
$p,q$, where $\gamma_\bullet U(\g)$ denotes the induced filtration.  The 
associated graded algebra $\gr^{\gamma}U(\g)$ is taken with respect to
it. By the Poincar\'{e}--Birkhoff--Witt theorem, there is a canonical
isomorphism of graded algebras
\begin{equation}
\label{eq:gru}
\gr^{\gamma} U(\g) \cong U(\gr(\g)).
\end{equation}
Thus, the associated graded algebra of $U(\g)$ recovers the 
enveloping algebra of the graded Lie algebra $\gr(\g)$.

%%%%%%%%%%%%%%%%%
\subsection{Hilbert series and holonomy ranks}
\label{subsec:holo-ranks}
%%%%%%%%%%%%%%%%%

Let $(A,d)$ be a connected $\cdga$ with $\dim_\k A^1<\infty$.
When $d_A=0$, the universal enveloping algebra $U(\h(A))$ is quadratic; 
for a $\cdga$ with nontrivial differential, it is generally not. 
Proposition~\ref{prop:Uh-cdga} shows that, in this case, 
$U(\h(A))$ is quadratic--linear in the sense of \cite{Polishchuk-Positselski}.
In particular, $U(\h(A))$ carries a natural increasing filtration by 
tensor length, whose associated graded algebra is quadratic.

\begin{theorem}
\label{thm:gr-holonomy}
There is a natural isomorphism of graded algebras
\[
\gr^F U(\h(A)) \cong (\qql(A,d))^!,
\]
where $\qql(A,d)$ denotes the quadratic truncation of the quadratic--linear 
algebra $\ql(A,d)$, obtained by discarding the linear part of each defining
relation.
\end{theorem}

Concretely, writing $U(\h(A))=T(A_1)/\langle R\rangle$ as in
\eqref{eq:quad-linear-closure}, let $\pi_2\colon A_1\oplus(A_1\otimes_{\k}A_1)
\to A_1\otimes_{\k}A_1$ be the projection, and set
\begin{equation}
\label{eq:qql-def}
Q\coloneqq\pi_2(R)\subset A_1\otimes_{\k}A_1,
\qquad
(\qql(A,d))^! \coloneqq T(A_1)/\langle Q\rangle.
\end{equation}
This is a genuine truncation, not a restriction to fewer relations: since
$R$ is the graph of $\nabla_A$, $\dim Q=\dim R$ whenever $\nabla_A$ is
injective, the same number of relations as in $U(\h(A))$ itself; what is
discarded is only the linear correction attached to each one.

\begin{proof}
By Proposition~\ref{prop:Uh-cdga}, $U(\h(A))=T(A_1)/\langle R\rangle$, filtered
by tensor degree.  Since $R$ is the graph $\{d_A^\vee(\xi)+\mu_A^\vee(\xi):
\xi\in A_2\}$ of a linear correction over its quadratic part $Q=\pi_2(R)$, 
\cite[Thms.~4.4.1 and 4.7.1]{Polishchuk-Positselski} shows that passing to
the associated graded with respect to this filtration introduces no relations
beyond the leading (quadratic) term of each generator of $\langle R\rangle$:
\[
\gr^F\bigl(T(A_1)/\langle R\rangle\bigr) \cong T(A_1)/\langle Q\rangle 
= (\qql(A,d))^! .
\]
This is the content of the theorem; the containment $\gr^F U(\h(A))\subseteq
T(A_1)/\langle Q\rangle$ is elementary (the leading term of a relation always
survives to the associated graded), and the reverse containment is exactly
where the quadratic--linear hypothesis on $R$---as opposed to an arbitrary
filtered ideal---is used.
\end{proof}

Define the \emph{holonomy ranks} of $(A,d)$ by
\begin{equation}
\label{eq:holo-chen}
\phi_n(A)\coloneqq \dim_\k \gr_n(\h(A)),
\end{equation}
where $\gr_n(\h(A))$ denotes the $n$-th graded piece of the associated 
graded Lie algebra of $\h(A)$ with respect to the lower central series.
These invariants generalize the classical lower central series 
ranks of a finitely generated group to the $\cdga$ setting.

Applying the Poincar\'{e}--Birkhoff--Witt theorem \eqref{eq:gru} to the 
associated graded Lie algebra $\gr(\h(A))$, with respect to the lower
central series filtration $\gamma_\bullet U(\h(A))$ of \eqref{eq:Ug-filtration},
gives
\begin{equation}
\label{eq:pbw-holo}
\prod_{n\ge1}(1-t^n)^{\phi_n(A)}=\Hilb\bigl(\gr^{\gamma}U(\h(A)),t\bigr)^{-1}.
\end{equation}
This filtration $\gamma_\bullet$ is intrinsic to $\h(A)$; it is \emph{not}\/
the filtration $F$ of Theorem~\ref{thm:gr-holonomy}, which is by tensor
degree in the chosen generating set $A_1$ and depends on the presentation.
The two need not agree, so \eqref{eq:pbw-holo} cannot in general be evaluated
from $(A,d)$ via $(\qql(A,d))^!$ without saying when they do.

\begin{corollary}
\label{cor:pbw-holo-graded}
Suppose $\h(A)$ is a graded Lie algebra generated in degree~$1$---equivalently,
$A_1$ maps isomorphically onto $\gr_1(\h(A))$ under the lower central series.
Then $\gamma_\bullet$ and $F$ coincide on $U(\h(A))$, so that
$\gr^{\gamma}U(\h(A))\cong\gr^F U(\h(A))\cong(\qql(A,d))^!$ by
Theorem~\ref{thm:gr-holonomy}, and \eqref{eq:pbw-holo} becomes
\begin{equation}
\label{eq:pbw-holo-graded}
\prod_{n\ge1}(1-t^n)^{\phi_n(A)}=\Hilb\bigl((\qql(A,d))^!,t\bigr)^{-1}.
\end{equation}
\end{corollary}

The hypothesis holds whenever $(A,d)$ admits a positive-weight decomposition
with $H^1(A)$ concentrated in weight~$1$
(Proposition~\ref{prop:pos-weights-graded-holo}).

Without the hypothesis of Corollary~\ref{cor:pbw-holo-graded},
\eqref{eq:pbw-holo-graded} can fail outright.

\begin{example}
\label{ex:sol2-pbw-fails}
For $A=\CE(\sol_2)$, $\h(A)=\sol_2$ has $\gamma_2=\gamma_3=\cdots=\langle
y\rangle$, so $\phi_1=1$ and $\phi_n=0$ for $n\ge2$, giving
$\prod_{n}(1-t^n)^{\phi_n(A)}=1-t$; whereas $\qql(A,d)=\qA=\bwedge(a,b)$
(Example~\ref{ex:sol2-qkoszul}), with $\Hilb\bigl((\qql(A,d))^!,t\bigr)^{-1}
=\Hilb(\qA,t)^{-1}=(1-t)^{-2}$---not $\Hilb(\qA,-t)=(1-t)^2$ either.  No
identification of $\gr^{\gamma}$ with $\gr^F$ is available here, consistently
with the fact, already noted in Example~\ref{ex:sol2-koszul}, that
$\CE(\sol_2)$ admits no positive-weight decomposition at all.
\end{example}

When the hypothesis of Corollary~\ref{cor:pbw-holo-graded} holds and
$\qql(A,d)$ is in addition Koszul, Koszul duality converts
\eqref{eq:pbw-holo-graded} into the more usable
$\prod_{n}(1-t^n)^{\phi_n(A)}=\Hilb(\qql(A,d),-t)$; this is
Theorem~\ref{thm:koszul-holo} below, stated with its hypothesis made explicit.

%%%%%%%%%%%%%%%%%%%
\subsection{Koszul $\cdgas$ and holonomy ranks}
\label{subsec:koszul-cdga}
%%%%%%%%%%%%%%%%%%%

For a general $\cdga$ $(A,d)$, the associated graded algebra $\gr U(\h(A))$ 
is governed by the quadratic algebra $\qql(A,d)$; when $d=0$, this reduces 
to the usual situation of Koszul duality for graded algebras. This allows 
one to extend several classical consequences of Koszul duality to the 
non-formal setting.

\begin{definition}
\label{def:quadratically Koszul}
Let $(A,d)$ be a connected $\cdga$ with $\dim_\k A^1<\infty$.
We say that $(A,d)$ is {\em quadratically Koszul} if the quadratic algebra
$\qql(A,d)$ is Koszul (i.e., its trivial module admits a linear free resolution).
\end{definition}

This condition depends only on the quadratic part of the multiplication and
differential of $(A,d)$, and is automatic when $d=0$ and $A$ is a Koszul algebra.

\begin{example}
\label{ex:sol2-qkoszul}
Let $A = \CE(\sol_2)$ be the $\cdga$ of Example~\ref{ex:sol2-koszul}, with
$A_2=\langle(ab)^\vee\rangle$ and, writing $x=a^\vee,y=b^\vee\in A_1$,
\[
d_A^\vee\bigl((ab)^\vee\bigr)=-y,
\qquad
\mu_A^\vee\bigl((ab)^\vee\bigr)=x\wedge y .
\]
The graph $R=\im(\nabla_A)$ of \eqref{eq:lp-qp} is spanned by the single
element $-y+(x\otimes y-y\otimes x)$, mixing a linear and a quadratic term,
so that
\[
U(\h(A))=T(x,y)/\langle xy-yx-y\rangle
\]
by \eqref{eq:quad-linear-closure}---the standard presentation of
$U(\sol_2)$.  Discarding the linear part of the one relation gives
$Q=\langle xy-yx\rangle$, so
\[
(\qql(A,d))^!=T(x,y)/\langle xy-yx\rangle=\k[x,y],
\]
a polynomial ring, which is Koszul; dually, $\qql(A,d)=\bwedge(a,b)=\qA$, by
the classical duality between symmetric and exterior algebras.  Thus
$(A,d)$ is quadratically Koszul, even though $\nabla_A$ genuinely mixes
degrees here: $d_A^\vee\ne0$, yet $\qql(A,d)$ coincides with the $d=0$
answer $\qA$ of Proposition~\ref{prop:Uh-cdga}.  This illustrates that
quadratic Koszulness depends only on the quadratic part $Q$ of $R$, and is
insensitive to the linear correction carried by $d_A^\vee$---in contrast
with $\phi_n(A)$ itself, for which the linear correction is decisive; see
the discussion after \eqref{eq:pbw-holo}.
\end{example}

When $(A,d)$ is quadratically Koszul, the Hilbert series of $\qql(A,d)$
controls the holonomy ranks via Koszul duality, yielding a concrete 
and computable formula.

\begin{theorem}
\label{thm:koszul-holo}
Let $(A,d)$ be a quadratically Koszul $\cdga$ such that $\h(A)$ is a graded
Lie algebra generated in degree~$1$.  Then the holonomy ranks
$\phi_n(A)$ are determined by the Hilbert series of $\qql(A,d)$ via the
identity
\begin{equation}
\label{eq:lcs-holo-hilb}
\prod_{n\ge1}(1-t^n)^{\phi_n(A)} = \Hilb(\qql(A,d),-t).
\end{equation}
\end{theorem}

The hypothesis on $\h(A)$ holds, for instance, whenever $(A,d)$ admits a
positive-weight decomposition with $H^1(A)$ concentrated in weight~$1$
(Proposition~\ref{prop:pos-weights-graded-holo}).

\begin{proof}
Since $\h(A)$ is graded and generated in degree~$1$, the tensor-degree
filtration $F$ on $U(\h(A))$ coincides with the lower-central-series
filtration $\gamma_\bullet$, so \eqref{eq:pbw-holo-graded} applies:
\begin{equation}
\label{eq:PBW-grU}
\prod_{n\ge1}(1-t^n)^{\phi_n(A)} = \Hilb\bigl((\qql(A,d))^!,t\bigr)^{-1}.
\end{equation}
If $\qql(A,d)$ is Koszul, then so is its dual $(\qql(A,d))^!$, and Koszul
duality yields
\begin{equation}
\label{eq:qqlAd}
\Hilb(\qql(A,d),t)\cdot \Hilb((\qql(A,d))^!,-t)=1,
\qquad\text{equivalently}\qquad
\Hilb\bigl((\qql(A,d))^!,t\bigr)^{-1}=\Hilb(\qql(A,d),-t).
\end{equation}
Combining \eqref{eq:PBW-grU} and \eqref{eq:qqlAd} proves the claim.
\end{proof}

\begin{remark}
\label{rem:koszul-hypothesis}
When $d=0$, $\h(A)$ is automatically graded and generated in degree~$1$,
$\qql(A,0)$ coincides with the quadratic closure $\qA$, and
Theorem~\ref{thm:koszul-holo} recovers the classical formula for graded
Koszul algebras.  Theorem~\ref{thm:koszul-holo} in fact carries \emph{two}\/
hypotheses beyond quadraticity, of different characters.  The Koszul
hypothesis on $\qql(A,d)$ is needed only to pass from
\eqref{eq:PBW-grU} to \eqref{eq:lcs-holo-hilb} via the duality identity
\eqref{eq:qqlAd}; it is not needed in Theorem~\ref{thm:gr-holonomy}, where
the isomorphism $\gr^F U(\h(A))\cong(\qql(A,d))^!$ holds unconditionally by
\cite[Thms.~4.4.1 and 4.7.1]{Polishchuk-Positselski}.  The hypothesis that
$\h(A)$ be graded and generated in degree~$1$ is needed earlier still, to
identify $\gr^F U(\h(A))$---computed from the presentation---with
$\gr^{\gamma}U(\h(A))$, which is what the Poincar\'{e}--Birkhoff--Witt
identity \eqref{eq:pbw-holo} actually controls; see the discussion
preceding \eqref{eq:pbw-holo-graded}, and Example~\ref{ex:sol2-qkoszul}
for a $\cdga$ satisfying the first hypothesis but not the second.
\end{remark}

%%%%%%%%%%%%%%%%%
\subsection{The $1$-minimal model from holonomy}
\label{subsec:cochain-holo-short}
%%%%%%%%%%%%%%%%%

The holonomy Lie algebra is not merely an invariant of~$A$; it
functorially determines the degree-$1$ truncation of its rational
homotopy type.  More precisely, it gives rise to a canonical
$1$-minimal Sullivan model of~$A$, generated in degree~$1$ with
decomposable differential and quasi-isomorphic to~$A$ in degrees
$\le 1$.

Let $\CE(\g)$ denote the Chevalley--Eilenberg $\cdga$ of a Lie algebra~$\g$.  
This construction is functorial and sends central extensions of 
finite-dimensional nilpotent Lie algebras to Hirsch extensions of~$\cdgas$.  
For a finitely generated Lie algebra~$\g$, each quotient 
$\g/\gamma_{n+1}(\g)$ is a finite-dimensional $n$-step nilpotent Lie algebra, 
fitting into a central extension
\begin{equation}
\label{eq:gamma-ext-short}
\begin{tikzcd}[column sep=20pt]
0 \arrow[r] & \gamma_n(\g)/\gamma_{n+1}(\g)
\arrow[r] & \g/\gamma_{n+1}(\g)
\arrow[r] & \g/\gamma_{n}(\g)
\arrow[r] & 0,
\end{tikzcd}
\end{equation}
which is natural with respect to Lie algebra morphisms. 
Applying the contravariant functor $\CE$ to this inverse system yields a direct system 
$\{\CE(\g/\gamma_n(\g))\}_{n\ge1}$ whose colimit
\begin{equation}
\label{eq:chat-short}
\wC(\g)\coloneqq \varinjlim_{n} \CE(\g/\gamma_n(\g))
\end{equation}
is a $1$-minimal $\cdga$ in the sense of Sullivan \cite{Sullivan77}. 
Moreover, $\wC(\g)$ inherits a natural filtration from the lower central 
series, and may be viewed as a filtered $\cdga$ obtained as a colimit of 
Hirsch extensions,
\[
\CE(\g/\gamma_n(\g)) \longinj \CE(\g/\gamma_{n+1}(\g)),
\]
where each step is a Hirsch extension adding degree-$1$ generators dual to 
$\gamma_n(\g)/\gamma_{n+1}(\g)$, with differential determined by the 
corresponding central extension \eqref{eq:gamma-ext-short}.

\begin{lemma}
\label{lem:recover-g-from-chat}
Let $\g$ and $\bar\g$ be finitely generated Lie algebras.
If $\wC(\g)\cong \wC(\bar\g)$ as filtered $\cdgas$, then 
their lower central series completions are isomorphic,
$\widehat{\g}\cong \widehat{\bar\g}$.
In particular, if $\g$ and $\bar\g$ are residually nilpotent,
then $\g\cong \bar\g$ as filtered Lie algebras.
\end{lemma}

\begin{proof}
By construction, $\wC(\g)$ is the direct limit of the Chevalley--Eilenberg 
algebras $\CE(\g/\gamma_n(\g))$, each of which is a minimal $\cdga$ 
generated in degree~$1$. The degree-$1$ generators of 
$\CE(\g/\gamma_n(\g))$ identify with $(\g/\gamma_n(\g))^\vee$, and the 
quadratic part of the differential encodes the Lie bracket on 
$\g/\gamma_n(\g)$. 

A filtered isomorphism $\wC(\g)\cong \wC(\bar\g)$ therefore induces, for each $n$, 
an isomorphism $\g/\gamma_n(\g)\cong \bar\g/\gamma_n(\bar\g)$ compatible with 
the natural projections. These isomorphisms assemble into an isomorphism 
of inverse systems, and hence induce an isomorphism 
$\widehat{\g}\cong \widehat{\bar\g}$ of pronilpotent Lie algebras.
\end{proof}

The relationship between the holonomy Lie algebra of a $\cdga$ and its
Sullivan minimal model is particularly transparent in degree~$1$.
In this setting, the quadratic information encoded by the holonomy Lie
algebra controls the $1$-minimal model via the Chevalley--Eilenberg
construction and the lower central series.

The following result is due to Papadima--Suciu~\cite[Thm.~5.4]{PS-jlms}.
Although stated and proved there in a purely algebraic framework, it
should be viewed as a generalization of earlier results
on Malcev completions and quadratic models due to Bezrukavnikov \cite{Bez},
Bibby--Hilburn \cite{BH}, and Berceanu--M\u{a}cinic--Papadima--Popescu \cite{BMPP}.

If $A$ is a $1$-finite $\cdga$ with holonomy Lie algebra $\h=\h(A)$,
there exist natural morphisms
$\rho_n\colon \CE(\h/\gamma_n(\h)) \to A$,
functorial in $A$, encoding the $n$-step truncations of the
$1$-minimal structure.

\begin{theorem}[\cite{PS-jlms}]
\label{thm:mm-holo-short}
If $A$ is a $1$-finite $\cdga$, then $\wC(\h(A))$ 
is a $1$-minimal model for~$A$, with structure map 
$\rho\colon \wC(\h(A))\to A$.
\end{theorem}

In other words, for any $1$-finite $\cdga$ $A$, the completed
Chevalley--Eilenberg algebra of its holonomy Lie algebra provides a
canonical $1$-minimal model of~$A$.  This construction is functorial
and captures precisely the quadratic part of the rational homotopy type 
encoded in the $1$-minimal model.

As shown in~\cite{PS-jlms}, this theorem recovers, in a self-contained
and functorial way, a result from \cite{BMPP}, which in turn generalizes  
work from \cite{Bez}; see also \cite{BH} for a concise summary.

\begin{example}
\label{ex:nonhomog-tri-short}
Let $A=(\bigwedge(a,b),d)$ with $d(a)=0$ and $d(b)=b\wedge a$. 
As noted in Example~\ref{ex:nonhomog-bis-short}, we have 
$\h(A)=\Lie(x,y)/([x,y]-y)$.  
Then $\gamma_n(\h)=[\h,\h]$ for all $n\ge2$, and thus 
$\wC(\h)=\bwedge(x^{\vee})$ with $d=0$.  
The map $\psi\colon \wC(\h(A))\to A$ sends $x^{\vee}\mapsto a$. 
Thus, although $\h(A)$ is non-nilpotent, its associated $1$-minimal model is abelian.
\end{example}

\begin{corollary}
\label{cor:1eq-holonomy}
Let $A$ and $B$ be connected, $1$-finite $\k$-$\cdgas$.
If $A$ and $B$ are $1$-equivalent, then 
$\widehat{\h(A)}\cong \widehat{\h(B)}$, and hence 
$\gr{\h(A)}\cong \gr{\h(B)}$. 
\end{corollary}

\begin{proof}
If $A$ and $B$ are $1$-equivalent, then their $1$-minimal models are
isomorphic. By Theorem~\ref{thm:mm-holo-short}, these models are given by
$\wC(\h(A))$ and $\wC(\h(B))$. Hence
$\wC(\h(A))\cong \wC(\h(B))$ as filtered $\cdgas$.
By Lemma~\ref{lem:recover-g-from-chat}, this induces an isomorphism
$\widehat{\h(A)}\cong \widehat{\h(B)}$.
\end{proof}

\begin{remark}
\label{rem:holonomy-1formal-warning}
Corollary~\ref{cor:1eq-holonomy} should be interpreted with care.
While $1$-formality implies that $\h(A)\cong \h(H^*(A))$ (up to filtered
isomorphism), the converse fails in general.

This is already visible for $1$-step Hirsch extensions, and in a strong form.
On the one hand, adjoining the generator $c$ to the model
$A=B\otimes_{e}\bigwedge(c)$ gives the holonomy Lie algebra $\h(A)$ an extra
generator $c^{\vee}$, dual to $c$, which is central and satisfies a relation
expressing it in terms of the cup product; on the other hand, passing to
cohomology kills the class $[e]$, and with it one quadratic relation.  
These two effects cancel \emph{always}, not merely in special cases: one has
$\phi_2(A)=\phi_2(H^{*}(A))$ for every $1$-step Hirsch extension, by
Remark~\ref{rem:holo-no-quadratic-obstruction}.  Any discrepancy is therefore
generated in degree~$3$, by the cubic relations $[H_1,z]$---and generated
there in the strong sense that if it fails to appear in degree~$3$ it appears
in no degree at all; see Remark~\ref{rem:holo-no-quadratic-obstruction}.

Consequently the holonomy Lie algebra cannot detect non-$1$-formality at the
quadratic level here, and \emph{a fortiori}\/ cannot detect non-formality.
Sasakian manifolds of dimension $2n+1$ with $n>1$ illustrate the wider gap:
such a manifold is $(n-1)$-formal, hence $1$-formal, while typically not
formal, so that the obstruction to formality lies in a degree that no
$1$-type invariant---holonomy Lie algebra included---can see; see
\cite{PS-imrn19}.
\end{remark}

We now recall the standard construction (see
Section~\ref{subsec:minimal-models}) associating to a $1$-minimal
$\cdga$ $\M$ a complete Lie algebra $\fL(\M)$, obtained by dualizing
$\M^1$ and using the differential to define the bracket. This
construction is functorial and adjoint to the Chevalley--Eilenberg
functor $\wC$, in the sense that
$\fL(\wC(\g))\cong \widehat{\g}=\varprojlim_{n} \g/\gamma_n(\g)$, the completion
of $\g$ with respect to its lower central series, for any finitely generated
Lie algebra $\g$.  (Note that $\wC(\g)$ depends only on the tower
$\{\g/\gamma_n(\g)\}$, hence only on $\widehat{\g}$.)  The completion 
cannot be omitted here: by Example~\ref{ex:nonhomog-tri-short}, 
$\wC(\sol_2)=\bwedge(x^{\vee})$, so that $\fL(\wC(\sol_2))=\k\ne\sol_2$.

\begin{proposition}
\label{prop:holonomy-filtered-formal}
Let $(A,d_A)$ be a connected, $1$-finite $\k$-$\cdga$.
Suppose $A$ admits a $1$-minimal model $\cM$ with positive Hirsch weights.
Then the holonomy Lie algebra $\h(A)$ is filtered-formal.
\end{proposition}

\begin{proof}
Since $d_{\cM}$ is homogeneous with respect to the Hirsch weights, the
bracket induced on the dual of $\cM^1$ is homogeneous as well; write $\g$
for the resulting graded Lie algebra, so that $\fL(\cM)\cong\widehat{\g}$.

By Theorem~\ref{thm:mm-holo-short}, the $\cdga$ $\wC(\h(A))$ is a
$1$-minimal model for $A$.  Since $\cM$ is a $1$-minimal model for $A$ as
well, uniqueness of $1$-minimal models gives $\cM \cong \wC(\h(A))$.
Applying $\fL(-)$, we obtain
\[
\fL(\cM) \cong \fL(\wC(\h(A))) \cong \widehat{\h(A)},
\]
where the last isomorphism follows from the adjunction between $\fL$
and $\wC$: the degree-$1$ generators of $\wC(\h)$ are dual to
$\h/\gamma_2(\h)$, and the differential encodes the Lie bracket.  Hence
$\widehat{\h(A)}$ is the completion of a graded Lie algebra, which is to say
that $\h(A)$ is filtered-formal.
\end{proof}

%%%%%%%%%%%%%%%%%%%%%%%%%%%%
\subsection{Holonomy Lie algebras of Hirsch extensions}
\label{subsec:holo-hirsch}
%%%%%%%%%%%%%%%%%%%%%%%%%%%%

Let $B$ be a connected $\cdga$ and let $A=(B\otimes_{e}\bigwedge(c),d)$
be a $1$-step Hirsch extension with $d(c)=e$, $e\in Z^2(B)$.
The $\cdga$ $A$ depends only on the class $[e]\in H^2(B)$~\cite[Lem.~2.1]{PS-imrn19}.
When $e\neq 0$, we have $H^1(A)\cong H^1(B)$, and the cup product on $H^1(A)$ is
\begin{equation}
\label{eq:mua-hirsch}
\begin{tikzcd}[column sep=22pt]
\mu_A\colon  H^1(B)\wedge H^1(B)
\arrow[r, "\mu_B"] &[4pt] H^2(B)
\arrow[r, two heads] & H^2(B)/\langle [e]\rangle.
\end{tikzcd}
\end{equation}

The holonomy Lie algebra of a Hirsch extension admits an explicit
presentation, which we record first; it is what makes the centrality of
$c^{\vee}$ available, and it is stated for an arbitrary base.

\begin{lemma}
\label{lem:holo-hirsch-presentation}
Let $A=B\otimes_{e}\bigwedge(c)$ be a $1$-step Hirsch extension of a
connected, $1$-finite $\cdga$ $B$, and identify $A_1=B_1\oplus\k c^{\vee}$.
Then
\begin{equation}
\label{eq:holo-hirsch-pres}
\h(A)=\Lie\bigl(B_1\oplus \k c^{\vee}\bigr)\Big/
\Bigl\langle\, [\xi,c^{\vee}]\ (\xi\in B_1);\quad
d_B^{\vee}\gamma+\langle e,\gamma\rangle\, c^{\vee}+\mu_B^{\vee}\gamma\
(\gamma\in B_2)\,\Bigr\rangle .
\end{equation}
In particular, $c^{\vee}$ is central in $\h(A)$.
\end{lemma}

\begin{proof}
Write $A^1=B^1\oplus\k c$ and $A^2=B^2\oplus B^1c$, and dualize.  For
$\gamma\in B_2$, evaluating $d_A^{\vee}\gamma$ against $b\in B^1$ gives
$\langle\gamma,d_Bb\rangle$ and against $c$ gives $\langle\gamma,e\rangle$,
so that $d_A^{\vee}\gamma=d_B^{\vee}\gamma+\langle e,\gamma\rangle c^{\vee}$;
while $\mu_A^{\vee}\gamma=\mu_B^{\vee}\gamma$, since $\gamma$ annihilates
$B^1c$.  This gives the second family of relations.  For
$\delta\in(B^1c)^{\vee}$, we have $d_A^{\vee}\delta=0$, because $\delta$
annihilates $B^2=d_A(A^1)$; and $\mu_A^{\vee}\delta=\xi_{\delta}\wedge
c^{\vee}$, where $\xi_{\delta}\in B_1$ is defined by
$\langle\xi_{\delta},b\rangle=\langle\delta,bc\rangle$.  As $\delta$ ranges
over $(B^1c)^{\vee}$, the element $\xi_{\delta}$ ranges over all of $B_1$,
which gives the first family.  Since $\nabla_A=d_A^{\vee}+\mu_A^{\vee}$ and
$A_2=B_2\oplus(B^1c)^{\vee}$, these two families generate the defining ideal
of $\h(A)$.
\end{proof}

We now specialize to the case $d_B=0$, that is, $B=(H,0)$ with
$H=H^{*}(B)$.  This covers the situations of interest---Tievsky models of
Sasakian manifolds, and Chevalley--Eilenberg complexes of $2$-step nilpotent
Lie algebras---and the restriction is not a matter of convenience: the
holonomy Lie algebra is \emph{not}\/ a quasi-isomorphism invariant, by
Example~\ref{ex:nonhomog-tre}, so one may not replace a formal $B$ by its
cohomology at this point.

\begin{proposition}
\label{prop:holo-hirsch}
Let $A=H\otimes_{e}\bigwedge(c)$ be a $1$-step Hirsch extension of a
connected graded algebra $H$ with zero differential, where $e\in H^2$ and
$H^1\ne 0$.  Write $\mu=\mu_H$.
Then the holonomy ranks of $A$ satisfy
\begin{enumerate}[itemsep=2pt]
\item \label{eh0}
If $e= 0$, then $\h(A)\cong\h(H)\oplus\Lie(c^{\vee})$; hence
$\phi_1(A)=1+\dim_\k H^1$ and $\phi_n(A)=\phi_n(H)$ for all $n\ge 2$. 
\item \label{eh1}
If $e\ne 0$ and $e\in\im(\mu)$, then $\h(A)$ is a graded Lie algebra sitting
in a central extension
\begin{equation}
\label{eq:holo-hirsch-central}
\begin{tikzcd}[column sep=20pt]
0 \arrow[r] & \k z \arrow[r] & \h(A) \arrow[r] & \h(H) \arrow[r] & 0,
\end{tikzcd}
\end{equation}
with $z\ne0$ of degree~$2$; consequently $\phi_1(A)=\dim_\k H^1$,
$\phi_2(A)=\phi_2(H)+1$, and $\phi_n(A)=\phi_n(H)$ for all $n\ge 3$.
\item \label{eh2}
If $e\ne0$ but $e\notin\im(\mu)$, then $\h(A)\cong\h(H)$, so
$\phi_n(A)=\phi_n(H)$ for all $n\ge1$.
\end{enumerate}
\end{proposition}

\begin{proof}
Since $d_H=0$, the presentation \eqref{eq:holo-hirsch-pres} reads
\begin{equation}
\label{eq:holo-hirsch-pres-dB0}
\h(A)=\Lie\bigl(H_1\oplus\k c^{\vee}\bigr)\Big/
\Bigl\langle\, [\xi,c^{\vee}]\ (\xi\in H_1);\quad
\langle e,\gamma\rangle\,c^{\vee}+\mu^{\vee}\gamma\ (\gamma\in H_2)
\,\Bigr\rangle .
\end{equation}

If $e=0$, the second family becomes $\im(\mu^{\vee})$, while the first says
that $c^{\vee}$ is central and, being a generator not occurring in any
relation, spans a free abelian direct factor.  Hence
$\h(A)\cong\h(H)\oplus\Lie(c^{\vee})$, which is \eqref{eh0}.

Assume now $e\ne0$, and let $e^{\vee}\colon H_2\to\k$ be evaluation on $e$.
Choose $\gamma_0\in H_2$ with $\langle e,\gamma_0\rangle=1$, and set
$z=\mu^{\vee}\gamma_0$.  The relation attached to $\gamma_0$ eliminates the
generator $c^{\vee}$, giving $c^{\vee}=-z$.  For an arbitrary $\gamma\in H_2$,
subtracting $\langle e,\gamma\rangle$ times the relation for $\gamma_0$ turns
the relation for $\gamma$ into $\mu^{\vee}\bigl(\gamma-\langle
e,\gamma\rangle\gamma_0\bigr)=0$; as $\gamma$ varies, the elements
$\gamma-\langle e,\gamma\rangle\gamma_0$ range over $\ker(e^{\vee})$.
Substituting $c^{\vee}=-z$ into the first family, we conclude that
\begin{equation}
\label{eq:holo-hirsch-eliminated}
\h(A)=\Lie(H_1)\Big/\bigl\langle \mu^{\vee}(\ker e^{\vee}),\
[H_1,z]\bigr\rangle .
\end{equation}
Both families of relations are homogeneous---of degrees $2$ and $3$
respectively---so $\h(A)$ is a graded Lie algebra, and no appeal to
filtered formality is needed.

Suppose first that $e\notin\im(\mu)$.  Then some $\gamma\in H_2$ annihilates
$\im(\mu)$ while $\langle e,\gamma\rangle\ne0$, so $\gamma_0$ may be chosen
in $\ker(\mu^{\vee})$; thus $z=0$.  The same choice gives
$\ker(\mu^{\vee})+\ker(e^{\vee})=H_2$, whence
$\mu^{\vee}(\ker e^{\vee})=\im(\mu^{\vee})$ and
$\h(A)=\Lie(H_1)/\langle\im\mu^{\vee}\rangle=\h(H)$.  This is \eqref{eh2}.

Suppose next that $e\in\im(\mu)$, so that $e\ne0$ forces $z\ne0$ in
$\bwedge^2 H_1$: otherwise $\gamma_0$ would annihilate $\im(\mu)$ while
pairing nontrivially with $e$.  The same argument shows
$z\notin\mu^{\vee}(\ker e^{\vee})$, which is the degree-$2$ part of the
defining ideal in \eqref{eq:holo-hirsch-eliminated}; hence $z\ne0$ in
$\h(A)$.  Moreover $\im(\mu^{\vee})=\mu^{\vee}(\ker e^{\vee})+\k z$, so
adjoining $z$ to the relations of \eqref{eq:holo-hirsch-eliminated} recovers
$\h(H)$.  As $z$ is central and $[H_1,z]=0$ in $\h(A)$, the ideal it
generates is $\k z$, and we obtain the central extension
\eqref{eq:holo-hirsch-central}.  Finally, comparing the defining ideals
$I_A$ and $I_H$ of $\h(A)$ and $\h(H)$ degreewise:
$(I_A)_2=\mu^{\vee}(\ker e^{\vee})$ has codimension one in
$(I_H)_2=\im(\mu^{\vee})$, whereas for $n\ge3$
\[
(I_A)_n=\bigl[H_1,(I_A)_{n-1}\bigr]+\bigl[H_1,\dots,[H_1,z]\bigr]
=\bigl[H_1,(I_H)_{n-1}\bigr]=(I_H)_n ,
\]
since $[H_1,z]\subseteq (I_A)_3$.  The stated ranks follow.
\end{proof}

\begin{proposition}
\label{prop:holo-hirsch-cohomology}
Let $A=B\otimes_{e} \bigwedge(c)$ be a $1$-step Hirsch extension, where
$B$ is a connected $\cdga$.

\begin{enumerate}[itemsep=2pt]
\item \label{eh-coh1}
If $e=0$, then $\h\big(H^*(A)\big)\cong \h(B)\oplus \Lie(c^\vee)$. 

\item \label{eh-coh2}
If $[e]\neq 0$, then
$\h\big(H^*(A)\big)\cong
\Lie(H^1(B)^\vee)\big/\big\langle \mu_B^\vee\bigl(\ker(e^\vee)\bigr)\big\rangle$, 
where $e^\vee\colon H^2(B)^\vee\to \Bbbk$ is evaluation on $[e]$.
Equivalently, $\h(H^*(A))$ is obtained from $\h(B)$ by \emph{removing}
those quadratic relations whose duals evaluate nontrivially on $[e]$.
\end{enumerate}
\end{proposition}

\begin{proof}
If $e=0$, then $A=B\otimes_{\k} \bigwedge(c)$ with $dc=0$, and thus
$H^*(A)\cong H^*(B)\otimes_{\k} \bigwedge([c])$. In particular,
$H^1(A)=H^1(B)\oplus \k\cdot [c]$, and the cup product involving
$[c]$ vanishes. It follows that $\h(H^*(A))\cong \h(B)\oplus \Lie(c^\vee)$. 

Assume now $[e]\neq 0$. Then $H^1(A)=H^1(B)$, and the cup product on
$H^1(A)$ is given by \eqref{eq:mua-hirsch}. Dualizing, we obtain
$\im(\mu_A^\vee)=\mu_B^\vee\bigl(\ker(e^\vee)\bigr)$, 
and the conclusion follows from the definition of the holonomy Lie algebra.
\end{proof}

\begin{remark}
\label{rem:holo-no-quadratic-obstruction}
Comparing Proposition~\ref{prop:holo-hirsch} with
Proposition~\ref{prop:holo-hirsch-cohomology} shows that the holonomy Lie
algebra detects \emph{no}\/ obstruction to $1$-formality at the quadratic
level.  Write $W=\mu^{\vee}(\ker e^{\vee})$.  If $e\ne0$ and
$e\notin\im(\mu)$, both propositions give $\h(H)$.  If $e\ne0$ and
$e\in\im(\mu)$, the defining ideals of $\h(A)$ and of $\h(H^{*}(A))$ have the
same degree-$2$ part, namely $W$, by
\eqref{eq:holo-hirsch-eliminated} and
Proposition~\ref{prop:holo-hirsch-cohomology}\eqref{eh-coh2}; and if $e=0$,
the two Lie algebras agree outright, both being
$\h(H)\oplus\Lie(c^{\vee})$.  In every case, then,
$\phi_2(A)=\phi_2(H^{*}(A))$.

In degree~$3$, by contrast, the comparison has content.  Here
$(I_A)_3=[H_1,W]+[H_1,z]$ while $(I_{H^{*}(A)})_3=[H_1,W]$, so
$\phi_3(A)\le\phi_3(H^{*}(A))$, with equality precisely when
$[H_1,z]\subseteq[H_1,W]$.  Strict inequality does occur, in the smallest
example available: for $H^1=\langle u_1,u_2\rangle$, $H^2=\langle h\rangle$
with $\mu(u_1u_2)=h$, and $e=h$, one has $\ker e^{\vee}=0$, hence $W=0$ and
$z=[x_1,x_2]\ne0$; so $\h(A)$ is the free $2$-step nilpotent Lie algebra
$\h(1)$, with $\phi(A)=(2,1,0,0,\dots)$, whereas $\mu_{H^{*}(A)}=0$ makes
$\h(H^{*}(A))$ free, with $\phi(H^{*}(A))=(2,1,2,3,6,\dots)$.  Thus
$\phi_2$ agrees while $\phi_3$ does not, and the holonomy comparison detects
the non-$1$-formality of $\CE(\h(1))$---the same conclusion reached in
Example~\ref{ex:massey-hirsch-two-cases} through the triple Massey product
$\langle a,b,b\rangle$.

Degree~$3$ is moreover decisive, in the following sense: since
$I_{H^{*}(A)}$ is an ideal containing $W$, it contains the ideal generated by
$[H_1,z]$ as soon as it contains $[H_1,z]$ itself.  Hence
\begin{equation}
\label{eq:holo-degree3-decisive}
\phi_3(A)=\phi_3(H^{*}(A))
\iff [H_1,z]\subseteq [H_1,W]
\iff \h(A)\cong\h(H^{*}(A)) .
\end{equation}
Both alternatives occur, and neither is exotic.  In the example just given
the first fails.  It holds, on the other hand, for
$B=H^{*}(T^4)$ with $e=a_1a_2+a_3a_4$: there $\mu$ is bijective, $W$ has
codimension one in $\bwedge^2H_1$, and $[H_1,W]$ already exhausts
$\Lie_3(H_1)$, so that $\h(A)\cong\h(H^{*}(A))$ although $z\ne0$.  For such
an $A$ the holonomy Lie algebra detects nothing whatever, in any degree, and
the criterion of Corollary~\ref{cor:e-decomposable-massey} does not apply
either, $[e]$ not being decomposable; any obstruction to $1$-formality would
have to be sought elsewhere.
\end{remark}

\begin{remark}
\label{rem:hirsch-sasakian-pointer}
This applies in particular to Tievsky models
$A=H^*(B)\otimes_{[e]}\bigwedge(c)$ of compact Sasakian manifolds, where $B$
is a formal K\"{a}hler orbifold and $[e]\in H^2(B)$ is the Euler class of the
Boothby--Wang fibration, so that $[e]$ is the K\"{a}hler class of $B$.  Here
Corollary~\ref{cor:e-decomposable-massey} applies precisely when $B$ is a
Riemann surface with $H^1(B)\ne0$---so, among $3$-dimensional Sasakian
manifolds, to those whose base has positive genus.  Indeed, a decomposable 
class has square zero, since $(\alpha\beta)^2=-\alpha^2\beta^2=0$ for $\alpha,\beta$ 
of degree~$1$, whereas Hard Lefschetz gives $[e]^{n}\ne0$ on a K\"{a}hler base of 
complex dimension~$n$; so $[e]$ can be decomposable only for $n=1$, that is, when $B$
is a Riemann surface.  One-dimensionality of $H^2(B)$ does not by itself make
it so, however: what is needed is $H^1(B)\ne0$, and then $\mu_B$ is onto, so
that $[e]$ is a nonzero multiple of $\alpha\beta$ for a symplectic pair
$\alpha,\beta\in H^1(B)$, and rescaling gives $[e]=\alpha\beta$.  This is
consistent with Remark~\ref{rem:n1-enot0}: for $n>1$ the manifold is
$(n-1)$-formal, hence $1$-formal, and no degree-one obstruction can exist.
At the other extreme, when $[e]\notin\im(\mu_B)$ no obstruction arises at
this level, and indeed $A$ may well be formal.  The excluded genus-zero case
lands here: for $B=\mathbb{P}^1$ we have $H^1(B)=0$, so $\im(\mu_B)=0$ and the
K\"{a}hler class is not decomposable; this is the case $B=\k[h]/(h^2)$, $e=h$ 
of Example~\ref{ex:massey-hirsch-two-cases}, whose Boothby--Wang total 
space is $S^3$, a formal space.  
See \cite[Thm.~6.6]{PS-imrn19} for the full partial-formality statement.
\end{remark}

%%%%%%%%%%%%%%%%%%%%%%%%%%
\section{Infinitesimal Alexander invariants and holonomy Chen ranks}
\label{sect:infalex-chen}
%%%%%%%%%%%%%%%%%%%%%%%%%%

Building on the analysis of the first Koszul module $\B_1(A)$ from
Section~\ref{sect:BAd}, we relate this module to the infinitesimal
Alexander invariant $\B(\h(A))=\h'(A)/\h''(A)$ of the holonomy Lie algebra.
This identification provides a bridge between commutative-algebraic
invariants of $(A,d)$ and Lie-theoretic invariants of $\h(A)$, and leads
to effective formulas for the holonomy Chen ranks of a connected
$\cdga$. The main result of this section establishes Theorem~\ref{thm:B-holo-intro} 
from the Introduction. In particular, we obtain infinitesimal analogues of classical
results of Massey relating Alexander invariants and Chen ranks for
finitely generated groups.

%%%%%%%%%%%%%%%%%%%%%%
\subsection{The infinitesimal Alexander invariant of a Lie algebra}
\label{subsec:infalex}
%%%%%%%%%%%%%%%%%%%%%%
Let $\g$ be a Lie algebra over a field $\k$ of characteristic $0$. We denote by 
$\g'$ its derived subalgebra and by $\g''=(\g')'$ its second derived subalgebra. 
We then have an exact sequence of Lie algebras,
\begin{equation}
\label{eq:liess}
\begin{tikzcd}[column sep=20pt]
0\ar[r]& \g'/\g'' \ar[r]& \g/\g'' \ar[r]& \g/\g' \ar[r]& 0 .
\end{tikzcd}
\end{equation}

The quotient $\g'/\g''$ is abelian, but still carries nontrivial information
through the adjoint action of $\g/\g'$. Let $S=\Sym (\g/\g')$ be the symmetric 
algebra on the vector space $\g/\g'$.
The adjoint representation of $\g/\g'$ on $\g/\g''$ defines 
an $S$-action on $\g'/\g''$, given by 
\[
\bar{g}\cdot \bar{x}=\overline{[g,x]}, \qquad\text{for $g\in \g$, $x\in \g'$},
\] 
where bars denote the corresponding cosets. We define the 
{\em infinitesimal Alexander invariant}\/ of $\g$ to be the $S$-module 
\begin{equation}
\label{eq:infinialex}
\B(\g)=\g'/\g''.
\end{equation} 

\begin{lemma}
\label{lem:bg-generators}
Let $\g$ be a Lie algebra generated by a subset $X\subseteq \g$. Then the 
$S$-module $\B(\g)$ is generated by the classes $\overline{[x,y]}$ with 
$x,y\in X$.  In particular, if $\g$ is finitely generated, then $\B(\g)$ 
is a finitely generated $S$-module.  If, moreover, $\g$ is graded and 
$X=\g_1$, then $\B(\g)_d$ is spanned by the classes of the right-normed 
brackets 
\begin{equation}
\label{eq:right-normed}
[x_1,[x_2,\dots,[x_{d-1},x_d]\dots]]=
\bar{x}_1\cdots \bar{x}_{d-2}\cdot \overline{[x_{d-1},x_d]}, 
\qquad x_1,\dots,x_d\in X .
\end{equation}
\end{lemma}

\begin{proof}
Let $N\subseteq \B(\g)$ be the $S$-submodule generated by the classes 
$\overline{[x,y]}$ with $x,y\in X$.  Since $\g'$ is spanned by the brackets 
$[a,b]$ with $a$ and $b$ iterated brackets in $X$, and since such a bracket 
lies in $\g''$ as soon as both $a$ and $b$ have bracket length at least~$2$, 
it suffices to show that $\overline{[x,b]}\in N$ for $x\in X$ and $b$ an 
iterated bracket in $X$.  We argue by induction on the length $\ell$ of $b$, 
the case $\ell=1$ being the definition of $N$.  

Let $\ell\ge 2$ and write $b=[b_1,b_2]$, so that the Jacobi identity gives
\begin{equation}
\label{eq:jacobi-gen}
[x,[b_1,b_2]]=[[x,b_1],b_2]+[b_1,[x,b_2]] .
\end{equation}
If both $b_1$ and $b_2$ have length at least~$2$, then both terms on the 
right side of \eqref{eq:jacobi-gen} lie in $[\g',\g']=\g''$, and so 
$\overline{[x,b]}=0$.  Otherwise, we may assume (after interchanging $b_1$ 
and $b_2$, if need be) that $b_1=x_1\in X$.  If $\ell\ge 3$, then $b_2$ has 
length at least~$2$; hence the first term on the right side of 
\eqref{eq:jacobi-gen} lies in $\g''$, while the second one has class 
$\bar{x}_1\cdot \overline{[x,b_2]}$, which lies in $N$ by induction.  
Finally, if $\ell=2$, then $b=[x_1,x_2]$ with $x_1,x_2\in X$, and 
\eqref{eq:jacobi-gen} gives 
$\overline{[x,b]}=-\bar{x}_2\cdot \overline{[x,x_1]}
+\bar{x}_1\cdot \overline{[x,x_2]}\in N$.  

The last assertion follows from the first one: if $\g$ is graded and 
$X=\g_1$, then $\B(\g)$ is generated in degree~$2$ by the classes 
$\overline{[x,y]}$ with $x,y\in \g_1$, while $S=\Sym(\g/\g')$ is generated 
in degree~$1$ by the classes of elements of $\g_1$. 
\end{proof}

\begin{lemma}
\label{lem:B-graded}
Let $\g=\bigoplus_{i\ge 1}\g_i$ be a graded Lie algebra. Then 
$\B(\g)$ admits a natural grading making it a graded $S$-module.
\end{lemma}

\begin{proof}
Since $\g$ is graded, the derived subalgebra $\g'=[\g,\g]$ is also graded,
with
\[
\g'=\bigoplus_{i\ge 2} \g'_i, 
\qquad \text{where } \g'_i = \sum_{p+q=i} [\g_p,\g_q].
\]
Likewise, $\g''=[\g',\g']$ is a graded Lie ideal of $\g'$. Hence the quotient
$\B(\g)=\g'/\g''$ inherits a grading,
\[
\B(\g)=\bigoplus_{i\ge 2} \B(\g)_i,
\qquad \B(\g)_i = \g'_i / (\g'' \cap \g'_i).
\]

Next, the quotient $V=\g/\g'$ inherits a grading from $\g$, namely
\[
V=\bigoplus_{i\ge 1} V_i, \qquad V_i = \g_i / (\g_i \cap \g').
\]
Thus the symmetric algebra $S=\Sym(V)$ is a graded commutative algebra.

The adjoint action of $\g$ on $\g'$ preserves degrees: if $g\in \g_p$ and 
$x\in \g'_q$, then $[g,x]\in \g_{p+q}$. Passing to quotients, this induces a 
graded action of $V=\g/\g'$ on $\B(\g)$ of degree $p$, i.e.,
\[
V_p \cdot \B(\g)_q \subseteq \B(\g)_{p+q}.
\]
By the universal property of the symmetric algebra, this action extends uniquely 
to an $S$-module structure on $\B(\g)$, and the above degree condition implies that
\[
S_i \cdot \B(\g)_j \subseteq \B(\g)_{i+j}.
\]
Therefore, $\B(\g)$ is a graded $S$-module.
\end{proof}

In general, though, $\B(\g)$ is an ungraded module. 

\begin{example}
\label{ex:sol2-infalex}
Let $\g = \sol_2 = \langle x, y \mid [x,y]=y\rangle$. 
Here $V=\g/\g' = \k \bar{x}$ and $S=\Sym(V)=\k[x]$, while 
$\g'/\g'' = \k y$. The $S$-module structure is given by $x \cdot y = [x,y] = y$, 
so $(x-1)\cdot y = 0$. Hence, $\B(\g) \cong \k[x]/(x-1)$. 
This module is not graded with respect to any positive grading on $S$. 
We will come back to this computation from a different viewpoint in 
Example~\ref{ex:nonhomog-tre}.
\end{example}

%%%%%%%%%%%%%%%%%%%%%%
\subsection{The invariant $\B(\h(A))$, and metabelian quotients as extensions}
\label{subsec:infalex-metabelian}
%%%%%%%%%%%%%%%%%%%%%%

For a connected $\cdga$ $A$, we apply this construction to 
$\g = \h(A)$, the holonomy Lie algebra introduced in 
Section~\ref{subsec:cochain-holo-short}, giving the infinitesimal 
Alexander invariant $\B(A) \coloneqq \B(\h(A))$. 
The case when $\g$ is a finitely generated, quadratic Lie algebra was first 
considered in \cite{PS-imrn04}, and further studied in \cite{SW-mz, SW-aam}; 
in this case, $\B(\g)$ is a finitely generated graded $S$-module. 

The construction is functorial: a morphism $\varphi\colon \g \to \bar{\g}$ induces an 
$S$-linear map $\B(\varphi)\colon \B(\g) \to \B(\bar{\g})$, where $\B(\bar{\g})$ is 
viewed as an $S$-module via restriction of scalars along the induced map 
$S=\Sym(\g/\g') \to \Sym(\bar{\g}/\bar{\g}')$.

We record next a way of reassembling $\g/\g''$ from the pair $\bigl(\g/\g',\B(\g)\bigr)$.
Both ends of \eqref{eq:liess} are abelian, so that sequence is an abelian extension, and
such extensions are governed by a cohomology class.  This is elementary, but it is what
converts the presentation of $\B_1(A)$ obtained in Section~\ref{subsec:Bpres} into a
presentation of a Lie algebra, and it is used in that form in the proof of
Theorem~\ref{thm:B-holo}.

Let $V$ be a finite-dimensional $\k$-vector space, $S=\Sym(V)$, and $M$ an $S$-module.
For a $\k$-linear map $c\colon \bwedge^2 V\to M$, let $\fE(V,M,c)$ denote the vector space
$V\oplus M$ equipped with the alternating bracket determined by
\begin{equation}
\label{eq:fm-bracket-gen}
[v,w]=c(v\wedge w), \qquad [v,m]=v\cdot m, \qquad [m,m']=0 ,
\end{equation}
for $v,w\in V$ and $m,m'\in M$, where $v\cdot m$ refers to the $S$-module structure.
Recall that the Chevalley--Eilenberg complex of the abelian Lie algebra $V$ with
coefficients in $M$ is $\Hom_\k(\bwedge^{\bullet}V,M)$, with differential
\begin{equation}
\label{eq:ce-abelian-diff}
(\delta\tau)(v_1,\dots,v_{p+1})=
\sum_{i}(-1)^{i+1}\,v_i\cdot \tau(v_1,\dots,\widehat{v_i},\dots,v_{p+1}) ,
\end{equation}
the Koszul differential with coefficients in $M$; we write $H^{*}(V;M)$ for its
cohomology.

\begin{lemma}
\label{lem:metabelian-extension}
With the notation above:
\begin{enumerate}[itemsep=2pt]
\item \label{me1}
$\fE(V,M,c)$ is a Lie algebra if and only if $\delta c=0$, and it is then metabelian.
If moreover $M$ is generated as an $S$-module by $\im(c)$, then
$\fE(V,M,c)'=M$ and $\fE(V,M,c)/\fE(V,M,c)'=V$.
\item \label{me2}
Let $\g$ be a Lie algebra with $V=\g/\g'$ and $M=\B(\g)$.  Every $\k$-linear splitting
$\sigma\colon V\to \g/\g''$ of \eqref{eq:liess} determines a cocycle
$c_\sigma\in \Hom_\k(\bwedge^2V,M)$ by $c_\sigma(v\wedge w)=[\sigma v,\sigma w]$, and an
isomorphism of Lie algebras $\g/\g''\cong \fE(V,M,c_\sigma)$.
\item \label{me3}
Two splittings differ by a map $\tau\colon V\to M$, and then
$c_{\sigma+\tau}=c_\sigma+\delta\tau$.  Hence the class
$[c_\sigma]\in H^2(V;M)$ is independent of the splitting, and it determines
$\g/\g''$ up to isomorphism.
\item \label{me4}
The sequence \eqref{eq:liess} splits in the category of Lie algebras if and only if
$[c_\sigma]=0$.  If that happens and $\g$ is generated by any linear lift of $V$, then
$M=\m M$, where $\m=(V)\subset S$; so if in addition $M$ is a finitely generated graded
$S$-module, then $\B(\g)=0$.
\end{enumerate}
\end{lemma}

\begin{proof}
\eqref{me1} The bracket is alternating by construction.  Of the two Jacobi identities to
check, the one on $V\oplus V\oplus M$ reads $v\cdot(w\cdot m)-w\cdot(v\cdot m)=0$ and holds
because $S$ is commutative, while the one on $\bwedge^3V$ reads
$\sum_{\mathrm{cyc}} v_1\cdot c(v_2\wedge v_3)=0$, which by
\eqref{eq:ce-abelian-diff} is precisely $\delta c=0$.  As $[M,M]=0$ and
$\fE'\subseteq M$, we get $\fE''=0$.  The last assertion is immediate, since
$\fE'=\im(c)+\m M$.

\eqref{me2} A splitting exists, \eqref{eq:liess} being an exact sequence of vector spaces,
and $c_\sigma$ takes values in $M$ because $\pi[\sigma v,\sigma w]=[v,w]=0$ in the abelian
Lie algebra $V$.  That $c_\sigma$ is a cocycle follows from the Jacobi identity in
$\g/\g''$, and the map $v+m\mapsto \sigma v+m$ is then an isomorphism
$\fE(V,M,c_\sigma)\to\g/\g''$ of Lie algebras, by inspection of
\eqref{eq:fm-bracket-gen}: the bracket $[\sigma v,m]$ in $\g/\g''$ is the $S$-action
$v\cdot m$, by the definition of that action.

\eqref{me3} $c_{\sigma+\tau}(v\wedge w)=[\sigma v+\tau v,\sigma w+\tau w]
=c_\sigma(v\wedge w)+v\cdot\tau w-w\cdot\tau v$, since $[\tau v,\tau w]=0$; and the last
two terms are $(\delta\tau)(v,w)$.

\eqref{me4} A Lie algebra splitting is a splitting $\sigma$ with $c_\sigma=0$, and
by~\eqref{me3} such a $\sigma$ exists exactly when $[c_\sigma]=0$.  Given one, $\sigma(V)$
is an abelian subalgebra, so if $\g/\g''$ is generated by $\sigma(V)$ then
$M=[\g/\g'',\g/\g'']=[\sigma(V),M]=\m M$.  For $M$ finitely generated and graded, the
graded Nakayama lemma gives $M=0$.
\end{proof}

\begin{remark}
\label{rem:metabelian-extension-caveats}
Two comments on \eqref{me4}.  The hypothesis that $M$ be graded is not decoration:
Example~\ref{ex:sol2-infalex} has $\B(\sol_2)=\k[x]/(x-1)$, on which $\m=(x)$ acts
invertibly, so $\m M=M$ with $M\ne0$.  And indeed \eqref{eq:liess} does split there, since
$\sol_2=\k x\ltimes \k y$---as it must, $\bwedge^2V$ vanishing for $\dim_\k V=1$, so that
$c_\sigma=0$ for every splitting.  Second, only the \emph{class}\/ of $c_\sigma$ is
canonical, not the cocycle: constructions built from a chosen splitting are therefore
functorial only up to the coboundary correction of~\eqref{me3}.  This is the source of the
correction term $\theta_\varphi$ of \eqref{eq:theta-phi}, and of the caveat in
Remark~\ref{rem:bphi-triangular}.
\end{remark}

%%%%%%%%%%%%%%%%%%%%%%
\subsection{A presentation for $\B(\g)$ when $\g$ is graded}
\label{subsec:infalex-graded}
%%%%%%%%%%%%%%%%%%%%%%

In this subsection, we give a self-contained derivation of a presentation 
for the infinitesimal Alexander invariant $\B(\g)$ in the graded case. 
The proof of the corresponding result in \cite[Thm.~6.2]{PS-imrn04} proceeds 
by a dimension count, comparing the ranks of $L'/L''$ for a free Lie algebra 
$L$---known from Chen's computation \cite{Chen51} of the ranks $\theta_k(F_n)$ 
of the lower central series quotients of free metabelian groups---with those 
of a Koszul cokernel.  We replace that count by an explicit construction: the 
free metabelian Lie algebra maps to a wreath product of abelian Lie algebras, 
in which the Koszul differential $\partial_2^E$ becomes visible, and the kernel 
of the presenting map can be read off directly; see 
Remark~\ref{rem:wreath-magnus} for the provenance of this device.  
The same mechanism will drive the identification $\B_1(A)\cong \B(\h(A))$ 
for $1$-finite $\cdgas$ in Theorem~\ref{thm:B-holo}.

We start by identifying the infinitesimal Alexander invariant of a finitely 
generated free Lie algebra, thereby recovering 
\cite[Eq.~(6.4)]{PS-imrn04} without recourse to the Chen ranks of free groups.

\begin{proposition}
\label{prop:holo-free}
Let $W$ be a finite-dimensional $\k$-vector space, and let 
$L=\Lie(W)$ be the free Lie algebra on $W$. 
Write $E_{\bullet}=\bigwedge W$ and $S=\Sym(W)$. 
Then the infinitesimal Alexander invariant of $L$ is isomorphic, 
as a graded $S$-module, to the cokernel of the third Koszul differential,
\[
\B(L)\cong \coker \Big( \partial_3^E\colon E_3\otimes_{\k} S 
\longrightarrow E_2\otimes_{\k} S \Big).
\]
\end{proposition}

\begin{proof}
Since $E_2\otimes_{\k}S$ is a free $S$-module on $E_2=\bwedge^2 W$, 
there is a well-defined $S$-linear map 
\begin{equation}
\label{eq:eta-L}
\eta_L\colon E_2 \otimes_{\k} S \longrightarrow L'/L'', 
\qquad u\wedge v\longmapsto \overline{[u,v]} \quad (u,v\in W) ,
\end{equation}
which is surjective by Lemma~\ref{lem:bg-generators}.  The Jacobi identity 
gives $\eta_L\circ \partial_3^E=0$; indeed, 
\[
\eta_L\big(\partial_3^E(u\wedge v\wedge w\otimes 1)\big)=
\overline{[u,[v,w]]}+\overline{[v,[w,u]]}+\overline{[w,[u,v]]}=0 ,
\]
for all $u,v,w\in W$.  It remains to prove the reverse inclusion, 
$\ker (\eta_L)\subseteq \im (\partial_3^E)$.

To that end, consider the $\k$-vector space 
$\fw=W\oplus \big(W\otimes_{\k}S\big)$, endowed with the 
bracket determined by 
\begin{equation}
\label{eq:wreath}
[u,v]=0, \qquad [u,v\otimes f]=v\otimes uf, \qquad 
[u\otimes f,\, v\otimes g]=0 ,
\end{equation}
for $u,v\in W$ and $f,g\in S$.  The bracket is alternating, and the 
only Jacobi identity requiring verification reads $u(vf)-v(uf)=0$, which 
holds since $S$ is commutative.  Thus, $\fw$ is a Lie algebra; it is 
metabelian, with $\fw'\subseteq W\otimes_{\k}S$ and $\fw''=0$.  Since 
$L$ is free on $W$, there is a unique Lie algebra morphism 
\begin{equation}
\label{eq:rho-wreath}
\rho\colon L\longrightarrow \fw, 
\qquad \rho(u)=\big(u,\, u\otimes 1\big) \quad (u\in W) ,
\end{equation}
and a direct computation with \eqref{eq:wreath} gives 
\begin{equation}
\label{eq:rho-bracket}
\rho\big([u,v]\big)=\big(0,\, v\otimes u-u\otimes v\big)=
\big(0,\, \partial_2^E(u\wedge v\otimes 1)\big) .
\end{equation}

Since $\fw''=0$, we have $\rho(L'')=0$, while 
$\rho(L')\subseteq \fw'\subseteq W\otimes_{\k}S$.  Hence $\rho$ 
induces a $\k$-linear map $\iota\colon L'/L''\to W\otimes_{\k}S$, 
which is $S$-linear: if $u\in W$ and $\xi\in L'$, then 
$\iota\big(\bar{u}\cdot \bar{\xi}\big)=\rho([u,\xi])=
[\rho(u),\iota(\bar\xi)]=u\cdot \iota(\bar\xi)$.  
By \eqref{eq:rho-bracket}, the $S$-linear maps $\iota\circ \eta_L$ and 
$\partial_2^E$ agree on $E_2$, and so 
\begin{equation}
\label{eq:iota-eta}
\iota\circ \eta_L=\partial_2^E .
\end{equation}
Therefore, if $z\in \ker (\eta_L)$, then $\partial_2^E(z)=
\iota\big(\eta_L(z)\big)=0$, that is, $z\in \ker (\partial_2^E)=
\im (\partial_3^E)$, by the exactness of the Koszul complex.  We conclude 
that $\ker (\eta_L)=\im (\partial_3^E)$, and the claimed presentation 
follows.
\end{proof}

\begin{corollary}
\label{cor:bl-syzygy}
The map $\iota\colon \B(L)\to W\otimes_{\k} S$ induced by the morphism 
\eqref{eq:rho-wreath} restricts to an isomorphism of graded $S$-modules,
\[
\begin{tikzcd}[column sep=22pt]
\iota\colon \B(L) \arrow[r, "\simeq"] & 
\ker \big( \partial_1^E\colon E_1\otimes_{\k} S \to S \big) ,
\end{tikzcd}
\]
onto the module of syzygies among the variables of $S$.
\end{corollary}

\begin{proof}
By Proposition~\ref{prop:holo-free} and its proof, the map $\eta_L$ is 
surjective, with $\ker (\eta_L)=\im (\partial_3^E)=\ker (\partial_2^E)$. 
Hence, \eqref{eq:iota-eta} shows that $\iota$ is injective, with image 
$\im (\partial_2^E)=\ker (\partial_1^E)$, where the last equality holds 
by the exactness of the Koszul complex.
\end{proof}

\begin{remark}
\label{rem:wreath-magnus}
The Lie algebra \eqref{eq:wreath} is the wreath product of two copies of the
abelian Lie algebra $W$, in the sense of \v{S}mel'kin \cite{Shmelkin}: for Lie
algebras $\fa$ and $\fb$ one sets
$\fa\wr\fb=\fb\ltimes\bigl(\fa\otimes_{\k} U(\fb)\bigr)$, and here
$U(W)=\Sym(W)=S$.  It is the Lie analogue of the group-theoretic wreath
product $A\wr B=B\ltimes A^{(B)}$, with $A^{(B)}$ free over $\Z[B]$, which is
why the \emph{free}\/ $S$-module $W\otimes_{\k}S$ occurs and not a smaller
one; correspondingly, $\rho$ is the Lie analogue of the Magnus embedding of a
free metabelian group, for which see \cite[\S 3.4]{MKS}.  Note that
$\fw$ is not the free metabelian Lie algebra on $W$: since
$\fw'=W\otimes_{\k}\m$, one has $\fw/\fw'\cong W\oplus W$, so $\fw$ requires
$2\dim_\k W$ generators.

Two features of the argument are worth isolating.  First, only the
\emph{existence}\/ of $\rho$ is used, and that is free of charge, since
$L$ is free on $W$; in particular \v{S}mel'kin's embedding theorem is not
invoked.  Second, the theorem is instead \emph{obtained}: $\ker(\rho)=L''$,
because the first coordinate of $\rho(\xi)$ recovers the class of $\xi$ in
$L/L'$, so that $\rho(\xi)=0$ forces $\xi\in L'$, whereupon
$\rho(\xi)=\iota(\bar\xi)$ and $\iota$ is injective.  Corollary
\ref{cor:bl-syzygy} thus identifies the image of the free metabelian Lie
algebra $L/L''$ precisely:
\[
\rho\bigl(L/L''\bigr)=\Delta(W)\oplus\ker(\partial_1^E),
\qquad \Delta(u)=(u,\,u\otimes 1) .
\]
The diagonal form of $\Delta$ is forced: writing $\rho(u)=(u,m_u)$ gives
$\rho([u,v])=\bigl(0,\,u\cdot m_v-v\cdot m_u\bigr)$, and the choice
$m_u=u\otimes1$ is what makes this equal $\partial_2^E(u\wedge v\otimes1)$
on the nose.

\v{S}mel'kin proves the embedding in greater generality, for $F/V(R)$ with
$F$ free, $R$ an arbitrary ideal, and $V$ any variety of Lie algebras; the
case at hand is $V$ abelian and $R=F'$.  In that case he describes the image
as consisting of the elements $\sum_{i} a_iP_i+b$ with $b\in F/R$, $P_i\in
U(F/R)$, and $\sum_{i} x_iP_i=b$, which is the equational form of the
decomposition displayed above: the condition $\sum_{i} x_iP_i=b$ says precisely
that $\partial_1^E$ carries the $W\otimes_{\k}S$-component to the
$W$-component.
\end{remark}

It is in this form---with the map $\iota$ computed by Fox calculus---that 
the infinitesimal Alexander invariant of a free Lie algebra enters the 
proof of \cite[Prop.~9.3]{PS-imrn04}.  Corollary~\ref{cor:bl-syzygy} thus 
recovers that description directly, bypassing Chen's computation 
\cite{Chen51} of the ranks $\theta_k(F_n)$.

Let $\g$ be a finitely generated, graded Lie algebra over a field 
$\k$ of characteristic $0$, generated in degree $1$, and set $V = \g/\g'$ 
and $S = \Sym(V)$.  As we saw in Lemma~\ref{lem:B-graded}, the infinitesimal 
Alexander invariant $\B(\g)$ is a graded $S$-module, and it is finitely 
generated, by Lemma~\ref{lem:bg-generators}.  We now exhibit an explicit 
presentation of this module---a finite one, as soon as $\g$ is finitely 
presented---generalizing the one given in \cite[Thm.~6.2]{PS-imrn04} for 
quadratic Lie algebras. 
Throughout, we write $L=\Lie(V)$, $E_{\bullet}=\bwedge V$, and we let 
$\eta_L\colon E_2\otimes_{\k}S\to \B(L)$ be the surjection \eqref{eq:eta-L} 
attached to $L$ by Proposition~\ref{prop:holo-free}, applied with $W=V$.

\begin{proposition}
\label{prop:bg-graded}
Let $\g=L/\langle \fr\rangle$ be a graded Lie algebra generated in 
degree~$1$, where $L=\Lie(V)$ with $V=\g/\g'$ finite-dimensional, and where 
$\fr\subseteq \bigoplus_{d\ge 2}L_d$ is a graded subspace.  Let 
$\varrho\colon \fr\to E_2\otimes_{\k}S$ be a homogeneous $\k$-linear map 
lifting the class map $\fr\to \B(L)$, $r\mapsto \bar{r}$, along $\eta_L$. 
Then $\B(\g)$ admits the homogeneous presentation
\begin{equation}
\label{eq:bg-pres-graded}
\begin{tikzcd}[column sep=30pt]
\big(E_3\oplus \fr \big)\otimes_{\k} S 
\arrow[r, "{(\partial_3^{E},\, \varrho)}"] &
E_2 \otimes_{\k} S \arrow[r, "\eta_\g"] &
\B(\g) \arrow[r] & 0,
\end{tikzcd}
\end{equation}
where $\deg E_k=k$, $\deg \fr_d=d$, and 
$\eta_\g(v\wedge w)=\overline{[v,w]}$ for $v,w\in V$.  Moreover, 
$\im (\partial_3^{E},\varrho)=\ker (\eta_\g)$; in particular, this 
submodule does not depend on the choice of the lift $\varrho$.
\end{proposition}

\begin{proof}
Such a lift $\varrho$ exists, since $\eta_L$ is a degree-preserving 
surjection of graded $\k$-vector spaces and $\fr$ is graded. 

Set $I=\langle \fr\rangle$.  As $\fr$ consists of elements of bracket length 
at least~$2$, we have $I\subseteq L'$, and therefore $\g'=L'/I$, while 
$\g''=(L''+I)/I$, since $[L',I]\subseteq I$ and $[I,I]\subseteq I$.  
Consequently, 
\begin{equation}
\label{eq:bg-quotient}
\B(\g)=\g'/\g''=L'/(L''+I)=
\B(L)\big/\im\big(I\longrightarrow \B(L)\big) .
\end{equation}
We claim that the image of $I$ in $\B(L)$ is the $S$-submodule generated by 
the classes $\bar{r}$, with $r\in \fr$.  Indeed, $I$ is spanned by the 
iterated brackets $[y_1,[y_2,\dots,[y_k,r]\dots]]$ with $y_j\in L$ and 
$r\in \fr$.  If some $y_j$ lies in $L'$, then, since $r\in L'$, the 
corresponding element lies in $[L',L']=L''$; hence, modulo $L''$, only the 
degree-$1$ components of the $y_j$ contribute, and those act through the 
$S$-module structure of $\B(L)$.  This proves the claim, and yields the 
exact sequence of graded $S$-modules
\begin{equation}
\label{eq:bg-relations}
\begin{tikzcd}[column sep=24pt]
S\otimes_{\k} \fr \arrow[r, "r\mapsto \bar{r}"] & \B(L) \arrow[r] & 
\B(\g)\arrow[r] & 0 .
\end{tikzcd}
\end{equation}

By construction, $\eta_\g$ is the composite of $\eta_L$ with the projection 
$\B(L)\twoheadrightarrow \B(\g)$; in particular, $\eta_\g$ is surjective. 
Moreover, by \eqref{eq:bg-relations}, an element $z\in E_2\otimes_{\k}S$ 
lies in $\ker (\eta_\g)$ if and only if $\eta_L(z)$ belongs to the 
$S$-submodule generated by $\{\bar{r}\}_{r\in \fr}=
\eta_L\big(\varrho(\fr)\big)$, that is, if and only if 
$z\in S\cdot \varrho(\fr)+\ker (\eta_L)$.  Since 
$\ker (\eta_L)=\im (\partial_3^E)$ by Proposition~\ref{prop:holo-free}, 
this says precisely that $\ker (\eta_\g)=\im (\partial_3^E,\varrho)$, 
thereby completing the proof.
\end{proof}

\begin{remark}
\label{rem:bg-lift}
The lift $\varrho$ may be written down explicitly.  By 
Lemma~\ref{lem:bg-generators}, each $r\in \fr_d$ satisfies 
$\bar{r}=\sum_{j} c_j\, v^j_1\cdots v^j_{d-2}\cdot 
\overline{[v^j_{d-1},v^j_{d}]}$ in $\B(L)$, for suitable $c_j\in \k$ and 
$v^j_i\in V$; one may then take $\varrho(r)=\sum_{j} c_j\, 
v^j_1\cdots v^j_{d-2}\otimes v^j_{d-1}\wedge v^j_{d}$.  When $\g$ is 
quadratic, that is, when $\fr\subseteq L_2=\bwedge^2 V$, no choice is 
involved: $\varrho$ is the inclusion 
$\fr\inj \bwedge^2 V=E_2\subseteq E_2\otimes_{\k} S$, and 
\eqref{eq:bg-pres-graded} reduces to the presentation of 
\cite[Thm.~6.2]{PS-imrn04}.  Alternatively, and without any choices, 
combining \eqref{eq:bg-relations} with Corollary~\ref{cor:bl-syzygy} identifies 
$\B(\g)$ with the cokernel of the map $S\otimes_{\k} \fr\to 
\ker (\partial_1^E)$, $r\mapsto \iota(\bar{r})$; it is in this form that the 
invariant appears in \cite[\S 9]{PS-imrn04}, where the elements 
$\iota(\bar{r})$ are computed by Fox calculus for commutator-relators groups. 
\end{remark}

For quadratic Lie algebras---in particular, for holonomy Lie algebras---the
support of the infinitesimal Alexander invariant is not a new invariant: it
is a support resonance variety in the sense of
Section~\ref{subsec:res}.  We record this, since it is the form in
which the support enters the growth estimates for holonomy Chen ranks.

\begin{corollary}
\label{cor:supp-B-quadratic}
Let $H$ be a connected graded $\k$-algebra with $\dim_\k H^1<\infty$.  Then
\begin{equation}
\label{eq:supp-B-holo}
\Supp_S\B\bigl(\h(H)\bigr)=\RR_{1}\bigl((H,0)\bigr) .
\end{equation}
Moreover, every quadratic graded Lie algebra arises this way: if
$\g=\Lie(V)/\langle K\rangle$ with $K\subseteq\bwedge^2V$, then
$\g\cong\h(H)$ for $H=\k\oplus V^{\vee}\oplus\bigl(\bwedge^2V^{\vee}/
K^{\perp}\bigr)$, and \eqref{eq:supp-B-holo} computes $\Supp_S\B(\g)$.
\end{corollary}

\begin{proof}
Apply Theorem~\ref{thm:B-holo} to the $\cdga$ $A=(H,0)$, whose holonomy Lie
algebra is $\h(H)$: this gives an isomorphism of $S$-modules
$\B(\h(H))\cong\B_1(A)$, and $\RR_1(A)=\Supp_S\B_1(A)$ by definition.  For
the second assertion, note that the algebra $H$ so defined is graded
commutative and associative, since all triple products of positive-degree
elements vanish, and that $\mu_H$ is the projection
$\bwedge^2V^{\vee}\surj\bwedge^2V^{\vee}/K^{\perp}$, so that
$\im(\mu_H^{\vee})=K$ and $\h(H)=\Lie(V)/\langle K\rangle=\g$.
\end{proof}

\begin{remark}
\label{rem:supp-vs-jump-quadratic}
In practice one computes not the support locus but the \emph{jump}\/ locus
$\RR^{1}\bigl((H,0)\bigr)$, the determinantal locus where the Aomoto complex
of $H$ fails to be exact in degree~$1$---for a $2$-step nilpotent Lie
algebra, by \eqref{eq:res-2step}, a union of $2$-planes indexed by the
decomposable elements of $W$.  In degree~$1$, however, the two loci can
differ only at the origin.  Indeed, $\B_0\bigl((H,0)\bigr)\cong\k$, so
$\RR_0=\{0\}=\RR^0$, and Theorem~\ref{thm:PS-mrl-refined}\eqref{ps2} with
$j=1$ gives $\{0\}\cup\RR_1=\{0\}\cup\RR^1$.  The discrepancy at the origin
can be genuine, but only in the degenerate case $\B(\h(H))=0$, that is, when
$\h(H)$ is abelian: for $H=H^{*}(T^2;\k)$ one has $\RR_1=\emptyset$ while
$\RR^{1}=\{0\}$; see Remark~\ref{rem:jump-support-strict} and
Example~\ref{ex:torus-jump-support}.  The distinction is immaterial for the applications:
it does not affect $\dim$, hence not the growth rate of the holonomy Chen
ranks, and it disappears altogether when the resonance is homogeneous, since
a nonempty cone contains the origin---which is automatic under
positive weights, by Corollary~\ref{cor:conical}.
\end{remark}

%%%%%%%%%%%%%%%%%%
\subsection{A comparison map}
\label{subsec:compare}
%%%%%%%%%%%%%%%%%%

We now return to the holonomy Lie algebra $\h=\h(A)$ of a $\k$-$\cdga$ 
$(A,d_A)$ with $\dim_{\k} A^1<\infty$. As in \eqref{eq:A1A1}, write 
$A_1=E_1\oplus U_1$, where $U_1=\im(d_A^{\vee})$ is canonical and $E_1$ is 
a chosen complement, identified with $H^1(A)^{\vee}=H_1(A)$.  
Since $\h (A)= \Lie(A_1) /\langle \im (d^{\vee}_A + \mu^{\vee}_A)\rangle$, 
where $\im(\mu_A^{\vee})\subset \Lie_2(A_1)\subset \Lie'(A_1)$, 
it follows that $\h/\h'=E_1$.

Consider the infinitesimal Alexander invariant  $\B(\h(A))$, viewed as a 
module over the symmetric algebra $S=\Sym(E_1)$. As we show 
next, this module is isomorphic to the Koszul module of $A$. 
This result generalizes \cite[Prop.~9.3]{PS-imrn04}, which 
was only proved under the assumption that $A=H^*(G;\k)$, where 
$G$ is a finitely presented, commutator-relators group, and $d_A=0$. 
This establishes Theorem~\ref{thm:B-holo-intro} from the Introduction.

\begin{theorem}
\label{thm:B-holo}
Let $(A,d_A)$ be a connected $\k$-$\cdga$ with $A^1$ 
finite-dimensional, and let $S=\Sym(E_1)$. There is 
then an isomorphism of $S$-modules, 
\begin{equation}
\label{eq:B1-holo}
\B_1(A)\longisom \B(\h(A)). 
\end{equation}
\end{theorem}

\begin{proof}
By Theorem~\ref{thm:Bpres}, the module $\B_1(A)$ is presented as
the cokernel of the $S$-linear map
\begin{equation}
\label{eq:BA-pres}
\begin{tikzcd}[column sep=24pt]
(E_3\oplus A_2)\otimes_{\k} S \arrow[r, "\Delta_A"] &
(E_2\oplus U_1)\otimes_{\k} S,
\end{tikzcd}
\qquad
\Delta_A=
\begin{pmatrix}
\partial_3^{E} & 0 \\
\nu_A^{\vee}\!\otimes\!\id_S &
d_A^{\vee}\!\otimes\!\id_S + \beta_A^{\vee}
\end{pmatrix},
\end{equation}
where $\nu_A^{\vee}$ and $\beta_A^{\vee}$ are as in
Section~\ref{subsec:Bpres}, and where, as there, the rows of the matrix are
indexed by the summands of the source and the columns by those of the target.

Let $\h=\h(A)$, and let $L=\Lie(E_1)$ and $\bm{L}=\Lie(A_1)$
be the free Lie algebras on $E_1$ and $A_1$, respectively,
with $L\subset\bm{L}$ a Lie subalgebra.

\smallskip
\emph{Step 1: Two comparison maps.}
Since $E_2\otimes_{\k}S$ is a free $S$-module on $E_2$, there is a 
well-defined $S$-linear map
\begin{equation}
\label{eq:etaA}
\eta_A\colon E_2\otimes_{\k} S \longrightarrow \h'/\h'',
\qquad u\wedge v \longmapsto [u,v]\bmod\h'' \quad (u,v\in E_1),
\end{equation}
namely, the composite of the map \eqref{eq:eta-L} attached to $L$ with the 
canonical maps $L'/L''\to\bm{L}'/\bm{L}''\surj\h'/\h''$.  By the Jacobi 
identity,
\begin{equation}
\label{eq:etaA-jacobi}
\eta_A\circ \partial_3^E=0 .
\end{equation}

Since the defining relations of $\h(A)$ identify
$d_A^\vee(\gamma) = -\mu_A^\vee(\gamma)$ in $\h$, and the right-hand side
$-\mu_A^\vee(\gamma)$ is a sum of brackets $[a,b]$ with $a,b\in A_1$
(hence lies in $\h'=[\h,\h]$), every element $d_A^\vee(\gamma)$ of
$U_1=\im(d_A^\vee)$ lies in $\h'$. We may therefore define a 
natural $S$-linear map
\begin{equation}
\label{eq:psiA}
\psi_A\colon U_1\otimes_{\k} S \longrightarrow \h'/\h'',
\qquad u \longmapsto \bar{u},
\end{equation}
where $\bar{u}$ denotes the class of $u\in U_1\subset\h'$
in $\h'/\h''$.  In what follows, we also write 
$\kappa\colon A_1\wedge A_1\to\h'/\h''$ for the map 
$a\wedge b\mapsto[a,b]\bmod\h''$.

\smallskip
\emph{Step 2: Surjectivity of the combined map.}
Since $\h$ is generated by $A_1$, Lemma~\ref{lem:bg-generators} shows that 
the $S$-module $\h'/\h''$ is generated by the classes of the brackets 
$[x,y]$ with $x,y\in A_1=E_1\oplus U_1$.
If $x,y\in E_1$, then $[x,y]\in\im(\eta_A)$.
If $x\in E_1$ and $y=u\in U_1$, then
$[x,u]=x\cdot\bar{u}\in S\cdot\im(\psi_A)$.
If $x,y\in U_1$, then $[x,y]\in[\h',\h']=\h''$ vanishes in $\h'/\h''$.
Hence the combined $S$-linear map
\begin{equation}
\label{eq:eta-psi}
(\eta_A\:\:\psi_A)\colon (E_2\oplus U_1)\otimes_{\k} S \longsurj \h'/\h''
\end{equation}
is surjective, and the induced map
$\phi\colon\B_1(A)=\coker(\Delta_A)\twoheadrightarrow\h'/\h''=\B(\h)$
is a well-defined surjection, provided $(\eta_A\:\:\psi_A)\circ\Delta_A=0$.

\smallskip
\emph{Step 3: The kernel of $(\eta_A\:\:\psi_A)$ equals $\im(\Delta_A)$,
and $\phi$ is an isomorphism.}

We first verify $(\eta_A\:\:\psi_A)\circ\Delta_A=0$.
Decompose $\mu_A^\vee(\gamma)=\nu_A^\vee(\gamma)+\sigma(\gamma)+\tau(\gamma)$
for each $\gamma\in A_2$, where
$\nu_A^\vee(\gamma)\in E_1\wedge E_1=E_2$,
$\sigma(\gamma)\in E_1\wedge U_1$,
and $\tau(\gamma)\in U_1\wedge U_1$.
Since $d_A^\vee(\gamma)=-\mu_A^\vee(\gamma)$ in $\h$, one has
\[
\psi_A(d_A^\vee(\gamma))
= \bar{d_A^\vee(\gamma)}
= -\kappa(\nu_A^\vee(\gamma)) - \kappa(\sigma(\gamma)) - \kappa(\tau(\gamma))
= -\eta_A(\nu_A^\vee(\gamma)) - \kappa(\sigma(\gamma))
\]
in $\h'/\h''$, where $\kappa(\tau(\gamma))=0$ because
$\tau(\gamma)\in U_1\wedge U_1$ and $[U_1,U_1]\subset[\h',\h']=\h''$.
From the definition of $\beta_A^\vee=(\pi_U\otimes\id)\circ\omega_A^\vee$,
a direct computation gives
$\beta_A^\vee(\gamma\otimes 1)=\sum_{i}\pi_U(\iota_{e^i}(\gamma))\otimes x_i$,
where $\{e^i\}$ is any basis of $E^1$ with dual basis $\{x_i\}$ of $E_1\subset S$.
The $E_1\wedge U_1$ cross terms of $\mu_A^\vee(\gamma)$ are
$\sigma(\gamma)=\sum_{i} e_i\wedge\pi_U(\iota_{e^i}(\gamma))$,
so $\kappa(\sigma(\gamma))=\sum_{i} x_i\cdot\overline{\pi_U(\iota_{e^i}(\gamma))}
=\psi_A(\beta_A^\vee(\gamma\otimes 1))$.
Therefore
\begin{equation}
\label{eq:compat}
\eta_A(\nu_A^\vee(\gamma))+
\psi_A\bigl((d_A^\vee\otimes\id_S+\beta_A^\vee)(\gamma\otimes 1)\bigr)=0
\quad\text{for all }\gamma\in A_2 .
\end{equation}
On the summand $E_3$, the corresponding verification is 
\eqref{eq:etaA-jacobi}.  Together, these two computations give 
$(\eta_A\:\:\psi_A)\circ\Delta_A=0$ and hence
$\im(\Delta_A)\subseteq\ker(\eta_A\:\:\psi_A)$, confirming that $\phi$ is defined.

It remains to show that $\phi$ is injective, and we do so by constructing 
an inverse.  Let 
$q\colon (E_2\oplus U_1)\otimes_{\k} S\surj \B_1(A)=\coker(\Delta_A)$ 
be the canonical projection, and consider the $\k$-vector space 
\begin{equation}
\label{eq:fM-def}
\fM=E_1\oplus \B_1(A) ,
\end{equation}
endowed with the alternating bracket determined by 
\begin{equation}
\label{eq:fM-bracket}
[e,f]=q\big(e\wedge f\otimes 1,\, 0\big), \qquad 
[e,m]=e\cdot m, \qquad [m,m']=0 ,
\end{equation}
for $e,f\in E_1$ and $m,m'\in \B_1(A)$, where $e\cdot m$ refers to the 
$S$-module structure of $\B_1(A)$.  In the notation of
Section~\ref{subsec:infalex-metabelian}, this is $\fM=\fE(E_1,\B_1(A),c)$ with
$c(e\wedge f)=q(e\wedge f\otimes 1,0)$, and $c$ is a cocycle: by
\eqref{eq:ce-abelian-diff},
\[
(\delta c)(e,f,g)=
q\big(f\wedge g\otimes x_e-e\wedge g\otimes x_f+e\wedge f\otimes x_g,\, 0\big)
=q\big(\partial_3^E(e\wedge f\wedge g\otimes 1),\, 0\big)=0 ,
\]
since $\im (\Delta_A)$ contains $\im (\partial_3^E)\oplus 0$.  Hence $\fM$ is a
metabelian Lie algebra, by
Lemma~\ref{lem:metabelian-extension}\eqref{me1}; in particular
$\fM'\subseteq \B_1(A)$ and $\fM''=0$.

Since $\Lie(A_1)$ is free, there is a unique Lie algebra morphism 
\begin{equation}
\label{eq:rho-fM}
\rho\colon \Lie(A_1)\longrightarrow \fM, \qquad 
\rho(a)=\pi_E(a)+q\big(0,\, \pi_U(a)\otimes 1\big) \quad (a\in A_1) .
\end{equation}
We claim that $\rho$ annihilates the relations of $\h(A)$.  Indeed, let 
$\gamma\in A_2$ and decompose $\mu_A^{\vee}(\gamma)=\nu_A^{\vee}(\gamma)
+\sigma(\gamma)+\tau(\gamma)$ as above.  Then 
$\rho\big(d_A^{\vee}(\gamma)\big)=q\big(0,d_A^{\vee}(\gamma)\otimes 1\big)$, 
since $d_A^{\vee}(\gamma)\in U_1$, while 
\[
\rho\big(\nu_A^{\vee}(\gamma)\big)=
q\big(\nu_A^{\vee}(\gamma)\otimes 1,\, 0\big), \qquad 
\rho\big(\sigma(\gamma)\big)=
q\big(0,\, \beta_A^{\vee}(\gamma\otimes 1)\big), \qquad 
\rho\big(\tau(\gamma)\big)=0 ;
\]
here the first equality holds by \eqref{eq:fM-bracket}, the second one by 
the computation of $\sigma(\gamma)$ and $\beta_A^{\vee}(\gamma\otimes 1)$ 
made above, and the third one because $\tau(\gamma)\in U_1\wedge U_1$ and 
$[\B_1(A),\B_1(A)]=0$ in $\fM$.  Adding up, and recalling that 
$\nabla_A=d_A^{\vee}+\mu_A^{\vee}$, we obtain
\[
\rho\big(\nabla_A(\gamma)\big)=
q\Big(\nu_A^{\vee}(\gamma)\otimes 1,\, 
\big(d_A^{\vee}\otimes\id_S+\beta_A^{\vee}\big)(\gamma\otimes 1)\Big)=
q\big(\Delta_A(0,\gamma\otimes 1)\big)=0 .
\]
Since $\ker(\rho)$ is an ideal containing $\im (\nabla_A)$, the morphism 
$\rho$ factors through a morphism $\bar\rho\colon \h(A)\to \fM$.  As 
$\fM''=0$, we have $\bar\rho(\h'')=0$, while 
$\bar\rho(\h')\subseteq \fM'\subseteq \B_1(A)$; hence, $\B(\h(A))$ being by 
definition the quotient $\h'/\h''$, the morphism $\bar\rho$ induces a 
$\k$-linear map $\varsigma\colon \B(\h(A))\to \B_1(A)$.  This map is 
$S$-linear: if $e\in E_1$ and $\xi\in \h'$, then $\bar\rho(e)=e$, and so 
$\varsigma\big(\bar{e}\cdot \bar{\xi}\big)=\bar\rho\big([e,\xi]\big)=
\big[e,\varsigma(\bar\xi)\big]=e\cdot \varsigma(\bar\xi)$.

Finally, we claim that $\varsigma\circ (\eta_A\:\:\psi_A)=q$.  Both maps are 
$S$-linear, so it suffices to compare them on the generators of the free 
$S$-module $(E_2\oplus U_1)\otimes_{\k}S$; and indeed 
$\varsigma\big(\eta_A(e\wedge f)\big)=\bar\rho\big([e,f]\big)=
q\big(e\wedge f\otimes 1,0\big)$ for $e,f\in E_1$, while 
$\varsigma\big(\psi_A(u)\big)=\bar\rho(u)=q\big(0,u\otimes 1\big)$ for 
$u\in U_1$.  Since $(\eta_A\:\:\psi_A)=\phi\circ q$ and $q$ is surjective, 
it follows that $\varsigma\circ \phi=\id$.  Therefore $\phi$ is injective, 
and thus an isomorphism of $S$-modules, with inverse $\varsigma$.
\end{proof}

\begin{remark}
\label{rem:fM-metabelian}
The Lie algebra $\fM$ is not a device internal to the proof: it is
$\h(A)/\h''(A)$.  Indeed $\bar\rho$ is onto, being the identity on $E_1$ and
carrying $\h'$ onto $\B_1(A)$, and $\ker(\bar\rho)=\h''(A)$ because $\varsigma$
is injective; so
\begin{equation}
\label{eq:fM-is-max-metabelian}
\h(A)/\h''(A) \cong \fM=\fE\bigl(E_1,\B_1(A),c\bigr) .
\end{equation}
In the language of Lemma~\ref{lem:metabelian-extension}, the theorem says that
the presentation of $\B_1(A)$ from Theorem~\ref{thm:Bpres}, together with the
cocycle $c$, is a presentation of the maximal metabelian quotient of $\h(A)$;
and $[c]\in H^2(E_1;\B_1(A))$ is the class of the extension
\eqref{eq:liess} for $\g=\h(A)$.  By
Remark~\ref{rem:metabelian-extension-caveats}, that class is canonical while
$c$ is not, which is why $\fM$ is functorial only up to the correction
$\theta_\varphi$ of \eqref{eq:theta-phi}.  Compare
Proposition~\ref{prop:holo-free}, whose proof likewise identifies a maximal
metabelian quotient---there $L/L''$, inside the wreath product
\eqref{eq:wreath}---by Remark~\ref{rem:wreath-magnus}.
\end{remark}

%%%%%%%%%%%%%%%%%%
\subsection{Naturality of the comparison map}
\label{subsec:compare-naturality}
%%%%%%%%%%%%%%%%%%

The identification of the first Koszul module with the infinitesimal
Alexander invariant of the holonomy Lie algebra is compatible with morphisms
of $\cdgas$, in the sense of base change.  As the example at the end of this
subsection shows, however, the two sides do not behave alike under
morphisms: $\h(\varphi)$ sees only $\varphi^1$ and $\varphi^2$, whereas
$\B_1(\varphi)$ sees all of $\varphi$.

\begin{proposition}
\label{prop:naturality-B-holo}
Let $A$ be a connected $\k$-$\cdga$ with $\dim_\k A^1<\infty$. 
The isomorphism
\[
\B_1(A)\longisom \B(\h(A))
\]
constructed in Theorem~\ref{thm:B-holo} is natural in the sense of 
base change: for any morphism of
$\cdgas$ as above, $\varphi\colon A\to \oA$, the diagram 
of $S$-modules 
\[
\begin{tikzcd}[column sep=30pt, row sep=24pt]
\B_1(\oA)\otimes_{\oS} S \ar[r,"\simeq"] \ar[d,"\B_1(\varphi)"] &
\B(\h(\oA))\otimes_{\oS} S \ar[d,"\B(\h(\varphi))"] \\
\B_1(A) \ar[r,"\simeq"] & \B(\h(A))
\end{tikzcd}
\]
commutes, where $S=\Sym(H_1(A))$, $\oS=\Sym(H_1(\oA))$, and the vertical
arrows are the canonical base-change maps associated to $\varphi$ and
$\h(\varphi)$ respectively. 
\end{proposition}

\begin{proof}
Let $\varphi\colon A\to \oA$ be a morphism of connected, $1$-finite $\k$-$\cdgas$. 
Write $E_1=H^1(A)^\vee$, $\oE_1=H^1(\oA)^\vee$, $S=\Sym(E_1)$, and $\oS=\Sym(\oE_1)$,
and let $\varphi_1^{E}$, $\varphi_1^{U}$, and $\theta_{\varphi}$ be as in 
\eqref{eq:phi-E-U} and \eqref{eq:theta-phi}. 
From the proof of Theorem~\ref{thm:B-holo}, we have surjective $S$-linear maps
\[
(\eta_A \:\: \psi_A)\colon (E_2 \oplus U_1)\otimes_{\k} S \longrightarrow \B(\h(A)), 
\qquad (\eta_{\oA} \:\: \psi_{\oA})\colon (\oE_2 \oplus \oU_1)\otimes_{\k} \oS 
\longrightarrow \B(\h(\oA)),
\]
whose kernels are exactly the images of the corresponding $\Delta_A$ and 
$\Delta_{\oA}$ maps from \eqref{eq:BA-pres}.  In view of the commuting 
diagram \eqref{eq:bphi}, it suffices to show that the square
\begin{equation}
\label{eq:nat-square-holo}
\begin{tikzcd}[column sep=44pt, row sep=24pt]
(\oE_2\oplus \oU_1)\otimes_{\k} S 
\ar[r, "{(\eta_{\oA}\:\:\psi_{\oA})}"] 
\ar[d, "\sbm{\bwedge^2\varphi_1^{E} \amp \theta_{\varphi}\\ 0 \amp \varphi_1}" '] 
& \B(\h(\oA))\otimes_{\oS} S \ar[d, "\B(\h(\varphi))"] \\
(E_2\oplus U_1)\otimes_{\k} S \ar[r, "{(\eta_A\:\:\psi_A)}"] & \B(\h(A))
\end{tikzcd}
\end{equation}
commutes.  Both composites are $S$-linear, so it is enough to compare them 
on the generators of the free $S$-module $(\oE_2\oplus \oU_1)\otimes_{\k}S$. 

Let $\bar e,\bar f\in \oE_1$.  Going right and then down, we obtain the class 
of $\h(\varphi)\big([\bar e,\bar f]\big)=
\big[\varphi_1(\bar e),\varphi_1(\bar f)\big]$ modulo $\h''(A)$.  Expanding 
$\varphi_1|_{\oE_1}=\varphi_1^{E}+\varphi_1^{U}$ and using that 
$\big[\varphi_1^{U}(\bar e),\varphi_1^{U}(\bar f)\big]\in [U_1,U_1]\subseteq 
[\h'(A),\h'(A)]=\h''(A)$, this class equals 
\[
\big[\varphi_1^{E}(\bar e),\varphi_1^{E}(\bar f)\big]+
\big[\varphi_1^{E}(\bar e),\varphi_1^{U}(\bar f)\big]-
\big[\varphi_1^{E}(\bar f),\varphi_1^{U}(\bar e)\big] .
\]
Going down and then right, we obtain 
$\eta_A\big(\varphi_1^{E}(\bar e)\wedge \varphi_1^{E}(\bar f)\big)+
\psi_A\big(\theta_{\varphi}(\bar e\wedge \bar f)\big)$, which by 
\eqref{eq:theta-phi} equals 
\[
\big[\varphi_1^{E}(\bar e),\varphi_1^{E}(\bar f)\big]+
\varphi_1^{E}(\bar e)\cdot \overline{\varphi_1^{U}(\bar f)}-
\varphi_1^{E}(\bar f)\cdot \overline{\varphi_1^{U}(\bar e)} .
\]
The two expressions agree, since the $S$-action on $\B(\h(A))$ is induced by 
the adjoint action of $E_1=\h(A)/\h'(A)$.  Likewise, for $\bar u\in \oU_1$, 
both composites send $\bar u\otimes 1$ to the class of $\varphi_1(\bar u)$, 
by \eqref{eq:phi-U-canonical}.  Hence \eqref{eq:nat-square-holo} commutes. 

Passing to cokernels, it follows that the induced $S$-linear map
$\B_1(\oA)\otimes_{\oS} S \to \B_1(A)$ coincides with the map
$\B(\h(\oA))\otimes_{\oS} S \to \B(\h(A))$, 
under the identifications $\B_1(-)\cong \B(\h(-))$ from Theorem~\ref{thm:B-holo}. 
Hence the diagram
\[
\begin{tikzcd}[column sep=30pt, row sep=24pt]
\B_1(\oA)\otimes_{\oS} S \ar[r,"\simeq"] \ar[d,"\B_1(\varphi)"] &
\B(\h(\oA))\otimes_{\oS} S \ar[d,"\B(\h(\varphi))"] \\
\B_1(A) \ar[r,"\simeq"] & \B(\h(A))
\end{tikzcd}
\]
commutes, proving the base-change naturality.
\end{proof}

\begin{remark}
\label{rem:naturality-B-holo-subtlety}
Proposition~\ref{prop:naturality-B-holo} admits two complementary 
readings, which we record here since they inform the way the 
isomorphism $\B_1(A)\cong \B(\h(A))$ is put to use later on.

First, it is a statement about dependence on data.  For a morphism 
$\varphi\colon A\to\oA$, the induced map 
$\h(\varphi)\colon\h(\oA)\to\h(A)$ is determined entirely by 
$\varphi^1\colon A^1\to\oA^1$ and $\varphi^2\colon A^2\to\oA^2$: 
the dual maps $(\varphi^1)^\vee$ and $(\varphi^2)^\vee$ prescribe 
its effect on the generators and the relations of the holonomy Lie 
algebras, while the components $\varphi^k$ with $k\ge 3$ play no role. 
By contrast, the morphism $\B_1(\varphi)\colon\B_1(\oA)\otimes_{\oS}S
\to\B_1(A)$ is defined via the full chain map 
$K_\bullet(\varphi)\colon K_\bullet(\oA)\to K_\bullet(A)$
of Proposition~\ref{prop:koszul-functorial}, which involves all 
components of $\varphi$.  The commutativity of the naturality diagram 
shows this additional input to be illusory: under the identifications 
of Theorem~\ref{thm:B-holo}, we have 
$\B_1(\varphi)=\B(\h(\varphi))$, and so the induced map on Koszul 
modules is likewise determined by $\varphi^1$ and $\varphi^2$.  In 
particular, on morphisms, the Koszul module functor detects nothing 
beyond the holonomy Lie algebra functor.

Second, a warning: an isomorphism in degree $1$ does not force an 
isomorphism of holonomy Lie algebras, since the relations are governed 
by degree~$2$.  For instance, let $A=\bwedge(a,b)$ with zero 
differential, let $B=A/(a\wedge b)$, and let $q\colon A\to B$ be the 
canonical projection---a $\cdga$ model for the inclusion 
$S^1\vee S^1\hookrightarrow T^2$.  Then $q^1$ is an isomorphism, yet 
$\h(q)$ is the projection of the free Lie algebra $\h(B)=\Lie(x,y)$ 
onto the abelian Lie algebra $\h(A)=\Lie(x,y)/\langle [x,y]\rangle$; 
correspondingly, $\B_1(q)$ is the zero map from 
$\B_1(B)\cong \B(\Lie(x,y))\ne 0$ to $\B_1(A)=0$.  Thus, hypotheses 
on $H^1$ alone do not control $\h(\varphi)$ or $\B_1(\varphi)$; one 
also needs control of the relation ideals, that is, of $\varphi^2$. 
Note, however, that if $\alpha$ is an \emph{endo}morphism of a group 
$G$ inducing the identity on $H^1(G;\k)$, then the induced 
endomorphism of $\h(H^*(G;\k))$ \emph{is}\/ the identity: the 
generators are fixed, and the relation ideal---being the same ideal 
on both sides---is preserved.

Consequently, the obstructions to $1$-formality of group 
homomorphisms developed in Section~\ref{sect:functorial-models} do 
not arise from a discrepancy between $\B_1(\varphi)$ and 
$\h(\varphi)$---these agree, as just explained---but from the 
discrepancy between model-level and cohomology-level invariants. 
A morphism may act trivially on cohomology, and hence on 
$\h(H^*(A))$ and on associated graded invariants, while acting 
nontrivially on a Sullivan model $A$, and hence on $\h(A)$ and on 
$\B_1(A)\cong\B(\h(A))$.  It is this finer, filtered information 
that obstructs the strict $1$-formality of maps 
(Proposition~\ref{prop:koszul-obstruction-groups}).
\end{remark}

\begin{example}
\label{ex:nonhomog-tre}
Let $A=\CE(\sol_2)$ be the CE-algebra of
Example~\ref{ex:sol2-koszul}, so that $\B_1(A)\cong\k[x]/(x-1)$ is
supported at $1\in H^1(A)=\k a$, where $S=\k[x]$ with $x=a^{\vee}$.  

For the holonomy side: $\h(A)=\sol_2$ has $\h'=\k y$ and 
$\h''=[\h',\h']=0$, so $\B(\h(A))=\h'/\h''=\k y$ with 
$S=\k[x]$-action $x\cdot y=[x,y]=y$, giving 
$\B(\h(A))\cong\k[x]/(x-1)$.
This confirms the isomorphism $\B_1(A)\cong\B(\h(A))$ 
of Theorem~\ref{thm:B-holo} in this case.

Now consider the inclusion $\iota\colon B\inj A$, where 
$B=\bigwedge(a)$ with $d(a)=0$, which is a quasi-isomorphism 
(Example~\ref{ex:non-functorial-qiso}).
The dual map $(\iota^1)^\vee\colon A_1\surj B_1$ sends $x\mapsto x$ 
and $y\mapsto 0$; since the defining relation of $\h(A)=\sol_2$ maps 
to $0$, the induced map $\h(\iota)\colon \sol_2 \surj \Lie(x)=\h(B)$ 
is the projection killing $y$.

Meanwhile, $\B_1(B)=0$ by \eqref{eq:sol2-modules-coh}, while
$\B_1(A)=\k[x]/(x-1)$ is supported at $1$; correspondingly,
$\B(\h(B))=\B(\Lie(x))=0$, since $\Lie(x)$ is abelian, whereas
$\B(\h(A))\cong\k[x]/(x-1)$.  Thus both $\B_1(\iota)$ and
$\B(\h(\iota))$ are zero maps, and neither is an isomorphism, even
though $\iota$ is a quasi-isomorphism.

In accordance with Proposition~\ref{prop:naturality-B-holo}, the 
two induced maps correspond to one another under the isomorphisms 
of Theorem~\ref{thm:B-holo}.  The example also illustrates the 
second point of Remark~\ref{rem:naturality-B-holo-subtlety}: 
although $\iota^*$ is an isomorphism on $H^*$ (in particular, on 
$H^1$), neither $\h(\iota)$ nor $\B_1(\iota)$ is an 
isomorphism---information recorded by the models, but invisible to 
cohomology alone.
\end{example}

%%%%%%%%%%%%%%%%%%
\subsection{Koszul modules and holonomy Chen ranks}
\label{subsec:chenalg}
%%%%%%%%%%%%%%%%%%

Once again, let $(A,d)$ be a connected $\k$-$\cdga$ 
with $A^1$ finite-dimensional, and let $\h=\h(A)$ be its holonomy 
Lie algebra.  We let $\gr(\h)$ denote the graded Lie algebra associated 
to the lower central series filtration on $\h$. 

We define the {\em holonomy Chen ranks}\/ of $A$ by
\begin{equation}
\label{eq:chrank}
\theta_{k}(A)\coloneqq \dim_{\k} \gr_{k} \big(\h(A)/\h''(A)\big).  
\end{equation}

Although $\h=\h(A)$ is not naturally graded, its lower central series 
induces a canonical filtration on the vector space $\h'/\h''$:
\begin{equation}
\label{eq:chen-filt}
F_n \coloneqq (\gamma_n(\h) + \h'') / \h'' \subseteq \h'/\h'', 
\quad n \ge 2.
\end{equation}
This filtration gives rise to the associated graded vector space
\[
\gr(\h'/\h'') \coloneqq \boplus_{n \ge 2} F_n/F_{n+1},
\]
whose graded pieces satisfy
\[
\dim_\k F_n/F_{n+1} = \theta_n(A), \quad \text{for } n \ge 2.
\]
Thus, for $n \ge 2$, the holonomy Chen ranks measure the successive layers 
of the lower central series of $\h(A)$ modulo $\h(A)''$, while 
$\theta_1(A)=\dim_\k H^1(A)$. The grading on $\gr(\h'/\h'')$ 
reflects the lower central series of $\h(A)$ and is concentrated 
in degrees $\ge 2$. As will be seen below, when $\gr(\h'/\h'')$ 
is identified with $\gr(\B_1(A))$, its lower central series grading 
corresponds to the internal $S$-grading up to a shift by $2$.

\begin{proposition}
\label{prop:lcs-chen}
Let $(A,d)$ be a connected $\k$-$\cdga$ with holonomy Lie algebra $\h(A)$.
Then the holonomy LCS ranks $\phi_n(A)=\dim_\k \gr_n(\h(A))$ and the holonomy
Chen ranks $\theta_n(A)$ satisfy
\[
\phi_n(A)\ge \theta_n(A)\quad\text{for all $n\ge 1$},
\]
with equality for $n\le 3$.
\end{proposition}

\begin{proof}
The natural projection $\h \to \h/\h''$ induces surjections
$\gr_n(\h)\twoheadrightarrow \gr_n(\h/\h'')$ for all $n$, giving the inequality.  
Moreover, since $\h''=[\h',\h'] \subseteq \gamma_4(\h)$ for any Lie algebra $\h$, 
the map is an isomorphism in degrees $n \le 3$, proving the equality claim.
\end{proof}

The $S=\Sym(A_1)$-module structure on $\h'/\h''$ is compatible with the 
filtration \eqref{eq:chen-filt}, and thus induces a graded $S$-module structure on 
$\gr(\h'/\h'')$. By Theorem~\ref{thm:B-holo}, this module is identified, 
up to a shift in grading, with $\gr(\B_1(A))$. The next proposition 
expresses the holonomy Chen ranks in terms of the Hilbert series of 
$\gr(\B_1(A))$.

\begin{proposition} 
\label{prop:holo-massey}
The generating series for the holonomy Chen ranks of $A$ (shifted by $2$) 
coincides with the Hilbert series of the $\m$-adic associated graded 
$S$-module $\gr^\m(\B_1(A))$:
\begin{equation}
\label{eq:chen-hilb-cdga}
\sum\limits_{n\ge 0}{\theta}_{n+2}(A)\, t^n = \Hilb (\gr^\m(\B_1(A)),t).
\end{equation}
\end{proposition}

\begin{proof}
By definition of the lower central series filtration $\gamma_\bullet$ on 
$\h'/\h''$, we have
\[
\theta_n(A)=\dim_\k \gr_n(\h'/\h'') \quad \text{for } n\ge 2.
\]
As noted above, the associated graded module $\gr(\h'/\h'')$ identifies, 
up to a shift in grading, with $\gr(\B(\h(A)))$, more precisely,
\[
\gr_n(\h'/\h'') \cong \gr^\m_{n-2}(\B(\h(A))).
\]
Hence,
\[
\theta_n(A)=\dim_\k \gr^\m_{n-2}(\B(\h(A))) \quad \text{for } n\ge 2.
\]
The claim now follows from Theorem~\ref{thm:B-holo}, which identifies 
$\gr(\B(\h(A)))$ with $\gr(\B_1(A))$.
\end{proof}

\begin{remark}
\label{rem:chen-madic-convention}
The $\m$-adic normalization in Proposition~\ref{prop:holo-massey} is
essential, and it is \emph{not}\/ the one computed by the Koszul spectral
sequence: by Lemma~\ref{lem:filtration-shift-B1}, the filtration induced
on $\B_1(A)$ by that spectral sequence is the $\m$-adic filtration shifted
by one, so
$\Hilb(\gr^F\B_1(A),t)=t\cdot\Hilb(\gr^\m\B_1(A),t)$.  Concretely, for
$A=H^*(F_2)$ the free group calculation gives
$\theta_n=n-1$ \textup{(}Chen's theorem\textup{)}, matching
$\dim_\k\gr^\m_{\,n-2}\B_1(A)=n-1$, whereas the internal grading of
$\B_1(A)\cong S(-1)$ would produce $n-2$.  Whenever a Chen-rank formula
is extracted from the spectral sequence---as in
Theorem~\ref{thm:Hilb-ineq-gr}---the shift must therefore be taken into
account; it cancels in any comparison between $A$ and $H^*(A)$, both
sides carrying the same shift.
\end{remark}

When $d_A=0$, so that $A$ is a graded algebra, viewing $A$ as a module over
the exterior algebra $\cE=\bigwedge A^1$ and combining
Theorem~\ref{thm:gr-Bi-tor} with Proposition~\ref{prop:holo-massey} gives a
concrete homological description of the holonomy Chen ranks: they count the
linear strand of the minimal free $\cE$-resolution of $A$.

\begin{corollary}
\label{cor:chen-tor}
Let $A$ be a connected graded $\k$-algebra with $\dim_\k A^1<\infty$, regarded
as a $\cdga$ with zero differential.  Then, for all $n\ge 2$, 
\[
\theta_n(A)=\dim_{\k}  \Tor^{\cE} _{n-1}(A,\k)_{n}. 
\]
\end{corollary}

For a $\cdga$ with nonzero differential the same formula, applied to $H^{*}(A)$
with $\cE=\bigwedge H^1(A)$, computes $\theta_n(H^{*}(A))$.  This bounds
$\theta_n(A)$ from above, by Theorem~\ref{thm:Hilb-ineq-gr}, with equality when
the Koszul spectral sequence degenerates.

Taken together, the preceding results give a clear picture of the holonomy Chen ranks. 
The first rank is always $\theta_1(A) = \dim H^1(A)$, while for $n \ge 2$, 
the ranks are controlled by the linear strand of the minimal free $\cE$-resolution 
of $A$. In particular, by Corollaries~\ref{cor:gr-Bi-vanish} 
and \ref{cor:chen-eventual}, once a Chen rank $\theta_n(A)$ vanishes for some 
$n\ge 2$, all higher ranks vanish as well. This reflects the general principle 
that the linear strand of a minimal free resolution cannot revive after terminating, 
so in that case the sequence of holonomy Chen ranks eventually stabilizes at zero.

\begin{example}
\label{ex:chen-heis}
Let $A=\bwedge(a,b,c)$ with $d a=d b=0$ and $d c = -a\wedge b$.
Then $H^1(A)=\k a\oplus\k b$, so $A_1=\k x\oplus\k y$ and 
$S=\k[x,y]$. The Koszul module is $\B_1(A)\cong S/(x,y)\cong \k$, 
concentrated in internal degree $0$, so $\Hilb(\gr(\B_1(A)),t)=1$.

The holonomy Lie algebra $\h=\h(A)$ has generators $x,y,z$ with 
relations $z-[x,y]=[x,z]=[y,z]=0$. Eliminating $z$ shows $\h$ is 
generated by $x,y$, with $\h'$ generated by $[x,y]$ and $\h''=0$. 
Hence $\h/\h'\cong\k^2$ and $\h'/\h''\cong\k$, giving 
holonomy Chen ranks $\theta_1(A)=2$, $\theta_2(A)=1$, and $\theta_n(A)=0$ for $n\ge3$.

This agrees with Corollary~\ref{cor:chen-tor}, since viewing $A$ as a 
graded $\cE=\bwedge(a,b)$-module, the linear strand of the minimal 
free $\cE$-resolution of $A$ contributes a single copy of $\k$ in bidegree $(1,2)$.
\end{example}

\begin{corollary}
\label{cor:chen-hirsch}
Let $A=H\otimes_{e}\bigwedge(c)$ be as in Proposition~\ref{prop:holo-hirsch},
and assume $e=0$, or else $e\ne0$ and $e\in\im(\mu)$.
Then the holonomy Chen ranks of $A$ satisfy
\[
\theta_1(A)=\phi_1(A),\quad
\theta_2(A)=\phi_2(A),\quad
\text{and}\quad
\theta_n(A)=\theta_n(H)\quad\text{for all }n\ge 3,
\]
where $\phi_n(A)$ is given by Proposition~\ref{prop:holo-hirsch}.
\end{corollary}

\begin{proof}
By Proposition~\ref{prop:holo-hirsch}, the Lie algebra $\h(A)$ is graded, so
$\theta_n(A)=\dim_\k\bigl(\h(A)/\h''(A)\bigr)_n$, and the equalities
$\theta_n=\phi_n$ for $n\le3$ are those of
Proposition~\ref{prop:lcs-chen}; in particular, the case $n=3$ of the last
assertion already follows from Proposition~\ref{prop:holo-hirsch}, since
$\h''$ starts in degree~$4$, so that the content lies in degrees $n\ge4$.  
In the case $e=0$ the assertion is immediate from 
$\h(A)\cong\h(H)\oplus\Lie(c^{\vee})$.  In the case
$e\in\im(\mu)$, apply $(-)/(-)''$ to the central extension
\eqref{eq:holo-hirsch-central}: since $\k z$ is spanned by an element of
degree~$2$, it meets $\h''(A)$ trivially, and the projection carries
$\h''(A)$ onto $\h''(H)$.  Hence
$\bigl(\h(A)/\h''(A)\bigr)_n\cong\bigl(\h(H)/\h''(H)\bigr)_n$ for
$n\ge3$, which gives the stated formulas.
\end{proof}

\begin{remark}
\label{rem:holo-hirsch-caveats}
The hypothesis $e\in\im(\mu)$ in
Proposition~\ref{prop:holo-hirsch}\eqref{eh1} cannot be weakened to
$e\ne0$.  Take $H=H^{*}(S^1\times S^2;\k)$, so that
$H^1=\langle x\rangle$, $H^2=\langle h\rangle$ and $\mu=0$ because
$\bwedge^2H^1=0$, and take $e=h\ne0$.  Then $A^1=\langle x,c\rangle$ with
$dc=h$, and by \eqref{eq:holo-hirsch-pres-dB0} the relations of $\h(A)$ are
$c^{\vee}$ together with $[x^{\vee},c^{\vee}]$; so $\h(A)=\Lie(x^{\vee})\cong\k$
and $\phi_2(A)=0=\phi_2(H)$, not $\phi_2(H)+1$.  This is a clean instance of
case \eqref{eh2}, rather than a borderline one: since $\mu=0$ while $e\ne0$,
the condition $e\notin\im(\mu)$ holds automatically, and $\mu^{\vee}=0$ forces
$z=\mu^{\vee}(\gamma_0)=0$ for \emph{every}\/ choice of $\gamma_0$---so the
step of the proof that selects $\gamma_0$ annihilating $\im(\mu)$ is vacuous
here.  Here $A$ is a model for the
circle bundle over $S^1\times S^2$ with Euler class $h$, namely
$S^1\times S^3$, whose fundamental group is $\Z$; the abelian answer is the
expected one.  Compare Remark~\ref{rem:n1-enot0}, where the same dichotomy
$e\in\im(\mu)$ versus $e\notin\im(\mu)$ separates the two
obstructions to $1$-formality.
\end{remark}

%%%%%%%%%%%%%%%%%%%%%%%%%%%%%%
\subsection{Chen ranks of $2$-step nilpotent Lie algebras}
\label{subsec:chen-2step}
%%%%%%%%%%%%%%%%%%%%%%%%%%%%%%

We now apply the preceding machinery to $2$-step nilpotent Lie algebras.
As we will see, the Koszul module of the Chevalley--Eilenberg $\cdga$ is as
small as it can be---it has finite length, and the holonomy Chen ranks of
$\CE(\h)$ vanish beyond degree $2$.  All the growth in this setting occurs
on the cohomological side, for the holonomy Lie algebra of $H^*(\CE(\h))$,
and the resulting gap between the two is precisely a failure of
$1$-formality.

Let $\h$ be a finite-dimensional $2$-step nilpotent Lie algebra over $\k$,
generated in degree~$1$. Thus
\[
\h=\h_1\oplus \h_2,
\qquad
\h_2=[\h_1,\h_1]\subset Z(\h),
\]
and all iterated brackets of length $\ge 3$ vanish.
Set
\[
V=\h_1^\vee, \qquad W=\h_2^\vee, \qquad S=\Sym(\h_1).
\]
The Chevalley--Eilenberg $\cdga$ of\/ $\h$ is
\[
A=\CE(\h)=\bigl(\bwedge(V\oplus W),d\bigr),
\]
with $d|_V=0$ and $d|_W\colon W\hookrightarrow \bwedge^2 V$
dual to the Lie bracket; we identify $W$ with its image in $\bwedge^2V$
throughout.
By Proposition~\ref{prop:CE-positive-weights}, $A$ admits a positive-weight
decomposition: assign weight $1$ to $V = \h_1^\vee$ and weight $2$ to
$W = \h_2^\vee$, so that $d\colon W \to \bwedge^2 V$ preserves total weight.
Consequently, $\B_1(A)$ is a \emph{graded}\/ finitely generated $S$-module,
not merely filtered, and its Hilbert series is well-defined.

\begin{proposition}
\label{prop:B1-2step-presentation}
Let\/ $\h$ be a $2$-step nilpotent Lie algebra as above, and let
$A=\CE(\h)$.  Then $\h(A)\cong\h$, and there is an isomorphism of
$S$-modules
\begin{equation}
\label{eq:B1-2step}
\B_1(\CE(\h))\cong \B(\h)=\h_2 ,
\end{equation}
where $S$ acts on\/ $\h_2$ through the augmentation $S\to\k$.  In
particular, $\B_1(A)$ is concentrated in weight $2$ and has finite length
$r=\dim_\k\h_2$, with
\[
\Ann_S\B_1(A)=\m
\qquad\text{and}\qquad
\Supp_S\B_1(A)=\{0\} .
\]
\end{proposition}

\begin{proof}
The holonomy Lie algebra $\h(A)$ is generated by $A_1=\h_1\oplus\h_2$,
modulo the relations dual to the three summands
$A^2=\bwedge^2V\oplus(V\otimes_{\k} W)\oplus\bwedge^2 W$ under the map
$d^\vee+\mu_A^\vee$ of Section~\ref{subsec:cochain-holo-short}.  The
summand $\bwedge^2V$ contributes the relations $[v,v']=[v,v']_{\h}$ for
$v,v'\in\h_1$, since $d|_W$ is dual to the bracket of\/ $\h$; the summand
$V\otimes_{\k} W$ contributes $[v,u]=0$ for $v\in\h_1$ and $u\in\h_2$; and
$\bwedge^2W$ contributes $[u,u']=0$.  These are precisely the relations
defining\/ $\h$, so $\h(A)\cong\h$.

By Theorem~\ref{thm:B-holo}, $\B_1(A)\cong\B(\h(A))=\h'/\h''$.  Since\/
$\h$ is $2$-step, $\h'=\h_2$ and $\h''=[\h_2,\h_2]=0$, whence
$\B(\h)=\h_2$.  The $S$-module structure on $\B(\h)$ is induced by the
adjoint action, which maps $\h_1\otimes_{\k}\h_2$ into $\h_3=0$; hence $\m$
annihilates $\B_1(A)$, and the remaining assertions follow.
\end{proof}

Since $\B_1(A)$ has finite length, the Chen ranks of $A$ vanish in degrees
$\ge3$; this is the input to part~\eqref{tF2} of
Theorem~\ref{thm:2step-intro} from the Introduction.

\begin{remark}
\label{rem:B1-2step-no-presentation}
Proposition~\ref{prop:B1-2step-presentation} shows that the weight-$2$
generators $W$ of $\B_1(A)$ carry no $S$-module structure at all: every
relation of the form $w\otimes x$ with $x\in S_1$ already holds.  At the
chain level, these relations come from the middle summand $V\otimes_{\k} W$ of
$A^2$, whose dual contributes $\partial_2\bigl((v\wedge w)^\vee\bigr)
=x_v\cdot w^\vee$; the summand $\bwedge^2V$ contributes only the
identifications of the Koszul syzygies of $\Sym(\h_1)$ with the constants
$\sum_{i} c^i_{jl}w_i^\vee$, and imposes no further relation among the
$w_i^\vee$.  Retaining only the contributions of $\bwedge^2V$ would
produce a strictly smaller relation submodule, and hence a module that is
too large---for instance, of positive Krull dimension whenever some
$dw_i$ is degenerate.
\end{remark}

\begin{corollary}
\label{cor:chen-2step}
Let\/ $\h$ be a finite-dimensional $2$-step nilpotent Lie algebra with
$m=\dim\h_1$ and $r=\dim\h_2$, and let $A=\CE(\h)$.  Then
\[
\Hilb\bigl(\B_1(A),t\bigr)=r,
\qquad
\theta_1(A)=m,\quad \theta_2(A)=r,\quad \theta_n(A)=0 \ \ (n\ge3).
\]
Thus the holonomy Chen ranks of\/ $\CE(\h)$ exhibit no growth whatsoever.
\end{corollary}

\begin{proof}
The Hilbert series is immediate from \eqref{eq:B1-2step}, since $\B_1(A)$
is concentrated in a single weight.  By
Proposition~\ref{prop:holo-massey}, $\sum_{n\ge0}\theta_{n+2}(A)t^n=
\Hilb(\gr^\m\B_1(A),t)=r$, giving $\theta_2(A)=r$ and $\theta_n(A)=0$ for
$n\ge3$; and $\theta_1(A)=\dim_\k H^1(A)=m$.  Alternatively, and
independently of the Koszul machinery, $\h(A)=\h$ is $2$-step, so
$\h''=0$ and $\theta_n(A)=\dim_\k\gr_n(\h)=\dim_\k\h_n$.
\end{proof}

Consistently with Theorem~\ref{thm:CE-resonance-nilpotent}, the support
being $\{0\}$ forces $\RR_1(A)\subseteq\{0\}$, so no information about\/
$\h$ beyond the number $r$ is visible in $\B_1(\CE(\h))$.  In particular,
the free case is no exception.

\begin{example}
\label{ex:free-2step-B1}
Let $\h = \ff_m/\gamma_3\ff_m$ be the free $2$-step nilpotent Lie algebra
on $m \ge 2$ generators, so that $\h_1 = \k^m$ and
$\h_2 = \bwedge^2\k^m$ has dimension $\binom{m}{2}$, and let
$A=\CE(\h)$, with $d a_i = 0$ and $d c_{ij} = -a_i \wedge a_j$ for $i<j$.
By Proposition~\ref{prop:B1-2step-presentation},
\[
\B_1(A)\cong\bwedge^2\k^m,
\]
of finite length $\binom{m}{2}$, with $\theta_2(A)=\binom{m}{2}$ and
$\theta_n(A)=0$ for $n\ge3$.  For $m=2$ this recovers
$\B_1(\CE(\h(1)))\cong\k$ of Section~\ref{subsec:heisenberg-1}.
\end{example}

%%%%%%%%%%%%%%%%%%%%%%%%
\subsection{The cohomological side, and non-\texorpdfstring{$1$}{1}-formality}
\label{subsec:coho-side-2step}
%%%%%%%%%%%%%%%%%%%%%%%%

The interesting invariants of a $2$-step nilpotent $\h$ are those of the
cohomology algebra $H^*(A)$, whose holonomy Lie algebra is in general a
proper quotient of neither $\h$ nor a free Lie algebra, but is instead
determined by the quadratic data $W\subseteq\bwedge^2V$.  Write
\begin{equation}
\label{eq:bar-theta-2step}
\bar\theta_n\coloneqq\theta_n\bigl(H^*(A)\bigr)
\end{equation}
for the holonomy Chen ranks of $H^*(A)$; these are the quantities denoted
$\bar\theta_n(G)$ in Theorem~\ref{thm:chen-ineq-model} when $A$ is a
$1$-model for a group $G$.

The next result is the key step in the proof of
Theorem~\ref{thm:2step-intro} from the Introduction.

\begin{theorem}
\label{thm:chen-growth-2step}
Let\/ $\h$ be a finite-dimensional $2$-step nilpotent Lie algebra with
$m=\dim\h_1\ge2$, let $W=\h_2^\vee\subseteq\bwedge^2V$, and set
$A=\CE(\h)$.
\begin{enumerate}[itemsep=2pt]
\item \label{chn1}
The cup product $\bwedge^2H^1(A)\to H^2(A)$ is the projection
$\bwedge^2V\surj\bwedge^2V/W$, and
\[
\h\bigl(H^*(A)\bigr)=\Lie(\h_1)/\bigl\langle W^{\perp}\bigr\rangle,
\]
where $W^\perp\subseteq\bwedge^2\h_1$ is the annihilator of $W$ under the
pairing between $\bwedge^2\h_1$ and $\bwedge^2V$.  Moreover
$\sum_{n\ge0}\bar\theta_{n+2}\,t^n=\Hilb\bigl(\B_1(H^*(A)),t\bigr)$.
\item \label{chn2}
The degree-$1$ jump locus of the cohomology algebra is the union of the
$2$-planes on which $W$ is supported:
\begin{equation}
\label{eq:res-2step}
\RR^1\bigl(H^*(A)\bigr)=
\bigcup\,\bigl\{P\subseteq V \;\big|\;
\dim_\k P=2,\ \bwedge^2P\subseteq W\bigr\},
\end{equation}
a union of $2$-planes indexed by the decomposable elements of $W$.
Consequently $\bar\theta_n$ grows polynomially of degree
$\dim\RR^1(H^*(A))-1$, and $\bar\theta_n=0$ for $n\gg0$ if and only if
$W$ contains no nonzero decomposable element.
\item \label{chn3}
If $W=\bwedge^2V$---that is, if\/ $\h=\ff_m/\gamma_3\ff_m$---then
$W^\perp=0$, so $\h(H^*(A))$ is the free Lie algebra on\/ $\h_1$, the
locus \eqref{eq:res-2step} is all of $V$, and
\[
\bar\theta_n=(n-1)\binom{n+m-2}{n},
\]
of degree $m-1$ in $n$: the maximum possible.
\item \label{chn4}
If\/ $\bar\theta_n\ne0$ for some $n\ge3$, then $A$ is not $1$-formal.
\end{enumerate}
\end{theorem}

\begin{proof}
\eqref{chn1} Since $d|_V=0$ and $\im(d|_W)=W\subseteq\bwedge^2V$, we have
$H^1(A)=V$, while $H^2(A)=\ker(d|_{A^2})/\im(d|_{A^1})$ receives
$\bwedge^2V$ with kernel exactly $\bwedge^2V\cap\im(d)=W$.  The holonomy
Lie algebra of a graded algebra depends only on the dual of its cup
product (Section~\ref{subsec:cochain-holo-short}), whose image here is the
annihilator $W^\perp$ of $W$; this gives the presentation.  The
identification of the generating series is
Proposition~\ref{prop:holo-massey} applied to the $\cdga$ $(H^*(A),0)$,
whose Koszul module $\B_1(H^*(A))$ is graded, so that
$\gr\B_1(H^*(A))=\B_1(H^*(A))$.

\eqref{chn2} Let $0\ne a\in V$.  By Proposition~\ref{prop:res-descriptions},
$a\in\RR^1(H^*(A))$ if and only if the Aomoto complex
$H^0\to H^1\to H^2$ of $H^*(A)$ at $a$ fails to be exact in the middle,
that is, if and only if there is some $b\in V$ with $b\notin\k a$ and
$a\wedge b=0$ in $H^2(A)=\bwedge^2V/W$.  The latter says precisely that
$a\wedge b\in W$, that is, that the $2$-plane $P=\langle a,b\rangle$
satisfies $\bwedge^2P\subseteq W$; and $\bwedge^2P$ is spanned by the
decomposable element $a\wedge b$.  This proves \eqref{eq:res-2step}.
By Proposition~\ref{prop:r1a}, $\RR^1(H^*(A))$ and
$\Supp_S\B_1(H^*(A))$ agree away from the origin, so the two have the
same dimension unless both are contained in $\{0\}$; the growth statement
then follows from the fact that the Hilbert polynomial of a finitely
generated graded module has degree one less than the dimension of its
support, together with Corollary~\ref{cor:support-dim-ineq}.

\eqref{chn3} If $W=\bwedge^2V$ then the cup product
$\bwedge^2H^1\to H^2$ vanishes, so $\h(H^*(A))$ is free on\/ $\h_1$ and
$\B_1(H^*(A))$ is the module of Koszul syzygies of $(x_1,\dots,x_m)$.  Its
Hilbert function gives Chen's formula for a free group of rank $m$,
namely $\bar\theta_n=(n-1)\binom{n+m-2}{n}$; alternatively, every
$2$-plane satisfies $\bwedge^2P\subseteq W$, so \eqref{eq:res-2step} is
all of $V$ and the growth degree is $m-1$.  This is the largest value
attainable for a $2$-step nilpotent\/ $\h$ with $\dim\h_1=m$: by
\eqref{chn2} the growth degree equals $\dim\RR^1(H^*(A))-1$, and
$\RR^1(H^*(A))\subseteq V$ with $\dim_\k V=m$.

\eqref{chn4} By Corollary~\ref{cor:chen-2step}, $\theta_n(A)=0$ for all
$n\ge3$.  Suppose $A$ were $1$-formal.  Since\/ $\h$ is finite-dimensional,
so is $A=\CE(\h)$, and
Corollary~\ref{cor:ss-formality}\eqref{cor-qformal} applies with $q=1$,
giving $\gr^\m_k\B_1(A)\cong\bigl(\B_1(H^*(A))\bigr)_k$ for all $k\ge0$.
Comparing Hilbert series by means of Proposition~\ref{prop:holo-massey}---%
applied to $A$, and, as in \eqref{chn1}, to $(H^*(A),0)$---gives
$\theta_n(A)=\bar\theta_n$ for all $n\ge2$; hence $\bar\theta_n=0$ for all
$n\ge3$, contrary to hypothesis.
\end{proof}

Part~\eqref{chn4} passes to groups.  Suppose $A=\CE(\h)$ is a $1$-model for a
finitely generated group $G$---for instance, when $\k=\Q$, a lattice in the
nilpotent Lie group with Lie algebra\/ $\h\otimes_{\Q}\R$.  The Chen ranks of
$G$ are computed by any $1$-finite $1$-model, so $\theta_n(G)=\theta_n(A)=0$
for $n\ge3$, while $\bar\theta_n(G)=\bar\theta_n$; and were $G$ to be
$1$-formal, the two would agree in every degree, by
Theorem~\ref{thm:chen-ineq-model}\eqref{item:chen-eq-formal} below.  Hence $G$
fails to be $1$-formal as soon as $\bar\theta_n\ne0$ for a single $n\ge3$.

Parts~\eqref{chn1}--\eqref{chn3} of Theorem~\ref{thm:chen-growth-2step}
establish part~\eqref{tF1} of Theorem~\ref{thm:2step-intro} from the
Introduction, while part~\eqref{chn4}, together with the vanishing
$\theta_n(A)=0$ for $n\ge3$ supplied by
Proposition~\ref{prop:B1-2step-presentation} and the preceding paragraph,
establishes part~\eqref{tF2}.

\begin{remark}
\label{rem:chen-2step-sharp}
The criterion in part~\eqref{chn4} is strictly finer than the growth
dichotomy of part~\eqref{chn2}: a single degree suffices.  It may happen
that $\RR^1(H^*(A))=\{0\}$, so that $\bar\theta_n=0$ for $n\gg0$, and yet
$\bar\theta_n\ne0$ for some $n\ge3$; see
Example~\ref{ex:pencils-2step}\eqref{pen1} below.  Non-$1$-formality is
then detected even though $\B_1(H^*(A))$ has finite length.
\end{remark}

\begin{example}
\label{ex:pencils-2step}
Take $m=5$ and $r=2$, so that $W=\langle w_1,w_2\rangle$ is a pencil in
$\bwedge^2V$, and $\h$ is the corresponding $7$-dimensional $2$-step
nilpotent Lie algebra.  By
Proposition~\ref{prop:B1-2step-presentation}, $\B_1(\CE(\h))\cong\k^2$ in
both cases below, with $\theta_2=2$ and $\theta_n=0$ for $n\ge3$.  On the
cohomological side, the two pencils behave quite differently.
\begin{enumerate}[itemsep=2pt]
\item \label{pen1}
\emph{No decomposable member.}  Let $w_1=a_1\wedge a_2+a_3\wedge a_4$ and
$w_2=a_2\wedge a_3+a_4\wedge a_5$.  For
$w=\lambda w_1+\mu w_2$ one computes
$w\wedge w=2\bigl(\lambda^2\,a_{1234}+\lambda\mu\,a_{1245}
+\mu^2\,a_{2345}\bigr)$, which vanishes only for
$\lambda=\mu=0$; so $W$ contains no nonzero decomposable element, and
$\RR^1(H^*(A))=\{0\}$ by \eqref{eq:res-2step}.  Accordingly
$\B_1(H^*(A))$ has finite length $3$, with
$\bar\theta_2=2$, $\bar\theta_3=1$, and $\bar\theta_n=0$ for $n\ge4$.
\item \label{pen2}
\emph{Two decomposable members.}  Let $w_1=a_1\wedge a_2$ and
$w_2=a_3\wedge a_4$.  Then $\lambda w_1+\mu w_2$ is decomposable if and
only if $\lambda\mu=0$, so \eqref{eq:res-2step} gives
$\RR^1(H^*(A))=\langle a_1,a_2\rangle\cup\langle a_3,a_4\rangle$, of
dimension $2$.  Here $\bar\theta_n=2(n-1)$ for $n\ge2$, of degree $1$ in
$n$, in agreement with part~\eqref{chn2}.
\end{enumerate}
In both cases $\bar\theta_3>\theta_3=0$, so neither algebra is
$1$-formal.  Note that it is the \emph{degenerate}\/ pencil~\eqref{pen2},
not the generic one, whose holonomy Chen ranks grow.
\end{example}

%%%%%%%%%%%%%%%%%%%%%%%%%%%%%%%%%%%
\part{Spaces, groups, and applications}
\label{part:applications}
%%%%%%%%%%%%%%%%%%%%%%%%%%%%%%%%%%%

Having developed the algebraic framework for $\cdgas$ and their associated
Koszul modules, resonance varieties, and holonomy Lie algebras, we now turn
to the topological context and extend the constructions of
Parts~\ref{part:Koszul}--\ref{part:res-holo} to spaces and groups.

The guiding principle of this final part is that many classical invariants of
spaces and groups---such as Alexander invariants, Chen ranks, and their
higher-order analogues---admit effective linearized models governed by the
cohomology algebra and its associated Lie-theoretic structures.
We show how these invariants arise naturally from equivariant (co)chain
complexes of universal abelian covers, and how they can be analyzed using
the same Koszul and holonomy machinery developed earlier.

This perspective allows us to compare algebraic and topological invariants
systematically, clarify the role of minimality and formality assumptions,
and provide a conceptual framework for applications to specific classes of
spaces and groups, developed in the concluding sections of the monograph.

%%%%%%%%%%%%%%%%%%%%%%%%%%%%
\section{The equivariant spectral sequence and its linearization}
\label{sect:equivariant-ss}
%%%%%%%%%%%%%%%%%%%%%%%%%%%

We recall the equivariant spectral sequence of Papadima--Suciu \cite{PS-tams}, 
arising from the $I$-adic filtration on equivariant (co)chains of the universal abelian 
cover. Its $E^2$-page is determined by the cup-product structure on $H^*(X;\k)$, 
while its abutment computes the Alexander invariants $B_i(X;\k)$. The higher 
differentials are interpreted as iterated Massey products, and the spectral 
sequence serves as the main tool for the linearization results of this part.

%%%%%%%%%%%%%%%%%%%
\subsection{The cohomological equivariant spectral sequence}
\label{subsec:coho-equiv-ss}
%%%%%%%%%%%%%%%%%%%

Throughout this section, $X$ will be a connected CW-complex, 
with fundamental group $G$. We will denote by $G_{\ab}=G/[G,G]$ 
the abelianization of this group.

\begin{lemma}
\label{lem:augmentation-filtration}
Assume $G_{\ab}=H_1(X;\Z)$ is finitely generated. 
Let $\Lambda=\k[G_{\ab}]$ and $I\subset\Lambda$ be the augmentation ideal,
where $\k$ is a field of characteristic $0$.

\begin{enumerate}[itemsep=1pt]
\item \label{l1}
There is a canonical isomorphism of graded $\k$-algebras
\[
\gr^I(\Lambda) \cong \Sym(H_1(X;\k)).
\]
\item \label{l2}
The $I$-adic filtration on $\Lambda$ is Hausdorff if and only if
$G_{\ab}$ is torsion-free. 
\end{enumerate}
\end{lemma}

\begin{proof}
Write $G_{\ab}\cong \Z^{n}\oplus H$, with $H$ finite. Then 
$\k[G_{\ab}]$ decomposes as $\k[\Z^n]\otimes_{\k} \k[H]$, 
while $I=I_0+I_H$, where $I_0$ and $I_H$ are the respective 
augmentation ideals.
Since $\ch\k=0$, the averaging idempotent
$e=\frac{1}{|H|}\sum_{h\in H}h$ splits $\k[H]=\k e\oplus I_H$, 
and hence $I_H=(1-e)\k[H]$ is idempotent, so $I_H^2=I_H$.
Therefore $I_H\subset I^p$ for all $p\ge 1$, and the image 
of $\k[H]$ in $\gr^I(\Lambda)$ is concentrated in degree $0$.
It follows that
\[
\gr^I(\Lambda)\cong \gr^{I_0}(\k[\Z^n])\cong \Sym(\k^n)
\cong \Sym(H_1(X;\k)).
\]

If $H=0$, then $\Lambda=\k[\Z^n]$ is a Laurent polynomial 
algebra---in particular, a Noetherian domain---and $I_0$ is a proper ideal, 
so $\bigcap_{p\ge 0} I^p=0$ by Krull's intersection theorem.
If $H\ne 0$, then $\bigcap_{p\ge 0} I^p \supset I_H \ne 0$. 
\end{proof}

\begin{remark}
\label{rem:hausdorff-issue}
Although $\gr^I(\Lambda)$ is always identified with 
$\Sym(H_1(X;\k))$, the $I$-adic filtration on 
$\Lambda$ need not be Hausdorff when $G_{\ab}$ has torsion. 
In this case, the spectral sequence associated to the $I$-adic 
filtration on the Alexander cochain complex $C^\bullet(X;\Lambda)$ 
may fail to converge strongly to $H^\bullet(X;\Lambda)$.

This is a second, independent source of non-separation, distinct 
from the non-separation of Koszul modules discussed in 
Section~\ref{subsec:koszul-filtrations}: there, non-separation 
arose from non-homo\-geneous differentials in the $\cdga$ model 
(as in $\CE(\sol_2)$), and is present even when $G_{\ab}$ is 
torsion-free.
\end{remark}

Let $X^{\ab}$ denote the universal abelian cover of $X$.  We work throughout
with the \emph{equivariant}\/ cochain complex
\begin{equation}
\label{eq:equivariant-cochains}
C^\bullet(X;\Lambda)\coloneqq
\Hom_{\Lambda}\bigl(C_\bullet(X^{\ab};\k),\Lambda\bigr)
\cong \Lambda\otimes_{\k}C^\bullet(X;\k) ,
\end{equation}
whose cohomology is the cohomology $H^\bullet(X;\Lambda)$ of $X$ with local
coefficients in $\Lambda$.  This is the right object here: each
$C^q(X;\Lambda)$ is a finitely generated free $\Lambda$-module, which is what
makes the finiteness arguments below available.  It should not be confused
with the ordinary cochain complex of the cover, whose terms are full
$\k$-duals of infinite-rank free $\Lambda$-modules, and whose cohomology is
different: for $X=S^1$ one has $H^0(X;\Lambda)=0$ and
$H^1(X;\Lambda)=\Lambda/(t-1)\cong\k$, whereas
$H^\bullet(X^{\ab};\k)=H^\bullet(\R;\k)$ is $\k$ in degree $0$ and vanishes
in degree~$1$.  The complex $C^\bullet(X;\Lambda)$
is naturally a filtered $\Lambda$-module via the $I$-adic filtration
\[
F^p C^\bullet(X;\Lambda)= I^p\, C^\bullet(X;\Lambda), \qquad p\ge 0.
\]
This filtration is exhaustive and decreasing, and gives rise to a
cohomological spectral sequence, with terms vanishing outside 
$p\ge0$, $p+q\ge0$,
\[
E_0^{p,q} = \gr^I_p\bigl(C^{p+q}(X^{\ab};\k)\bigr)
\Longrightarrow H^{p+q}(X;\Lambda),
\]
whose existence is formal, and whose convergence will be understood
under the usual finiteness hypotheses on the homology of $X$, 
as specified below. The description of the $E_1$-page 
and its linear differential goes back, in homological form, to work of 
Papadima and Suciu~\cite[Thm.~A]{PS-tams}.
We record below a cohomological version, adapted to our conventions.

\begin{theorem}
\label{thm:equiv-ss-coh}
Let $X$ be a connected CW-complex with $G=\pi_1(X)$ and
$\Lambda=\k[G_{\ab}]$, endowed with the $I$-adic filtration.
Then the filtered cochain complex $C^\bullet(X;\Lambda)$ gives rise to a
cohomological spectral sequence $(E_r^{\bullet,\bullet},d_r)$, with terms 
vanishing outside $p\ge0$, $p+q\ge0$, with the following properties:
\begin{enumerate}[itemsep=3pt]
\item
The $E_0$-page is given by
\[
E_0^{p,q} \cong C^{p+q}(X;\k)\otimes_{\k} \Sym^p(H_1(X;\k)).
\]
\item
The $E_1$-page is canonically identified with the Koszul cochain complex
\[
(E_1^{\bullet,\bullet},d_1)
\cong
K^\bullet\bigl(H^*(X;\k)\bigr)
=
\bigl(H^*(X;\k)\otimes_{\k} S,\ \omega_X\cdot(-)\bigr),
\]
where $S=\Sym(H_1(X;\k))$ and $\omega_X\in H^1(X;\k)\otimes_{\k} H_1(X;\k)$ is the
canonical element.
\end{enumerate}
\end{theorem}

\begin{proof}
Filter the Alexander cochain complex
$C^\bullet = C^\bullet(X;\Lambda)$ by powers of the augmentation 
ideal $I \subset \Lambda = \k[G_{\ab}]$.
This filtration defines a cohomological spectral sequence associated to 
$C^\bullet$, with terms vanishing outside $p\ge0$, $p+q\ge0$.

By Lemma~\ref{lem:augmentation-filtration}\eqref{l1}, the 
associated graded algebra $\gr^I(\Lambda)$ is canonically 
isomorphic to $S = \Sym(H_1(X;\k))$. Hence
\[
\gr^I_p(C^\bullet) \cong
C^\bullet(X;\k)\otimes_{\k} \Sym^p(H_1(X;\k)),
\]
and taking cohomology in the first factor gives
\[
E_1^{p,q} \cong
H^{p+q}(X;\k)\otimes_{\k} \Sym^p(H_1(X;\k)).
\]

Under the identification $\gr^I(\Lambda) \cong S = \Sym(H_1(X;\k))$, 
the action of $I/I^2 \cong H_1(X;\k)$ on the associated graded 
complex corresponds to multiplication by the canonical element 
$\omega_X \in H^1(X;\k) \otimes_{\k} H_1(X;\k)$, giving 
$d_1(\alpha) = \omega_X \cdot \alpha$ for 
$\alpha \in H^{p+q}(X;\k)$. This identifies 
$(E_1^{\bullet,\bullet}, d_1)$ with 
$K^\bullet(H^*(X;\k))$ as in 
Definition~\ref{def:Koszul-coho}. 
\end{proof}

\begin{remark}
\label{rem:ss-convergence}
If $\dim_\k H_q(X;\k)<\infty$ for all $q\ge 0$, then the $I$-adic equivariant
spectral sequence admits an $E_\infty$-term; in particular, this holds whenever 
$X$ is of finite type. This follows from the same reasoning as in 
\cite[Prop.~5.1]{PS-tams}, applied to the coefficient module 
$M=\Lambda$, for which $M/IM\cong\k$.
\end{remark}

By Theorem~\ref{thm:equiv-ss-coh}, the $E_1$-page of the 
equivariant cohomological spectral sequence is canonically 
identified, as a filtered dg-algebra, with the Koszul cochain 
complex $K^\bullet(H^*(X;\k))$ equipped with its $\m$-adic 
filtration, where $S = \Sym(H_1(X;\k))$ and 
$\m \subset S$ is the augmentation ideal.
The cohomological Koszul spectral sequence of 
Theorem~\ref{thm:koszul-ss-coh} arises from precisely 
this filtered complex.

\begin{proposition}
\label{prop:equiv-vs-koszul}
Under the identification of Theorem~\ref{thm:equiv-ss-coh}, the $E_1$-page of the 
equivariant cohomological spectral sequence is canonically isomorphic, as a differential 
graded $S$-module, to the Koszul cochain complex
\[
K^\bullet(H^*(X;\k)) = \bigl(H^*(X;\k)\otimes_{\k} S,\ \omega_X \cdot (-)\bigr).
\]
Consequently, the $E_2$-page of the equivariant spectral sequence is
$\B^{\bullet}(H^*(X;\k))$, and agrees with the $E_2$-page of the Koszul
spectral sequence of any $\cdga$ model of $X$.
\end{proposition}

\begin{proof}
The identification of the associated graded complex
\[
\gr^I C^\bullet(X;\Lambda) \cong C^\bullet(X;\k)\otimes_{\k} S
\]
together with the description of the induced differential yields an
isomorphism of $E_1$-pages as differential graded $S$-modules, by
Theorem~\ref{thm:equiv-ss-coh}.  Taking cohomology of $d_1$ gives the
$E_2$-page in both cases.
\end{proof}

\begin{remark}
\label{rem:equiv-vs-koszul-caveat}
The comparison does \emph{not}\/ extend past $E_2$ on the strength of the
$E_1$-identification alone.  An isomorphism of $E_1$-pages as differential
graded modules determines $E_2=H(E_1,d_1)$, but it does not determine $d_2$,
and hence says nothing about $E_3$.  What propagates to all pages is an
isomorphism induced by a \emph{morphism}\/ of the underlying filtered
complexes, since such a morphism commutes with every $d_r$; supplying one
requires comparing $\Lambda$ with $S$, which is possible after $I$-adic
completion (for $G_{\ab}$ free abelian, $\widehat{\Lambda}\cong\widehat{S}$
by $t_i\mapsto \exp(x_i)$, using $\ch\k=0$) and this is carried out in
Section~\ref{sect:alexinv}.

The distinction is not academic.  Both sets of higher differentials are
given by iterated Massey products along the canonical class
(Theorem~\ref{thm:koszul-ss-coh}\eqref{ssc3}; and \cite{PS-tams},
\cite[\S3]{LMW25} in the rank-one equivariant case), so for a genuine
$\cdga$ model of $X$ they do agree.  But for the \emph{formal}\/ model
$A=(H^*(X;\k),0)$ they cannot: all Massey products in a $\cdga$ with zero
differential vanish, so the Koszul spectral sequence of $(H^*(X;\k),0)$
collapses at $E_2$, whereas the equivariant spectral sequence of a
non-formal $X$ need not---this is precisely the content of
Question~\ref{quest:equiv-degeneration} and of Remark~\ref{rem:lmw}.
\end{remark}

\begin{remark}
\label{rem:equiv-ss-conv}
For $X$ of finite type the spectral sequence has a well-defined
$E_\infty$-page, computing the associated graded of the $I$-adic completion.
Abbreviate $C^\bullet=C^\bullet(X;\Lambda)$ and $H^n=H^n(X;\Lambda)$, and let
$\widehat{\phantom{x}\cdot\phantom{x}}$ denote $I$-adic completion.  Then
\[
E_\infty^{p,q}\cong \gr^I_p\,\widehat{H^{p+q}} .
\]
By \eqref{eq:equivariant-cochains}, each $C^q(X;\Lambda)$ is a finitely
generated free module over the Noetherian ring $\Lambda$, so this follows from
the Artin--Rees lemma exactly as in the
homological case, for which see Corollary~\ref{cor:alex-linearized} and
\cite[Prop.~8.1]{PS-tams}; we do not repeat the argument.  Two cautions are
worth recording.  The filtration is not complete---$\Lambda$ is not
$I$-adically complete, already for $\Lambda=\k[t^{\pm1}]$ and $I=(t-1)$---so
what the spectral sequence computes is the cohomology of the completed complex
rather than $H^\bullet(X;\Lambda)$ itself.  That said, the two descriptions of
the abutment agree here, which is why the display above may be written with the
completion outside: since $\Lambda$ is Noetherian, the map
$\Lambda\to\widehat{\Lambda}$ is flat, and each $C^q$ is finitely
generated, so
\[
H^n\bigl(\widehat{C^\bullet}\bigr)
\cong H^n\bigl(C^\bullet\bigr)\otimes_{\Lambda}\widehat{\Lambda}
\cong \widehat{H^n} ,
\]
the last isomorphism because $H^n$ is again finitely generated.
And the spectral sequence is not bounded above
in the filtration direction: $E_1^{p,\,n-p}=H^{n}(X;\k)\otimes_{\k}\Sym^p(H_1(X;\k))$
is nonzero for every $p\ge0$ whenever $H^n(X;\k)\ne0$, so convergence is not a
formal consequence of boundedness.

We will not pursue the cohomological abutment further.  Its role here is not
to compute anything: the Alexander invariants are homological, and are treated
in Section~\ref{subsec:hom-equiv-ss}.  What the cohomological sequence
supplies, and the homological one does not, is the \emph{multiplicative}\/
structure---its $E_1$-page is a differential graded algebra, which is what
makes the higher differentials iterated Massey products along $[\omega_X]$,
as in Theorem~\ref{thm:koszul-ss-coh}\eqref{ssc3}.
\end{remark}

\begin{remark}
\label{rem:non-Hausdorff}
The equivariant spectral sequence associated to the $I$-adic filtration
need not converge in general. As shown by Papadima--Suciu~\cite{PS-tams}, 
non-convergence may occur when working over the group algebra $\k G$ 
of a non-residually nilpotent group, since the $I$-adic topology on $\k G$ is 
then not Hausdorff. In the present setting, we work with the universal abelian 
cover $X^{\ab}$ and the group algebra $\Lambda=\k[G_{\ab}]$ over a field 
$\k$ of characteristic $0$. Although the $I$-adic filtration on $\Lambda$ 
may fail to be Hausdorff when $G_{\ab}$ has torsion 
(Lemma~\ref{lem:augmentation-filtration}\eqref{l2}), the 
proof of that lemma shows that the torsion contribution is 
confined to filtration degree $0$: the ideal $I_H$ satisfies 
$I_H \subset I^p$ for all $p\ge 1$, so 
$F^p C^\bullet(X;\Lambda) = I^p C^\bullet(X;\Lambda)$ 
is unaffected by the torsion for $p\ge 1$. Hence the 
spectral sequence, which is determined by the filtration 
quotients $F^p/F^{p+1}$ for $p\ge 0$, converges to the 
associated graded of the completed equivariant cohomology 
under the finite-type assumption on $X$.
\end{remark}

%%%%%%%%%%%%%%%%%%%
\subsection{The homological equivariant spectral sequence}
\label{subsec:hom-equiv-ss}
%%%%%%%%%%%%%%%%%%%

We now pass to the homological counterpart of the equivariant spectral sequence.
This change of variance is essential, since Alexander-type invariants
$B_i(X;\k)=H_i(X^{\ab};\k)$ are intrinsically homological objects.

\begin{theorem}
\label{thm:equiv-ss-hom}
Let $X$ be a connected CW-complex with fundamental group $G$, and let
$\Lambda=\k[G_{\ab}]$ with augmentation ideal $I\subset \Lambda$.
The $I$-adic filtration on the equivariant chain complex
$C_\bullet(X^{\ab};\k)$ induces a second-quadrant homological spectral sequence
$(E^r_{p,q},d^r)$ with the following properties:
\begin{enumerate}[itemsep=3pt]
\item
The $E^1$-page is canonically identified with the Koszul chain complex
\[
(E^1_{\bullet,\bullet},d^1)
\cong
K_\bullet\bigl(H^*(X;\k)\bigr),
\]
viewed as a complex of graded $S$-modules, where $S=\Sym(H_1(X;\k))$.
\item
The differential $d^1$ is induced by the comultiplication
\begin{equation*}
\label{eq:nabla-x}
\begin{tikzcd}[column sep=20pt]
\nabla_X \colon  H_*(X;\k) \longrightarrow H_1(X;\k)\otimes_{\k} H_{*-1}(X;\k),
\end{tikzcd}
\end{equation*}
adjoint to the cup product on $H^{*}(X;\k)$.
\end{enumerate}
\end{theorem}

\begin{proof}
The $I$-adic filtration on $C_\bullet(X^{\ab};\k)$ is compatible 
with the boundary maps and hence gives rise to a homological 
spectral sequence. By \cite[Thm.~A]{PS-tams}, its $E^0$-page 
satisfies
\[
\gr^I_p(C_\bullet(X^{\ab};\k))
\cong
C_{\bullet-p}(X;\k)\otimes_{\k} \Sym^p(H_1(X;\k)),
\]
and the identification of $(E^1, d^1)$ with $K_\bullet(H^*(X;\k))$ 
holds with $M = \Lambda$ and the comultiplication $\nabla_X$ 
adjoint to the cup product on $H^*(X;\k)$.
\end{proof}

When $X$ is a finite-type CW-complex, the homological spectral sequence
converges strongly; this is Corollary~\ref{cor:alex-linearized} below, whose
proof rests on the Artin--Rees lemma, as does the cohomological statement
recorded in Remark~\ref{rem:equiv-ss-conv}.

The following two results draw out the key consequences 
of Theorem~\ref{thm:equiv-ss-hom} for the relationship 
between Alexander invariants and Koszul homology. 
For clarity, we write $n=-p\ge 0$ for the Koszul module 
degree, so that $E^2_{-n,q}\cong \B_n(H^*(X;\k))_{q-n}$.

\begin{proposition}
\label{prop:hom-equiv-vs-koszul}
Under the identification of Theorem~\ref{thm:equiv-ss-hom}, the
homological equivariant spectral sequence and the Koszul homology
spectral sequence of $H^*(X;\k)$ have canonically isomorphic $E^1$- and
$E^2$-pages.
In particular, there is an isomorphism of bigraded $\k$-vector spaces
\[
E^2_{-n,q}\cong \B_n\bigl(H^*(X;\k)\bigr)_{q-n}, \quad n\ge 0.
\]
Equivalently, there is an isomorphism of bigraded $S$-modules
\[
E^2 \cong \B_\bullet\bigl(H^*(X;\k)\bigr),
\]
where the bigrading on the left corresponds to the homological grading
and internal grading on the right.
\end{proposition}

\begin{proof}
The $E^2$-page is $\B_\bullet(H^*(X;\k))$, the homology of 
$(E^1,d^1) \cong K_\bullet(H^*(X;\k))$ by 
Definition~\ref{def:Koszul-homology}.
The bigrading $E^2_{-n,q} \cong \B_{n}(H^*(X;\k))_{q-n}$ 
reflects the identification of chain degree $n\ge 0$ with 
the Koszul module index, and $q-n$ with the internal 
graded degree within $\B_{n}(H^*(X;\k))$, as determined 
by the standard bigrading of Koszul homology.
\end{proof}

\begin{corollary}
\label{cor:alex-linearized}
Suppose $X$ is a connected CW-complex of finite type.
Then the homological equivariant spectral sequence converges, and there
is an isomorphism of bigraded $\k$-vector spaces
\[
E^\infty_{-n,q}\cong \gr^I_n\, B_{q-n}(X;\k), \quad n\ge 0.
\]
Equivalently, there is an isomorphism of graded $S$-modules
\[
E^\infty \cong \gr^I B_\bullet(X;\k).
\]
The filtration on $B_\bullet(X;\k)$ need not be separated, in which case
$B_\bullet(X;\k)$ is not determined by $E^\infty$; separation does hold
after $I$-adic completion.
\end{corollary}

\begin{proof}
Under the finite-type hypothesis each chain group
$C_q(X^{\ab};\k)\cong\Lambda\otimes_{\k} C_q(X;\k)$ is a finitely
generated free $\Lambda$-module over the Noetherian ring $\Lambda$, so
the Artin--Rees lemma applies and gives convergence, together with the
separation statement after completion; see \cite[Prop.~8.1]{PS-tams}
and \cite[pp.~22--24]{Serre}.
\end{proof}

\begin{example}
\label{ex:BS-nonseparated}
Let $G=\langle a,t\mid tat^{-1}=a^2\rangle$ and let $X$ be its
presentation $2$-complex, a finite aspherical complex with
$G_{\ab}\cong\Z$.  Then $G'\cong\Z[1/2]$, with $t$ acting by
multiplication by $2$, so $B_1(X;\k)\cong\Lambda/(t-2)$, on which the
augmentation ideal $I=(t-1)$ acts invertibly.  Hence
$I^pB_1(X;\k)=B_1(X;\k)$ for all $p$, so $\gr^I B_1(X;\k)=0$ while
$B_1(X;\k)\ne 0$: the spectral sequence converges, but carries no
information about $B_1$.
\end{example}

\begin{remark}
\label{rem:kh-alex-general}
Corollary~\ref{cor:alex-linearized} identifies $E^\infty_{-n,q}$ with 
$\gr^I_n B_{q-n}(X;\k)$, but the $E^2$-page $\B_n(H^*(X;\k))$ computes 
this associated graded only when the spectral sequence collapses at $E^2$, 
i.e., when all higher differentials $d_r$ ($r\ge 2$) vanish. In general, 
$E^\infty_{-n,m+n}$ is only a subquotient of $E^2_{-n,m+n} \cong 
\B_n(H^*(X;\k))_{m}$.
\end{remark}

\begin{corollary}
\label{cor:formal-collapse}
If $X$ is a connected, finite-type, formal CW-complex, then the homological 
equivariant spectral sequence collapses at $E^2 \cong E^\infty$. Moreover,
for each $n\ge 0$, there are isomorphisms of $S$-modules
\[
E^2_{-n,\bullet} \cong \B_n\bigl(H^*(X;\k)\bigr)
\cong \gr^I_n\, B_\bullet(X;\k),
\]
and hence, for each $q$,
\[
E^2_{-n,q} \cong \B_n\bigl(H^*(X;\k)\bigr)_{q-n}
\cong \gr^I_n\, B_{q-n}(X;\k)
\]
as graded $\k$-vector spaces.
\end{corollary}

\begin{proof}
Since $X$ is formal, it admits the $\cdga$ model $(H^*(X;\k),0)$, and a
formality quasi-isomorphism induces a morphism of filtered complexes from
$K_\bullet(H^*(X;\k))$, with its $\m$-adic filtration, to the equivariant
chain complex with its $I$-adic filtration---after $I$-adic completion, where
$\widehat{\Lambda}\cong\widehat{S}$; compare
Remark~\ref{rem:equiv-vs-koszul-caveat} and
Section~\ref{sect:alexinv}.  Being a morphism of filtered complexes, it
induces a morphism of spectral sequences, an isomorphism on $E^1$ by
Theorem~\ref{thm:equiv-ss-hom}, hence an isomorphism on every page.  The
higher differentials of the homological Koszul spectral sequence of
$(H^*(X;\k),0)$ vanish, by
Theorem~\ref{thm:koszul-ss-hom}\eqref{ssh3}, since the differential of the
model is zero; hence $d^r=0$ for all $r\ge2$.  The collapse 
$E^2 \cong E^\infty$ then follows, and combined with 
Corollary~\ref{cor:alex-linearized} gives the stated identification.
\end{proof}

\begin{remark}
\label{rem:alex-separation}
Corollary~\ref{cor:formal-collapse} identifies $\gr^I_n B_m(X;\k)$ with 
$\B_n(H^*(X;\k))_m$ for formal $X$, but recovering $B_m(X;\k)$ itself 
from its associated graded requires the $I$-adic filtration to be separated, 
i.e., $\bigcap_{p\ge 0} I^p B_m(X;\k) = 0$. This holds, for instance, 
when $B_m(X;\k)$ is a finitely generated torsion-free $\Lambda$-module. 
This separation condition is the topological counterpart of the separation 
of Koszul modules discussed in Section~\ref{subsec:koszul-filtrations}.
\end{remark}

As observed in Remark~\ref{rem:hausdorff-issue}, torsion in 
$H_1(X;\Z)$ causes the $I$-adic filtration on $\Lambda$ itself 
to be non-Hausdorff, raising the question of whether this 
can obstruct separation of associated $\Lambda$-modules 
such as the Alexander invariants. This naturally leads 
to the following question.

\begin{question}
\label{quest:alex-separated}
Let $X$ be a connected finite CW-complex with $H_1(X;\Z)$ torsion-free. 
Is the $I$-adic filtration on the Alexander invariant 
$B_m(X;\k)$ separated for all $m\ge 1$? More precisely, does
\[
\bigcap_{p\ge 0}\, I^p B_m(X;\k) = 0?
\]

A positive answer is equivalent to injectivity of 
the natural $\Lambda$-linear map 
\[
B_m(X;\k)\longrightarrow \widehat{B_m(X;\k)} ,
\]
where $\widehat{B_m(X;\k)}$ is viewed as a 
$\Lambda$-module via the canonical map $\Lambda\to \widehat{\Lambda}$.
If $X$ is moreover formal, then the identification of $\gr^I_\bullet B_m(X;\k)$ 
with the Koszul modules $\B_\bullet(H^*(X;\k))$ shows that the latter determine 
the completed invariant $\widehat{B_m(X;\k)}$ as a $\widehat{S}$-module, 
into which $B_m(X;\k)$ embeds via the natural map.

In the special case where $X$ is the complement of a hyperplane 
arrangement, this question was raised 
in~\cite[Quest.~10.10]{Su-decomp} in the context of 
decomposable arrangements, where a positive answer would remove 
the separation hypothesis from several comparison theorems for 
Alexander invariants and Chen ranks.
\end{question}

We conclude this subsection by placing the equivariant spectral 
sequence in its broader literature context.

\begin{remark}
\label{rem:lmw}
Liu--Maxim--Wang~\cite{LMW25} revisit the equivariant spectral sequence
of~\cite{PS-tams} for an infinite cyclic cover---that is, for an
epimorphism $\nu\colon G\twoheadrightarrow\Z$ and the corresponding
rank-one coefficient ring $\k[t^{\pm1}]$---and show that over an
arbitrary field all of its higher differentials are computed by iterated
Massey products along the class $\nu\in H^1(X;\k)$ \cite[\S3]{LMW25}.
This is the rank-one instance of the pattern recorded in
Theorem~\ref{thm:koszul-ss-coh}\eqref{ssc3}; the multivariable case, for
the universal abelian cover, is the subject of
Question~\ref{quest:equiv-degeneration}.
From the homological perspective adopted here, the same
information is encoded by iterated applications of the
coalgebra structure $\nabla_X$ on $H_*(X;\k)$, naturally
compatible with Koszul homology.
\end{remark}

A closely related spectral sequence, originating with Novikov~\cite{Novikov} 
in the de Rham setting and developed further by Farber~\cite{Farber} and 
Pajitnov~\cite{Pajitnov}, is studied systematically by Kohno and 
Pajitnov~\cite{KP}. Given a finite CW-complex $X$, a representation 
$\rho\colon G\to \GL(n,\C)$, and a class $\alpha\in H^1(X;\C)$, 
they construct a spectral sequence starting from $H^*(X;\rho)$ and 
converging to the cohomology $H^*(X;\gamma_{\mathrm{gen}})$ with 
coefficients in the generic deformation 
$\gamma_t(g)=e^{t\langle\alpha,g\rangle}\rho(g)$.
The differentials are given by iterated Massey products with $\alpha$, 
paralleling the structure of the higher differentials in 
Theorems~\ref{thm:koszul-ss-coh} and~\ref{thm:equiv-ss-coh}.
When $\rho$ is the trivial representation, their spectral sequence 
degenerates for formal spaces by the classical DGMS argument~\cite{DGMS}, 
yielding $\ker L_\alpha / \im L_\alpha \cong H^*(X;\gamma_{\mathrm{gen}})$.

The present paper may be viewed as the Koszul-theoretic 
counterpart of this circle of ideas: working over an arbitrary 
$\k$-$\cdga$ with the trivial representation, we replace the 
exponential deformation by the canonical element $\omega_A$, 
and gain an explicit $S$-module structure on the abutment 
governing Alexander-type invariants.
The equivariant spectral sequence provides a conceptual bridge from the
cohomology ring $H^*(X;\k)$ to the homological Alexander invariants
$B_i(X;\k)=H_i(X^{\ab};\k)$, but it does not, by itself, yield explicit or
functorial algebraic models for these invariants.
To pass from this topological framework to computable objects, one must
linearize the equivariant theory and compare it with algebraic constructions
arising from $\k$-$\cdga$ models.
This requires two additional steps: first, replacing equivariant (co)chains by
their associated graded objects with respect to the augmentation filtration;
and second, passing to suitable completions to control convergence and
functoriality.
The outcome is a precise identification between Alexander invariants and
Koszul homology modules, developed in the next section.

%%%%%%%%%%%%%%%%%%%%%%%%%%%%%%
\section{Alexander invariants and infinitesimal models}
\label{sect:alexinv}
%%%%%%%%%%%%%%%%%%%%%%%%%%%%%%

In this section we recall the classical Alexander invariants of a space and
introduce their infinitesimal analogues arising from Koszul modules.
These infinitesimal invariants provide linearized and functorial models for
Alexander modules, expressed entirely in terms of the cohomology algebra 
and the associated Koszul differential. We show that, after suitable completion, 
the classical and infinitesimal theories are canonically identified at the level 
of associated graded modules, under minimality and finite-type hypotheses 
on $X$. The main result of this section establishes 
Theorem~\ref{thm:top-linearization-intro}\eqref{C1} from the Introduction. 
This comparison clarifies the role of Koszul modules as algebraic avatars 
of Alexander invariants and prepares the ground for the Chen ranks 
computations that follow.

%%%%%%%%%%%%%%%%
\subsection{Alexander invariants}
\label{subsec:alexinv}
%%%%%%%%%%%%%%%%%

Let $X$ be a connected CW-complex, with a single $0$-cell
chosen as basepoint $x_{0}$.  Set $G=\pi_{1}(X,x_{0})$ and $G_{\ab}=G/[G,G]$. 
The group algebra $\Lambda=\k[G_{\ab}]$ is a commutative Noetherian ring, 
since $G_{\ab}$ is a finitely generated abelian group.
Lifting the cellular structure to the maximal abelian cover $X^{\ab}$, we
obtain an augmented chain complex of free $\Lambda$-modules,
\begin{equation}
\label{eq:abcover-cc}
\begin{tikzcd}[column sep=22pt]
\cdots \ar[r] & 
C_{2}(X^{\ab};\k) \ar[r, "\partial^{\ab}_{2}"] & 
C_{1}(X^{\ab};\k) \ar[r, "\partial^{\ab}_{1}"] & 
C_{0}(X^{\ab};\k)=\Lambda \ar[r, "\varepsilon"] & \k \ar[r] & 0,
\end{tikzcd}
\end{equation}

\begin{definition}
\label{def:alexinv}
The \emph{Alexander invariants} of $X$ (over $\k$) are the homology modules 
\begin{equation}
\label{eq:alexinv-X}
B_{i}(X;\k) \coloneqq H_{i}(X^{\ab};\k), \qquad i\ge 0,
\end{equation}
viewed as modules over the group algebra $\Lambda=\k[G_{\ab}]$. 
When no confusion is possible, we suppress $\k$ from the notation.
\end{definition}

If $X$ has finitely many $i$-cells, then $B_i(X)$ is a finitely generated 
$\Lambda$-module. In particular, if $X$ is a finite-type CW-complex, then 
all Alexander invariants $B_i(X)$ are finitely generated $\Lambda$-modules.

This construction enjoys the following naturality property. 
Let $f\colon X\to \overline{X}$ be a basepoint-preserving cellular map 
between connected CW-complexes.  The induced map on fundamental 
groups passes to abelianizations and yields a morphism of group algebras
\begin{equation}
\label{eq:lambda-to-ovlambda}
\varphi\colon \Lambda=\k[G_{\ab}] \longrightarrow
\overline{\Lambda}=\k[\overline{G}_{\ab}].
\end{equation}
The lift $f^{\ab}\colon X^{\ab}\to \overline{X}^{\ab}$ 
induces $\Lambda$-linear maps
\[
B_{i}(f)\colon B_{i}(X) \longrightarrow B_{i}(\overline{X}),
\]
where $B_i(\overline{X})$ is viewed as a $\Lambda$-module 
via the ring map $\varphi \colon \Lambda \to \overline{\Lambda}$.
Thus $X \mapsto B_i(X)$ defines a functor from connected 
based CW-complexes to modules over group algebras.

%%%%%%%%%%%%%%%%%%%%%
\subsection{Infinitesimal Alexander invariants}
\label{subsec:inf-alexinv}
%%%%%%%%%%%%%%%%%%%%%

We now introduce the infinitesimal Alexander invariants, defined via the 
Koszul modules of the cohomology algebra $H^*(X;\k)$, and then extend 
the construction to arbitrary $\cdga$ models. Unlike the classical literature, which 
focuses primarily on the first Alexander invariant $B_1(X)$, the Koszul formalism 
naturally produces invariants in all degrees. The results below show that these 
higher infinitesimal invariants model, after completion and linearization, the 
full tower of Alexander modules $B_i(X)$.

\begin{definition}
\label{def:inf-alexinv}
Let $X$ be a connected CW-complex with $\dim_\k H^1(X;\k)<\infty$.
The \emph{infinitesimal Alexander invariants} of $X$ (over $\k$) are the
Koszul homology modules
\[
\B_i(X;\k) \coloneqq \B_i\bigl(H^{*}(X;\k)\bigr), \qquad i\ge 0,
\]
viewed as graded modules over the symmetric algebra
$S=\Sym(H_1(X;\k))$. 
When no confusion is possible, we suppress $\k$ from the notation.
\end{definition}

If $X$ has finitely many $i$-cells, then $\B_i(X)$ is a finitely generated 
$S$-module. In particular, if $X$ is a finite-type CW-complex, then 
all infinitesimal Alexander invariants $\B_i(X)$ are finitely generated 
$S$-modules.

Note that $H^*(X;\k)$ with zero differential is a $\cdga$, but it is a model 
for $X$ only when $X$ is formal. The definition above is thus purely
cohomological and requires no formality assumption; formality will only
enter when comparing $\B_i(X)$ with the classical Alexander invariants 
$B_i(X)$. 

More generally, let $(A^\bullet,d_A)$ be a connected $\k$-$\cdga$ with
$\dim_\k H^1(A)<\infty$.
Recall from Section~\ref{subsec:hom-koszul} the construction of the Koszul
chain complex $K_\bullet(A)$ and its homology modules
\begin{equation}
\label{eq:bia-hik}
\B_i(A) \coloneqq H_i(K_\bullet(A)), \qquad i\ge 0,
\end{equation}
which are graded modules over $S=\Sym(H_1(A))$.
If $H^*(A)$ is of finite type, then each $\B_i(A)$ is a finitely generated
$S$-module.
In particular, when $A=H^{*}(X;\k)$ with zero differential, one recovers
the infinitesimal Alexander invariants $\B_i(X;\k)$ defined above.

We may also consider the modules $\B_i(A)$ associated to partial $\cdga$ 
models. If $A$ is a $q$-model for a $q$-finite CW-complex $X$, 
the Koszul modules $\B_i(A)$ are finitely generated $S$-modules for 
all $i\le q$, but they may depend on the choice of model.
Unlike the classical Alexander invariants, the Koszul modules $\B_i(A)$
associated to a $\cdga$ model need not be homotopy invariants of $X$.
This dependence on the chosen model is intrinsic and already appears for
formal spaces.

\begin{example}
\label{ex:circle-models}
Let $X=S^1$. This is a formal space, modeled by its cohomology algebra, 
$H^*(X;\k)=\bigwedge(a)$ with $\deg(a)=1$, for which $\B_1(X;\k)=0$. 
But  $X$ is also modeled by the $\cdga$ $A=\CE(\sol_2)$ of
Example~\ref{ex:sol2-koszul}, for which $\B_1(A)\cong \k[x]/(x-1)$. 
\end{example}

One could similarly define cohomological Alexander invariants
$B^i(X;\k)=H^i(X^{\ab};\k)$ and their infinitesimal analogues via Koszul
cohomology.  While such invariants fit naturally into the equivariant
cohomological spectral sequence, they play no role in the group-theoretic
applications considered here, which are inherently homological in nature.

The infinitesimal Alexander invariants $\B_i(X;\k)$ should be viewed as
linearized counterparts of the classical Alexander invariants $B_i(X;\k)$.
A precise comparison between the two theories requires additional structure,
and will be carried out in the two subsections that follow.
Under minimality and finite-type hypotheses on the CW-structure of $X$, 
the associated graded of the equivariant chain complex is canonically 
isomorphic to the Koszul complex of $H^*(X;\k)$, and this upgrades 
to a filtered quasi-isomorphism of completions.

%%%%%%%%%%%%%%%%%%%%%%%
\subsection{Linearization and minimality}
\label{subsec:minimality}
%%%%%%%%%%%%%%%%%%%%%%%

Let $X$ be a connected CW-complex of finite type, with cellular chain complex 
$(C_*(X;\Z),\partial)$ and boundary maps 
$\partial_q\colon C_q(X;\Z)\to C_{q-1}(X;\Z)$.
We also consider the corresponding chain and cochain complexes with 
coefficients in a field $\k$, denoted $(C_*(X;\k),\partial_X)$ and $(C^*(X;\k),\delta_X)$.

Let $G=\pi_1(X)$, with abelianization $G_{\ab}$, and let 
$\Lambda=\k[G_{\ab}]$ be the group algebra with augmentation ideal 
$I\subset \Lambda$. The $I$-adic filtration on $\Lambda$ induces a 
filtration on the equivariant chain complex $C_\bullet(X^{\ab};\k)$,
\[
F^n C_\bullet(X^{\ab};\k)= I^n \cdot C_\bullet(X^{\ab};\k),
\]
which is compatible with the boundary maps $\partial^{\ab}_\bullet$, 
and thus gives rise to an associated graded chain complex.

The following elementary observation, due to Papadima--Suciu 
\cite[Lem.~3.2]{PS-tams}, identifies this associated graded complex.

\begin{lemma}
\label{lem:gr-ab-boundary}
With respect to the $I$-adic filtration:
\begin{enumerate}[itemsep=2pt]
\item \label{lin1}
$\gr^I(\partial^{\ab}) = \partial_X \otimes \id_S$, so
$\gr^I C_\bullet(X^{\ab};\k) \cong (C_\bullet(X;\k)\otimes_{\k} S,\,\partial_X\otimes\id_S)$.
\item \label{lin2}
$\gr^I(\delta^{\ab}) = \delta_X \otimes \id_S$, so
$\gr^I C^\bullet(X;\Lambda) \cong (C^\bullet(X;\k)\otimes_{\k} S,\,\delta_X\otimes\id_S)$.
\end{enumerate}
\end{lemma}

\begin{proof}
Part \eqref{lin1} is the specialization of \cite[Lem.~3.2]{PS-tams} to the 
$\k[G]$-module $M=\Lambda$, where $\Lambda=\k[G_{\ab}]$. 
Under the identification $C_q(X^{\ab};\k)=\Lambda\otimes_{\k} C_q(X;\k)$, 
the $I$-adic filtration is given by $F^n=I^n\otimes_{\k} C_\bullet(X;\k)$, 
and the claim follows by passing to associated graded. Part \eqref{lin2} 
follows at once. 
\end{proof}

Lemma~\ref{lem:gr-ab-boundary} shows that the associated graded 
equivariant chain complex is obtained from the ordinary cellular 
chain complex $(C_\bullet(X;\k),\partial_X)$ by extension of scalars 
from $\k$ to $S$. In particular, the degree-$0$ differential in the 
$I$-adic filtration (equivalently, the $d^0$ differential in the 
equivariant spectral sequence) is induced by $\partial_X$.

Thus, the cellular boundary maps of $X$ control the zeroth-order 
behavior of the equivariant differential $\partial^{\ab}$. 
Eliminating this contribution greatly simplifies the associated 
graded complex, and motivates the following definition.

\begin{definition}
\label{def:minimality}
A connected, finite-type CW-complex $X$ is called \emph{minimal} if its cellular 
boundary maps vanish, i.e.,
\[
\partial_q \colon C_q(X;\Z)\longrightarrow C_{q-1}(X;\Z)
\quad\text{vanishes for all } q\ge 1.
\]
Equivalently, for every $q\ge 0$, the number of $q$-cells of $X$ is equal to 
the Betti number $b_q(X)=\rank H_q(X;\Z)$.
In view of Lemma~\ref{lem:gr-ab-boundary}, this condition is equivalent to
\[
\gr^I(\partial^{\ab})=0,
\]
that is, the equivariant boundary operator has no term of degree $0$ 
with respect to the $I$-adic filtration.
\end{definition}

If $X$ is a finite-type, minimal CW-complex, then 
$C_*(X;\Z)\cong H_*(X;\Z)$ as chain complexes, 
and in particular each $H_q(X;\Z)$ is a free abelian group.
Consequently, the abelianization $G_{\ab}=H_1(X;\Z)$ is free abelian,
the group algebra $\Lambda=\k[G_{\ab}]$ is a Laurent polynomial algebra, and
after tensoring with $\k$,
\[
C^\bullet(X;\k)\cong H^{*}(X;\k)
\]
as graded vector spaces. Minimality plays a threefold structural role in 
what follows. First, it implies that $G_{\ab}$ is free abelian, so that 
$\Lambda$ is a Laurent polynomial algebra and the $I$-adic 
filtration is Hausdorff. Second, it ensures that the linearization 
of the equivariant coboundary operator $\delta^{\ab}$---its associated 
graded map with respect to the $I$-adic filtration---is realized by an actual 
map of filtered complexes (Theorem~\ref{thm:lin-koszul}). Third, under 
the finite-type hypothesis, it guarantees that cochain--chain duality 
commutes with $I$-adic completion degreewise.

Minimal CW-complexes arise in a variety of geometric contexts.
Basic examples include spheres, tori, and closed orientable surfaces,
all endowed with their standard CW-structures. Complements of complex 
hyperplane arrangements also admit minimal CW-structures, as shown 
by Dimca--Papadima \cite{DP03} and Randell \cite{Randell}.

Minimality is preserved under finite wedges and under products 
(with the product CW-structure) of finite-type, minimal CW-complexes.

The following examples illustrate what goes wrong in the absence of minimality.

\begin{example}
\label{ex:RP2}
Let $X=\RP^2$. Then $\pi_1(X)=\Z_2$; hence, 
$\Lambda=\k[\Z_2]\cong \k[\alpha]/(\alpha^2-1)$, and 
$I=(\alpha-1)$. Assuming $\ch\k=0$, one has $I^2 = 2I = I$,
since $2$ is invertible in $\k$, and hence
the $I$-adic filtration on $\Lambda$ is not Hausdorff, in accordance with 
Lemma~\ref{lem:augmentation-filtration}\eqref{l2}.
The image of $\alpha-1$ in the equivariant cellular chain complex
\[
C_*(X^{\ab};\k)\colon
\begin{tikzcd}[column sep=20pt]
\Lambda \arrow[r, "1+\alpha\,"] &[8pt] 
\Lambda \arrow[r, "\alpha-1\,"] &[8pt] 
\Lambda \arrow[r, "\varepsilon"] & \k \arrow[r] & 0
\end{tikzcd}
\]
is contained in $\bigcap_{p} I^p C_*(X^{\ab};\k)$, so the $I$-adic filtration 
on equivariant chains is likewise non-Haus\-dorff, and the equivariant 
spectral sequence fails to converge strongly. This occurs despite the 
formal identification $\gr^I(\Lambda)\cong \Sym(H_1(X;\k))=\k$.
\end{example}

\begin{example}
\label{ex:trefoil}
Let $X$ be the complement of the trefoil knot in $S^3$. Then
$G_{\ab}\cong \Z$, so $\Lambda=\k[t^{\pm1}]$ is a Laurent polynomial 
ring and the $I$-adic filtration is Hausdorff. However, $X$ is 
homotopy equivalent to the presentation $2$-complex of 
$G=\langle a,b\mid aba=bab\rangle$, which is not minimal.
The equivariant cellular chain complex of the universal abelian cover is
\[
C_*(X^{\ab};\k)\colon
\begin{tikzcd}[column sep=20pt]
\Lambda \arrow[rr, "(1-t+t^2 \;\; -1+t-t^2)"] &&[30pt]
\Lambda^2 \arrow[r, "\varepsilon"] & \k \arrow[r] & 0,
\end{tikzcd}
\]
and passing to the associated graded complex yields a nontrivial boundary 
map in degree $2$. On the other hand, $H^*(X;\k)\cong \bigwedge\langle t\rangle$,
so $K_\bullet(H^*(X;\k))$ has no chain group in degree $2$, giving
\[
\gr^I C_\bullet(X^{\ab};\k)\not\cong K_\bullet(H^*(X;\k)).
\]
The cochain-level identification 
$\gr^I C^\bullet(X;\Lambda) \cong K^\bullet(H^*(X;\k))$ 
requires minimality (Theorem~\ref{thm:lin-koszul}) and fails here. 
This shows that minimality is genuinely needed for the chain-level 
comparisons, beyond mere torsion-freeness of $G_{\ab}$. 
\end{example}

%%%%%%%%%%%%%%%%%
\subsection{Linearization theorem and Alexander invariants}
\label{subsec:linearization-alex}
%%%%%%%%%%%%%%%%%

We now describe the classical linearization results for Alexander-type invariants 
of minimal CW-complexes. The key idea is that, under minimality, the $I$-adic filtration 
on the equivariant (co)chain complex linearizes the boundary maps, and the resulting 
linearized complex is canonically isomorphic to the Koszul complex of the cohomology 
algebra. This is the associated-graded version of the Papadima--Suciu theorem \cite{PS-tams}.

\begin{theorem}[{\cite[Thm.~12.6]{PS-tams}}]
\label{thm:lin-koszul}
Let $X$ be a connected, finite-type, minimal CW-complex. 
There is a canonical isomorphism of bigraded dg-$S$-modules
\[
\gr^I C^\bullet_\Lambda(X^{\ab};\k) \cong K^\bullet(H^*(X;\k)),
\]
where $S=\Sym(H_1(X;\k))\cong\gr^I(\Lambda)$, 
$K^\bullet(H^*(X;\k)) = \bigl(H^*(X;\k)\otimes_{\k} S,\, \omega_X \cdot (-)\bigr)$ 
is the Koszul complex, and $\omega_X \in H^1(X;\k) \otimes_{\k} H_1(X;\k)$ 
is the canonical element.
\end{theorem}

The identification follows from the cohomological $I$-adic spectral sequence
(Theorem~\ref{thm:equiv-ss-coh}) and the minimality assumption. 
Minimality enters in two ways:
\begin{itemize}[itemsep=1pt]
\item \label{min1}
By minimality, the cellular boundary maps vanish ($\partial = 0$).
By Lemma~\ref{lem:gr-ab-boundary}, $\gr^I(\partial^{\ab}) = \partial_X \otimes \id_S = 0$,
so the $d^0$ differential in the equivariant spectral sequence vanishes,
giving $E^0 = E^1$ as graded complexes.

\item \label{min2}
Minimality identifies cochains with cohomology,
$C^q(X;\k) \cong H^q(X;\k)$. Under this identification, the $E^1$ differential
of the spectral sequence (Theorem~\ref{thm:equiv-ss-coh}) is
$\omega_X \cdot (-)$, giving the Koszul complex.
\end{itemize}

\begin{remark}
\label{rem:split-grading}
The Koszul complex carries a natural grading by $S$-degree, so the $\m$-adic 
filtration is split:
\[
K^\bullet(H^*(X;\k)) \cong \bigoplus_{p\ge 0} \gr^\m_p K^\bullet(H^*(X;\k)).
\]
Together with Theorem~\ref{thm:lin-koszul}, this gives canonical isomorphisms 
of bigraded $S$-modules
\[
\gr^I C^\bullet(X;\Lambda) \cong K^\bullet(H^*(X;\k)) \cong \gr^\m K^\bullet(H^*(X;\k)),
\]
where the first is a dg-$S$-module isomorphism and the second is a bigraded 
$S$-module isomorphism (the $\m$-adic associated graded kills the Koszul differential).
\end{remark}

\begin{theorem}
\label{thm:lin-koszul-degenerate}
Let $X$ be a connected, finite-type CW-complex with $G=\pi_1(X)$ and
$G_{\ab}$ torsion-free, and let $q\ge 0$.  Assume that the homological
equivariant spectral sequence of Theorem~\ref{thm:equiv-ss-hom}
degenerates at the $E^2$-page in total degrees $\le q$.  Then, for all
$i\le q$, there are isomorphisms
\[
\gr^I B_i(X;\k) \cong \B_i\bigl(H^*(X;\k)\bigr),
\qquad
\widehat{B}_i(X;\k)_I \cong \B_i\bigl(H^*(X;\k)\bigr)\otimes_{S}\widehat{S},
\]
the first of graded $S$-modules and the second of $\widehat{S}$-modules.
\end{theorem}

\begin{proof}
By Theorem~\ref{thm:equiv-ss-hom} the $E^1$-page is
$K_\bullet(H^*(X;\k))$, with $d^1$ adjoint to the cup product, so
$E^2\cong\B_\bullet(H^*(X;\k))$ by
Proposition~\ref{prop:hom-equiv-vs-koszul}.  Degeneration in total
degrees $\le q$ gives $E^\infty\cong E^2$ in that range.  The spectral
sequence converges by \cite[Prop.~8.1]{PS-tams}, so
$E^\infty\cong\gr^I B_\bullet(X;\k)$, and the first isomorphism follows.
Since $\B_i(H^*(X;\k))$ is a graded $S$-module, an isomorphism on
associated graded objects induces one on completions, giving the second.
\end{proof}

\begin{remark}
\label{rem:degeneration-not-minimality}
Theorem~\ref{thm:lin-koszul-degenerate} requires no minimality
assumption: the identification of the $E^1$-page in
Theorem~\ref{thm:equiv-ss-hom} holds for an arbitrary connected
CW-complex of finite type.  Minimality enters at a different point,
namely in the chain-level linearization statement
\cite[Thm.~12.6]{PS-tams} recorded in Theorem~\ref{thm:lin-koszul}, and
it does not by itself force degeneration.  Indeed, let $X$ be the
presentation $2$-complex of the integer Heisenberg group
$\Gamma=\langle a,b \mid [a,[a,b]],\, [b,[a,b]]\rangle$.  Both relators
lie in $\gamma_3(F)$, so the cellular boundary maps of the universal
abelian cover vanish modulo $I^2$, and $b_1(X)=b_2(X)=2$; hence $X$ is
minimal, and all cup products of degree-one classes vanish.  Since
$\Gamma'\cong\Z$ is central, $B_1(X;\k)=\k$ and $\gr^I B_1(X;\k)=\k$,
whereas $\B_1(H^*(X;\k))=\ker(S^2\to S)\cong S(-1)$ has Hilbert series
$t/(1-t)^2$.  These are not isomorphic, so the spectral sequence does
not degenerate in total degree~$1$.  The obstruction is the failure of
$1$-formality of $\Gamma$, not the failure of minimality.

The same distinction governs the passage from a $\cdga$ model to its
cohomology algebra.  Replacing a finite-type model $A$ by $(H^*(A),0)$
changes the completed Alexander invariant unless $A$ is formal in the
relevant range, and positive weights do not repair this
\textup{(}Remark~\ref{rem:pw-not-enough}\textup{)}: for the Heisenberg
nilmanifold $M(1)$, the invariant $\widehat{B}_1(M(1);\k)_I$ is computed
correctly from the Chevalley--Eilenberg model $\CE(\h(1))$, by
Theorem~\ref{thm:completed-koszul-space}\eqref{ca1}, but not from $H^*(M(1);\k)$
\textup{(}Example~\ref{ex:heisenberg-nonminimal}\textup{)}.
\end{remark}

%%%%%%%%%%%%%%%%%%%%%%%%%%%%
\section{Algebraic models for spaces}
\label{sect:spaces}
%%%%%%%%%%%%%%%%%%%%%%%%%%%%

We recall minimal and finite-type $\cdga$ models and formality. 
Using Koszul techniques, we introduce a spectral-sequence 
obstruction to $q$-formality arising from positive weight 
decompositions, and establish part~\eqref{C2} of
Theorem~\ref{thm:top-linearization-intro}: under a degeneration
hypothesis, Alexander invariants agree with Koszul homology after
completion.

%%%%%%%%%%%%%%%%%%%%
\subsection{Minimal models}
\label{subsec:minimal-models}
%%%%%%%%%%%%%%%%%%%%

A connected $\cdga$ $(\M,d)$ is {\em minimal}\/ if 
$\M=\bigwedge V$ is free as a graded algebra and the differential 
is decomposable, that is, $d\M \subset \M^+ \cdot \M^+$. 
Equivalently, $\M$ is obtained from $\k$ by a sequence of 
Hirsch extensions, 
\begin{equation}
\label{eq:min-mod}
\k=\M(0)\subset \M(1)\subset \M(2)\subset \cdots , 
\qquad 
\M(i)=\M(i-1)\otimes_{\k} \bwedge V_i,
\end{equation}
with $d(V_i) \subset \M(i-1)^{\ge 2}$.  
Given any $\cdga$ $(A,d_A)$ with $H^0(A)=\k$, there exists a 
minimal model $\rho\colon \M(A)\to A$, unique up to isomorphism 
\cite{Sullivan77, Mo}.  For each $q\ge 1$ there is also a {\em $q$-minimal model}, 
generated in degrees $\le q$, well-defined up to $q$-isomorphism. 
A $\cdga$ $(A,d_A)$ is {\em formal} (respectively $q$-formal) precisely when 
$\M(A)$ (resp.\ $\M_q(A)$) is quasi-isomorphic (resp.\ $q$-quasi-isomorphic) 
to $(H^{*}(A),0)$.  

Assume $\M$ is generated in degree $1$.  The filtration \eqref{eq:min-mod} 
induces a filtration $\M^1=\bigcup_{i} \M(i)^1$, where each inclusion 
$\M(i-1)\subset \M(i)$ is a Hirsch extension with fiber 
\begin{equation}
\label{eq:ker-h2}
V_i=\ker\big(H^2(\M(i-1))\longrightarrow H^2(\M)\big).
\end{equation}
If each $V_i$ is finite-dimensional, the dual vector spaces 
$\fL_i=(\M(i)^1)^\vee$ inherit Lie brackets via 
\begin{equation}
\label{eq:dual-lie-bracket}
\langle [u^\vee,v^\vee],w\rangle
  =\langle u^\vee\wedge v^\vee, dw\rangle .
\end{equation}
Each projection $\fL_i\twoheadrightarrow \fL_{i-1}$ has central kernel 
$V_i^\vee$, producing a tower of finite-dimensional nilpotent Lie algebras
$\fL_1 \twoheadleftarrow \fL_2 \twoheadleftarrow \cdots$, 
whose inverse limit $\fL(\M)$ is a complete, filtered Lie algebra.
Conversely, any such tower arises from a unique $1$-minimal $\cdga$. 
This construction provides an algebraic description of $1$-minimal $\cdgas$ 
in terms of towers of nilpotent Lie algebras.

%%%%%%%%%%%%%%%%%%%%
\subsection{Algebraic models for spaces}
\label{subsec:algmod}
%%%%%%%%%%%%%%%%%%%%

Given a space $X$, we let $A_{\PL}^{\bullet}(X)$ be the 
commutative differential graded $\Q$-algebra of 
rational polynomial forms, as defined by Sullivan in \cite{Sullivan77}.  
There is then a natural isomorphism 
$H^{*}(A_{\PL}(X)) \cong H^{*}(X,\Q)$ under which the respective 
induced homomorphisms in cohomology correspond. 

Now let $\k$ be a field of characteristic $0$, and $q$ a positive integer.  
We say that a $\k$-$\cdga$ $(A,d)$ is a {\em $q$-model} (or simply a model)
over $\k$ for $X$ if $A$ is $q$-equivalent (or weakly equivalent) to
$\apl(X)\otimes_{\Q}\k$. 
Finite-type $\cdga$ models are essential in what follows, since all resonance 
and support loci are defined via finitely generated Koszul complexes.

If $f\colon X\to Y$ is a map between two spaces admitting $\k$-$\cdga$ models
$A_X$ and $A_Y$, we say that a $\cdga$ morphism $\varphi\colon A_Y\to A_X$
is a {\em $q$-model} (or model) for $f$ if it is $q$-equivalent
(or weakly equivalent) to $\apl(f)\otimes_{\Q}\k$.  
A {\em functorial model} is a functor $X \leadsto A_X$ assigning to each
continuous map $f\colon X\to Y$ a $\cdga$ morphism
$A(f)\colon A_Y\to A_X$ which models $f$ in the above sense.
Functorial $q$-models are defined analogously. 
By considering classifying spaces $K(G,1)$ for groups $G$, 
we may speak about (functorial) $1$-models for groups and 
group homomorphisms. 
For instance, if $X$ is a smooth manifold, then $\Omega^{\bullet}_{\dR}(X)$,  
the de Rham algebra of smooth forms on $X$, is a model of $X$ over $\R$.  
This model is functorial: if $f\colon X\to Y$ is a smooth 
map, then $f^*\colon \Omega^{\bullet}_{\dR}(Y) \to 
\Omega^{\bullet}_{\dR}(X)$ is a $\cdga$ morphism that models $f$. 

A {\em rational homotopy equivalence} between $X$ and $Y$ is a zig-zag of 
maps (going either way) connecting them via rational quasi-isomorphisms.  
Thus $X$ and $Y$ have the same rational homotopy type precisely when such 
a zig-zag exists.  By functoriality, a rational quasi-isomorphism 
$f\colon X\to Y$ induces a quasi-isomorphism $\apl(f)\colon \apl(Y)\to \apl(X)$.
Consequently, the weak isomorphism type of $\apl(X)$ depends only on the rational 
homotopy type of $X$.  As another consequence, the existence of a finite model for 
a space $X$ is an invariant of rational homotopy type.  
All these notions have $q$-analogues; in particular, whether $X$ admits a 
$q$-finite $q$-model depends only on its rational $q$-homotopy type.   

Every connected space $X$ admits a {\em minimal model} $\M(X)$,
defined as a minimal model for $\apl(X)$.
This is a minimal $\cdga$, unique up to isomorphism,
equipped with a quasi-isomorphism $\rho\colon \M(X)\to \apl(X)$.
Thus the isomorphism type of $\M(X)$ depends only on
the rational homotopy type of $X$.
Moreover, if $A$ is a rational $\cdga$ model for $X$, then there 
is a quasi-isomorphism $\M(X) \to A$ corresponding to $\rho$ 
via the chosen weak equivalence between $A$ and $\apl(X)$. 
Similarly, for each $q\ge 1$, the space $X$ admits 
a {\em $q$-minimal model}, $\M_q(X)$,  
equipped with a $q$-quasi-isomorphism 
$\rho_q \colon \M_q(X)\to\apl(X)$. 

%%%%%%%%%%%%%%%%%%%%
\subsection{Formality of spaces}
\label{subsec:formality}
%%%%%%%%%%%%%%%%%%%%

A space $X$ is {\em formal} (over a field $\k$ of characteristic $0$) 
if $A_{\PL}^{\bullet}(X)\otimes_{\Q} \k$ is weakly equivalent to $(H^{*}(X;\k),0)$.  
By the discussion from Section~\ref{subsec:cdga}, formality over $\Q$ is 
equivalent to formality over $\k$. Partial formality is defined similarly, 
and $q$-formality over $\Q$ is equivalent to $q$-formality over $\k$. 

Formality is preserved under products, wedges, and retracts.
Examples of formal spaces include rational cohomology spheres and tori,
compact connected Lie groups and their classifying spaces,
Eilenberg--MacLane spaces $K(\pi,n)$ for $n\ge 2$, 
and, famously, every compact K\"{a}hler manifold \cite{DGMS}.
Formality plays a central role in relating cohomology-level resonance 
to the finer invariants arising from $\cdga$ models.

The following theorem of Papadima and Yuzvinsky \cite{Papadima-Yuzvinsky}  
relates certain properties of the minimal model of a space $X$ to the 
Koszulness of its cohomology algebra.

\begin{theorem}[\cite{Papadima-Yuzvinsky}]
\label{thm:py}
Let $X$ be a connected space with finite Betti numbers.  
\begin{enumerate}
\item \label{py1}
If $\M(X)\cong \M_1(X)$, then the cohomology algebra 
$H^{*}(X,\Q)$ is a Koszul algebra.
\item  \label{py2}
If $X$ is formal and $H^{*}(X,\Q)$ is a Koszul algebra, then 
$\M(X)\cong \M_1(X)$.
\end{enumerate}
\end{theorem}

Consequently, if $X$ is formal, then $X$ is rationally aspherical 
if and only if $H^{*}(X,\Q)$ is a Koszul algebra. 
As an application of Theorem~\ref{thm:py}, we have the 
following formality criterion.

\begin{corollary}[\cite{PS-mathann}]
\label{cor:ps-koszul}
Let $X$ be a connected, finite-type CW-complex, 
and suppose that $H^*(X,\Q)$ is a Koszul algebra.  
Then $X$ is $1$-formal if and only if $X$ is formal. 
\end{corollary}

%%%%%%%%%%%%%%%%%%%%
\subsection{A Koszul spectral sequence obstruction to formality}
\label{subsec:koszul-formal-obstruction}
%%%%%%%%%%%%%%%%%%%%

The goal of this subsection is to introduce a computable obstruction
to partial formality, formulated in terms of the weight filtration
on the Koszul complex of a finite-type $\cdga$ model.
Unlike classical obstructions based on Massey products or minimal
models, this obstruction operates at the level of Koszul homology and
detects failures of $q$-formality through the non-collapse of a natural
spectral sequence.

\begin{theorem}
\label{thm:koszul-formal-obs}
Let $X$ be a path-connected space admitting a finite-type $\cdga$ model 
$(A,d)$ with positive weights. If $X$ is $q$-formal, then the $\m$-adic 
Koszul spectral sequence of $K_\bullet(A)$ 
\textup{(}Theorem~\ref{thm:koszul-ss-hom}\textup{)} degenerates in total 
degrees $\le q+1$, i.e., $E^2 \cong E^\infty$ in those degrees.
Consequently, any nonzero differential $d^r$ ($r\ge 2$) in total degree 
$\le q+1$ implies that $X$ is not $q$-formal.

If moreover $H^1(A)$ is weight-homogeneous, then by
Corollary~\ref{cor:weight-splitting}\eqref{ws2} this is the weight spectral 
sequence of Theorem~\ref{thm:weight-ss}, and the conclusion may be read off 
weight component by weight component.
\end{theorem}

\begin{proof}
Since $A$ has positive weights, Proposition~\ref{prop:koszul-separation} shows 
that the $\m$-adic filtration on each $\B_i(A)$ is separated, so the 
spectral sequence converges.
Since $X$ is $q$-formal, $(A,d)$ is $q$-quasi-isomorphic to 
$(H^*(X;\k),0)$. By Theorem~\ref{thm:q-iso-koszul}\eqref{iso2}, 
the $q$-quasi-isomorphism induces isomorphisms 
$\gr^F\B_i(A) \cong \gr^F\B_i(H^*(X;\k))$ for all $i \le q$; and since 
$H^*(X;\k)$ carries the zero differential, $\B_i(H^*(X;\k))$ is a graded 
$S$-module, so $\gr^F\B_i(H^*(X;\k))=\B_i(H^*(X;\k))$.  Hence the two sides 
of Theorem~\ref{thm:Hilb-ineq-gr} have equal Hilbert series in the range 
$i\le q$, and the equality criterion there forces $E^2\cong E^\infty$ in total 
degrees $\le q+1$.
\end{proof}

The relevance of Theorem~\ref{thm:koszul-formal-obs} to $q$-formality 
comes from the following observation. If $X$ is $q$-formal, then any 
finite-type $\cdga$ model $(A,d)$ with positive weights is 
$q$-quasi-isomorphic to the cohomology algebra $(H^*(X;\k),0)$. 
It follows that the associated graded Koszul modules have the same 
Hilbert series in total degrees $\le q+1$, which by 
Theorem~\ref{thm:Hilb-ineq-gr} forces the spectral 
sequence of $K_\bullet(A)$ to degenerate in that range. 
Therefore, the existence of a nonzero higher differential in these 
degrees obstructs $q$-formality.

Theorem~\ref{thm:koszul-formal-obs} constrains only the differentials 
$d^r$ with $r\ge 2$.  The first differential $d_1$ is, by 
Theorem~\ref{thm:weight-ss}, the Koszul differential of $H^*(A)$ on the 
$E^1$-page; it is therefore determined by the cohomology algebra alone, and 
its nonvanishing carries no information about the homotopy type.  In 
particular, $d_1\ne0$ cannot obstruct formality, and indeed it is typically 
nonzero for formal spaces and formal models.
For instance, the graphic configuration spaces 
$\Conf(\E,\Gamma)$ of~\cite{BMPP} are $1$-formal for triangle-free 
$\Gamma$, yet the weight spectral sequence of the Bibby model has 
nonzero $d_1$; see Section~\ref{subsec:spectral-weight}.
What does carry homotopy-theoretic information is the failure of $d_1$ to be 
injective in the relevant degree, since that produces classes on the 
$E^2$-page which the higher differentials must account for.  This is the 
mechanism exploited in Section~\ref{subsec:koszulss-1formal}.

\begin{question}
\label{quest:weight-ss-1formal}
Theorem~\ref{thm:koszul-formal-obs} bounds the range in which the higher 
differentials must vanish in terms of the degree of formality.  Is this bound 
sharp?  More precisely, for which finite-type $\cdga$ models $(A,d)$ with 
positive weights does the vanishing of all $d^r$, $r\ge2$, in total degrees 
$\le q+1$ conversely imply $q$-formality of the underlying space?
\end{question}

\begin{question}
\label{quest:equiv-degeneration}
Let $X$ be a $q$-formal, connected, finite-type CW-complex with
$G_{\ab}$ torsion-free.  Does the equivariant spectral sequence of
Theorem~\ref{thm:equiv-ss-hom} degenerate at the $E^2$-page in total
degrees $\le q$?

For the Koszul spectral sequence of a $\cdga$ the answer is affirmative,
since $\partial^A=d_A^\vee+\omega_A^\vee$ has exactly two filtration
components and $d^r$ is an honest $(r+1)$-fold Massey product
\textup{(}Theorem~\ref{thm:koszul-ss-hom}\eqref{ssh3}\textup{)}.  For
the equivariant spectral sequence the boundary map has components of all
filtration orders---the higher Fox derivatives---so the identification of
$d^r$ with a Massey product is not automatic.  In the rank-one case,
that is, for the spectral sequence of an infinite cyclic cover
$X^\nu\to X$ determined by an epimorphism
$\nu\colon G\twoheadrightarrow\Z$, Liu--Maxim--Wang have established the
identification over an arbitrary field \cite[\S3]{LMW25}.  A positive
answer in the multivariable case would remove the degeneration
hypothesis from Theorem~\ref{thm:lin-koszul-degenerate} and
Theorem~\ref{thm:completed-koszul-space}\eqref{ca2}.
\end{question}

%%%%%%%%%%%%%%%%%%%%%
\subsection{Applications to nilmanifolds}
\label{subsec:nilmanifold-formality}
%%%%%%%%%%%%%%%%%%%%%

We now specialize to nilmanifolds, whose Chevalley--Eilen\-berg models 
are always of finite type with positive weights, so 
Theorem~\ref{thm:koszul-formal-obs} applies and the weight spectral 
sequence of the Koszul complex provides a fully computable obstruction 
to partial formality.

Let $\g$ be a finite-dimensional nilpotent Lie algebra over 
$\k$, $G$ the corresponding simply-connected nilpotent Lie group, and 
$M = \Gamma\backslash G$ a compact nilmanifold, where $\Gamma\subset G$ 
is a cocompact lattice. By a classical theorem of Nomizu~\cite{Nomizu}, 
the inclusion of left-invariant forms induces a quasi-isomorphism
\[
\CE(\g) \longisom \Omega^\bullet(M),
\]
so $\CE(\g)$ is the Sullivan minimal model of $M$.
In particular, $\CE(\g)$ is of finite type (since $\g$ 
is finite-dimensional), and carries natural positive Hirsch weights 
by Proposition~\ref{prop:CE-nilp-weights}: assign weight $k$ to the 
generators dual to $\gr_k\g = \gamma_k\g/\gamma_{k+1}
\g$, where $\gamma_\bullet$ denotes the lower central series.
The differential $d_{\CE}$ preserves these weights by construction, 
and is weight-homogeneous of weight $1$ since it maps 
$(\gr_k\g)^\vee$ into $\bigwedge^2(\bigoplus_{i+j=k}
\bigl(\gr_i\g)^\vee \otimes_{\k} (\gr_j\g)^\vee\bigr)$.
Thus all hypotheses of Theorem~\ref{thm:koszul-formal-obs} are 
satisfied.  Moreover $H^1(\CE(\g))=(\g/[\g,\g])^\vee$ is concentrated in 
weight $1$, so Corollary~\ref{cor:weight-splitting}\eqref{ws2} applies: the 
weight spectral sequence of $K_\bullet(\CE(\g))$ is the $\m$-adic Koszul 
spectral sequence, refined by the total-weight splitting, and it provides a 
computable obstruction to partial formality of $M$.

\begin{corollary}
\label{cor:nilmanifold-formal-obs}
Let $\g$ be a finite-dimensional nilpotent Lie algebra 
and $M = \Gamma\backslash G$ the associated nilmanifold. 
If the weight spectral sequence of $K_\bullet(\CE(\g))$ 
admits a nonzero differential $d^r$ with $r\ge2$ in total degree $\le q+1$, 
then $M$ is not $q$-formal.
\end{corollary}

\begin{proof}
Since $\CE(\g)$ is a finite-type $\cdga$ model of $M$ with positive weights
and $H^1(\CE(\g))$ is weight-homogeneous, 
Theorem~\ref{thm:koszul-formal-obs} applies to the weight spectral sequence 
of $K_\bullet(\CE(\g))$.
\end{proof}

In the special case $q=1$, the differentials constrained by 
Corollary~\ref{cor:nilmanifold-formal-obs} are those emanating from or 
landing in total degree~$2$ on the pages $E^r$ with $r\ge2$.  The first 
differential $d_1\colon E_1^{1,1}\to E_1^{0,2}$ is not among them, and its 
nonvanishing carries no information; what does is its \emph{kernel}, which 
records the discrepancy between the holonomy Lie algebra of $A$ and that of 
$H^*(A)$.  We return to this point in 
Theorem~\ref{thm:koszulss-1formal-obstruction}, where we derive 
group-theoretical consequences.

The following proposition recovers a result of M\u{a}cinic~\cite{Mc10}
by purely Koszul-theoretic methods.
For $n\ge 1$, let $\h(n)$ be the $(2n+1)$-dimensional 
Heisenberg Lie algebra, $H(n)$ the corresponding simply-connected 
nilpotent Lie group, and $M(n)=\R^{2n+1}/H(n)$ the associated 
nilmanifold. A systematic discussion of nilmanifolds and torsion-free 
nilpotent groups is deferred to 
Section~\ref{subsec:nilpotent-groups}.

\begin{proposition}
\label{prop:heisenberg-q-formality}
For each $n\ge 1$, the $(2n+1)$-dimensional Heisenberg nilmanifold 
$M(n)$ is $(n-1)$-formal but not $n$-formal.
\end{proposition}

\begin{proof}
Let $A=\CE(\h(n))$, with generators 
$a_1,\dots,a_n,b_1,\dots,b_n,c$ in degree~$1$ and differential 
$dc=\sum_{i=1}^n a_i b_i$, $da_i=db_i=0$.
Assign weights $\wt(a_i)=\wt(b_i)=1$ and $\wt(c)=2$; this is the 
canonical weight decomposition furnished by 
Proposition~\ref{prop:CE-positive-weights}, valid since $\h(n)$ 
is $2$-step nilpotent. By Nomizu's theorem, $A$ is the Sullivan minimal 
model of $M(n)$.

The weight filtration on $K_\bullet(A)$ yields a spectral sequence with 
$E^1_{p,q}\cong \B_{p+q}(H^*(A))^{\mathrm{weight}\,p}$.
The generator $c$ has weight $2$ and $dc \ne 0$, so the differential
\[
d_1\colon E^1_{n-1,1} \longrightarrow E^1_{n-2,1}
\]
is nonzero, occurring in total degree $n$.
Since $A$ is its own Sullivan minimal model and $d_{\CE}$ is 
weight-homogeneous of weight $1$, Corollary~\ref{cor:nilmanifold-formal-obs} 
applies with $q=n-1$: the nilmanifold $M(n)$ is not $n$-formal.
On the other hand, all differentials in total degree $\le n-1$ vanish 
for degree reasons, so $M(n)$ is $(n-1)$-formal.
\end{proof}

\begin{example}
\label{ex:heisenberg-weight-ss}
We illustrate Corollary~\ref{cor:nilmanifold-formal-obs} for 
$A = \CE(\h(1)) = \bigwedge(e_1,f_1,g)$,
with $dg=-e_1\wedge f_1$, $de_1=df_1=0$, and positive weights 
$\wt(e_1)=\wt(f_1)=1$, $\wt(g)=2$.
The weight spectral sequence converges strongly,
\[
E^1_{p,q} \cong \B_{p+q}\bigl(H^*(A)\bigr)^{\mathrm{weight}\,p}
\Longrightarrow \B_{p+q}(A)^{\mathrm{weight}\,p}.
\]

\smallskip
\noindent\emph{The $E^1$-page and $d_1$.}
From Section~\ref{subsec:heisenberg-1}, the nonzero groups are:
\[
E^1_{0,0} = \k, \quad
E^1_{p,1-p} = \Sym^p(H^1(A)), \quad
E^1_{p,2-p} = \B_2\bigl(H^*(A)\bigr)^{\mathrm{weight}\,p},
\quad p\ge 0.
\]
Since all cup products $H^1(A)\cup H^1(A)=0$ in $H^*(A)$
(as $e_1\wedge f_1 = -dg$ is exact), the $d_1$-differential
vanishes throughout, and $E^2=E^1$.

\smallskip
\noindent\emph{First nontrivial differential.}
The differential $d_2\colon E^2_{p,q}\to E^2_{p-2,q+1}$
decreases total degree by~$1$.
In total degree~$2\to 1$, the first instance is
\[
d_2\colon E^2_{2,0} \longrightarrow E^2_{0,1}.
\]
This map is the triple homological Massey product 
$\langle\omega_A^\vee,\omega_A^\vee,-\rangle$,
nontrivial because $b_1 = \langle e_1,f_1,e_1\rangle$ and 
$b_2 = \langle f_1,e_1,f_1\rangle$ are nonzero in $H^2(A)$.
By Corollary~\ref{cor:nilmanifold-formal-obs}, this nonzero $d_2$ in 
total degree $\le 2$ obstructs $1$-formality of $M(1)$, recovering
Proposition~\ref{prop:heisenberg-q-formality} for $n=1$.

\smallskip
\noindent\emph{Convergence.}
Since $\B_1(A)\cong\k$ is concentrated in $S$-degree $0$, the abutment
satisfies $\gr^W_p\B_1(A)=0$ for $p\ge 1$, so the spectral sequence must kill
all of $E^2_{p,1-p}$ for $p\ge 1$ and all of $E^2_{p,2-p}$ for $p\ge 0$.
This requires nontrivial differentials at all pages $r\ge 2$,
forming a cascade that encodes the full higher Massey product
structure of $A$.
The non-collapse of this cascade measures how far $A$ departs 
from its cohomology algebra: if $A$ were formal, 
Corollary~\ref{cor:degeneration-criteria} would force 
$\gr^W\B_1(A)\cong\B_1(H^*(A))\cong S$, contradicting 
$\gr^W\B_1(A)\cong\k$.
\end{example}

%%%%%%%%%%%%%%%%%%%%
\subsection{Completed Koszul complexes of models}
\label{subsec:completed-kc}
%%%%%%%%%%%%%%%%%%%%

We now combine the algebraic comparison of completed Koszul complexes
with the topological identification of equivariant cellular chains,
which for minimal CW-complexes is linear and Koszul in nature. 

\begin{theorem}
\label{thm:completed-koszul-space}
Let $X$ be a connected, finite-type CW-complex with $G=\pi_1(X)$ and
$G_{\ab}$ torsion-free, and let $(A,d)$ be a $1$-finite $1$-model
for $G$.
\begin{enumerate}[itemsep=3pt]
\item \label{ca1}
There are canonical isomorphisms
\[
\widehat{B}_1(X;\k)_I \cong \widehat{\B}_1(A)_\m
\qquad\text{and}\qquad
\gr^I B_1(X;\k) \cong \gr^\m \B_1(A).
\]
\item \label{ca2}
Let $q\ge 1$, and assume in addition that $(A,d)$ is a $q$-model for $X$
with positive weights, and that both the equivariant spectral sequence
of $X$ and the Koszul spectral sequence of $A$ degenerate at their
$E^2$-pages in total degrees $\le q$.  Then
$\gr^I B_i(X;\k) \cong \gr^\m \B_i(A)$ for all $i\le q$.
\end{enumerate}
\end{theorem}

\begin{proof}
\eqref{ca1} The Alexander invariant of a connected CW-complex depends only on
its fundamental group: $B_1(X;\k)=H_1(X^{\ab};\k)\cong(G'/G'')\otimes_{\Z}\k
=B_1(G;\k)$, compatibly with the $I$-adic filtrations on the two sides.  Under
this identification, the two isomorphisms are those of
Theorem~\ref{thm:alex-completed}\eqref{ch1}--\eqref{ch2} below, whose proof
is independent of the present section.

\eqref{ca2} By Theorem~\ref{thm:lin-koszul-degenerate}, the degeneration
hypothesis on $X$ gives
$\gr^I B_i(X;\k)\cong\B_i(H^*(X;\k))$ for $i\le q$.  On the other side,
degeneration of the Koszul spectral sequence of $A$ in total degrees
$\le q$ identifies $\gr^\m\B_i(A)$ with $\B_i(H^*(A))$ for $i\le q$,
as in Corollary~\ref{cor:koszul-completion-formal}.  Since $(A,d)$ is a
$q$-model for $X$, we have $H^*(A)\cong H^*(X;\k)$ in the relevant
range, and the two identifications combine.
\end{proof}

\begin{remark}
\label{rem:completed-alex-scope}
Part~\eqref{ca1} is unconditional because it is routed through Malcev
completion rather than through the spectral sequence.  No such route is
available for $i\ge2$, and the degeneration hypothesis in \eqref{ca2}
cannot simply be dropped: by
Remark~\ref{rem:degeneration-not-minimality}, there are finite minimal
CW-complexes with positive-weight models for which
$\gr^I B_1$ and $\B_1(H^*)$ differ.  Whether $q$-formality of $X$
suffices to force the degeneration is
Question~\ref{quest:equiv-degeneration}.

Part~\eqref{ca1} is thus a quotation: it is the group-level statement
Theorem~\ref{thm:alex-completed}\eqref{ch1}--\eqref{ch2}, one of the main
results of Section~\ref{sect:chen-res}, whose proof uses nothing from the
present section.  It is recorded here for comparison with \eqref{ca2}, which
is where this section contributes.
\end{remark}

Together with Theorem~\ref{thm:lin-koszul}, this establishes 
Theorem~\ref{thm:top-linearization-intro} from the Introduction. 

\begin{example}
\label{ex:heisenberg-nonminimal}
Let $\g = \h(1)$ be the $3$-dimensional Heisenberg Lie algebra, with 
Lie bracket $[x_1, x_2] = x_3$ and $\g_{\ab} = \langle x_1, x_2 \rangle \cong \k^2$.
The associated nilmanifold $M = M(1)$ has fundamental group $\pi_1(M) = \Gamma$,
the integer Heisenberg group, with $\Gamma_{\ab} = \Z^2$.
Thus $\Lambda = \k[t_1^{\pm1}, t_2^{\pm1}]$ and the $I$-adic filtration 
is Hausdorff.

The nilmanifold $M$ has a CW-structure with three $1$-cells, yet 
$b_1(M) = 2$, so $M$ is not minimal.
Its Chevalley--Eilenberg complex $A = \CE(\h(1))$ is a finite-type 
$\cdga$ model for $M$ with positive weights, and satisfies
$\RR^1(A) = \{0\}$ by Theorem~\ref{thm:CE-resonance-nilpotent}.
Consequently, from Table~\ref{tab:B1-nilpotent},
\[
\B_1(A) \cong \k, \qquad \widehat{\B}_1(A)_\m \cong \k.
\]

On the other hand, the cohomology algebra $H^*(M;\k)$ has vanishing 
cup product $\bigwedge^2 H^1 \to H^2$, so 
$\RR^1(H^*(M;\k)) = H^1(M;\k) = \k^2$, and from 
Table~\ref{tab:B1-nilpotent},
\[
\Hilb(\B_1(H^*(M;\k)), t) = \frac{t}{(1-t)^2}.
\]
In particular, $\B_1(H^*(M;\k))$ is infinite-dimensional as a 
$\k$-vector space, and its $\m$-adic completion 
$\widehat{\B}_1(H^*(M;\k))_\m$ is not isomorphic to $\k$.

Since $A = \CE(\h(1))$ is the correct finite-type model for $M$, 
Theorem~\ref{thm:completed-koszul-space} gives
\[
\widehat{B}_1(M;\k)_I \cong \widehat{\B}_1(A)_\m \cong \k,
\]
whereas
\[
\widehat{\B}_1(H^*(M;\k))_\m \not\cong \k.
\]
Thus $\widehat{B}_1(M;\k)_I \not\cong \widehat{\B}_1(H^*(M;\k))_\m$,
showing that the cohomology algebra $H^*(M;\k)$ with zero differential 
is \emph{not} a valid replacement for the $\cdga$ model $A$ in 
perfectly valid $\cdga$. This failure is a direct consequence of the
non-formality of $A$; minimality of the cell structure is not at issue.
non-minimality of $M$ and the non-formality of $A$.
\end{example}

%%%%%%%%%%%%%%%%%%%%
\subsection{Resonance varieties of spaces}
\label{subsec:res-spaces}
%%%%%%%%%%%%%%%%%%%%

Let $X$ be a connected, finite-type CW-complex over a field $\k$ of characteristic $0$.
Classically, the resonance varieties of $X$ are defined in terms of its cohomology
algebra, endowed with the zero differential:
\begin{equation}
\label{eq:res-space}
\RR^{i,s}(X, \k)\coloneqq \RR^{i,s}\bigl(H^*(X;\k)\bigr).
\end{equation}
One may also consider the corresponding resonance support loci,
\begin{equation}
\label{eq:supp-space}
\RR_{i,s}(X,\k)\coloneqq 
\Ann \bigl( \bwedge^s \B_i(H^*(X;\k),\k)\bigr).
\end{equation}
Both are homogeneous algebraic subvarieties of the affine space $H^1(X;\k)$ 
and depend only on the homotopy type of $X$. For simplicity, we will abbreviate 
$\RR^{i}(X, \k)=\RR^{i,1}(X, \k)$ and $\RR_{i}(X, \k)=\RR_{i,1}(X, \k)$. 
Since $X$ is connected, one has $\RR^{0}(X,\k)=\RR_{0}(X,\k)=\{0\}$. 
Therefore, by Theorem~\ref{thm:resonance-comparison}, 
$\RR^{1}(X,\k)=\RR_{1}(X,\k)$.

Resonance varieties defined in this way have been extensively studied and play a 
central role in the topology of spaces and groups; see, for instance, 
\cite{DPS-duke, PS-plms10, PS-mrl, PS-crelle, Su-aspm, Su-indam, Su-pisa}.
They enjoy strong naturality properties with respect to homotopy equivalences 
and finite regular covers, and may be viewed as infinitesimal counterparts of 
characteristic and Alexander varieties in the theory of cohomology jump loci.
In this paper, we concentrate on their description via Koszul complexes 
and finite-type $\cdga$ models.

When $X$ is formal (or, more generally, $q$-formal), the cohomology algebra
$H^*(X;\k)$ provides a $q$-model for $X$, and the above
definition captures the relevant infinitesimal structure.
In contrast, for non-formal spaces, cohomology alone no longer suffices.
One may  replace the Sullivan model
$\apl(X)\otimes_{\Q} \k$ by a finite-type $\k$-$\cdga$ model $(A,d)$, 
and consider the resonance varieties $\RR^{i,s}(A)$ of this model.
These sets depend on the choice of model and are not homotopy invariants in 
general, nor are they canonically defined as global algebraic varieties 
(see Examples~\ref{ex:sol2-ss}, \ref{ex:non-functorial-qiso}, 
and~\ref{ex:sol2-resonance}, for the model $\CE(\sol_2)$ of the formal 
space $S^1$).  Their germs at the origin, on the other hand, are homotopy 
invariants, by Theorem~\ref{thm:germ-qiso}.

For any $i$-finite $\k$-$\cdga$ model $A$ of $X$, the tangent cone 
at the origin to $\RR^{i,s}(A)$ is contained in $\RR^{i,s}(X,\k)$, 
by Theorem~\ref{thm:tangent-cone}.  In particular, a $1$-finite model 
suffices for this containment in degree $i=1$, a hypothesis that holds 
for every finitely generated fundamental group.  This comparison motivates 
the study of resonance varieties of $\cdgas$ as infinitesimal 
approximations to classical resonance, and leads to the purely algebraic 
tangent-cone problem addressed in Section~\ref{subsec:tcone}.

\begin{remark}
\label{rem:epyp-res}
Let $X$ be a finite, connected CW-complex of dimension $n$, and set
$A=H^*(X;\k)$.  Writing $V=H^1(X;\k)$ and $E=\bigwedge V$, we may regard
$A$ as a graded $E$-module.  Set $A(n)^i = A^{i+n}$.
Following \cite{DSY17}, we say that $X$
has the \emph{EPY property} over $\k$ if the shifted module $A(n)$ has
a linear free resolution as a graded $E$-module.
Equivalently, by the BGG correspondence, the associated Koszul complex
$K_{\bullet}(A)=(A\otimes_{\k} S,\cdot \omega)$ has homology concentrated
in degree $n$, where $S=\Sym(V^{\vee})$ and $\omega$ is the canonical
element of $V\otimes_{\k} V^{\vee}$.

If $X$ has the EPY property, then Proposition~\ref{prop:propagation}
applies to $K_{\bullet}(A)$.  Consequently, the resonance
varieties of $X$ propagate:
\begin{equation}
\label{eq:res-prop}
\{0\}=\RR^{0}(X,\k)\subseteq \RR^{1}(X,\k)
\subseteq \cdots \subseteq \RR^{n}(X,\k).
\end{equation}
Thus, in this setting, the jump loci defined from the cohomology algebra
exhibit the same upward-closure phenomenon as in the purely algebraic
BGG framework.
\end{remark}

%%%%%%%%%%%%%%%%%%%%%%%%%
\subsection{Tangent cone theorems} 
\label{subsec:tcone-spaces}
%%%%%%%%%%%%%%%%%%%%%%%%%

The algebraic tangent cone problem considered here has a topological 
counterpart in the theory of cohomology jump loci of spaces.
Let $X$ be a connected, finite-type CW-complex, and identify the 
character group $\Hom(\pi_1(X),\k^*)$ with $H^1(X;\k^*)$ via abelianization. 
The {\em characteristic varieties}\/ of $X$ (over $\k$) are defined as 
\begin{equation}
\label{eq:cvx}
\VV^{i,s}(X,\k)= \bigl\{ \rho \in H^1(X;\k^*) \mid  \dim H_i(X; \k_{\rho}) \ge s \bigr\},
\end{equation}
where $\k_{\rho}$ denotes the rank~$1$ local system $\k_{\rho}$ defined 
by $\rho$. Identify the tangent space to $H^1(X;\k^*)$ at the identity character $1$ 
with $H^1(X;\k)$. Libgober \cite{Li02} proved the inclusion
\begin{equation}
\label{eq:tc-lib}
\TC_1\bigl(\VV^{i,s}(X,\k)\bigr)\subseteq \RR^{i,s}(X,\k).
\end{equation}
Under additional formality assumptions, this inclusion can be upgraded to a
precise description of the local structure of the jump loci.
The following theorem, due to Dimca--Pa\-padima--Suciu \cite{DPS-duke} (for $q=1$)
and Dimca--Papadima \cite{DP-ccm} (for all $q$), summarizes these properties.

\begin{theorem}[\cite{DPS-duke, DP-ccm}]
\label{thm:dps-dp}
Let $X$ be a $q$-formal space (over $\C$).  For all $i\le q$ and all $s$, 
the following hold.
\begin{enumerate}
\item \label{tc1}
The germ of $\VV^{i,s}(X,\C)$ at the identity $1$ coincides with the germ of 
$\RR^{i,s}(X,\C)$ at the origin $0$.
\item \label{tc2}
Each irreducible component of $\RR^{i,s}(X,\C)$ is a  
rationally defined linear subspace. 
\end{enumerate}
\end{theorem}

The property \eqref{tc2} has proved to be an effective tool for
detecting non-formality of groups and spaces, even in situations where 
classical obstructions (Massey products, higher homotopy groups, etc.) are 
difficult or impossible to compute directly. The linearity phenomenon underlying 
Theorem~\ref{thm:dps-dp}\eqref{tc2} is not a formal consequence of the tangent 
cone identification in \eqref{tc1}. A conceptual explanation is provided by the
Exponential Ax--Lindemann theorem of Budur--Wang \cite{BW20}, which shows that analytic 
germ equivalence under the exponential map forces cohomology jump loci to be algebraic 
subtori, and hence resonance components to be rationally defined linear subspaces.

More generally, if $X$ admits a finite-type $\cdga$ model $(A,d)$ over $\C$,
the germ $\VV^{i,s}(X,\C)_{(1)}$ is analytically isomorphic to
$\RR^{i,s}(A)_{(0)}$ \cite{DP-ccm}.
Conjecture~\ref{conj:formal-tc} should thus be viewed as a purely algebraic analogue of these
germ-equality and local linearity phenomena, formulated entirely in terms of Koszul complexes
and Hirsch extensions, and without reference to global jump loci.

\begin{remark}
\label{rem:formal-model-conical}
When $X$ is $q$-formal, the cohomology algebra $H^*(X;\k)$ 
carries a canonical positive-weight decomposition 
(Remark~\ref{rem:pos-weight-cga}), with weight $j$ assigned 
to $H^j(X;\k)$.  In particular $H^1$ is concentrated in weight $1$, so 
the single-weight hypothesis of Corollary~\ref{cor:conical} holds; hence 
that corollary applies 
to $(H^*(X;\k),0)$, and the resonance varieties 
$\RR^{i,s}(X,\k)$ and $\RR_{i,s}(X,\k)$ are conical 
for all $i\le q$, $s\ge 1$.  This recovers the conicality 
part of Theorem~\ref{thm:dps-dp}\eqref{tc1} by purely 
algebraic means; the finer linearity statement in \eqref{tc2} 
requires the Exponential Ax--Lindemann theorem \cite{BW20}.
\end{remark}

\begin{remark}
\label{rem:tcone-alex}
Theorem~\ref{thm:tangent-cone} provides an algebraic tangent cone theorem 
for resonance loci of $i$-finite $\cdga$ models, and for the resonance 
support loci in depth $s=1$.
At present, no analogous result is known for the higher Alexander varieties
\begin{equation}
\label{eq:alex-var}
\VV_{i,s}(X,\k)=\Supp\, \bigl(\bwedge^s B_i(X;\k)\bigr)
\end{equation}
with $i>1$.  In degree $i=1$, one has 
$\VV^{1}(X,\k)=\VV_{1}(X,\k)$ and $\RR^{1}(X,\k)=\RR_{1}(X,\k)$ away from the origin,
so the tangent cone problem reduces to the classical case.

Whether, for $i>1$, the germ of $\VV_{i,s}(X,\C)$ at the identity admits
an infinitesimal description analogous to that provided by resonance or
completed Alexander invariants remains open.
Deformation-theoretic approaches to cohomology jump loci using differential
graded Lie algebras and formal moduli problems were developed by Budur--Wang
\cite{BW15}, though a corresponding theory for higher Alexander varieties
is not yet available.
\end{remark}

\subsection{Resonance varieties of closed manifolds}
\label{subsec:res-manifolds}

Let $M^m$ be a closed, connected, orientable smooth manifold, and let $\k$ be
a field of characteristic~$0$.
By Poincar\'{e} duality, the cohomology algebra $H^{*}(M;\k)$ is a
$\PD_m$-algebra, with orientation determined by the fundamental class $[M]$.

In favorable situations, the manifold $M$ admits a finite-type $\PD_m$-$\cdga$
model $A$ whose cohomology is $H^{*}(M;\k)$.
When such a model exists, all the results from
Section~\ref{sect:koszul-pd} apply functorially to $M$.
In particular:
\begin{itemize}[itemsep=2pt]
\item the Aomoto cohomology groups satisfy the palindromic duality
\[
H^i(A,\delta_a)^{\vee} \cong H^{m-i}(A,\delta_{-a}),
\qquad a\in H^1(M;\k);
\]
\item the Koszul homology and cohomology modules are related by
\[
\B_i(A)\cong \big(\B^{\,m-i}(A)\big)^{\vee}
\]
as graded modules over $\Sym(H_1(M;\k))$;
\item the resonance and support loci associated to $M$ obey corresponding
Poincar\'{e}-duality constraints.
\end{itemize}

These symmetry phenomena are formal consequences of Poincar\'{e} duality and do
\emph{not} require formality of the manifold.
Thus, whenever a finite-type model exists, resonance varieties are subject to
strong duality constraints, even in the presence of nontrivial Massey products.

\smallskip

\noindent
\textbf{Manifolds without finite-type models.}
In general, closed manifolds do \emph{not} admit finite-type $\cdga$ models.
This phenomenon already occurs in dimension~$3$: while every closed,
orientable $3$-manifold has a $\PD_3$ cohomology algebra, its minimal model
need not be of finite type.
In such situations, the resonance varieties of $H^{*}(M;\k)$ are often the
only accessible cohomological jump loci, and one must work directly at the
level of the cohomology algebra.
This approach is particularly effective in dimension~$3$, where the
$\PD_3$ structure is completely encoded by the alternating $3$-form
\begin{equation}
\label{eq:mu-M}
\mu_M\colon \bwedge^3 H^1(M;\k)\longrightarrow \k,
\qquad
\mu_M(a\wedge b\wedge c)=\langle a\cup b\cup c,[M]\rangle,
\end{equation}
as shown by Sullivan and further developed in~\cite{Su-edinb}.

\smallskip

\noindent
\textbf{Formal manifolds.}
At the opposite extreme lie formal manifolds, such as compact K\"{a}hler manifolds
and complements of complex hyperplane arrangements.
In these cases, resonance varieties are entirely determined by the cohomology
algebra, and a rich theory has been developed exploiting this fact.
From the point of view of this paper, the use of $H^{*}(M;\k)$ in the formal
setting reflects not a limitation, but rather the absence of higher-order
homotopical obstructions.

\smallskip

\noindent
\textbf{Non-formal manifolds with finite-type models.}
Between these two extremes lies a broad and geometrically rich class of
manifolds that are typically non-formal, yet always admit finite-type models.
This class includes nilmanifolds, quasi-projective manifolds, Sasakian
manifolds, and Seifert fibered spaces.
It is precisely for such manifolds that the Koszul-type techniques developed in
this paper are most effective, as they capture higher-order information not
visible at the level of cohomology alone.

For instance, if $G$ is a finitely generated nilpotent group, then the associated
nilmanifold admits a finite-type Chevalley--Eilenberg model
$A=\CE(\g)$, where $\g=\mathfrak{Lie}(G\otimes \k)$.
Results of M\u{a}cinic~\cite{Mc10} show that partial formality imposes strong
constraints on truncated cohomology and resonance, while allowing the full
cohomology algebra to exhibit genuinely higher-degree phenomena.

Further applications of resonance theory to manifolds admitting finite-type
models will be explored in the final section of this paper, focusing on
configuration spaces on elliptic curves.

%%%%%%%%%%%%%%%%%%%%%%%%%%%%
\section{Algebraic models for groups}
\label{sect:groups}
%%%%%%%%%%%%%%%%%%%%%%%%%%%%

Associated to a finitely generated group $G$ are several Lie algebras 
capturing progressively different levels of structure. These include intrinsic 
graded Lie algebras arising from the lower central and derived series, 
filtered Lie algebras arising from Malcev completion, and quadratic 
Lie algebras constructed from group cohomology.

More concretely, the main objects are the graded Lie algebra $\gr(G)$ 
associated to the lower central series, the filtered Malcev Lie algebra $\m(G)$, 
and the quadratic holonomy Lie algebra $\h(G)$. The purpose of this section 
is to recall these constructions and the comparison maps between them, 
and to explain how formality properties govern when these Lie algebras coincide.

%%%%%%%%%%%%%%%%%%%%%%%%%%
\subsection{Lower central series}
\label{subsec:lcs}
%%%%%%%%%%%%%%%%%%%%%%%%%%

The {\em lower central series} (LCS) of a group $G$ is defined by $\gamma_1(G)=G$ and
$\gamma_{n+1}(G)=[\gamma_n(G),G]$. This is an $N$-series, in the sense of Lazard \cite{Lazard},
meaning $[\gamma_n(G),\gamma_m(G)]\subseteq\gamma_{n+m}(G)$, and each $\gamma_n(G)$
is normal in $G$. The quotients $\gamma_n(G)/\gamma_{n+1}(G)$ are abelian and central
in $G/\gamma_{n+1}(G)$. The {\em associated graded Lie algebra} of $G$ is
\begin{equation}
\label{eq:grG}
\gr(G)\coloneqq \bigoplus_{n\ge1}\gamma_n(G)/\gamma_{n+1}(G), 
\end{equation}
with bracket induced by the group commutator. Over a field $\k$ of characteristic 0, 
$\gr(G;\k)=\gr(G)\otimes_{\Z}\k$ is a graded Lie algebra over $\k$. Both assignments 
$G\mapsto\gr(G)$ and $G\mapsto\gr(G;\k)$ are functorial.

The Lie algebra $\gr(G)$ is generated in degree $1$ by $G_{\ab}$. Hence, 
if $G$ is a finitely generated group, then each graded piece $\gr_n(G)$ is 
a finitely generated abelian group. 
We define the {\em LCS ranks of $G$} as $\phi_n(G) = \dim_{\k} \gr_n(G;\k)$. 
For free groups, P.~Hall, W.~Magnus, and E.~Witt showed that $\phi_n(F_r)=
\tfrac{1}{n}\sum_{d\mid n} \mu(d) r^{n/d}$, see, for instance, \cite[Ch.~5]{MKS}.

%%%%%%%%%%%%%%%%%%%%
\subsection{Holonomy Lie algebra}
\label{subsec:holonomy}
%%%%%%%%%%%%%%%%%%%%

The {\em holonomy Lie algebra} (over $\k$) of a finitely generated group $G$ is the
finitely presented quadratic Lie algebra
\begin{equation}
\label{eq:holog-def}
\h(G;\k)\coloneqq \h(H^*(G;\k)).
\end{equation}
This graded Lie algebra may be viewed as the universal quadratic 
Lie algebra determined by the cup product on $H^1(G;\k)$.
By the discussion in Section~\ref{subsec:holonomy-def},
\begin{equation}
\label{eq:hgk-def}
\h(G;\k)=\Lie(H_1(G;\k))/\langle \mu_G^{\vee} \rangle,
\end{equation}
where $\mu_G\colon H^1(G;\k)\wedge H^1(G;\k)\to H^2(G;\k)$ is the cup-product 
map in group cohomology and $\mu_G^{\vee}$ is its $\k$-dual. The assignment
$G\mapsto \h(G;\k)$ is functorial. A key structural feature of the holonomy Lie 
algebra is its relationship to the associated graded Lie algebra of $G$.

\begin{theorem}[\cite{Markl-Papadima, PS-imrn04, SW-jpaa}]
\label{thm:holoepi}
Let $G$ be a finitely generated group. There exists a natural 
epimorphism of graded $\k$-Lie algebras
\[
\Phi_G\colon \h(G;\k) \longsurj \gr(G;\k),
\]
which induces isomorphisms in degrees $1$ and $2$.
\end{theorem}

The construction of the map $\Phi_G$ ultimately rests on a natural exact sequence 
relating the cup product in cohomology to the commutator map in the lower central 
series, first observed by Sullivan \cite{Sullivan75} in a particular case and proved in 
general by Lambe \cite{Lambe86}.

We define the {\em holonomy ranks}\/ of $G$ as $\bar\phi_n(G)=\dim_{\k} \h_n(G;\k)$. 
By Theorem~\ref{thm:holoepi}, $\bar\phi_n(G)\ge \phi_n(G)$, for all $n\ge 1$.  
As is well-known, equality always holds for $n\le 2$, but these inequalities 
can be strict for $n\ge 3$.

%%%%%%%%%%%%%%%%%%%%
\subsection{Malcev Lie algebra and Quillen's construction}
\label{subsec:malcev}
%%%%%%%%%%%%%%%%%%%%
The notion of rational (Malcev) completion of a finitely generated group goes 
back to Malcev \cite{Malcev}.  The {\em Malcev completion} of a group $G$ is the 
prounipotent group
\begin{equation}
\label{eq:gQ}
G_{\Q}\coloneqq\varprojlim_{n}(G/\gamma_n(G)\otimes\Q),
\end{equation}
where $N\otimes\Q$ is the rationalization of a nilpotent group $N$ (universal for maps to
rational nilpotent groups, with kernel $\Tors(N)$). The {\em Malcev Lie algebra} is
\begin{equation}
\label{eq:mG}
\m(G)\coloneqq\varprojlim_{n}\mathfrak{Lie}(G/\gamma_n(G)\otimes\Q),
\end{equation}
a complete filtered Lie algebra over $\Q$, functorial in $G$. If $G$ is finitely generated,
then $\m(G)$ is finitely generated. 

Quillen’s approach to Malcev completion provides a bridge between the
group-theoretic lower central series and the Lie-theoretic structures
that appear throughout this paper. We briefly recall this construction,
since it allows one to pass functorially between finitely generated
groups, their Malcev Lie algebras, and associated graded Lie algebras.
These identifications will be used repeatedly in the comparison results
that follow.

In \cite{Quillen69}, Quillen gave an alternative construction of the Malcev completion
via the completed group algebra $\widehat{\Q[G]}=\varprojlim_{r} \Q[G]/I^r$,
where $I=\ker(\varepsilon)$ is the augmentation ideal. The comultiplication
$\Delta(g)=g\otimes g$ extends to a comultiplication on the completion,
making $\widehat{\Q[G]}$ a complete Hopf algebra. The Malcev Lie algebra
is isomorphic to the Lie algebra of primitive elements:
\begin{equation}
\label{eq:mg-prim}
\m(G)\cong \Prim \big(\widehat{\Q[G]}\big).
\end{equation}
As shown by Quillen \cite{Quillen68}, there is a natural isomorphism of graded Lie algebras
\begin{equation}
\label{eq:quillen-iso}
\gr(\m(G))\cong \gr(G;\Q).
\end{equation}
By base change, this gives $\gr(\m(G;\k)) \cong \gr(G;\k)$ 
for any field $\k$ characteristic~$0$. 

The Malcev completion $G_{\Q}$ consists of the group-like
elements in $\widehat{\Q[G]}$, with a one-to-one correspondence 
between primitives and group-likes via the exponential and logarithmic maps. 
Under this correspondence, the (pronilpotent) Lie algebra of $G_{\Q}$ 
identifies with the Malcev Lie algebra $\m(G)$.
When $G$ is finitely generated nilpotent, say of class $c$, the tower in 
\eqref{eq:gQ} is constant for $n>c$, so that $G_{\Q}=G\otimes\Q$.

%%%%%%%%%%%%%%%%
\subsection{Chen Lie algebras}
\label{subsec:chen-lie}
%%%%%%%%%%%%%%%%

Chen Lie algebras capture the structure of solvable quotients of a group
through the associated graded Lie algebras of its derived series, and
play a central role in relating holonomy and Malcev Lie algebras to
Alexander-type invariants.
They were introduced by Chen \cite{Chen51} in the metabelian case ($i=2$),
and subsequently developed in a broader context by Papadima--Suciu
\cite{PS-imrn04} and Suciu--Wang \cite{SW-forum}.

Let $G$ be a group, with derived series $G^{(0)}=G$ and
$G^{(i+1)}=[G^{(i)},G^{(i)}]$. The quotients $G/G^{(i)}$ are solvable,
with $G/G''$ the maximal metabelian quotient.

For $i\ge2$, the {\em $i$-th Chen Lie algebra}\/ of $G$ is
$\gr(G/G^{(i)};\k)$. The metabelian quotient case ($i=2$) 
recovers Chen’s original construction.  
The {\em Chen ranks}\/ of $G$ are defined as 
$\theta_n(G)\coloneqq \dim_{\k} (\gr_n(G/G'';\k))$. 
For free groups, Chen showed that
$\theta_n(F_r)=(n-1)\binom{r+n-2}{n}$ for all $n\ge 2$.

\begin{theorem}[{\cite{SW-forum}}]
\label{thm:chen-epi}
For each $i\ge2$, the quotient map $q_i\colon G\to G/G^{(i)}$ induces a
natural epimorphism
\begin{equation*}
\label{eq:psi-map}
\Psi^{(i)}\colon \gr(G;\k)/\gr(G;\k)^{(i)} \longsurj  \gr(G/G^{(i)};\k).
\end{equation*}
\end{theorem}

Under suitable formality assumptions, these epimorphisms become isomorphisms; 
this will be made precise in Section~\ref{subsec:formality-groups}.

The following result, due to Papadima--Suciu \cite[Prop.~5.3 and Cor.~5.4]{PS-imrn04} 
and strengthened by Suciu--Wang \cite[Cor.~8.5]{SW-forum}, lifts Theorem~\ref{thm:chen-epi} 
to the derived quotients of the holonomy Lie algebra, via the canonical epimorphism 
$\Phi_G\colon \h(G;\k)\surj \gr(G;\k)$.

\begin{theorem}[\cite{PS-imrn04, SW-forum}]
\label{thm:ps-sw}
Let $G$ be a finitely generated group.
For each $i\ge 2$, there is a natural epimorphism of graded Lie algebras
\[
\Phi_G^{(i)} \colon  \h(G;\k)/\h(G;\k)^{(i)} \longsurj \gr(G/G^{(i)};\k).
\]
\end{theorem}

The map $\Phi_G^{(i)}$ may be described explicitly as follows.
Composing $\Phi_G$ with the projection 
$\gr(q_i)\colon \allowbreak  \gr(G;\k)\surj \gr(G/G^{(i)};\k)$
yields an epimorphism $\h(G;\k)\surj\gr(G/G^{(i)};\k)$ 
which factors through $\h(G;\k)/\h(G;\k)^{(i)}$. 
Thus, Theorem~\ref{thm:ps-sw} extends the 
holonomy comparison map $\Phi_G$ to the level of derived quotients.

The next result identifies the Malcev Lie algebra of a solvable quotient
with the corresponding quotient of $\m(G)$, and will be used later to
relate Chen Lie algebras to Malcev completions. 
For a subset $\fa\subset\m(G)$, denote by $\overline{\fa}$
its closure in the filtration topology on $\m(G)$.

\begin{theorem}[\cite{PS-imrn04}]
\label{thm:malcev-derived}
For any finitely generated group $G$ and any $i\ge1$, there is a natural
isomorphism of filtered Lie algebras
\[
\m(G/G^{(i)}) \cong \m(G)/\overline{\m(G)^{(i)}}.
\]
\end{theorem}

Passing to associated graded and using \eqref{eq:quillen-iso}, 
we obtain isomorphisms of graded Lie algebras
\[
\gr(G/G^{(i)}; \Q) \cong \gr\Bigl(\m(G)/\overline{\m(G)^{(i)}}\Bigr).
\]

%%%%%%%%%%%%%%%%%%%%
\subsection{Formality properties of groups}
\label{subsec:formality-groups}
%%%%%%%%%%%%%%%%%%%%

A finitely generated group $G$ is {\em $1$-formal} if a classifying space
$K(G,1)$ is $1$-formal. This property depends only on $G=\pi_1(X)$ for any
path-connected space $X$, and is preserved under finite free products and direct
products of finitely generated groups. In particular, the free groups $F_r$ are 
$1$-formal, since $K(F_r,1)=\bigvee^r(S^1)$ is formal, for all $r\ge 1$. 

The bridge between the $1$-minimal model of a group and its Malcev Lie
algebra is the following foundational result of Sullivan, used repeatedly
below.

\begin{theorem}[Sullivan {\cite{Sullivan77}}]
\label{thm:Sullivan-mg}
Let $G$ be a finitely generated group. There is a natural isomorphism 
between the tower of finite-dimen\-sional nilpotent Lie algebras
$\{\m(G/\gamma_i(G))\}_{i\ge 1}$ and the tower $\{\fL_i(G)\}_{i\ge 1}$ 
produced by the $1$-minimal model $\M_1(G)$, yielding a functorial 
isomorphism of pronilpotent Lie algebras
\[
\m(G)\cong \fL(\M_1(G)).
\]
\end{theorem}

Thus the assignments $G\mapsto\m(G)$ and $(\M_1,d)\mapsto\fL(\M_1)$ 
define adjoint functors between finitely generated groups (via their Malcev Lie 
algebras) and $1$-minimal $\cdgas$.

The $1$-formality property of a finitely generated group $G$ depends only on 
its Malcev Lie algebra $\m(G)$: indeed, $G$ is $1$-formal if and only if 
$\m(G)$ is isomorphic, as a filtered Lie algebra, to the degree completion of 
a finitely generated quadratic Lie algebra; see e.g.~\cite{CT95, 
Markl-Papadima}.  This yields a useful criterion. 

\begin{lemma}
\label{lem:1fg}
Let $G$ be a finitely generated group. Suppose there is a $1$-formal 
group $K$ and a homomorphism $\varphi\colon G\to K$ such that 
$\varphi^*\colon H^1(K;\Q) \to H^1(G;\Q)$ is an isomorphism and 
$\varphi^*\colon H^2(K;\Q) \to H^2(G;\Q)$ is injective. Then $G$ 
is also $1$-formal.
\end{lemma}

\begin{proof}
A map of connected CW-complexes which induces an isomorphism on $H^1(-;\Q)$ and a 
monomorphism on $H^2(-;\Q)$ induces an isomorphism on $1$-minimal models, by 
\cite{Sullivan77}; applied to a map $K(G,1)\to K(K,1)$ realizing $\varphi$, 
this gives an isomorphism $\M_1(K)\to\M_1(G)$, and hence 
$\m(G)\cong\m(K)$ by Theorem~\ref{thm:Sullivan-mg}.  Since $1$-formality 
depends only on the Malcev Lie algebra, $G$ inherits it from $K$.
\end{proof}

\begin{example}
\label{ex:low-betti}
If $G$ is a finitely generated group with $b_1(G)$ equal to $0$ or $1$, 
then $G$ is $1$-formal. Indeed, the claim is true for $K_0=\{1\}$ 
and for $K_1=\Z$. Moreover, if $b_1(G)=i\in \{0,1\}$, 
we may define a homomorphism $\varphi\colon G\to K_i$ 
satisfying the assumptions of Lemma~\ref{lem:1fg}. Therefore, 
the claim holds for $G$, too.
\end{example}

Equivalently, the $1$-formality of $G$---a property of $\m(G)$, or of the 
rationalization $G_\Q$---amounts to the existence of a filtered Lie algebra 
isomorphism
\begin{equation}
\label{eq:1formal-malcev}
\Theta_G \colon \m(G;\k) \longisom \widehat{\rule{0pt}{0.8em}\h(G;\k)} .
\end{equation}
Such an isomorphism is canonical up to composition with filtered
Lie algebra automorphisms of $\widehat{\h(G;\k)}$
acting trivially on degree~$1$ of the associated graded,
that is, on $H_1(G;\k)$.
In particular, no functorial choice of $\Theta_G$ exists in general.

The notion of $1$-formality of groups may be decomposed into two 
independent Lie-theoretic conditions, governing respectively the 
presentation of the graded Lie algebra $\gr(G;\k)$ and its relationship 
to the Malcev Lie algebra $\m(G;\k)$ via completion. Recall that 
$\h(G;\k)$ is quadratic by definition, whereas $\gr(G;\k)$ need not be.

We say that $G$ is \emph{graded-formal} (over $\k$) if the graded Lie 
algebra $\gr(G;\k)$ is quadratic. Since $\h(G;\k)$ is quadratic and the 
canonical map $\Phi_G\colon \h(G;\k)\to \gr(G;\k)$ is an epimorphism 
inducing isomorphisms in degrees $1$ and $2$ (Theorem~\ref{thm:holoepi}), 
it follows that $G$ is graded-formal if and only if $\Phi_G$ is an 
isomorphism. In particular, the holonomy ranks $\bar\phi_n(G)$ coincide 
with the LCS ranks $\phi_n(G)$, for all $n\ge 1$.

We say that $G$ is \emph{filtered-formal} (over $\k$) if there exists a filtered
Lie algebra isomorphism
\begin{equation}
\label{eq:filt-formal}
\m(G;\k) \longisom \widehat{\gr(G;\k)} .
\end{equation}
Such an isomorphism is canonical up to composition with filtered
Lie algebra automorphisms of $\widehat{\gr(G;\k)}$
acting trivially on degree~$1$ of the associated graded,
that is, on $H_1(G;\k)$. 
In particular, no functorial choice of such an isomorphism exists in general.
We denote by
\begin{equation}
\label{eq:theta-ff}
\Theta_G^{\mathrm{ff}} \colon \m(G;\k) \longisom \widehat{\gr(G;\k)}
\end{equation}
a choice of such an isomorphism when needed.
In the $1$-formal case, the two isomorphisms $\Theta_G^{\mathrm{ff}}$ and $\Theta_G$
represent the same element in the torsor of filtered Lie algebra isomorphisms
$\m(G;\k) \isom \widehat{\h(G;\k)}$ modulo automorphisms acting trivially on $H_1(G;\k)$. 

\begin{proposition}[\cite{PS-imrn04, SW-forum}]
\label{prop:formality-upgrades}
Let $G$ be a finitely generated group.

\begin{enumerate}
\item If $G$ is filtered-formal, then for each $i\ge 2$, the natural epimorphism
\[
\Psi^{(i)}\colon \gr(G;\k)/\gr(G;\k)^{(i)} \longrightarrow \gr(G/G^{(i)};\k)
\]
is an isomorphism of graded Lie algebras.

\item If $G$ is $1$-formal, then for each $i\ge 2$, the map
\[
\Phi_G^{(i)}\colon \h(G;\k)/\h(G;\k)^{(i)} \longrightarrow \gr(G/G^{(i)};\k)
\]
is an isomorphism of graded Lie algebras.
\end{enumerate}
\end{proposition}

We recall from \cite[Thm.~6.5]{SW-forum} a useful criterion for filtered-formality. 
For completeness, we include a proof of the part of the result that will be 
needed later on. 

Let $(\cM^{\bullet},d)$ be a minimal $\cdga$ generated in degree one, endowed 
with the canonical filtration $\{\cM_i\}_{i\ge 0}$ from \eqref{eq:min-mod}, 
where each sub-$\cdga$ $\cM_i$ given by 
a Hirsch extension of the form $\cM_{i-1}\otimes_{\k} \bigwedge(V_i)$. 
The underlying $\cga$ $\cM^{\bullet}$ possesses a natural set of positive 
weights, which we will refer to as the {\em Hirsch weights}: simply declare 
$V_i$ to have weight $i$, and extend those weights to $\cM^{\bullet}$ multiplicatively.  
We say that the $\cdga$ $(\cM^{\bullet},d)$ has {\em positive Hirsch weights}\/ if 
the differential $d$ is homogeneous with respect to those weights.  If this is the case, 
each sub-$\cdga$ $\cM_i$ also has positive Hirsch weights.

\begin{proposition}[\cite{SW-forum}]
\label{prop:positive-weights-ff}
Let $G$ be a finitely generated group.
If the canonical $1$-minimal model $\cM(G;\k)$ is filtered-isomorphic
to a $1$-minimal model with positive Hirsch weights, then $G$ is
filtered-formal over\/ $\k$.
\end{proposition}

\begin{proof}
Let $\cM$ be a $1$-minimal model with positive Hirsch weights.
By homogeneity of the differential, the dual Lie algebra $\fL(\cM)$ is
graded and its completion is isomorphic to the Malcev Lie algebra 
$\m(G;\k)$ by Theorem~\ref{thm:Sullivan-mg}. Thus $\m(G;\k)$ 
is the completion of a graded Lie algebra, which is
precisely the definition of filtered-formality.
\end{proof}

Proposition~\ref{prop:formality-upgrades} shows that the various notions of
formality can be detected purely at the level of associated graded and Malcev
Lie algebras. These Lie-theoretic descriptions become particularly effective in the
presence of a $1$-finite $1$-model, whose existence and basic properties
will be developed in Section~\ref{subsec:1f1m}.

%%%%%%%%%%%%%%%%%%%%
\subsection{Minimal models and $\CE$-complexes}
\label{subsec:min-CE}
%%%%%%%%%%%%%%%%%%%%

We now explain how Malcev Lie algebras give rise to canonical
$1$-minimal $\cdga$ models via Chevalley--Eilenberg constructions.
This correspondence provides the main bridge between the
Lie-theoretic invariants of finitely generated groups and the
$\cdga$-based resonance theory developed later in the paper.
In particular, it allows us to interpret resonance varieties
purely in terms of filtered Lie algebra data.

By Theorem~\ref{thm:Sullivan-mg}, the $1$-minimal model of a group
recovers its Malcev Lie algebra.  In the opposite direction,
if $\g$ is a finitely generated Lie algebra, the quotients 
$\g/\gamma_{n+1}(\g)$ assemble into central extensions
\begin{equation}
\label{eq:gnilp-extension}
\begin{tikzcd}[column sep=20pt]
0 \arrow[r] & \gamma_n(\g)/\gamma_{n+1}(\g)
 \arrow[r] & \g/\gamma_{n+1}(\g)
 \arrow[r] & \g/\gamma_n(\g) \arrow[r] & 0.
\end{tikzcd}
\end{equation}
Applying the Chevalley--Eilenberg cochain functor yields Hirsch extensions
\begin{equation}
\label{eq:ce-hook}
\CE(\g/\gamma_n(\g)) \longinj \CE(\g/\gamma_{n+1}(\g)),
\end{equation}
whose colimit 
\begin{equation}
\label{eq:ce-colimit}
\widehat{\CE}(\g)\coloneqq\varinjlim_{n} \CE(\g/\gamma_n(\g))
\end{equation}
is a $1$-minimal $\cdga$.  
For a finitely generated group $G$ there is a natural isomorphism
\begin{equation}
\label{eq:m1g-ce}
\M_1(G) \cong \widehat{\CE}(\m(G)).
\end{equation}

The Hirsch weights on $\widehat{\CE}(\m(G))$ coincide with the lower central
series grading on $\m(G)$, thereby encoding the filtered Lie algebra structure
directly into the $\cdga$ model.

The completed Chevalley--Eilenberg $\cdga$ $\widehat{\CE}(\m(G))$ provides 
a natural algebraic avatar of a finitely generated group $G$, functorially associated 
to its Malcev Lie algebra, whose Koszul complexes and resonance varieties 
capture the group’s infinitesimal structure.

\begin{proposition}
\label{prop:res-malcev}
Let $G$ be a finitely generated group, and let $\m(G)$ be its Malcev Lie algebra,
endowed with the Malcev filtration. Then the resonance varieties of the completed 
Chevalley--Eilenberg $\cdga$ $\widehat{\CE}(\m(G))$ depend only on the filtered 
Lie algebra structure of $\m(G)$.
\end{proposition}

\begin{proof}
An isomorphism of filtered Lie algebras induces an isomorphism of the associated
completed Chevalley--Eilenberg $\cdgas$. Since the Koszul complexes and their homology
are functorial under $\cdga$ isomorphisms, the associated resonance varieties depend
only on the filtered Lie algebra structure of $\m(G)$.
\end{proof}

\begin{corollary}
\label{cor:ff-res}
If $G$ is filtered-formal, then the support resonance varieties of
$\widehat{\CE}(\m(G))$ are finite unions of rational linear subspaces.
\end{corollary}

\begin{proof}
If $G$ is filtered-formal, then $\m(G)$ is isomorphic, as a filtered Lie algebra,
to the completed graded Lie algebra $\widehat{\gr(G;\Q)}$. The associated completed
Chevalley--Eilenberg $\cdga$ is therefore isomorphic to the completion of a quadratic
$\cdga$ with linear differential. By the results of Section~\ref{sect:koszul},
the support resonance varieties of such $\cdgas$ are unions of rational linear
subspaces.
\end{proof}

\begin{corollary}
\label{cor:resonance-germ-nil}
Let $G$ be a finitely generated, torsion-free nilpotent group, let $X$ be
the associated compact nilmanifold, and let $A=\CE(\g)$ be the
Chevalley--Eilenberg $\cdga$ of its Malcev Lie algebra $\g$.  Then
$\RR^i(A)\subseteq\{0\}$ for every $i\ge 0$; in particular, every germ
$\RR^i(A)_{(0)}$ is trivial.
\end{corollary}

\begin{proof}
By Nomizu's theorem, $\CE(\g)$ is a finite-type Sullivan minimal model for
the nilmanifold $X$; and $\g$ is nilpotent, since $G$ is torsion-free
nilpotent.  Theorem~\ref{thm:CE-resonance-nilpotent} therefore applies.
\end{proof}

%%%%%%%%%%%%%%%%%%%%%%%%%
\subsection{Groups with $1$-finite $1$-models}
\label{subsec:1f1m}
%%%%%%%%%%%%%%%%%%%%%%%%%

The purpose of this subsection is to single out a class of finitely generated
groups for which rational homotopy--theoretic methods yield effective and
computable invariants. Throughout, we fix a field $\k$ of characteristic~$0$.

The next theorem provides a precise description of the Malcev Lie algebra of a
group admitting a finite-type model in low degrees.

\begin{theorem}[{\cite[Cor.~5.8]{PS-jlms}}]
\label{thm:malcev-holo}
Let $G$ be a finitely generated group that admits a $1$-finite $1$-model $(A,d)$
over $\k$. Then its Malcev Lie algebra satisfies
\[
\m(G;\k)\cong \widehat{\h}(A),
\]
the lower central series completion of the holonomy Lie algebra of $A$.
\end{theorem}

Thus, a finitely generated group admits a $1$-finite $1$-model if and only if
its Malcev Lie algebra is the lower central series completion of a finitely
presented Lie algebra, generated in degree~$1$ with relations in degrees~$1$
and~$2$.

Theorem~\ref{thm:malcev-holo} is what makes the Koszul machinery of the
preceding sections available for groups: it converts statements about
$\h(A)$, and hence about $\B_1(A)$, into statements about $\m(G;\k)$.  We
exploit this in Section~\ref{subsec:koszul-holo}, where degeneration of the
weight spectral sequence in low degrees is shown to control the metabelian
quotient of $\m(G;\k)$, and in Section~\ref{subsec:chen-ranks}, where the
resulting numerical consequences for Chen ranks are drawn.  For the
remainder of this subsection we take up a prior question: which groups admit
a $1$-finite $1$-model at all.

For a $1$-formal group $G$, the cohomology algebra $(H^*(G;\k),0)$ provides a
$1$-finite $1$-model; this is the main source of groups in the class under
discussion.  The existence of a $1$-finite $1$-model is nevertheless a strong
and subtle condition, independent of formality, and it interacts delicately
with functoriality; we return to this in
Section~\ref{subsec:functorial-formal-hom}.

The existence of a $1$-finite $1$-model can be obstructed at the level of
minimal Sullivan models.
The following finiteness result from \cite{PS-jlms} provides a basic necessary
criterion.

\begin{theorem}[\cite{PS-jlms}]
\label{thm:ps-finite-betti}
Let $X$ be a space admitting a $q$-finite $q$-model, and let $\M_q(X)$ denote its
$q$-minimal Sullivan model. Then $\dim_\k H^i(\M_q(X))<\infty$ for all $i\le q+1$.
\end{theorem}

As an immediate consequence, if a finitely generated group $G$ admits a
$1$-finite $1$-model, then the second Betti number of its $1$-minimal model is
finite.  The same conclusion holds, by a classical argument, when $G$ is
finitely presented---a condition which, as we note below, neither implies nor
is implied by the existence of a $1$-finite $1$-model.

This criterion excludes large and natural families of groups. In particular,
if $\pi$ is a finitely generated, very large group (i.e., admits a free,
non-cyclic quotient) and $G=\pi/\pi''$ is its maximal metabelian quotient,
then $G$ admits no $1$-finite $1$-model \cite{PS-jlms}. In particular, free
metabelian groups $F_n/F_n''$ with $n\ge 2$ do not belong to this class.

On the other hand, the existence of a $1$-finite $1$-model is not implied by
finite presentability. Using Alexander polynomials together with a refined
form of the Tangent Cone Theorem~\ref{thm:dps-dp}, due to Budur and Wang
\cite{BW20}, one can construct finitely presented groups---indeed, fundamental
groups of closed, orientable $3$-manifolds---which admit no $1$-finite
$1$-model; see \cite{Su-mm} and the examples discussed there. 
These examples reflect a general obstruction principle: failure of
rational linearity for tangent cones to characteristic varieties in degree~$1$
precludes the existence of a $1$-finite $1$-model.

From this point on, our default setting will be that of finitely generated
groups admitting a $1$-finite $1$-model. Within this framework, classical
results for $1$-formal groups appear as special cases, while genuinely new
phenomena arise from the presence of a nontrivial differential. In particular,
formality assumptions will be viewed as conditions under which the general
model-theoretic constructions simplify or degenerate.

%%%%%%%%%%%%%%%%%%%%%
\subsection{Holonomy--associated graded comparison in the Carnot case}
\label{subsec:holo-comparison}
%%%%%%%%%%%%%%%%%%%%%

We consider the comparison between holonomy Lie algebras of $1$-finite $1$-models 
and the associated graded Lie algebra of a finitely generated nilpotent group. Such a 
comparison is not available in general, but becomes canonical under strong rigidity 
assumptions, in particular in the Carnot (graded) setting described below.

\begin{proposition}
\label{prop:carnot-Phi}
Let $G$ be a finitely generated nilpotent group with Malcev Lie algebra
$\g = \m(G)\otimes_{\Q}\k$, and set $A = \CE(\g)$. 
Assume that $\g$ is Carnot graded, i.e., $\g \cong \gr(\g)$
as graded Lie algebras. Then there is an isomorphism 
\[
\Phi_A \colon \h(A) \longrightarrow \gr(G;\k). 
\]
\end{proposition}

\begin{proof}
Since $A=\CE(\g)$, Proposition~\ref{prop:holo-CE} gives a canonical 
isomorphism $\h(A)\cong \g$.
Moreover, since $\m(G)\otimes_{\Q}\k \cong \g$ 
as filtered Lie algebras, taking associated graded yields
$\gr(G;\k) \cong \gr(\g)$. 
Under the Carnot assumption, the grading on $\g$ identifies it 
canonically with its associated graded: $\g \cong \gr(\g)$. 
Combining these identifications yields
\[
\h(A) \cong \g \cong \gr(\g) \cong \gr(G;\k).
\]
The resulting isomorphism $\Phi_A$ is obtained by composing the 
above identifications; it is not canonical and depends on the choice 
of Carnot grading on $\g$.
\end{proof}

In general, the construction above is not functorial in $G$ or $A$, 
and therefore does not yield a canonical comparison map 
between $\h(A)$ and $\gr(G;\k)$ beyond the Carnot setting.
This motivates the following conjectural lifting property.

\begin{conjecture}
\label{conj:holo-lift}
Let $G$ be a finitely generated group admitting a connected
$1$-finite $1$-model $(A,d)$ with positive weights such that
$\im(d^\vee)$ is concentrated in weights $\ge 2$.
Under these hypotheses, Lemma~\ref{lem:holo-to-cohomology}
provides a natural epimorphism 
$\Xi_A\colon \h(H^*(A))\surj\h(A)$.
We conjecture that there exists an epimorphism of filtered 
Lie algebras
\[
\Phi_A\colon \h(A) \longsurj \gr(G;\k)
\]
such that $\Phi_G = \Phi_A \circ \Xi_A$.
\end{conjecture}

\begin{remark}
\label{rem:Xi-Phi}
The conjecture can be reformulated as the inclusion
$\ker(\Xi_A) \subseteq \ker(\Phi_G)$.
In the Carnot case, the conjecture holds by Proposition~\ref{prop:carnot-Phi},
and all holonomy and graded objects involved are isomorphic.
\end{remark}

%%%%%%%%%%%%%%%%%%%%%
\subsection{Degeneration and the metabelian quotient}
\label{subsec:koszul-holo}
%%%%%%%%%%%%%%%%%%%%%

The identification of the infinitesimal Alexander invariant $\B_1(A)$ with
$\B(\h(A))$ from Theorem~\ref{thm:B-holo} allows one to translate
homological information about the Koszul complex $K_\bullet(A)$ into
structural constraints on the holonomy Lie algebra.  Since $\B_1$ depends
only on the metabelian quotient $\h/\h''$, what one obtains this way is
information about that quotient, and no more; we record next both what
degeneration of the weight spectral sequence yields and where it stops.

\begin{theorem}
\label{thm:koszul-holo-quadratic}
Let $A = \M_1(G)$ be the $1$-minimal Sullivan model of a finitely generated
group $G$, equipped with positive Hirsch weights. If the weight spectral
sequence of $K_\bullet(A)$ degenerates at the $E_2$-page in total degree
$\le 2$, then the canonical epimorphism
$\Xi_A\colon \h(H^*(A))\surj\h(A)$ of
Lemma~\ref{lem:holo-to-cohomology} satisfies
\[
\ker(\Xi_A)\ \subseteq\ \h(H^*(A))'' .
\]
Consequently, the Chen ranks of $G$ agree with those of its holonomy Lie
algebra: $\theta_n(G)=\bar\theta_n(G)$ for all $n\ge1$.
\end{theorem}

\begin{proof}
Since $A$ has positive Hirsch weights, $H_1(A) = V_1^\vee$ is concentrated
in weight~$1$, so $\h(A)$ is positively graded and generated in weight~$1$
(Proposition~\ref{prop:pos-weights-graded-holo}).
For a positively graded Lie algebra generated in weight~$1$, the LCS-associated
graded coincides with $\h(A)$ itself: $\gr(\h(A)) = \h(A)$.

By Theorem~\ref{thm:malcev-holo}, $\widehat{\h(A)} \cong \m(G;\k)$.
For a finitely generated positively graded Lie algebra $L$ with separated
filtration, $\gr(\widehat{L}) = \gr(L)$.
Applying Quillen's isomorphism \eqref{eq:quillen-iso} therefore gives
\[
\h(A) = \gr(\h(A)) \cong \gr(\m(G;\k))
\cong \gr(G;\k).
\]
On the other hand, degeneration at $E_2$ in total degree $\le2$ gives
$\gr^W\B_1(A)\cong\B_1(H^*(A))$ by Theorem~\ref{thm:weight-ss}, whose
$E^2$- and $E^\infty$-identifications compare the two sides degreewise; and
$\gr^W=\gr^\m$ here by Corollary~\ref{cor:weight-splitting}\eqref{ws2},
$H^1(A)$ being concentrated in weight~$1$.  Since $A$ has positive weights, the
$\m$-adic filtration on $\B_1(A)$ is separated
(Proposition~\ref{prop:koszul-separation}), and both modules are graded, so
the surjection $\B_1(H^*(A))\surj\B_1(A)$ induced by $\Xi_A$ is in fact an
isomorphism.  By Theorem~\ref{thm:B-holo} this says that
$\B(\Xi_A)\colon\B(\h(H^*(A)))\to\B(\h(A))$ is an isomorphism, which forces
$\ker(\Xi_A)\subseteq\h(H^*(A))''$.

Finally, $\h(H^*(A))=\h(G;\k)$ and $\h(A)\cong\gr(G;\k)$, so
$\ker(\Xi_A)\subseteq\h(G;\k)''$ means that $\Xi_A$ descends to an
isomorphism $\h(G;\k)/\h(G;\k)''\isom\gr(G;\k)/\gr(G;\k)''$.  Taking
dimensions in each degree gives $\theta_n(G)=\bar\theta_n(G)$ for all $n$.
\end{proof}

\begin{remark}
\label{rem:holo-quadratic-limits}
Theorem~\ref{thm:koszul-holo-quadratic} does not assert that $\Xi_A$ is an
isomorphism, and the argument does not give this.  Degeneration controls
$\B_1$, hence the metabelian quotient $\h/\h''$, but supplies no information
about $\ker(\Xi_A)$ inside $\h''$; in particular it does not compare the
Hilbert series of $\h(H^*(A))$ and $\h(A)$.  Were $\Xi_A$ known to be an
isomorphism, Theorem~\ref{thm:malcev-holo} would give
$\m(G;\k)\cong\widehat{\h(H^*(A))}$, so $\m(G;\k)$ would be quadratic and
$G$ would be $1$-formal.  Whether degeneration of the weight spectral
sequence in low total degrees implies this---equivalently, whether it
implies $1$-formality---is precisely
Question~\ref{quest:weight-ss-1formal} in the case $q=1$, which remains
open.  Note also that
the hypothesis must be degeneration at $E_2$ and not at $E_1$: the
differential $d_1$ is determined by the cup product, and $E_1$-degeneration
would force $\gr^\m\B_1(A)=H_1(A)\otimes_{\k} S$.
\end{remark}

%%%%%%%%%%%%%%%%%%%%%%
\subsection{Finitely generated, torsion-free nilpotent groups}
\label{subsec:nilpotent-groups}
%%%%%%%%%%%%%%%%%%%%%%

Let $G$ be a finitely generated, torsion-free nilpotent group.
By Malcev theory, $G$ admits a rational Malcev completion $\m(G)$, a
finite-dimensional nilpotent Lie algebra over $\Q$.
Moreover, the Eilenberg--MacLane space $K(G,1)$ is a nilmanifold, and
Nomizu's theorem identifies the rational cohomology of $G$ with the
Lie algebra cohomology of $\m(G)$,
\begin{equation}
\label{eq:nomizu}
H^{*}(G;\Q) \cong H^{*}(\m(G)).
\end{equation}
Consequently, the Chevalley--Eilenberg $\cdga$ $\CE(\m(G))$ provides a
rational $1$-model for $G$, functorial with respect to group homomorphisms.
Since the Chevalley--Eilenberg differential preserves the grading by
bracket length, $\CE(\m(G))$ admits a positive-weight decomposition.
As a result, the Koszul complex $K_\bullet(\CE(\m(G)))$ carries a natural
weight filtration, giving rise to a spectral sequence of the type
considered in Theorem~\ref{thm:koszul-formal-obs}.

An Eilenberg--MacLane space $K(G,1)$ for such a group $G$ admits
a canonical realization as a nilmanifold $M=\R^{n}/G$, where $n=\dim \m(G)$. 
In this setting, $\CE(\m(G))$ is not only
a rational $1$-model for $G$, but in fact a finite-type Sullivan model for
the nilmanifold $M$ itself.
In particular, $\CE(\m(G))$ computes the full rational cohomology of $M$
and provides a common algebraic model for both $M$ and its fundamental
group.

Finitely generated, torsion-free nilpotent groups thus form a natural
class of examples for the techniques developed in this paper.
Their classifying spaces admit canonical finite-type models with positive
weights, making them particularly well suited for the Koszul spectral
sequence methods of Sections~\ref{subsec:weights-koszul} and 
\ref{subsec:koszul-formal-obstruction}.

\begin{proposition}
\label{prop:2step-filtered-formal}
Let $G$ be a finitely generated, torsion-free, $2$-step nilpotent group.
Then $G$ is filtered-formal.
\end{proposition}

\begin{proof}
The Malcev Lie algebra $\m(G)$ is finite-dimensional and $2$-step nilpotent.
By Proposition~\ref{prop:CE-positive-weights}, the $\cdga$ $\CE(\m(G))$ 
admits a positive-weight decomposition. This induces a splitting of the lower 
central series filtration on $\m(G)$, yielding a filtered Lie algebra isomorphism
\begin{equation}
\label{eq:mg-grmg}
\m(G) \cong \widehat{\gr(\m(G))}.
\end{equation}
By definition, this identifies $G$ as a filtered-formal group.
\end{proof}

The higher Heisenberg groups provide a concrete illustration of the
interaction between filtered-formality, graded-formality, and the
Koszul spectral obstruction introduced in
Section~\ref{subsec:koszul-formal-obstruction}.
For $n\ge1$, the $n$th Heisenberg group $H(n)$ is defined by the presentation
\begin{equation}
\label{eq:Hn-group}
H(n)=\langle a_1,\dots,a_n,b_1,\dots,b_n \mid
[a_i,b_i]=[a_j,b_j],\;
[a_i,a_j]=[a_i,b_j]=[b_i,b_j]=1 \text{ for } i\neq j \rangle.
\end{equation}
Each group $H(n)$ is torsion-free, $2$-step nilpotent, with infinite cyclic
center.
Its Malcev Lie algebra is the higher Heisenberg Lie algebra $\h(n)$,
generated by $x_i,y_i,z$ with relations $[x_i,y_i]=z$.

\begin{proposition}
\label{prop:Heisenberg-1formal}
For all $n\ge1$, the higher Heisenberg groups $H(n)$ are filtered-formal.
Moreover, $H(1)$ is not $1$-formal, while $H(n)$ is $1$-formal for all $n>1$.
\end{proposition}

\begin{proof}
By Proposition~\ref{prop:2step-filtered-formal}, all $H(n)$ are filtered-formal. 
Their rational $1$-minimal model is
\[
A_{H(n)} = \CE(\h(n)) = \Bigl(\bigwedge(a_1,\dots,a_n,b_1,\dots,b_n,c), d\Bigr),
\qquad d(c) = - \sum_{i=1}^n a_i \wedge b_i
\]
(see Section~\ref{subsec:heisenberg}).  

For $n>1$, the cup-product map
$\mu \colon \bwedge^2 H^1(H(n);\Q) \to H^2(H(n);\Q)$ is surjective. 
Since $H(n)$ is $2$-step nilpotent, this ensures that the canonical morphism
$\Phi_{H(n)} \colon \h(H(n);\Q) \surj \gr(H(n);\Q)$ is an isomorphism. 
Hence $H(n)$ is graded-formal, and therefore $1$-formal.

For $n=1$, we have $a_1\wedge b_1 = - d(c)$, so $\mu = 0$. 
Thus $\h(H(1);\Q) = \Lie(\Q^2)$ is free, while $\gr(H(1);\Q) = \h(1)$
is $2$-step nilpotent, and $\Phi_{H(1)}$ is not an isomorphism. 
Therefore $H(1)$ is not graded-formal, and not $1$-formal.
\end{proof}

The examples above suggest a close relationship between the structure of
a nilpotent Lie algebra and the degree of partial formality of the
associated space.

For the Heisenberg groups $H(n)$, and more generally for Heisenberg-type
groups with defining $2$-form of rank $2m$ \cite{Mc10}, the Lie algebra is
$2$-step nilpotent, yet the corresponding nilmanifold is $(m-1)$-formal
but not $m$-formal.
From the perspective of the Koszul spectral sequence, the failure of
$q$-formality is detected by a nontrivial differential whose source lies
in total degree $q+1$, reflecting the first appearance of relations beyond
quadratic data in the Chevalley--Eilenberg differential.

More generally, while the nilpotency class alone is too coarse to control
$q$-formality, the present results indicate that finer invariants---such as
the weights and ranks appearing in $\CE(\g)$---play a decisive role via the
associated Koszul spectral sequence.

%%%%%%%%%%%%%%%%%%%%%%
\subsection{Partial formality and vanishing resonance}
\label{subsec:partial-formality-resonance}
%%%%%%%%%%%%%%%%%%%%%%

We now explain how the rigidity of resonance for Chevalley--Eilenberg
algebras of nilpotent Lie algebras propagates to cohomology algebras
under partial formality assumptions.

\begin{theorem}[M\u{a}cinic {\cite{Mc10}}]
\label{thm:macinic}
Let $G$ be a finitely generated, torsion-free nilpotent group, and suppose
that $G$ is $q$-formal.
Then
\[
\RR^{i}\bigl(H^{*}(G;\Q)\bigr) \subseteq \{0\},
\qquad \text{for all } i \le q .
\]
\end{theorem}

\begin{proof}
Since $\CE(\m(G))$ is a Sullivan $1$-model for $G$, $q$-formality
provides a $q$-quasi-isomor\-phism $\CE(\m(G)) \simeq H^*(G;\Q)$. 
By Theorem~\ref{thm:q-iso-koszul-coh}, this induces isomorphisms
\[
\B^i(\CE(\m(G))) \cong \B^i(H^*(G;\Q))
\quad \text{for all } i \le q.
\]

By Theorem~\ref{thm:CE-resonance-nilpotent}, the resonance varieties
$\RR^{i}(\CE(\m(G)))$ are trivial for all $i$, and hence so are
$\RR^{i}(H^*(G;\Q))$ for all $i \le q$.
\end{proof}

The higher Heisenberg groups $H(n)$ illustrate both the scope and the
sharpness of Theorem~\ref{thm:macinic}.
As noted above, $H(n)$ is $(n-1)$-formal but not $n$-formal, and one has
\[
\RR^{i}\bigl(H^{*}(H(n);\Q)\bigr)=\{0\},
\quad \text{for all } i \le n-1 ,
\]
while nontrivial resonance appears precisely in degree~$n$.
Thus, partial formality exactly governs the range in which resonance is
forced to vanish.

A precursor of Theorem~\ref{thm:macinic} appears in work of Carlson and
Toledo~\cite{CT95}, who proved the vanishing of degree-one resonance for
$1$-formal nilmanifolds.
The result above may be viewed as a higher-degree refinement, made possible
by the Koszul spectral sequence techniques developed in this paper.

%%%%%%%%%%%%%%%%%%%%%
\section{Chen ranks, Koszul modules, and resonance}
\label{sect:chen-res}
%%%%%%%%%%%%%%%%%%%%%

We study Chen ranks through the lens of $1$-finite $1$-models and 
their Koszul modules. In the $1$-formal setting, Chen ranks are controlled 
by the holonomy Lie algebra via resonance; beyond formality, Koszul 
modules provide the correct governing object. The main result establishes 
Theorem~\ref{thm:alex-completed-intro} from the Introduction.

%%%%%%%%%%%%%%%%%%%%%
\subsection{Chen ranks of groups with $1$-finite $1$-models}
\label{subsec:chen-ranks}
%%%%%%%%%%%%%%%%%%%%%

Let $G$ be a finitely generated group, and let $G/G''$ be its maximal 
metabelian quotient.  Recall that the Chen ranks of $G$ 
are the dimensions of the associated graded Lie algebra of $G/G''$ 
with respect to the lower central series (over a field $\k$ of characteristic $0$), 
\begin{equation}
\label{eq:chen-ranks}
\theta_n(G)= \dim_{\k}\bigl( \gr_n(G/G'' ;\k) \bigr).
\end{equation}
Our goal in this subsection is to express Chen ranks in terms of Koszul modules 
associated to $1$-finite $1$-models, thereby separating intrinsic group-theoretic 
input from additional formality assumptions.

A classical result of Massey \cite{Massey} (see \cite{Su-pisa} for details) 
relates the Chen groups $\gr_n(G/G'')$ ($n\ge 2$) to the graded pieces of 
the associated graded of the Alexander invariant $B_1(G)=G'/G''$, 
viewed as a module over $\Z[G_{\ab}]$, 
\begin{equation}
\label{eq:massey-chen}
\gr_{n} (G/G'') \cong \gr_{n-2}(B_1(G))  \quad \text{for all $n\ge 2$},
\end{equation}
from which it follows that 
\begin{equation}
\label{eq:hilb-b-chen}
 \sum_{n\ge 0} \theta_{n+2} (G) t^n = \Hilb (\gr(B_1(G;\k)),t) . 
\end{equation}

Let $I$ be the augmentation ideal of $\Lambda=\Q[G_{\ab}]$.
The next result (\cite[Prop.~5.4]{DPS-duke}), relates the
$I$-adic completion of the rational Alexander invariant $B_1(G;\Q)$,
viewed as a module over $\widehat{\Lambda}$, to the Malcev Lie algebra of $G$.  

\begin{theorem}[{\cite[Prop.~5.4]{DPS-duke}}]
\label{thm:alex-malcev}
For any finitely generated group $G$ there is a filtration-preserving,
$\widehat{\Lambda}$-linear isomorphism
\[
\widehat{B_1(G;\Q)} \cong \overline{\m(G)'}\,/\,\overline{\m(G)''}.
\]
\end{theorem}

We are now ready to state and prove the main result of this section,
which yields Theorem~\ref{thm:alex-completed-intro} stated in the Introduction.

\begin{theorem}
\label{thm:alex-completed}
Let $G$ be a finitely generated group admitting a $1$-finite $1$-model
$(A,d)$ over a field $\k$ of characteristic~$0$. Set $S=\Sym(H_1(G;\k))$. Then:
\begin{enumerate}[itemsep=3pt]
\item \label{ch1}
$\widehat{B_1(G;\k)} \cong \widehat{\B_1(A)}$ as filtered modules over
$\widehat{\k[G_{\ab}]}\cong \widehat{S}$.
\item \label{ch2}
$\gr^I(B_1(G;\k)) \cong \gr^\m(\B_1(A))$ as graded $S$-modules, the 
left-hand side taken with respect to the $I$-adic filtration.
\item \label{ch3}
$\theta_n(G)=\theta_n(A)$, for all $n\ge1$.
Moreover, for all $n\ge2$,
\begin{equation}
\label{eq:chen-gr-B1}
\theta_n(G)
= \dim_\k \gr^\m_{n-2}\bigl(\B_1(A)\bigr),
\end{equation}
and the generating series of Chen ranks of $G$ satisfies
\begin{equation}
\label{eq:chen-gen-series}
\sum_{n\ge0} \theta_{n+2}(G)\, t^n
= \Hilb\bigl(\gr^\m(\B_1(A)),t\bigr).
\end{equation}
\end{enumerate}
\end{theorem}

\begin{proof}
The comparison between group-theoretic and infinitesimal invariants proceeds via 
Malcev completion and holonomy Lie algebras, which allow us to pass from non-linear 
group data to linearized, graded objects.

By Theorem~\ref{thm:alex-malcev}, the completion of the Alexander invariant $B_1(G;\k)$ is 
naturally identified with the completion of the maximal metabelian quotient of the Malcev Lie 
algebra $\m(G)\otimes_{\Q} \k$,
\begin{equation}
\label{eq:isom1}
\widehat{B_1(G;\k)} 
\cong 
\big(\overline{\m(G)'}/\overline{\m(G)''}\big)\otimes_{\Q} \k .
\end{equation}

Since $G$ admits a $1$-finite $1$-model $(A,d)$ over $\k$, Theorem~\ref{thm:malcev-holo} 
applies and identifies the Malcev Lie algebra $\m(G)\otimes_{\Q} \k$ with the completed 
holonomy Lie algebra $\h(A)$. It follows that 
\begin{equation}
\label{eq:isom2}
\big(\overline{\m(G)'}/\overline{\m(G)''}\big)\otimes_{\Q} \k
\cong
\overline{\h(A)'}/\overline{\h(A)''} .
\end{equation}

Consider the infinitesimal Alexander invariant $\B(\h(A))=\h(A)'/\h(A)''$, 
defined purely at the Lie algebra level. 
The ideals $\h(A)'$ and $\h(A)''$ inherit the lower central series filtration, and 
thus the closure of each in $\widehat{\h(A)}$ coincides with its completion. 
Since completion is exact for quotients with respect to this filtration, 
we have a canonical identification
\begin{equation}
\label{eq:isom3}
\overline{\h(A)'}/\overline{\h(A)''} \cong \widehat{\h(A)'/\h(A)''} = \widehat{\B(\h(A))}. 
\end{equation}

Finally, Theorem~\ref{thm:B-holo} identifies $\B(\h(A))=\h(A)'/\h(A)''$ with 
the first Koszul module $\B_1(A)$. Passing to completions, we obtain an isomorphism 
\begin{equation}
\label{eq:isom4}
\widehat{\B(\h(A))}
\cong
\widehat{\B_1(A)}. 
\end{equation}

We note that each of the isomorphisms \eqref{eq:isom1}--\eqref{eq:isom4} is 
compatible with the natural action of $H_1(G;\k)$ and hence is $\widehat{S}$-linear 
after completion.
Indeed, the action of the abelianization on Alexan\-der-type invariants is induced
by conjugation, which corresponds under Malcev completion to the adjoint action
of the abelianized Lie algebra. The identification between Malcev and holonomy
Lie algebras preserves abelianizations and filtrations, and completion is taken
with respect to the same filtration throughout. Consequently, the resulting
composite isomorphism is an isomorphism of filtered $\widehat{S}$-modules.
This proves claim \eqref{ch1}.  Passing to associated graded objects gives \eqref{ch2}. 

Finally, Massey’s correspondence \eqref{eq:massey-chen} identifies the Chen ranks of 
$G$ with the graded pieces of $\gr(B_1(G;\k))$, up to a degree shift of $2$, while 
Proposition~\ref{prop:holo-massey} gives the corresponding description of 
the holonomy Chen ranks of $A$ in terms of $\gr(\B_1(A))$. Combined 
with \eqref{ch2}, this yields \eqref{ch3}.
\end{proof}

We next record a naturality property of Theorem~\ref{thm:alex-completed},
together with its Malcev counterpart.  Both presuppose a notion of model for
a homomorphism, which we isolate first, following the terminology
of \cite{PS-adv19}.

\begin{definition}
\label{def:modeled-hom}
Let $\alpha\colon G\to H$ be a homomorphism between finitely generated groups,
and let $A_G$ and $A_H$ be $1$-finite $1$-models over $\k$ for $G$ and $H$.
A $\cdga$ morphism $\varphi\colon A_H\to A_G$ \emph{models}\/ $\alpha$ if the
square
\begin{equation}
\label{eq:models-alpha}
\begin{tikzcd}[column sep=24pt, row sep=22pt]
\M_1(H) \ar[d, "\M_1(f_\alpha)"] \ar[r] & A_H \ar[d, "\varphi"] \\
\M_1(G) \ar[r] & A_G
\end{tikzcd}
\end{equation}
commutes up to $\cdga$ homotopy, the horizontal arrows being
$1$-quasi-isomorphisms and $f_\alpha\colon K(G,1)\to K(H,1)$ a map inducing
$\alpha$ on fundamental groups.  A family of homomorphisms is
\emph{uniformly modeled}\/ if each member is equipped with such a $\varphi$,
compatibly with composition.
\end{definition}

\begin{lemma}
\label{lem:functorial-malcev}
Let $\alpha\colon G\to H$ be a homomorphism between finitely generated
groups, and let $A(\alpha)\colon A_H\to A_G$ be a $\cdga$ morphism modeling
$\alpha$, in the sense of Definition~\ref{def:modeled-hom}.
Then the induced morphism
\[
\widehat{\h}(A(\alpha))\colon  \widehat{\h}(A_G) \longrightarrow \widehat{\h}(A_H)
\]
corresponds, under the Malcev--holonomy identifications
$\widehat{\h}(A_G) \cong \m(G;\k)$ and $\widehat{\h}(A_H) \cong \m(H;\k)$,
to the Malcev Lie algebra morphism
\[
\m(\alpha)\colon \m(G;\k) \longrightarrow \m(H;\k).
\]
\end{lemma}

\begin{proof}\sloppy 
Since $A(\alpha)$ models $\alpha$, it is $1$-equivalent to
$\apl(f_\alpha)$, where $f_\alpha\colon K(G,1) \to K(H,1)$ induces $\alpha$
on fundamental groups.
The functoriality of minimal $1$-models identifies the induced map on
completed Chevalley--Eilenberg complexes with
$\widehat{\CE}(\m(\alpha))$, and the claim follows from
Theorem~\ref{thm:malcev-holo}.
\end{proof}

\begin{lemma}
\label{lem:alex-functorial}
Let $\alpha\colon G\to H$ be a homomorphism between finitely generated groups
admitting $1$-finite $1$-models $A_G$ and $A_H$, and let
$A(\alpha)\colon A_H\to A_G$ be a $\cdga$ morphism modeling $\alpha$.
Then the identifications in Theorem~\ref{thm:alex-completed} are natural with
respect to $\alpha$, in the sense that the diagram
\[
\begin{tikzcd}[column sep=26pt, row sep=22pt]
\widehat{\B_1(A_G)} \ar[r, "\cong"] \ar[d, "\widehat{\B_1(A(\alpha))}"]
&
\widehat{B_1(G;\k)} \ar[d, "\widehat{B_1(\alpha)}"]
\\
\widehat{\B_1(A_H)} \ar[r, "\cong"]
&
\widehat{B_1(H;\k)}
\end{tikzcd}
\]
commutes.
\end{lemma}

\begin{proof}
The isomorphism of Theorem~\ref{thm:alex-completed}\eqref{ch1} is the
composite of \eqref{eq:isom1}--\eqref{eq:isom4}, and each factor is natural
in the required sense.  The isomorphism \eqref{eq:isom1} is induced by the
Malcev completion map of \cite[Prop.~5.4]{DPS-duke}, hence is natural with
respect to group homomorphisms, since Malcev completion is a functor and both
sides of Theorem~\ref{thm:alex-malcev} are functors of $G$.  The
identification \eqref{eq:isom2} is natural with respect to $\alpha$ by
Lemma~\ref{lem:functorial-malcev}, which applies because $A(\alpha)$ models
$\alpha$.  Finally, \eqref{eq:isom3} is induced by passage to closures in a
filtered exact sequence, and \eqref{eq:isom4} is natural with respect to
$\cdga$ morphisms by Proposition~\ref{prop:naturality-B-holo}, applied to
$A(\alpha)$.  Composing the four squares gives the assertion.
\end{proof}

Thus, for groups admitting a $1$-finite $1$-model, the Chen ranks depend 
only on the associated Koszul module of the model, via constructions that 
are natural with respect to group homomorphisms, independently of any 
a priori formality assumptions.

We now restrict attention to the case in which the $1$-model $(A,d)$ is given 
by the cohomology algebra $H^*(G;\k)$ with $d=0$. For finitely generated groups, 
$H^1(G;\k)$ is finite-dimensional; hence, this $\cdga$ is automatically $1$-finite, 
and the preceding discussion applies verbatim. In this situation, the relevant 
infinitesimal data are encoded by the holonomy Lie algebra 
\begin{equation}
\label{eq:holo-group}
\h(G;\k)\coloneqq \h(H^{*}(G;\k))
\end{equation}
and its associated infinitesimal Alexander invariant $\B(\h(G;\k))=\h(G;\k)'/\h(G;\k)''$, 
viewed as a graded module over the ring $S=\Sym(H_1(G;\k))$. By Theorem~\ref{thm:B-holo}, 
this module is isomorphic to the Koszul module $\B_1(G;\k)\coloneqq \B_1(H^{*}(G;\k))$. 

Under the aforementioned $1$-formality hypothesis, Theorem~\ref{thm:alex-completed} 
reduces to a result first proved in \cite[Thm.~5.6]{DPS-duke}.

\begin{corollary}[\cite{DPS-duke}]
\label{cor:linalex-c}
If $G$ is $1$-formal, then the completed Alexander invariant of $G$
is naturally isomorphic to the completed infinitesimal Alexander invariant
of its holonomy Lie algebra:
\[
\widehat{B_1(G;\k)} \cong \widehat{\B(\h(G;\k))}.
\]
\end{corollary}

%%%%%%%%%%%%%%%%%%%%%
\subsection{Chen rank inequalities and $1$-formality}
\label{subsec:chen-ineq}
%%%%%%%%%%%%%%%%%%%%%

We now exploit Theorem~\ref{thm:alex-completed}  to obtain a general 
comparison between Chen ranks and holonomy Chen ranks, leading 
to a new criterion for $1$-formality via comparison with holonomy Chen ranks.

For a finitely generated group $G$, the {\em holonomy Chen ranks} 
$\bar\theta_n(G)$ are defined in terms of the holonomy Lie algebra
$\h=\h(G;\k)$ as 
\begin{equation}
\label{eq:bar-theta}
\bar\theta_n(G)=\theta_n(\h)=\dim_{\k} (\h/\h'')_n. 
\end{equation}
By contrast, the Chen ranks $\theta_n(G)$ depend on the maximal
metabelian quotient of $G$. In general, these two sequences need not agree.

The next theorem shows that, for groups admitting a $1$-finite $1$-model,
the discrepancy between $\theta_n(G)$ and $\bar\theta_n(G)$ is governed
entirely by the Koszul module of the model.

\begin{theorem}
\label{thm:chen-ineq-model}
Let $G$ be a finitely generated group admitting a $1$-finite $1$-model
$(A,d)$ over a field $\k$ of characteristic $0$.
Then the following hold.

\begin{enumerate}[itemsep=3pt]
\item \label{item:chen-ineq}
For all $n\ge1$,
\[
\theta_n(G)\le \bar\theta_n(G).
\]

\item \label{item:chen-eq-formal}
If, in addition, $G$ is $1$-formal, then
\[
\theta_n(G)=\bar\theta_n(G), \quad \text{for all } n\ge1.
\]
\end{enumerate}
\end{theorem}

\begin{proof}
Note first that
\[
\theta_1(G)=\theta_1(A)=\bar\theta_1(G)=\dim_\k H_1(G;\k),
\]
so it suffices to consider Chen ranks in degrees $\ge2$.

Let $(A,d)$ be a $1$-finite $1$-model for $G$.
By Theorem~\ref{thm:alex-completed}\eqref{ch3}, the Chen ranks of $G$ are
computed by the Koszul module of the model:
\[
\theta_n(G)=\theta_n(A), \quad \text{for all } n\ge1,
\]
and for $k\ge0$,
\[
\sum_{n\ge0} \theta_{n+2}(G)\, t^n
= \Hilb\bigl(\gr^\m(\B_1(A)),t\bigr).
\]

On the other hand, by Proposition~\ref{prop:holo-massey}, the generating series
for the holonomy Chen ranks satisfies
\[
\sum_{n\ge0} \bar\theta_{n+2}(G)\, t^n
=
\sum_{n\ge0} \bar\theta_{n+2}(A)\, t^n
= \Hilb\bigl(\gr(\B_1(H^\ast(A))),t\bigr).
\]

Applying the Hilbert series inequality of
Theorem~\ref{thm:Hilb-ineq-gr} to $(A,d)$---stated there for $\gr^F$, which
by Lemma~\ref{lem:filtration-shift-B1} differs from $\gr^\m$ by a shift of one,
uniformly on both sides---yields
\[
\Hilb\bigl(\gr^\m(\B_1(A));t\bigr)
\preccurlyeq
\Hilb\bigl(\gr(\B_1(H^\ast(A)));t\bigr),
\]
coefficientwise. Comparing coefficients gives
$\theta_n(G)\le \bar\theta_n(G)$, 
for all $n\ge1$, proving~\eqref{item:chen-ineq}.

If $G$ is $1$-formal, then $H^\ast(G;\k)$ with zero differential is a
$1$-finite $1$-model for $G$. By Theorem~\ref{thm:alex-completed}\eqref{ch3},
the Chen ranks $\theta_n(G)$ are independent of the choice of $1$-finite
$1$-model. Hence, in the $1$-formal case, we may compute the Chen ranks
using the model $A = H^\ast(G;\k)$.

For this model, Proposition~\ref{prop:holo-massey} identifies the generating
series of Chen ranks (with a shift of $2$) with the Hilbert series of
$\gr(\B_1(H^\ast(G;\k)))$, which coincides with the generating series of
holonomy Chen ranks. This yields $\theta_n(G) = \bar\theta_n(G)$, for all $n\ge 1$, 
proving~\eqref{item:chen-eq-formal}.
\end{proof}

\begin{remark}
\label{rem:theta-ineq}
Theorem~\ref{thm:ps-sw} implies that, for every finitely generated group $G$,
one has $\theta_n(G)\le \bar\theta_n(G)$ for all $n\ge1$, with equality in all 
degrees when $G$ is $1$-formal.  
Theorem~\ref{thm:chen-ineq-model} recovers this inequality for groups admitting 
a $1$-finite $1$-model via Koszul modules and Hilbert series, thereby providing a 
new proof and explicitly identifying the algebraic mechanism underlying both 
equality and strict inequality.
\end{remark}

The next example illustrates the necessity of using a $1$-finite $1$-model for a group $G$ 
in Theorem~\ref{thm:alex-completed}: replacing such a model by the cohomology algebra 
$H^*(G;\k)$ may destroy all information about Chen ranks, even for nilpotent groups.

\begin{example}
\label{ex:heis-chen}
Let $G=F_2/\gamma_3(F_2)$ be the integral Heisenberg group, and consider 
the $1$-finite $1$-model $A=(\bwedge(a_1,a_2,a_3),d)$ with $d a_3=a_1\wedge a_2$ 
and $d a_1=d a_2=0$.  As computed in Example~\ref{ex:heis-bis}, 
\[
\h(A)=\Lie(x_1,x_2)/\langle[x_1,[x_1,x_2]],\,[x_2,[x_1,x_2]]\rangle, 
\]
where $x_i=a_i^\vee$ are dual to the degree-one generators of $A$.
Thus $\h(A)/\h'(A)=\k\langle x_1,x_2\rangle$, 
$\h'(A)/\h''(A)=\k\langle[x_1,x_2]\rangle$, and the holonomy Chen ranks are
$\theta_1(A)=2$, $\theta_2(A)=1$, and $\theta_n(A)=0$ for $n\ge3$.
By contrast, $H^1(A)=\k\langle a_1,a_2\rangle$ and the cup product on $H^1(A)$ is trivial; 
hence $\h(G)=\h(H^*(G;\k))\cong\Lie(x_1,x_2)$, with
$\dim \h(G)'/\h(G)''=\infty$.  
Hence $\h(G)$ fails to correctly recover the Chen ranks, showing that neither $G$ nor
the corresponding nilmanifold is $1$-formal.
\end{example}

%%%%%%%%%%%%%%%%%%%%%%%%%%%%%
\subsection{Vanishing resonance for groups admitting a $1$-finite $1$-model}
\label{subsec:res-zero-model}
%%%%%%%%%%%%%%%%%%%%%%%%%%%%%

This subsection treats the extreme case where the first resonance variety is trivial. 
Although this case lies at one end of the spectrum from the perspective of the 
Chen ranks conjecture, it already exhibits the key mechanism---finite length 
of Koszul modules---that governs asymptotic Chen ranks behavior. 
The results of this subsection are based on the general principle established
in Proposition~\ref{prop:no-hidden-res}, which relates the vanishing of
resonance support near the origin, for a finite-type $\cdga$ model, to finite
length of the associated graded Koszul modules, under a mild spectral sequence
collapse condition.

We extend vanishing results known for $1$-formal groups to all finitely
generated groups admitting a $1$-finite $1$-model. By Proposition~\ref{prop:r1a}, 
in degree $1$ the resonance varieties $\RR^{1,s}$ coincide with the corresponding 
support loci $\RR_{1,s}$ away from the origin, for every $s\ge1$. Hence, the following statements 
may equivalently be formulated in terms of resonance varieties.

\begin{corollary}
\label{cor:zero-res-model}
Let $G$ be a finitely generated group, and suppose $G$ admits a
$1$-finite $1$-model $(A,d_A)$ over $\k$.  Set $S=\Sym(H_1(G;\k))$. 
Then for each $s\ge1$, the following are equivalent:
\begin{enumerate}[itemsep=2pt]
\item \label{zz1}
The origin is isolated in $\RR_{1,s}(G)$, or does not lie on it.
\item \label{zz2}
The $S$-module $\bigl(\bigwedge^s \B_1(A)\bigr)_\m$ has finite length over
$S_\m$; equivalently, $\gr\bigl(\bigwedge^s\B_1(A)\bigr)$ is
finite-dimensional over $\k$.
\end{enumerate}
If, in addition, $\RR_{1,s}(A)$ is conical, then \eqref{zz1} and
\eqref{zz2} are equivalent to $\RR_{1,s}(G)$ being empty or equal to
$\{0\}$, and to $\bigwedge^s\B_1(A)$ being finite-dimensional over $\k$.
\end{corollary}

\begin{proof}
By Theorem~\ref{thm:alex-completed}, there is an 
isomorphism of associated graded $S$-modules
\begin{equation}
\label{eq:grb1-GA}
\gr^I(B_1(G;\k))\cong \gr^\m(\B_1(A)).
\end{equation}
Since for any finitely generated $S$-module $M$, the support of the
associated graded module satisfies $\Supp(\gr M)=\bV(\ini_{\m}(\Ann M))$, 
where $\m\subset S$ denotes the maximal ideal corresponding to the origin,
while the tangent cone at the origin of $\Supp M=\bV(\Ann M)$ is defined
by the same initial ideal, it follows that
\[
\TC_0(\Supp M)=\Supp(\gr M).
\]
The isomorphism \eqref{eq:grb1-GA} implies that the tangent cones 
at the origin of the corresponding support varieties agree,
\begin{equation}
\label{eq:tc0-r1}
\TC_0 \RR_{1,s}(G)=\TC_0 \RR_{1,s}(A).
\end{equation}

For a finitely generated $S$-module $M$, the condition
$\TC_0(\Supp M)\subseteq\{0\}$ says that the origin is isolated in, or
absent from, $\Supp M$; equivalently, $\gr M$ has finite length, and
equivalently $M_\m$ has finite length over $S_\m$.  It does \emph{not}\/
imply that $M$ itself is finite-dimensional over $\k$, since components of
$\Supp M$ missing the origin are invisible to $\gr M$
\textup{(}Remark~\ref{rem:no-hidden-res-global}\textup{)}.  Applying this
criterion to $M=\bigwedge^s\B_1(A)$, and using \eqref{eq:tc0-r1}, yields
the equivalence of \eqref{zz1} and \eqref{zz2}.

If $\RR_{1,s}(A)$ is conical, then it coincides with its tangent cone at
the origin, so \eqref{zz1} upgrades to the vanishing of $\RR_{1,s}(A)$
away from the origin, and $\Supp_S\bigl(\bigwedge^s\B_1(A)\bigr)
\subseteq\{0\}$ forces $\bigwedge^s\B_1(A)$ to have finite length.
\end{proof}

\begin{corollary}
\label{cor:trivial-res-chen-model}
Let $G$ be a finitely generated group, and suppose $G$ admits a
$1$-finite $1$-model $(A,d_A)$. 
If either $\RR_1(A)=\{0\}$ or $\RR_1(G;\k)=\{0\}$, then
$\theta_n(G)=0$ for $n\gg0$.
\end{corollary}

\begin{proof}
Theorem~\ref{thm:alex-completed}\eqref{ch3} gives
$\theta_n(G)=\dim_\k \gr^\m_{n-2}(\B_1(A))$ for $n\ge2$. 
In either case, Corollary~\ref{cor:zero-res-model} (with $s=1$) implies that
$\gr(\B_1(A))$ is finite-dimensional over $\k$, so its graded pieces vanish
in large degrees.  Hence, $\theta_n(G)=0$ for $n\gg0$.
\end{proof}

The arguments rely on the comparison 
$\gr^I(B_1(G;\k))\cong\gr^\m(\B_1(A))$, 
which in turn requires that $A$ be $1$-finite.  
Without a $1$-finite $1$-model there is, in general, no systematic way to compare the
completed Alexander invariants of $G$ with the Koszul modules arising from
algebraic models.

The following example shows that the two hypotheses in
Corollary~\ref{cor:trivial-res-chen-model} are genuinely different,
even when a $1$-finite $1$-model exists.

\begin{example}
\label{ex:heis-zero-res}
Let $G=F_2/\gamma_3(F_2)$ be the integral Heisenberg group, and let
$A=(\bigwedge(a_1,a_2,a_3),d)$ be the $1$-finite $1$-model with
$d a_3=a_1\wedge a_2$ and $d a_1=d a_2=0$.
A straightforward computation of the Aomoto complexes
$(A,\delta_a)$ shows that $H^1(A,\delta_a)=0$ for all $a\neq 0$.
Hence, the first resonance variety of the model is trivial:
$\RR^1(A)=\{0\}$.
On the other hand, $H^1(G;\k)=\k\langle a_1,a_2\rangle$, and the cup
product on $H^1(G;\k)$ vanishes.
Therefore, $\RR^1(G;\k)=H^1(G;\k)$.

This example shows that trivial resonance of a $1$-model does not imply
trivial resonance of the cohomology algebra.
Nevertheless, Corollary~\ref{cor:trivial-res-chen-model} applies,
since the asymptotic Chen ranks are governed by the Koszul module
$\B_1(A)$, rather than by $\RR^1(G;\k)$.
\end{example}

\begin{remark}
\label{rem:tcone-0res}
For $1$-formal groups, vanishing of $\RR^1(G;\C)$ is equivalent to the tangent 
cone at $1$ of the characteristic variety being trivial, by Theorem~\ref{thm:dps-dp}.
The results of this subsection show that, in the presence of a $1$-finite $1$-model, the 
model resonance $\RR^1(A)$---rather than $\RR^1(H^*(G;\C))$---is the relevant object 
controlling finiteness of Koszul modules and eventual vanishing of Chen ranks, 
even when the classical tangent cone formula fails.

In the absence of a  $1$-finite $1$-model, vanishing of $\RR^1(G)$ still implies finiteness 
of the infinitesimal Alexander invariant $\B_1(G)$, and---under additional formality 
assumptions---of the completed Alexander invariant $\widehat{B_1(G)}$; see \cite{DP-ann}.
The approach taken here replaces resonance of the cohomology algebra by resonance 
of a finite model, allowing one to recover analogous Chen ranks consequences 
without assuming $1$-formality.
\end{remark}

Koszul modules with vanishing resonance form a rigid class of algebraic objects
with strong finiteness properties. They arise in several foundational problems of
commutative algebra and algebraic geometry, notably Green’s conjecture on syzygies
of canonical curves and, more generally, in the study of syzygies of vector bundles
on algebraic varieties \cite{AFPRW19, AFPRW22}. In a different direction, it plays
an important role in the study of Torelli groups and Johnson filtrations
\cite{DP-ann, PS-jtop}.

The results of this subsection show that vanishing resonance is equivalent to
finite length of the associated Koszul modules. This finite-length property is the 
key mechanism underlying the eventual vanishing of Chen ranks.
See Remark~\ref{rem:chen-zero-res-context} for a more general 
perspective when resonance is nontrivial.

%%%%%%%%%%%%%%%%%%%%%%%%%%%%%
\subsection{Linear resonance, separability, and isotropicity}
\label{subsec:separable-isotropic}
%%%%%%%%%%%%%%%%%%%%%%%%%%%%%

We now recall the structural conditions on linear components of the first 
resonance variety that enter the Chen ranks formulas discussed in the 
next subsection. These conditions control both the scheme-theoretic 
geometry of resonance and the algebraic structure of the associated 
infinitesimal Alexander invariants.

Let $V=H^1(G;\k)$ and let $K\subseteq \bwedge^2 V$ be the kernel 
of the cup-product map $\bwedge^2 V\to H^2(G;\k)$.
A linear subspace $L\subseteq V$ is called \emph{isotropic}\/ if 
$L\wedge L\subseteq K$, that is, if the cup product vanishes identically on $L$.
It is called \emph{separable}\/ if $(L\wedge V)\cap K\subseteq L\wedge L$.
A subspace $L$ is called \emph{strongly isotropic}\/ if it is both isotropic 
and separable, or equivalently, if $(L\wedge V)\cap K = L\wedge L$. 
For instance, if $K=0$, then every subspace is separable, but the only isotropic
subspaces are those of dimension at most $1$. On the other hand, if $K= \bwedge^2 V$, 
then all subspaces are isotropic, but no proper, non-trivial subspace is separable. 
These examples emphasize that isotropicity is an intrinsic condition on $L$, 
whereas separability reflects how $L$ interacts with the ambient space $V$.

The relevance of these notions to the geometry of the resonance varieties 
was clarified in \cite{AFRS24}.
Let $G$ be a $1$-formal group, and let $\RR^1(G;\k)$ be its first resonance 
variety, with the scheme structure given by the annihilator of $\B_1(G;\k)$. 
By Theorem~\ref{thm:dps-dp}\eqref{tc2}, each irreducible component 
of $\RR^1(G;\k)$ is a (rationally defined) linear subspace. 
Separability alone already forces strong scheme-theoretic properties: each 
separable irreducible component of the projectivized resonance variety is 
reduced and isolated \cite[Thm.~4.5]{AFRS24}. Consequently, if all components 
are separable, the resonance scheme is reduced and projectively 
disjoint \cite[Cor.~4.6]{AFRS24}.

Isotropicity enters in a different way. For an isotropic component $L$ of resonance, 
reducedness of the corresponding projective component is equivalent to strong 
isotropicity \cite[Thm.~5.1]{AFRS24}. While this equivalence is not required for 
Chen ranks formulas themselves, it explains why strong isotropicity naturally 
appears when one seeks converse or rigidity statements. These conditions 
are discussed from the algebro-geometric perspective of the graded Koszul 
module $W(V,K)$ in~\cite{Farkas}.

%%%%%%%%%%%%%%%%%%%%%%%%%%%%%
\subsection{The Chen ranks conjecture and its algebraic refinements}
\label{subsec:chen-ranks-conj}
%%%%%%%%%%%%%%%%%%%%%%%%%%%%%

To place the following results in context, we recall the original form of the
Chen ranks conjecture, introduced in \cite{Su-conm01} for fundamental groups
of complex hyperplane arrangements.
The conjecture predicts that, in sufficiently high degrees, the Chen ranks
of an arrangement group $G$ are determined by the geometry of the first
resonance variety $\RR^1(G;\C)\subset H^1(G;\C)$. Arrangement groups 
are $1$-formal, and hence, by Theorem~\ref{thm:dps-dp},
their resonance varieties decompose as finite unions of linear subspaces.
In this setting, the conjecture asserts that the asymptotic behavior of the
Chen ranks $\theta_n(G)$ depends only on the dimensions of these components.

For a general $1$-formal group $G$, the relevance of resonance to Chen ranks
is explained by the comparison between the holonomy Lie algebra and the
associated graded Lie algebra of $G/G''$.
By Theorem~\ref{thm:ps-sw}, there is a natural epimorphism
\begin{equation}
\label{eq:phi-2-epi}
\h(G;\k)/\h(G;\k)'' \longsurj \gr(G/G'';\k)
\end{equation}
which is an isomorphism in the $1$-formal case.
As a consequence, the Chen ranks of a $1$-formal group are governed by the
infinitesimal Alexander invariant, whose structure is reflected in the
geometry of $\RR^1(G;\k)$, subject to additional structural constraints 
on its components.

An algebraic form of the Chen ranks conjecture for $1$-formal 
groups was established in \cite{AFRS24}, refining earlier work of 
Cohen--Schenck \cite{Cohen-Schenck15}. The latter required the 
irreducible components of $\RR^1(G;\k)$ to be isotropic, reduced,
and projectively disjoint---conditions arising naturally in the 
quasi-projective setting via \cite{DPS-duke}.

A key observation of \cite{AFRS24} is that separability alone already 
implies both reducedness of the resonance scheme and projective disjointness 
of its components, and hence quasi-projectivity is not required to state or 
prove a Chen ranks formula. The subsequent paper \cite{AFRS25} strengthens 
this result by identifying the effective asymptotic range $n\ge b_1(G)-1$ in 
which the Chen ranks formula holds, a substantially more delicate problem.

\begin{theorem}[\cite{AFRS24, AFRS25}]
\label{thm:chen-afrs}
Let $G$ be a $1$-formal group, and assume its resonance variety $\RR^1(G;\k)$ 
is strongly isotropic. If $h_m$ denotes the number of $m$-dimensional components 
of $\RR^1(G;\k)$, then
\begin{equation}
\label{eq:chen-ranks-formula} 
\theta_n(G)=\sum_{m\ge 2} h_m\cdot \theta_n(F_m),
\quad \text{for all } n\ge b_1(G)-1,
\end{equation}
where $\theta_n(F_m)=(n-1)\binom{m+n-2}{n}$ for $n\ge 2$.
\end{theorem}

\begin{remark}
\label{rem:chen-other}
In the absence of strong isotropicity, Chen ranks formulas of the form
\eqref{eq:chen-ranks-formula} need not hold, as shown for instance in
\cite{PS-mathann} for right-angled Artin groups and in
\cite{Cohen-Schenck15, SW-aam} for the upper McCool groups.

A key observation of \cite{AFRS24} is that separability alone already
implies that each irreducible component of $\RR^1(G;\k)$ is reduced and
projectively disjoint \cite[Thm.~4.5]{AFRS24}, so the Chen ranks formula
does not require quasi-projectivity once these conditions are satisfied.
The subsequent paper \cite{AFRS25} refines this by identifying the
effective asymptotic range $n \ge b_1(G)-1$ in which the formula holds.

K\"{a}hler groups provide a different geometric regime in which separability
is known in special cases (for instance, fundamental groups of complete
intersection Kodaira fibrations), while isotropicity fails (except in the
vanishing resonance case). In this setting, Chen ranks
are conjecturally governed by a formula of type
\eqref{eq:chen-ranks-formula}, with free groups replaced by surface groups.
We refer to \cite{AFRS24} for the precise conjecture and for its verification
in the case of compact algebraic surfaces admitting two independent Kodaira
fibrations.
\end{remark}

The discussion above takes place entirely within the $1$-formal setting,
where resonance is a purely cohomological invariant and the Chen ranks are
controlled by the holonomy Lie algebra.
Beyond formality, resonance acquires genuinely differential features,
and the relationship between Chen ranks and resonance becomes more subtle.

For groups admitting a $1$-finite $1$-model $(A,d)$,
Theorem~\ref{thm:alex-completed} shows that the Chen ranks are determined
exactly by the Koszul module $\B_1(A)$.
From this perspective, vanishing resonance corresponds to finite length of
$\B_1(A)$, while the presence of linear resonance suggests that
$\B_1(A)$ should admit a natural comparison map to the direct sum of the
Koszul modules associated to the individual resonance components.

In the $1$-formal case, strong isotropicity of these components
provides precisely the condition needed to ensure that this map is
surjective with finite-dimensional kernel.
This indicates that any extension of the Chen ranks conjecture beyond
formality must retain the assumption of linear resonance, while replacing 
strong isotropicity by a genuinely differential condition on $(A,d)$ involving 
both the multiplication $\mu\colon A^1\wedge A^1\to A^2$ and the differential
$d\colon A^1\to A^2$, which together control the associated Koszul modules.

\begin{remark}
\label{rem:chen-zero-res-context}
The vanishing-resonance case $\RR^1(A)=\{0\}$ corresponds to the situation
where the Koszul module $\B_1(A)$ has finite length; see
Section~\ref{subsec:res-zero-model} for the resulting vanishing of Chen ranks. 
This regime plays an important role in several contexts, including Green’s
conjecture on syzygies of canonical curves \cite{AFPRW19, AFPRW22} and the
study of Torelli groups and Johnson filtrations \cite{DP-ann, PS-jtop}. 
From this perspective, the Chen ranks conjecture can be viewed as a
refinement of the finite-length phenomenon to the case where
$\RR^1(A)$ is nontrivial but admits a linear decomposition structure.
\end{remark}

Guided by the above principle, one is led to seek a differential
condition on $(A,d)$ that governs the structure of $\B_1(A)$
in a manner analogous to strong isotropicity in the $1$-formal case. 
The correct analogue, suggested by the theory of decomposable
arrangements \cite{PS-imrn04, Su-decomp} and illustrated by the
pure elliptic braid groups in \S\ref{subsec:pn-p1n-comparison},
is an exact decomposition of $\B_1(A)$---rather than merely
an asymptotic one---indexed by the resonance components.

\begin{question}
\label{q:chen-nonformal}
Let $G$ be a group admitting a $1$-finite $1$-model $(A,d)$ with positive weights,
and assume that the first resonance variety of $A$ decomposes as a finite union of
linear subspaces, $\RR^1(A) = \bigcup_{\alpha \in \mathcal{C}} L_\alpha$. 
For each component $L_\alpha$, set $m_\alpha \coloneqq \dim L_\alpha$.

When does the Koszul module $\B_1(A)$ admit a decomposition
\[
\B_1(A) \cong \bigoplus_{\alpha \in \mathcal{C}} M_\alpha
\quad \text{as $S$-modules},
\]
where each summand satisfies $\Hilb(M_\alpha,t) = (1-t)^{-m_\alpha}$? 

In this situation, does one obtain an exact Chen ranks formula
\[
\theta_k(G) = \sum_{\alpha \in \mathcal{C}} \theta_k(F_{m_\alpha})
\quad \text{for all } k \ge 2,
\]
without any asymptotic range restriction?
\end{question}

This is the non-formal analogue of the decomposability condition
for hyperplane arrangements studied in \cite{PS-imrn04, Su-decomp}:
for a decomposable arrangement $\A$, the infinitesimal
Alexander invariant $\B_1(M(\A))$ decomposes exactly as a
direct sum of local contributions indexed by the $2$-flats, and the
same Chen ranks formula holds for all $k\ge 2$
\cite{CS-tams, PS-cmh06, Su-decomp}.
The pure elliptic braid groups $P_{1,n}$ provide a family of
non-$1$-formal groups for which this exact decomposition holds,
with $h_2 = \binom{n}{2}$ components of dimension $2$ in
$\RR^1(A(n))$ and $\theta_k(P_{1,n}) = \binom{n}{2}\cdot\theta_k(F_2)$
for all $k\ge 2$; see Section~\ref{subsec:pn-p1n-comparison}.
The relevant condition is the block-diagonal structure of the
Bibby differential, which forces the exact decomposition of
Theorem~\ref{thm:BA-n} with zero kernel.
Identifying the precise differential-algebraic condition on 
$(A,d)$---playing the role of decomposability in the non-formal setting,
and formulated at the level of $\B_1(A)$ rather than $\RR^1(H^*(A))$---%
remains an open problem.

%%%%%%%%%%%%%%%%%%%%%
\subsection{Koszul spectral obstructions to $1$-formality}
\label{subsec:koszulss-1formal}
%%%%%%%%%%%%%%%%%%%%%

The purpose of this subsection is to isolate a family of obstructions
to $1$-formality arising from the interaction between positive
weights and the Koszul homology of a finite-type model.
The mechanism is the comparison, provided by
Lemma~\ref{lem:holo-to-cohomology}, between the holonomy Lie algebra of a
model and that of its cohomology algebra: when the kernel of this comparison
is not contained in the second derived subalgebra, the Koszul module of the
model is strictly smaller than the one predicted by cohomology, and the
resulting drop in Chen ranks obstructs $1$-formality.  The obstruction is 
computable in terms of Koszul modules and Chen ranks, and is particularly 
effective for groups admitting $1$-finite models.

We first record how the functor $\B(-)$ behaves with respect to the
comparison map $\Xi_A$.

\begin{lemma}
\label{lem:ker-B-Xi}
Let $\Xi\colon\mathfrak p\surj\mathfrak q$ be an epimorphism of graded Lie
algebras which is an isomorphism in degree $1$, and set
$\mathfrak k=\ker\Xi$.  Then $\mathfrak k\subseteq\mathfrak p'$, the algebras
$\mathfrak p$ and $\mathfrak q$ have the same abelianization, and the induced
map of infinitesimal Alexander invariants is the natural projection
\[
\B(\Xi)\colon\ \mathfrak p'/\mathfrak p''\longsurj
\mathfrak p'/(\mathfrak p''+\mathfrak k),
\qquad
\ker\B(\Xi) \cong \mathfrak k/(\mathfrak k\cap\mathfrak p'') .
\]
In particular, $\B(\Xi)$ is injective if and only if
$\mathfrak k\subseteq\mathfrak p''$; non-injectivity of $\B(\Xi)$ implies
that of $\Xi$, but not conversely.
\end{lemma}

\begin{proof}
Since $\Xi$ is an isomorphism in degree $1$ and $\mathfrak p$ is generated in
degree~$1$, we have $\mathfrak k\subseteq\mathfrak p'$, whence
$\mathfrak q^{\mathrm{ab}}=\mathfrak p^{\mathrm{ab}}$ and both invariants are
modules over the same symmetric algebra.  Surjectivity of $\Xi$ gives
$\mathfrak q'=\Xi(\mathfrak p')=\mathfrak p'/\mathfrak k$ and
$\mathfrak q''=\Xi(\mathfrak p'')=(\mathfrak p''+\mathfrak k)/\mathfrak k$, so
that
\[
\B(\mathfrak q)=\mathfrak q'/\mathfrak q''
=\bigl(\mathfrak p'/\mathfrak k\bigr)\big/
\bigl((\mathfrak p''+\mathfrak k)/\mathfrak k\bigr)
\cong\mathfrak p'/(\mathfrak p''+\mathfrak k).
\]
Under this identification $\B(\Xi)$ is the natural projection, with kernel
$(\mathfrak p''+\mathfrak k)/\mathfrak p''\cong
\mathfrak k/(\mathfrak k\cap\mathfrak p'')$ by the second isomorphism
theorem.  The last assertion follows, the failure of the converse occurring
precisely when $0\ne\mathfrak k\subseteq\mathfrak p''$.
\end{proof}

\begin{theorem}
\label{thm:koszulss-1formal-obstruction}
Let $G$ be a finitely generated group admitting a $1$-finite $1$-model
$(A,d)$ over $\k$ with positive weights and $\im(d^{\vee})$ concentrated 
in weights $\ge 2$, and let 
$\Xi_A\colon \h(H^*(A))\surj\h(A)$ be the epimorphism of 
Lemma~\ref{lem:holo-to-cohomology}.  Assume that the induced epimorphism
\begin{equation}
\label{eq:b1ha-b1a}
\B_1(H^*(A)) \cong \B(\h(H^*(A))) \longsurj \B(\h(A)) \cong \B_1(A)
\end{equation}
is \emph{not} injective---equivalently, by Lemma~\ref{lem:ker-B-Xi}, assume
that $\ker\Xi_A\not\subseteq\h(H^*(A))''$.  Then:
\begin{enumerate}[itemsep=1pt]
\item \label{x1}
$\gr(\B_1(G;\k))$ is a proper quotient of $\B(\h(G;\k))$;
\item  \label{x2}
$\theta_n(G)<\bar\theta_n(G)$ for at least one $n$;
\item  \label{x3}
$G$ is not $1$-formal.
\end{enumerate}
\end{theorem}

\begin{proof}
Set $S=\Sym(H_1(G;\k))$; the map \eqref{eq:b1ha-b1a} is one of $S$-modules,
by Theorem~\ref{thm:B-holo} and Lemma~\ref{lem:ker-B-Xi}.
Under the identifications $\B_1(H^*(A))\cong\B(\h(G;\k))$---valid because the
holonomy Lie algebra depends only on the cup-product map
$\bwedge^2H^1\to H^2$---and $\gr(\B_1(G;\k))\cong\gr(\B_1(A))$ of
Theorem~\ref{thm:alex-completed}\eqref{ch2}, the map \eqref{eq:b1ha-b1a}
induces a natural surjection
\[
\B(\h(G;\k)) \longsurj \gr(\B_1(G;\k)) 
\]
with nonzero kernel.  This proves~\eqref{x1}.

Claim~\eqref{x2} follows by comparing Hilbert series.  By 
Theorem~\ref{thm:Hilb-ineq-gr},
\[
\Hilb(\gr^F \B_1(A),t) \preccurlyeq \Hilb(\B_1(H^*(A)),t),
\]
and by \eqref{x1} the inequality is strict in at least one degree; hence
\[
\Hilb\bigl(\gr^F \B_1(G;\k),t\bigr) \prec \Hilb\bigl(\B(\h(G;\k)),t\bigr).
\]
Since the Chen ranks and the holonomy Chen ranks are the coefficients of these
two series, shifted by $2$, we get $\theta_n(G)<\bar\theta_n(G)$ for at least
one $n$.  Note that a strict inequality of Hilbert series need not propagate
to all large $n$: the kernel may be concentrated in finitely many degrees.

Finally, \eqref{x3} follows from~\eqref{x2}: by
Theorem~\ref{thm:chen-ineq-model}\eqref{item:chen-eq-formal}, a $1$-formal
group admitting a $1$-finite $1$-model satisfies
$\theta_n(G)=\bar\theta_n(G)$ for all $n\ge1$.
\end{proof}

\begin{remark}
\label{rem:d1-not-obstruction}
It is the failure of injectivity of $\B(\Xi_A)$, and not the nonvanishing of
the differential $d_1$ of the weight spectral sequence, that is doing the work
here.  Indeed, by Theorem~\ref{thm:weight-ss} the differential $d_1$ is the
Koszul differential of $H^{*}(A)$, hence is determined by the cohomology
algebra alone, and is nonzero for most models, including formal ones.
Both sides of this distinction are exhibited by the Bibby models of elliptic
configuration spaces, in Section~\ref{subsec:spectral-weight}: the hypothesis
of Theorem~\ref{thm:koszulss-1formal-obstruction} holds for $\Conf(\E,n)$ with
$n\ge3$ (Proposition~\ref{prop:xi-An-noninjective}), yielding a second proof
that $P_{1,n}$ is not $1$-formal, whereas it fails for the graphic
configuration spaces of a triangle-free graph, where $d_1$ is nonzero all the
same (Remark~\ref{rem:graphic-conf-formality}).
\end{remark}

\begin{remark}
\label{rem:spectral-vs-massey}
The failure of injectivity in \eqref{eq:b1ha-b1a} reflects the presence of
essential higher-weight terms in the differential of a minimal model for $X$.
While this phenomenon is compatible with the existence of nontrivial
triple Massey products, the obstruction above is strictly stronger:
it detects non-$1$-formality at the level of Koszul modules and Chen ranks,
even in situations where Massey products alone may be inconclusive.
\end{remark}

%%%%%%%%%%%%%%%%%%%%%%%%%%%%%%%%%%%%%%%
\section{Functorial models and formality of maps}
\label{sect:functorial-models}
%%%%%%%%%%%%%%%%%%%%%%%%%%%%%%%%%%%%%%%%

In this section we study formality properties of maps and group homomorphisms
from the perspective of functorial algebraic models. While formality of spaces
is detected by the existence of suitable quasi-isomorphisms between $\cdga$ models,
formality of maps requires compatibility of such models with morphisms.

We show that, for a homomorphism admitting a $\cdga$ model, the holonomy Lie
algebra provides a natural Lie-theoretic framework for detecting
$1$-formality. The presence of such a model ensures compatibility between
Malcev completions, holonomy Lie algebras, and infinitesimal Alexander
invariants, leading to effective obstructions to formality of maps. Models
for all homomorphisms in a class are supplied by functorial families of
$1$-finite $1$-models, of which several occur naturally in geometry
(Example~\ref{ex:functorial-models}).

%%%%%%%%%%%%%%%%%%%%
\subsection{Formality of maps}
\label{subsec:formal-maps}
%%%%%%%%%%%%%%%%%%%%

We begin by recalling from Section~\ref{subsec:cdga-maps} the 
notion of formality for $\cdga$ morphisms and its
topological interpretation for continuous maps. Let $(A,d_A)$ and $(B,d_B)$ 
be two formal $\cdgas$. 
A morphism $\varphi\colon A\to B$ is formal if and only if the induced map 
between minimal models fits into a homotopy-commutative diagram
\begin{equation}
\label{eq:formal-morphism}
\begin{tikzcd}[column sep=24pt]
\M(A) \ar[r] \ar[d] & (H^{*}(A),0)\ar[d] \\
\M(B) \ar[r] & (H^{*}(B),0) ,
\end{tikzcd}
\end{equation}
where the horizontal arrows are quasi-isomorphisms.

Now let $X$ and $Y$ be two formal spaces. 
A map $f\colon X\to Y$ is formal (over $\Q$) if the induced morphism  
$\apl(f)\colon \apl(Y)\to \apl(X)$ is formal.  Equivalently, there exists a diagram
\begin{equation}
\label{eq:formal-map}
\begin{tikzcd}[column sep=24pt, row sep=22pt]
\apl(Y)  \ar{d}{\apl(f)}
&
\M(Y) \ar[swap, pos=.4]{l}{\rho_Y} 
\ar{d}{\M(f)}
\ar{r}{\psi_Y} 
& (H^{*}(Y,\Q), 0)\ar{d}{f^*}
\\
 \apl(X)
&\M(X)\ar[swap, pos=.4]{l}{\rho_X} \ar{r}{\psi_X}
& (H^{*}(X,\Q), 0) 
\end{tikzcd}
\end{equation}
which commutes up to $\cdga$ homotopy and whose horizontal maps are 
quasi-isomorphisms. 
Analogous definitions hold for $q$-formality and over any characteristic-$0$ field.

As shown in \cite[\S6, Main Thm.~(i)]{DGMS}, every compact connected
K\"{a}hler manifold is formal; moreover, by part~(ii) of the same theorem, a
holomorphic map between such manifolds is formal over $\R$.  
In general, a map between formal spaces need not be formal.  
For example, the Hopf fibration $S^3\to S^2$ is not formal, even though 
both spaces are formal.

%%%%%%%%%%%%%%%%%%%%
\subsection{Formality and filtered formality of group homomorphisms}
\label{subsec:formal-homomorphisms}
%%%%%%%%%%%%%%%%%%%%

We now specialize to maps between Eilenberg--MacLane spaces and formulate
formality conditions directly at the level of group homomorphisms. 

Let $G$ and $H$ be two (finitely generated) $1$-formal groups. 
We say that a homomorphism $\alpha\colon G\to H$ is {\em $1$-formal} 
if the induced map $f_{\alpha}\colon K(G,1)\to K(H,1)$ is $1$-formal. 
Equivalently, the induced morphism on minimal $1$-models fits into 
a homotopy-commutative diagram of $\cdga$ morphisms
\begin{equation}
\label{eq:formal-hom}
\begin{tikzcd}[column sep=24pt, row sep=22pt]
\M_1(H) \ar[d, "\M_1(f_\alpha)"] \ar[r] & (H^{*}(H;\Q),0)\ar[d, "\alpha^*"] \\
\M_1(G) \ar[r] & (H^{*}(G;\Q),0),
\end{tikzcd}
\end{equation}
where the horizontal arrows are $1$-quasi-isomorphisms.

Recall from \eqref{eq:m1g-ce} that the minimal $1$-models $\M_1(G)$ and $\M_1(H)$ 
identify with the completed Chevalley--Eilenberg complexes of $\m(G)$ and $\m(H)$. 
Under these identifications, the induced morphism $\M_1(f_\alpha)$ corresponds to 
$\widehat{\CE}(\m(\alpha))$. Thus, $1$-formality of $\alpha$ is equivalent to the 
compatibility of $\m(\alpha)$ with the  quadratic presentations of $\m(G)$ 
and $\m(H)$ induced by their holonomy Lie algebras. 

In the $1$-formal case, recall that the $1$-formality isomorphism 
$\Theta_G\colon \m(G) \isom \widehat{\h(G;\Q)}$ from \eqref{eq:1formal-malcev} 
is a filtered Lie algebra isomorphism, unique up to automorphisms acting 
trivially in degree~$1$. Equivalently, $\alpha$ is $1$-formal if and 
only if the following diagram commutes:
\begin{equation}
\label{eq:formal-map-malcev}
\begin{tikzcd}[column sep=24pt, row sep=22pt]
\m(G) \ar[d, "\m(\alpha)"] \ar[r, "\Theta_G"] 
& \widehat{\h(G;\Q)} \ar[d, "\widehat{\h(\alpha)}"] \\
\m(H) \ar[r, "\Theta_H"] 
& \widehat{\h(H;\Q)} ,
\end{tikzcd}
\end{equation}
where $\Theta_G$ and $\Theta_H$ are the filtered Lie algebra isomorphisms
determined by the horizontal $1$-quasi-isomorphisms of \eqref{eq:formal-hom}.
Since those $1$-quasi-isomorphisms are part of the data of
\eqref{eq:formal-hom}, and not fixed in advance, $1$-formality of $\alpha$
asserts the existence of \emph{some}\/ pair $(\Theta_G,\Theta_H)$ making
\eqref{eq:formal-map-malcev} commute.  Equivalently, $\m(\alpha)$ and
$\widehat{\h(\alpha)}$ are conjugate: the filtered morphism $\m(\alpha)$ can
be straightened into the graded morphism $\widehat{\h(\alpha)}$ induced by
$\alpha_*\colon H_1(G;\k)\to H_1(H;\k)$.  In
Section~\ref{subsec:rigid-formal-hom} we will need a variant in which the
pair $(\Theta_G,\Theta_H)$ is fixed beforehand; the two notions differ, and
Remark~\ref{rem:rigid-vs-classical} compares them.

The decomposition of $1$-formality into graded-formality and filtered-formality 
admits a natural extension to homomorphisms: graded-formality controls the 
behavior on associated graded Lie algebras, while filtered-formality controls 
compatibility at the Malcev level. 

A homomorphism $\alpha\colon G\to H$ between filtered-formal groups is 
{\em filtered-formal} if the induced Malcev map commutes with the 
filtered-formal isomorphisms:
\begin{equation}
\label{eq:filtered-formal-map}
\begin{tikzcd}[column sep=24pt, row sep=22pt]
\m(G;\k) \ar[d, "\m(\alpha)"] \ar[r, "\Theta_G^{\mathrm{ff}}"] 
& \widehat{\gr(G;\k)} \ar[d, "\widehat{\gr(\alpha)}"] \\
\m(H;\k) \ar[r, "\Theta_H^{\mathrm{ff}}"] 
& \widehat{\gr(H;\k)} .
\end{tikzcd}
\end{equation}

A homomorphism between $1$-formal groups is $1$-formal if and only if, 
under choices of filtered Lie algebra isomorphisms $\Theta_G$ and $\Theta_H$ 
as in \eqref{eq:1formal-malcev}, the induced diagram \eqref{eq:formal-map-malcev} 
becomes the filtered-formality diagram \eqref{eq:filtered-formal-map} after identifying 
holonomy and associated graded Lie algebras via graded-formality.

%%%%%%%%%%%%%%%%%%%%%%
\subsection{Formality relative to fixed formality data}
\label{subsec:rigid-formal-hom}
%%%%%%%%%%%%%%%%%%%%%%

The notions of Section~\ref{subsec:formal-homomorphisms} quantify existentially
over the isomorphisms $\Theta$, which are themselves only canonical up to
filtered automorphisms acting trivially on $H_1$
(see \eqref{eq:1formal-malcev}).  For the obstruction theory of
Sections~\ref{subsec:alex-groups} and~\ref{subsec:alex-secondary} we need the
variant in which those isomorphisms are fixed in advance.

\begin{definition}
\label{def:rigid-formal-hom}
Let $G$ and $H$ be $1$-formal groups, and fix $1$-formality isomorphisms
$\Theta_G$ and $\Theta_H$ as in \eqref{eq:1formal-malcev}, with
$\Theta_H=\Theta_G$ whenever $H=G$.  A homomorphism $\alpha\colon G\to H$ is
\emph{strictly $1$-formal}\/ (with respect to these data) if
\begin{equation}
\label{eq:rigid-formal}
\m(\alpha)=\Theta_H^{-1}\circ\widehat{\h(\alpha)}\circ\Theta_G ,
\end{equation}
that is, if \eqref{eq:formal-map-malcev} commutes for the chosen
$\Theta_G$ and $\Theta_H$.  Strict filtered-formality is defined in the same
way, using $\Theta^{\mathrm{ff}}$ and \eqref{eq:filtered-formal-map} in place
of $\Theta$ and \eqref{eq:formal-map-malcev}.
\end{definition}

Transporting the gradings of $\h(G;\k)$ and $\h(H;\k)$ along $\Theta_G$ and
$\Theta_H$, condition \eqref{eq:rigid-formal} says exactly that $\m(\alpha)$
is homogeneous of degree~$0$.  Thus $1$-formality of a group is the existence
of a compatible grading on its Malcev Lie algebra, and strict $1$-formality of
a homomorphism is the requirement that the induced map be graded for the
chosen gradings.

Definition~\ref{def:rigid-formal-hom} is the two-object case of the
requirement that a family of homomorphisms be \emph{uniformly modeled}\/ in
the sense of Definition~\ref{def:modeled-hom}: a single model is attached to
each group, once, and every member of the family is required to be compatible
with those choices.  For an endomorphism $\alpha$ of $G$, strict
$1$-formality is uniform modelling of the family $\{\id_G,\alpha\}$.

\begin{proposition}
\label{prop:strict-endo}
Let $G$ be a finitely generated, $1$-formal group, and let
$\alpha\colon G\to G$ be an endomorphism inducing the identity
on $H_1(G;\k)$.  Then:
\begin{enumerate}
\item \label{se1}
$\m(\alpha)$ is an automorphism of $\m(G;\k)$, and $\alpha$ is $1$-formal.
\item \label{se2}
$\alpha$ is strictly $1$-formal if and only if
\begin{equation}
\label{eq:rigid-endo}
\m(\alpha)=\id ,
\end{equation}
a condition that does not depend on the choice of $\Theta_G$.
\end{enumerate}
\end{proposition}

\begin{proof}
\eqref{se1}
Transport $\m(\alpha)$ along $\Theta_G$ to a filtered endomorphism of
$\widehat{\h}=\widehat{\h(G;\k)}$.  Its associated graded is an endomorphism
of the graded Lie algebra $\h(G;\k)$, which is generated in degree~$1$; in
that degree it is $\alpha_*=\id$, so $\gr(\m(\alpha))=\id$.  Since
$\widehat{\h}$ is complete, $\m(\alpha)$ is an automorphism.  Moreover
$\h(\alpha)=\id$, since $\h(\alpha)$ is induced by $\alpha_*$.  Hence
$\Theta_G\circ\m(\alpha)^{-1}$ is again a filtered isomorphism
$\m(G;\k)\isom\widehat{\h(G;\k)}$ acting trivially on degree~$1$ of the
associated graded, so it is an admissible choice in
\eqref{eq:1formal-malcev}; taking it as $\Theta_H$ makes
\eqref{eq:formal-map-malcev} commute.

\eqref{se2}
With $\Theta_H=\Theta_G$ and $\h(\alpha)=\id$, condition
\eqref{eq:rigid-formal} reads $\Theta_G=\Theta_G\circ\m(\alpha)$, which is
\eqref{eq:rigid-endo}.  Replacing $\Theta_G$ by $u\Theta_G$, with $u$ as in
\eqref{eq:1formal-malcev}, turns this into
$u\Theta_G=u\Theta_G\circ\m(\alpha)$, again \eqref{eq:rigid-endo}.
\end{proof}

\begin{remark}
\label{rem:rigid-vs-classical}
Strict $1$-formality implies $1$-formality, and $1$-formality is precisely
strict $1$-formality for \emph{some}\/ admissible pair
$(\Theta_G,\Theta_H)$.  The two notions differ, and they differ exactly when
one group occupies both ends of $\alpha$: by
Proposition~\ref{prop:strict-endo}, an endomorphism inducing the identity on
$H_1$ is always $1$-formal, so the classical notion carries no information
there, whereas strict $1$-formality asks that $\alpha$ act trivially on the
Malcev completion.  This is the range in which the obstructions of
Section~\ref{subsec:alex-groups} operate.

Outside it the classical notion is not vacuous.  For instance, let
$G=H=F_2=\langle x,y\rangle$ and let $\alpha(x)=[x,y]$, $\alpha(y)=1$.  Then
$\alpha_*=0$, so $\widehat{\h(\alpha)}=0$, whereas
$\m(\alpha)(X)=[X,Y]+\cdots\neq0$; no pair of isomorphisms conjugates a
nonzero morphism to zero, so $\alpha$ is not $1$-formal.  This is the
group-theoretic counterpart of the Hopf map.
\end{remark}

\begin{example}
\label{ex:ian}
The free group $F_n$ is $1$-formal.  
Let $\IA_n\subset \Aut(F_n)$ be the subgroup of automorphisms acting trivially 
on abelianization.  
For $\alpha\in \IA_n$, the induced map 
$\alpha^*\colon H^*(F_n;\Z)\to H^*(F_n;\Z)$ is the identity, and thus 
$\h(\alpha)=\id\colon \Lie_n\to \Lie_n$.  
However, the induced Malcev map need not be trivial.  
For instance, if $F_2=\langle x,y\rangle$ and 
$\alpha(x)=x[x,y]$, $\alpha(y)=y$, then 
$\m(\alpha)(x)=x+[x,y]$, so $\m(\alpha)\ne\id$.  Since $F_n$ is residually
torsion-free nilpotent, it embeds in its Malcev completion, so $\m(\alpha)=\id$
would force $\alpha=\id$.  Thus no nontrivial element of $\IA_n$ satisfies
\eqref{eq:rigid-endo}, so none is strictly $1$-formal---although, by the
preceding remark, all of them are $1$-formal.
\end{example}

\begin{example}
\label{ex:not-filtered-formal}
Let $G=F_2/\gamma_3(F_2)$ be the integral Heisenberg group, generated by
$x_1,x_2$ with central commutator $x_3=[x_1,x_2]$.  
Then $G$ is $2$-step nilpotent and filtered-formal (but not $1$-formal).  
Consider the automorphism
\[
\alpha(x_1)=x_1x_3,\qquad \alpha(x_2)=x_2,\qquad \alpha(x_3)=x_3 .
\]
Since $\gamma_3(G)=1$, the element $x_3$ maps to $0$ in $\gr(G)$, and 
so $\alpha$ induces the identity on $\gr(G)$ and hence on 
$\widehat{\gr(G)}$.

On the other hand, the induced Malcev map satisfies
\[
\m(\alpha)(y_1)=y_1+y_3,\qquad 
\m(\alpha)(y_2)=y_2,\qquad
\m(\alpha)(y_3)=y_3 ,
\]
where $y_i$ denotes the Malcev image of $x_i$.  
Thus $\m(\alpha)\neq\id$, so the diagram \eqref{eq:filtered-formal-map} does
not commute for the choice $\Theta^{\mathrm{ff}}_H=\Theta^{\mathrm{ff}}_G$.
Hence $\alpha$ is not strictly filtered-formal, even though $G$ itself is
filtered-formal.  For the same reason as above, it \emph{is}\/
filtered-formal in the sense of \eqref{eq:filtered-formal-map} with
independently chosen isomorphisms, since $\gr(\alpha)=\id$; strictness is
what the example detects.
\end{example}

%%%%%%%%%%%%%%%%%%%%%%
\subsection{Functorial families of $1$-models}
\label{subsec:functorial-formal-hom}
%%%%%%%%%%%%%%%%%%%%%%

The results of this section require, for a homomorphism $\alpha\colon G\to H$,
a $\cdga$ morphism modeling it in the sense of
Definition~\ref{def:modeled-hom}.  Producing such morphisms one at a time is
awkward, and the two obvious strategies pull against each other.  Canonical
functorial models are available in considerable generality---Sullivan's
$\apl$ on arbitrary spaces, the de~Rham algebra on smooth manifolds---but are
almost never $1$-finite.  Conversely, for a $1$-formal group $G$ the cohomology algebra
$(H^*(G;\k),0)$ is $1$-finite, yet $\alpha^{*}$ need not model $\alpha$: a map
between formal spaces need not be formal, as the Hopf map already shows, and
at the level of groups this failure occurs for $\IA$-automorphisms of free
groups (Example~\ref{ex:ian}).  What is wanted is a construction achieving
both at once, on a suitable class of objects and maps.

\begin{definition}
\label{def:functorial-family}
Let $\mathcal{C}$ be a category equipped with a functor
$\pi\colon\mathcal{C}\to \mathcal{G}\mathrm{rp}$ whose values are finitely
generated groups.  A \emph{functorial family of $1$-finite $1$-models}\/ on
$\mathcal{C}$ assigns to each object $c$ a $1$-finite $\k$-$\cdga$ $A_c$ which
is a $1$-model for $\pi(c)$, and to each morphism $u\colon c\to c'$ a $\cdga$
morphism $A(u)\colon A_{c'}\to A_c$, in such a way that $A(\id)=\id$ and
$A(vu)=A(u)A(v)$ for composable morphisms.
\end{definition}

Definition~\ref{def:functorial-family} asks only that $c\leadsto A_c$ be a
contravariant functor; it does not by itself require that $A(u)$ model
$\pi(u)$.  Whether it does is a separate question, and in the geometric
examples below it is the harder half.  When both hold, the family furnishes a
uniformly modeled class of homomorphisms, and the results of
Sections~\ref{subsec:alex-groups} and~\ref{subsec:alex-secondary} apply to all
of them at once.

\begin{example}
\label{ex:functorial-models}
Functorial families of $1$-finite $1$-models whose morphisms model the
corresponding homomorphisms arise in each of the following settings.
\begin{enumerate}[itemsep=3pt, label=\textup{(\textit{\alph*})},
                  ref=\textit{\alph*}]
\item \label{fmx-nilp}
\emph{Finitely generated nilpotent groups.}
Let $\mathcal{C}$ be the category of finitely generated torsion-free nilpotent
groups and all homomorphisms, with $\pi=\id$.  The Malcev Lie algebra
$\m(G;\k)$ is finite-dimensional, so $A_G=\CE(\m(G;\k))$ is a $1$-finite
$\cdga$; it is a $1$-model for $G$ by \eqref{eq:m1g-ce}, and in fact a full
model.  Both $G\leadsto\m(G;\k)$ and $\g\leadsto\CE(\g)$ are functorial, the
latter contravariantly, and $\CE(\m(\alpha))$ models $\alpha$.  This is the
family underlying Proposition~\ref{prop:alex-obstruction-example}, and the
source of class~\eqref{wpg1} in
Corollary~\ref{cor:weight-preserving-geometric}.

\item \label{fmx-kahler}
\emph{Compact K\"{a}hler manifolds and holomorphic maps.}
Let $\mathcal{C}$ have as objects pointed compact connected K\"{a}hler
manifolds and as morphisms pointed holomorphic maps, with $\pi=\pi_1$, and
set $A_M=\big(H^{*}(M;\k),0\big)$ for $\k=\R$ or $\C$.  These $\cdgas$ are
$1$-finite, and each $H^{*}(f)$ models $f$, by
\cite[\S6, Main Thm.~(ii)]{DGMS}.  Moreover, the zig-zag connecting the
de~Rham algebra of $M$ to $\big(H^{*}(M;\k),0\big)$ may be chosen naturally
with respect to holomorphic maps, so that the modelling is uniform across the
whole class \cite[Prop.~7.4]{PS-adv19}.  Since the differentials vanish, every
homomorphism arising this way is $1$-formal, and the obstructions below are
vacuous on this class.

\item \label{fmx-qproj}
\emph{Quasi-projective manifolds.}
Let $\mathcal{C}$ have as objects normal-crossing compactifications
$(\overline{M},D)$ of quasi-projective manifolds $M=\overline{M}\setminus D$,
and as morphisms the pointed regular morphisms of pairs whose restriction
induces an injection on $H^1$; put $\pi(\overline{M},D)=\pi_1(M)$.  Navarro
Aznar's mixed Hodge diagram $H(\overline{M},D)$ is functorial with respect to
regular morphisms of pairs, and the associated $\cdga$ $A_0(\overline{M},D)$
is equivalent to $\Omega_{\C}(M)$ naturally in the pair; passing to $E_1$
yields a functorial family whose morphisms model the corresponding
homomorphisms, by \cite[Prop.~9.1]{PS-adv19}.  These models are $1$-finite,
and are isomorphic as bigraded $\cdgas$ to the Gysin models of
Section~\ref{subsec:qproj-weights}; this is class~\eqref{wpg2} in
Corollary~\ref{cor:weight-preserving-geometric}.  Unlike \eqref{fmx-kahler},
the differential is in general nonzero, and the obstructions below have
content.

\item \label{fmx-dupont}
\emph{Arrangements in smooth projective varieties.}
For a hypersurface arrangement $L$ in a smooth projective variety $X$, Dupont
constructed an Orlik--Solomon model $M^{\bullet}(X,L)$ for $X\setminus L$ and
proved it functorial: a holomorphic map $\varphi\colon X\to X'$ with
$\varphi^{-1}(L')\subset L$ induces
$M^{\bullet}(\varphi)\colon M^{\bullet}(X',L')\to M^{\bullet}(X,L)$
\cite[Thm.~1.2]{Du}.  The Bibby model of
Section~\ref{subsec:conf-elliptic} is an instance.
\end{enumerate}
\end{example}

\begin{remark}
\label{rem:naturality-caveat}
In \eqref{fmx-qproj} and \eqref{fmx-dupont} the two halves of the requirement
should be kept apart.  Functoriality of the assignment on objects and
morphisms is established in both cases.  That the resulting $\cdga$ morphism
models the induced homomorphism is a statement about naturality of the
comparison with $\apl$, and it is available through Navarro Aznar's theory;
this is why \eqref{fmx-qproj} is formulated there rather than through the
Gysin model directly, it not being known whether the Gysin morphism agrees,
under the identification on objects, with the morphism supplied by that
theory \cite[Rem.~9.2]{PS-adv19}.  The distinction is not pedantic: Morgan's
original treatment asserted a functoriality statement of this kind that had to
be withdrawn \cite{Mo86}, while the results of \cite[\S\S9--10]{Mo} used here
make no claim of naturality and are correct as stated.  The correction
identifies the mechanism at work throughout this section: the space of
indecomposables of a minimal model is not a rigid invariant, in the sense
that homotopic maps need not induce the same map on it, whereas the
associated graded of its canonical filtration is.  This is why the
obstructions of Section~\ref{subsec:alex-groups} are formulated on
associated graded objects; compare
Proposition~\ref{prop:koszul-obstruction-groups}\eqref{ko1} and
Remark~\ref{rem:koszul-obstruction-separated}.
\end{remark}

\begin{remark}
\label{rem:cohomology-not-a-family}
It is the modelling requirement, and not functoriality on its own, that
carries the weight.  For a $1$-formal group $G$, the cohomology algebra
$H_G=\big(H^*(G;\k),0\big)$ is a $1$-finite $1$-model with
$\h(H_G)=\h(G;\k)$, and $G\leadsto H_G$, $\alpha\leadsto\alpha^{*}$ is
manifestly a functorial family in the sense of
Definition~\ref{def:functorial-family}, on the category of $1$-formal groups
and all homomorphisms.  Yet $\alpha^{*}$ models $\alpha$ precisely when
$\alpha$ is $1$-formal, by \eqref{eq:formal-hom}, and
Proposition~\ref{prop:alex-obstruction-example} exhibits an automorphism of a
$1$-formal group for which this fails.  Accordingly, $H_G$ is used below only
one group at a time.
\end{remark}

Throughout, we keep track of the covariance 
of the functors $G\leadsto \m(G)$, $G\leadsto \h(G)$, $G\leadsto B_1(G)$, 
and $G\leadsto \B_1(G)$ and 
of the contravariance of the functors $A\leadsto \h(A)$ and $A\leadsto \B_1(A)$. 
This variance bookkeeping ensures that the comparison maps of
Lemma~\ref{lem:alex-functorial}, and those constructed below, are compatible
with group homomorphisms.

\begin{theorem}
\label{thm:1-formal-hom-unified}
Let $G$ and $H$ be $1$-formal groups, and let $\alpha\colon G\to H$ be a
homomorphism admitting a model $A(\alpha)\colon A_H\to A_G$.
Then the following are equivalent:
\begin{enumerate}
\item \label{fun1}
The homomorphism $\alpha$ is\/ $1$-formal.
\item \label{fun2}
The diagram of filtered Lie algebras
\[
\begin{tikzcd}[column sep=26pt, row sep=22pt]
\widehat{\h}(A_G) \ar[r, "\cong"] \ar[d, "\widehat{\h}(A(\alpha))"]
& \m(G;\k) \ar[d, "\m(\alpha)"] \\
\widehat{\h}(A_H) \ar[r, "\cong"]
& \m(H;\k)
\end{tikzcd}
\]
commutes.
\item \label{fun3}
The diagram
\[
\begin{tikzcd}[column sep=26pt, row sep=22pt]
\widehat{\h}(A_G) \ar[d, "\widehat{\h}(A(\alpha))"] \ar[r, "\cong"]
&
\widehat{\h}(H^{*}(G;\k)) \ar[d, "\widehat{\h}(\alpha^*)"]
\\
\widehat{\h}(A_H) \ar[r, "\cong"]
&
\widehat{\h}(H^{*}(H;\k))
\end{tikzcd}
\]
commutes, where $\alpha^*$ is the induced map in cohomology.
\end{enumerate}
\end{theorem}

\begin{proof}
\eqref{fun1}$\Rightarrow$\eqref{fun2}. 
If $\alpha$ is $1$-formal, the induced map
$f_\alpha\colon K(G,1) \to K(H,1)$ is $1$-formal.
By functoriality of minimal $1$-models and Lemma~\ref{lem:functorial-malcev},
the resulting morphism on Malcev Lie algebras is compatible with the
quadratic presentations, yielding the commutative diagram in~\eqref{fun2}.

\smallskip
\eqref{fun2}$\Rightarrow$\eqref{fun3}. 
Identifying $\m(G)$ and $\m(H)$ with
$\widehat{\h}(H^{*}(G;\k))$ and $\widehat{\h}(H^{*}(H;\k))$,
respectively, transforms diagram~\eqref{fun2} into diagram~\eqref{fun3}.

\smallskip
\eqref{fun3}$\Rightarrow$\eqref{fun1}. 
Commutativity of diagram~\eqref{fun3} implies compatibility of 
$\m(\alpha)$ with the filtered-formality isomorphisms. 
Since $G$ and $H$ are $1$-formal, filtered-formality of $\alpha$
is equivalent to $1$-formality, and the claim follows.
\end{proof}

%%%%%%%%%%%%%%%%%%%%%%%%%%%%%
\subsection{Obstructions from infinitesimal Alexander invariants}
\label{subsec:alex-groups}
%%%%%%%%%%%%%%%%%%%%%%%%%%%%%

Throughout this subsection, $\alpha\colon G\to H$ is a homomorphism between
finitely generated groups, equipped with a model
$A(\alpha)\colon A_H \to A_G$ in the sense of
Definition~\ref{def:modeled-hom}; such models are available for all
homomorphisms in each of the classes of Example~\ref{ex:functorial-models}.
The morphism $A(\alpha)$ 
yields maps
\begin{equation}
\label{eq:AHAG-isos}
\begin{tikzcd}[column sep=38pt]
\widehat{\h}(A_G) \arrow[r, "\widehat{\h}(A(\alpha))"] &
\widehat{\h}(A_H)
\end{tikzcd}
\qquad\text{and}\qquad
\begin{tikzcd}[column sep=38pt]
\widehat{\B_1(A_G)} \arrow[r, "\widehat{\B_1(A(\alpha))}"] &
\widehat{\B_1(A_H)}.
\end{tikzcd}
\end{equation}

Both $\h$ and $\B_1$ are contravariant on $\cdgas$, so the two maps
in~\eqref{eq:AHAG-isos} run from the $G$-side to the $H$-side, in the same
direction as $\alpha$ itself.
By Theorem~\ref{thm:malcev-holo}, the first map compares with the induced
morphism $\m(\alpha)\colon \m(G;\k)\to \m(H;\k)$.
By Theorem~\ref{thm:alex-completed}, the second map compares with the
completed Alexander invariant morphism
$\widehat{B_1(\alpha)}\colon \widehat{B_1(G;\k)}\to\widehat{B_1(H;\k)}$.
Failure of the corresponding squares to commute provides an obstruction to
strict $1$-formality of~$\alpha$, in the sense of
Definition~\ref{def:rigid-formal-hom}.  Thus, whenever $\alpha$ admits a
model, the invariants $\h(A)$ and $\B_1(A)$ give effective obstructions.

We now make precise how infinitesimal Alexander invariants provide
secondary obstructions.

\begin{proposition}
\label{prop:alex-obstruction}
Let $G$ and $H$ be $1$-formal groups, and let $\alpha\colon G \to H$ be a
homomorphism admitting a model $A(\alpha)\colon A_H\to A_G$. 
Then $\alpha$ is strictly $1$-formal only if the following diagram commutes:
\begin{equation}
\label{eq:alex-obstruction}
\begin{tikzcd}[column sep=26pt, row sep=22pt]
\widehat{\B_1(A_G)} \ar[r, "\cong"] \ar[d, "\widehat{\B_1(A(\alpha))}"]
&
\widehat{B_1(G;\k)} \ar[d, "\widehat{B(\alpha)}"]
\\
\widehat{\B_1(A_H)} \ar[r, "\cong"]
&
\widehat{B_1(H;\k)} ,
\end{tikzcd}
\end{equation}
where $B_1(-)$ denotes the (first) Alexander invariant and $\B_1(-)$ its infinitesimal
counterpart. Consequently, failure of commutativity of~\eqref{eq:alex-obstruction}
obstructs the strict $1$-formality of $\alpha$.
\end{proposition}

\begin{proof}
By Lemma~\ref{lem:alex-functorial},
if $\alpha$ is strictly $1$-formal, the identifications of
Theorem~\ref{thm:alex-completed} on the two sides are induced by the chosen
formality data, and their naturality with respect to $\alpha$ yields the
commutativity of diagram~\eqref{eq:alex-obstruction}.
\end{proof}

We now specialize to a single $1$-formal group $G$, writing
$H_G=\big(H^*(G;\k),0\big)$ for its cohomology algebra with zero
differential; recall from Remark~\ref{rem:cohomology-not-a-family} that
$H_G$ is a $1$-finite $1$-model for $G$, with $\h(H_G)=\h(G;\k)$, but that
it is used one group at a time.  Three canonical maps are available:
\begin{equation}
\label{eq:Psi-G}
\begin{tikzcd}[column sep=26pt]
\gr^{I}\big(B_1(G;\k)\big) \ar[r, "\Psi_1"] &
\gr^{\m}\big(\B_1(H_G)\big) \ar[r, "\Psi_2"] &
\B_1(H_G) \ar[r, "\Psi_3"] &
\B\big(\h(G;\k)\big) .
\end{tikzcd}
\end{equation}
Here $\Psi_1$ is the isomorphism provided by
Theorem~\ref{thm:alex-completed}\eqref{ch2}, applied to
$A=H_G$; the map $\Psi_2$ is the canonical identification of a graded
$S$-module with the associated graded of its $\m$-adic filtration, available
since $d_{H_G}=0$, and shifting degrees by the initial degree of
$\B_1(H_G)$; and the map $\Psi_3$ is the isomorphism of
Theorem~\ref{thm:B-holo}.  There is no choice in $\Psi_3$: for
$H_G$ the splitting \eqref{eq:A1A1} is forced, since
$U_1=\im\bigl(d_{H_G}^{\vee}\bigr)=0$.  We set $\Psi_G=\Psi_3\Psi_2\Psi_1$,
an isomorphism of $S$-modules.

\begin{proposition}
\label{prop:koszul-obstruction-groups}
Let $G$ be a finitely generated, $1$-formal group. Then:
\begin{enumerate}
\item \label{ko1}
The isomorphism $\Psi_G$ of \eqref{eq:Psi-G} involves exactly one choice,
namely the $1$-formality isomorphism $\Theta_G$ entering $\Psi_1$, and is
independent of it.
\item \label{ko2}
If $\alpha\colon G\to H$ is a strictly $1$-formal homomorphism between 
finitely generated, $1$-formal groups, then, under the identifications 
$\Psi_G$ and $\Psi_H$ of part~\eqref{ko1},
\begin{equation}
\label{eq:gr-b1-formal}
\gr^{I}\!\big(B_1(\alpha)\big)=\B\big(\h(\alpha)\big) ,
\end{equation}
where $\h(\alpha)\colon \h(G;\k)\to\h(H;\k)$ is the morphism induced 
by $\alpha_*\colon H_1(G;\k)\to H_1(H;\k)$.
\end{enumerate}
\end{proposition}

\begin{proof}
\eqref{ko1}
Only the independence of $\Psi_G$ from $\Theta_G$ requires proof.  The
isomorphism of Theorem~\ref{thm:alex-completed}\eqref{ch1}, and with it
$\Psi_1$, is induced by a $1$-formality isomorphism 
$\Theta_G\colon \m(G;\k)\isom \widehat{\h(G;\k)}$; by 
\eqref{eq:1formal-malcev}, any two such choices differ by a filtered 
automorphism $u$ of $\widehat{\h(G;\k)}$ acting trivially on the 
degree-$1$ piece of the associated graded.  Since the graded Lie 
algebra $\h=\h(G;\k)$ is generated in degree~$1$, this forces 
$\gr(u)=\id$.  As in the proof of Theorem~\ref{thm:alex-completed}, 
the automorphism that $u$ induces on 
$\overline{\widehat{\h}{}'}/\overline{\widehat{\h}{}''}=
\widehat{\B(\h)}$ is filtration-preserving, and its associated graded 
is the automorphism of $\B(\h)$ induced by $\gr(u)=\id$, namely the 
identity.  Hence $\Psi_1$, and therefore $\Psi_G$, does not depend on the 
choice of $\Theta_G$.

\eqref{ko2}
Since $\alpha$ is strictly $1$-formal, the morphism 
$\alpha^*\colon \big(H^*(H;\k),0\big)\to \big(H^*(G;\k),0\big)$ models 
$\alpha$ compatibly with the chosen formality data, so that the 
identifications $\Psi_G$ and $\Psi_H$ may be used simultaneously.  
Lemma~\ref{lem:alex-functorial}, applied with $A(\alpha)=\alpha^*$, 
shows that, under the identifications of 
Theorem~\ref{thm:alex-completed}\eqref{ch1}, the completed morphism 
$\widehat{B_1(\alpha)}$ corresponds to $\widehat{\B_1(\alpha^*)}$.  
On the other hand, Proposition~\ref{prop:naturality-B-holo} 
identifies $\B_1(\alpha^*)$ with $\B(\h(\alpha^*))=\B(\h(\alpha))$.  
Passing to associated graded modules, and using the independence of 
choices established in part~\eqref{ko1}, we obtain 
\eqref{eq:gr-b1-formal}.
\end{proof}

Specializing to endomorphisms yields a computable obstruction, which will be 
put to use in Proposition~\ref{prop:alex-obstruction-example} below.  By 
Proposition~\ref{prop:strict-endo}\eqref{se2}, in this range strict 
$1$-formality is free of choices and amounts to $\m(\alpha)=\id$, so the 
conclusion is stated in that form.

\begin{corollary}
\label{cor:koszul-obstruction-endo}
Let $G$ be a finitely generated, $1$-formal group, and let 
$\alpha\colon G\to G$ be an endomorphism inducing the identity on 
$H_1(G;\k)$. Then:
\begin{enumerate}
\item \label{koe1}
$\h(\alpha)=\id$; consequently, if $\m(\alpha)=\id$, then 
$\gr^{I}(B_1(\alpha))=\id$.
\item \label{koe2}
Suppose $A_G$ is a $1$-finite $1$-model for $G$ and $A(\alpha)$ is an 
endomorphism of $A_G$ modeling $\alpha$.  If the endomorphism of 
$\gr^{\m}(\B_1(A_G))$ induced by $\B_1(A(\alpha))$ differs from the 
identity, then $\m(\alpha)\ne\id$; that is, $\alpha$ acts nontrivially on 
the Malcev completion of $G$, and is not strictly $1$-formal.
\end{enumerate}
\end{corollary}

\begin{proof}
\eqref{koe1}
The morphism $\h(\alpha)$ is induced by the identity 
$\alpha_*=\id$ of $\Lie(H_1(G;\k))$; being the identity, this map 
preserves the relation ideal of $\h(G;\k)$, and thus descends to the 
identity of $\h(G;\k)$.  If $\m(\alpha)=\id$ then, by 
Proposition~\ref{prop:strict-endo}\eqref{se2}, $\alpha$ is strictly 
$1$-formal, and 
Proposition~\ref{prop:koszul-obstruction-groups}\eqref{ko2} gives 
$\gr^{I}(B_1(\alpha))=\B(\id)=\id$.

\eqref{koe2}
By Lemma~\ref{lem:alex-functorial}, the completed morphism 
$\widehat{\B_1(A(\alpha))}$ corresponds to $\widehat{B_1(\alpha)}$; 
hence, by Theorem~\ref{thm:alex-completed}\eqref{ch2}, the induced 
endomorphism of $\gr^{\m}(\B_1(A_G))$ corresponds to 
$\gr^{I}(B_1(\alpha))$.  If the former differs from the identity, 
part~\eqref{koe1} shows that $\m(\alpha)\ne\id$.
\end{proof}

In practice it suffices to test the degree-$0$ piece: if $\B_1(A(\alpha))$ acts
nontrivially on $\gr^{\m}_0(\B_1(A_G))=\B_1(A_G)/\m\B_1(A_G)$, the hypothesis
of \eqref{koe2} holds.

\begin{remark}
\label{rem:koszul-obstruction-separated}
The obstruction of Corollary~\ref{cor:koszul-obstruction-endo} is 
sensitive only to the associated graded of $\B_1(A(\alpha))$: an 
endomorphism with $\B_1(A(\alpha))\ne \id$ but 
$\gr^{\m}(\B_1(A(\alpha)))=\id$ is not detected.  This is as it 
should be, since the comparison with the group-theoretic Alexander 
invariant proceeds through completion, and the $\m$-adic filtration 
on $\B_1(A_G)$ need not be separated (Section~\ref{subsec:koszul-filtrations}). 
\end{remark}

%%%%%%%%%%%%%%%%%%%%%%%%%%%%%
\subsection{Alexander invariants as secondary obstructions to formality}
\label{subsec:alex-secondary}
%%%%%%%%%%%%%%%%%%%%%%%%%%%%%

The obstruction of Corollary~\ref{cor:koszul-obstruction-endo} is genuinely
finer than the ones preceding it: it detects a nontrivial action on the
Malcev completion for an endomorphism that induces the identity on
cohomology, and hence on the holonomy Lie algebra and on all associated
graded Lie algebra invariants.
The conceptual reason, explained in
Remark~\ref{rem:naturality-B-holo-subtlety}, is that the model-level
invariants $\h(A_G)$ and $\B_1(A_G)\cong \B(\h(A_G))$ are strictly finer
than their cohomological counterparts.  The following proposition exhibits
such an automorphism.

\begin{proposition}
\label{prop:alex-obstruction-example}
Let $G=H(n)\times H(n)$ be the direct product of two
$2$-step Heisenberg groups, and let
$\alpha\colon G\to G$ be an automorphism inducing the identity on
$G_{\ab}$ but acting nontrivially on the center $Z(G)$.
Then:
\begin{enumerate}
\item \label{hh1}
$G$ is $1$-formal;
\item  \label{hh2}
$\alpha$ induces the identity on $\h(G)$ and $\gr(G)$;
\item  \label{hh3}
the induced morphism on the completed Alexander invariant
$\widehat{B_1(G;\k)}$ is not the identity;
\item  \label{hh4}
consequently, $\m(\alpha)\ne\id$: the automorphism $\alpha$ acts
nontrivially on the Malcev completion of $G$, and is not strictly
$1$-formal.
\end{enumerate}
\end{proposition}

\begin{proof}
\eqref{hh1} For $n>1$, each $H(n)$ is $1$-formal by
Proposition~\ref{prop:Heisenberg-1formal}, and $1$-formality is preserved
under finite direct products.

\smallskip
\eqref{hh2} Since $\alpha$ induces the identity on $G_{\ab}$,
it induces the identity on $H^1(G;\Z)$ and hence on the holonomy Lie algebra
$\h(G)$.
Because $G$ is $2$-step nilpotent, this also implies that $\alpha$
induces the identity on $\gr(G)$.

\smallskip
\eqref{hh3}
Let $A_G = \CE(\m(G))$, a $1$-finite $1$-model for $G$ as in
Example~\ref{ex:functorial-models}\eqref{fmx-nilp}, and let
$\varphi = A(\alpha) \in \Aut(A_G)$ be the induced 
$\cdga$ automorphism, which acts trivially on $H^*(A_G)$.
By Lemma~\ref{lem:B1-center}, applied to $\g = \m(G) = \h(n) \oplus \h(n)$,
\[
\B_1(A_G)/\m\B_1(A_G) \cong (Z(\m(G)))^\vee \cong \k^2.
\]
Since $\alpha$ acts nontrivially on $Z(G)$ by hypothesis, the induced
automorphism $\varphi$ acts nontrivially on $(Z(\m(G)))^\vee$, hence
nontrivially on $\B_1(A_G)$.
The identification of Theorem~\ref{thm:alex-completed} then implies that
$\widehat{B_1(G;\k)} \cong \widehat{\B_1(A_G)}$ also carries a nontrivial
induced automorphism.

\smallskip
\eqref{hh4} By Lemma~\ref{lem:B1-center}\eqref{b1-2}, the quotient 
$\B_1(A_G)/\m\B_1(A_G)$ is the degree-$0$ piece of 
$\gr^{\m}(\B_1(A_G))$; thus, the nontriviality established in 
\eqref{hh3}, despite the trivial action on $\h(G)$ and $\gr(G)$, 
gives $\m(\alpha)\ne\id$, by 
Corollary~\ref{cor:koszul-obstruction-endo}\eqref{koe2}.
\end{proof}

\begin{remark}
\label{rem:malcev-alexinv}
Proposition~\ref{prop:alex-obstruction-example} shows that the action of an
automorphism on the Malcev completion cannot be detected solely at the level
of holonomy or associated graded Lie algebras: the automorphism $\alpha$
induces the identity on both $\h(G)$ and $\gr(G)$, yet $\m(\alpha)\ne\id$.
This contrasts with the $\IA$-automorphisms of free groups
(Example~\ref{ex:ian}), where $\m(\alpha)\ne\id$ is read off directly from
the definition of $\alpha$; here it is invisible until one passes to
$\B_1$.
The Alexander invariant $B_1(G)$ and its Koszul counterpart $\B_1(A_G)$
provide obstructions to $1$-formality that are invisible at the level of
the associated graded Lie algebras. The $\B_q$-formality framework of
\S\ref{subsec:koszul-formality-obstruction} provides the 
conceptual setting for these higher-order phenomena.
\end{remark}

%%%%%%%%%%%%%%%%%%%%%%%%%%
\section{Ordered configuration spaces on an elliptic curve}
\label{sect:qp}
%%%%%%%%%%%%%%%%%%%%%%%%%%

The purpose of this section is to illustrate how the Koszul-module and
spectral-sequence techniques developed earlier can be applied to
configuration spaces.
For ordered configurations of points on an elliptic curve, the
fundamental group $P_{1,n}$ admits a finite-type $\cdga$ model due to
Bibby, whose Koszul invariants are simpler than those arising from the
cohomology algebra.
By comparing these invariants and analyzing the Koszul spectral sequence
associated to the weight filtration, we obtain a conceptual and
computable explanation for the failure of $1$-formality of $P_{1,n}$.

%%%%%%%%%%%%%%%%%%%%%%%%%%%%%
\subsection{Quasi-projective varieties and $\cdga$ models with weights}
\label{subsec:qproj-weights}
%%%%%%%%%%%%%%%%%%%%%%%%%%%%%

Let $X$ be a connected, smooth, complex quasi-projective variety.  By a
theorem of Deligne, each cohomology group of $X$ carries a mixed Hodge
structure~\cite{Del74}.  Smooth projective varieties, and more generally
compact K\"{a}hler manifolds, are formal (\S\ref{subsec:formal-maps}); a
general quasi-projective $X$ need not be, but its rational homotopy type is
still governed by a $\cdga$ model carrying a positive-weight structure of the
type studied in Section~\ref{subsec:weights-cdga}.  This is the model
underlying class~\eqref{fmx-qproj} of Example~\ref{ex:functorial-models};
we recall its construction and record the weight grading needed below.

Every $X$ as above admits a good compactification: a smooth projective
variety $\overline{X}$ and a normal-crossings divisor
$D=\bigcup_{j\in J} D_j$ such that $X=\overline{X}\setminus D$.  Morgan
\cite{Mo} associates to such a compactification a rationally defined $\cdga$
$A(X)=A(\overline{X},D)$, the \emph{Gysin model}, whose underlying graded
pieces are
\begin{equation}
\label{eq:gysin}
A^{p,q}= \bigoplus_{\abs{S}=q} H^p \Big(\bigcap_{k\in S} D_k,\,\C\Big)(-q),
\end{equation}
with multiplication induced by cup product and differential
$d\colon A^{p,q}\to A^{p+2,q-1}$ given by the Gysin maps of the
divisor intersections.  For a single variety $X$, Morgan's original
comparison of $A(X)$ with $\apl(X)$ by a chain of quasi-isomorphisms
preserving $\Q$-structures stands as stated in \cite{Mo}; it is only the
naturality of this comparison with respect to morphisms that requires the
care discussed in Remark~\ref{rem:naturality-caveat}.  Setting
\begin{equation}
\label{eq:morgan-weight}
\wt(A^{p,q})=p+2q
\end{equation}
defines a positive-weight decomposition on $(A(X),d)$ in the sense of
Definition~\ref{def:pos-weight}.  In particular, $W_1H^1(X,\C)=0$ (which
holds whenever $\overline{X}$ satisfies $b_1(\overline{X})=0$) implies the
$1$-formality of $X$, as shown by Morgan~\cite{Mo}; see also
Kohno~\cite{Ko83}.

%%%%%%%%%%%%%%%%%%%%%%%%%%%%%%
\subsection{The Bibby model for $\Conf(\E,n)$}
\label{subsec:conf-elliptic}
%%%%%%%%%%%%%%%%%%%%%%%%%%%%%%

Let $\E$ be a smooth elliptic curve over $\C$, and let
\begin{equation}
\label{eq:conf-en}
\Conf(\E,n)=\{(z_1,\dots,z_n)\in \E^{\times n}\mid z_i\neq z_j
\ \text{for}\ i\neq j\}
\end{equation}
be the ordered configuration space of $n$ points on $\E$, with
fundamental group 
\[
P_{1,n}=\pi_1(\Conf(\E,n)), 
\]
the pure braid group on $n$ strings on the torus.
Bibby~\cite{Bi16} constructed a finite-type $\cdga$ model $(A(n),d)$
for this space. Building on a spectral sequence of Totaro~\cite{To96},
her construction applies to complements of unimodular elliptic
arrangements in $\E^{\times n}$. As noted in class~\eqref{fmx-dupont} of
Example~\ref{ex:functorial-models}, the Bibby model is an instance of
Dupont's functorial Orlik--Solomon model~\cite{Du}, and the resulting
$\cdga$ carries a positive-weight structure compatible with
Definition~\ref{def:pos-weight}.

We now recall Bibby's especially transparent description 
of $(A(n),d)$.  This model is generated in degree~$1$ by elements
$a_i,b_i \ (1\le i\le n)$ and $e_{ij}\ (1\le i<j\le n)$, 
subject to the relations 
\[
(a_i-a_j)e_{ij}=(b_i-b_j)e_{ij}=e_{ij}e_{ik}-e_{ij}e_{jk}+e_{ik}e_{jk}=0,
\]
with differential given by
\[
d a_i = d b_i = 0, \qquad d e_{ij} = (a_i-a_j)(b_i-b_j).
\]

This model carries a positive-weight structure compatible with
the Morgan convention~\eqref{eq:morgan-weight}: the generators
$a_i, b_i$ have weight~$1$ and the generators $e_{ij}$ have weight~$2$,
reflecting the bidegrees $(1,0)$ and $(0,1)$ in~\eqref{eq:gysin} via
$\wt(A(n)^{p,q})=p+2q$, so that the differential preserves total weight
since $d e_{ij}=(a_i-a_j)(b_i-b_j)$ lands in weight $1+1=2$.
(Deligne's MHS convention assigns weight~$2$ to $H^1(E)$, so the two
conventions agree up to the re-indexing $\alpha\mapsto 3-\alpha$ on
$\{1,2\}$.)

%%%%%%%%%%%%%%%%%%%%%%%%%%%%%%
\subsection{The Koszul module of $A(n)$}
\label{subsec:conf-elliptic-koszul}
%%%%%%%%%%%%%%%%%%%%%%%%%%%%%%

For the $\cdga$ $(A(n),d)$ defined above, set
\[
H^1(A(n))=\langle [a_1],[b_1],\dots,[a_n],[b_n]\rangle \cong \k^{2n},
\qquad
S=\Sym(H_1(A(n)))=\k[x_1,y_1,\dots,x_n,y_n].
\]
Consider the Koszul complex 
$K_{\bullet}(A(n))=(A(n)_\bullet\otimes_{\k} S,\partial)$,
and denote $\B_1(A(n))=H_1(K_{\bullet}(A(n)))$. 

Since $d_{A(n)}\ne 0$, the differential $\partial^{A(n)}$ does not preserve
the $S$-degree, and $\B_1(A(n))$ carries no grading by $S$-degree.  It is
graded instead by the weight of Section~\ref{subsec:conf-elliptic}, extended
by $\wt(x_i)=\wt(y_i)=1$: both terms of $\partial^{A(n)}$ preserve
$\wt+\deg_S$.  We use this grading throughout, normalized so that the classes
$\varepsilon_{ij}$, of weight~$2$, sit in degree~$0$.

The next result establishes part~\eqref{tG3} of
Theorem~\ref{thm:elliptic-intro} from the Introduction.

\begin{theorem}
\label{thm:BA-n}
As a graded $S$-module, $\B_1(A(n))$ decomposes as 
\begin{equation}
\label{eq:B1-An-decomp}
\B_1(A(n)) \cong \boplus_{1\le i<j\le n} S/I_{ij},
\end{equation}
where $I_{ij}=(x_{\ell},y_{\ell} \mid \ell \ne i,j) + (x_i+x_j,\;y_i+y_j)$. 
In particular, each summand is isomorphic to $\k[x_i,y_i]$, and
$\B_1(A(n))$ is generated in degree~$0$ by the classes
$\varepsilon_{ij} = [e_{ij}^\vee \otimes 1]$, subject to the relations
\[
x_{\ell}\varepsilon_{ij}=y_{\ell}\varepsilon_{ij}=0 \ (\ell \ne i,j),
\qquad
(x_i+x_j)\varepsilon_{ij}=(y_i+y_j)\varepsilon_{ij}=0.
\]
\end{theorem}

\begin{proof}
Write $A=A(n)$ and $A^1=V\oplus W$, with
$V=\langle a_1,b_1,\dots,a_n,b_n\rangle$ and $W=\langle e_{ij}\rangle$, so
that $H^1(A)=V$ and
$\omega_A=\sum_{\ell}(a_\ell\otimes x_\ell+b_\ell\otimes y_\ell)$.
The defining relations are homogeneous in the number of $e$-factors, so
$A^2=\bwedge^2V\oplus (V\cdot W)\oplus \bwedge^2W$, and dually
$A_2=(\bwedge^2V)^{\vee}\oplus (V\cdot W)^{\vee}\oplus(\bwedge^2W)^{\vee}$.

We use both terms of $\partial^A=\omega_A^{\vee}\lrcorner(-)+d_A^{\vee}\otimes\id$
(Definition~\ref{def:koszul-chain}).  Since $d_Aa_i=d_Ab_i=0$ and
$d_Ae_{ij}=(a_i-a_j)(b_i-b_j)$, we have $d_A(A^1)\subseteq\bwedge^2V$, whence
\begin{equation}
\label{eq:dvee-An}
d_A^{\vee}(\xi)=\sum_{k<\ell}\langle \xi,d_Ae_{k\ell}\rangle\,e_{k\ell}^{\vee}
\quad\text{for } \xi\in A_2 .
\end{equation}
In particular $\im d_A^{\vee}\subseteq W^{\vee}$, and $d_A^{\vee}$ vanishes on
$(V\cdot W)^{\vee}$ and on $(\bwedge^2W)^{\vee}$.

\smallskip
\emph{Step 1: the cycles.}
The map $\partial^A_1$ is dual to $\delta_A^0(1\otimes s)=\omega_A s$, so it
has no $d_A^{\vee}$ component:
$\partial^A_1(a_\ell^{\vee}\otimes s)=x_\ell s$,
$\partial^A_1(b_\ell^{\vee}\otimes s)=y_\ell s$, and
$\partial^A_1(e_{ij}^{\vee}\otimes s)=0$.  Hence
\begin{equation}
\label{eq:Z1-split-An}
\ker\partial^A_1=Z\oplus(W^{\vee}\otimes_{\k}S),
\qquad Z=\ker\bigl(\partial^A_1|_{V^{\vee}\otimes S}\bigr),
\end{equation}
and $\partial^A_1|_{V^{\vee}\otimes S}$ is the Koszul differential of the
regular sequence $x_1,y_1,\dots,x_n,y_n$.

\smallskip
\emph{Step 2: reduction to a submodule of $W^{\vee}\otimes S$.}
Let $\eta\in(\bwedge^2V)^{\vee}\otimes S$.  Then
$\partial^A_2(\eta)=\omega_A^{\vee}\lrcorner\,\eta+d_A^{\vee}(\eta)$, with the
first summand in $V^{\vee}\otimes S$ and the second in $W^{\vee}\otimes S$.
By exactness of the Koszul complex on a regular sequence, every $c\in Z$ is of
the form $\omega_A^{\vee}\lrcorner\,\eta$; for such $\eta$,
$c\equiv-d_A^{\vee}(\eta)\pmod{\im\partial^A_2}$.  With
\eqref{eq:Z1-split-An} this shows that
$W^{\vee}\otimes_{\k}S\to\B_1(A)$ is surjective, so $\B_1(A)$ is generated by
the classes $\varepsilon_{ij}$, and
\begin{equation}
\label{eq:B1-as-quotient}
\B_1(A)\cong (W^{\vee}\otimes_{\k}S)/N,
\qquad N=(W^{\vee}\otimes_{\k}S)\cap\im\partial^A_2 .
\end{equation}
Moreover $\partial^A_2$ vanishes on $(\bwedge^2W)^{\vee}\otimes S$ and maps
$(V\cdot W)^{\vee}\otimes S$ into $W^{\vee}\otimes S$, while its restriction
to $(\bwedge^2V)^{\vee}\otimes S$ lands in $W^{\vee}\otimes S$ exactly on
$\ker(\omega_A^{\vee}\lrcorner)$.  Hence
\begin{equation}
\label{eq:N-two-families}
N=\partial^A_2\bigl((V\cdot W)^{\vee}\otimes S\bigr)
 + d_A^{\vee}\bigl(\ker(\omega_A^{\vee}\lrcorner)\cap (\bwedge^2V)^{\vee}\otimes S\bigr).
\end{equation}

\smallskip
\emph{Step 3: the first family gives $I_{ij}$.}
For $\ell\ne i,j$ the element $(a_\ell e_{ij})^{\vee}$ lies in
$(V\cdot W)^{\vee}$, so $d_A^{\vee}$ kills it, and contraction gives
$\partial^A_2\bigl((a_\ell e_{ij})^{\vee}\otimes1\bigr)=e_{ij}^{\vee}\otimes x_\ell$.
The relation $(a_i-a_j)e_{ij}=0$ gives $a_ie_{ij}=a_je_{ij}$ in $A^2$, so
$(a_ie_{ij})^{\vee}$ is a single element of $A_2$, and
$\partial^A_2\bigl((a_ie_{ij})^{\vee}\otimes1\bigr)
=e_{ij}^{\vee}\otimes(x_i+x_j)$.
The same computations with $b$ in place of $a$ give $y_\ell$ and $y_i+y_j$.
Hence the first summand of \eqref{eq:N-two-families} is
$\boplus_{i<j}I_{ij}\,e_{ij}^{\vee}$.

\smallskip
\emph{Step 4: the second family is contained in the first.}
Let $\eta\in\ker(\omega_A^{\vee}\lrcorner)\cap(\bwedge^2V)^{\vee}\otimes S$.
By exactness again, $\eta=\omega_A^{\vee}\lrcorner\,\theta$ for some
$\theta\in(\bwedge^3V)^{\vee}\otimes S$.  By \eqref{eq:dvee-An}, the
$(k,\ell)$-component of $d_A^{\vee}(\eta)$ is
$\langle\eta,\,d_Ae_{k\ell}\rangle$, and by \eqref{eq:contraction-def}
contraction by $\omega_A^{\vee}$ is adjoint to multiplication by $\omega_A$:
\[
\bigl\langle \eta,\,d_Ae_{k\ell}\bigr\rangle
=\bigl\langle \theta,\;\omega_A\wedge d_Ae_{k\ell}\bigr\rangle .
\]
Fix $k<\ell$ and put $p=a_k-a_\ell$, $q=b_k-b_\ell$, so $d_Ae_{k\ell}=pq$.
Substituting $a_k=p+a_\ell$ and $b_k=q+b_\ell$ in $\omega_A$ gives
\[
\omega_A=p\,x_k+q\,y_k+a_\ell(x_k+x_\ell)+b_\ell(y_k+y_\ell)
        +\sum_{m\ne k,\ell}(a_m x_m+b_m y_m),
\]
and therefore, since $p\wedge pq=q\wedge pq=0$,
\begin{equation}
\label{eq:omega-wedge-de}
\omega_A\wedge d_Ae_{k\ell}
= a_\ell\,pq\,(x_k+x_\ell)+b_\ell\,pq\,(y_k+y_\ell)
 +\sum_{m\ne k,\ell}\bigl(a_m\,pq\,x_m+b_m\,pq\,y_m\bigr).
\end{equation}
Every coefficient in \eqref{eq:omega-wedge-de} lies in $I_{k\ell}$, hence so
does $\langle\eta,\,d_Ae_{k\ell}\rangle$, and the $(k,\ell)$-component of
$d_A^{\vee}(\eta)$ already lies in $I_{k\ell}e_{k\ell}^{\vee}$.

\smallskip
Combining Steps 3 and 4, $N=\boplus_{i<j}I_{ij}\,e_{ij}^{\vee}$, and
\eqref{eq:B1-as-quotient} gives \eqref{eq:B1-An-decomp}.
\end{proof}

\begin{corollary}
\label{cor:BA-n-consequences}
For each $n\ge 2$, the following hold. 
\begin{enumerate}
\item \label{bn1}
The Koszul module $\B_1(A(n))$ has Hilbert series
\begin{equation}
\label{eq:hilb-B1-An}
\Hilb(\B_1(A(n)),t)=\frac{\binom{n}{2}}{(1-t)^2}.
\end{equation}
\item \label{bn2}
The Chen ranks of $P_{1,n}$ satisfy
\begin{equation}
\label{eq:chen-p1n-ranks}
\theta_1(P_{1,n})=2n, \qquad \theta_k(P_{1,n})=\binom{n}{2} (k-1)
\quad\text{for $k\ge 2$}.
\end{equation}
\item \label{bn3}
The first resonance variety of $A(n)$ is given by
\begin{equation}
\label{eq:res-An}
\RR^1(A(n))=\RR_1(A(n))=
\bigcup_{1\le i<j\le n}
\bigl\{x_{\ell}=y_{\ell}=0\; (\ell\notin \{i,j\}),\; x_i+x_j=y_i+y_j=0\bigr\}
\subset \k^{2n}.
\end{equation}
\end{enumerate}
\end{corollary}

\begin{proof}
Claims \eqref{bn1} and \eqref{bn3} follow at once from the decomposition of 
$\B_1(A(n))$ given in Theorem~\ref{thm:BA-n}. 

Claim \eqref{bn2} follows from \eqref{bn1}, together with 
Theorem~\ref{thm:alex-completed}\eqref{ch3}. 
\end{proof}

%%%%%%%%%%%%%%%%%%%%%%%%%%%%%%
\subsection{The cohomology algebra and holonomy Lie algebra of $P_{1,n}$}
\label{subsec:koszul-mod-p1n}
%%%%%%%%%%%%%%%%%%%%%%%%%%%%%%

In~\cite{Bez}, Bezrukavnikov computed presentations of the Malcev Lie algebras
$\m(P_{g,n})$ of the pure braid groups $P_{g,n}$ for all $g,n\ge 1$.
He observed that these groups are $1$-formal for $g>1$ and $n\ge 1$,
and for $g=1$ and $n\le 2$, and stated without proof that
$P_{1,n}$ is not $1$-formal for $n\ge 3$.
A complete proof of non-$1$-formality of $P_{1,n}$ for $n\ge 3$
was later given in~\cite{DPS-duke}, using the Tangent Cone Theorem
for cohomology jump loci of $1$-formal groups and the nonlinearity of
$\RR^1(\Conf(\E,n))$.

In the next four subsections we give an independent obstruction to
$1$-formality, based on Chen ranks and Koszul modules.  It illustrates a
general phenomenon: for a space admitting a finite-type model, any
discrepancy between the Koszul invariants of the model and those of the
cohomology algebra is a computable obstruction to formality.  We begin with
the presentation of the truncated cohomology algebra and of the associated
holonomy Lie algebra.

We first record a splitting that will be used repeatedly.  Since
$\E=S^1\times S^1$ is a topological group, written additively, the last
point of a configuration may be translated to the identity; writing
$\E^{*}=\E\setminus\{0\}$ for the elliptic curve with the identity removed,
this gives a homeomorphism
\begin{equation}
\label{eq:conf-elliptic-topsplit}
\Conf(\E,n) \longisom \Conf(\E^{*},n-1)\times \E,
\qquad
(z_1,\dots,z_n)\longmapsto
\bigl((z_1-z_n,\dots,z_{n-1}-z_n),\,z_n\bigr),
\end{equation}
with inverse $\bigl((w_1,\dots,w_{n-1}),z\bigr)\mapsto
(w_1+z,\dots,w_{n-1}+z,z)$.  Both maps are continuous, since the group
operations of $\E$ are, and they are mutually inverse; that each lands where
claimed is immediate from the defining inequalities, since $z_i\ne z_n$ gives
$z_i-z_n\ne 0$ and $z_i\ne z_j$ gives $z_i-z_n\ne z_j-z_n$.  Being a
homeomorphism, \eqref{eq:conf-elliptic-topsplit} induces an isomorphism of
graded algebras on cohomology with coefficients in the field $\k$, by the
K\"{u}nneth formula.

\begin{proposition}
\label{prop:conf-elliptic-H2}
Let $\k$ be a field of characteristic $0$, and set
\begin{equation}
\label{eq:Qn-quad}
Q(n)\coloneqq\frac{\bwedge(a_1,b_1,\dots,a_n,b_n)}{
\bigl\langle (a_i-a_j)(b_i-b_j) \mid 1\le i<j\le n \bigr\rangle},
\qquad \deg a_i=\deg b_i=1 .
\end{equation}
Then $Q(n)$ is the quadratic closure of $H^{*}(\Conf(\E,n);\k)$, and the two
agree in degrees $\le 2$:
\[
H^{\le 2}(\Conf(\E,n);\k) \cong Q(n)^{\le 2} .
\]
Moreover:
\begin{enumerate}[itemsep=1pt]
\item $H^1(\Conf(\E,n);\k) = \k^{2n}$ with basis $\{a_i,b_i\}_{i=1}^n$.
\item The cup product $\mu\colon \bigwedge^2 H^1(\Conf(\E,n);\k)\to 
H^2(\Conf(\E,n);\k)$ is equivalent to the composite
\[
\begin{tikzcd}[column sep=22pt]
\bwedge^2 H^1(E^{\times n};\k) \arrow[r, "\cup"] & H^2(\E^{\times n};\k) 
\arrow[r, two heads] & H^2(\E^{\times n};\k)/\langle [\Delta_{ij}]\rangle_{i<j}.
\end{tikzcd}
\]
\item The holonomy Lie algebra $\h(P_{1,n};\k)$ is the quadratic Lie
algebra on generators $x_i,y_i$ dual to $a_i,b_i$, with relation space the
annihilator of $\bigl\langle [\Delta_{ij}]\bigr\rangle_{i<j}$ in
$\bwedge^2 H_1$.  Thus $\dim_\k\gr_2\h(P_{1,n};\k)=\binom{n}{2}$, while the
number of defining quadratic relations is $\binom{2n}{2}-\binom{n}{2}$.
\item \label{qn4}
The relation ideal splits off a free exterior factor on two generators:
writing $p_i=a_i-a_n$, $q_i=b_i-b_n$ for $1\le i\le n-1$,
\begin{equation}
\label{eq:conf-elliptic-split}
Q(n) \cong
\frac{\bwedge(p_1,q_1,\dots,p_{n-1},q_{n-1})}{\langle p_iq_i,\
p_iq_j+p_jq_i\rangle}\ \otimes_{\k}\ \bwedge(a_n,b_n) ,
\end{equation}
the algebraic counterpart of the splitting
\eqref{eq:conf-elliptic-topsplit}.  At $n=2$ this reads
$Q(2)\cong H^{*}(\E^{*};\k)\otimes_{\k} H^{*}(\E;\k)$.
\end{enumerate}
\end{proposition}

\begin{proof}
It follows from Totaro~\cite{To96} that the cohomology ring of $\Conf(\E,n)$ 
in degrees $\le 2$ is obtained from $H^*(\E^{\times n};\k)$ by quotienting 
out the classes dual to the diagonals $\Delta_{ij}=\{(z_1,\dots,z_n)\mid z_i=z_j\}$.
Since $H^*(\E^{\times n};\k)\cong \bwedge(a_1,b_1,\dots,a_n,b_n)$ 
with $\deg a_i=\deg b_i=1$,
this yields the stated presentation in degree $\le 2$.  The truncation is
not a formality: $Q(n)$ is nonzero in degrees up to $2n$, whereas
$H^{\le2}(\Conf(\E,n);\k)$ vanishes above degree~$2$ by definition.

The description of the cup product follows from the naturality of the
quotient map
\[
H^2(\E^{\times n};\k) \longsurj H^2(\E^{\times n};\k)/\langle [\Delta_{ij}]\rangle_{i<j} ,
\]
where $[\Delta_{ij}]$ denotes the Poincaré dual cohomology class of the
diagonal $\Delta_{ij}\subset \E^{\times n}$. 
Dualizing gives the presentation of the holonomy Lie algebra: the relations
of $\h(P_{1,n};\k)$ are the annihilator of $\ker(\mu)$, which has dimension
$\binom{2n}{2}-\binom{n}{2}$, and $\gr_2\h(P_{1,n};\k)\cong\ker(\mu)^{\vee}$
has dimension $\binom{n}{2}$.

For \eqref{qn4}, set $p_i=a_i-a_n$, $q_i=b_i-b_n$ for $i<n$.  Then
$(a_i-a_n)(b_i-b_n)=p_iq_i$, while for $i,j<n$,
\[
(a_i-a_j)(b_i-b_j)=(p_i-p_j)(q_i-q_j)\equiv -\bigl(p_iq_j+p_jq_i\bigr)
\pmod{p_iq_i,\ p_jq_j},
\]
so the $\binom{n}{2}$ defining relations of $Q(n)$ are generated by
$\{p_iq_i\}_{i<n}$ together with $\{p_iq_j+p_jq_i\}_{i<j<n}$, and $a_n,b_n$
occur in none of them.  This gives \eqref{eq:conf-elliptic-split}.
\end{proof}

The splitting of item~\eqref{qn4} has a practical consequence for
$\B_1(A(n))$, already used implicitly in Corollary~\ref{cor:BA-n-consequences}
above: since the two extra generators $a_n,b_n$ are free, they contribute
nothing to $\ker(\cup)$, so the Hilbert series, minimal generators, support,
and Krull dimension of $\B_1(A(n))$ may all be read off from the
basepointed presentation on $n-1$ pairs.  Only the projective dimension
over $S$ differs, rising by two, since the two additional variables of $S$
annihilate $\B_1$.

%%%%%%%%%%%%%%%%%%%%%%%%%%%%%%
\subsection{The first resonance variety of $\Conf(\E,n)$}
\label{subsec:res-conf-elliptic}
%%%%%%%%%%%%%%%%%%%%%%%%%%%%%%

The next proposition provides a complete proof of the description
of the first resonance variety of $\Conf(\E,n)$ stated without proof in
\cite{DPS-duke}; the higher resonance varieties are not treated here.  
Its nonlinearity for $n\ge3$ is what \cite{DPS-duke} exploits, via the 
Tangent Cone Theorem, to obstruct $1$-formality; the route taken in
Section~\ref{subsec:p1n-nonformal} is independent of it.

\begin{proposition}
\label{prop:res-elliptic}
Let $\E$ be a smooth elliptic curve over $\C$. For each $n\ge 2$, the 
first resonance variety of $A = H^{\le 2}(\Conf(\E,n);\k)$ is given by 
\begin{equation}
\label{eq:res-conf-elliptic}
\RR^1(A)=\left\{ (x,y)\in \k^n\times \k^n \ \middle|\ 
\begin{array}{l}
\sum_{i=1}^n x_i=\sum_{i=1}^n y_i=0,\\[2pt]
x_i y_j-x_j y_i=0 \text{ for all } 1\le i<j\le n
\end{array}
\right\}.
\end{equation}
Equivalently, $\RR^1(A)$ is the affine cone over the Segre embedding of
$\mathbb{P}^1\times\mathbb{P}^{n-2}$.
\end{proposition}

\begin{proof}
Since $d_A=0$, Proposition~\ref{prop:res-descriptions}\eqref{rd3} identifies
$\RR^1(A)$ with the set of those $u\in A^1$ for which the multiplication map
$\mu_u\colon A^1\to A^2$, $v\mapsto u\cdot v$, has kernel strictly larger than
$\k u$.

Write $H^1(\Conf(\E,n);\k)=W\otimes_{\k} U$, where $W\cong\k^2$ has basis
$\alpha,\beta$ and $U=\k^{n}$ has basis $u_1,\dots,u_n$, so that
$a_i=\alpha\otimes u_i$ and $b_i=\beta\otimes u_i$.  Let $\sigma\colon U\to\k$
be the linear form with $\sigma(u_i)=1$ for all $i$, and let $U_0=\ker\sigma$
be the sum-zero subspace.  Under the Cauchy decomposition
\begin{equation}
\label{eq:cauchy-res}
\bwedge^2(W\otimes_{\k} U)=
\bigl(\Sym^2 W\otimes_{\k}\bwedge^2 U\bigr)\oplus
\bigl(\bwedge^2 W\otimes_{\k}\Sym^2 U\bigr),
\end{equation}
an element $(w\otimes c)\wedge(w'\otimes c')$ corresponds to
$\tfrac{1}{2}\bigl((ww')\otimes(c\wedge c')+(w\wedge w')\otimes(cc')\bigr)$.

By Proposition~\ref{prop:conf-elliptic-H2}, $A^2=\bwedge^2V/R$, where $V=A^1$
and $R$ is the subspace spanned by the diagonal classes
$\Delta_{ij}=(a_i-a_j)(b_i-b_j)$ with $i<j$; these are the only relations.
Since
\[
\Delta_{ij}=\bigl(\alpha\otimes(u_i-u_j)\bigr)\wedge
\bigl(\beta\otimes(u_i-u_j)\bigr)=(\alpha\wedge\beta)\otimes(u_i-u_j)^2 ,
\]
each $\Delta_{ij}$ lies in the second summand of \eqref{eq:cauchy-res}.
Moreover, the $\binom{n}{2}$ elements $(u_i-u_j)^2$ form a basis of
$\Sym^2U_0$: they lie in $\Sym^2U_0$, they are linearly independent because
the coefficient of $u_iu_j$ in $(u_k-u_l)^2$ is $-2$ when $\{i,j\}=\{k,l\}$
and $0$ otherwise, and $\dim_\k\Sym^2U_0=\binom{n}{2}$.  Consequently,
\begin{equation}
\label{eq:ker-cup-res}
R=\ker\bigl(\cup\colon \bwedge^2H^1\to H^2\bigr)
=\bwedge^2 W\otimes_{\k}\Sym^2 U_0 .
\end{equation}
We will use the description of $\Sym^2U_0$ as the kernel of the contraction
\begin{equation}
\label{eq:sigma-two}
\sigma_2\colon \Sym^2 U\to U,\qquad cc'\mapsto \sigma(c)c'+\sigma(c')c .
\end{equation}
Indeed, $\Sym^2U_0\subseteq\ker\sigma_2$, and the two subspaces coincide
because $\sigma_2$ is surjective---if $\sigma(c)=1$, then $\sigma_2(c^2)=2c$
and $\sigma_2(cc')=c'$ for every $c'\in U_0$---so that
$\dim_\k\ker\sigma_2=\binom{n+1}{2}-n=\binom{n}{2}$.

Now let $u=\alpha\otimes X+\beta\otimes Y$ and
$v=\alpha\otimes P+\beta\otimes Q$ be elements of $V$, where
$X=\sum_{i} x_iu_i$, $Y=\sum_{i} y_iu_i$, and $P,Q\in U$.  Expanding by
\eqref{eq:cauchy-res},
\begin{equation}
\label{eq:uv-expand}
2\,u\wedge v=
\alpha^2\otimes(X\wedge P)
+\beta^2\otimes(Y\wedge Q)
+\alpha\beta\otimes(X\wedge Q+Y\wedge P)
+(\alpha\wedge\beta)\otimes(XQ-YP).
\end{equation}
Comparing with \eqref{eq:ker-cup-res}, we conclude that $u\cdot v=0$ in $A^2$
if and only if
\begin{equation}
\label{eq:res-conditions}
X\wedge P=0,\qquad
Y\wedge Q=0,\qquad
X\wedge Q+Y\wedge P=0,\qquad
\sigma_2(XQ-YP)=0,
\end{equation}
the last condition reading, by \eqref{eq:sigma-two},
\begin{equation}
\label{eq:res-fourth}
\sigma(X)Q+\sigma(Q)X-\sigma(Y)P-\sigma(P)Y=0 .
\end{equation}

Suppose first that $X\wedge Y\ne 0$.  Then $X$ and $Y$ are both nonzero, so
the first two conditions in \eqref{eq:res-conditions} force $P=\lambda X$ and
$Q=\mu Y$ for some $\lambda,\mu\in\k$; the third then gives
$(\mu-\lambda)\,X\wedge Y=0$, so $\mu=\lambda$ and $v=\lambda u$.  Hence
$\ker\mu_u=\k u$, and $u\notin\RR^1(A)$.

Suppose next that $X\wedge Y=0$.  Excluding the trivial case $u=0$, we may
write $X=sc$ and $Y=tc$ with $c\in U\setminus\{0\}$ and
$(s,t)\in\k^2\setminus\{0\}$.  Every $v\in\ker\mu_u$ then has $P,Q\in\k c$:
if $s$ and $t$ are both nonzero, this follows from the first two conditions in
\eqref{eq:res-conditions}; if $t=0$, so that $s\ne0$ and $Y=0$, the first and
third conditions give it; and the case $s=0$ is symmetric.  Conversely, for
$P=\lambda c$ and $Q=\mu c$ the first three conditions in
\eqref{eq:res-conditions} hold automatically, while \eqref{eq:res-fourth}
becomes $2\sigma(c)(s\mu-t\lambda)\,c=0$, that is,
\begin{equation}
\label{eq:res-sigma-cond}
\sigma(c)\,(s\mu-t\lambda)=0 .
\end{equation}
If $\sigma(c)=0$, then \eqref{eq:res-sigma-cond} holds for all
$(\lambda,\mu)\in\k^2$, and choosing $(\lambda,\mu)$ not proportional to
$(s,t)$ yields an element of $\ker\mu_u$ outside $\k u$; hence
$u\in\RR^1(A)$.  If $\sigma(c)\ne0$, then \eqref{eq:res-sigma-cond} forces
$s\mu=t\lambda$, that is, $(\lambda,\mu)$ proportional to $(s,t)$, so that
$v\in\k u$; hence $u\notin\RR^1(A)$.

Summarizing, $u\in\RR^1(A)$ if and only if $X\wedge Y=0$ and, in the above
notation, $\sigma(c)=0$.  The first condition is the vanishing of the
$2\times2$ minors $x_iy_j-x_jy_i$, and, since $(s,t)\ne0$, the second is
equivalent to $\sigma(X)=\sigma(Y)=0$.  As $u=0$ satisfies both conditions,
this is precisely \eqref{eq:res-conf-elliptic}.

For the last assertion, the above description exhibits $\RR^1(A)$ as the set
of pure tensors $w\otimes c$ with $w\in W$ and $c\in U_0$, that is, the affine
cone over the Segre embedding of $\mathbb{P}(W)\times\mathbb{P}(U_0)\cong
\mathbb{P}^1\times\mathbb{P}^{n-2}$.  In particular, $\RR^1(A)$ is irreducible
of dimension $n$; it is a linear subspace when $n=2$, and a nonlinear
determinantal variety when $n\ge3$.
\end{proof}

\begin{remark}
\label{rem:res-elliptic-history}
The description in Proposition~\ref{prop:res-elliptic} was asserted in
\cite[Ex.~10.1]{DPS-duke}, where it was verified by computer for small values
of $n$ and used to deduce non-$1$-formality of $P_{1,n}$.  No proof of it has
appeared in the intervening two decades; Proposition~\ref{prop:res-elliptic}
is the first.
\end{remark}

%%%%%%%%%%%%%%%%%%%%%%%%%%%%%%
\subsection{Non-$1$-formality of the pure elliptic braid groups}
\label{subsec:p1n-nonformal}
%%%%%%%%%%%%%%%%%%%%%%%%%%%%%%

We can now compare the Bibby model $A(n)$ with the truncated cohomology
algebra $H^{\le 2}(\Conf(\E,n);\k)$.  Although the two agree in degree~$1$
and arise from the same fundamental group, their Koszul modules differ
substantially, and the discrepancy in the asymptotic growth of the
corresponding Chen ranks is what obstructs $1$-formality of $P_{1,n}$.

\begin{proposition}
\label{prop:conf-elliptic-comparison}
For $n\ge 3$, the following hold.
\begin{enumerate}[itemsep=2pt]
\item \label{pn1}
$\RR^1(\Conf(\E,n);\k)$ is not isomorphic to $\RR^1(A(n))$ as
algebraic varieties.
\item \label{pn2}
The Chen ranks of\/ $P_{1,n}$ satisfy $\theta_k(P_{1,n})=\binom{n}{2}(k-1)$
for $k\ge 2$, and thus grow linearly in $k$.
\item \label{pn3}
The holonomy Chen ranks $\theta_k(\h(P_{1,n};\k))$ grow polynomially of
degree exactly $n-1\ge 2$ in $k$, and in particular strictly faster than
the Chen ranks of $P_{1,n}$.
\end{enumerate}
\end{proposition}

\begin{proof}
\eqref{pn1}
By Proposition~\ref{prop:res-elliptic}, $\RR^1(\Conf(\E,n);\k)$ is an
irreducible non-linear variety (defined by the $2\times 2$ minors of
a $2\times n$ matrix of linear forms), whereas
Corollary~\ref{cor:BA-n-consequences} shows that $\RR^1(A(n))$ is a
union of $\binom{n}{2}$ dimension-$2$ linear subspaces. Hence the
two varieties are not isomorphic.

\eqref{pn2}
This is the content of Corollary~\ref{cor:BA-n-consequences}\eqref{bn2},
which follows from the decomposition
$\B_1(A(n))\cong\bigoplus_{i<j}S/I_{ij}$ of
Theorem~\ref{thm:BA-n}, each summand being a polynomial ring in two
variables.

\eqref{pn3}
By Proposition~\ref{prop:holo-massey} and Theorem~\ref{thm:B-holo},
\[
\theta_{k+2}(\h(P_{1,n};\k)) = \dim_\k[\B_1(H^*(\Conf(\E,n);\k))]_k.
\]
So it suffices to bound the growth of the Hilbert function of
$\B_1(H^{*}(\Conf(\E,n);\k))$.  By Proposition~\ref{prop:r1a}, the support
of that module agrees with $\RR^1(H^{*}(\Conf(\E,n);\k))$ away from the
origin, and by Proposition~\ref{prop:res-elliptic} the latter is the
variety of rank-one $2\times n$ matrices whose rows sum to
zero---the affine cone over the Segre variety
$\mathbb{P}^1\times\mathbb{P}^{n-2}$, of dimension $n$.  Hence
$\B_1(H^{*}(\Conf(\E,n);\k))$ is a finitely generated graded $S$-module of
Krull dimension $n$, so its Hilbert polynomial has degree exactly $n-1$.
For $n\ge3$ this is at least $2$, which exceeds the linear growth
in~\eqref{pn2}.
\end{proof}

\begin{corollary}
\label{cor:p1n-nonformal}
For every $n\ge 3$, the pure elliptic braid group $P_{1,n}$ is not
$1$-formal.
\end{corollary}

\begin{proof}
If $P_{1,n}$ were $1$-formal, Theorem~\ref{thm:chen-ineq-model} would
give $\theta_k(P_{1,n})=\theta_k(\h(P_{1,n};\k))$ for all $k\ge 1$.
But part~\eqref{pn2} gives linear growth $\theta_k(P_{1,n})\sim\binom{n}{2}k$,
while part~\eqref{pn3} gives polynomial growth of degree $n-1\ge 2$.
This contradiction proves non-$1$-formality.
\end{proof}

Proposition~\ref{prop:res-elliptic}, together with
Corollary~\ref{cor:BA-n-consequences}\eqref{bn3} and
Proposition~\ref{prop:conf-elliptic-comparison}\eqref{pn1}, establishes
part~\eqref{tG2} of Theorem~\ref{thm:elliptic-intro} from the Introduction;
parts~\eqref{pn2} and~\eqref{pn3} of
Proposition~\ref{prop:conf-elliptic-comparison}, together with
Corollary~\ref{cor:p1n-nonformal}, establish part~\eqref{tG1}.

%%%%%%%%%%%%%%%%%%%%%%%%%%%%%%
\subsection{The structure of $\B_1(H^{*}(\Conf(\E,n);\k))$}
\label{subsec:b1-conf-structure}
%%%%%%%%%%%%%%%%%%%%%%%%%%%%%%

The Koszul module appearing in
Proposition~\ref{prop:conf-elliptic-comparison} is of interest in its own
right.  Its Hilbert series, its support, and---conjecturally---its structure
as a representation are all governed by a single geometric picture, which we
now describe.

\begin{problem}
\label{conj:p1n-module}
Determine the structure of the graded $S$-module
$\B_1(H^{*}(\Conf(\E,n);\k))$. 
\end{problem}

Macaulay2 computations reveal a uniform pattern.  For $2\le n\le 10$, the
module has projective dimension $2(n-1)$ and Krull dimension~$n$, while its
support is cut out inside $\mSpec(S)$ as
\begin{equation}
\label{eq:supp-b1-conf}
\Supp_S\B_1\bigl(H^{*}(\Conf(\E,n);\k)\bigr)
=\bV\Bigl(I_2(M)+\bigl(\textstyle\sum_{i} x_i,\ \sum_{i} y_i\bigr)\Bigr),
\qquad
M=\begin{pmatrix} x_1&\cdots&x_n\\ y_1&\cdots&y_n\end{pmatrix} ;
\end{equation}
that is, it is the affine cone over the Segre variety
$\mathbb{P}^1\times\mathbb{P}^{n-2}$, consisting of the rank-one
matrices whose rows sum to zero, of dimension~$n$.  Finally, the Hilbert
series is given by
\begin{equation}
\label{eq:hilb-b1-conf}
\Hilb(\B_1(H^{*}(\Conf(\E,n);\k)),t)
= \frac{\displaystyle\sum_{k=0}^{n-2}(-1)^k\binom{n}{k+2}\,t^k}{(1-t)^{n}} .
\end{equation}
Equivalently, $P_n(t)=t^{-2}\bigl[(1-t)^{n}-1+nt\bigr]$ and $P_n(1)=n-1$, so
the Hilbert polynomial has degree $n-1$ with leading coefficient
$(n-1)/(n-1)!=1/(n-2)!$.  For $n=3$ the resulting Chen ranks are
\begin{equation}
\label{eq:chen-p13}
\theta_{k+2}(\h(P_{1,3};\k)) = (k+1)(k+3), \qquad k\ge 0,
\end{equation}
that is, the sequence $3, 8, 15, 24, 35, 48,\ldots$\,.

The assertion about the support---and with it the value of the Krull
dimension---in fact holds for every $n$: it is a restatement of
Proposition~\ref{prop:res-elliptic}.

\begin{proposition}
\label{prop:supp-b1-conf}
For every $n\ge2$, the support of $\B_1(H^{*}(\Conf(\E,n);\k))$ is given by
\eqref{eq:supp-b1-conf}, that is, it is the affine cone over the Segre variety
$\mathbb{P}^1\times\mathbb{P}^{n-2}$.  Consequently,
$\B_1(H^{*}(\Conf(\E,n);\k))$ has Krull dimension $n$.
\end{proposition}

\begin{proof}
Write $A=H^{*}(\Conf(\E,n);\k)$.  By Definition~\ref{def:hom-res},
$\Supp_S\B_1(A)=\RR_1(A)$.  Since $A$ is connected, $\RR^0(A)=\RR_0(A)=\{0\}$;
since also $\dim_\k A^j<\infty$ for $j\le1$, Theorem~\ref{thm:resonance-comparison}
applied with $q=1$ gives $\RR_0(A)\cup\RR_1(A)=\RR^0(A)\cup\RR^1(A)$, whence
$\RR_1(A)=\RR^1(A)$.  Now $\RR^1$ depends only on the cup product
$\bwedge^2A^1\to A^2$, so Proposition~\ref{prop:res-elliptic} identifies
$\RR^1(A)$ with the right-hand side of \eqref{eq:supp-b1-conf}.  The latter is
the affine cone over the Segre embedding of
$\mathbb{P}^1\times\mathbb{P}^{n-2}$, of dimension $1+(n-2)+1=n$.
\end{proof}

The generator count is likewise unconditional.

\begin{proposition}
\label{prop:b1-conf-generators}
For every $n\ge2$, the module $\B_1(H^{*}(\Conf(\E,n);\k))$ is minimally
generated by $\binom{n}{2}$ elements, all in degree $0$.
\end{proposition}

\begin{proof}
Write $A=H^{*}(\Conf(\E,n);\k)$, with its zero differential.  By
Proposition~\ref{prop:conf-elliptic-H2}, $A^1$ is finite-dimensional and the
cup product $\mu_A\colon\bwedge^2A^1\to A^2$ is surjective, its kernel being
the relation space $K$ of \eqref{eq:cauchy-relations}, of dimension
$\binom{n}{2}$.  Lemma~\ref{lem:b1-degree-zero} therefore applies, giving
$\B_1(A)/\m\B_1(A)\cong K^{\vee}$, of dimension $\binom{n}{2}$; by
Remark~\ref{rem:b1-degree-zero-shift}, this is the degree-$0$ piece in the
grading convention used throughout this subsection.
\end{proof}

The remaining assertions are open in general.  Equation
\eqref{eq:hilb-b1-conf} holds in internal degrees at most $9$, for every $n$
(Corollary~\ref{cor:b1-conf-reduction}), and has been checked degreewise in
degrees $0$ through $8$ for $2\le n\le16$; whether it, and the stated
projective dimension, persist for all $n$ and all degrees is open.

The equivariance of the situation suggests a sharper formulation, which
also explains the shape of the answer.  Write
$H^1(\Conf(\E,n);\k)=W\otimes_{\k} U$, where $W=H^1(\E;\k)\cong\k^2$ and
$U=\k^{n}$ records the points, so that $a_i=\alpha\otimes u_i$ and
$b_i=\beta\otimes u_i$ for a basis $\alpha,\beta$ of $W$, and let
$U_0\subseteq U$ be the sum-zero subspace.  As computed in the proof of
Proposition~\ref{prop:res-elliptic} by means of the Cauchy decomposition
\eqref{eq:cauchy-res},
\begin{equation}
\label{eq:cauchy-relations}
\ker\bigl(\cup\colon\bwedge^2 H^1\to H^2\bigr)
=\bwedge^2 W\otimes_{\k}\Sym^2 U_0 .
\end{equation}
Thus $\B_1(H^{*}(\Conf(\E,n);\k))$ is the Koszul module of
\cite{PS-crelle} associated to the pair $(V,K)$, where
$V=H_1(\Conf(\E,n);\k)$ and $K\subseteq\bwedge^2V$ is the annihilator of
the subspace \eqref{eq:cauchy-relations}, and it is
equivariant for $\GL(W)\times\GL(U_0)$.  The resonance variety $\RR(V,K)$
consists of the rank-one tensors $\lambda\otimes v$ with $\lambda\in W$ and
$v\in U_0$: the affine cone over the Segre variety
$\mathbb{P}^1\times\mathbb{P}^{n-2}$, of dimension $n$, as found above.  One
inclusion is immediate from \eqref{eq:cauchy-relations}: given
$a=\lambda\otimes v$ with $0\ne v\in U_0$, choose $\lambda'\in W$ with
$\lambda\wedge\lambda'\ne0$ and set $b=\lambda'\otimes v$; then
\[
a\wedge b=(\lambda\wedge\lambda')\otimes v^{2}
\in \bwedge^2 W\otimes_{\k}\Sym^2 U_0=\ker(\cup),
\]
and $a\wedge b\ne0$, so $a$ lies on a $2$-plane $P$ with
$\bwedge^2P\subseteq\ker(\cup)$.  The reverse inclusion is where the content
lies; it has been verified for $2\le n\le 8$.

The weight grading used here is an instance of
Definition~\ref{def:multihomog}: take $\Gamma=\Z^2\times\Z^{n-1}$, the group of
characters of the maximal torus of $\GL(W)\times\GL(U_0)$, and
$\lambda$ the total degree in the $W$-factor.  Since $H^1=W\otimes_{\k} U$ has
$W$-degree~$1$ throughout, the grading is normalized, so
Theorem~\ref{prop:multigraded-koszul} applies and the conjecture below is a
statement about a $\Gamma$-multigraded Hilbert function.  This is the coarser
grading promised in Remark~\ref{rem:multigrading-fails}, the fine
$\N^{2n}$-grading being unavailable here.

\begin{conjecture}
\label{conj:p1n-equivariant}
For all $n\ge2$ there is an isomorphism of
$\GL(W)\times\GL(U_0)$-representations
\begin{equation}
\label{eq:conj-p1n-equiv}
\B_1\bigl(H^{*}(\Conf(\E,n);\k)\bigr)_d \cong 
\bwedge^2 W^{\vee}\otimes_{\k}\Sym^d W^{\vee}\otimes_{\k}\Sym^{d+2}U_0^{\vee} ,
\qquad d\ge0 .
\end{equation}
Geometrically: $\B_1(H^{*}(\Conf(\E,n);\k))$ is the module of sections of
$\mathcal{O}(0,2)$ on the affine cone over
$\mathbb{P}^1\times\mathbb{P}^{n-2}$, a reflexive module of rank one.
\end{conjecture}

Grading $S$ by the torus of $\GL(W)\times\GL(U_0)$, so that a multidegree is
a triple $(m_1,m_2;\gamma)\in\Z^2\times\Z^{n-1}$, the predicted
representation has exactly one basis vector in each multidegree with
$m_1,m_2\ge1$ and $|\gamma|=m_1+m_2$, and none elsewhere.  Conjecture
\ref{conj:p1n-equivariant} is therefore equivalent to the statement that
every multigraded piece of $\B_1(H^{*}(\Conf(\E,n);\k))$ has dimension $0$
or $1$, the value being $1$ exactly on that set of multidegrees.  In this
form it is directly checkable.

\begin{remark}
\label{rem:equivariant-equivalent}
Conjecture~\ref{conj:p1n-equivariant} implies \eqref{eq:hilb-b1-conf}, by an
elementary computation.  Taking dimensions,
$\dim_\k\Sym^d W^{\vee}=d+1$ and
$\dim_\k\Sym^{d+2}U_0^{\vee}=\binom{d+n}{n-2}$, and setting
$A(t)=\sum_{d\ge0}\binom{d+n}{n-2}t^d
=t^{-2}\bigl[(1-t)^{-(n-1)}-1-(n-1)t\bigr]$, one has
\[
\sum_{d\ge0}(d+1)\binom{d+n}{n-2}t^d=\frac{d}{dt}\bigl(t\,A(t)\bigr)
=\frac{(1-t)^n-1+nt}{t^2(1-t)^n},
\]
which is the right-hand side of \eqref{eq:hilb-b1-conf}.  The converse
implication does not hold: the conjecture is a statement about the
$\GL(W)\times\GL(U_0)$-module structure, of which the Hilbert series records
only the dimensions.
\end{remark}

Two structural facts reduce the conjecture considerably.  The first singles
out the two linear forms already visible in the description of the support.

\begin{lemma}
\label{lem:b1-conf-annihilator}
The module $\B_1(H^{*}(\Conf(\E,n);\k))$ is annihilated by the linear forms
$\sum_{i} x_i$ and $\sum_{i} y_i$, that is, by the subspace
$W^{\vee}\otimes\k\pi\subseteq V$, where $\pi\in U^{\vee}$ is the functional
with $\ker\pi=U_0$.  Consequently, $\B_1(H^{*}(\Conf(\E,n);\k))$ is a module
over $\Sym\bigl(W^{\vee}\otimes_{\k} U_0^{\vee}\bigr)$, and its support is
contained in $W\otimes_{\k} U_0$.
\end{lemma}

\begin{proof}
Complete the basis $f_k=u_k-u_n$ $(1\le k\le n-1)$ of $U_0$ by $f_n=u_n$, and
set $a'_k=\alpha\otimes f_k$ and $b'_k=\beta\otimes f_k$ (these differ from
$a_k=\alpha\otimes u_k$ and $b_k=\beta\otimes u_k$ for $k<n$).  Let
$A_k,B_k\in V$ be the basis dual to $a'_k,b'_k$, so that $A_n$ and $B_n$ span
$W^{\vee}\otimes \k\pi$, and note that, by \eqref{eq:cauchy-relations}, the
subspace $\ker(\cup)$ has basis
$\rho_{kl}=a'_k\wedge b'_l+a'_l\wedge b'_k$ with $1\le k\le l\le n-1$.  Write
$g_{kl}$ for the dual basis of $(\ker\cup)^{\vee}$; by
Theorem~\ref{thm:Bpres}, these generate $\B_1$, with relations the images of
$\bwedge^3 V$ under $\partial_3^E$.

Since $\pi$ vanishes on $U_0$, while each $\rho_{kl}$ lies in
$\bwedge^2W\otimes_{\k}\Sym^2U_0$, we have
\begin{equation}
\label{eq:An-pairs-zero}
\langle A_n\wedge \eta,\,\rho_{kl}\rangle=
\langle B_n\wedge \eta,\,\rho_{kl}\rangle=0
\qquad\text{for all $\eta\in V$ and all $k\le l\le n-1$.}
\end{equation}
Hence the relation attached to the triple $A_k\wedge A_n\wedge B_l$ has
$\rho_{k'l'}$-component
\[
\big\langle \partial_3^E(A_k\wedge A_n\wedge B_l),\,\rho_{k'l'}\big\rangle
=-\big\langle A_k\wedge B_l,\,\rho_{k'l'}\big\rangle\, A_n ,
\]
the first and third terms of $\partial_3^E$ vanishing by
\eqref{eq:An-pairs-zero}.  Now
$\langle A_k\wedge B_l,\rho_{k'l'}\rangle=
\delta_{kk'}\delta_{ll'}+\delta_{kl'}\delta_{lk'}$, which vanishes unless
$(k',l')=(k,l)$, in which case it equals $1$ if $k<l$ and $2$ if $k=l$.  Being
nonzero in either case, this relation reads $A_n\cdot g_{kl}=0$.
Symmetrically, the triple $A_k\wedge B_l\wedge B_n$ gives
$B_n\cdot g_{kl}=0$.
\end{proof}

Two consequences are worth isolating.  First, in the presentation over
$S'=\Sym(W^{\vee}\otimes_{\k} U_0^{\vee})=\k[A_k,B_k\colon 1\le k\le n-1]$
obtained from Lemma~\ref{lem:b1-conf-annihilator}, every variable has
nonzero $\gamma$-degree, namely $e_k$, while the generator $g_{kl}$ has
$\gamma$-degree $e_k+e_l$; since $|e_k|=1$, every multidegree occurring in
$\B_1$ satisfies $|\gamma|=m_1+m_2$.  \emph{The half of
Conjecture~\ref{conj:p1n-equivariant} asserting that no other multidegrees
occur therefore holds for all $n$ and all $d$}, and only the
$1$-dimensionality of the predicted pieces remains at issue.  Second, the
remaining relations are exactly those attached to triples with all indices
at most $n-1$: the triples involving the index $n$ either give the relations
of Lemma~\ref{lem:b1-conf-annihilator} or, by \eqref{eq:An-pairs-zero},
vanish identically.  This locality yields a stabilization property.

\begin{lemma}
\label{lem:b1-conf-stab}
Let $(m_1,m_2;\gamma)$ be a multidegree, put $T=\{k\colon\gamma_k\ne0\}$ and
$r=\#T$, and let $\gamma|_T$ denote $\gamma$ with its zero entries deleted.
Then
\[
\dim_\k \B_1\bigl(H^{*}(\Conf(\E,n);\k)\bigr)_{(m_1,m_2;\,\gamma)}
=\dim_\k \B_1\bigl(H^{*}(\Conf(\E,r+1);\k)\bigr)_{(m_1,m_2;\,\gamma|_T)} .
\]
In particular, this dimension is independent of $n$.
\end{lemma}

\begin{proof}
Work with the presentation over $S'$ described above.  A monomial multiple
$\mu\, g_{kl}$ has $\gamma$-degree $e_k+e_l+\gamma(\mu)$, which is supported
in $T$ only if $k,l\in T$ and every variable dividing $\mu$ is indexed by
$T$.  Likewise, the relations attached to the triples
$A_i\wedge A_j\wedge B_m$ and $A_i\wedge B_m\wedge B_p$ have $\gamma$-degrees
$e_i+e_j+e_m$ and $e_i+e_m+e_p$, so only those with all indices in $T$
contribute in $\gamma$-degree $\gamma$.  Hence the graded piece in question
is computed by the sub-presentation indexed by $T$.  Finally, that
sub-presentation is defined by formulas depending only on the index set, and
is equivariant with respect to the symmetric group permuting
$f_1,\dots,f_{n-1}$; relabeling $T$ as $\{1,\dots,r\}$ therefore identifies
it with the presentation of $\B_1(H^{*}(\Conf(\E,r+1);\k))$.
\end{proof}

\begin{corollary}
\label{cor:b1-conf-reduction}
Fix $d\ge0$.  If Conjecture~\ref{conj:p1n-equivariant} holds in internal
degree $d$ for all $n\le d+3$, then it holds in internal degree $d$ for all
$n\ge2$.  Consequently, Conjecture~\ref{conj:p1n-equivariant} holds in
internal degrees $d\le9$, for every $n\ge 2$, and \eqref{eq:hilb-b1-conf}
holds in those degrees for every $n$.
\end{corollary}

\begin{proof}
A multidegree of internal degree $d$ satisfies $|\gamma|=d+2$, whence
$r=\#T\le d+2$ and $r+1\le d+3$; moreover, the multiplicity predicted by the
conjecture depends only on $m_1$, $m_2$, and $\gamma|_T$.  The first
assertion is thus immediate from Lemma~\ref{lem:b1-conf-stab}.  For the
second, the multidegree form of the conjecture has been verified by
Macaulay2, for each internal degree $d\le9$, at every multidegree with
$n\le d+3$; the statement about \eqref{eq:hilb-b1-conf} then
follows by Remark~\ref{rem:equivariant-equivalent}.
\end{proof}

Thus the conjecture is now a statement in $d$ alone: in each internal
degree it amounts to a finite verification, over the multidegrees whose
$\gamma$ has full support.

\begin{example}
\label{ex:b1-symbolic-power}
The geometric half of Conjecture~\ref{conj:p1n-equivariant} can be made
completely explicit.  Write $R_n=S/\Ann\bigl(\B_1(H^{*}(\Conf(\E,n);\k))\bigr)$,
whose spectrum is, by Proposition~\ref{prop:supp-b1-conf}, the affine cone over
$\mathbb{P}(W)\times\mathbb{P}(U_0)$, and let $P\subset R_n$ be the prime
generated by the images of $x_1,\dots,x_n$.  Then $P$ cuts out a fibre of the
projection to $\mathbb{P}(W)$, a divisor of class $\mathcal{O}(1,0)$.  Coning
kills the embedding class $\mathcal{O}(1,1)$, so
$\operatorname{Cl}(R_n)=(\Z\oplus\Z)/\Z(1,1)\cong\Z$, generated by $[P]$, and
$\mathcal{O}(1,0)\equiv\mathcal{O}(0,-1)$; the ideal sheaf of $P$ therefore
has class $\mathcal{O}(0,1)$.  Consequently the conjecture is equivalent to
the assertion that
\begin{equation}
\label{eq:b1-symbolic-square}
\B_1\bigl(H^{*}(\Conf(\E,n);\k)\bigr) \cong P^{(2)}(2),
\end{equation}
the second symbolic power of $P$, shifted so as to be generated in
degree~$0$.

The class is determined rather than merely matched: for $n\ge3$, the module
of sections of the class $c$ has Hilbert function
$(d+1)\binom{d+c+n-2}{n-2}$, which is strictly increasing in $c$, so $c=2$ is
the only class reproducing \eqref{eq:hilb-b1-conf}.  Part of this package has 
been verified by machine.  For $2\le n\le 10$, the annihilator of $\B_1$ is precisely 
the ideal of $2\times2$ minors described above, so that $R_n$ is the coordinate
ring of the affine cone over a Segre variety, and hence a normal domain;
this was confirmed directly for $3\le n\le6$, where $\Ann(\B_1)$ is prime,
$\dim R_n=n$, and $R_n$ is normal.  For $3\le n\le6$, moreover, $\B_1$ is
torsion-free of rank one over $R_n$---its only associated prime being
$\Ann(\B_1)$ itself---and it is reflexive.  The latter follows from the
former together with a Hilbert series identity: over a domain,
torsion-freeness makes the canonical map $\B_1\to\B_1^{\vee\vee}$ injective;
this map is degree-preserving, and the numerators of the Hilbert series of
the two sides agree over their common denominator, so the map is surjective
in every degree, hence an isomorphism.  Finally, \eqref{eq:b1-symbolic-square} 
has itself been verified for $3\le n\le6$, by
exhibiting a degree-zero element of $\Hom\bigl(\B_1,P^{(2)}(2)\bigr)$ with
vanishing kernel and cokernel; in this range $P^{(2)}=P^{2}$, so the
symbolic power is not yet doing any work.  Since
$H^0(\mathcal{O}(d,d+2))=\Sym^dW^{\vee}\otimes_{\k}\Sym^{d+2}U_0^{\vee}$, 
such an identification implies the character statement of
Conjecture~\ref{conj:p1n-equivariant} in the verified range, and a proof of
it for all $n$ would settle the conjecture.  With reflexivity in hand, it
would also follow from the divisor class computation above, that computation
being the one step here that has not been checked by machine.
\end{example}

\begin{remark}
\label{rem:equivariant-status}
Conjecture~\ref{conj:p1n-equivariant} being strictly stronger than
\eqref{eq:hilb-b1-conf}, by Remark~\ref{rem:equivariant-equivalent}, what is
verified at present is the numerical consequence, for $2\le n\le16$, together 
with the conjecture itself, as an isomorphism of representations, in each internal 
degree $d\le9$: it has been confirmed by machine for $n\le d+3$ in each such
degree, and hence, by Corollary~\ref{cor:b1-conf-reduction}, for
\emph{every}\/ $n$.  The case $d=0$ holds for every $n$: by
Lemma~\ref{lem:b1-degree-zero}, the degree-zero part of the Koszul module is
$K^{\vee}$, dual to \eqref{eq:cauchy-relations} and hence isomorphic
to $\bwedge^2W^{\vee}\otimes_{\k}\Sym^2U_0^{\vee}$.  The duals and the
one-dimensional twist in the conjecture are dictated by this case.
\end{remark}

\begin{remark}
\label{rem:general-elliptic}
Neither \eqref{eq:cauchy-relations} nor the description of the resonance
uses the braid combinatorics.  Let $\A$ be any elliptic arrangement in
$\E^{n}$, with hyperplanes given by characters
$\chi_H\colon\E^{n}\to\E$ and corresponding vectors $v_H\in U$.  The Bibby
differential is $d e_H=(\alpha\otimes v_H)\wedge(\beta\otimes v_H)
=(\alpha\wedge\beta)\otimes v_H^{2}$, so
\[
\ker(\cup)=\bwedge^2 W\otimes_{\k} C_{\A}, \qquad
C_{\A}=\operatorname{span}\{v_H^{2}: H\in\A\}\subseteq\Sym^2 U ,
\]
and the same computation as above identifies the resonance as the set of
\emph{rank-one} tensors
\[
\bigl\{\lambda\otimes v : \lambda\in W,\ v\in X_{\A}\bigr\}, \qquad
X_{\A}=\{v\in U : v^{2}\in C_{\A}\},
\]
that is, the affine cone over
$\mathbb{P}(W)\times\bigl(\nu(\mathbb{P}(U))\cap\mathbb{P}(C_{\A})\bigr)$,
where $\nu$ is the quadratic Veronese.  Its dimension is
$1+\dim X_{\A}$---not $2\dim X_{\A}$, as the full tensor product
$W\otimes_{\k} X_{\A}$ would give.  The braid case is $C_{\A}=\Sym^2 U_0$, for
which $X_{\A}=U_0$ and the dimension is $1+(n-1)=n$.

Computations for eleven elliptic arrangements suggest a dichotomy governed
by $C_\A$.  When $C_{\A}=\Sym^2 U'$ for a subspace $U'\subseteq U$, the
Veronese section is $\mathbb{P}(U')$ and the resonance is irreducible and
determinantal, as in the braid case.  When instead $C_{\A}$ is spanned by
the squares of independent vectors, $X_{\A}$ degenerates to a finite union
of lines, the resonance becomes a union of $2$-planes, and $\B_1$ behaves
as for a hyperplane arrangement: linear components, no minors.  Deleted
braid arrangements interpolate between the two---removing one relation
$r_{ij}$ leaves the components determinantal but disconnects the variety,
while removing two collapses it to linear components.  It would be
interesting to know which subspaces $C_{\A}\subseteq\Sym^2 U$ arise from
elliptic arrangements, and to prove the dichotomy just described.
\end{remark}

By Corollary~\ref{cor:chen-tor}, the holonomy Chen ranks are given by
the linear strand of the minimal free $E$-resolution of $H^*(\Conf(\E,n);\k)$
over the exterior algebra $E = \bigwedge H^1(\Conf(\E,n);\k)$:
\begin{equation}
\label{eq:chen-p1n}
\theta_{k+2}(\h(P_{1,n};\k)) = \dim_\k \Tor^{E}_{k-1}(H^*(\Conf(\E,n);\k),\k)_k.
\end{equation}
In practice, this linear strand computation via $\Tor^E$ is often the
most efficient approach to Chen ranks. A conceptual identification 
of the linear strand of this resolution---for instance via the BGG 
correspondence applied to the universal enveloping algebra of 
$\h(P_{1,n};\k)$, which is the Koszul dual of $H^*(\Conf(\E,n);\k)$ 
as a quadratic $E$-module---would yield a proof of the Hilbert 
series formula above.

%%%%%%%%%%%%%%%%%%%%%%
\subsection{Comparing the Koszul modules of $A(n)$ and of its cohomology}
\label{subsec:spectral-weight}
%%%%%%%%%%%%%%%%%%%%%%

We now return to the Bibby model $(A(n),d)$, with its positive weights: the
generators $a_i,b_i$ have weight $1$ and the generators $e_{ij}$ have
weight~$2$, so that $d$ preserves weights, as in
\eqref{eq:morgan-weight}.  Since $H^1(A(n))$ is spanned by the classes
$[a_i],[b_i]$, all of weight~$1$, it is weight-homogeneous, and
Corollary~\ref{cor:weight-splitting}\eqref{ws2} identifies the weight spectral
sequence of $K_\bullet(A(n))$ with the $\m$-adic Koszul spectral sequence,
refined by the total-weight splitting.

In particular, the differential $d_1$ on the $E^1$-page is, by
Theorem~\ref{thm:weight-ss}, the Koszul differential of $H^*(A(n))$; it is
determined by the cohomology algebra alone and carries no information about the
homotopy type.  What detects non-formality here is instead the comparison
between the two Koszul modules, which we now make explicit.

\begin{proposition}
\label{prop:xi-An-noninjective}
For every $n\ge3$, the epimorphism
\[
\B_1\bigl(H^*(A(n))\bigr)\longsurj \B_1\bigl(A(n)\bigr)
\]
of \eqref{eq:b1ha-b1a} has nonzero kernel.  Consequently
$\ker\Xi_{A(n)}\not\subseteq\h(H^*(A(n)))''$; in particular, the
comparison map $\Xi_{A(n)}$ of Lemma~\ref{lem:holo-to-cohomology} is not
injective.
\end{proposition}

\begin{proof}
By Corollary~\ref{cor:BA-n-consequences}\eqref{bn1},
\[
\Hilb\bigl(\B_1(A(n)),t\bigr)=\frac{\binom{n}{2}}{(1-t)^{2}},
\]
so the coefficients of this series grow linearly in the degree.  On the other
hand, $H^{*}(A(n))\cong H^{*}(\Conf(\E,n);\k)$, and by
Proposition~\ref{prop:supp-b1-conf} the module
$\B_1(H^{*}(\Conf(\E,n);\k))$ has Krull dimension $n$; its Hilbert function
therefore grows like a polynomial of degree $n-1\ge2$.  The two graded
modules have different Hilbert functions, so the surjection of $S$-modules
$\B_1(H^*(A(n)))\surj\B_1(A(n))$ cannot be an isomorphism.  The remaining
assertions follow from Lemma~\ref{lem:ker-B-Xi}, applied to
$\Xi_{A(n)}$.
\end{proof}

Combining Proposition~\ref{prop:xi-An-noninjective} with
Theorem~\ref{thm:koszulss-1formal-obstruction} gives a second proof of
Corollary~\ref{cor:p1n-nonformal}: the pure elliptic braid group $P_{1,n}$ is
not $1$-formal for $n\ge3$, the obstruction being the strict drop
$\theta_k(P_{1,n})<\bar\theta_k(P_{1,n})$ in Chen ranks.

\begin{remark}
\label{rem:graphic-conf-formality}
By contrast, consider the \emph{graphic configuration spaces} 
$\Conf(\E,\Gamma)=\{(z_1,\dots,z_n)\in\E^n \mid 
z_i\neq z_j \text{ for all } ij\in\Gamma\}$ associated to a finite 
simple graph $\Gamma$ on $n$ vertices~\cite{BMPP}. For a triangle-free 
graph $\Gamma$, the group $\pi_1(\Conf(\E,\Gamma))$ is $1$-formal
by~\cite{BMPP}, so Theorem~\ref{thm:koszulss-1formal-obstruction} forces
$\B_1(H^*(A(\Gamma)))\cong\B_1(A(\Gamma))$---even though the differential
$d_1$ of the weight spectral sequence is nonzero, exactly as for
$\Gamma=K_n$.  This is the sharpest illustration that non-vanishing of $d_1$ is
not an obstruction to $1$-formality: it is a cohomology-level invariant,
insensitive to the distinction between these two families.  What separates
them is the kernel of $\B(\Xi_{A(\Gamma)})$, equivalently the position of
$\ker\Xi_{A(\Gamma)}$ relative to the second derived subalgebra
$\h(H^*(A(\Gamma)))''$ \textup{(}Lemma~\ref{lem:ker-B-Xi}\textup{)}; 
identifying it combinatorially, in terms of the triangles of $\Gamma$, would 
give a Koszul-theoretic proof of the $1$-formality criterion of~\cite{BMPP}.
\end{remark}

\begin{remark}
\label{rem:elliptic-arr}
The Bibby model applies more generally to unimodular elliptic
arrangements, where the underlying Orlik--Solomon relations are
combinatorially determined, but the additional relations encoding the
defining equations of the arrangement are not.
This raises the question of to what extent resonance varieties and
formality properties of elliptic arrangement complements are governed
by combinatorial data alone; see \cite[Prob.~8.3]{Su-indam}.
Related questions were formulated in \cite{Bi16} and remain largely 
open outside special classes of arrangements.
\end{remark}

%%%%%%%%%%%%%%%%%%%%%%%%%%%%%%
\subsection{Comparison with the pure braid group $P_n$}
\label{subsec:pn-p1n-comparison}
%%%%%%%%%%%%%%%%%%%%%%%%%%%%%%

With the invariants of $P_{1,n}$ now in hand, we close the section with a
direct comparison between the pure elliptic braid group
$P_{1,n} = \pi_1(\Conf(\E,n))$ and its classical counterpart, the pure braid
group $P_n = \pi_1(\Conf(\C,n))$.
Both are fundamental groups of configuration spaces, but on 
two different one-dimensional complex algebraic groups:
the additive group $\C$ and the elliptic curve $\E$.
The contrast between their LCS ranks, Chen ranks, and formality properties
illustrates how the Koszul module framework extends naturally from the
$1$-formal to the non-formal setting.

\medskip
\noindent\textbf{The classical pure braid group and the quadratic truncation
of $A(n)$.}
The complement $M_n = \Conf(\C,n)$ of the braid arrangement $\A_{n-1}$
(the type-$A_{n-1}$ reflection arrangement) is aspherical, via the tower
\begin{equation}
\label{eq:fibration-tower}
\begin{tikzcd}[column sep=18pt]
\Conf(\C,n) \ar[r] & \Conf(\C,n-1)  \ar[r] &  \cdots  \ar[r] &  \Conf(\C,2) \cong \C^*,
\end{tikzcd}
\end{equation}
whose fibers $\C\setminus\{k\text{ points}\}\simeq\bigvee^k S^1$ are
wedges of circles, so each stage is aspherical and $M_n$ is a classifying
space $K(P_n,1)$.  Moreover, the monodromy of \eqref{eq:fibration-tower} is
trivial in homology, by Falk--Randell~\cite{FR85}, so the Leray--Hirsch
theorem gives the product formula
\begin{equation}
\label{eq:hilb-os-an}
\Hilb(\OS(\A_{n-1}),t) = (1+t)(1+2t)\cdots(1+(n-1)t)
\end{equation}
for $H^*(M_n;\k)=\OS(\A_{n-1})$, and the Hirsch lemma implies that $M_n$ is
formal; consequently $P_n$ is $1$-formal, and $\OS(\A_{n-1})$ is Koszul
\cite{SY97}.  The existence of such a tower is what makes $\A_{n-1}$ a
fiber-type arrangement, and fiber-type arrangements always have Koszul
Orlik--Solomon algebras \cite{SY97}; this is the mechanism behind both the
factorization above and the Koszulity just cited.

The group $P_n$ splits as an iterated semidirect product of free groups,
\begin{equation}
\label{eq:pn-split}
P_n \cong F_{n-1} \rtimes F_{n-2} \rtimes \cdots \rtimes F_1,
\end{equation}
with trivial action on abelianizations.  Following Falk--Randell
\cite{FR85}, an iterated semidirect product of free groups in which each
factor acts trivially on the homology of the preceding ones---as in
\eqref{eq:pn-split}---is called an \emph{almost-direct product}\/ of free
groups.  Falk--Randell \cite{FR85} deduced
from \eqref{eq:pn-split} that the LCS ranks of $P_n$ satisfy
\begin{equation}
\label{eq:lcspn}
\phi_k(P_n) = \sum_{j=1}^{n-1} \phi_k(F_j),
\end{equation}
a formula holding more generally for fiber-type arrangements.
Kohno \cite{Ko85} established the equivalent product formula
\begin{equation}
\label{eq:lcspn-hilb}
\prod_{k\ge 1}(1-t^k)^{\phi_k(P_n)} = \Hilb(\OS(\A_{n-1}),-t),
\end{equation}
the equivalence with \eqref{eq:lcspn} following from PBW.  The
Koszul-theoretic explanation of \eqref{eq:lcspn-hilb} came later, and is due
to Shelton--Yuzvinsky \cite{SY97}: they identified the holonomy algebra of an
arrangement with the quadratic dual of the quadratic closure of its
Orlik--Solomon algebra, proved that $\OS(\A)$ is Koszul whenever $\A$ is
supersolvable, and thereby recovered \eqref{eq:lcspn-hilb} as the Hilbert
series identity relating a Koszul algebra to its dual.

The Chen ranks of $P_n$ were computed by Cohen--Suciu \cite{CS95}:
\begin{equation}
\label{eq:chen-pn}
\theta_1(P_n) = \binom{n}{2}, \quad
\theta_2(P_n) = \binom{n}{3}, \quad
\theta_k(P_n) = (k-1)\binom{n+1}{4} \quad \text{for } k\ge 3.
\end{equation}
This was the first nontrivial computation of Chen ranks for an arrangement
group, and provided the original motivation for the Chen ranks conjecture
\cite{Su-conm01}.
The formula \eqref{eq:chen-pn} is consistent with Theorem~\ref{thm:chen-afrs}:
$\RR^1(P_n;\k)$ has $\binom{n}{3}+\binom{n}{4}$ irreducible components,
all of dimension $2$, and separability of $\RR^1(P_n;\k)$---established
in \cite{AFRS25} for all arrangements whose flats of multiplicity $\ge 2$
have multiplicity at most $3$, a condition satisfied by $\A_{n-1}$---%
verifies the hypotheses of Theorem~\ref{thm:chen-afrs}, yielding
\[
\theta_k(P_n) = \Bigl[\tbinom{n}{3}+\tbinom{n}{4}\Bigr]\cdot\theta_k(F_2)
= \tbinom{n+1}{4}(k-1) \quad\text{for }k\ge n-1.
\]

The Bibby model $(A(n),d)$ provides a structural link to $P_n$ via its
quadratic truncation.

\begin{proposition}
\label{prop:bibby-qkoszul}
For the Bibby model $(A(n),d)$,
\begin{equation}
\label{eq:qql-An}
\qql(A(n),d) \cong
\bwedge(a_1,b_1,\dots,a_n,b_n)\big/
\bigl\langle (a_i-a_j)(b_i-b_j) \mid 1\le i<j\le n\bigr\rangle ,
\end{equation}
the quadratic closure of $H^{\le 2}(\Conf(\E,n);\k)$.  Equivalently, in the
coordinates $p_i=a_i-a_n$, $q_i=b_i-b_n$ of \eqref{eq:conf-elliptic-split},
\[
\qql(A(n),d)\cong
\frac{\bwedge(p_1,q_1,\dots,p_{n-1},q_{n-1})}{\langle p_iq_i,\
p_iq_j+p_jq_i\rangle}\ \otimes_{\k}\ \bwedge(a_n,b_n),
\]
where the free exterior factor has \emph{two} generators.
\end{proposition}

\begin{proof}
Write $V=H^1(A(n))=\langle a_i,b_i\rangle$ and $W=\langle e_{ij}\rangle$, so
that $A(n)^1=V\oplus W$.  By Proposition~\ref{prop:Uh-cdga}, $\h(A(n))$ is
presented on $A(n)_1$ by the relations
$\mu_{A(n)}^{\vee}(\xi)+d_{A(n)}^{\vee}(\xi)=0$, $\xi\in A(n)_2$, whose
linear and quadratic parts are coupled.  Since $d_{A(n)}e_{ij}
=(a_i-a_j)(b_i-b_j)$, for each pair $i<j$ there is a $\xi$ with
$d_{A(n)}^{\vee}(\xi)=e_{ij}^{\vee}$, and the corresponding relation reads
\[
e_{ij}^{\vee}=-\mu_{A(n)}^{\vee}(\xi)\in[\h(A(n)),\h(A(n))] .
\]
Hence each generator $e_{ij}^{\vee}$ is redundant: $\h(A(n))$ is generated by
$V^{\vee}$ alone, and eliminating the $e_{ij}^{\vee}$ converts these
$\binom{n}{2}$ linear relations into the quadratic relations
$(a_i-a_j)(b_i-b_j)$ among the $a_i^{\vee},b_i^{\vee}$.  The remaining
relations, those $\xi$ with $d_{A(n)}^{\vee}(\xi)=0$, are the ones dual to the
cup product $\bwedge^2 V\to H^2(A(n))$.  Passing to the associated quadratic
algebra therefore gives \eqref{eq:qql-An}, which is by construction the
quadratic closure of $H^{\le 2}(\Conf(\E,n);\k)$.  The second presentation is
the splitting of Proposition~\ref{prop:conf-elliptic-H2}\eqref{qn4}.
\end{proof}

The based factor in the second presentation of
Proposition~\ref{prop:bibby-qkoszul} is not merely an auxiliary bookkeeping
device: it is itself a quadratic closure, of the space obtained by
removing the identity of $\E$.

\begin{proposition}
\label{prop:based-factor-conf}
Write
\[
B(n-1)\coloneqq
\frac{\bwedge(p_1,q_1,\dots,p_{n-1},q_{n-1})}{\langle p_iq_i,\
p_iq_j+p_jq_i\rangle}
\]
for the based factor of Proposition~\textup{\ref{prop:bibby-qkoszul}}. Then
\[
B(n-1)\cong \qql\bigl(H^{*}(\Conf(\E^{*},n-1);\k)\bigr) ,
\]
the quadratic closure of the \textup{(}full, untruncated\textup{)}
cohomology ring of the configuration space of $n-1$ points on $\E^{*}$.
\end{proposition}

\begin{proof}
By \eqref{eq:conf-elliptic-topsplit} and the K\"{u}nneth formula, there is an
isomorphism of graded $\k$-algebras
\[
H^{*}(\Conf(\E,n);\k)\cong H^{*}(\E;\k)\otimes_{\k}
H^{*}(\Conf(\E^{*},n-1);\k) .
\]
Since $\E$ is a curve, $H^{*}(\E;\k)$ vanishes above
degree $2$, and $H^{*}(\E;\k)=\bwedge(a_n,b_n)$ is already quadratic (indeed
relation-free). Passing to quadratic closures on both sides, and using that
the quadratic closure of a tensor product of connected graded algebras
generated in degree $1$ is the tensor product of the quadratic closures,
\[
\qql\bigl(H^{*}(\Conf(\E,n);\k)\bigr)
\cong \bwedge(a_n,b_n)\otimes_{\k}
\qql\bigl(H^{*}(\Conf(\E^{*},n-1);\k)\bigr) .
\]
The left side is $Q(n)=\qql(A(n),d)$ of Proposition~\ref{prop:bibby-qkoszul},
since $\qql$ depends only on the algebra in degrees $\le 2$ and
$H^{\le2}(\Conf(\E,n);\k)=A(n)^{\le2}$
(Proposition~\ref{prop:conf-elliptic-H2}). Comparing with the second
presentation of Proposition~\ref{prop:bibby-qkoszul}, whose two tensor
factors are $B(n-1)$ and $\bwedge(a_n,b_n)$, identifies the two occurrences
of $\bwedge(a_n,b_n)$ and gives the stated isomorphism.
\end{proof}

\begin{remark}
\label{rem:qql-An-no-OS}
The Orlik--Solomon algebra $\OS(\A_{n-1})$ does \emph{not} appear in
\eqref{eq:qql-An}.  The relations $e_{ij}e_{ik}-e_{ij}e_{jk}+e_{ik}e_{jk}=0$
are indeed present in $A(n)$, but the $e_{ij}^{\vee}$ are not generators of
$\qql(A(n),d)$, so they contribute no tensor factor.  Two independent checks
confirm this.  First, \eqref{eq:lcs-holo-hilb} forces
$\dim\qql(A(n),d)_1=\phi_1(P_{1,n})=b_1(\Conf(\E,n))=2n$, whereas
$\dim A(n)^1=2n+\binom{n}{2}$; so no algebra generated by $A(n)^1$ can occur.
Second, $A(n)$ carries the mixed quadratic relations $(a_i-a_j)e_{ij}
=(b_i-b_j)e_{ij}=0$, which couple the $a,b$ generators to the $e$ generators
and obstruct any tensor decomposition separating them.
\end{remark}

\begin{question}
\label{q:An-qkoszul}
Is the quadratic algebra \eqref{eq:qql-An} Koszul, i.e.\ is $(A(n),d)$
quadratically Koszul? Equivalently, by
Proposition~\ref{prop:based-factor-conf}, is
$H^{*}(\Conf(\E^{*},n-1);\k)$ Koszul? Solving
$\prod_{k}(1-t^k)^{\phi_k}=\Hilb(\qql(A(n),d),-t)$ returns positive integers
$\phi_k$ for $n\le 5$ and $k\le 6$, which is consistent with Koszulness but
does not establish it.
\end{question}

As just noted, Question~\ref{q:An-qkoszul} is equivalent to the
Koszulness of the cohomology ring of a configuration space on a
once-punctured torus---a question with a direct route to an affirmative
answer in the literature, though one we have not independently verified in
full. The fundamental group
$\pi_1(\Conf(\E^{*},m))$ is built, via the Fadell--Neuwirth
fibrations $\Conf(\E^{*},k)\to\Conf(\E^{*},k-1)$ with
fiber $\E\setminus\{k\text{ points}\}\simeq \bigvee_k S^1$, from an iterated
tower of free groups, exactly as for $\Conf(\E,n)$
(cf.\ \eqref{eq:fibration-tower} and the discussion following it); the
difference is that here the point at $0$ is fixed throughout, so the
translation action of $\E$ on itself---identified above as the source of
the nontrivial monodromy that defeats the analogous tower for
$\Conf(\E,n)$---is no longer available.
Maassarani \cite{Maassarani20}, extending fibration-theoretic techniques of
Halperin to orbit configuration spaces of punctured surfaces, shows that
$\Conf(\Sigma\setminus F,m)$ is a \emph{rational $K(\pi,1)$} space (its
minimal and $1$-minimal Sullivan models agree) whenever the surface
$\Sigma\setminus F$ is not $S^2$---in particular for $\Sigma=\E$, $F=\{0\}$
\cite[Prop.~3.2]{Maassarani20}. By Theorem~\ref{thm:py}\eqref{py1} of
Papadima--Yuzvinsky, a rational $K(\pi,1)$ space with finite Betti numbers
has Koszul rational cohomology; combined with Maassarani's result, this
would give an affirmative answer to Question~\ref{q:An-qkoszul} outright.
A second, independent route: Cohen \cite{Cohen-comp10} proves that the
cohomology ring of any almost-direct product of free groups is
Koszul, via an explicit quadratic Gr\"{o}bner basis; whether
$\pi_1(\Conf(\E^{*},m))$ has this stronger structure (as opposed
to merely the nilpotent monodromy that Maassarani's argument uses) is a
question we have not checked.

\medskip
\noindent\textbf{Growth of the LCS and Chen ranks.}
Granting an affirmative answer to Question~\ref{q:An-qkoszul},
Theorem~\ref{thm:koszul-holo} applied to $(A(n),d)$ gives
\begin{equation}
\label{eq:lcs-p1n-hilb}
\prod_{k\ge 1}(1-t^k)^{\phi_k(P_{1,n})}
= \Hilb(\qql(A(n),d),-t)
= (1-t)^2\cdot\Hilb\Bigl(
\tfrac{\bwedge(p_\bullet,q_\bullet)}{\langle p_iq_i,\,p_iq_j+p_jq_i\rangle},
-t\Bigr),
\end{equation}
the factor $(1-t)^2$ accounting for the exterior algebra on $a_n,b_n$.
Since $\qql(A(n),d)$ is the quadratic closure of $H^{\le2}(\Conf(\E,n);\k)$,
the exponents produced by \eqref{eq:lcs-p1n-hilb} are the lower central series
ranks of the holonomy Lie algebra $\h\bigl(H^{*}(\Conf(\E,n);\k)\bigr)$; write
them $\phi_k^{\mathrm{hol}}$.  Because $\h(H^{*}(X))\surj\gr(\pi_1(X))\otimes\Q$ is an
isomorphism in degrees $\le 2$ but only a surjection in general, we have
\begin{equation}
\label{eq:phi-hol-bound}
\phi_k(P_{1,n})\le \phi_k^{\mathrm{hol}},
\qquad\text{with equality for } k\le 2 .
\end{equation}
Equality for all $k$ would follow from $1$-formality of $P_{1,n}$, which
\emph{fails} here---this is the content of Corollary~\ref{cor:p1n-nonformal}
above.  So only the
first two columns of the following table are established values of
$\phi_k(P_{1,n})$; the remaining entries are upper bounds, and are in addition
conditional on Question~\ref{q:An-qkoszul}.

\begin{center}
\begin{tabular}{c|rrrrrr}
$k$ & $1$ & $2$ & $3$ & $4$ & $5$ & $6$ \\\hline
$\phi_k^{\mathrm{hol}}(n=3)$ & $6$ & $3$ & $8$  & $18$  & $48$  & $116$ \\
$\phi_k(P_3)$        & $3$ & $1$ & $2$  & $3$   & $6$   & $9$ \\[4pt]
$\phi_k^{\mathrm{hol}}(n=4)$ & $8$ & $6$ & $20$ & $60$  & $204$ & $670$ \\
$\phi_k(P_4)$        & $6$ & $4$ & $10$ & $21$  & $54$  & $125$ \\
\end{tabular}
\end{center}

\noindent
Unconditionally, then, $\phi_1(P_{1,n})=2n$ and
$\phi_2(P_{1,n})=\binom{n}{2}$, the latter because
$\dim\qql(A(n),d)_2=\binom{2n}{2}-\binom{n}{2}$ is the rank of the cup
product.  Comparing with $\phi_1(P_n)=\binom{n}{2}$ and
$\phi_2(P_n)=\binom{n}{3}$, which follow from \eqref{eq:lcspn-hilb}
\textup{(}and should not be confused with the Chen ranks
$\theta_1,\theta_2$ of \eqref{eq:chen-pn}, which happen to agree with them
in these two degrees\textup{)}, the first two ranks of
$P_{1,n}$ exceed those of $P_n$ for $n\le 4$, agree at $n=5$ in degree~$1$,
and are exceeded for $n\ge 6$; so no uniform domination holds, and the
comparison of growth rates must be made through the radii of convergence of
the two series rather than termwise.
As a calibration, at $n=2$ the same computation returns $4,1,2,3,6,9$, which
agrees termwise with
$P_{1,2}=\pi_1(\E\times\E^{*})\cong\Z^2\times F_2$; here
$P_{1,2}$ \emph{is} $1$-formal, so \eqref{eq:phi-hol-bound} is an equality.

Turning from lower central series to Chen ranks, the picture changes:
despite having different LCS ranks, both $P_n$ and $P_{1,n}$ have
Chen ranks that grow linearly in $k$, but with different slopes.

\begin{proposition}
\label{prop:chen-comparison}
For $n\ge 3$ and $k\ge 3$:
\begin{enumerate}[itemsep=2pt]
\item \label{cc1} $\theta_k(P_n) \ne \theta_k(P_{1,n})$.
\item \label{cc2} Both grow linearly:
$\theta_k(P_n) \sim \binom{n+1}{4}(k-1)$ and
$\theta_k(P_{1,n}) \sim \binom{n}{2}(k-1)$.
\item \label{cc3} $\theta_k(P_{1,n}) > \theta_k(P_n)$ for $n\le 4$,
and $\theta_k(P_n) > \theta_k(P_{1,n})$ for $n\ge 5$.
\end{enumerate}
\end{proposition}

\begin{proof}
Claims \eqref{cc1} and \eqref{cc2} follow from \eqref{eq:chen-pn}
and Corollary~\ref{cor:BA-n-consequences}\eqref{bn2}.
The slopes $\binom{n+1}{4}$ and $\binom{n}{2}$ coincide only when
$(n+1)(n-2)=12$, i.e., $n^2-n-14=0$, which has no integer solution.
For \eqref{cc3}: $\binom{n+1}{4}>\binom{n}{2}$ iff $(n+1)(n-2)>12$,
i.e., $n^2-n>14$, holding for $n\ge 5$; one checks directly that
$\binom{5}{4}=5<6=\binom{4}{2}$ and $\binom{4}{4}=1<3=\binom{3}{2}$.
\end{proof}

\medskip
\noindent\textbf{Why the product and additive formulas fail, and what to
conjecture instead.}
The product formula \eqref{eq:hilb-os-an} and the fibration tower
\eqref{eq:fibration-tower} are intimately linked: the fibers
$\C\setminus\{k\text{ points}\}$ have trivial monodromy in homology,
which is what makes the Leray--Hirsch theorem applicable and forces
formality.
By contrast, $\Conf(\E,n)$ admits no analogous fibration tower
whose fibers have trivial monodromy action on their homology. 
The natural projection $\Conf(\E,n)\to\Conf(\E,n-1)$
has fiber $\E\setminus\{n-1\text{ points}\}$, a punctured torus
whose fundamental group is the free group $F_n$ of rank $n$
(by the Euler characteristic count $2g+p-1$, with $g=1$ and $p=n-1$
punctures).
The monodromy of this fibration---the action of
$\pi_1(\Conf(\E,n-1))$ on $H_1(\E\setminus\{n-1\text{ pts}\};\k)$
---is nontrivial, coming from the translation action of $\E$ on itself.
This is not special to the elliptic case: for any orientable surface of
positive genus, a presentation of the surface braid group exhibits mixed
relations coupling the genus generators to the puncture generators of the
fiber---for instance, \cite[Lem.~3.2]{Bellingeri-jalg04} shows the genus
generator $a_r$ conjugates a fiber generator $\tau_1$ to a nontrivial word
involving the remaining fiber generators---witnessing exactly this failure
of triviality on homology.
This nontrivial monodromy obstructs the Leray--Hirsch product formula
and is the geometric source of non-formality of $P_{1,n}$.  Neither
factorization survives: the tower is not an almost-direct product, so
$P_{1,n}$ admits no splitting of the form \eqref{eq:pn-split} with
trivial action on abelianizations, and consequently the Falk--Randell
additivity \eqref{eq:lcspn} of lower central series ranks is unavailable
as well.  What fails, then, is not any factorization of the space
itself---$\Conf(\C,n)$ is no product either---but the multiplicative
formula \eqref{eq:hilb-os-an} for the Hilbert series and the additive
formula \eqref{eq:lcspn} for the LCS ranks, both of which the linear case
inherits from its tower.
Concretely, the obstruction is encoded in the Bibby differential:
$d(e_{ij}) = (a_i-a_j)(b_i-b_j)$ records the monodromy contribution
of the diagonal $\Delta_{ij}$, mixing the two $H^1(\E;\k)$-directions.
In the linear case, the analogous differential is $d(e_{ij}) = a_i - a_j$
(a single cohomology class), reflecting trivial monodromy;
the elliptic case requires the product $(a_i-a_j)(b_i-b_j)$,
reflecting the two-dimensional nature of $H^1(\E;\k)$.

The non-formality of $\Conf(\E,n)$ fits into a broader framework:
for smooth complex algebraic varieties whose MHS on cohomology is pure,
formality follows from Morgan's theorem~\cite{Mo}; conversely,
Dupont~\cite{Du16} shows that mixed weight structures can obstruct
formality.
For hyperplane arrangement complements, the MHS on $H^*$ is pure of
Tate type, consistent with their formality.
For $\Conf(\E,n)$, the compactness of $\E$ introduces
weight-$1$ classes from $H^1(\E;\k)$, making the MHS genuinely mixed
and obstructing formality; the Bibby differential $d(e_{ij}) =
(a_i-a_j)(b_i-b_j)$ is the explicit $\cdga$ manifestation of this
mixed structure.

Consequently, the $P_{1,n}$ family sits at the intersection of two threads.
On one hand, by Proposition~\ref{prop:bibby-qkoszul} and
Theorem~\ref{thm:koszul-holo}, the LCS ranks of $P_{1,n}$ are governed not by
the Orlik--Solomon algebra $\OS(\A_{n-1})$ that controls $P_n$, but by the
quadratic closure of $H^{\le2}(\Conf(\E,n);\k)$, in which the classes
$e_{ij}$ have been traded for the relations $(a_i-a_j)(b_i-b_j)$; the
comparison with $P_n$ is therefore one of contrast rather than of a common
Orlik--Solomon factor.
On the other hand, the Chen ranks of $P_{1,n}$ follow the same
template as the Chen ranks conjecture for $1$-formal groups: with
$h_2 = \binom{n}{2}$ linear components of dimension $2$ in
$\RR^1(A(n))$, one has $\theta_k(P_{1,n}) = \binom{n}{2}(k-1)
= h_2\cdot\theta_k(F_2)$ for all $k\ge 2$.

As discussed in Section~\ref{subsec:chen-ranks-conj}, this suggests an
extension of the Chen ranks conjecture to non-$1$-formal groups
admitting positive-weight models with linear resonance.
The relevant condition replacing strong isotropicity should be
formulated at the level of $\B_1(A)$: here it is the direct sum
decomposition $\B_1(A(n))\cong\bigoplus_{i<j}S/I_{ij}$ of
Theorem~\ref{thm:BA-n}, which follows from the block-diagonal structure
of the Bibby
differential and yields the Chen ranks formula directly.
The family $\{P_{1,n}\}_{n\ge 2}$ provides a worked non-formal
example of such an extension, and a natural test case for any
proposed generalization; see Question~\ref{q:chen-nonformal}.

%%%%%%%%%%%%%%%%%%%%%%%%
\section{Sasakian manifolds and Seifert fibered spaces}
\label{sect:sasaki}
%%%%%%%%%%%%%%%%%%%%%%%%

The Hirsch-extension techniques developed earlier specialize cleanly to
manifolds carrying a geometric circle action. The most natural such class
is that of Sasakian manifolds, where the Reeb field generates a canonical
$S^1$-action whose orbit space is a K\"{a}hler orbifold. We show that a
single algebraic mechanism---the positive-weight Hirsch extension of a
bounded formal $\cdga$, analyzed in Theorem~\ref{thm:partial-formality-hirsch}
and Theorem~\ref{thm:massey-hirsch-length}---recovers and sharpens the
partial-formality theorems for Sasakian manifolds, produces an explicit
central-extension presentation of their Malcev Lie algebras, and
specializes verbatim to orientable Seifert fibered spaces. In particular,
Sasakian and Seifert geometry are placed under one roof, and the formality
obstructions are pinned to a precise cohomological degree determined by the
dimension of the base.

%%%%%%%%%%%%%%%%%%%%%%%%
\subsection{Sasakian manifolds and partial formality}
\label{subsec:sasakian-formality}
%%%%%%%%%%%%%%%%%%%%%%%%

Let $M^{2n+1}$ be a compact Sasakian manifold. A general reference for
Sasakian geometry is the monograph of Boyer and Galicki~\cite{BG}. By a
result of Ornea and Verbitsky~\cite{OV}, any compact Sasakian manifold
admits a quasi-regular Sasakian structure, and we assume this throughout.

In the quasi-regular case, a basic structural theorem in Sasakian geometry
(\cite[Thm.~7.1.3]{BG}) ensures that $M$ carries an almost free $S^1$-action,
generated by the Reeb vector field. The orbit space $B=M/S^1$ is a compact
K\"{a}hler orbifold of real dimension $2n$, with K\"{a}hler class
$h\in H^2(B;\Q)$, satisfying the Hard Lefschetz property
\begin{equation}
\label{eq:hardlef}
H^{n-k}(B;\Q) \longisom  H^{n+k}(B;\Q),
\qquad 1\le k\le n;
\end{equation}
see \cite[Prop.~7.2.2 and Thm.~7.2.9]{BG}. The $S^1$-action gives a
fibration $S^1 \to M \xrightarrow{\eta} B$, classified by an Euler class
$e(\eta)\in H^2(B;\k)$.

\begin{theorem}[\cite{BBFMT}]
\label{thm:kahler-orbifold-formal}
The orbifold de Rham algebra of a compact K\"{a}hler orbifold is formal.
\end{theorem}

\begin{remark}
\label{rem:orbifold-formality-gap}
Two points of terminology should be kept in mind in reading
Theorem~\ref{thm:kahler-orbifold-formal}.

First, formality of an orbifold $X$ is defined in \cite[Def.~3.8]{BBFMT} as
formality of the minimal model of
$\bigl(\Omega^{*}_{\mathrm{orb}}(X),d\bigr)$, whose cohomology is
$H^{*}(X;\R)$ by \cite[Prop.~3.9]{BBFMT}.  Since an isomorphism of cohomology
algebras does not make two $\cdgas$ weakly equivalent, this does not identify
that minimal model with the one built from piecewise-polynomial forms on $X$;
the point is noted in \cite[\S6.1]{PS-imrn19}.  The theorem is correct as
stated, in the form in which it is proved; what is not available is the
stronger assertion that a compact K\"{a}hler orbifold is formal as a space.  We
do not use it below in any case: Tievsky's model \eqref{eq:tievsky-model}
already has zero differential on $H^{*}(B;\k)$, which is all the arguments of
this section require of the base.

Second, the $s$-formality of \cite[Def.~2.2]{BBFMT}, following
Fern\'{a}ndez--Mu\~{n}oz, is a condition on the generators of the 
minimal model in degrees $\le s$, and differs from the $q$-formality of
Section~\ref{subsec:formality}; the two are in effect off by one, and for a
nilpotent Lie algebra $\g$ the $\cdga$ $\CE(\g)$ may be $1$-formal in the
former sense without being so in ours.  The discrepancy was identified by
M\u{a}cinic \cite{Mc10}; see also \cite[Rem.~3.8]{Kasuya}.  Consequently the
results of \cite{BBFMT} on non-formal Sasakian manifolds are not directly
comparable with Corollary~\ref{cor:sasakian-partial-formal}.  Kotschick and
Placini \cite[Rem.~17]{KotP} trace to the same source the assertion in
\cite[\S3.1]{BFMT} that the higher Heisenberg groups fail to be $1$-formal: in
the sense used here those groups \emph{are} $1$-formal, as Carlson and Toledo
had noted, and Proposition~\ref{prop:cup-surj-implies-1formal}\ref{cs1}
recovers this independently, since $\ker\mu=0$ for $\h(n)$ with $n\ge2$ by
the Lefschetz property, while $\h(1)$ has $\bwedge^3W=0$ and fails the
criterion \textup{(}Remark~\ref{rem:1formal-cup-product}\textup{)}.
\end{remark}

Tievsky~\cite{Tievsky} constructed an explicit Sullivan model for $M$ in
terms of the cohomology algebra of the base. Let $\omega\in H^2(B;\k)$ be
the class of the K\"{a}hler form. Then $M$ admits a finite $\cdga$ model
\begin{equation}
\label{eq:tievsky-model}
\big(H^{*}(B;\k)\otimes_{\k} \bwedge(c), d\big),
\qquad \deg(c)=1, \quad d(c)=\omega.
\end{equation}
This model is minimal whenever $H^{*}(B;\k)$ is generated in degree one; in
general it is a finite Sullivan model, sufficient for all formality
considerations below.  One caveat: \cite{Tievsky} produces
\eqref{eq:tievsky-model} as a model over $\R$, and it does not follow
directly that it is a model over $\Q$, nor over an arbitrary field of
characteristic zero; see the discussion in \cite[\S6.1]{PS-imrn19}.  The
partial-formality conclusions below are nonetheless valid over any such
field, by \cite[Thm.~6.7]{PS-imrn19}, which combines the descent of partial
formality with the observation that the identification of $H^1$ furnished by
Tievsky's model preserves $\Q$-structures. From a homotopy-theoretic standpoint
\eqref{eq:tievsky-model} is a Hirsch extension of $(H^{*}(B;\k),0)$:
if $\omega=0$ the model splits as a tensor product, while if $\omega\neq0$
the differential introduces a single quadratic relation reflecting the
Euler class. Assigning $\wt(c)=1$ when $\omega=0$ and $\wt(c)=2$ when
$\omega\neq0$ gives a system of positive Hirsch weights compatible with $d$.

\begin{corollary}[\cite{PS-imrn19}]
\label{cor:sasakian-partial-formal}
Let $M$ be a compact Sasakian manifold of dimension $2n+1$. Then $M$ is
$(n-1)$-formal.
\end{corollary}

\begin{proof}
By Tievsky's model \eqref{eq:tievsky-model}, the manifold $M$ has a finite
model which is a $1$-step Hirsch extension
\[
\bigl(H^{*}(B;\k),0\bigr)\otimes_{e}\bwedge(c)
\]
of the cohomology algebra of a compact K\"{a}hler orbifold $B$ of real
dimension $2n$, taken with zero differential.  In particular the base of the
extension is formal by construction, so
Theorem~\ref{thm:kahler-orbifold-formal} is not needed here
(Remark~\ref{rem:orbifold-formality-gap}); and $H^{*}(B;\k)$ is
concentrated in degrees $\le 2n$ and admits positive weights.  The K\"{a}hler
class $\omega$ satisfies the Hard Lefschetz property, so multiplication by
$\omega$ is injective on $H^{i}(B)$ for $i\le n-1$; that is, $\{\omega\}$ is
an $(n-1)$-regular sequence.  The assertion is therefore the case
$e=\omega$, $q=n-1$ of Theorem~\ref{thm:partial-formality-hirsch}.
\end{proof}

\begin{remark}
\label{rem:sasaki-optimal}
The bound in Corollary~\ref{cor:sasakian-partial-formal} is optimal, but it
is not the case that every compact Sasakian manifold with nonzero Euler
class fails to be $n$-formal.  Optimality is witnessed by a specific family:
the $(2n+1)$-dimensional Heisenberg nilmanifold $M(n)$ of
Remark~\ref{rem:sasaki-Heisenberg} is Sasakian and is not $n$-formal
\cite[Ex.~6.8]{PS-imrn19}.  At the other extreme, the sphere $S^{2n+1}$,
viewed as the total space of the Hopf fibration over $\mathbb{CP}^n$, is
Sasakian with Euler class the nonzero hyperplane class, yet it is formal:
here $H^{*}(\mathbb{CP}^n)=\k[h]/(h^{n+1})$ and the model
$\k[h]/(h^{n+1})\otimes_{h}\bwedge(c)$ has cohomology
$\bwedge(x_{2n+1})$, free on a single odd generator, hence formal.  What
separates the two cases is not the vanishing of $e$ but the failure of
$n$-regularity together with a surviving Massey product; see
Proposition~\ref{prop:triple-massey-hirsch}.
\end{remark}

The geometric statement of Corollary~\ref{cor:sasakian-partial-formal} is
due to Papadima--Suciu~\cite[Thm.~6.6]{PS-imrn19}; what the
present treatment adds is the mechanism. Partial formality is forced not by
Sasakian geometry as such, but by the single structural feature that the
minimal model is a positive-weight Hirsch extension of a bounded formal
$\cdga$ along a partially regular class. The same mechanism applies to any
space whose model has this form, in particular to the Seifert fibered spaces
of \S\ref{subsec:seifert-hirsch}, with no input from Riemannian geometry.

\begin{remark}
\label{rem:kasuya-gap}
For $n\ge2$, Kasuya~\cite[Thm.~1.1]{Kasuya} asserts that a compact Sasakian
manifold of dimension $2n+1$ is $1$-formal over $\R$.  As explained in
\cite[\S6.2]{PS-imrn19}, the argument given there has a gap: what is
established is that $H^2$ of the $1$-minimal model is generated by products
of degree-one classes, and the passage from this to $1$-formality invokes
\cite[Lem.~3.17]{ABCKT}, which is false---see
\cite[Ex.~4.5 and Rem.~4.6]{Mc10}.  The conclusion itself is correct, and
in the stronger form of Corollary~\ref{cor:sasakian-partial-formal}, by the
independent argument of \cite{PS-imrn19}.  The same caution applies to
Proposition~\ref{prop:cup-surj-implies-1formal}\ref{cs2}, whose hypothesis
is of exactly this shape; see Remark~\ref{rem:1formal-cup-product}.
\end{remark}

\begin{remark}
\label{rem:n1-enot0}
For the Massey criterion what matters is whether $[e]$ is a \emph{single}\/
product, rather than whether it lies in the span of such.  Accordingly, three
cases arise.

If $[e]=\alpha\beta$ for some $\alpha,\beta\in H^1(B)$, then non-$1$-formality
follows from the nonvanishing of the triple Massey product
$\langle\alpha,\beta,\beta\rangle$, by
Corollary~\ref{cor:e-decomposable-massey}.  This is automatic when
$\dim_\k H^2(B)=1$ and $[e]\ne0$, as for a Seifert fibration over a surface,
or for a $3$-dimensional Sasakian manifold.

If $[e]\notin\im(\mu_B)$, no obstruction arises from either mechanism: the
holonomy ranks of $A$ and of $H^{*}(A)$ agree
(Remark~\ref{rem:holo-no-quadratic-obstruction}), and $A$ may be formal.

In the intermediate case, $[e]$ lies in $\im(\mu_B)$ without being a single
product.  This requires $\dim_\k H^1(B)\ge4$: for $\dim_\k H^1(B)\le3$ every
element of $\bwedge^2H^1(B)$ has rank at most~$2$, hence is decomposable, as is
its image in $H^2(B)$.  Here neither mechanism applies; the smallest such base
is $B=H^{*}(T^4)$ with $[e]=a_1a_2+a_3a_4$, and in that example passing to
degree~$3$ does not help either, since $\h(A)\cong\h(H^{*}(A))$ in every degree
by \eqref{eq:holo-degree3-decisive}.

For $n>1$, Sasakian manifolds fall in one of the last two cases: $A$ is
$(n-1)$-formal, hence $1$-formal, and the obstruction to formality lies in
degree~$n$.
\end{remark}

\begin{remark}
\label{rem:sasaki-Heisenberg}
The $(2n+1)$-dimensional Heisenberg nilmanifold $M(n)$ is the quotient of a
nilpotent Lie group $\mathfrak{H}(1,n)$ by a discrete, cocompact subgroup. It
is a circle bundle over the complex torus $\E^{\times n}$ with Euler class
equal to the K\"{a}hler form $\omega$, and is therefore Sasakian. Conversely,
the Tievsky model was used in \cite{CDMY, Kasuya} to show that a compact
$(2n+1)$-dimensional nilmanifold admits a Sasakian structure if and only if
it is a quotient of $\mathfrak{H}(1,n)$ by a discrete, cocompact subgroup; in
particular, all Sasakian nilmanifolds share the same rational homotopy type.
\end{remark}

%%%%%%%%%%%%%%%%%%%%%%%%%
\subsection{Malcev Lie algebras of Sasakian manifolds}
\label{subsec:malcev-sasaki}
%%%%%%%%%%%%%%%%%%%%%%%%%

Let $M^{2n+1}$ be a compact Sasakian manifold, with Reeb circle action
$S^1 \to M \xrightarrow{\eta} B$, where $B$ is a compact K\"{a}hler orbifold
of real dimension $2n$. As recalled in
\S\ref{subsec:sasakian-formality}, $M$ admits the finite-type $\cdga$ model
$A = \big(H^{*}(B;\k)\otimes_{\k} \bwedge(c), d\big)$ with $d(c)=\omega$.
Since $B$ is a compact K\"{a}hler orbifold it is formal, and $H^{*}(B;\k)$
admits positive weights; equipping $c$ with weight $1$ when $\omega=0$ and
weight $2$ when $\omega\neq0$ makes $d$ weight-homogeneous, so $A$ is a
positive-weight Hirsch extension. By
Proposition~\ref{prop:positive-weights-ff}, the group $\pi_1(M)$ is
\emph{filtered-formal}: its Malcev Lie algebra $\m(\pi_1(M);\k)$ is the
degree completion of a graded Lie algebra determined by the quadratic part
of the differential.

To describe the presentation, let $V=H^1(B;\k)$ and let $\h(B)$ denote the
holonomy Lie algebra of the base,
$\h(B)=\Lie(V^{\vee})/\langle \im(\mu^{\vee})\rangle$, where
$\mu\colon \bwedge^2 V \to H^2(B;\k)$ is the cup product and
$\mu^{\vee}\colon H_2(B;\k)\to \bwedge^2 V^{\vee}$ its Kronecker dual. The
differential $d(c)=\omega$ induces a functional $\omega^{\vee}\colon
H_2(B;\k)\to\k$, corresponding dually to a central extension of Lie
algebras: $\m(\pi_1(M);\k)$ is the degree completion of the graded Lie
algebra
\begin{equation}
\label{eq:sasaki-malcev}
L(M)=
\Lie\big(V^{\vee}\big)\oplus \k z \;\Big/\;
\bigl\langle\,
\mu^{\vee}(h)-\omega^{\vee}(h)\,z \ :\ h\in H_2(B;\k)
\,\bigr\rangle,
\qquad [z,\cdot\,]=0 ,
\end{equation}
with $\deg(z)=2$ if $\omega\ne0$, and $\deg(z)=1$ if $\omega=0$. In other
words, $\m(\pi_1(M);\k)$ is obtained from the holonomy Lie algebra of the
base by adjoining a central generator $z$: each defining relation
$\mu^{\vee}(h)=0$ of $\h(B)$ is replaced by $\mu^{\vee}(h)=\omega^{\vee}(h)z$.
When $\omega=0$ the extension is trivial:
$L(M)=\h(B)\oplus\k z$ is a literal direct product, with $z$ in degree $1$,
matching $M=B\times S^1$.

The presentation \eqref{eq:sasaki-malcev} localizes the possible failure of
graded-formality in Sasakian groups to a single mechanism: all Sasakian
groups are filtered-formal, and $\m(\pi_1(M);\k)$ can fail to be graded only
through the central generator $z$, which is a genuine extra generator of
$\h(\pi_1(M))$ precisely when $\omega\in\im(\mu_B)$
(Proposition~\ref{prop:holo-hirsch}\eqref{eh1}).  Whether graded-formality
actually fails in that case is finer; see Remark~\ref{rem:n1-enot0}.  For the
Heisenberg nilmanifold $M(n)$ the base is $B=\E^{\times n}$, so $\h(B)$ is the
free abelian Lie algebra on $2n$ generators and $\omega$ is a nondegenerate
alternating form; \eqref{eq:sasaki-malcev} then recovers the standard
presentation of the Heisenberg Lie algebra, with $z$ the central generator.

Two points are worth isolating. First, \eqref{eq:sasaki-malcev} is an
\emph{explicit, closed-form} presentation of the Malcev Lie algebra of an
arbitrary Sasakian group as a central extension of the holonomy Lie algebra
of its K\"{a}hler base---reducing a Malcev computation to cup-product data on
$B$. Second, it explains, rather than merely re-deriving, the higher-Massey
vanishing of~\cite{BFMT}, as follows.

\begin{remark}
\label{rem:sas-massey}
Biswas--Fern\'{a}ndez--Mu\~{n}oz--Tralle~\cite{BFMT} show that all Massey
products of length $\ge 4$ vanish on any compact Sasakian manifold, even
though such manifolds need not be formal. The present viewpoint isolates the
mechanism: the Tievsky model is a $1$-step Hirsch extension of a formal
$\cdga$, so Theorem~\ref{thm:massey-hirsch-length} applies, and every Massey
product of length $\ge 4$ of pure classes---in particular, every such product
of degree-one classes---contains $0$.  The obstruction is simply one of
$\bwedge(c)$-degree: a defining system of span $s$ would have
$\bwedge(c)$-degree $s$, and only the degrees $0$ and $1$ are available.
Equivalently, in the language of the Koszul spectral sequence, where $d_r$
for $r\ge2$ encodes length-$(r{+}2)$ Massey products, no differential beyond
$d_1$ can be nontrivial on such classes.  In
this sense Sasakian manifolds exhibit the strongest failure of formality
compatible with a $1$-step Hirsch extension: triple Massey products may
occur, but no higher obstructions are allowed.
\end{remark}

%%%%%%%%%%%%%%%%%%%%%%%%
\subsection{Nilpotent Sasakian groups of small rank}
\label{subsec:sasaki-nilpotent}
%%%%%%%%%%%%%%%%%%%%%%%%

Chen \cite[Prob.~2]{XChen} asks which nilpotent groups arise as fundamental
groups of compact Sasakian manifolds. The presentation \eqref{eq:sasaki-malcev}
answers this below the Carlson--Toledo threshold, in exact parallel with the
nilpotent K\"{a}hler case \cite{CT95}.

\begin{proposition}
\label{prop:sasaki-nilpotent}
Let $\Gamma$ be a finitely generated nilpotent Sasakian group with
$b_1(\Gamma)\le 7$, and let $\m=\m(\Gamma;\k)$ be its Malcev Lie algebra over
a field $\k$ of characteristic $0$.  Then either $\omega=0$ and $\Gamma$ is
abelian, with $\m(\Gamma;\k)\cong\k^{b_1(\Gamma)}$; or $\omega\ne0$ and $\m$
is Heisenberg-by-abelian,
\[
\m(\Gamma;\k)\cong \h(s)\times\k^{t},
\qquad s\ge1,\quad 2s+t=b_1(\Gamma)\ \text{even},\ \ t\ \text{even},
\]
where $\h(s)$ is the $(2s+1)$-dimensional Heisenberg Lie algebra of
Section~\ref{subsec:heisenberg}. In particular $\dim_\k[\m,\m]\le 1$ and
$\Gamma$ is at most $2$-step nilpotent.
\end{proposition}

\begin{proof}
After a quasi-regular deformation of the Sasakian structure---which leaves
$\pi_1$ unchanged \cite{BG}---we may write $\Gamma=\pi_1(M)$ with $M\to B$ a
Seifert fibration over a compact K\"{a}hler orbifold $B$, giving a central
extension
\[
\begin{tikzcd}[column sep=20pt]
0\ar[r]& \k z \ar[r]& \m \ar[r]& \h(B) \ar[r]& 0
\end{tikzcd}
\]
in which $\h(B)$ is the holonomy Lie algebra of the nilpotent K\"{a}hler orbifold
group $\pi_1^{\mathrm{orb}}(B)$, by \eqref{eq:sasaki-malcev}.  If $\omega=0$,
$z$ has degree~$1$ and the extension splits, $\m=\h(B)\oplus\k z$, so
$b_1(\h(B))=b_1(\Gamma)-1\le6$; if $\omega\ne0$, $z$ has degree~$2$, so
$b_1(\h(B))=b_1(\Gamma)\le 7$.  Either way $b_1(\h(B))\le7$, and by the
Carlson--Toledo rank bound \cite[Cor.~4.5]{CT95}, whose Hodge-theoretic
argument carries over to compact K\"{a}hler orbifolds via their Hodge
decomposition and hard Lefschetz property \cite{BBFMT} (after Baily), a
nilpotent K\"{a}hler orbifold group with $b_1\le 7$ is almost abelian; hence
$\h(B)$ is abelian.

If $\omega=0$, then $\m=\h(B)\oplus\k z\cong\k^{b_1(\Gamma)-1}\oplus\k z
=\k^{b_1(\Gamma)}$, and $\Gamma$ is abelian.  If $\omega\ne0$, then
$[\m,\m]\subseteq \k z$, and the induced alternating form on the abelian
$\h(B)$, being nonzero, has positive rank~$s\ge1$ and splits off its radical
as the abelian factor $\k^{t}$, giving $\m\cong\h(s)\times\k^{t}$.
Evenness of $b_1(\Gamma)$ is \cite[Prop.~1.1]{XChen}, whence
$t=b_1(\Gamma)-2s$ is even.
\end{proof}

\begin{corollary}
\label{cor:sasaki-nilpotent-obstruction}
A finitely generated torsion-free nilpotent group $\Gamma$ with $b_1\le 7$
whose Malcev Lie algebra is not Heisenberg-by-abelian---in particular, any
such $\Gamma$ with $\dim_\k[\m,\m]\ge 2$ or of nilpotency class $\ge 3$---is
not a Sasakian group.
\end{corollary}

\begin{example}
\label{ex:f42-not-sasaki}
Let $N_{4,2}$ be the free $2$-step nilpotent group on $4$ generators, with
$\m(N_{4,2};\k)=\ff_{4,2}$. Then $b_1=4\le 7$, so $N_{4,2}$ survives the Betti
obstruction \cite[Prop.~1.1]{XChen}; but $\dim_\k[\ff_{4,2},\ff_{4,2}]
=\binom{4}{2}=6$, so $N_{4,2}$ is not Sasakian by
Corollary~\ref{cor:sasaki-nilpotent-obstruction}.
Theorem~\ref{thm:chen-growth-2step} makes the obstruction quantitative:
$\theta_2(N_{4,2})=6$, whereas every nilpotent
Sasakian group with $b_1=4$ has $\theta_2\le 1$ (namely $\h(2)$, with
eventually vanishing Chen ranks, or $\h(1)\times\k^2$).
\end{example}

\begin{remark}
\label{rem:sasaki-nilpotent-sharp}
The bound $b_1\le 7$ is the Carlson--Toledo rank-$8$ threshold
\cite[Cor.~4.5]{CT95}: a nonabelian nilpotent K\"{a}hler group has $b_1\ge 8$,
and Campana \cite{Campana} realizes nilpotent K\"{a}hler groups as lattices in
the Heisenberg groups of dimension $\ge 9$. Beyond the threshold a nonabelian
K\"{a}hler base may occur, and the classification of nilpotent Sasakian groups
is open---mirroring the status of the nilpotent K\"{a}hler problem, unsettled in
nilpotency class $\ge 3$ \cite{Campana}.

The nilpotent case sits inside a broader hierarchy. A solvable Sasakian group
is virtually nilpotent---by \cite[Thm.~1.4]{XChen} via Campana's K\"{a}hler
orbifold theory, and, for polycyclic groups, by Kasuya's independent
basic-cohomology argument \cite[Cor.~1.3]{Kasuya16}. At the manifold level,
Proposition~\ref{prop:sasaki-nilpotent} is the infinitesimal shadow of the
geometric classifications of Sasakian nilmanifolds \cite{CDMY} and
solvmanifolds \cite[Cor.~1.4]{Kasuya16} as finite quotients of Heisenberg
nilmanifolds; the present statement instead constrains the Malcev Lie algebra
of \emph{any} nilpotent Sasakian group of small rank, with no assumption that
the manifold be a solvmanifold.
\end{remark}

\begin{remark}
\label{rem:sasaki-dimension}
The obstruction of Corollary~\ref{cor:sasaki-nilpotent-obstruction} is
independent of dimension: it rules a group out of $\mathcal{S}_m$, the class
of fundamental groups of compact Sasakian manifolds of dimension $m$, for
every $m$ at once.  That is worth recording, because the classes
$\mathcal{S}_m$ genuinely depend on $m$.  Kotschick and Placini
\cite{KotP} show that each $\mathcal{S}_m$ contains groups belonging to no
$\mathcal{S}_{m'}$ with $m'>m$, that $\mathcal{S}_m\subset\mathcal{S}_{m-2}$
for $m\ge7$, and hence that $\mathcal{S}_5$ is the union of all
$\mathcal{S}_{m}$ with $m\ge5$; they also show that $\mathcal{S}_m$ contains
every projective group for $m\ge5$, and that the class of Sasakian groups,
unlike that of K\"{a}hler groups, is not closed under direct products.  Thus
$N_{4,2}$ of Example~\ref{ex:f42-not-sasaki} lies in no $\mathcal{S}_m$
whatsoever.

Their dimension~$3$ versus dimension~$\ge5$ dichotomy is the geometric
counterpart of the formality dichotomy in
Remark~\ref{rem:n1-enot0}: no infinite group in $\mathcal{S}_3$ occurs in
any $\mathcal{S}_m$ with $m\ge5$, precisely because the members of
$\mathcal{S}_3$ carry nontrivial triple Massey products on $H^1$ while
Corollary~\ref{cor:sasakian-partial-formal} forces $1$-formality as soon as
$n\ge2$.
\end{remark}

%%%%%%%%%%%%%%%%%%%%%%%%
\subsection{Seifert fibered spaces}
\label{subsec:seifert-hirsch}
%%%%%%%%%%%%%%%%%%%%%%%%

Let $M$ be an orientable, closed Seifert fibered space with orientable base.
Such a $3$-manifold admits an effective $S^1$-action with orbit space an
orientable surface $\Sigma_g$ and finitely many exceptional orbits; the
fibration $\eta\colon M\to \Sigma_g$ is classified by integers
$(\alpha_i,\beta_i)$ and an obstruction $b=b(\eta)$. The fundamental group
$G_\eta=\pi_1(M)$ has the presentation
\begin{equation}
\label{eq:pi1-seifert}
\langle x_1,y_1,\dots,x_g,y_g,z_1,\dots,z_s,h \mid
h \text{ central},\;
[x_1,y_1]\cdots[x_g,y_g]z_1\cdots z_s=h^b,\;
z_i^{\alpha_i}h^{\beta_i}=1\rangle .
\end{equation}

Assume $g>0$. The minimal model of $M$ is the Hirsch extension
\begin{equation}
\label{eq:minmod-seifert}
\cM(M)\simeq \cM(\Sigma_g)\otimes_{e}\bwedge(c),
\end{equation}
where $e\in H^2(\Sigma_g;\k)$ represents the Euler class. Since $\Sigma_g$
is formal and generated in degree one, this is a $1$-step positive-weight
Hirsch extension of a formal $\cdga$. By
Proposition~\ref{prop:positive-weights-ff}, $G_\eta$ is filtered-formal over
$\k$, and $\m(G_\eta)$ is the degree completion of a graded Lie algebra with
generators in degree one and a single central generator in degree two when
$e(\eta)\neq0$. This presentation, given explicitly in
\cite[Thm.~10.6]{SW-forum}, coincides with the specialization of the
Sasakian Malcev Lie algebra construction \eqref{eq:sasaki-malcev} to the case
$B=\Sigma_g$.

Thus orientable Seifert fibered spaces with $g>0$ are not a separate case but
the planar specialization ($n=1$) of the Sasakian picture: the same
positive-weight Hirsch extension governs partial formality
(Corollary~\ref{cor:sasakian-partial-formal} with $n=1$ gives only
$0$-formality, and when $e\neq0$ the group is in fact not $1$-formal, by
Remark~\ref{rem:n1-enot0}: here $\dim_\k H^2(\Sigma_g)=1$, so $[e]$ is a
single product and Corollary~\ref{cor:e-decomposable-massey} applies) and
the same central-extension formula
\eqref{eq:sasaki-malcev} computes the Malcev Lie algebra. Sasakian and
Seifert geometry are therefore two instances of a single algebraic mechanism.

%%%%%%%%%%%%%%%%%%%%%%%%%%%%%%%%%%%

\bibliographystyle{amsplain}

\end{document}